\documentclass[a4paper,11pt]{amsart}

\usepackage{mathptmx}
\usepackage{helvet}
\usepackage{courier}
\usepackage{graphicx}
\usepackage{multicol}
\usepackage{footmisc}
\usepackage[T1]{fontenc}
\usepackage{euscript}
\usepackage{rotating}
\usepackage{epic}
\usepackage{tikz-cd}
\usepackage{multicol}
\usepackage{tikz}
\tikzstyle{every picture}+=[remember picture]
\usepackage{comment}
\usepackage{amsfonts,amsthm,amsmath}
\usepackage{amssymb}
\usepackage{amscd}
\usepackage{enumerate}
\usepackage{enumitem}
\usepackage{soul}
\usepackage{array}
\usepackage{array,tabularx}
\usepackage[all,cmtip]{xy}
\usepackage{pgfplots}
\usetikzlibrary{arrows}
\usepackage{lscape}
\usepackage{float}

\usepackage{pict2e}
\usepackage[dvipsnames]{xcolor}
\usepackage{quiver}

\numberwithin{equation}{section}
\numberwithin{equation}{subsection}

\theoremstyle{plain}

\newtheorem{theorem}[equation]{Theorem}
\newtheorem{lemma}[equation]{Lemma}
\newtheorem{proposition}[equation]{Proposition}

\newtheorem{corollary}[equation]{Corollary}

\newtheorem{question}[equation]{Question}

\newtheorem{thm}[equation]{Theorem}
\newtheorem{cor}[equation]{Corollary}

\newtheorem{prop}[equation]{Proposition}

\newtheorem{claim}[equation]{Claim}

\newtheorem{obs}[equation]{Observation}

\theoremstyle{definition}

\newtheorem{example}[equation]{Example}
\newtheorem{nonexample}[equation]{Nonexample}
\newtheorem{examples}[equation]{Examples}
\newtheorem{remark}[equation]{Remark}

\newtheorem{define}[equation]{Definition}
\newtheorem{problem}[equation]{Problem}

\newtheorem{nota}[equation]{Notation}

\newtheorem{ass}[equation]{Assumption}
\newtheorem{conj}[equation]{Conjecture}

\newtheorem{bekezdes}[equation]{}

\newcommand{\bC}{{\mathbb C}}

\newcommand{\bQ}{{\mathbb Q}}
\newcommand{\bR}{{\mathbb R}}

\newcommand{\bZ}{{\mathbb Z}}
\newcommand{\bN}{{\mathbb N}}

\newcommand{\cV}{{\mathcal V}}

\newcommand{\cI}{{\mathcal I}}

\newcommand{\cO}{{\mathcal O}}
\newcommand{\calO}{{\mathcal O}}
\newcommand{\cF}{{\mathcal F}}

\newcommand{\cL}{{\mathcal L}}

\newcommand{\tX}{\widetilde{X}}

\newcommand{\C}{{\calc}}

\newcommand{\ix}{\index}

\def\blfootnote{\xdef\@thefnmark{}\@footnotetext}

\newcommand{\co}{\cO}

\newcommand{\frakv}{\mathfrak{v}}

\newcommand{\bH}{{\mathbb H}}

\newcommand{\calF}{{\mathcal F}}

\newcommand{\hh}{\mathfrak{h}}

\newcommand{\calv}{\mathcal{V}}

\newcommand{\sw}{{\mathfrak{sw}}}

\newcommand{\Z}{\mathbb{Z}}
\newcommand{\Q}{\mathbb{Q}}
\newcommand{\R}{\mathbb{R}}
\newcommand{\g}{\gamma}

\begingroup\makeatletter\ifx\SetFigFontNFSS\undefined%
\gdef\SetFigFontNFSS#1#2#3#4#5{%
  \reset@font\fontsize{#1}{#2pt}%
  \fontfamily{#3}\fontseries{#4}\fontshape{#5}%
  \selectfont}%
\fi\endgroup%

\newcommand{\frh}{\mathfrak{h}}

\newcommand{\frv}{\mathfrak{v}}

\newcommand{\calt}{{\mathcal T}}

\def\C{\mathbb C}
\def\Q{\mathbb Q}
\def\R{\mathbb R}
\def\bH{\mathbb H}

\def\Z{\mathbb Z}
\def\N{\mathbb N}
\def\I{\mathbb I}

\newcommand{\cdstack}[2]{\begin{matrix}#1\\#2\end{matrix}}
\newcommand{\dvdstack}[3]{\begin{matrix}#1\\#2\\#3\end{matrix}}
\definecolor{rred}{RGB}{214, 92, 92}
\tikzset{
    invisible/.style={opacity=0},
    visible on/.style={alt={#1{}{invisible}}},
    alt/.code args={<#1>#2#3}{%
      \alt<#1>{\pgfkeysalso{#2}}{\pgfkeysalso{#3}}%
  }
}

\newenvironment{sizeddisplay}[1]
 {\par\nopagebreak#1\noindent\ignorespaces}
 {\nopagebreak\ignorespacesafterend}

\pgfplotsset{compat=1.18}

\usepackage{hyperref}

\author{Andr\'as N\'emethi}
\thanks{The authors are partially supported by  ``\'Elvonal (Frontier)'' Grant KKP 144148}
\address{Alfr\'ed R\'enyi Institute of Math.,
    Re\'altanoda utca 13-15, H-1053, Budapest, Hungary \newline
    \hspace*{3mm} ELTE - Univ. of Budapest, Dept. of Geo.,
    P\'azm\'any P\'eter s\'et\'any 1/A, 1117, Budapest, Hungary \newline \hspace*{3mm}
    BBU - Babe\c{s}-Bolyai Univ., Str, M. Kog\u{a}lniceanu 1, 400084 Cluj-Napoca, Romania
    \newline \hspace*{3mm}
    BCAM - Basque Center for Applied Math.,
    Mazarredo, 14 E48009 Bilbao, Basque Country, Spain}
\email{nemethi.andras@renyi.hu }

\author{Gergő Schefler}
\address{HUN-REN Alfr\'ed R\'enyi Institute of Mathematics,
Re\'altanoda utca 13-15, H-1053, Budapest, Hungary \newline
 \hspace*{4mm} ELTE - University of Budapest, Dept. of Geometry, Budapest, Hungary }
\email{schefler@renyi.hu}

\vspace*{-1cm}

\title{Lattice homology
of integrally closed submodules and Artin algebras }

\begin{document}

\begin{abstract}
The general construction of \emph{lattice (co)homology} assigns to a lattice $\mathbb{Z}^r$ and a weight function $w:\mathbb{Z}^r \to \mathbb{Z}$ a bigraded $\mathbb{Z}[U]$-module $\mathbb{H}_*$.
The weight function $w$ is often obtained from some geometric data as the difference of two `height functions'. In this paper we consider the case when these height functions are Hilbert functions of valuative multifiltrations on a Noetherian $k$-algebra $\mathcal{O}$ and a finitely generated $\mathcal{O}$-module $M$. We introduce the notion of \textit{`realizable submodules'} in $M$, the prime example of which are finite codimensional integrally closed submodules in the sense of Rees (or integrally closed ideals when $M=\mathcal{O}$). We prove, that whenever two sets of `extended' discrete valuations \textit{`realize'} the same submodule $N \leq M$, then, although the corresponding lattices and weight functions might be different, the resulting lattice homology modules are isomorphic and have Euler characteristic $\dim_k(M/N)$. In this way, we associate a well-defined lattice homology to any quotient of type $M/N$, where $N$ is a realizable submodule of $M$. We also present some structural and computational results: e.g., we geometrically characterize the (lattice) homological dimension of integrally closed monomial ideals of $k[x,y]$.

The main upshot of the paper, however, is the possibility of categorifying numerical invariants defined as codimensions of realizable submo\-dules or integrally closed ideals. The geometric applications include: the \emph{delta invariant}  $\delta(C, o)$ of a reduced curve singularity; the \emph{geometric genus} $p_g(X, o)$, the \emph{irregularity} $q(X, o)$ and the various \emph{plurigenera} of higher dimensional isolated normal singularities. The corres\-ponding categorifications generalize the analytic lattice homologies of Ágoston and the first author. 

We also discuss relations with singularity deformation theory: we initiate the study of \textit{`conductor ideal constant'} Gorenstein deformations, which automatically preserve the analytic lattice homology. 

The manuscript presents computations
for several key families of algebras and singularities, including, e.g., 
Newton nondegenerate singularities
and deformations which does not preserve the lattice homology.
\end{abstract}

\maketitle

\thispagestyle{empty}

\vspace*{-1cm}

\tableofcontents

\newpage

\section{Introduction}

\subsection{Lattice homology constructions}\,

The theory of lattice (co)homology is a connecting bridge between several mathematical fields. In its current form, it encompasses a family of closely related invariants of isolated singularities of complex analytic varieties, making connections to analytic invariants, commutative algebra, and, most prominently, to low-dimensional topology, where it also gained independent development. In fact, one can think of (the analytic) lattice homology as a complex analytic analogue of Heegaard Floer theory.

The original topological version of lattice homology \cite{Nlattice, NBook} is defined combinatorially from the dual resolution graph of a normal surface singularity, hence,  it is strongly related to the low-dimensional topology  of the link. It is, in fact, equivalent to its Heegaard Floer homology (see \cite{Nlattice, OSSzSpectral} for some special cases, and \cite{Zemke} for the most general statement), which is, a priori, much more difficult to compute.
Moreover, from  the analytic viewpoint, this topological lattice homology {has the ability} to bound (or, in certain cases, compute) analytic invariants as well, making it a useful tool for the continuation and generalization of the Artin--Laufer program.

Parallel to the topological {version}, the analytic lattice homology \cite{AgArr, AgNe1, AgNeHigh, NBook} is also defined, which reveals even deeper information about the analytic type: it {can help} to distinguish different analytic structures supported on a fixed topological type and categorifies the geometric genus. It comes equipped with a natural morphism from the topological version.

The analytic lattice homology of reduced curve singularities was also recently introduced by Ágoston and the first author in \cite{AgNeCurves}. This is a categorification of the delta invariant, 
{moreover, it conjecturally has the following `functoriality' property: flat singularity}
deformations induce nontrivial graded morphisms {between} the analytic lattice homology {modules}. Its low-dimensional topology counterparts were defined by Ozsváth, Stipsicz and Szabó for knots \cite{OSSzKnots}, and Borodzik, Liu and Zemke for links \cite{BLZ}. They are all strongly connected  with the Heegaard Floer Knot\,/\,Link homology (see, e.g., \cite{NFilt}, where a connecting spectral sequence was presented).

The goal of this manuscript is to show that the {domain and scope} of lattice homology can be broadly {extended, it can be} used in {several} different  situations, providing new, computable and interesting invariants. 
We present a simple algebraic language which unifies all the existing analytic lattice homology theories and generalizes them to a much more abstract setting. This viewpoint enables us to rely on the very same, so-called Independence Theorem, when proving well-definedness and structural results, using combinatorial techniques instead of the original case-by-case analytic arguments. 

Moreover, our methods extend to the purely commutative algebraic setting of finite codimensional integrally closed ideals of $k$-algebras (or, more generally, integrally closed submodules in a finitely generated module), obtaining a well-defined categorification of their codimension. This level of generality holds much promise both in singularity theory and in the more general realm of commutative algebra as well. 
Our results open the gates to multiple research directions and raise interesting questions, we present a large list of these at the end of the manuscript.

Below we present a map of existing lattice homology theories and their connections to low dimensional topology in order to provide better context for our new construction introduced in this article.

\begin{figure}[H]
\centering{
\resizebox{15cm}{!}{
\begin{tikzcd}[ampersand replacement=\&, row sep=large, column sep=tiny]
	{\textit{\Huge Present construction}} \&  \& {\textit{\Huge Analytic LH's}} \& {\cdstack{\textit{\Huge Filtered}}{\textit{\Huge analytic LH's}}}\& {\textit{\Huge Topological LH's}} \&  \& \cdstack{ \textit{\Huge Low dimensional}}{\textit{\Huge topology invariants}} \&  \&  \\
	\&\& \dvdstack{\textit{\Huge Analytic LH}}{\textit{\Huge of isolated singularities}}{\text{\Huge [Ágoston--Némethi \cite{AgNeHigh}]}} \&\& \textit{\Huge ?} \&\&\&\& {} \textit{\Huge $\dim_{\mathbb{C}} > 2$} \\
	\&\&\& {\dvdstack{\textit{\Huge Filtered LH of}}{\textit{\Huge surface singularities}}{\text{\Huge [Némethi \cite{NFilt2}]}}}\&\&\&\&\&  {} \\
	\cdstack{\textit{\Huge LH of realizable}}{\textit{\Huge submodules}} \arrow[uurr, bend left=8, leftrightarrow] \arrow[rr, leftrightarrow]  \arrow[ddrr, bend right=20, leftrightarrow] \&\& \dvdstack{\textit{\Huge Analytic LH}}{\textit{\Huge of normal surface singularities}}{\text{\Huge [Ágoston--Némethi \cite{AgNe1}]}} \arrow[uu, rightsquigarrow, "\textit{\LARGE generalization}" ']   \&\& \dvdstack{\textit{\Huge Topological LH}}{\textit{\Huge of normal surface singularities}}{\text{\Huge [Némethi \cite{Nlattice, NBook}]}} \arrow[rr, leftrightarrow, "\textit{\LARGE $\cong$}", "\text{\LARGE \cite{Zemke}}" '] \arrow[ll,  "\textit{\LARGE $\exists$  morphism}"] \arrow[ul, rightsquigarrow, bend left=10] \&\& \cdstack{\textit{\Huge Heegaard Floer}}{\textit{\Huge homology}} \arrow[dd, rightsquigarrow, "\textit{\LARGE filtration}" '] \&\& {} \textit{\Huge $\dim_{\mathbb{C}} = 2$} \\
	\&  \&\& {\dvdstack{\textit{\Huge Filtered LH of}}{\textit{\Huge curve singularities}}{\text{\Huge [Némethi \cite{NFilt}]}}} \arrow[rd,  bend right=10]\&\&\&\&\& {} \\
	\&\& \dvdstack{\textit{\Huge Analytic LH of}}{ \textit{\Huge curve singularities}}{\text{\Huge [Ágoston--Némethi \cite{AgNeCurves}]}} \ar[ru, rightsquigarrow, bend right=10] \&\& \dvdstack{\textit{\Huge Knot\,/\,Link lattice homology}}{\text{\Huge [Ozsváth--Stipsicz--Szabó \cite{OSSzKnots}]}}{ \text{\Huge [Borodzik--Liu--Zemke \cite{BLZ}]}}  \arrow[rr, leftrightarrow, "\textit{\LARGE $\cong$}"] \&\& \cdstack{\textit{\Huge Heegaard Floer}}{\textit{\Huge Knot\,/\,Link homology}} \&\& {} \textit{\Huge $\dim_{\mathbb{C}} =1$} \\
    \&  \&\& \&\&\&\&\& {}
\end{tikzcd}}
\caption{A map of lattice homology (LH) theories and the present construction} \label{figure1}} 
\end{figure}
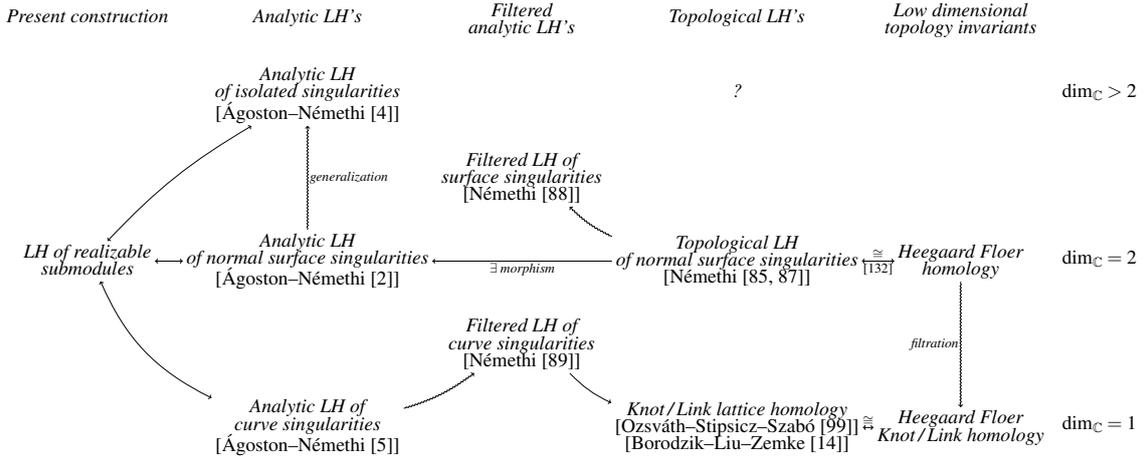

\vspace{-1cm}

\subsection{The lattice homology module associated with a lattice and a weight function}\,

The lattice homology construction starts with a lattice $\mathbb{Z}^r$ 
and a weight function $w:\Z^r\to\Z$, which has some finiteness properties and is bounded from below. In fact, this weight function is the main geometric ingredient, the one which differentiates the distinct theories, and ideally carries important information about the studied object.  
{For example,} in the topological lattice homology of normal surface singularities the lattice is freely generated by the 
irreducible components of the exceptional divisor of a good resolution, whereas the weight function is given by the Riemann--Roch expression on the lattice points considered as homological cycles \cite{Nlattice}. 
In the analytic case, the weight function (defined in the first quadrant of the same lattice) is given by the difference of the Hilbert function corresponding to the divisorial valuations and the rank of the first sheaf cohomology of the structure sheaves of these divisors \cite{AgArr,AgNe1}.  
The level sets and critical points of these weight functions determine a fundamentally new type of geometry of these singularities, which lies at the heart of the definition of lattice homology. 

The construction itself is based on the {following} filtration of $\mathbb{Z}^r \otimes \mathbb{R}=\mathbb{R}^r$ (equipped with the standard cubical decomposition), the `tower of spaces'
$\emptyset=S_{m_w-1}\subset S_{m_w}\subset S_{m_w+1}\subset \cdots \subset S_n \subset \cdots \subset \mathbb{R}^r$, where $S_n$ is the union of those cubes, which have all vertices with $w$-weight at most $n\in \mathbb{Z}$ and $m_w$ is the minimal value of the weight function. Then the bigraded lattice homology $\mathbb{Z}[U]$-module is 

\begin{center}
    $\mathbb{H}_{\ast}(\mathbb{R}^r, w)=\oplus_{q \geq 0}\mathbb{H}_q(\mathbb{R}^r, w)=\oplus_{q \geq 0} \oplus_{n \geq m_w}\mathbb{H}_{q,-2n}(\mathbb{R}^r, w)$, with $\mathbb{H}_{q,-2n}=H_q(S_n, \mathbb{Z})$. 
\end{center}
The $q$-grading is called the \textit{homological grading}, whereas the $n$-grading is the \textit{weight grading}.  This latter is --- by definition --- \emph{even}, in order to emphasize the compatibility with the Heegaard Floer theories. The bigraded $U$-action is given by the natural inclusion maps 
 $H_q(S_{n-1},\Z)\to H_q(S_n,\Z)$. 
  The {\it reduced}
 lattice homology 
also  appears naturally: $\bH_{\rm red, *}(\mathbb{R}^r, w)=\oplus _q \bH_{{\rm red}, q}$, where $\bH_{{\rm red}, q}=\oplus _n
 \widetilde{H}_q(S_n,\bZ)$.

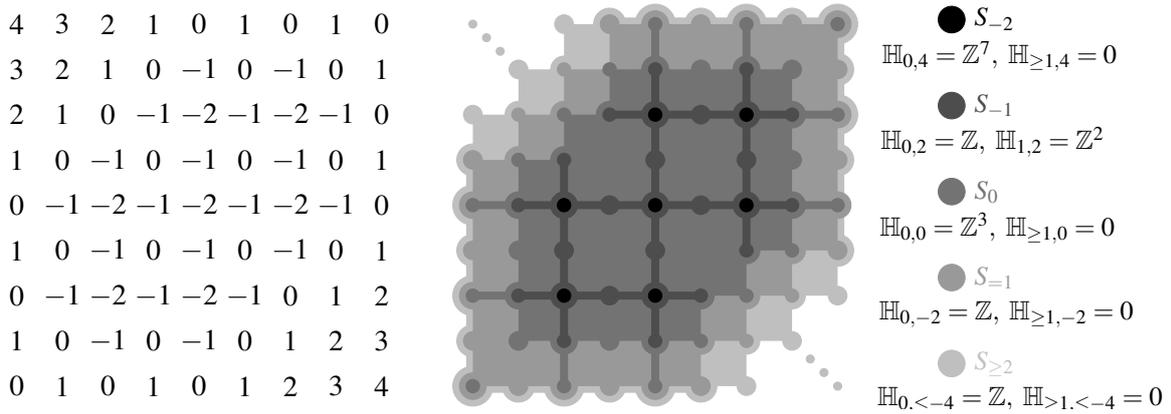
\begin{figure}
\definecolor{wqwqwq}{rgb}{0.3,0.3,0.3}
\definecolor{cqcqcq}{rgb}{0.75,0.75,0.75}
\definecolor{aqaqaq}{rgb}{0.6,0.6,0.6}
\definecolor{yqyqyq}{rgb}{0.45,0.45,0.45}
    \centering
\begin{tikzpicture}[line cap=round,line join=round,>=triangle 45,x=0.60cm,y=0.60cm]
\clip(-10.5,-0.5) rectangle (24,9);

\draw  (-10,0) node {$0$};
\draw  (-9,0) node {$1$};
\draw  (-8,0) node {$0$};
\draw  (-7,0) node {$1$};
\draw  (-6,0) node {$0$};
\draw  (-5,0) node {$1$};
\draw  (-4,0) node {$2$};
\draw  (-3,0) node {$3$};
\draw  (-2,0) node {$4$};

\draw  (-10,1) node {$1$};
\draw  (-9,1) node {$0$};
\draw  (-8,1) node {$-1$};
\draw  (-7,1) node {$0$};
\draw  (-6,1) node {$-1$};
\draw  (-5,1) node {$0$};
\draw  (-4,1) node {$1$};
\draw  (-3,1) node {$2$};
\draw  (-2,1) node {$3$};

\draw  (-10,2) node {$0$};
\draw  (-9,2) node {$-1$};
\draw  (-8,2) node {$-2$};
\draw  (-7,2) node {$-1$};
\draw  (-6,2) node {$-2$};
\draw  (-5,2) node {$-1$};
\draw  (-4,2) node {$0$};
\draw  (-3,2) node {$1$};
\draw  (-2,2) node {$2$};

\draw  (-10,3) node {$1$};
\draw  (-9,3) node {$0$};
\draw  (-8,3) node {$-1$};
\draw  (-7,3) node {$0$};
\draw  (-6,3) node {$-1$};
\draw  (-5,3) node {$0$};
\draw  (-4,3) node {$-1$};
\draw  (-3,3) node {$0$};
\draw  (-2,3) node {$1$};

\draw  (-10,4) node {$0$};
\draw  (-9,4) node {$-1$};
\draw  (-8,4) node {$-2$};
\draw  (-7,4) node {$-1$};
\draw  (-6,4) node {$-2$};
\draw  (-5,4) node {$-1$};
\draw  (-4,4) node {$-2$};
\draw  (-3,4) node {$-1$};
\draw  (-2,4) node {$0$};

\draw  (-10,5) node {$1$};
\draw  (-9,5) node {$0$};
\draw  (-8,5) node {$-1$};
\draw  (-7,5) node {$0$};
\draw  (-6,5) node {$-1$};
\draw  (-5,5) node {$0$};
\draw  (-4,5) node {$-1$};
\draw  (-3,5) node {$0$};
\draw  (-2,5) node {$1$};

\draw  (-10,6) node {$2$};
\draw  (-9,6) node {$1$};
\draw  (-8,6) node {$0$};
\draw  (-7,6) node {$-1$};
\draw  (-6,6) node {$-2$};
\draw  (-5,6) node {$-1$};
\draw  (-4,6) node {$-2$};
\draw  (-3,6) node {$-1$};
\draw  (-2,6) node {$0$};

\draw  (-10,7) node {$3$};
\draw  (-9,7) node {$2$};
\draw  (-8,7) node {$1$};
\draw  (-7,7) node {$0$};
\draw  (-6,7) node {$-1$};
\draw  (-5,7) node {$0$};
\draw  (-4,7) node {$-1$};
\draw  (-3,7) node {$0$};
\draw  (-2,7) node {$1$};

\draw  (-10,8) node {$4$};
\draw  (-9,8) node {$3$};
\draw  (-8,8) node {$2$};
\draw  (-7,8) node {$1$};
\draw  (-6,8) node {$0$};
\draw  (-5,8) node {$1$};
\draw  (-4,8) node {$0$};
\draw  (-3,8) node {$1$};
\draw  (-2,8) node {$0$};

\fill[line width=0pt,color=cqcqcq,fill=cqcqcq,fill opacity=1] (0,0) -- (6,0) -- (6,1) -- (7,1) -- (7,2) -- (8,2) -- (8,8) -- (2,8) -- (2,7) -- (1,7) -- (1,6) -- (0,6) -- cycle;

\filldraw [cqcqcq] (0,0) circle (7.5pt);
\filldraw [cqcqcq] (1,0) circle (6pt);
\filldraw [cqcqcq] (2,0) circle (7.5pt);
\filldraw [cqcqcq] (3,0) circle (6pt);
\filldraw [cqcqcq] (4,0) circle (7.5pt);
\filldraw [cqcqcq] (5,0) circle (6pt);
\filldraw [cqcqcq] (0,1) circle (6pt);
\filldraw [cqcqcq] (0,2) circle (7.5pt);
\filldraw [cqcqcq] (0,3) circle (6pt);
\filldraw [cqcqcq] (0,4) circle (7.5pt);
\filldraw [cqcqcq] (0,5) circle (6pt);
\filldraw [cqcqcq] (8,8) circle (7.5pt);
\filldraw [cqcqcq] (7,8) circle (6pt);
\filldraw [cqcqcq] (6,8) circle (7.5pt);
\filldraw [cqcqcq] (5,8) circle (6pt);
\filldraw [cqcqcq] (4,8) circle (7.5pt);
\filldraw [cqcqcq] (3,8) circle (6pt);
\filldraw [cqcqcq] (8,7) circle (6pt);
\filldraw [cqcqcq] (8,6) circle (7.5pt);
\filldraw [cqcqcq] (8,5) circle (6pt);
\filldraw [cqcqcq] (5,4) circle (7.5pt);
\filldraw [cqcqcq] (8,3) circle (6pt);
\filldraw [cqcqcq] (6,0) circle (3.5pt);
\filldraw [cqcqcq] (7,1) circle (3.5pt);
\filldraw [cqcqcq] (8,2) circle (3.5pt);
\filldraw [cqcqcq] (0,6) circle (3.5pt);
\filldraw [cqcqcq] (1,7) circle (3.5pt);
\filldraw [cqcqcq] (2,8) circle (3.5pt);
\filldraw [cqcqcq] (6,1) circle (5pt);
\filldraw [cqcqcq] (7,2) circle (5pt);
\filldraw [cqcqcq] (1,6) circle (5pt);
\filldraw [cqcqcq] (2,7) circle (5pt);

\filldraw [cqcqcq] (0,8) circle (1.5pt);
\filldraw [cqcqcq] (0.3,7.7) circle (1.5pt);
\filldraw [cqcqcq] (0.6,7.4) circle (1.5pt);
\filldraw [cqcqcq] (8,0) circle (1.5pt);
\filldraw [cqcqcq] (7.7,0.3) circle (1.5pt);
\filldraw [cqcqcq] (7.4,0.6) circle (1.5pt);

\draw [line width=5pt,color=cqcqcq] (0,0)-- (5,0);
\draw [line width=5pt,color=cqcqcq] (0,0)-- (0,5);
\draw [line width=5pt,color=cqcqcq] (3,8)-- (8,8);
\draw [line width=5pt,color=cqcqcq] (8,3)-- (8,8);

\fill[line width=0pt,color=aqaqaq,fill=aqaqaq,fill opacity=1] (0,0) -- (5,0) -- (5,1) -- (6,1) -- (6,2) -- (7,2) -- (7,3) -- (8,3) -- (8,8) -- (3,8) -- (3,7) -- (2,7) -- (2,6) -- (1,6) -- (1,5) -- (0,5) -- cycle;

\filldraw [aqaqaq] (0,0) circle (5pt);
\filldraw [aqaqaq] (1,0) circle (3.5pt);
\filldraw [aqaqaq] (2,0) circle (5pt);
\filldraw [aqaqaq] (3,0) circle (3.5pt);
\filldraw [aqaqaq] (4,0) circle (5pt);
\filldraw [aqaqaq] (5,0) circle (3.5pt);
\filldraw [aqaqaq] (0,1) circle (3.5pt);
\filldraw [aqaqaq] (6,1) circle (2.5pt);
\filldraw [aqaqaq] (7,2) circle (2.5pt);
\filldraw [aqaqaq] (0,3) circle (3.5pt);
\filldraw [aqaqaq] (8,3) circle (3.5pt);
\filldraw [aqaqaq] (0,5) circle (3.5pt);
\filldraw [aqaqaq] (8,5) circle (3.5pt);
\filldraw [aqaqaq] (1,6) circle (2.5pt);
\filldraw [aqaqaq] (2,7) circle (2.5pt);
\filldraw [aqaqaq] (8,7) circle (3.5pt);
\filldraw [aqaqaq] (3,8) circle (3.5pt);
\filldraw [aqaqaq] (5,8) circle (3.5pt);
\filldraw [aqaqaq] (7,8) circle (3.5pt);
\filldraw [aqaqaq] (5,1) circle (5pt);
\filldraw [aqaqaq] (6,2) circle (5pt);
\filldraw [aqaqaq] (7,3) circle (5pt);
\filldraw [aqaqaq] (1,5) circle (5pt);
\filldraw [aqaqaq] (2,6) circle (5pt);
\filldraw [aqaqaq] (0,2) circle (5pt);
\filldraw [aqaqaq] (0,4) circle (5pt);
\filldraw [aqaqaq] (8,4) circle (5pt);
\filldraw [aqaqaq] (8,6) circle (5pt);
\filldraw [aqaqaq] (4,8) circle (5pt);
\filldraw [aqaqaq] (6,8) circle (5pt);
\filldraw [aqaqaq] (8,8) circle (5pt);

\fill[line width=0pt,color=yqyqyq,fill=yqyqyq,fill opacity=1] (1,1) -- (5,1) -- (5,2) -- (6,2) -- (6,3) -- (7,3) -- (7,7) -- (3,7) -- (3,6) -- (2,6) -- (2,5) -- (1,5) -- cycle;

\draw [line width=3pt,color=yqyqyq] (4,7)-- (4,8);
\draw [line width=3pt,color=yqyqyq] (6,7)-- (6,8);
\draw [line width=3pt,color=yqyqyq] (7,6)-- (8,6);
\draw [line width=3pt,color=yqyqyq] (7,4)-- (8,4);
\draw [line width=3pt,color=yqyqyq] (4,1)-- (4,0);
\draw [line width=3pt,color=yqyqyq] (2,1)-- (2,0);
\draw [line width=3pt,color=yqyqyq] (1,2)-- (0,2);
\draw [line width=3pt,color=yqyqyq] (1,4)-- (0,4);

\filldraw [yqyqyq] (0,0) circle (2.5pt);
\filldraw [yqyqyq] (2,0) circle (2.5pt);
\filldraw [yqyqyq] (4,0) circle (2.5pt);
\filldraw [yqyqyq] (1,1) circle (3.5pt);
\filldraw [yqyqyq] (2,1) circle (5pt);
\filldraw [yqyqyq] (3,1) circle (3.5pt);
\filldraw [yqyqyq] (4,1) circle (5pt);
\filldraw [yqyqyq] (5,1) circle (2.5pt);
\filldraw [yqyqyq] (0,2) circle (2.5pt);
\filldraw [yqyqyq] (1,2) circle (5pt);
\filldraw [yqyqyq] (5,2) circle (5pt);
\filldraw [yqyqyq] (6,2) circle (2.5pt);
\filldraw [yqyqyq] (1,3) circle (3.5pt);
\filldraw [yqyqyq] (3,3) circle (3.5pt);
\filldraw [yqyqyq] (5,3) circle (3.5pt);
\filldraw [yqyqyq] (6,3) circle (5pt);
\filldraw [yqyqyq] (7,3) circle (2.5pt);
\filldraw [yqyqyq] (0,4) circle (2.5pt);
\filldraw [yqyqyq] (1,4) circle (5pt);
\filldraw [yqyqyq] (7,4) circle (5pt);
\filldraw [yqyqyq] (8,4) circle (2.5pt);
\filldraw [yqyqyq] (1,5) circle (2.5pt);
\filldraw [yqyqyq] (2,5) circle (5pt);
\filldraw [yqyqyq] (3,5) circle (3.5pt);
\filldraw [yqyqyq] (5,5) circle (3.5pt);
\filldraw [yqyqyq] (7,5) circle (3.5pt);
\filldraw [yqyqyq] (2,6) circle (2.5pt);
\filldraw [yqyqyq] (3,6) circle (5pt);
\filldraw [yqyqyq] (7,6) circle (5pt);
\filldraw [yqyqyq] (8,6) circle (2.5pt);
\filldraw [yqyqyq] (3,7) circle (2.5pt);
\filldraw [yqyqyq] (4,7) circle (5pt);
\filldraw [yqyqyq] (5,7) circle (3.5pt);
\filldraw [yqyqyq] (6,7) circle (5pt);
\filldraw [yqyqyq] (7,7) circle (3.5pt);
\filldraw [yqyqyq] (4,8) circle (2.5pt);
\filldraw [yqyqyq] (6,8) circle (2.5pt);
\filldraw [yqyqyq] (8,8) circle (2.5pt);

\draw [line width=3pt,color=wqwqwq] (1,2)-- (5,2);
\draw [line width=3pt,color=wqwqwq] (2,1)-- (2,5);
\draw [line width=3pt,color=wqwqwq] (1,4)-- (7,4);
\draw [line width=3pt,color=wqwqwq] (6,7)-- (6,3);
\draw [line width=3pt,color=wqwqwq] (3,6)-- (7,6);
\draw [line width=3pt,color=wqwqwq] (4,1)-- (4,7);

\begin{scriptsize}
\filldraw [wqwqwq] (2,1) circle (2.5pt);
\filldraw [wqwqwq] (4,1) circle (2.5pt);
\filldraw [wqwqwq] (1,2) circle (2.5pt);
\filldraw [wqwqwq] (2,2) circle (5pt);
\filldraw [wqwqwq] (3,2) circle (3.5pt);
\filldraw [wqwqwq] (4,2) circle (5pt);
\filldraw [wqwqwq] (5,2) circle (2.5pt);
\filldraw [wqwqwq] (2,3) circle (3.5pt);
\filldraw [wqwqwq] (4,3) circle (3.5pt);
\filldraw [wqwqwq] (6,3) circle (2.5pt);
\filldraw [wqwqwq] (1,4) circle (2.5pt);
\filldraw [wqwqwq] (2,4) circle (5pt);
\filldraw [wqwqwq] (3,4) circle (3.5pt);
\filldraw [wqwqwq] (4,4) circle (5pt);
\filldraw [wqwqwq] (5,4) circle (3.5pt);
\filldraw [wqwqwq] (6,4) circle (5pt);
\filldraw [wqwqwq] (7,4) circle (2.5pt);
\filldraw [wqwqwq] (2,5) circle (2.5pt);
\filldraw [wqwqwq] (4,5) circle (3.5pt);
\filldraw [wqwqwq] (6,5) circle (3.5pt);
\filldraw [wqwqwq] (3,6) circle (2.5pt);
\filldraw [wqwqwq] (4,6) circle (5pt);
\filldraw [wqwqwq] (5,6) circle (3.5pt);
\filldraw [wqwqwq] (6,6) circle (5pt);
\filldraw [wqwqwq] (7,6) circle (2.5pt);
\filldraw [wqwqwq] (4,7) circle (2.5pt);
\filldraw [wqwqwq] (6,7) circle (2.5pt);

\filldraw [black] (2,2) circle (2.5pt);
\filldraw [black] (4,2) circle (2.5pt);
\filldraw [black] (2,4) circle (2.5pt);
\filldraw [black] (4,4) circle (2.5pt);
\filldraw [black] (6,4) circle (2.5pt);
\filldraw [black] (4,6) circle (2.5pt);
\filldraw [black] (6,6) circle (2.5pt);
\end{scriptsize}

\begin{small}
\filldraw [black] (10.5,8.1) circle (5pt) node [anchor=west]{\ \ $S_{-2}$};
\draw (11.55, 7.3) node {$\mathbb{H}_{0,4}=\mathbb{Z}^7, \ \mathbb{H}_{\geq 1,4}=0$};
\filldraw [wqwqwq] (10.5,6.2) circle (5pt) node [anchor=west]{\ \ $S_{-1}$};
\draw (11.4, 5.4) node {$\mathbb{H}_{0,2}=\mathbb{Z}, \ \mathbb{H}_{1,2}=\mathbb{Z}^2$};
\filldraw [yqyqyq] (10.5,4.3) circle (5pt) node [anchor=west]{\ \ $S_{0}$};
\draw (11.5, 3.5) node {$\mathbb{H}_{0,0}=\mathbb{Z}^3, \ \mathbb{H}_{\geq 1,0}=0$};
\filldraw [aqaqaq] (10.5,2.4) circle (5pt) node [anchor=west]{\ \ $S_{= 1}$};
\draw (11.7, 1.6) node {$\mathbb{H}_{0,-2}=\mathbb{Z}, \ \mathbb{H}_{\geq 1, -2}=0$};
\filldraw [cqcqcq] (10.5,0.5) circle (5pt) node [anchor=west]{\ \ $S_{\geq 2}$};
\draw (12, -0.3) node {$\mathbb{H}_{0,\leq-4}=\mathbb{Z}, \ \mathbb{H}_{\geq 1, \leq -4}=0$};
\end{small}

\end{tikzpicture}
\caption{The $S_n$ spaces associated with a concrete weight function on $R((0,0), (8,8))$}
    \label{fig:pelda}
\end{figure}

In some cases one only considers the cubes contained in a rectangle $R(0, d):=\{x \in \mathbb{R}^r \,:\, 0 \leq x \leq 0\}$, which provides $\mathbb{H}^{\ast}(R(0, d), w)$. An example with a rank $2$ lattice is presented in Figure \ref{fig:pelda}.  

 The requirements on $w$ ensure that the $\{S_n\}_n$ spaces are finite CW complexes, which are usually empty for $n\ll 0$ and contractible for $n \gg 0$.
 Therefore, in such a case
 $\bH_{\rm red, *}$ has finite $\bZ$-rank and the \emph{Euler characteristic} is defined as $ eu(\bH_*):= -m_w+\textstyle{\sum_q}\, (-1)^q\ {\rm rank}_\bZ\
 \bH_{{\rm red}, q}$. Notice that 
 the smallest weight $m_w$ 
  (the analogue of the $d$-invariant used in Heegaard Floer homology \cite{OSzF})
 and 
 $\bH_{\rm red, *}$ together determine the whole  $\bH_*$.

An important observation is that the general construction\,/\,definition  of lattice homo\-logy  is very flexible: by providing different weight functions one
obtains different lattice homologies.
On the other hand, the following question arises naturally: when do two different weight functions yield the same lattice homology module? The main result of this article, the Independence Theorem, gives a sufficient condition for this in a rather large class of naturally occurring weight functions. 

\subsection{Categorification with lattice homology}\,

In mathematical classification procedures one uses invariants to distinguish different objects. If a certain
invariant is not sufficiently `strong', then one tries to endow it
with some additional structure. A {\it categorification}
of an invariant is a (co)homology theory whose Euler characteristic is the invariant considered. It is clearly a finer invariant: the Betti numbers {and the torsion part provide} more subtle information than the Euler characteristic itself.

In the last decades several famous categorifications were introduced.
E.g., in knot theory, the Khovanov invariant was introduced as
the categorification of the Jones polynomial \cite{Kho},
the Link Heegaard Floer
homology as the categorification of the Alexander polynomial \cite{OSzLink}.
Or, in the 3-manifold theory, the Heegaard Floer homology of
Ozsv\'ath and Szab\'o is the categorification of the Seiberg--Witten invariant
\cite{OSz,OSz7}.

The existing lattice homology theories fit naturally into this trend. The topological lattice homo\-logy $\bH_{top, *}$
associated with links of normal surface singularities 
is a  categorification of the normalized Seiberg--Witten invariant \cite{NJEMS}, whereas the
 analytic version 
gives a categorification of the geometric genus 
\cite{AgNe1}. 
In the reduced curve singularity case, the Euler characteristic of the analytic lattice homology is the delta invariant \cite{AgNeCurves}. 

In fact, Ágoston and the first author described a simple recipe to obtain lattice homological categorifications. 
They studied the properties of lattice homology modules associated with special, so called `split' weight functions of the form $w(\ell)=h(\ell)+h^\circ(\ell)-h^\circ(0)$,
where $h$ is some increasing and $h^\circ$ some decreasing function along the lattice. 
They realized, that if the pair $(h, h^\circ)$   
satisfies certain combinatorial properties, then the Euler characteristic of the lattice homology can be computed from $h^\circ$ (cf. \cite[Theorem 5.2.1]{AgNe1} or the Categorification Theorem \ref{th:comblattice} here).  
{This theorem provides a powerful tool for categorifying numerical invariants via lattice homology.}
In fact, Ágoston and the first author built their analytic lattice homology theories on this cornerstone observation.

In this article we will elaborate on this construction, considering the case when the height functions $h$ and $h^\circ$ come as Hilbert functions of some valuative multifiltrations on a $k$-algebra $\mathcal{O}$ and a module $M$ over it. In this {setting} the {required combinatorial properties can be {forced to hold} by a slight adjustment in the definition. In fact, this is the multiplication by $2$ in (\ref{filtrations-intro}), which marks a shift in philosophy: moving from searching height functions with right combinatorial properties to deforming interesting\,/\,relevant ones to meet the requirements. 
The Independence Theorem then claims that the resulting lattice homology is a well-defined categorification of the $k$-codimension of a submodule $N \leq M$ `realized' by the multifiltration. This construction enables us to 
give new categorifications of other important analytic invariants of normal surface singularities such as the irregularity or several plurigenera, as well as to extend the definition of the analytic lattice homology to singularities with not necessarily $\mathbb{Q}HS^3$ link (a necessary requirement in the old versions). Even more, we certainly hope that the generality of the construction will help us produce new invariants in many geometric\,/\,algebraic situations where similar valuative multifiltrations occur naturally. We firmly believe that,  through our results, purely commutative algebraic problems will become accessible, too.

\subsection{The Independence Theorem and the lattice homology of realizable submodules}\label{ss:main}\,

Consider the following setup: let $k$ be a field, $\mathcal{O}$ a Noetherian $k$-algebra and $M$ a finitely generated $\mathcal{O}$-module.  Given a discrete (semi-)valuation $\mathfrak{v}:\mathcal{O} \rightarrow \mathbb{N} \cup \{\infty\}$ on $\mathcal{O}$, we call a map $\mathfrak{v}^M: M \rightarrow \overline{\mathbb{Z}}=\mathbb{Z}\cup \{\infty\}$  its \textit{`extension'}, if it satisfies the following (see Definition \ref{def:edv}):
\begin{align*}
&\text{(a)}\ \frv^M(fm) = \frv(f) + \frv^M(m) &\ \text{ for all }f \in \cO, \ m \in M, \text{ and} \\
 &\text{(b)}\ \frv^M(m+m') \geq \min \{ \frv^M(m), \frv^M(m') \} &\ \text{ for all } m, m' \in  M. \hspace{14mm}
\end{align*}
 If $M$ is not torsion, such extensions always exist, {in fact}, there are infinitely many (cf. Remark \ref{rem:translating}). 

    To a finite collection $\mathcal{D}=\{ (\mathfrak{v}_v, \mathfrak{v}_v^M)\}_{v=1}^r$ of such extended discrete valuations on $\mathcal{O}$ and $M$ {(of \emph{cardinality} $|\mathcal{D}|=r$)} we associate the following {descending} multifiltrations of $\mathcal{O}$ and $M$:
\begin{equation}\label{filtrations-intro}
    \mathbb{Z}^r \ni \ell = \sum_v \ell_v e_v \mapsto  \begin{cases}
        \, \mathcal{F}_{\mathcal{D}}(\ell) :=\{ \ f \in \mathcal{O}\,:\, 2\cdot\mathfrak{v}_v(f) \geq \ell_v \ \forall v\} \ \ \, \triangleleft \ \mathcal{O}\,; \\
        \mathcal{F}_{\mathcal{D}}^M(\ell) :=\{ m \in M\,:\, 2\cdot\mathfrak{v}_v^M(m) \geq \ell_v \ \forall v\} \leq M.
    \end{cases}
\end{equation}
(The multiplication by $2$ {serves technical purposes}, see subsection \ref{ss:CDPforD}.) In fact, we will only consider such extended valuations for which the following \emph{`height'} functions are finite {(cf. Assumption \ref{ass:fincodim})}:
\begin{equation*}
    \mathfrak{h}_{\mathcal{D}}:  \ell \mapsto \dim_k\mathcal{O}/\mathcal{F}_{\mathcal{D}}(\ell) \text{ and }\mathfrak{h}^\circ_{\mathcal{D}}: \ell \mapsto \dim_kM/\mathcal{F}_{\mathcal{D}}^M(-\ell).
\end{equation*} 
Finally, we consider the weight function $w_{ \mathcal{D}}(\ell)=\mathfrak{h}_{\mathcal{D}}(\ell)+\mathfrak{h}^\circ_{\mathcal{D}} (\ell)-\mathfrak{h}^\circ_{\mathcal{D}} (0)$ on the lattice $\mathbb{Z}^r$. The main theorem of this manuscript states that the corresponding bigraded lattice homology $\mathbb{Z}[U]$-module only depends on what the submodule $\mathcal{F}_{\mathcal{D}}^M(0)\leq M$ is, not on the collection $\mathcal{D}$ itself:

\begin{theorem}[= Independence Theorem \ref{th:IndepMod}]\label{th:Indep-intro} Let $\mathcal{D}=\{(\frv_v, \frv^M_v)\}_{v}$ and $\mathcal{D}'=\{(\frv'_{v'}, {\frv'_{v'}}^{M})\}_{v'}$ be two collections of extended discrete valuations.  
       Suppose that $\mathcal{F}_{\mathcal{D}}^M(0) = \mathcal{F}_{\mathcal{D}'}^M(0)=:N \leq M$.
        Then the spaces $S_{n, \mathcal{D}}$ and $S_{n, \mathcal{D}'}$ 
        associated with the corresponding lattices and weight functions 
        are homotopy equivalent for every $n \in \mathbb{Z}$.
        Even more, the homotopy equivalences commute (up to homotopy) with the inclusions
        $S_{n, \mathcal{D}}\hookrightarrow S_{n+1, \mathcal{D}}$ and $S_{n, \mathcal{D}'}\hookrightarrow S_{n+1, \mathcal{D}'}$, hence
        $\mathbb{H}_*(\mathbb{R}^{|\mathcal{D}|}, w_{\mathcal{D}}) \cong  \mathbb{H}_*(\mathbb{R}^{|\mathcal{D}'|}, w_{\mathcal{D}'}) \text{ as bigraded } \mathbb{Z}[U]\text{-modules.}$
        We call it the \emph{ lattice homology of the submodule $N\leq M$} and denote it by $\boxed{\mathbb{H}_*(N \hookrightarrow_\mathcal{O} M)}$.
\end{theorem}

\noindent(For the precise formulation regarding the homotopy type of the $\{S_n\}_n$ filtration see subsection \ref{ss:filteredhom}.) The idea of the proof is that we associate the `missing' extended valuations one by one to $\mathcal{D}$ and $\mathcal{D}'$, to get to the `join' realization  $\mathcal{D} \cup \mathcal{D}'$, and prove that along any single step the homotopy types of the $S_n$-spaces are preserved. This last part is achieved by relying on Quillen's Fiber Theorem \cite[Theorem A]{Quillen}, or, alternatively, the technique of quasifibrations by Dold and Thom \cite{DoldThom}.

\begin{define}[= Definition \ref{def:REAL}]\label{def:Real-intro}
    A finite $k$-codimensional submodule $N \leq M$ is called \textit{`realizable'} if some finite collection of 
    extended discrete valuations 
    $\mathcal{D}$ satisfies 
    $N=\mathcal{F}^M_{\mathcal{D}}(0)$. In this case $\mathcal{D}$ is called a \textit{`realization'} of $N$.
\end{define}

Therefore, the main consequence of the Independence Theorem \ref{th:Indep-intro} is that to any realizable submodule we can assign a well-defined lattice homology computable through any arbitrary rea\-lization.
Its main properties are as follows:
\begin{theorem}[= Theorem \ref{th:properties}]\label{th:propertiesbc-intro}\,

\noindent (b) The Euler characteristic is well-defined and satisfies $eu(\mathbb{H}_*(N \hookrightarrow M))=\dim_k (M/N)$, i.e., the lattice homology $\mathbb{H}_*(N \hookrightarrow M)$ categorifies the codimension ${\rm codim}_k(N \hookrightarrow M)$.  \\
\noindent (c) The homotopy type of the $S_n$ spaces (hence $\mathbb{H}_*(N \hookrightarrow M)$ itself) only depend on the quotient module $M/N$ over the algebra $\mathcal{O}/{\rm Ann}_{\mathcal{O}}(M/N)$, where ${\rm Ann}_{\mathcal{O}}(M/N)$ denotes the annihilator in $\mathcal{O}$ of $M/N$. 
\end{theorem}

If we choose $M = \mathcal{O}$, then a finite codimensional ideal $\mathcal{I} \triangleleft \mathcal{O}$  is realizable if and only if $\mathcal{I} = \mathcal{F}_{\mathcal{D}}(d)$ for some collection $\mathcal{D}$ of discrete (semi-)valuations and lattice point $d \in \mathbb{Z}^{|\mathcal{D}|}$ (see Corollary \ref{cor:idealREAL}). One can also see, that in this case the weight function is automatically symmetric with respect to the lattice point $2d$ (i.e., $w_{\mathcal{D}}(\ell)=w_{\mathcal{D}}(2d-\ell)$ for any $\ell \in \mathbb{Z}^{|\mathcal{D}|}$, cf. Proposition \ref{prop:symforideals}), implying an inherent $\mathbb{Z}_2$-symmetry of {all the $S_n$-spaces}. We claim that the homotopy equivalences of the Independence Theorem respect this $\mathbb{Z}_2$-symmetry up to homotopy, hence it is inherited by the lattice homology module as well (cf. Theorem \ref{th:Indep}). Therefore, we denote the `symmetric' lattice homology of realizable ideals by $\boxed{\mathbb{SH}_{\ast}(\mathcal{I} \triangleleft \mathcal{O})}$, which is a categorification of the codimension ${\rm codim}_k(\mathcal{I} \hookrightarrow \mathcal{O})$. By the above Theorem \ref{th:propertiesbc-intro}, it only depends on the Artin algebra $\mathcal{O}/\mathcal{I}$.

These main results are rather surprising at first. Indeed, a realizable submodule $N\leq M$ can have multiple very different realizations (with different cardinality and different types of valuations --- see, e.g., Examples \ref{ex:uj1}, \ref{bek:10.2.1} and \ref{bek:10.6.3}), yet according to the Independence Theorem \ref{th:Indep-intro}, each of these yields the very same lattice homology. We also emphasize that $\mathbb{H}_*(N\hookrightarrow M)$ is often highly nontrivial (see, e.g., non-trivial $\mathbb{H}_1$ in Example \ref{ex:2552}), it recovers deep information about the input algebras and modules. On a more philosophical note, we could say that this new invariant sensitively captures the geometry in some broad and vague sense of some `universal multiplicative ordering'  of the quotient $\mathcal{O}$-module $M/N$, which is respected by every extended discrete valuation.

\subsection{Realizable and integrally closed submodules and Artin algebras}\,

The Independence Theorem \ref{th:Indep-intro} naturally raises the question of which submodules and ideals are realizable in the sense of Definition \ref{def:Real-intro}. The valuative nature of the concept suggests connections with Rees' integrally closed submodules (cf. \cite{ReesMod}) and ideals. Indeed, we have the following result:

\begin{theorem}[= Theorem \ref{th:REES}]
    Let $k$ be a field, $\mathcal{O}$ a Noetherian $k$-algebra, $N \leq M$ an inclusion of finitely generated $\mathcal{O}$-modules with ${\rm codim}_k(N\hookrightarrow M)<\infty$. Then $N$ is integrally closed in $M$ (in the sense of Rees) if and only if there exists a realization $\mathcal{D}$ of $N$ (in the sense of Definition \ref{def:Real-intro}) such that for every (semi-)valuation $\mathfrak{v}_v \in \mathcal{D}$ the prime ideal $\mathfrak{p}_v=\{ f \in \mathcal{O}\,:\, \mathfrak{v}_v(f)=\infty\}$ is \emph{minimal}. Hence, finite codimension and integral closedness implies realizability.
\end{theorem}

\noindent(See also Corollary \ref{cor:integralviaextendedval} for the equivalent reformulation of the fact that an element $m \in M$ is integral over $N$ in the sense of Rees  using extensions of (semi-)valuations.)

In the case of finite codimensional ideals the notions of realizability and integral closedness agree due to the equational criterion of integral closedness (cf. Corollary \ref{cor:intersection}, see also the valuative criterion of integral closedness or the existence of Rees valuations in \cite{SH}).

On the other hand, part \textit{(c)} of Theorem \ref{th:propertiesbc-intro} implies that our new lattice homology construction is an invariant of {finite $k$-dimensional quotients $M/N$ over Artin $k$-algebras $\mathcal{O}/{\rm Ann}_{\mathcal{O}}(M/N)$, where, by Proposition \ref{prop:pullbackring} \textit{(2)}, this ideal ${\rm Ann}_{\mathcal{O}}(M/N)$ is realizable.} Therefore, the lattice homology theory of realizable submodules concerns the representation theory of `integrally reduced' Artin algebras:

\begin{define}[= Definition \ref{def:intreduced}]
     We call a Noetherian $k$-algebra $A$ \emph{`integrally reduced'} if there exists a presentation $A \cong \mathcal{O}/\mathcal{I}$ with $\mathcal{O}$ a Noetherian $k$-algebra and $\mathcal{I} =\overline{\mathcal{I}}$ an integrally closed ideal. 
\end{define}

We prove that the property of a Noetherian $k$-algebra $A$ that it is form $k[x_1, \ldots, x_m]/\mathcal{I}$ {with $\mathcal{I}=\overline{\mathcal{I}}$} is well-defined, i.e. in any presentation of an integrally reduced algebra with generators and relations, the ideal of relations is inherently integrally closed (cf. Corollary \ref{cor:intredwdef}). In the case when $A$ is additionally Artinian, then, by Theorem \ref{th:propertiesbc-intro} \textit{(c)}, $\boxed{\mathbb{SH}_*(A)}:=\mathbb{SH}_*(\mathcal{I} \triangleleft k[x_1, \ldots, x_m])$ is a well-defined and rather strong invariant categorifying the $k$-dimension of $A$ (nevertheless, not complete: see Example \ref{ex:nonisoArtin}). We also prove that an Artin $k$-algebra is integrally reduced if and only if all the local components of its direct product decomposition are {so} (cf. Proposition \ref{prop:proddecompintred}), and provide Künneth-type relations for the corresponding symmetric lattice homology {modules} (cf. Proposition \ref{prop:chainproduct} and Corollary \ref{prop:Künneth})

\bekezdes \textbf{Realizability in the local complex analytic case.} For $\mathcal{O}=\mathcal{O}_{X, o}$, the local algebra of a complex analytic spacegerm $(X, o)$, $M=\mathcal{O}^p$ a free module  and $N \leq M$ a submodule, Gaffney defined integral closedness via holomorphic map germs of form $\psi:(\mathbb{C}, 0) \rightarrow (X, o)$ (cf. \cite{Gaff, GaffKl}). In fact, his definition is equivalent to Rees' and we also show how it induces realizability (cf. Theorem \ref{th:Gaffney}).

Moreover, in this complex analytic setting Lejeune-Jalabert and Teissier characterized integral reducedness of $A$ by the following internal property (cf. \cite[2.1 Théorème]{LT}): 
\begin{equation*}
\mbox{
 if $a\in A$ and $\gamma(a)=0$
 for any homomorphism $\gamma:A\to \bC\{t\}/(t^m)$, then $a=0$.}
 \end{equation*}
A similar internal characterization in the general case would be highly desirable.

In conclusion, to any finite codimensional submodule (or ideal), which is integrally closed in Rees' or Gaffney's sense, we can associate a well-defined lattice homology with our Independence Theorem.

\subsection{Structural results and {conceptual} meaning}\,

In order to understand what our new lattice homology construction measures, it is fruitful to understand its structure better. First, looking at the homological grading, we see that, by virtue of the construction (a compact sub{complex} of $\mathbb{R}^r$ has trivial $H_{\geq r}$), the \emph{`(lattice) homological deg\-ree'} ${\rm homdim}(N \hookrightarrow M):=\max \{ q\,:\, \mathbb{H}_{q}(N \hookrightarrow M)\not=0\}$ is strictly smaller than the cardinality {$|\mathcal{D}|$ of any  realization $\mathcal{D}$} of $N \lneq M$  (cf. Proposition \ref{prop:homdegupperbound}).
[If $N=M$ we define ${\rm homdim}(N \hookrightarrow M)$ to be $-1$.] 
We observe that in the case of integrally closed ideals  {the minimum of all such cardinalities} might not be obtained by the set of Rees valuations (cf. Example \ref{bek:10.2.2}), i.e., it is a different invariant of the ideal. Moreover, we conjecture the following:
\begin{conj}[= Conjecture \ref{conj:elso}]\label{conj:elso-intro}
    Suppose that $\mathcal{I} \triangleleft k[x_1, x_2]$ is a finite codimensional, integrally closed ideal. If $\mathbb{SH}_{\geq 1}(\mathcal{I} \triangleleft k[x_1, x_2])=0$, then $\mathcal{I}$ can be realized by a single (semi-)valuation.
\end{conj}

In fact, in the case of \emph{monomial} ideals $\mathcal{M} \triangleleft k[x_1, x_2]$ the conjecture holds true, even more, we have the following characterization:

\begin{theorem}[= Theorem \ref{th:homdim}]\label{th:homdim-intro}
    Let $\mathcal{M}$ be a finite codimensional integrally closed monomial ideal in $k[x_1, x_2]$. Then we have the following:
    \begin{align*}
    {\rm homdim}(\mathcal{M} \triangleleft k[x_1,x_2]) =&\, \min \{|\mathcal{D}|\,:\, \mathcal{D} \text{ a realization of } \mathcal{M}\} -1\\ =&\, \min \{|\mathcal{D}|\,:\, \mathcal{D} \text{ a realization with monomial valuations of } \mathcal{M}\} -1.
\end{align*}
\end{theorem}

\noindent (For a more geometric picture see Proposition \ref{prop:linedim-1=kerbdim}.) The proof
relies on finding a special combination of ring elements {(using the Newton polytope), the existence of} which helps to bound both the homological dimension and the cardinality of any realization (cf. Proposition \ref{prop:htopnem0}). 
We expect Theorem \ref{th:homdim-intro} to also be true in higher dimensional polynomial rings, cf. Remark \ref{rem:homdimhighdim}.

Secondly, we also study the weight grading. {It turns out} that the set  $\{n\,: \, (\bH_{{\rm red}, *}(N \hookrightarrow M))_{-2n} \not=0\} $ (with $\mathcal{O}, \ M$ and $N$ {general} as above) is bounded both from above and below (this follows, e.g., from Theorem \ref{th:properties} \textit{(a)}). The lower bound $m_w=-\max\{ n \in \mathbb{Z} \,:\, \mathbb{H}_{0,2n}(N \hookrightarrow M) \neq 0\}$ (cf. Corollary \ref{cor:minw}) is a challenging new invariant, which, in the special case of the analytic lattice homology of reduced curve singularities, is strongly related with the{ir} Cohen--Macauley type (see \cite{HofN}). On the other hand, by a slight generalization of \cite[Theorem 6.1.1]{KNS2}, the upper bound (under some conditions) is universal:

\begin{theorem}[= Nonpositivity Theorem  \ref{th:upperbound}]\label{th:nonpos-intro}
    If $k$ is algebraically closed and $\mathcal{O}$ is \emph{local}, then the $S_n$-spaces are contractible for all $n>0${, while, if $N \neq M$, then $S_0$ is not even connected}. Thus, 
    \begin{center}
        $(\bH_{{\rm red}, q}(N \hookrightarrow M))_{-2n} = 0$ for all $n>0, \ q\geq 0$.
    \end{center}
\end{theorem}

\noindent (For the strongest formulation see Theorem \ref{th:contr}. We also give an adapted version for integrally reduced Artin $k$-algebras based on their product decomposition into local algebras in Theorem \ref{th:NonposArtin}.)
The proof relies on the careful study of \emph{`generalized local minimum points'} of the weight function. For their definition and main properties see subsection \ref{ss:glm} (compare also with \cite[section 3]{KNS2}).

\subsection{Categorification of numerical invariants in complex analytic singularity theory.}\,

A main consequence of the Independence Theorem is the possibility of categorifying important numerical invariants. In algebraic geometry and, in particular, in complex analytic singularity theory numerous important invariants are defined as codimensions of realizable submo\-dules or integrally closed ideals. In most of these cases, already  the corresponding quotient modules or algebras are invariants of the singularities (often they can also be defined intrinsically, without the need of a resolution), hence so are their lattice homology modules. In this article we collect numerous such categorifications, moreover, we give new (more general) descriptions of all existing analytic lattice homologies. 

We hope that the new theory will have similar power and applicability as the topological version
(or, as the $HF$-theory),
with the difference that  in this case
its  applications will be mostly in the analytic theory of singularities, or, more generally, in commutative algebra.

\bekezdes \textbf{Reformulation of the analytic lattice homology of curve singularities.}\,

The analytic lattice homology of reduced  curve singularities was introduced by Ágoston and the first author in \cite{AgNeCurves}.
It connects the analytic and topological properties of these singularities in a natural way. 
 In a certain sense it can be viewed as an algebraic analogue of the Heegaard Floer Link theory applicable to non-planar singularities as well (see the connections in \cite{GorNem2015, NFilt}). 

 The original construction assigns to a reduced curve singularity $(C, o)$ with irreducible decomposition $(C, o)=\bigcup_{v=1}^{r} (C_v, o)$ a lattice $\mathbb{Z}^r$, while the weight function is obtained from the Hilbert function associated with the valuations corresponding to the normalizations of the components. 
 The resulting analytic lattice homology module $\boxed{\bH_{an,*}(C,o)}$ has Euler characteristic the delta invariant $\delta(C,o)$ and is conjectured to be functorial with respect to flat deformations \cite{AgNeCurves}. By previous work of Kubasch and the authors, the $\mathbb{Z}[U]$-module structure of $\bH_{an,0}(C,o)$ detects the Gorenstein property \cite{KNS2}; for plane curve singularities it detects the multiplicity \cite{KNS1, KNS2}, moreover, in the irreducible case it is a complete embedded topological invariant \cite{KS3}.

Our first result is a reformulation of the definition of this analytic lattice homology theory using the framework of realizable submodules. 
Indeed, if we denote by $n_{C}: \overline{C} \rightarrow C$ the normalization map (corresponding to the integral closure $n_C^*: \mathcal{O}_{C, o} \hookrightarrow \overline{\mathcal{O}_{C, o}}$), by $\omega_{C}^R$ the module of Rosenlicht's regular differential forms on $\overline{C}$ and by $\Omega^1_{\overline{C}}$ the submodule in $\omega_{C}^R$ of germs of holomorphic differential forms on $\overline{C}$, then, using  the  perfect duality between $\omega^R_{C}/ \Omega^1_{\overline{C}}$ and
 $\overline{\calO_{C, o}}/\calO_{C,o}$ (cf. \cite[Ch 2 Sect 9]{Serre}),  we get

\begin{theorem}[= Theorem \ref{th:equivforcurves}]
    For any reduced curve singularity $\mathbb{H}_{\mathrm{an}, *}(C, o) \cong \mathbb{H}_*(\Omega^1_{\overline{C}} \hookrightarrow_{\mathcal{O}_{C, o}} \omega_{C}^R)$.
\end{theorem}

Since $\mathcal{O}_{C, o}$ is local we can apply the Nonpositivity Theorem \ref{th:nonpos-intro} as well --- in fact, this was the setting in which it was originally proved by Kubasch and the authors in \cite{KNS2}. On another note, the module $\omega_C^R$ is free of rank one over $\mathcal{O}_{C, o}$ if and only if the curve singularity $(C, o)$ is Gorenstein. Therefore, in this case the weight function is symmetric and, in fact $\mathbb{H}_{an, *}(C, o) \cong \mathbb{SH}_*(\mathcal{C} \triangleleft \mathcal{O}_{C, o})$ (cf. Corollary \ref{cor:curveArtin}), where $\mathcal{C}$ is the conductor ideal ${\rm Ann}_{\mathcal{O}_{C, o}}\big(\overline{\mathcal{O}_{C, o}}/\mathcal{O}_{C, o}\big)$.

\bekezdes \textbf{Generalizations of the analytic lattice homology of normal surface singularities.}\,

The analytic lattice homology of normal surface singularities of Ágoston and the first author \cite{AgNe1} arose as the analogue of the topological version \cite{Nlattice}, with the aim of capturing analytic instead of to\-po\-logical data from the local algebra and various associated sheaf cohomologies. In \cite{AgNe1} it is defined for normal singularities with \emph{rational homology sphere link} and categorifies the geometric genus.  

The original construction of $\boxed{\mathbb{H}_{{an}, *}(X, o)}$ considers a good resolution $\phi: \widetilde{X} \rightarrow X$ of the normal surface singularity $(X, o)$, and assigns to it a lattice formally generated by the irreducible components of the exceptional divisor $E$ (as prime divisors
). The weight function (defined on the positive orthant) is obtained from the Hilbert function associated with the divisorial valuations and the {dimensions} of the first sheaf cohomologies of the corresponding structure sheaves. Well-definedness is proved by careful analysis of the behaviour of these numerical invariants under blow-ups. Our main observation is that, due to Laufer's duality and geometric genus formula \cite{Laufer72}
\begin{equation*}
    p_g(X, o)= \dim_{\mathbb{C}} {H^0(\widetilde{X}\setminus E, \Omega^{2}_{\widetilde{X}})}/{H^0(\widetilde{X}, \Omega^{2}_{\widetilde{X}})},
\end{equation*}
the analytic weight function agrees with the one associated with the realization of $H^0(\widetilde{X}, \Omega^{2}_{\widetilde{X}})$ inside $H^0(\widetilde{X}\setminus E, \Omega^{2}_{\widetilde{X}})$ (in the sense of Definition \ref{def:Real-intro}) by the natural extensions to meromorphic $2$-forms of the divisorial valuations. This already implies that $\mathbb{H}_*\big(H^0(\widetilde{X}, \Omega^{2}_{\widetilde{X}})\hookrightarrow_{H^0(\widetilde{X}, \mathcal{O}_{\widetilde{X}})} H^0(\widetilde{X}\setminus E, \Omega^{2}_{\widetilde{X}})\big)$ makes sense, and by relying on work of Yau \cite{Yau0, Yau1} and the Independence Theorem \ref{th:Indep-intro}, we can prove that it is independent of the resolution and Stein representative chosen (cf. Corollary \ref{cor:irreg&p_gcategorification} \textit{(b)}). We denote the output of this new construction by $\boxed{\mathbb{H}_*((X, o); \Omega^2)}$ and call it the `\textit{analytic lattice homology of $2$-forms}'. In fact, contrary to the original definition of Ágoston and the first author, this new formulation gives a well-defined invariant for any normal surface singularity (without restrictions on the link or analytic structure), nevertheless, in the rational homology sphere link case it agrees with that: 
\begin{theorem}[= Theorem \ref{th:newvsoldhighdim}]\label{th:newandoldhighdim-intro}
    Let $(X, o)$ be a normal surface singularity with rational homology sphere link. Then
        $   \mathbb{H}_*((X, o); \Omega^2) \cong\mathbb{H}_{{an},*}(X, o)$,
    i.e., our `new' construction generalizes the `old' one.
\end{theorem}

The Gorenstein property for normal surface singularities means that 
the sheaf of $2$-forms  $\Omega^2_{\widetilde{X} \setminus E}$ is holomorphically trivial, hence, the module $H^0(\widetilde{X}\setminus E, \Omega^{2}_{\widetilde{X}})$ is free of rank one above $
H^0(\widetilde{X},\calO_{\widetilde{X}})$. This implies, that the weight function corresponding to any realization is symmetric, moreover, in  this case the submodule $H^0(\widetilde{X},  \Omega^{2}_{\widetilde{X}})$ corresponds to the ideal $H^0(\widetilde{X}, \mathcal{O}_{\widetilde{X}}(-Z_K))$, hence we have $$\mathbb{H}_*((X, o); \Omega^2) \cong \mathbb{SH}_*\big(\phi_*(\cO_{\widetilde{X}}(-Z_K))_o \triangleleft \mathcal{O}_{X, o}\big),$$ 
where $Z_K$ denotes the anticanonical cycle (cf. Corollary \ref{cor:GORCOND}). Analogously to the curve case, we call $\phi_*(\mathcal{O}_{\widetilde{X}}(-Z_K))_o$ the \emph{`conductor ideal'} of $\mathcal{O}_{X, o}$ (cf. Definition \ref{def:CONDINDEF}, see also the \emph{conducteur d'adjonction} of Merle and Teissier \cite{MT}) and prove its independence from the resolution $\phi$ (cf. Proposition \ref{prop:CONDINDEP}).

We also highlight that we can naturally generalize the analytic lattice homology of $2$-forms in the following sense: if we consider instead $1$-forms we can get a categorification of the irregularity
\begin{equation*}
    q(X, o)=\dim_{\mathbb{C}} {H^0(\widetilde{X}\setminus E, \Omega^{1}_{\widetilde{X}} )}/{H^0(\widetilde{X}, \Omega^{1}_{\widetilde{X}})}.
\end{equation*}
Indeed, the natural extensions of divisorial valuations give a realization of ${H^0(\widetilde{X}, \Omega^{1}_{\widetilde{X}})}$ and 
$$\boxed{\mathbb{H}_*((X, o); \Omega^1)}:=\mathbb{H}_*\big( H^0(\widetilde{X}, \Omega^{1}_{\widetilde{X}})\hookrightarrow_{H^0(\widetilde{X}, \mathcal{O}_{\widetilde{X}})} H^0(\widetilde{X}\setminus E, \Omega^{1}_{\widetilde{X}}) \big)$$
is independent of the resolution and Stein representative chosen (cf. Corollary \ref{cor:irreg&p_gcategorification} \textit{(a)}). We call it the `\textit{analytic lattice homology of $1$-forms}' and, by Theorem \ref{th:propertiesbc-intro} \textit{(b)}, it has Euler characteristic $q(X, o)$. We stress that these new invariants are highly nontrivial, we give concrete examples in subsection \ref{ex:whq}.

For both $p=1$ and $2$, we can reformulate the defining modules (by using results from \cite{Yau0, Yau1}) so as to be able to app\-ly the Nonpositivity Theorem \ref{th:nonpos-intro} to the setting of  $\mathbb{H}_*((X, o); \Omega^p)$: Corollary \ref{cor:upperboundcurves} claims that the corresponding $S_n$-spaces are contractible for every $n > 0$. We recall that a similar statement holds in the setting of the topological lattice homology (cf. \cite[Corollary 7.3.13]{NBook}), which asserts that the $S_n$-spaces are connected for every $n > 0$.

An interesting consequence of the Independence Theorem \ref{th:Indep-intro} is that in this analytic setting we immediately get a lattice reduction--type theorem (in the spirit of \cite[Reduction Theorem 3.7]{LN1}), which is especially useful for accelerating computation by reducing the rank of the lattice. Indeed, an equivalent reformulation of the Independence Theorem claims that
if $\mathcal{D}' \subset \mathcal{D}$ are both realizations of the same submodule $N \leq M$, then we can use $\mathcal{D}'$ instead of $\mathcal{D}$ to compute the lattice homology $\mathbb{H}_*(N \hookrightarrow M)$, which, translated to the language of analytic lattice homology of $p$-forms, gives the following result (which is slightly stronger than the original Analytic Reduction Theorem from \cite{AgNe1}):

\begin{theorem}[= Reduction Theorem \ref{prop:redhigh} and Remark \ref{rem:badvertices} (b) -- (d)]\label{th:RedThm-intro}
    Let $(X, o)$ be a normal surface singularity with good resolution $\phi: \widetilde{X} \rightarrow X$ and reduced exceptional divisor $E=\cup_{v \in \mathcal{V}}E_v$. In order to compute $\mathbb{H}_*((X, o); \Omega^p)$ (for $p=1, 2$) it is enough to consider the divisorial valuations corresponding to a subset $\overline{\mathcal{V}} \subset \mathcal{V}$ of irreducible components such that all the connected components of the  full subgraph of the resolution graph $\Gamma_\phi$ supported by $\mathcal{V} \setminus \overline{\mathcal{V}}$ are rational graphs (in the sense of Artin \cite{Artin62, Artin66}). Such a subset $\overline{\mathcal{V}}$ can be, for example, the set of nodes, but usually it is even smaller.
\end{theorem}

Another important  remark is that all of these results remain true for any isolated normal singularity $(X, o)$ of dimension $n\geq 2$. Indeed, we can consider the module of meromorphic $p$-forms (with $p=n-1$ or $p=n$) to obtain lattice homological categorifications $\mathbb{H}_*((X, o);\Omega^p)$ of the irregularity and geometric genus. In fact, a large part of section \ref{s:dnagy} is already written in this more general setting.

\bekezdes \textbf{Lattice homological categorification of the plurigenera of normal surface singularities.}\,

Once we have such a construction for the categorification of the irregularity and the geometric genus of a normal
surface singularity $(X, o)$, it is natural to ask: is there an analogous
theory for the several plurigenera? Recall that in the literature there are different versions:
the $m$-th $L^2$-plurigenus of Watanabe \cite{Watanabe80},
the $m$-th plurigenus of Kn\"oller \cite{Kno73} and
the $m$-th log-plurigenus of Morales \cite{Mor83}.
In this manuscript we consider their common generalizations (cf. \cite[Definition 6.8.57]{NBook})
$$\Delta_{m,n}(X, o)=\dim_{\mathbb{C}} {H^0(\widetilde{X}\setminus E, \mathcal{O}_{\widetilde{X}}(mK_{\widetilde{X}}))}/{H^0(\widetilde{X}, \mathcal{O}_{\widetilde{X}}(m(K_{\widetilde{X}}+E)+nE))} \hspace{5mm} (m \in \mathbb{Z}_{>0},\ n \in \mathbb{Z}\cap[-m, 0]).$$
In fact, the defining submodule is realizable and the corresponding lattice homology, denoted by $\mathbb{H}_*((X, o), \Delta_{m, n})$, is a well-defined categorification of the plurigenus (cf. Corollary \ref{cor:Delta_m,ncategorification}). Moreover, similarly to the previous cases, we can apply the Nonpositivity Theorem \ref{th:nonpos-intro} (cf. Corollary \ref{cor:nonposforDeltamn}) and have symmetric weight function in the Gorenstein case (cf. Remark \ref{rem:gorsymforplurig}). On the other hand, the Reduction Theorem does not have such a geometric reinterpretation as in \ref{th:RedThm-intro} (see Remark \ref{rem:RedThmforplurig}).

We recall, that the authors already had an attempt of categorifying the plurigenera of Gorenstein normal surface singularities in \cite{NSplurig}. The new construction is significantly different from that, seems more natural and has more desirable properties (see Remark \ref{rem:newplurivsold} and part (c) of Example \ref{ex:eltolasplurig}).

We also remark that the results for the $m$-th $L^2$-plurigenus of Watanabe (i.e., for $\Delta_{m,-1}$) extend to higher dimensional isolated normal singularities as well (cf. Corollaries \ref{cor:delta_mcategorification} and \ref{cor:nonposfordeltam})

\subsection{Applications to singularity deformations and Newton nondegenerate singularities}\,

In our new formulation, it is apparent that the analytic lattice homology (or the analytic lattice homology of $n$-forms for isolated normal singularities of dimension $n\geq2$) has deep connections to deformation theory: if the defining quotient module remains the same along a singularity deformation, then the lattice homology will stay the same as well. (By our convention, the studied curve singularities are always reduced, whereas the higher dimensional singularities are isolated and normal.) In the embedded Gorenstein case we can formulate this condition more explicitly using the (previously introduced) conductor ideal:

\begin{define}[= Definition \ref{def:deformations} (1)]
    Let $F=\{(X_t, 0)\subset (\mathbb{C}^m, 0)\}_{t \in (\mathbb{C},0)}$ be a deformation of embedded Gorenstein singularities and let  us denote by $q_t:\cO_{\C^m,0}\to \cO_{X_t,0}$ the natural quotient morphism. $F$  is a {\it `conductor ideal constant deformation'} if the pullback
$q_t^{-1}(\mathcal{C}_{X_t,0})\subset \cO_{\C^m,0}$ is independent of the parameter $t$,
where $\mathcal{C}_{X_t,0}$ is the conductor ideal of $(X_t,0)$.
\end{define}

Applying part \textit{(c)} of Theorem \ref{th:propertiesbc-intro} we automatically obtain the following: 

\begin{theorem}[= Theorem \ref{th:DEFCOND}] 
Let $F=\{(X_t, 0)\subset (\mathbb{C}^m, 0)\}_{t \in (\mathbb{C},0)}$ be a conductor ideal constant deformation of the $n$-dimensional Gorenstein germs $(X_t,0)$ as above. Then
\begin{itemize}
    \item if $n=1$, then the analytic lattice homology $\mathbb{H}_{an,*}(X_t, 0)$ stays constant;
    \item if $n\geq 2$, then  the analytic 
lattice homology of $n$-forms 
$\bH_{*}((X_t,0);\Omega^n)$  stays constant.
\end{itemize}
\end{theorem}

Some prime examples of this class are the `naive'  and `non-naive' deformations of superisolated surface singularities (see subsections \ref{ss:superiso} -- \ref{bek:SIdef} and Proposition \ref{prop:11.5.5}). 

In comparison with other types of deformations, conductor ideal constant deformations are automatically $\delta$-- or $p_g$-constant.
The converse, for curve singularities is false (see, e.g., Example \ref{ex:LIST} (d)), whereas, in the surface case, it is the subject of \cite[Conjecture 11.9.48]{NBook}, which expects that the analytic lattice homology $\mathbb{H}_{an,*}$ is stable along $p_g$-constant deformations whenever the link is a $\mathbb{Q}HS^3$. If we, however, omit the rational homology sphere link condition and work with the more general invariant $\mathbb{H}_*(\ . \ ;\Omega^n)$, than this statement is also false (see part (e) of Example \ref{ex:LIST}, where the link of the generic fiber is not a $\mathbb{Q}HS^3$).

In fact, most of these examples rely on the theory of hypersurface singularities $f:(\mathbb{C}^m, 0) \rightarrow (\mathbb{C}, 0)$ with Newton nondegenerate principal part. These were introduced by Kouchnirenko and Varchenko \cite{Kou, Var} and, by a theorem of Oka (cf. \cite[Theorem 2.1]{Oka}), their embedded topological type depends only on their Newton boundary. Moreover, by a result of Merle and Teissier (cf. \cite[2.1.1 Théorème]{MT}, see also our formulation in Theorem \ref{th:MerleTeissier}), the pullback of the conductor ideal $q^{-1}(\mathcal{C})\triangleleft\mathcal{O}_{\mathbb{C}^m, 0}$ of the germ $f$  can be expressed as the adjoint of the monomial ideal $\mathcal{M}\triangleleft \mathcal{O}_{\mathbb{C}^m, 0}$ generated by monomials with exponent vector in the Newton polytope of $f$.

\begin{cor}[= Theorem \ref{th:comparison} and Corollary \ref{cor:adj=combforNewton}]
    Let $f:(\mathbb{C}^m, 0) \rightarrow (\mathbb{C}, 0)$ be an isolated hypersurface singularity with Newton nondegenerate principal part. Then, its analytic lattice homology (in the $m>2$ case we can think of the analytic lattice homology of $(m-1)$-forms as well) agrees with the symmetric lattice homology $\mathbb{SH}_*({\rm adj}(\mathcal{M}) \triangleleft \cO_{\mathbb{C}^m, 0})$ of the monomial ideal ${\rm adj}(\mathcal{M})$. In particular, $\mathbb{H}_{an,*}(\{f=0\}, 0)$ can be combinatorially computed in terms of the Newton boundary of $f$.
\end{cor}

This observation can be used to construct topologically distinct singularities with identical lattice homology modules (see, e.g., Examples \ref{bek:10.2.1} and \ref{ex:NNSurfaces}) and, also, to answer positively a  question proposed by the first author in \cite[Problem 11.9.52]{NBook}:

\begin{cor}[= Corollary \ref{cor:HanstabforNN}]
    The analytic lattice homology (of $(m-1)$-forms) is constant along an equisingular family of hypersurface singularities (in $(\mathbb{C}^m, 0)$) with Newton nondegenerate principal part associated with a fixed Newton boundary.  
\end{cor}

We also highlight that, by Theorem \ref{th:homdim-intro}, the homological dimension of the analytic lattice homo\-logy of Newton nondegenerate plane curve singularities can be characterized as follows:

\begin{cor}[= Corollary \ref{cor:htopforNNcurves}]
 Let $(C,0)\subset (\mathbb{C}^2, 0)$ be a plane curve singularity with Newton nondegenerate principal part and convenient Newton boundary.
Then we have the following:
   \begin{align*}
      {\rm homdim}_{an}(C, 0):&=\begin{cases}
      -1 & \text{ if }(C, 0) \text{ is regular}\,;\\
          \max\{s\,:\,\mathbb{H}_{an,\, s}(C, 0)\neq 0\} & \text{ if }(C, 0) \text{ is singular}.
      \end{cases}\\
      &=\min\{|\mathcal{D}|\,:\, \mathcal{D} \text{ is a realization of }q^{-1}(\mathcal{C})\triangleleft\mathcal{O}_{\mathbb{C}^2, 0}\}-1\\
      &=\min\{|\mathcal{D}|\,:\, \mathcal{D} \text{ is a realization with  monomial valuations of }q^{-1}(\mathcal{C})\triangleleft\mathcal{O}_{\mathbb{C}^2, 0}\}-1
  \end{align*}
  where $q^{-1}(\mathcal{C})$ is the pullback of the conductor ideal $\mathcal{C}$ of $(C,o)$ along the quotient $q:\mathcal{O}_{\mathbb{C}^2, 0} \rightarrow \mathcal{O}_{C, 0}$.
\end{cor}

\subsection{Generalizations of the construction}\,

Remarkably, the conditions of the Independence Theorem \ref{th:Indep-intro} can be weakened. Theorem \ref{th:propertiesbc-intro} \textit{(c)} already implies that we do not necessarily need $M$ to be finitely generated, it is enough to have $M/N$ finite dimensional over $k$. Similarly, the Noetherian condition on $\mathcal{O}$ is also required more for convenience (we stress, however, that Rees defined integrally closed submodules in this exact setting). Moreover, instead of counting $k$-dimensions, we could also use the length (of Artinian  quotient modules) over (a Noetherian ring) $\mathcal{O}$ to define the height functions and, hence, obtain a well-defined lattice homology theory (though, this might give a different outcome for $k$-algebras than the given construction).  

On the other hand, we can also leave the realm of rings and modules and consider only a $k$-bilinear `multiplication' map $\bullet: M_1 \times M_{-1} \rightarrow M_0$ on three (possibly) different $k$-vector spaces. (Our standard setting corresponds to the substitutions $M_1=\mathcal{O}$, $M_{-1}=M_0=M$ and $\bullet$ being the ring action.)  The corresponding versions of the Independence Theorem are given in Theorems \ref{th:Indep-MGUF} and \ref{th:Indep-MGF}. They are used to produce a reformulation of the definition of the \emph{equivariant} analytic lattice homology after \cite{AgNeIII} of normal surface singularities with rational homology sphere link (cf. Corollary \ref{cor:equivariant}) and the corresponding Reduction Theorem (cf. Theorem \ref{th:REDANh}).

Yet another approach to generalization arises from the following property of our new construction:

\begin{theorem}[= Theorem \ref{th:properties} \textit{(a)}]
    Given  a Noetherian $k$-algebra $\mathcal{O}$, a finitely generated $\mathcal{O}$-module $M$ and $N\leq M$ realizable, $\mathbb{H}_*(N \hookrightarrow_{\mathcal{O}} M)$ can be computed on a finite rectangle: for any  realization $\mathcal{D}$ of $N$ and  large enough lattice point $0\ll d \leq \infty$ (depending on $\mathcal{D}$) all the 
inclusions $S_n \cap R(0, d) \hookrightarrow S_n$ are homotopy equivalences. (In fact, this result holds in each of the previously presented settings.)
\end{theorem}

Thus,
from the purely computational point of view, it is enough to consider the multifiltrations $\mathcal{F}_\mathcal{D}$ and $\mathcal{F}_{\mathcal{D}}^M$ of (\ref{filtrations-intro}) only on this finite rectangle $R(0, d)$. Consequently, instead of the discrete valuations $\{\mathfrak{v}_v\}_v$ of $\mathcal{D}$ one can consider valuative maps, so called \emph{`quasi-valuations'} $\mathfrak{q}_v:\mathcal{O} \rightarrow \mathbb{N}_{d_v}$ (cf. Definition \ref{def:quasi-val}) with values in the ordered commutative monoid $(\mathbb{N}_{d_v}, \ast_{d_v})=\big(\{0, 1, \ldots, d\}, \ast_{d_v}\big)$  (with $d_v\in \mathbb{Z}_{>0}$) having the usual ordering and the unusual multiplication $a \ast_{d_v} b=\min \{ a+b, d_v\}$ (cf. Definition \ref{def:Nd}). 

Such quasi-valuations can be obtained from discrete valuations by simple post-composition with the $\min_{d_v}:\mathbb{N} \rightarrow \mathbb{N}_{d_v}, \ a \mapsto \min\{a, d_v\}$ map (cf. Example \ref{ex:levagottv}), or, more interestingly, by restricting to quotients: if $\mathfrak{v}_v(\mathcal{I})=d_v$ for some discrete valuation $\mathfrak{v}_v:\mathcal{O} \rightarrow \overline{N}$ and ideal $\mathcal{I} \triangleleft \mathcal{O}$, then the well-defined map $\mathfrak{v}_v\big|_{\mathcal{O}/\mathcal{I}}:\mathcal{O}/\mathcal{I} \rightarrow \mathbb{N}_{d_v}, \ f + \mathcal{I} \mapsto \mathfrak{v}_v(f+\mathcal{I})$ is a quasi-valuation (cf. Example \ref{ex:quasirest}).

Analogously to the case described in subsection \ref{ss:main}, we can define \emph{extended} quasi-valuations (cf. Definition \ref{def:extquasi}), quasi-realizability (cf. Definition \ref{def:REAL-q}) and the adjusted Independence Theorem holds true in this context as well (cf. Theorem \ref{th:IndepMod-q}). There are two main advantages of working with this quasi-valuative setup (cf. Remark \ref{rem:robustness-q}):
\begin{itemize}
    \item quasi-realizability is a more robust property with respect to module operations than realizability: for example, if $N$ is quasi-realizable in $M$, then so is $N/N$ in $M/N$ (which is false for realizability or integral closedness, see, e.g., the example in Remark \ref{rem:imageofintclosed});
    \item quasi-valuations are, in a sense, finite objects, they can be determined by finite bases and the quasi-values of their elements (cf. Remarks \ref{rem:constructionofq-val} and \ref{rem:constructionofextq-val}), thus, they are much easier to construct than (non-monomial) discrete valuations (cf. Propositions \ref{prop:quasivalcond} and \ref{prop:extquasivalcond}).
\end{itemize}
 Consequently, they provide a more flexible 
 setting for further exploration, and, especially, to test  Conjecture \ref{conj:elso-intro} (see, e.g.,  Example \ref{ex:curvewith1qv}).

\subsection{Structure of the article}\,

The structure of the article 
is the following.
Section \ref{s:Prem1} contains the general definition of lattice homology and graded root
 associated with a lattice and weight function, here we follow \cite{NOSz,NGr,Nlattice}.
 In section \ref{ss:CombLattice} we present the combinatorial backbone of the theory: several  useful reduction arguments, which depend only on
combinatorial properties, are separated here. We especially highlight the requirements of the Categorification Theorem of Ágoston and the first author. 

The main results are discussed in section \ref{s:4}, where we introduce the concept of \emph{extended} discrete valuations, \emph{realizable} submodules and state the Independence Theorem. We present the main structural properties of the lattice homology of realizable submodules, but the more technical proofs are postponed in order to help easier understanding of the concept.
In section \ref{s:ideals} we discuss the case of ideals, where realizability agrees with integral closedness and the weight function is automatically symmetric. In the last part of this section we present the first concrete examples and computations. We introduce `integrally reduced' algebras in section \ref{ss:ARTIN}, discuss their well-definedness and study their symmetric lattice homology modules.

In sections \ref{s:deccurves} and \ref{s:dnagy} we present concrete applications of our construction to provide ca\-tegorifications of numerical invariants of complex analytic singularities. First, in the curve case, we show a reinterpretation of the analytic lattice homology as corresponding to a realizable submodule. Secondly, we provide a natural categorification of the geometric genus and irregularity of isolated singularities with dimension $\geq 2$. This extends the ori\-ginal definition of the analytic lattice homology of Ágoston and the first author and comes with strong Reduction-- and Nonpositivity Theorems. Finally, we treat the various plurigenera of normal surface singularities: since they can be defined as codimensions of realizable submodules, we can associate with them natural lattice homological categorifications. We compare these with the lattice \emph{co}homological invariants introduced previously by the authors  in \cite{NSplurig}. 

The following section \ref{s:ND} discusses the case of monomial ideals and Newton nondegenerate singularities. 
It turns out that for this latter class, the analytic lattice homology can be computed combinatorially from the Newton boundary. More precisely, for such analytic germs there are two naturally ocurring weight functions (one based on the analytic data, the other one on combinatorial valuations of the Newton diagram), which, in general, do not agree, but give the same lattice homology. In fact, the question of their equivalence was the main starting point of our research.

Section \ref{s:homdim} discusses the newly introduced \emph{(lattice) homological dimension} of a realizable submodule $N \leq M$ and its relation to the number of (extended) valuations required to obtain a realization of $N$. We present a general statement providing bounds for both values, and, in the case of integrally closed monomial ideals in the two variable polynomial ring, we prove that these determine each other.

In section \ref{s:NNSI} we discuss the relation of our analytic invariants with singularity deformations. The case of superisolated surface germs is particularly highlighted. We introduce the \emph{conductor ideal constant} Gorenstein singularity deformations and present multiple examples and non-examples. 

The proof of the Independence Theorem is located in section \ref{s:proof}. Its generalized forms are also presented here, which are able to describe the equivariant analytic lattice homology for surfaces of Ágoston and the first author \cite{AgNeIII} as well. Section \ref{s:nonpositivity} contains the proof of the Nonpositivity Theorem.

In section \ref{ss:quasi} we describe a more pragmatic approach towards our lattice homology contruction via the so-called \emph{quasi-valuations}. These are in some sense finite quotients of discrete valuations, which provide a more flexible setting and, a priori, a larger domain for our definitions.

In the Appendix we compare relevant algebraic properties of submodules: realizability and integral closedness in the sense of Rees (respectively, in the sense of Gaffney). Throughout the manuscript we also
present several examples and problems\,/\,conjectures regarding the new theory. In the final subsection \ref{s:QOP} we list and collect further related questions and open problems.

\subsection{Acknowledgements.} Parts of this manuscript constitute the backbone of the second author's PhD dissertation. He is especially grateful to his referees Tamás László and András Stipsicz, who read through his thesis and thus parts of a previous version of this paper as well. The authors are also thankful to Alexander A. Kubasch for the inspiring previous work, moreover to András Szűcs for suggesting the use of Quillen's Fiber Theorem in the proof of the main statement.

\section{Basic properties of lattice homology and the graded root}\label{s:Prem1}

\subsection{The lattice homology module associated with  a weight function \cite{NOSz,Nlattice}}\, \label{ss:latweight}

The lattice homology construction 
associates
 a bigraded $\mathbb{Z}[U]$-module $\bH_*$ to a
triple $(\mathbb{Z}^r, \{e_v\}_v, w)$ consisting of a lattice, a fixed basis $\{e_v\}_v$ 
of it and a weight function $w$,
see the original sources \cite{ NOSz,NGr,Nlattice,NBook} or paragraphs \ref{bek:211}--\ref{9complexb} below.   
[Here $U$ is a free variable; hence, in fact, 
$\bH_*$ is a bigraded $\mathbb{Z}$-module  with a specified bihomogeneous endomorphism, denoted by $U$.] For the construction in a more general setting (which will not be used in the present manuscript) see \cite{AgNePoset}.

$\bH_*$ has a direct sum decomposition
$\oplus_{q\geq 0}\bH_q$ of graded ${\mathbb Z}[U]$-modules, 
each $\bH_q$  is $2\mathbb{Z}$-graded, with the $U$-action being homogeneous of degree $(-2)$;
we use this  convention to  sustain
compatibility with the Heegaard Floer (Link) theory (for more on these connections see  \cite{NOSz,NFilt,NFilt2,Zemke}). For the sake of notational simplicity we introduce the following distinguished graded $\Z[U]$-modules.

\bekezdes\label{9zu1} {\bf $\Z[U]$-modules.} We modify the usual grading of the polynomial
ring $\Z[U]$ in such a way that the degree of $U^m$ is $-2m$ ($m\geq 0$).
Besides $\calt^-_0:=\Z[U]$, considered as a $\Z[U]$-module with this grading,  we will consider for $n\geq 1$ the modules
$\calt_0(n):=\Z[U]/(U^n)$, too, with the induced grading. Moreover, we can consider the $\Z[U]$-module
$\Z[U,U^{-1}]$ with the above grading and  denote by $\calt_0^+$
its quotient by the submodule  $U\cdot \Z[U]$.

More generally, for any graded $\Z[U]$-module $P$ with
$d$-homogeneous part $P_d$, and  for any  $k\in\Z$,   we
denote by $P[k]$ the same module graded in such a way
that $P[k]_{d+k}=P_{d}$. Then set $\calt^-_k:=\calt^-_0[k]$ and
$\calt_k(n):=\calt_0(n)[k]$. Hence, for $m\in \Z$, we have 
$\calt_{-2m}^-\cong\Z\langle U^{m}, U^{m+1},\ldots\rangle$ as a $\Z$-module, 
whereas
$\calt_{2m}^+\cong\Z\langle U^{-m}, U^{-m-1},\ldots\rangle$.\\

The lattice homology module can be defined in two different but equivalent ways, one highlighting its geometric, the other its algebraic nature. The first construction,  via 
the filtration of subspaces $\{S_n\}_n$ of $\R^r$, will be reviewed in the next paragraphs.
For the other one, via chain complexes, see  subsection \ref{ss:chain}.

\bekezdes\label{bek:211} {\bf The lattice and the cubical decomposition.}
First we consider a finitely generated free $\mathbb{Z}$-module, with a fixed basis
$\{e_v\}_{v\in\mathcal{V}}$, denoted by $\mathbb{Z}^r$
(hence $r=|\mathcal{V}|$).
There is a natural
partial ordering  $\leq $ of $\bZ^r$ (and of  $\mathbb{R}^r$), defined coordinatewise,
induced by the fixed basis:  for lattice points $\ell_1,\ell_2\in \mathbb{Z}^r$ with $\ell_j=\sum_v \ell_{jv}e_v$ ($j\in\{1,2\}$)
 one says that $\ell_1\leq \ell_2$ if $\ell_{1v}\leq \ell_{2v}$ for all $v\in \mathcal{V}$. We define the minima and the maxima of two lattice points $\ell_1= \sum_v \ell_{1, v} e_v$ and $\ell_2=\sum_v \ell_{2, v} e_v$ as $\min\{\ell_1, \ell_2\} = \sum_v \min\{\ell_{1, v}, \ell_{2, v}\} e_v$ and $\max\{\ell_1, \ell_2\} = \sum_v \max\{\ell_{1, v}, \ell_{2, v}\} e_v$.
The support of a cycle $\ell=\sum \ell_ve_v$ is defined as  ${\rm supp}(\ell)=\cup_{\ell_v\not=0}e_v$. 

The space 
$\mathbb{Z}^r\otimes \mathbb{R}=\mathbb{R}^r$ has a natural cellular decomposition into cubes. The
set of zero dimensional cubes consists of  the lattice points of
$\mathbb{Z}^r$. Any $\ell\in \mathbb{Z}^r$ and subset
$I\subset \mathcal{V}$ of
cardinality $q$  define  a $q$-dimensional cube  $(\ell, I)\subset \mathbb{R}^r$, which has its
vertices in the lattice points $\{\ell+e_{I'}\}_{I'}$, where
$I'$ runs over all subsets of $I$ and $e_{I'}$ denotes $\sum_{v\in I'}e_v$. We regard these cubes as {\it closed} cubes of ${\mathbb R}^r$, which contain their boundaries as well. We will often use the notation $\square_q$ for cubes of dimension $q$, whereas $\square_q^\circ$ will denote their relative interior.

\bekezdes \textbf{The weight function.} \label{bek:wf}
A `{\it weight function}' 
$w_0:\mathbb{Z}^r\to \mathbb{Z}$ is an integer valued function 
on the lattice  with the property  that 
\begin{equation}
    \text{ the set } w_0^{-1}(\,(-\infty,k]\,) \text{ is finite 
for any } k\in\mathbb{Z}. 
\end{equation}\label{weightfctn}
Once $w_0$ is fixed, 
we endow each cube with a weight: for any $q$-cube $\square_q$  we set  
\begin{equation}\label{eq:9weight}
\mbox{$ \ w_q(\square_q)=\max \{ w_0(\ell) \, : \, \ell$ is a vertex of $\square_q\}$}.\end{equation}

\noindent [For the more general notion of {`a set of compatible weight functions'} see \cite[3.1.4 Definition]{Nlattice}.] \\
When the dimension of the cube  is clear from the context we will omit the index $q$ from our notation. We also often use the notation $w$ for the set $\{w_q\}_q$ of weight functions as well.

\bekezdes\label{9complexb} {\bf The $S_n$ spaces and the lattice homology.}
For each $n\in \mathbb{Z}$ we
define $S_n\subset \mathbb{R}^r$ as the union of all
those closed cubes $\square$ (of any dimension), which have $w(\square)\leq
n$.
These spaces give an increasing closed cubical complex filtration: $\ldots \subset S_{n-1} \subset S_n \subset S_{n+1} \subset \ldots \subset \mathbb{R}^r$. 
Clearly, $S_n=\emptyset$, whenever $n<m_w:=\min\{w_0(\ell)\,:\, \ell\in \mathbb{Z}^r\}$ (we will sometimes also denote this by $\min w_0$). 

Then for any  $q\geq 0$, set
$$\mathbb{H}_q(\mathbb{R}^r,w):=\oplus_{n\geq m_w}\, H_q(S_n,\mathbb{Z})\ \ \mbox{and}\ \
\mathbb{H}_{{\rm red},q}(\mathbb{R}^r,w):=\oplus_{n\geq m_w}\, \widetilde{H}_q(S_n,\mathbb{Z}).$$
We write $\bH_*$ for $\oplus_{q\geq 0}\bH_q$, we call $q$ the 
\textit{`homological grading'}. On the other hand, 
each $\mathbb{H}_q$ is $2\mathbb{Z}$-graded, the summand of 
$(-2n)$-homogeneous elements $\mathbb{H}_{q, -2n}:=(\mathbb{H}_q)_{-2n}$ is  $H_q(S_n,\mathbb{Z})$.
This is called the \textit{`weight grading'}.  

Also, $\mathbb{H}_q$ is a $\mathbb{Z}[U]$-module: the $U$-action is given by
the homological maps $H_q(S_{n},\mathbb{Z})\to H_q(S_{n+1},\mathbb{Z})$ induced by the inclusions $S_n \hookrightarrow S_{n+1}$.
The same is true for $\mathbb{H}_{{\rm red},*}$.
 Moreover, for
$q=0$, any fixed `base-point' $\ell_w\in S_{m_w}\subset S_n$ provides a
simultaneous augmentation, i.e., for each $n \geq m_w$ a splitting of form
 $H_0(S_n,\mathbb{Z})=
\mathbb{Z}\oplus \widetilde{H}_0(S_n,\mathbb{Z})$; hence an augmentation of the
$\mathbb{Z}[U]$-modules
$$\mathbb{H}_0\simeq
\calt^-_{-2m_w}
\oplus \mathbb{H}_{{\rm red},0}=(\oplus_{n\geq m_w}\mathbb{Z})\oplus (
\oplus_{n\geq m_w}\widetilde{H}_0(S_n,\mathbb{Z}))\ \ \mbox{and} \ \
\mathbb{H}_*\simeq
\calt^-_{-2m_w}
\oplus \mathbb{H}_{{\rm red},*}.$$

Although
$\mathbb{H}_{{\rm red},*}(\mathbb{R}^r,w)$ has finite $\mathbb{Z}$-rank in any fixed
weight degree, in general --- without certain additional properties of $w_0$ --- it is not
finitely generated over $\mathbb{Z}$, in fact, not even over $\mathbb{Z}[U]$.

\bekezdes\label{9SSP} {\bf Restrictions.} Assume that $T\subset \mathbb{R}^r$ is a subspace
of $\mathbb{R}^r$ consisting of some union of (closed) cubes.  For any $q\geq 0$ define $\mathbb{H}_q(T,w)$ as
$\oplus_{n\geq\min{\{w_0|_T\}}} H_q(S_n\cap T,\mathbb{Z})$. It naturally has a graded $\mathbb{Z}[U]$-module
structure.  Moreover, the inclusion maps $S_n\cap T\hookrightarrow S_n$  induce a natural bigraded
$\mathbb{Z}[U]$-module homomorphism
$$\mathbb{H}_*(T,w)\to  \mathbb{H}_*(\mathbb{R}^r,w) \ \ \ \text{(of bidegree $(0,0)$)}.$$
\noindent In some cases it can happen that the weights are only defined for cubes belonging to $T$.
In our applications,  $T$ (besides the trivial $T=\mathbb{R}^r$ case) will be one of the following:\\
(i)  the first quadrant $(\mathbb{R}_{\geq 0})^r$;\\
(ii) a rectangle $R(0,d)=\{x\in \mathbb{R}^r\,:\, 0\leq x\leq d\}$ for some lattice point $d \in (\mathbb{Z}_{\geq 0})^r$, or \\
(iii) a path of composed edges in the lattice starting from $0$: a sequence $\gamma=\{x_i\}_{i=0}^t$ of lattice points is called a \textit{path from $x_0$ to $x_t$} if 
$x_i\not=x_j$ for $i\not=j$ and $x_{i+1}=x_i + \epsilon_i e_{v(i)}$, with $\epsilon_i \in \{\pm1\}$ and $v(i)\in \mathcal{V}$ for every
$0\leq i<t$ ($\gamma$ is called \textit{increasing} if $\epsilon_i=1$ for all $0 \leq i <t$); in this case  $T$ is the
union of $0$-cubes marked by the points $\{x_i\}_i$ and of the
segments of type  $[x_i,x_{i+1}]$ (i.e., $1$-cubes of form $(x_i, \{e_{v(i)}\})$). We often also assume that $x_0=0$.

\bekezdes \label{bek:eu}{\bf The Euler characteristic of $\mathbb{H}_*$.}
Let $T$ be  as in \ref{9SSP} and assume that each $\mathbb{H}_{{\rm red},q}(T,w)$ has finite $\mathbb{Z}$-rank for every $q \in \mathbb{Z}_{\geq 0}$.
(This happens, e.g., when $T$ is a finite complex.)
{Whenever the following formula makes sense,} we define the Euler characteristic of $\mathbb{H}_*(T,w)$ as
$$eu(\mathbb{H}_*(T,w)):=-\min\{w_0|_T\} +
\sum_q(-1)^q{\rm rank}_{\mathbb{Z}}\left(\mathbb{H}_{{\rm red},q}(T,w)\right).$$
By \cite[Theorem 2.3.7]{NJEMS}, for $T=R(0, d)$ (for some $d\in(\mathbb{Z}_{\geq 0})^r$): 
$eu(\bH_*(T,w))=
\sum_{\square_q\subset T} (-1)^{q+1}w(\square_q)$.

\bekezdes \label{bek:LC}{\bf Cohomological formulation.} Given a lattice $\mathbb{Z}^r$ and a weight function $w$, we could consider the singular \textit{co}homologies of the cubical complexes $S_n$ (instead of their singular homologies), which would yield the \textit{lattice cohomology} module $\mathbb{H}^*(\mathbb{R}^r, w)=\oplus_{q \geq 0} \oplus_{n \geq m_w} H^q(S_n, \mathbb{Z})$. Obviously, this reverses the direction of the $U$-action, which we treat by changing the weight grading to its negative, to ensure that the action is homogeneous of degree $(-2)$ similarly to the case of Heegaard Floer theory. On the other hand, most of the results remain the same in the cohomological and homological formulation (e.g. the Euler characteristic formulae, well-definedness and lattice reduction results), specifically those, which rely only on  the homotopy types of the $S_n$-spaces. Most of the original literature (e.g. \cite{Nlattice, NBook}) is written in the cohomological convention with the hope of using its extra ring structure, however, some of the recent results have nicer formulation in the homological setting (e.g. \cite{NFilt, NFilt2}). We opted to use this latter to be able to use a non-vanishing tensor product operation in subsection \ref{ss:KUNNETH}.

\subsection{Example: the topological lattice homology of a complex normal surface singularity $(X,o)$ \cite{NOSz, Nlattice, NBook}}\label{toplc} [We remark, that the original sources use the cohomological convention.] 

 We consider a   good
  resolution $\phi: \widetilde{X} \rightarrow X$ of $(X,o)$ and we assume that 
  {\it its link $L_X$ is a
rational homology sphere}. Let $\mathcal{V}$ be the indexing set of irreducible exceptional curves
$\{E_v\}_{v \in \mathcal{V}}$ and set  $r:=|\mathcal{V}|$.
Then the lattice $L:=H_2(\widetilde{X}, \bZ)$ is isomorphic to $\bZ^r$ with
a fixed basis consisting of the homological fundamental classes of the exceptional curves $\{E_v\}_v$ (also denoted by $\{E_v\}_v$). It carries a negative definite intersection form $(\,,\,)$,  which naturally extends to $L\otimes \bQ$.
Consider the anticanonical cycle $Z_K:=-c_1(\Omega_{\widetilde{X}}^2)\in H^2(\widetilde{X},\bZ)\simeq{\rm Hom}_{\bZ}(L,\bZ)
\simeq
\{\ell\in L\otimes \bQ\,:\, (\ell, L)\subset  \bZ\}$.
  Then, the Riemann--Roch expression
 $\chi: \ell \mapsto -(\ell, \ell-Z_K)/2$ defines a weight function
$w_{top,0}(\ell)=\chi(\ell)$ on $\mathbb{Z}^r$, hence 
a  weight for any cube via identity (\ref{eq:9weight}).

 The  $\bZ[U]$-modules $\bH_*(\bR^r,  w_{top}  )$ and
$\bH_{red,*}(\bR^r,  w_{top}  )$ obtained in this way are
called the {\em topological lattice homologies}.
One proves that they
depend only on the (diffeomorphism type of the) link $L_X$ and, thus, they are  independent of the choice of the good resolution $\phi$   (see \cite[Proposition 3.4.2]{Nlattice} or \cite[Proposition 11.1.24]{NBook}).
As  $\bH_{red,*}(\bR^r,  w_{top}  )$
 has finite $\bZ$-rank, its Euler characteristic is well-defined. By \cite[Theorem B]{NJEMS},\- it satisfies
$eu(\bH_*(\bR^r,  w_{top}  ))=
\sw_{\sigma_{can}}(L_X)-(Z_K^2+|\mathcal{V}|)/8$, where  $\sw$ denotes the Seiberg--Witten invariant ${\rm Spin}^c(L_X)\to \bQ$, which associates a rational number
$\sw _{\sigma}(L_X)$  to each spin$^c$-structure $\sigma$ of the link. In other words, {\it 
$\bH_*(\bR^r,  w_{top}  )$
is the categorification of $\sw_{\sigma_{can}}(L_X)$}    (normalized by $(Z_K^2+|\mathcal{V}|)/8$), where $\sigma_{can}$ is the canonical spin$^c$ structure
of $L_X$.

It turns out that in this construction the set of cubes of $\bR^r$ can be replaced by
those from the first quadrant $(\bR_{\geq 0})^r$ providing the very  same homology.
In the sequel we will write $\bH_{top, *}(X, o)$} for \ $\bH_*(\bR^r,  w_{top}  )=\bH_*((\bR_{\geq0})^r,  w_{top}  )$.

We must also mention that, in fact, there are more variants of the topological lattice homology according to different spin$^c$-structures. These correspond to different choices of the characteristic element 
$k \in Char:=\{\ell'\in L\otimes \bQ\,|\, (\ell'+\ell,\ell)\in 2\bZ \ \text{for all $l\in L$}\}$
used in the  weight function $\chi_k(\ell) := -(\ell, \ell+k)/2$ of these theories.  They provide bigraded  $\bZ[U]$-modules with similar properties, cf. 
\cite{NOSz, Nlattice, NBook}. 

\subsection{The homological graded root}\label{ss:gradedroot}\

Consider the setting of subsection \ref{ss:latweight}. 
The graded $\bZ[U]$-module $\bH_0$ has an improvement, which can also be used as a 
pictorial presentation of $\bH_0$. 
It is the {\it homological graded root} $\mathfrak{R}$, a $\bZ$-graded connected graph, which describes the containment relations of the different connected components in the tower of spaces $\{S_n\}_{n \in \mathbb{Z}}$. We will use it mainly for presenting examples in a concise yet geometric way.

\begin{define}\label{bek:gradedroot}
 The vertex set $\mathfrak{V}$ of $\mathfrak{R}$ is $\bZ$-graded, 
 $\mathfrak{V}=\sqcup _n
\mathfrak{V}_n$. The vertices  $\{v_{-n,i}\}_i\in \mathfrak{V}_{-n}$ correspond bijectively to the connected components $\{S_{n,i}\}_i$ of $S_n$. The edges are defined as follows: if two components $S_{n,i}$ and $S_{n+1,j}$
of $S_n$ and $S_{n+1}$  satisfy $S_{n,i}\subset S_{n+1,j}$ then we connect $v_{-n,i}$ and $v_{-(n+1),j}$ by an edge. 
\end{define}
From $\mathfrak{R}$ one can read  $\bH_0$ as follows: each vertex of degree $(-n)$
of $\mathfrak{R}$
generates a free $\bZ$-summand in $\mathbb{H}_{0, -2n}=(\bH_0)_{-2n}=H_0(S_n,\bZ)\cong\oplus_i H_0(S_{n, i},\bZ)$, and the $U$-action is marked by the edges. On the other hand, $\mathbb{H}_0$ in general does not characterize $\mathfrak{R}$, see e.g. \cite[3.6 Examples (b)]{NOSz}. 
For more see, e.g., \cite{NGr,Nlattice,NBook} for the cohomological version, and 
\cite{NFilt,NFilt2} for the homological version.

\section{Combinatorial properties of lattice homology with split  weight function} \label{ss:CombLattice}

The weight function $w_0$ plays a crucial role in the definition of lattice homology. Indeed, a priori we are working on the contractible space $\mathbb{R}^r$, yet, the weight function determines the nontrivial level sets $\{S_n\}_n$. 
In the case of the topological lattice homology (cf. Example \ref{toplc}) $w_0$ is given by a Riemann-Roch expression,
however in several analytic and combinatorial situations it has a different, specific structure:
it is a sum of two  (geometrically defined) functions, one of them  increasing, the other one  decreasing and they satisfy some key additional properties, too, see below. 
This decomposition usually reflects (or, is the consequence of)  some kind of
(e.g., Serre--, Gorenstein--) duality.   In this section  we present  the general combinatorial setup  of such weight functions.
In the first subsection we follow \cite{AgNe1,NBook} (where the \emph{co}homological version is treated), while the second subsection contains some new technical results whose proofs the reader might skip at first reading to get to the description of the main constructions faster.

\subsection{The combinatorial setup.}

\bekezdes\textbf{Split weight function.} \label{bek:comblattice}
Consider the lattice  $\bZ^r$ with a fixed basis $\{e_v\}_{v\in\mathcal{V}}$, $|\mathcal{V}|=r$.
Fix also  an element
$d\in (\bZ_{\geq 0})^r\cup\{ \infty\}$ and
consider the rectangle $R:=R(0,d)$.
Furthermore,  assume that
to each $\ell\in R\cap \bZ^r$ we assign

(i)   an integer $h(\ell)$ such that $h(0)=0$ and $h(\ell+e_v)\geq h(\ell)$
for any $v \in \mathcal{V}$, $\ell, \ell+e_v\in R \cap \bZ^r$\\
\indent \hspace{3mm} (i.e., $h$ is increasing); 

(ii) and an integer $h^\circ (\ell)$ such that $h^\circ (\ell+e_v)\leq  h^\circ (\ell)$
for any $v \in \mathcal{V}$ and $\ell, \ell+e_v\in R \cap \bZ^r$\\
\indent \hspace{4mm} (i.e., $h^\circ$ is decreasing).

\noindent
E.g., once  $h$ and $d$ are fixed with (i) and $d\in(\mathbb{Z}_{\geq 0})^r$ finite,
a possible choice for $h^\circ $ is
 $h^{sym}_d$, the symmetrization of $h$ with respect to $d$, i.e., $h^{sym}_d(\ell):=h(d-\ell)$ for every $\ell \in R \cap \mathbb{Z}^r$. 
 We sometimes call $h$ and $h^\circ$ the \textit{`height functions'}.

We  consider the set of cubes of $R$ as in \ref{bek:211} and
the weight function
\begin{equation} \label{eq:w_0def}
   w_0:R\cap \mathbb{Z}^r\rightarrow \bZ\ \ \mbox{defined by} \ \ w_0(\ell):=h(\ell)+h^\circ (\ell)-h^\circ (0).
\end{equation}
\noindent Clearly  $w_0(0)=0$. Moreover, we define  the other functions $\{w_q\}_{q \geq 1}$ by   (\ref{eq:9weight}).
 They define the lattice homology module $\bH_*(R,w)$.
 In particular, if $\mathbb{H}_{{\rm red},*}(R, w)$ has finite rank (which in the $d < \infty$ case is automatic), we have a well-defined
Euler characteristic  $eu(\bH_*(R,w))$ as well (see paragraph \ref{bek:eu}).

\bekezdes \textbf{Key properties.} We will focus on pairs  $(h,h^\circ)$ which satisfy
certain additional combinatorial properties.

\begin{define}\label{def:matroidh}
We say that $h$ satisfies the  {\it `matroid rank inequality'}\, if
 \begin{equation}\label{eq:matroid}
 h(\ell_1)+h(\ell_2)\geq h(\min\{\ell_1,\ell_2\})+h(\max\{\ell_1,\ell_2\}), \ \ \ell_1,\ell_2\in R \cap \bZ^r.
 \end{equation}
 Note that the matroid rank inequality and the increasing property  of
 $h$ implies the following so-called  {\it `stability property'}
 \begin{equation}\label{eq:stability}
 h(\ell)=h(\ell+e_v)\ \ \Rightarrow\  \ h(\ell+\bar{\ell})=h(\ell+\bar{\ell}+e_v);
 \end{equation}
valid for any $\ell, \ \bar{\ell} \in R \cap \mathbb{Z}^r$ such that $\bar{\ell}\geq 0$,
 $\bar{\ell}\ngeq e_v$ (i.e., $e_v \notin {\rm supp}(\bar{\ell})$) and $\ell+ \overline{\ell} + e_v \in R \cap \bZ^r$.
\end{define}

\begin{example} \label{ex:multifilt}
Assume that $V$ is a  $k$-vector space (with $k$ some underlying field) endowed with  a rank $r$ decreasing  multifiltration $\mathcal{F}$ on it, i.e., for any lattice point $\ell \in \Z^r$ there is assigned a linear subspace $\mathcal{F}(\ell) \subset  V$, such that $\mathcal{F}(\ell) \subset \mathcal{F}(\ell')$ if $\ell \geq \ell'$. We assume that $\mathcal{F}(0)=V$.
If we require that
$\frh(\ell) := \dim_k(V/\mathcal{F}(\ell))<\infty$ for any $\ell \in \Z^r$, then we can define  the \textit{Hilbert function} of $\mathcal{F}$ as $\ell\mapsto \frh(\ell)$.
It satisfies $\frh(0)=0$ and it is increasing: $\frh(\ell) \geq \frh(\ell')$ if $\ell \geq \ell'$.
However, in this generality the matroid rank inequality is usually not satisfied.
\end{example}

On the other hand, we have the following simple criterion:

\begin{lemma}\label{lemma:matroid} 
    Suppose that $V$ and  $\{\mathcal{F}(\ell)\}_{\ell\in \mathbb{Z}^r}$ are as in Example \ref{ex:multifilt}. If the
  multifiltration satisfies the identity
    \begin{equation}\label{eq:FaFb=Fmax}
        \mathcal{F}(\max\{\ell_1,\ell_2\}) = \mathcal{F}(\ell_1)\cap \mathcal{F}(\ell_2), \ \ \mbox{for all} \  \ell_1, \ell_2 \in \bZ^r,
    \end{equation}
 then the Hilbert function
 $\ell \mapsto \frh(\ell)=\dim_k( V/\mathcal{F}(\ell))$
satisfies  the matroid rank inequality.
\end{lemma}

\begin{proof}
    Left to the reader.
\end{proof}

\begin{example}\label{ex:Z^rgrading}
    ($V$, $\mathcal{F}$ as in Example \ref{ex:multifilt}.) If the multifiltration $\mathcal{F}$ corresponds to a  $\bZ^r$-grading  $V \cong \oplus_{a\in\bZ^r}V_{a}$ of the vector space , i.e., $\mathcal{F}(\ell):= \oplus _{a\geq \ell}V_a$ for all $\ell \in (\mathbb{Z}_{\geq 0})^r$, then it naturally satisfies condition (\ref{eq:FaFb=Fmax}).
\end{example}

\begin{define}\label{def:COMPGOR}
 We say that the pair $(h,h^\circ)$      satisfies the {\it 
 `Combinatorial Duality  Property'} in $R(0, d)$ (often abbreviated as CDP) if
$h(\ell+e_v)-h(\ell)$ and $h^\circ (\ell+e_v)-h^\circ (\ell)$  are not   simultaneously  nonzero
for any $\ell,\, \ell+e_v\in R\cap\bZ^r$. Furthermore, whenever $d<\infty$, 
 we say that the height function $h$  satisfies  the CDP  in $R(0,d)$ if
 the pair $(h,h^{sym}_d)$ satisfies  it.
\end{define}

Examples of pairs $(h, h^\circ)$ satisfying the CDP can be found in \cite{AgNe1, AgNeCurves} or in several parts of this  note (see e.g., Example \ref{ex:pairs} or 
Lemma \ref{lem:CDPcomb}). 

\begin{theorem}[Categorification Theorem]\label{th:comblattice}  \cite[Theorem 5.2.1]{AgNe1}
Assume that   $d<\infty$, the height function   $h$ satisfies the stability property, and the pair $(h,h^\circ)$
satisfies the Combinatorial Duality  Property. Then
\begin{center}
 $eu(\bH_*(R,w))=h^\circ (0)-h^\circ (d)  $, 
 \end{center}
 i.e., $\mathbb{H}_*$ categorifies the difference $h^\circ (0)-h^\circ (d)$.
\end{theorem}

This theorem provides us a simple recipe to obtain lattice homological categorifications of certain numerical invariants. Our results in this manuscript specialize this recipe even further.

\begin{example}\label{ex:pairs}
T. \'Agoston and the first author already defined numerous lattice (co)homo\-logy theories in the realm of complex singularity theory using the decomposition of the weight function into pairs  $(h,h^\circ)$ satisfying the above properties. These constructions have the common key property that the resulting lattice homology modules categorify important geometric invariants. For example, the analytic lattice homology of reduced complex curve singularities \cite{AgNeCurves} categorifies the delta invariant, whereas the analytic lattice homology of normal surface singularities (with $\mathbb{Q}HS^3$ link) \cite{AgNe1} categorifies the geometric genus. We will review (and highly generalize in an algebraic context) these constructions from a unified viewpoint in sections \ref{s:deccurves} and \ref{s:dnagy}.

We note that the weight function  of the topological lattice homology of normal surface singularities (cf. Example \ref{toplc}) does not have a splitting into 
a pair $(h,h^\circ)$ with the required properties, at least known by the authors.
\end{example}
 
\subsection{Weight functions satisfying the matroid rank inequality}\,

In this subsection we will provide  some new technical results  (combinatorial in nature) regarding weight functions satisfying the matroid rank inequality. On first reading the reader might skip the proofs of these statements in order to reach faster the main construction of this manuscript --- the lattice homology of integrally closed submodules --- , however, these results will be needed to prove its main properties. Our setting is the following: we work over the lattice $\mathbb{Z}^r$ with standard basis $\{e_v\}_{v \in \mathcal{V}}$ and consider a weight function $w_0:\mathbb{Z}^r \rightarrow \mathbb{Z}$ defined on some rectangle $R(d^-, d^+)$ (here we allow $d^-=-\infty$ and\,/\,or $d^+=\infty$ as well) satisfying the matroid rank inequality as described in (\ref{eq:matroid}), 
or, equivalently, satisfying:
\begin{equation*}
    w_0(\ell + e_v)-w_0(\ell) \geq w_0(\ell + \bar{\ell} + e_v)-w_0(\ell + \bar{\ell}) \hspace{5mm} \forall\, \ell, \bar{\ell} \in \mathbb{Z}^r, \text{ such that } \bar{\ell}\geq 0, \bar{\ell}\ngeq e_v,
\end{equation*}
and $R(\ell, \ell+\bar{\ell} +e_v) \subset R(d^-, d^+)$.

\begin{examples}\label{ex:hhcircmatroid}
    (i) If $h$ and $h^\circ$ of paragraph \ref{bek:comblattice} 
    are defined on the rectangle $R(0, d)$ and 
    both satisfy the matroid rank inequality, then so does the weight function $w_0$ defined by formula (\ref{eq:w_0def}). Moreover, this is also true if $h$ satisfies the matroid rank inequality and $h^\circ:=h^{sym}_d$ is defined by symmetrization   with respect to $d$. \\
    \noindent (ii) A simple computation shows that 
    the weight function $w_{top, 0}$ of the topological lattice homology of complex normal surface singularities (see Example \ref{toplc})  satisfies the matroid rank inequality. 
\end{examples}

\begin{prop}\label{prop:finrectcomb}
    Let $w_0: \mathbb{Z}^r \rightarrow \mathbb{Z}$ be any function defined on the whole lattice $\mathbb{Z}^r$ satisfying the matroid rank inequality. Suppose that the following conditions are satisfied:
    \begin{enumerate}
        \item there exist lattice points $d^-, d^+ \in \mathbb{Z}^r, \  d^- \leq  d^+,$ such that for any $\ell^+\geq d^+, \ \ell^- \leq d^-$ and any $v \in \mathcal{V}$ we have $w_0(\ell^++e_v) \geq w_0(\ell^+)$ and $w_0(\ell^- - e_v) \geq w_0(\ell^-)$;
        \item for any lattice point $\ell \in \mathbb{Z}^r$ and index $v \in \mathcal{V}$ we have $\lim_{k \to \infty}w_0(\ell \pm  ke_v) =\infty$.
    \end{enumerate}
    Then $w_0$ is a weight function, i.e., for any integer $n \in \mathbb{Z}$, the set $w_0^{-1}((-\infty, n])$ is finite  (cf. (\ref{weightfctn})). Moreover, the inclusion $R(d^-, d^+) \hookrightarrow \mathbb{R}^r$ of cubical complexes induces a bigraded $\mathbb{Z}[U]$-module isomorphism $\mathbb{H}_*(R(d^-, d^+), w) \cong \mathbb{H}_*(\mathbb{R}^r, w)$ between the lattice homology modules associated with the weight function $w_0$ (defined in section \ref{ss:latweight}). In fact, this follows from the fact that the inclusions $S_n \cap R(d^-, d^+) \hookrightarrow S_n$ are all homotopy equivalences for every $n \geq m_w=\min w_0$ (they are homotopy inverses of strong deformation retractions $S_n \searrow S_n \cap R(d^-, d^+)$).
\end{prop}

\begin{proof}
For the first part, it is enough to prove that for every $n \in \mathbb{Z}$ the cubical complex $S_n$, defined in paragraph \ref{9complexb}, is finite, i.e., it is contained in a finite rectangle. With this aim, we will construct inductively a sequence of finite rectangles
\begin{equation*}
    R(d^-, d^+) = R(d^-_0, d^+_0) \subset  R(d^-_1, d^+_1) \subset \cdots 
    \subset R(d^-_k, d^+_k)  \subset  R(d^-_{k+1}, d^+_{k+1}) \subset \cdots
\end{equation*} such that for any $n \in \mathbb{Z}$ there exists some $k_n \in \mathbb{N}$ such that $S_n \subset R(d^-_{k_n}, d^+_{k_n})$.

Once a pair $(d_k^-, d_k^+)$ is constructed let us denote by $b_k:=\min \left\{w_0\big|_{\partial R\left(d_k^-, d_k^+\right)}\right\}$. Then we wish to construct the sequences of lattice points $\{d^+_{k}\}_{k \geq 0}$ and $\{ d^-_{k}\}_{k \geq 0}$ (marking the endpoints of our finite rectangles) such that they also satisfy the following:
    \begin{itemize}
        \item $d^+_{k+1}>d^+_{k}$ and $d^-_{k+1}<d^-_{k}$ for all $k$;
        \item $\lim_{k \to \infty}d^\pm_{k,v}=\pm \infty$ for every $v \in \mathcal{V}$;
        \item $b_{k+1}\geq b_{k}+1$ and 
        \item $w_0 \Big|_{R\left(d^-_{k+1}, d^+_{k+1}\right) \setminus R\left(d^-_{k}, d^+_{k}\right)}\geq b_k$;
    \end{itemize}
     hence making sure $S_{b_k} \subset R(d^-_{k+1}, d^+_{k+1})$ and $\lim_{k \to \infty}b_k=\infty$.
     
     For the sake of convenience we introduce a total ordering of the index set $\mathcal{V}$. The idea now is to choose $d^\pm_{0}=d^\pm$ and then $d^+_{k+1}=d^+_{k}+\sum_{v=1}^rl_{k,v}e_v$, where $l_{k,v}$ are recursively defined by being the smallest integer $l$ such that 
     \begin{center}
     $w_0(d^+_{k}+\sum_{w=1}^{v-1}l_{k,w}e_w + le_v) > w_0(d^+_{k}+\sum_{u=1}^{v-1}l_{k,w}e_w)$.
     \end{center}
     Such an $l_{k,v}$ must exist due to condition \textit{(2)}, whereas by condition  \textit{(1)} we have
     \begin{center}
         $w_0(d^+_{k}+\sum_{w=1}^{v-1}l_{k,w}e_w + le_v) = w_0(d^+_{k}+\sum_{w=1}^{v-1}l_{k,w}e_w)$ if $0 \leq l < l_{k, v}$.
     \end{center} 
     Then, by the matroid rank inequality (\ref{eq:matroid}) used for lattice points 
     \begin{center}
         $\ell_1 \in R\left(d_k^-, d_k^+ + \sum_{w=1}^{v-1}l_{k,w}e_w + l_{k,v}e_v\right) \setminus R\left(d_k^-, d_k^+ +\sum_{w=1}^{v-1}l_{k,w}e_w\right)$
     \end{center} and $\ell_2=d_k^++\sum_{w=1}^{v-1}l_{k,w}e_w$ we get (inductively for every $v \in \mathcal{V}$) that 
     \begin{equation*}\begin{split}
     w_0 \big|_{R\left(d_k^-, d_k^+ + \sum_{w=1}^{v-1}l_{k,w}e_w + l_{k,v}e_v\right) \setminus R\left(d_k^-, d_k^+ +\sum_{w=1}^{v-1}l_{k,w}e_w\right)} & \geq b_k, \ \
    \text{ and }
         \\
         w_0 \big|_{R\left(d_k^-, d_k^+ + \sum_{w=1}^{v-1}l_{k,w}e_w + l_{k,v}e_v\right) \cap \lbrace x \in \mathbb{R}^r\,:\, x_v = d^+_{k,v} + l_{k,v} \rbrace} & \geq b_k+1.
     \end{split}\end{equation*}
     By detailed analysis of these extensions, we can verify through the matroid rank inequality, that (after doing the symmetric construction to obtain $d_{k+1}^-$ from $d_k^-$) indeed $b_{k+1} > b_k$. Thus $w_0$ is truly a weight function.

    For the second assertion, by the finiteness of the $S_n$ complexes, it is enough to prove that for any $k \in \mathbb{N}$ and $n \in \mathbb{Z}$, the inclusion $\iota_{k,n}:R(d^-_k, d^+_k) \cap S_n \hookrightarrow  R(d^-_{k+1}, d^+_{k+1}) \cap S_n$ is a homotopy equivalence. Notice that $\iota_{k,n}$ is just the composition of inclusion maps of type
    \begin{equation*}
        R\Big(d_k^-, d_k^+ + \sum_{w=1}^{v-1}l_{k,w}e_w + le_v\Big) \cap S_n \hookrightarrow R\Big(d_k^-, d_k^+ + \sum_{w=1}^{v-1}l_{k,w}e_w + (l+1)e_v\Big) \cap S_n, \ 0 \leq l < l_{k,v}. 
    \end{equation*}
    Let us denote temporarily the lattice point $d_k^+ + \sum_{w=1}^{v-1}l_{k,w}e_w + le_v$ by $c$. We claim that these inclusions admit a strong deformation retraction homotopy inverse: 
    \begin{center}
        $[0, 1]\ni t \mapsto r_t:R(d^-_k, c+e_v) \cap S_n \rightarrow R(d^-_k, c)\cap S_n$.
    \end{center} 
    The key observation is that for any lattice point $x\in R(d^-_k, c)\cap \Z^r$ with 
    $x_{v}=c_{v}$, we have, by matroid rank inequality, 
$w_0(x+e_{v})- w_0(x)\geq w_0(c+e_v)-w_0(c) \geq 0$, hence for every lattice point $x+e_v \in S_n$ we have $x \in S_n$, too. Then the retraction $r_1: R(d^-_k, c+e_v) \cap S_n \rightarrow R(d^-_k, c)\cap S_n$  can be defined as follows. On  $R(d^-_k, c)\cap S_n$ the map $r_1$  is the identity, while if $\square =(x,I)\subset S_n\cap R(d^-_k, c+e_v) )$
is not included  in   $R(d^-_k, c)$, then we retract it parallel to the $v^{\text{th}}$ coordinate axis to
$(x, I\setminus \{v\})$. The strong deformation retraction can now be defined simply, extending this retraction in a linear way. Furthermore, these retractions defined on individual cubes glue together. The composition of these retractions will yield the desired homotopy inverse for $\iota_{k,n}$ for any $k$ and $n$.
\end{proof}

\begin{remark}
    Conditions \textit{(1)} and \textit{(2)} of Proposition \ref{prop:finrectcomb} and the matroid rank inequality imply the following:
    \textit{
    \begin{enumerate}
        \setcounter{enumi}{2}
        \item for any (or, equivalently, for some) 
     sequences $\{ \ell^+_k\}_{k \in \mathbb{N}}$ and $\{ \ell^-_k\}_{k \in \mathbb{N}}$ of lattice points, such that $\lim_{k \to \infty} \ell^\pm_{k,v}=\pm \infty \ \  \forall \,v \in \mathcal{V}$, we have $\lim_{k \to \infty}w_0(\ell^\pm_k)=\infty$.
    \end{enumerate}}
   \noindent  However, conditions \textit{(1)} and \textit{(3)} do not imply condition \textit{(2)}, see, e.g., $w_0: \mathbb{Z}^2 \rightarrow \mathbb{Z}, \ (l_1, l_2) \mapsto |l_2|$.
\end{remark}

\begin{cor}\label{cor:finrectcomb}
    Let $h, h^\circ: \mathbb{Z}^r \rightarrow \mathbb{Z}$ be an increasing and, respectively, 
    a decreasing function as in paragraph \ref{bek:comblattice}, both defined on the full lattice and both satisfying the matroid rank inequality, and set $w_0=h+h^\circ - h^{\circ}(0)$ as in (\ref{eq:w_0def}). Suppose that the following conditions hold:
    \begin{itemize}
        \item[(1')]there exist lattice points $d^-, d^+ \in \mathbb{Z}^r, \ d^- 
        \leq  d^+$, such that for any $\ell^+\geq d^+, \ \ell^- \leq d^-$ we have $h^\circ (\ell^+)=h^\circ (d^+)$ and $h(\ell^-) =h(d^-)$; 
        \item[(2')] for any index $v \in \mathcal{V}$ we have $\lim_{k \to \infty}h(ke_v) =\lim_{k \to \infty}h^\circ (-ke_v) =\infty$. 
    \end{itemize}
    Then $w_0$ is a weight function, i.e., for any integer $n \in \mathbb{Z}$, the set $w_0^{-1}((-\infty, n])$ is finite. Moreover, the inclusion $R(d^-, d^+) \hookrightarrow \mathbb{R}^r$ induces a $\mathbb{Z}[U]$-module isomorphism $\mathbb{H}_*(R(d^-, d^+), w) \cong \mathbb{H}_*(\mathbb{R}^r, w)$. In fact, this follows from the fact that the inclusions $S_n \cap R(d^-, d^+) \hookrightarrow S_n$ are all homotopy equivalences for every $n \geq m_w$ (they are homotopy inverses of strong deformation retractions $S_n \searrow S_n \cap R(d^-, d^+)$).
\end{cor}

\begin{proof}
    We only have to check the conditions of Proposition \ref{prop:finrectcomb}. By part (i) of Example \ref{ex:hhcircmatroid}, $w_0$ satisfies the matroid rank inequality, whereas Condition \textit{(1)} clearly follows from the property  {\it (1')} of $h$ and $h^\circ$, with the same $d^\pm$. Note also, that by the matroid rank inequality and the decreasing property of $h^\circ$, condition {\it (1')} implies that
    \begin{equation}\label{eq:0d-stabilization}
        \forall\, \ell \in \mathbb{Z}^r \text{ we have } h^\circ(\ell) = h^\circ(\min \{ \ell, d^+\}) \ (\text{similarly, }h(\ell) = h(\max \{ \ell, d^-\})).
    \end{equation}
    [These identities do not depend on each other, they remain true without the $d^- 
    \leq d^+$ condition as well.]

    Considering condition \textit{(2)} of Proposition \ref{prop:finrectcomb}, we will prove that for any lattice point $\ell \in \mathbb{Z}^r$ and index $v \in \mathcal{V}$ the limit $\lim_{k\to \infty}w_0(\ell + ke_v)=\infty$ (the  $k\to -\infty$ case is completely analogous). First let us suppose that $\ell\geq {\max\{d^+, 0\}}$. Then for every $k \geq 0$ we have $w_0(\ell + ke_v)=h(\ell + ke_v) + h^\circ(d^+)-h^\circ(0)$, hence it is enough to prove that $\lim_{k\to \infty}h(\ell + ke_v)=\infty$. But this follows from $\lim_{k \to \infty}h(ke_v)=\infty$ by the increasing property of $h$ (since $h(\ell + ke_v) \geq h(ke_v) $ for all $k$).  

    Now let us consider an arbitrary  $\ell \in \mathbb{Z}^r$. Notice that for any integer $k \geq k_0:=-\ell_v+d^+_v$ we have $w_0(\ell + ke_v)=h(\ell + ke_v) + h^\circ(\min\{\ell + k_0e_v, d^+\}) + h^\circ (0)$, hence once again it is enough to prove $\lim_{k\to \infty}h(\ell + ke_v)=\infty$. {Now, we can choose the lattice point $\ell_0^+ :=\max\{\ell + k_0 e_v, d^+, d_v^+e_v\}\in \mathbb{Z}^r$ satisfying $\ell^+_{0, v}=(\ell+k_0e_v)_v=\ell_v+k_0 =d^+_v$, and then, by the matroid rank inequality, we have the bound $h(\ell + ke_v)\geq h(\ell + k_0e_v) + h(\ell^+_0 + (k-k_0)e_v)  - h(\ell^+_0) $ for any $k \geq k_0$, with the first and last term on the right hand side fixed.  Since $\ell^+_0 + (\max\{k_0, -\ell_v\}-k_0)e_v\geq \max \{d^+, 0\}$ we are done by the previous paragraph.}
\end{proof}

We also have the following extension property 
of the  Combinatorial Duality Property under the assumption of the height functions satisfying the matroid rank inequality. 

\begin{prop}\label{prop:rectCDP->fullCDP}
    Let $h, h^\circ: \mathbb{Z}^r \rightarrow \mathbb{Z}$ be an increasing and, respectively, a decreasing function as in paragraph \ref{bek:comblattice}, both defined on the full lattice, satisfying the matroid rank inequality, and suppose that there exist lattice points $d^-, d^+ \in \mathbb{Z}^r, \ d^- \leq d^+$ such that for any $\ell^+\geq d^+, \ \ell^- \leq d^-$ we have $h^\circ (\ell^+)=h^\circ (d^+)$ and $h(\ell^-) =h(d^-)$.
    If the pair $(h, h^\circ)$ satisfies the Combinatorial Duality Property (cf. Definition \ref{def:COMPGOR}) in the rectangle  $R(d^-, d^+)\cap \mathbb{Z}^r$, then it also satisfies it in the whole lattice $\mathbb{Z}^r$.
\end{prop}

\begin{proof}
    We have to prove that $h(\ell+e_v)-h(\ell)$ and $h^\circ (\ell+e_v)-h^\circ (\ell)$  cannot be    simultaneously nonzero for any $\ell \in \mathbb{Z}^r, \ v \in \mathcal{V}$. By identity (\ref{eq:0d-stabilization}), if $\ell_v <d^-_v$, then $h(\ell) = h(\ell +e_v)$, whereas if $\ell_v \geq d^+_v$, then $h^\circ(\ell) = h^\circ(\ell + e_v)$, hence the only nontrivial case is $d_v^- \leq \ell_v < d^+_v$. 
    
    Fix then such an $\ell$ and coordinate $v\in \mathcal{V}$. Let us denote by $\ell^-:=\max \{ \ell, d^-\}$, by $\ell^+:=\min\{ \ell,d^+\}$ and by $\ell^R:=\max\{\ell^+, d^-\}=\min\{ \ell^-, d^+\} \in R(d^-, d^+-e_v)$. Now, from  the increasing property of $h$, the decreasing property of $h^\circ$ and the matroid rank inequality (cf. (\ref{eq:matroid}), see also the stability property (\ref{eq:stability}) for the increasing $h$), we get the following: $h(\ell^R+e_v)=h(\ell^R)$ would imply $h(\ell^-+e_v)=h(\ell^-)$, whereas $h^\circ(\ell^R+e_v)=h^\circ(\ell^R)$ would imply $h^\circ(\ell^++e_v)=h^\circ(\ell^+)$. Using identity (\ref{eq:0d-stabilization}) for $\ell$ and $\ell +e_v$, we see that the identity $h(\ell^R+e_v)=h(\ell^R)$ would imply $h(\ell+e_v)=h(\ell)$, whereas the identity $h^\circ(\ell^R+e_v)=h^\circ(\ell^R)$ would imply $h^\circ(\ell+e_v)=h^\circ(\ell)$. However, from the assumptions we know that at least one of $h(\ell^R+e_v)-h(\ell^R)$ and $h^\circ (\ell^R+e_v)-h^\circ (\ell^R)$  is zero. 
\end{proof}
\newpage
\section{Lattice homology of realizable submodules}\label{s:4}

In this section we present the main conceptual result of the manuscript.
First we set some notations and describe the general constructions involved. In fact, we will develop further the construction of Example  \ref{ex:multifilt}, although, contrary to that case, this time  we will consider two filtrations.
More specifically, we will replace the vector space $V$ of Example \ref{ex:multifilt}
by   a Noetherian  ring $\cO$, with a $k$-algebra structure,  and the first 
$\bZ^r$-filtration will be associated with $r$ discrete valuations on $\cO$.
On the other hand, the second multifiltration will be supported on an $\cO$-module $M$. 

Throughout this section $k$ will always denote a field (we do not impose any restriction on it, though, in all of our applications this will be $\mathbb{C}$), $\mathcal{O}$ will be a $k$-algebra, usually assumed to  be Noetherian (only for convenience and stronger relations with commutative algebra), and $M$ a finitely generated $\mathcal{O}$-module (once again, finite generacy is assumed only for convenience --- for the more general setup see subsection \ref{ss:MOSTGEN}).
In several of our applications the corresponding $k$-algebra $\cO$ will be the local analytic ${\mathbb{C}}$-algebra $\cO_{X,o}$ of a  complex analytic spacegerm $(X,o)$, whereas $M$ an $\mathcal{O}$-module of some kind of differential forms.

\subsection{Extended discrete valuations and the associated filtrations}\label{ss:discval+filt}\,

We will use the notion of discrete valuations in the following (generalized) sense
(cf.  e.g., \cite{FJ,S}). (Such valuations have different names in the literature, as 
semivaluations, or pseudo-valuations.) 
\begin{define} \label{def:gdv}
A function $\frv:\cO \rightarrow \overline{\N}= \N \cup \{ \infty \}$
is a discrete valuation if  the following properties  hold:
\begin{itemize}
    \item[(a)] $\frv(fg) = \frv(f) + \frv(g)$ \ for all  $f, g \in \cO$\
    (where $\infty + n = n+ \infty = \infty +\infty=\infty$ for every $n \in \N$);
    \item[(b)] $\frv(f+g) \geq \min \{ \frv(f), \frv(g) \}$ \ for all $f, g \in \cO\ $;
    \item[(c)] 
    $\frv$ is nonconstant, i.e., there exists some $f \in \mathcal{O}$ with $\mathfrak{v}(f) \in \mathbb{N}_{>0}$ \
    (hence, the value of its powers can get arbitrarily high).
\end{itemize}

\end{define}
From the definition one can deduce (see e.g., the argument in \cite[1.2]{FJ})
that $\frv(0)=\infty$, $\frv|_{k^*}=0$ and the set $\mathfrak{p}_{\mathfrak{v}}:=\{f\in\cO\,:\, \frv(f)=\infty\}$ (which we will sometimes call the \textit{`core'} of the valuation) is a prime ideal of $\cO$.  Moreover,  in property (b) we necessarily have equality if $\mathfrak{v}(f) \neq \mathfrak{v}(g)$. 

\begin{remark}
Note that Definition \ref{def:gdv}
is weaker than the `classical' definition of the discrete valuation: we do not require 
neither that $\cO$ is a domain, nor that $\mathfrak{p}_{\mathfrak{v}}=(0)$.  Also, we see immediately that for any nilpotent element $f \in \mathcal{O}$, we must have $\mathfrak{v}(f)= \infty$. Also, any such discrete valuation descends to a `classical' discrete valuation $\mathfrak{v}_{\mathcal{O}/\mathfrak{p}_{\mathfrak{v}}}$ on $\mathcal{O}/\mathfrak{p}_{\mathfrak{v}}$ with $\mathfrak{v}_{\mathcal{O}/\mathfrak{p}_{\mathfrak{v}}}(f+\mathfrak{p}_{\mathfrak{v}}):=\mathfrak{v}(f)$ being independent of the representative chosen.
\end{remark}

\begin{define} \label{def:edv}
For a fixed $\cO$-module $M$ we will call a pair $(\mathfrak{v}, \mathfrak{v}^M)$ of maps an `\textit{extended discrete valuation}' if 
\begin{itemize}
    \item  $\mathfrak{v}:\cO \rightarrow \overline{\N}$ is a discrete valuation in the sense of Definition \ref{def:gdv};
    \item $\mathfrak{v}^M: M \rightarrow \overline{\Z}= \mathbb{Z} \cup \{ \infty\}$ satisfies the following properties:
    \begin{itemize}
    \item[(a)] $\frv^M(fm) = \frv(f) + \frv^M(m)$ \ for all  $f \in \cO, m \in M$\\
    (where, similarly as above,  $\infty + n = n+ \infty = \infty +\infty=\infty$ for every $n \in \mathbb{Z}$);
    \item[(b)] $\frv^M(m+m') \geq \min \{ \frv^M(m), \frv^M(m') \}$ \ for all $m, m' \in  M$;
    \item[(c)]  and $\frv^M$ is nonconstant, meaning that there exists some $m \in M$ with $\mathfrak{v}^M(m)\neq \infty$.
\end{itemize}
\end{itemize}

\end{define}
Similarly as above, in (b) one has  equality if $\frv^M(m) \neq \frv^M(m')$). 
For ease of language we will call $\mathfrak{v}$ `\textit{extended}' and $\mathfrak{v}^M$ its `\textit{extension}'. Moreover, for the pair $(\frv,\frv^M)$
we  might just  write  $\frv $.

Similar (although real valued) objects, called `$v$-valuative functions', were recently studied in \cite{Hada}, see especially \cite[Definition 3.1]{Hada}.

\begin{example}\label{ex:trivext}
    If $\mathfrak{v}: \mathcal{O} \rightarrow \overline{\mathbb{N}}$ is a discrete valuation and $M=\mathcal{O}^p$, then we can construct an extension e.g., as $\mathfrak{v}^M: (f_1, \ldots, f_p) \mapsto \min \{\mathfrak{v}(f_1), \ldots, \mathfrak{v}(f_p)\}$.  Moreover, any extension can be restricted to a submodule $N \leq M$.
\end{example}

\begin{remark}\label{rem:translating}
    (i)  If $(\mathfrak{v}, \mathfrak{v}^M)$ is an extended discrete valuation, then so are $(n\mathfrak{v}, n\mathfrak{v}^M)$ and $(\mathfrak{v}, \mathfrak{v}^M+n')$ for any $n \in \mathbb{N}\setminus \{0\}$ and $n' \in \mathbb{Z}$.
    
    \noindent (ii) With $\mathfrak{p}_{\mathfrak{v}}:=\{f\in\cO\,:\, \frv(f)=\infty\}$ we have $\mathfrak{v}^M\big|_{\mathfrak{p}_{\mathfrak{v}}M} \equiv \infty$. Moreover, $\mathfrak{v}^M$ descends to the $\mathcal{O}/\mathfrak{p}_{\mathfrak{v}}$-module $M/\mathfrak{p}_{\mathfrak{v}}M$, which sends the torsion part of $M/\mathfrak{p}_{\mathfrak{v}}M$ to $\infty$ as well. 

    \noindent (iii) 
In fact, the extensions $\mathfrak{v}^M$ to $M$ correspond to extensions to the $\mathcal{O}_{\mathfrak{p}_{\mathfrak{v}}}/(\mathfrak{p}_{\mathfrak{v}})_{\mathfrak{p}_{\mathfrak{v}}}={\rm Frac}(\mathcal{O}/\mathfrak{p}_{\mathfrak{v}})$-vector space $M_{\mathfrak{p}_{\mathfrak{v}}}/\mathfrak{p}_{\mathfrak{v}}M_{\mathfrak{p}_{\mathfrak{v}}}$, and,  thus, for any given valuation $\mathfrak{v}$, if $M_{\mathfrak{p}_{\mathfrak{v}}}/\mathfrak{p}_{\mathfrak{v}}M_{\mathfrak{p}_{\mathfrak{v}}}\neq 0$,  there are infinitely many extensions $\mathfrak{v}^M$, see Remark \ref{rem:exttoMandW}.
\end{remark}

\begin{define}
    To any extended valuation $(\frv,\frv^M)$ on a fixed $\mathcal{O}$-module $M$ we associate the following descending filtrations: 
    \begin{equation}\label{eq:fvMl}
    \text{for any } n\in\Z: \ \ \ \cF_{\frv}(n) := \{ f \in \cO \,:\, \frv(f) \geq n \}\triangleleft\mathcal{O} \ \ \text{ and } \ \ \cF^M_{\frv}(n) := \{ m \in M \,:\, \frv^M(m) \geq n \}\leq M .
\end{equation}
Moreover, for any finite (multi)set $\mathcal{D}=\{\frv_1, \ldots, \frv_r\}$
of extended discrete valuations and tuple\,/\,lattice point $\ell=(\ell_1,\ldots, \ell_r)\in \mathbb{Z}^r$ we also define
\begin{equation}\label{eq:fdMl}
    \mathcal{F}_{\mathcal{D}}:\ \Z^r \ni \ell \mapsto \bigcap_{v =1}^{r}\mathcal{F}_{\frv_v}(\ell_v) \triangleleft \mathcal{O} \ \ \text{ and }\ \ \mathcal{F}^M_{\mathcal{D}}:\ \Z^r \ni \ell \mapsto \bigcap_{v =1}^{r}\mathcal{F}^M_{\frv_v}(\ell_v) \leq M.
\end{equation}
\end{define}

\begin{ass}\label{ass:fincodim}
    In the sequel we will always assume that for any $\ell \in \mathbb{Z}^r$ the ideal $\mathcal{F}_{\mathcal{D}}(\ell)$ (respectively the submodule $\mathcal{F}_{\mathcal{D}}^M(\ell)$) has finite $k$-codimension in $\cO$ (respectively in $M$). Moreover, a \textit{`finite collection of extended (discrete) valuations'}  will always mean a finite (multi)set of extended discrete valuations satisfying this assumption.
\end{ass}

\begin{define}
    For any finite collection $\mathcal{D}=\{\frv_1, \ldots, \frv_r\}$
of extended discrete valuations 
we  introduce the following Hilbert functions:
\begin{equation}\label{eq:handhM}
    \mathfrak{h}_{\mathcal{D}}, \mathfrak{h}^M_{\mathcal{D}}: \mathbb{Z}^r \rightarrow \mathbb{Z}, \hspace{5mm} \mathfrak{h}_{\mathcal{D}}(\ell) = \dim_k\mathcal{O}/\mathcal{F}_{\mathcal{D}}(\ell) \hspace{5mm} \text{and} \hspace{5mm} \mathfrak{h}^M_{\mathcal{D}}(\ell) = \dim_kM/\mathcal{F}_{\mathcal{D}}^M(\ell).
\end{equation}
Since $\mathcal{F}_\mathcal{D}(0)=\mathcal{O}$, we have $\mathfrak{h}_{\mathcal{D}}(0)=0$ and both $\mathfrak{h}_{\mathcal{D}}$ and $\mathfrak{h}^M_{\mathcal{D}}$ are nonnegative and increasing.
\end{define}

\begin{obs} \label{obs:hhMmatroid}
    As the multifiltrations satisfy condition (\ref{eq:FaFb=Fmax}), both Hilbert functions satisfy the matroid rank inequality.
\end{obs}

\subsection{The lattice homology associated with a finite collection of extended valuations}\, 

Let us recall our setup: $k$ is a field, $\cO$ a Noetherian $k$-algebra, $M$ a finitely generated $\cO$-module and $\mathcal{D}=\{\frv_1, \ldots, \frv_r\}$ a finite collection of extended discrete valuations (satisfying Assumption \ref{ass:fincodim}). 

In view of paragraph \ref{bek:comblattice}, we can associate to this setting the following lattice homology theory: set the lattice $\mathbb{Z}^r$ (where $r=|\mathcal{D}|$ denotes the number of extended valuations) with standard basis $\{e_v\}_v$ and consider the functions 
\begin{equation*}
    \mathfrak{h}_{\mathcal{D}}, \mathfrak{h}^\circ_{\mathcal{D}}: \mathbb{Z}^r \rightarrow \mathbb{Z}, \hspace{4mm} \mathfrak{h}_{\mathcal{D}}(\ell) = \dim_k\mathcal{O}/\mathcal{F}_{\mathcal{D}}(\ell) \hspace{4mm} \text{and} \hspace{4mm} \mathfrak{h}_{\mathcal{D}}^\circ(\ell) = \mathfrak{h}^M_{\mathcal{D}}(-\ell)=\dim_kM/\mathcal{F}_{\mathcal{D}}^M(-\ell).
\end{equation*}
Set the weight $w_{\mathcal{D}, 0}(\ell)=\frh_{\mathcal{D}}(\ell)+\frh^\circ_{\mathcal{D}} (\ell)-\frh^\circ_{\mathcal{D}} (0)$ for every $\ell \in \mathbb{Z}^r$ as in (\ref{eq:w_0def}). Next, extend the weight function $w_{\mathcal{D}, 0}$ to higher dimensional cubes  following (\ref{eq:9weight}). Then run the construction of paragraph \ref{9complexb} to obtain the cubical complex $S_{n, \mathcal{D}}$, for every $n \in \mathbb{Z}$, and get the lattice homology $\mathbb{Z}[U]$-module $\mathbb{H}_*(\mathbb{R}^r, w_{\mathcal{D}})$. 

\begin{lemma}\label{lem:w0matroid}
    The function $w_{\mathcal{D}, 0}$ defined above satisfies the matroid rank inequality. Moreover, for any $k \in \mathbb{N}$, the set $w_{\mathcal{D}, 0}^{-1}((-\infty, k])$ is finite, i.e., $w_{\mathcal{D}, 0}$ is indeed a weight function (in the sense of paragraph \ref{bek:wf}).
\end{lemma}
We defer the proof to subsection \ref{ss:ddep}.
\vspace{2mm}

As we will see in the next section this general construction (applied for different modules and collections $\mathcal{D}$) recovers several previously defined lattice homology constructions  as particular cases, moreover, we will list several new applications as well.

The construction, which  relies on the collection of extended valuations $\mathcal{D}$,
suggests that the output $\bH_*(\mathbb{R}^r, w_{\mathcal{D}})$ depends  essentially on the choice of $\mathcal{D}$. However, our main theorem, stated in the following subsection, says that (under certain additional assumption regarding $\mathcal{D}$, which is
`almost automatically' satisfied)
 the lattice homology module $\bH_*(\mathbb{R}^r, w_{\mathcal{D}})$
 depends only on the pair $\mathcal{F}_{\mathcal{D}}^M(0)\leq M$. 

\subsection{The main theorem}\,

The precise statement is the following:

\begin{theorem}[Independence Theorem]\label{th:IndepMod} Let $k$ be any field, $\mathcal{O}$ a (Noetherian) $k$-algebra and $M$ a finitely generated module over it. Let $\mathcal{D}=\{\frv_1, \ldots, \frv_r\}$ and $\mathcal{D}'=\{\frv'_1, \ldots, \frv'_{r'}\}$ be two collections of extended discrete valuations (satisfying Assumption \ref{ass:fincodim}).  
        Suppose that the following conditions hold:
        \begin{itemize}
            \item $\mathcal{F}_{\mathcal{D}}^M(0) = \mathcal{F}_{\mathcal{D}'}^M(0)$ and
            \item both pairs $(\mathfrak{h}_{\mathcal{D}}, \mathfrak{h}_{\mathcal{D}}^\circ)$ and $(\mathfrak{h}_{\mathcal{D}'}, \mathfrak{h}_{\mathcal{D}'}^{\circ})$ satisfy the Combinatorial Duality Property (cf. Definition \ref{def:COMPGOR}).
        \end{itemize}
        Then the spaces $S_{n, \mathcal{D}}$ and $S_{n, \mathcal{D}'}$ 
        associated with the corresponding lattices and weight functions 
        are homotopy equivalent for every $n \in \mathbb{Z}$.
        Even more, the homotopy equivalences respect the inclusions
        $S_{n, \mathcal{D}}\hookrightarrow S_{n+1, \mathcal{D}}$ and $S_{n, \mathcal{D}'}\hookrightarrow S_{n+1, \mathcal{D}'}$, hence
        \begin{equation*}
        \mathbb{H}_*(\mathbb{R}^r, w_{\mathcal{D}}) \cong  \mathbb{H}_*(\mathbb{R}^{r'}, w_{\mathcal{D}'}) \text{ as bigraded } \mathbb{Z}[U]\text{-modules.}
        \end{equation*}
\end{theorem}

 That is, the lattice homology and the homotopy types of the cubic complexes $S_{n, \mathcal{D}}$ are independent of the collection $\mathcal{D}$ of extended discrete valuations, they are invariants merely of the module $M$ and 
 the finite codimensional submodule $\mathcal{F}_{\mathcal{D}}^M(0) \leq M$. 
 
The proof of the Independence Theorem will be given in section \ref{s:proof}, here we will discuss its main consequences.
For a more general formulation see subsection \ref{ss:MOSTGEN}.

The Independence Theorem implies that, once we fixed the module $M$,
we can associate a well-defined lattice homology module to 
any submodule $N$ of $M$, whenever $N$ can be `realized' as  $\mathcal{F}_{\mathcal{D}}^M(0) $
for some collection $\mathcal{D}$  for which the pair
 $(\mathfrak{h}_{\mathcal{D}}, \mathfrak{h}_{\mathcal{D}}^\circ)$ satisfies CDP. 
We will investigate this class of submodules and give the first properties.
We proceed in two steps, first we test the `realizability' condition
$N=\mathcal{F}_{\mathcal{D}}^M(0)$, then the additional  CDP property.

In the sequel we will omit  $\mathcal{D}$ from the notation of the Hilbert functions, weight functions  and  $S_n$ spaces whenever the finite collection of extended discrete valuations we  use  is clear from the context.

\subsection{Lattice homology of realizable submodules and its main properties}\label{ss:lathomofmods}\,

Set again some field $k$, a Noetherian $k$-algebra $\mathcal{O}$ and a finitely generated $\mathcal{O}$-module $M$.

\begin{define}\label{def:REAL}
    Let $N$ be a finite codimensional submodule of $M$. $N$ is called \textit{`realizable'} if some finite collection $\mathcal{D}$ of 
    extended discrete valuations 
    satisfies 
    $N=\mathcal{F}^M_{\mathcal{D}}(0)$. In this case $\mathcal{D}$ is called a \textit{`realization of $N$'}. 
      We say that the realization is a `{\it CDP realization}'
if the associated pair of functions $(\mathfrak{h}_{\mathcal{D}}, \mathfrak{h}^\circ_{\mathcal{D}})$
satisfies the Combinatorial Duality Property. 
\end{define}

These realizability conditions on $N \leq M$ are highly nontrivial, such a realization by some $\mathcal{D}$ might not exist. For example if we choose $M=\mathcal{O}$, then the condition for the existence of such a realization is equivalent to the ideal $N$ being integrally closed, see  Corollary \ref{cor:intersection} below. In the general finitely generated module setting we have the following result:

\begin{prop}[=Theorem \ref{th:REES}]\label{prop:intC}
    If the finitely generated, finite $k$-codimensional submodule $N$ is integrally closed in $M$ in the sense of Rees \cite{ReesMod} (see also \cite[Section 16]{SH} or Definition \ref{def:intclosedmod} below), then it is realizable.
\end{prop}

The next natural question is the following:  is the CDP assumption in the realizability of $N$ (and hence in the Independence Theorem) an `essential' additional restriction? We have the following answer:

\begin{obs}[=Lemma \ref{lem:duplatrukk}]\label{obs:4.3.3}
    If $\mathcal{D}$ is a realization of $N$, then we can easily construct from it a CDP realization $\mathcal{D}^{\natural}$ just by postcomposing every extended discrete valuation with the `doubling' map $\overline{\mathbb{Z}} \xrightarrow{\cdot2}\overline{\mathbb{Z}}$ (the new valuations are denoted simply by $\{2\frv_v\}_v$).  
\end{obs}

Thus, using the Independence Theorem \ref{th:IndepMod}, Proposition \ref{prop:intC} and 
Observation \ref{obs:4.3.3},
for any  {\it realizable\,/\,finite codimensional integrally closed} submodule $N$  we can define its lattice homology.

\begin{define}\label{def:LCofMOD}
    Fix a Noetherian $k$-algebra $\mathcal{O}$ and a finitely generated $\mathcal{O}$-module $M$. Consider a (finitely generated) finite $k$-codimensional realizable submodule $N \leq M$.  From all the possible ones choose an arbitrary  
    realization $\mathcal{D}=\{\mathfrak{v}_v\}_{v\in \mathcal{V}}$ 
    in the sense of Definition \ref{def:REAL}. Fix the lattice $\mathbb{Z}^r$, with rank $r=|\mathcal{V}|$ (with its standard basis). Take the finite collection $\mathcal{D}^\natural=\{2\mathfrak{v}_v\}_{v\in \mathcal{V}}$ of extended discrete valuations and consider the corresponding multifiltrations $\mathbb{Z}^r \ni \ell \mapsto \mathcal{F}_{\mathcal{D}^\natural}(\ell), \mathcal{F}_{\mathcal{D}^\natural}^M(\ell)$ (introduced in (\ref{eq:fdMl})). Then define the functions
\begin{equation}\label{eq:handhcirc0}
    \mathfrak{h}, \mathfrak{h}^\circ: \mathbb{Z}^r \rightarrow \mathbb{Z}, \hspace{5mm} \mathfrak{h}(\ell) = \dim_k\mathcal{O}/\mathcal{F}_{\mathcal{D}^\natural}(\ell) \hspace{5mm} \text{and} \hspace{5mm} \mathfrak{h}^\circ(\ell) = \dim_kM/\mathcal{F}_{\mathcal{D}^\natural}^M(-\ell).
\end{equation}
Set the weight function
$w_0(\ell)=\frh(\ell)+\frh^\circ (\ell)-\frh^\circ (0)$ on the lattice points $\ell \in \mathbb{Z}^r$ (see Lemma \ref{lem:w0matroid}) and extend it to every cube in the cubical decomposition of $\mathbb{R} \otimes \mathbb{Z}^r$ as in (\ref{eq:9weight}). Finally, consider the sublevel sets $S_n$  and define through them the lattice homology module $\mathbb{H}_*(\mathbb{R}^r, w)$ as in paragraph \ref{9complexb}. 

Then, by the above discussion, the lattice homology module is independent of the choice of $\mathcal{D}$ (justifying the name of `Independence Theorem' of Theorem \ref{th:IndepMod}),
it depends only on the pair $N\leq M$.
We will denote it by $\mathbb{H}_*(N \hookrightarrow_{\mathcal{O}} M)$. When the base ring $\mathcal{O}$ is clear from the context, we will omit it from the notation.
\end{define}

\begin{remark}\label{rem:h-h'}
    Note that (using the notions of (\ref{eq:handhcirc0}) $\mathfrak{h}^\circ(0)=\dim_kM/N$, hence, for the lattice points $\ell \in \mathbb{Z}^r, \ \ell \geq 0$,  the weight function
$w_0(\ell)=\frh(\ell)+\frh^\circ (\ell)-\frh^\circ (0)$ agrees with
\begin{equation}\label{eq:w0}
w_0(\ell)  =
\dim_k\mathcal{O}/\mathcal{F}_{\mathcal{D}^\natural}(\ell)-
\dim_k\mathcal{F}_{\mathcal{D}^\natural}^M(-\ell)/N.\end{equation}
\end{remark}

\begin{remark}\label{rem:nodoublingifCDP}
    By the Independence Theorem \ref{th:IndepMod}, if the chosen collection 
    $\mathcal{D}$ is already a CDP realization of $N$, then
    we do not need to convert it to the `doubled' version $\mathcal{D}^\natural$, the construction with $\mathcal{D}$ will yield the same result.
\end{remark}

We state here the first properties of $\mathbb{H}_*(N \hookrightarrow M)$. They will be discussed and proved in later (sub)sections. The first algebraic examples will be given in subsection \ref{ss:NEWEX}.

\begin{theorem}\label{th:properties}
    Let $\mathcal{O}$ be a Noetherian $k$-algebra, $N \leq M$ finitely generated modules with $N$ realizable (hence, by definition, $\dim_k(M/N) < \infty$). \\
\noindent (a) $\mathbb{H}_*(N \hookrightarrow_{\mathcal{O}} M)$ can be computed on a finite rectangle: for any fixed CDP realization $\mathcal{D}=\{\frv_1, \ldots, \frv_r\}$ of $N$ and lattice point $0\ll d \leq \infty$ large enough (the precise meaning of `large enough' depends on the chosen realization), the inclusion 
$R(0,d)\subset \mathbb{R}^r$ induces a  bigraded $\Z[U]$-module 
isomorphism $\bH_*(R(0,d),w)\to \bH_*(\mathbb{R}^r,w)\cong \mathbb{H}_*(N \hookrightarrow M)$ (where $w$ is as defined in (\ref{eq:w0})). In fact, the inclusions $S_n \cap R(0, d) \hookrightarrow S_n$ are all homotopy equivalences for every $n \geq m_w$ (they are homotopy inverses of strong deformation retractions $S_n \searrow S_n \cap R(0, d)$). \\
\noindent (b) The Euler characteristic is well-defined and satisfies $eu(\mathbb{H}_*(N \hookrightarrow M))=\dim_k (M/N)$, i.e., the lattice homology $\mathbb{H}_*(N \hookrightarrow M)$ categorifies the codimension ${\rm codim}_k(N \hookrightarrow M)$.  \\
\noindent (c) The $\mathbb{Z}[U]$-module (and the homotopy type of the $S_n$ spaces) only depends on the quotient module $M/N$ over the algebra $\mathcal{O}/{\rm Ann}_{\mathcal{O}}(M/N)$, where ${\rm Ann}_{\mathcal{O}}(M/N)$ denotes the annihilator of the quotient $\mathcal{O}$-module $M/N$. 
\end{theorem}

This theorem shows that we can construct a well-defined categorification of any numerical invariant, whenever this number can be defined as the codimension of some realizable submodule.
It turns out that in the theory of complex analytic singularities many invariants are of this nature: the delta invariant of reduced curve singularities (see section \ref{s:deccurves}),  the geometric genus and irregularity of isolated singularities of dimension $\geq 2$ (see section \ref{s:dnagy}) and the plurigenera of normal surface singularities (see subsections \ref{s:pluri}--\ref{ss:plurig1}) are some highlighted examples.

\subsection{Realizability of $N$}\label{ss:realofN}\,

In this subsection we discuss the realizability property of $N$ (introduced in Definition \ref{def:REAL}) and compare it with the integral closedness property of Rees. 
Recall that a finite $k$-codimensional submodule $N$ of a finitely generated $\mathcal{O}$-module $M$ (where $\mathcal{O}$ is a Noetherian $k$-algebra) is realizable, if there exists a finite collection $\mathcal{D}=\{(\frv_v,\frv_v^M)\}_{v \in \mathcal{V}}$ of extended discrete valuations such that $$N= \mathcal{F}_{\mathcal{D}}^M(0) = \{ m \in M\,:\, \mathfrak{v}^M_v(m) \geq 0 \ \mbox{for all} \ v \in \mathcal{V}\}.$$

Let us start with the following immediate observations
regarding realizability.

\begin{lemma}\label{lem:contractionREAL}
(a) Any finite collection $\mathcal{D}=\{(\frv_v,\frv_v^M)\}_{v \in \mathcal{V}}$ of extended discrete valuations and lattice point $d=(d_1, \ldots, d_r)\in \Z^r $ defines a realizable submodule  $N:= \mathcal{F}^M_{\mathcal{D}}(d)$. 

(b)     Let $\varphi:M' \rightarrow M$ be a finitely generated $\mathcal{O}$-module morphism and suppose that $N \leq M$ is a finitely generated, finite $k$-codimensional realizable submodule. Then $\varphi^{-1}(N) \leq M'$ is also finite $k$-codimensional  and realizable (whenever $M'$ admits nontrivial extended valuations). 
\end{lemma}

\begin{proof} (a)
By Remark \ref{rem:translating} (i), a possible realization can be given by $\{(\frv_v,\frv_v^M-d_v)\}_{v \in \mathcal{V}}$.

(b) 
    For any realization $\mathcal{D}=\{(\frv_v,\frv_v^M)\}_{v \in \mathcal{V}}$ of $N$, the compositions $\mathfrak{v}^{M'}_v=\mathfrak{v}^M_v \circ \varphi$ will also give extensions of the discrete valuations $\mathfrak{v}_v$ on $\mathcal{O}$, except when $\mathfrak{v}^M_v \big|_{\varphi(M')} \equiv \infty$ . Hence, if the index set $\mathcal{V}':= \mathcal{V} \setminus \{ v\,:\, \mathfrak{v}^M_v \big|_{\varphi(M')} \equiv \infty \}$ is nonempty, then  $\mathcal{D}':=\{(\frv_v,\frv_v^{M'})\}_{v \in \mathcal{V}'}$ gives a realization of $\varphi^{-1}(N) = \mathcal{F}_{\mathcal{D}'}^{M'}(0)$. On the other hand, if $\mathcal{V}' = \emptyset$, then $\varphi^{-1}(N) = M'$, which is realizable in $M'$ (if $M'$ admits nontrivial extended valuations). Indeed, by Assumption \ref{ass:fincodim}, for any extended valuation $(\mathfrak{v}, \mathfrak{v}^{M'})$ there exists some $d_{\mathfrak{v}} \in \mathbb{Z}$ such that $\mathfrak{v}^{M'}-d_{\mathfrak{v}} \geq 0$ on $M'$.
\end{proof}

Let us now recall 
Rees's definition
of {\it integral dependence over a submodule} in the setting of  an
inclusion $N \leq M$ of finitely generated $\mathcal{O}$-modules (with $\mathcal{O}$ Noetherian). 
\begin{nota}
For any prime ideal $\mathfrak{p} \triangleleft\mathcal{O}$ let us denote the fraction field of $\mathcal{O}/\mathfrak{p}$ by $\kappa(\mathfrak{p})$ and
 the $\kappa(\mathfrak{p})$-vector space $M \otimes_{\mathcal{O}}\kappa(\mathfrak{p}) \cong M_{\mathfrak{p}}/\mathfrak{p} M_{\mathfrak{p}}$ by $W(\mathfrak{p})$. We denote by $\iota_\mathfrak{p}$ the composition $M \rightarrow M/\mathfrak{p}M \rightarrow W(\mathfrak{p})$ of the canonical maps.
 \end{nota}

\begin{define}\label{def:intclosedmod}\cite{ReesMod}, see also \cite[Section 16]{SH}.   
An element $m \in M$ is said to be
\textit{`integral over $N$'} if for every \textit{minimal} prime ideal $\mathfrak{p}$ in $\mathcal{O}$ and every discrete valuation
ring $\mathcal{O}_{\mathfrak{v}}$ (in the classical sense) bet\-ween 
 $\mathcal{O}/\mathfrak{p} \leq \mathcal{O}_{\mathfrak{v}}\leq \kappa(\mathfrak{p})$ (i.e., the corresponding valuation $\mathfrak{v}$ satisfies $\mathfrak{v} \geq 0$ on $\mathcal{O}/\mathfrak{p}$ and $\mathfrak{v}(a)=\infty$ if and only if $a=0$ in $\mathcal{O}/\mathfrak{p}$) we have $\iota_{\mathfrak{p}}(m) \in \iota_{\mathfrak{p}}(N)\mathcal{O}_{\mathfrak{v}}$.

The submodule $N$ is called \textit{`integrally closed in $M$'}, if every element $m \in M$, which is integral over $N$, already lies in $N$.
\end{define}

\begin{theorem}\label{th:REES}
    Let $k$ be a field, $\mathcal{O}$ a Noetherian $k$-algebra, $N \leq M$ an inclusion of finitely generated $\mathcal{O}$-modules with ${\rm codim}_k(N\hookrightarrow M)<\infty$. Then $N$ is integrally closed in $M$ if and only if there exists a finite collection $\mathcal{D}=\{ \mathfrak{v}_v\}_{v \in \mathcal{V}}$ of extended discrete valuations, with their cores $\mathfrak{p}_v=\{ f \in \mathcal{O}\,:\, \mathfrak{v}_v(f)=\infty\}$ being \emph{minimal} prime ideals,
    such that $N = \mathcal{F}_{\mathcal{D}}^M(0)$. Hence, finite codimension and integral closedness implies realizability.
\end{theorem}

The proof of Theorem \ref{th:REES} and further  algebraic discussions 
are separated in the Appendix (see especially Corollary \ref{cor:Rees=Real}).
We do not have a proof for the opposite 
statement however (that is, whether the realizability by a certain $\mathcal{D}$ implies
that $N$ is integrally closed in $M$).
Nevertheless, by Theorem \ref{th:REES}
we can apply our results from subsection \ref{ss:lathomofmods} safely  for any finite codimensional integrally closed submodule (in the sense of Rees). 

\bekezdes \textbf{The local analytic case.}\label{bek:locan}\,

In the case when $\cO=\cO_{X,o}$, the local $\C$-algebra of a complex analytic spacegerm $(X, o)$, and $M=\cO^p$ for some $p\geq 1$, integrally closed submodules were considered by Gaffney in \cite{Gaff,GaffKl} as well. By his  definition, $N$ is integrally closed 
(`in the sense of Gaffney')  if $\overline{N}=N$, where $\overline{N}$ --- the integral closure of $N$ in the sense of Gaffney --- consists of elements $x\in 
\cO_{X,o}^p$ such that for all holomorphic map germs $\phi:(\C,0)\to (X,o)$ the induced morphism $\phi^{*}_M:M=\mathcal{O}^p \rightarrow \mathcal{O}_{\mathbb{C}, 0}^p=\mathbb{C}\{t\}^p$ maps $x$ into the submodule $(\phi_M^*(N))\cO_{\C,0}$ generated over $\cO_{\C,0}$ by the image of $N$.
Notice, that any such map germ $\phi$ induces a  discrete valuation in the sense of Definition \ref{def:gdv}:
\begin{equation}\label{eq:curvemapval}
    \mathfrak{v}_\phi: \mathcal{O} \xrightarrow{\phi^*}\cO_{\C,0} \cong \mathbb{C}\{t\} \xrightarrow{{\rm ord}_t} \overline{\mathbb{N}},
\end{equation} which could be extended to $\mathcal{O}^p$ in numerous ways (see, e.g., Example \ref{ex:trivext}). However, it is not clear 
whether each such extension $\mathfrak{v}_\phi^{\mathcal{O}^p}$ would satisfy the inequality $\mathfrak{v}_\phi^{\mathcal{O}^p}(x) \geq \mathfrak{v}_\phi^{\mathcal{O}^p}(N)$. More explicitly, it is not clear whether any extension $\mathfrak{v}_\phi^{\mathcal{O}^p}$ would be compatible with the `extension of scalars' $N \rightarrow (\phi^*_M(N))\mathcal{O}_{\mathbb{C}, o}$ along the valued ring extension $(\mathcal{O}/{\rm ker} \phi^*, \mathfrak{v}_\phi) \hookrightarrow (\mathbb{C}\{t\}, {\rm ord}_t)$ (for even more on this topic see Question \ref{q:Gaffney}). On the other hand, if we start from an extension ${\rm ord}_t^{\mathbb{C}\{t\}^p}:\mathbb{C}\{t\}^p \rightarrow \overline{\mathbb{Z}}$ of ${\rm ord}_t:\mathbb{C}\{t\} \rightarrow \overline{\mathbb{N}}$ and consider $\mathfrak{v}_\phi^{\mathcal{O}^p}:={\rm ord}_t^{\mathbb{C}\{t\}^p} \circ \phi^*_M$, then we 
 get $\mathfrak{v}_\phi^{\mathcal{O}^p}(x) \geq \mathfrak{v}_\phi^{\mathcal{O}^p}(N)$. 
In fact, the property that $ \phi^*_M(x)\in(\phi^*_M(N))\cO_{\C,0}$ can be detected with such extensions to obtain the following characterization:

\begin{theorem}\label{th:Gaffney}
    Let $\mathcal{O}=\mathcal{O}_{X, o}$ be the local algebra of a complex analytic spacegerm $(X, o)$, $M=\mathcal{O}^p$ a free module of rank $p \geq 1$  and $N \leq M$ a finite codimensional finitely generated submodule. Then $N$ is integrally closed (in the sense of Gaffney) in $M$ if and only if there exists a finite set of holomorphic map germs $\{\phi_{v}:(\mathbb{C}, 0) \rightarrow (X, o)\}_{v=1}^r$ and extensions $\{{\rm ord}_t^v:\mathbb{C}\{t\}^p\rightarrow\overline{\mathbb{Z}}\}_{v=1}^r$ of the standard order function ${\rm ord}_t:\mathbb{C}\{t\}\rightarrow\overline{\mathbb{N}}$, such that the collection $\mathcal{D}=\{ (\mathfrak{v}_{\phi_v}:={\rm ord}_t \circ \phi^*, \mathfrak{v}_{\phi_v}^M:={\rm ord}_t^v\circ \phi^*_M)\}_{v =1}^r$ satisfies $N = \mathcal{F}_{\mathcal{D}}^M(0)$. Hence, in this setting finite codimension and integral closedness  (in the sense of Gaffney) implies realizability.
\end{theorem}

We will give a proof  in paragraph \ref{par:locan}. Nevertheless, it is definitely false that every discrete valuation on the local algebra $\mathcal{O}$ (as in Definition \ref{def:intclosedmod})
is of this form:  the Hilbert function of such a valuation (cf. (\ref{eq:handhM})) jumps at most one at a step (compare with Remark \ref{rem:arcvaluations}).

It is noteworthy, however, that, even though their characterizations with our language differ, Gaffney showed in \cite[p. 305]{Gaff}, that in this analytic setting with $M=\mathcal{O}^p$, the two definitions of integral closure (in the sense of Rees respectively Gaffney) agree.

\subsection{A realization induces a CDP realization}\label{ss:CDPforD}\ 

 In this subsection we will show that if $N$ admits some realization by some finite collection 
 $\mathcal{D}$, then from $\mathcal{D}$ one can produce in a natural way several CDP realizations of $N$.

\begin{lemma}\label{lem:duplatrukk}
Let $\cO$, $M$  and $N$ be as in subsection \ref{ss:lathomofmods}.  Assume that
 $\mathcal{D}=\{(\frv_v, \frv^M_v)\}_{v=1}^r$  is a finite collection of extended discrete valuations on $M$, such that
$\mathcal{F}^M_{\mathcal{D}}(0)= N$. Then the following facts hold. 

(1)  The multiset
\begin{center}
    $\mathcal{D}^\sharp:=\{ (\frv_1,\frv^M_1),  (\frv_2,\frv^M_2), 
\ldots , (\frv_r,\frv^M_r), 
(\frv_{r+1},\frv^M_{r+1})=(\frv_1, \frv^M_1),
,\ldots , (\frv_{2r},\frv^M_{2r})=(\frv_r, \frv^M_r)\}$
\end{center}
 consisting of\, $(2r)$ extended discrete valuations 
is a CDP realization of $N$.

(2)
The collection $\mathcal{D}^\natural:=\{ (\frv_v^\natural,\frv^{\natural,M}_{v})\}_{v=1}^r$  
of extended discrete valuations,
    defined as 
    \begin{center} $(\frv_v^{\natural},\frv^{\natural,M}_{v}):=(2\frv_v,2\frv^M_v)$ for all $v \in \{1, \ldots, r\}$ (cf. Remark \ref{rem:translating} \text{(i)}),
    \end{center}
    gives a CDP realization of $N$.
\end{lemma}

\begin{proof} \textit{(1)} 
    Let us suppose indirectly, that there exist a lattice point $\ell^\sharp \in
    \Z^{2r}$
    and a coordinate direction $v^\sharp \in \mathcal{V}^\sharp=\{1, 2, \ldots, 2r\}$ such that
    \begin{equation*}
        \frh_{\mathcal{D}^\sharp}(\ell^\sharp + e_{v^\sharp}) > \frh_{\mathcal{D}^\sharp}(\ell^\sharp) \ \ \text{ and } \ \  \frh^{\circ}_{\mathcal{D}^{\sharp}}(\ell^\sharp+e_{v^{\sharp}}) < \frh^{\circ}_{\mathcal{D}^{\sharp}}(\ell^\sharp).
    \end{equation*}
This means, that there exists some $f \in \mathcal{F}_{\mathcal{D}^\sharp}(\ell^\sharp)\triangleleft \mathcal{O}$ such that 
 $\frv_{v^\sharp}(f)= \ell^\sharp_{v^\sharp}$, and, similarly,
there exists 
$m\in \mathcal{F}_{\mathcal{D}^\sharp}^M(-\ell^\sharp - e_{v^\sharp}) \leq M$ such that $\frv^{M}_{v^\sharp}(m)= -\ell^\sharp_{v^\sharp}-1$.
This implies that  
 $\frv^{M}_{w^{\sharp}}(fm)\geq 0$ for all $w^\sharp\not=v^\sharp$  and $\frv^{M}_{v^\sharp}(fm)=-1$.
However, this statement,  by the repetition of the extended valuations, leads to a contradiction, 
since $\frv^{M}_{v^\sharp}=\frv^{M}_{w^\sharp}$ for $w^\sharp=\begin{cases}
    v^\sharp+r, \text{ if } v^\sharp \leq r, \\
    v^\sharp-r, \text{ if } v^\sharp > r.
\end{cases}$

 \textit{(2)} 
 If the $e_v$-coefficient $\ell^{\natural}_v$ of $\ell^{\natural}=\sum_{w\in \mathcal{V}}\ell^{\natural}_w e_w$ is odd, then $\mathfrak{h}_{\mathcal{D}^{\natural}}(\ell^{\natural}+e_v)= \mathfrak{h}_{\mathcal{D}^{\natural}}(\ell^{\natural})$. 
 On the other hand, if $\ell^{\natural}_v$ is even, then $\mathfrak{h}^{\circ}_{\mathcal{D}^\natural}(\ell^{\natural}+e_v)= \mathfrak{h}^{\circ}_{\mathcal{D}^\natural}(\ell^{\natural})$. The CDP then automatically holds.
\end{proof}

We will sometimes denote $\mathfrak{h}_{\mathcal{D}^{\natural}}$  by $\mathfrak{h}^{\natural}$, and  $\mathfrak{h}^{\circ}_{\mathcal{D}^\natural}$ by $\mathfrak{h}^{\circ\,\natural}$.

\subsection{Reduction to a finite rectangle and Euler characteristic}\label{ss:ddep}\,

In this subsection we will prove that, given any CDP realization $\mathcal{D}$ of $N$, the lattice homology module $\mathbb{H}_*(N \hookrightarrow M)$ can be computed on a finite rectangle $R(0, d)$ (where $d$ depends on the CDP realization chosen) and it categorifies the codimension ${\rm codim}_k(N \hookrightarrow M)$.

First, notice that, by  Definition \ref{def:edv} of extended discrete valuations and Assumption \ref{ass:fincodim}, for any lattice point $d \gg0$ large enough $\mathcal{F}^M_{\mathcal{D}}(-d)=M$.

\begin{nota}\label{not:d_D}
    Associated with a (CDP) realization $\mathcal{D}$ of $N$ we will denote the unique smallest lattice point $d\geq 0$ satisfying $\mathcal{F}^M_{\mathcal{D}}(-d)=M$ by $d_\mathcal{D}\in \Z^r$ (for well-definedness use the matroid rank inequality (\ref{eq:matroid}) for $\mathfrak{h}^\circ$, cf. Observation \ref{obs:hhMmatroid}). Then for all $d \geq d_\mathcal{D}$ we have $\mathfrak{h}^\circ (d) \equiv 0$. Notice that, although this lattice point $d_{\mathcal{D}}$ is unique, in most geometric situations it is not straightforward to compute, see e.g., subsection \ref{ss:zcoh}.
\end{nota}

We now prove  that the function $w_{0}$ satisfies the matroid rank inequality, which will allow us to use Corollary \ref{cor:finrectcomb} to reduce the full lattice  to the finite rectangle $R(0, d_{\mathcal{D}})$.

\begin{proof}[Proof of Lemma \ref{lem:w0matroid} and Theorem \ref{th:properties} (a)]
    From Observation \ref{obs:hhMmatroid} and part (i) of Example \ref{ex:hhcircmatroid}, it is clear, that $w_0$ satisfies the matroid rank inequality. 

    On the other hand, 
    for any  standard basis element $e_v$ of $\mathbb{Z}^r$, 
    by the nontriviality conditions in the Definitions \ref{def:gdv} (c) and \ref{def:edv} (c) of the extended discrete valuations (applied for $\frv_v$), 
     we have  $\lim_{k \to \infty}\mathfrak{h}(ke_v) =\infty$ and $\lim_{k\to -\infty}\mathfrak{h}^\circ(ke_v) = \infty$. Moreover, for any lattice points $\ell^+\geq d_\mathcal{D}$ and $\ell^- \leq 0$ we have $\frh^\circ (\ell^+)=\frh^\circ (d_\mathcal{D})=0$ and $h(\ell^-) =h(0)=0$. Therefore, we can use Corollary \ref{cor:finrectcomb} to prove that $w_0$ is a weight function and that the lattice homology module $\mathbb{H}_*(N \hookrightarrow M)$ can be computed on any rectangle $R(0, d)$ with $d_\mathcal{D} \leq d \leq \infty$. 
    \end{proof}

\begin{proof}[Proof of Theorem \ref{th:properties} (b)] Fix any CDP realization $\mathcal{D}$ of $N$ and a lattice point $d$ with $d_\mathcal{D}\leq d\leq \infty$. Then 
by part \textit{(a)} of this Theorem \ref{th:properties}  and Theorem \ref{th:comblattice} applied for the rectangle $R(0, d)$ we get
    $eu(\bH_*(N \hookrightarrow M))=eu(\bH_*(R(0, d), w))=\frh^\circ (0)-\frh^\circ(d)$.
Hence
$eu(\bH_*(N \hookrightarrow M))=\dim_k M/N$.
\end{proof}

\begin{cor}\label{cor:rectCDP->fullCDP}
    In the setting of Theorem \ref{th:properties}, a realization $\mathcal{D}$ of $N$ is a CDP realization if and only if the corresponding height functions satisfy the CDP condition in the rectangle $R(0, d_{\mathcal{D}})$.
\end{cor}

\begin{proof}
    The nontrivial direction follows from the definition of the lattice point $d_{\mathcal{D}}$ and Proposition \ref{prop:rectCDP->fullCDP}, since the height functions satisfy the matroid rank inequality (cf. Observation \ref{obs:hhMmatroid}).
\end{proof}

\subsection{Reduction to  $M/N$. The lattice homology of $M/N$}\label{ss:MperN}\,

In this subsection we will prove part {\it (c)} of Theorem \ref{th:properties}, namely
that the lattice homology of a realizable $\mathcal{O}$-submodule $N \leq M$ only depends on the quotient module $M/N$ over the quotient algebra $\mathcal{O}/{\rm Ann}_{\mathcal{O}}(M/N)$. 

First, we fix the setting.
 Let $\mathcal{O}$ be a Noetherian $k$-algebra and assume that
 $Q$ is a finite $k$-dimensional  $\mathcal{O}$-module such that it can be realized as
$M/N$, where $M$ is a
finitely generated $\mathcal{O}$-module 
 and $N$ is realizable submodule (e.g., an integrally closed submodule).
 Since $M$ is finitely generated, there exists a surjective homomorphism $\alpha:\cO^p\to M$ from a free module of finite rank $p \in \mathbb{N}$. 
 By Lemma \ref{lem:contractionREAL} \textit{(b)} $\mathcal{N}:=\alpha ^{-1}(N)$ is realizable in $\cO^p$, hence the quotient module
 $Q$ can also be represented as $\cO^p/\mathcal{N}$, with $\mathcal{N}$ realizable.

\begin{prop}\label{prop:quotient} (1) Consider the situation as above: $N$ a finite $k$-codimensional realizable $\cO$-submodule in $M$, $\alpha:\cO^p\to M$ a surjective $\mathcal{O}$-module homomorphism, and
$\mathcal {N}:=\alpha^{-1}(N)$ in $\cO^p$. Then 
$\bH_*(N \hookrightarrow M) \cong \bH_*(\mathcal{N} \hookrightarrow \cO^p)$.

(2)
Assume that $Q$ is a finite $k$-dimensional $\mathcal{O}$-module 
such that it can be realized as
$M/N$, where $N$ is a realizable submodule in $M$.
 Then the lattice homology $\bH_*(N \hookrightarrow_{\mathcal{O}} M)$ is independent of the choice of the realization $M/N$ of $Q$, it depends merely on the quotient module $Q$. 
\end{prop}

\begin{proof}
{\it (1)} Consider a CDP realization $\mathcal{D}=\{(\frv_v, \mathfrak{v}_v^M)\}_{v=1}^r$ of $N$. Then define the collection of extended discrete valuations $\alpha^*\mathcal{D}=\{(\frv_v, \mathfrak{v}_v^M\circ \alpha) \}_{v=1}^r$ of $\cO^p$. It clearly gives a realization of $\mathcal{N} \leq \mathcal{O}^p$.
Moreover, as $\alpha$ is surjective, the  height functions and, hence, the weight functions associated with the two setups agree on the infinite rectangle $R(0, \infty)$ (hence, by Corollary \ref{cor:rectCDP->fullCDP},  $\alpha^*\mathcal{D}$ is a CDP realization of $\mathcal{N}$, too). Therefore, by the Independence Theorem \ref{th:IndepMod} and Remark \ref{rem:nodoublingifCDP}, we can use the same lattice to compute the lattice homology modules, which, by Theorem \ref{th:properties} \textit{(a)}, agree.

{\it (2)}
By part {\it (1)} we can assume that the larger module $M$ is free of finite rank. Suppose now that we have $Q=\cO^{p'}/\mathcal{N}'=\cO^{p}/\mathcal{N}$ two different realizations. Then 
the identity map $Q\to Q$
lifts to a homomorphism $\beta: \cO^{p'}\to \cO^{p}$ such that $\beta^{-1}(\mathcal{N})=\mathcal{N}'$ (we stress, however, that $\beta$ does not have to be surjective).
Now if $\mathcal{D}=\{(\frv_v, \mathfrak{v}_v^{\mathcal{O}^p})\}_v$ is a CDP realization of $\mathcal{N}$ in $\mathcal{O}^p$, then the collection
$\mathcal{D}':=\beta^*\mathcal{D}=\{(\mathfrak{v}_v, \frv_v^{\mathcal{O}^p}\circ \beta)\}_v$ is a realization of $\mathcal{N}'$ in $\mathcal{O}^{p'}$.  
Moreover, we claim, that the corresponding height and weight functions agree on the
infinite rectangle $R(0, \infty)$, which is really the essential part of the lattice (since this already implies that $\mathcal{D}'$ is a CDP realization and the lattice homology modules agree as in part \textit{(1)}). Indeed, $\beta$ identifies the corresponding filtrations and Hilbert functions meaning that for any $\ell \in (\mathbb{Z}_{\geq 0})^r$ we have $\mathfrak{h}_{\mathcal{D}}(\ell) = \mathfrak{h}_{\mathcal{D}'}(\ell)$ and
\begin{align*}
    \mathfrak{h}^\circ_{\mathcal{D}'}(\ell) =&\ \dim_k \dfrac{\mathcal{O}^{p'}}{\beta^{-1}(F_{\mathcal{D}}^{\mathcal{O}^{p}}(-\ell))} = \dim_k \dfrac{(\mathcal{O}^{p'}/\mathcal{N}')}{\beta^{-1}(F_{\mathcal{D}}^{\mathcal{O}^{p}}(-\ell))/\beta^{-1}(F_{\mathcal{D}}^{\mathcal{O}^{p}}(0))}=\\
    =&\ \dim_k \dfrac{\mathcal{O}^p/\mathcal{N}}{\mathcal{F}_{\mathcal{D}}^{\mathcal{O}^p}(-\ell)/\mathcal{F}_{\mathcal{D}}^{\mathcal{O}^p}(0)}=\mathfrak{h}^\circ_{\mathcal{D}}(\ell).
\end{align*}
Thus, just as before, 
we get $\bH_*(\mathcal{N}'\hookrightarrow \cO^{p'})\cong\bH_*(\mathcal{N}\hookrightarrow \cO^{p})$.
\end{proof}

On the other hand, as $\mathcal{O}$ is a Noetherian $k$-algebra, there exists a surjective morphism $\gamma:\cO_m\to \cO$ for some $m \in \mathbb{N}$,
 where $\cO_m$ denotes the polynomial ring $k[x_1,\ldots, x_m]$ over the field $k$. Then clearly every $\mathcal{O}$-module $M$ is also a $\cO_m$-module via $\gamma$.
  [In the analytic\,/\,complete case the same statements remain true with $\mathcal{O}_m=\mathbb{C}\{x_1, \ldots, x_m\}$ or $\mathcal{O}_m=k[[x_1, \ldots, x_m]]$.]

\begin{proposition}\label{prop:pullbackring}
    (1) 
    Consider the above situation: let $N$ be a finite $k$-codimensional realizable $\mathcal{O}$-submodule in $M$ and $\gamma: \mathcal{O}_m \rightarrow \mathcal{O}$ a surjective ring morphism. Then the lattice homology of $N$ in $M$ is independent of the base ring, that is,  $\mathbb{H}_*(N \hookrightarrow_{\mathcal{O}} M) \cong \mathbb{H}_*(N \hookrightarrow_{\mathcal{O}_m} M)$.\\
\noindent [Notice also, that in this setting $\gamma$ induces an isomorphism $\mathcal{O}_m/{\rm Ann}_{\mathcal{O}_m}(M/N) \cong \mathcal{O}/{\rm Ann}_{\mathcal{O}}(M/N)$.]

    (2)
    Given any realization $\mathcal{D}$ of $N$, we have the following 
    (for notation see \ref{not:d_D}): 
    \begin{equation*}
        \mathcal{F}_{\mathcal{D}}(d_{\mathcal{D}})={\rm Ann}_{\mathcal{O}}(M/N)= (N:_\mathcal{O} M)=\{ f \in \mathcal{O}\,:\, fm \in N \ \ \forall\, m \in M\} \triangleleft \mathcal{O}.
    \end{equation*}
\end{proposition}

\begin{proof}
    {\it (1)} Consider a CDP realization $\mathcal{D}=\{(\frv_v, \mathfrak{v}_v^M)\}_{v=1}^r$ of $N\leq M$ over $\mathcal{O}$, i.e., $\mathcal{F}_{\mathcal{D}}^M(0)=N$. Then define the collection of extended discrete valuations $\gamma^*\mathcal{D}:=\{(\frv_v \circ \gamma, \mathfrak{v}_v^M) \}_{v=1}^r$ of the $\mathcal{O}_m$-module $M$. They satisfy the axioms of Definition \ref{def:edv} and give a realization of the $\mathcal{O}_m$-submodule $N$.
Moreover, as $\gamma$ is surjective, the  Hilbert functions and, hence, the weight functions associated with the two setups agree (and thus $\gamma^*\mathcal{D}$ is a CDP realization) and we can use the Independence Theorem \ref{th:IndepMod}. 

{\it (2)} Recall that, for a given realization $\mathcal{D}=\{(\frv_v, \mathfrak{v}_v^M)\}_{v\in \mathcal{V}}$ of $N$, the element $d_{\mathcal{D}}$ is defined as the smallest lattice point satisfying $\mathcal{F}_{\mathcal{D}}^M(-d_{\mathcal{D}})=M$. By the definition of the filtration $\mathcal{F}_{\mathcal{D}}^M$ (cf. (\ref{eq:fdMl})) this means that for every $v \in \mathcal{V}$ the coordinate $d_{\mathcal{D},v}$ is defined uniquely by the containments
\begin{equation*}
    \mathcal{F}_v^M(-d_{\mathcal{D}, v}) =M, \text{ whereas } \mathcal{F}_v^M(-d_{\mathcal{D}, v}+1) \subsetneq M, \text{ i.e., } -d_{\mathcal{D}, v}=\min \{ \mathfrak{v}_v(m)\,:\, m \in M\}.
\end{equation*}
From Definition \ref{def:edv} of the extended discrete valuations and the fact that $N=\mathcal{F}_{\mathcal{D}}^M(0)$, it is now clear that $\mathcal{F}_{\mathcal{D}}(d_{\mathcal{D}})={\rm Ann}_{\mathcal{O}}(M/N)$. 
\end{proof}

We are now ready to prove part \textit{(c)} of Theorem \ref{th:properties} stating that the lattice homology $\mathbb{H}_*(N \hookrightarrow_{\mathcal{O}} M)$ only depends on the  $\mathcal{O}/{\rm Ann}_{\mathcal{O}}(M/N)$-module structure of the quotient module  $M/N$.

\begin{proof}[Proof of Theorem \ref{th:properties} (c)]
Suppose that we have two different realizations of the $\mathcal{O}/{\rm Ann}_{\mathcal{O}}(M/N)$-module $M/N$ with $N$ realizable. By the previous Propositions \ref{prop:quotient} and \ref{prop:pullbackring} it is enough to consider the follo\-wing simplified case: let
$\mathcal{O}_m$ and $\mathcal{O}_{m'}$ be two polynomials rings over $k$, and consider two finite $k$-codimensional realizable submodules $\mathcal{N} \leq \mathcal{O}_m^p$ and $\mathcal{N}'\leq \mathcal{O}_{m'}^{p'}$, such that we are given
\begin{equation*}\begin{split}
   \text{ a ring isomorphism }& \phi:\mathcal{O}_{m'}/(\mathcal{N}' :_{\mathcal{O}_{m'}}\mathcal{O}_{m'}^{p'}) \rightarrow \mathcal{O}_{m}/(\mathcal{N}:_{\mathcal{O}_{m}}\mathcal{O}_{m}^{p}), \ \ \text{ and }\\
  \text{ a }\phi\text{-equivariant module isomorphism } & \psi:\mathcal{O}_{m'}^{p'}/\mathcal{N}' \rightarrow \mathcal{O}_{m}^{p}/\mathcal{N}.
\end{split}\end{equation*}
    We have to prove $\mathbb{H}_*\big(\mathcal{N} \hookrightarrow_{\mathcal{O}_m} \mathcal{O}_m^p\big) \cong \mathbb{H}_*\big(\mathcal{N}' \hookrightarrow_{\mathcal{O}_{m'}} \mathcal{O}_{m'}^{p'}\big)$ in this setting.

    First we choose a lift $\Phi: \mathcal{O}_{m'} \rightarrow \mathcal{O}_{m}$ of $\phi$ with $\Phi^{-1}\big( (\mathcal{N} :_{\mathcal{O}_m}\mathcal{O}_m^p) \big)= (\mathcal{N'}:_{\mathcal{O}_{m'}}\mathcal{O}_{m'}^{p'})$ and then a $\Phi$-equivariant lift $\Psi: \mathcal{O}_{m'}^{p'} \rightarrow \mathcal{O}_m^{p}$ of $\psi$ with $\Psi^{-1}(\mathcal{N}) = \mathcal{N}'$. These exist due to the projectivity of polynomial rings and free modules. Then if we consider any CDP realization $\mathcal{D}=\{(\mathfrak{v}_v, \mathfrak{v}_v^{\mathcal{O}_m^p})\}_{v \in \mathcal{V}}$ of $\mathcal{N}$, then $\mathcal{D}':=\{(\mathfrak{v}_v\circ \Phi, \mathfrak{v}_v^{\mathcal{O}_m^p}\circ\Psi)\}_{v \in \mathcal{V}}$ will give a realization of $\mathcal{N}'$. Therefore, provided that we can prove that $\mathcal{D}'$ is a CDP realization,  we can use the same lattice $\mathbb{Z}^{|\mathcal{V}|}$ to compute the lattice homology modules. 
    Moreover, since $\psi$ is an isomorphism, we have $d_{\mathcal{D}}=d_{\mathcal{D}'}$, so we can use, in fact, the same rectangle $R(0, d_{\mathcal{D}})$ as well (see Theorem \ref{th:properties} \textit{(a)}). Now, similarly to the proof of part \textit{(2)} of Proposition \ref{prop:quotient} we get that for any $\ell \in \mathbb{Z}^{|\mathcal{V}|}, \ \ell \leq d_{\mathcal{D}}$ we have $\mathfrak{h}_{\mathcal{D}}(\ell) = \mathfrak{h}_{\mathcal{D}'}(\ell)$ (here we use part \textit{(2)} of Proposition \ref{prop:pullbackring} to be able to rely on the isomorphism $\phi$) and for any $\ell \in (\mathbb{Z}_{\geq 0})^{|\mathcal{V}|}$ we have $\mathfrak{h}_{\mathcal{D}}^\circ(\ell) = \mathfrak{h}_{\mathcal{D}'}^\circ(\ell)$. Thus, $\mathcal{D}'$ is a CDP realization (by Corollary \ref{cor:rectCDP->fullCDP}) and on the essential rectangle $R(0, d_{\mathcal{D}})$ the weight functions $w_{\mathcal{D}}$ and $w_{\mathcal{D}'}$ agree, hence the lattice homology modules corresponding to the different descriptions of the quotient module  $\mathcal{O}_{m}^{p}/\mathcal{N}$, with $\mathcal{N}$ realizable, agree.
\end{proof}

\begin{define}
    Let $\mathcal{O}$ be a Noetherian $k$-algebra and $Q$ a finite $k$-dimensional $\mathcal{O}$-module such that $Q$ can be written as $M/N$ with $M$ finitely generated and $N$ realizable. Then the lattice homology $\mathbb{H}_*(N \hookrightarrow M)$ only depends on the $\mathcal{O}/{\rm Ann}_{\mathcal{O}}(Q)$-module $Q$ (see the previously proved
    Theorem \ref{th:properties} \textit{(c)}). We refer to it as the {\it 
    lattice homology of $Q$}, 
    and it will be denoted by $\bH_*(Q)$.
\end{define}

 By Theorem \ref{th:properties} \textit{(b)}, $eu(\bH_*(Q))=\dim(Q)$ holds.
    In other words, the matching $Q\rightsquigarrow \bH_*(Q)$ is the categorification of the numerical invariant $Q\rightsquigarrow \dim_k Q$.
 
Notice, that, in the setting of the previous definition,  the ideal ${\rm Ann}_{\mathcal{O}}(Q)$ is integrally closed (use Proposition \ref{prop:pullbackring} \textit{(2)} and Corollary \ref{cor:intersection}). Therefore, the Artin $k$-algebra $\mathcal{O}/{\rm Ann}_{\mathcal{O}}(Q)$ is of form $\mathcal{O}_m/\mathcal{I}$ with $\overline{\mathcal{I}}=\mathcal{I}$. For more on these types of algebras see subsection \ref{ss:ARTIN}.

\subsection{Bounds on the gradings}\label{ss:bounds} \,

In this subsection we will present some natural bounds on the gradings of the nonzero elements of the lattice homology of realizable submodules.
Fix a finite $k$-codimensional realizable submodule $N$ in the finitely generated $\mathcal{O}$-module $M$, where $\mathcal{O}$ is a Noetherian $k$-algebra. Recall that  its lattice homology 
has the following structure: $\mathbb{H}_*(N \hookrightarrow M) =\oplus_{q \geq 0}\oplus_{n \in \mathbb{Z}}\mathbb{H}_{q, 2n}$, where $q$ denotes the homological grading and $2n$ denotes the weight grading. 

\bekezdes\textbf{Homological grading.} \label{bek:homdeg} It is convenient and natural to introduce the following terminology: the {\it (lattice) homological dimension } of $N\leq M$ (or of 
$M/N$) is
\begin{equation}\label{eq:defhomdim}
    {\rm homdim}(N \hookrightarrow M):=\begin{cases}
        \max \{ q\,:\, \mathbb{H}_{q}(N \hookrightarrow M)\not=0\} & \text{if }N\neq M;\\
        -1 & \text{if }N= M.
    \end{cases}
\end{equation}

\begin{remark}
    Equivalently, ${\rm homdim}(N \hookrightarrow M)=
        \max \big(\{ q\,:\, \mathbb{H}_{{\rm red},q}(N \hookrightarrow M)\not=0\}\cup\{-1\})$. Indeed, if $N\neq M$, then $\mathbb{H}_{{\rm red}, *}\neq 0$, since $eu(\mathbb{H}_*)\neq 0$ and, using any realization, the lattice point $0$, with $w(0)=0$, gives a distinct connected component of $S_0$ (by a similar argument to Remark \ref{rem:w(R(0,E)<kappa}).
\end{remark}

If $\mathcal{D}$ is a 
realization of $N$ with cardinality $r\geq1$, then by Lemma \ref{lem:duplatrukk} {\it (2)} we have a CDP realization $\mathcal {D}^\natural$ with $|\mathcal{D}^\natural|=r$.
Then each $S_n$ space associated with $\mathcal{D}^\natural$  is a compact subspace of $\R^r$, hence  ${H}_{\geq r}(S_n,\Z)=0$ for all $n \in \mathbb{Z}$, that is, 
$\oplus_{q\geq r}\bH_{q}(N \hookrightarrow M)=0$. [The case of $\mathcal{D}=\emptyset$ corresponds to $N=M$ and ${\rm homdim}(N \hookrightarrow M)=-1$.] We therefore have the following  bound:

\begin{prop}\label{prop:homdegupperbound}
    ${\rm homdim}(N \hookrightarrow M) \leq \min \{|\mathcal{D}|\,:\, \mathcal{D} \text{ a realization of } N\} -1$.
\end{prop}
In particular, if we can find a 
realization $\mathcal{D}$ consisting of only one element, then $\bH_*$ is concentrated merely in $\bH_0$ (this is the case, e.g., for the analytic lattice homology of irreducible curve singularities or almost rational surface singularities, see sections \ref{s:deccurves} and \ref{s:dnagy}). 

From the other point of view, this inequality can be also seen as a lower bound for the numerical invariant $ \min \{ |\mathcal{D}|\,:\, \mathcal{D} \text{ a realization of } N\} $ associated with a pair $N\leq M$. This invariant is hardly nontrivial. In the case of realizable (or, equivalently, integrally closed, cf. Corollary \ref{cor:intersection}) ideals $\mathcal{I}=N\leq M=\cO$, 
it is usually strictly smaller than the number of Rees valuations (see Example \ref{bek:10.2.2}), which are uniquely associated to the ideal $\mathcal{I}$ (cf. \cite{Rees}). 

Regarding the inequality in Proposition \ref{prop:homdegupperbound}, it is not clear how sharp it actually is.
In the case of finite codimensional integrally closed ideals, in all our computations with  ${\rm homdim}(\mathcal{I} \hookrightarrow \cO)=0$ we succeeded to find  a single valuation which realizes $\mathcal{I}$. This encourages us to make the following conjecture:

\begin{conj}\label{conj:elso}
    Let us suppose that $\mathcal{I} \triangleleft \mathcal{O}_2=k[x_1, x_2]$ is a finite codimensional, integrally closed ideal. We conjecture that if ${\rm homdim}(\mathcal{I} \hookrightarrow \mathcal{O})=0$, then $\mathcal{I}$ can be realized by a single valuation, i.e., $\mathcal{I}=\{ f \in \mathcal{O}_2\,:\, \mathfrak{v}(f) \geq d_{\mathcal{I}}\}$ for some discrete valuation $\mathfrak{v}$ and integer $d_{\mathcal{I}}$.
\end{conj}

In fact, in the case of monomial ideals $\mathcal{M} \triangleleft \mathcal{O}_2=k[x_1, x_2]$ the conjecture holds true, even more, the following much stronger converse of Proposition \ref{prop:homdegupperbound} is also true:

\begin{theorem}\label{th:homdimmon}
    For a finite codimensional integrally closed monomial ideal $\mathcal{M} \triangleleft \mathcal{O}_2$:
    \begin{equation*}
        {\rm homdim}(\mathcal{M} \hookrightarrow \mathcal{O}_2) = \min \{|\mathcal{D}|\,:\, \mathcal{D} \text{ a realization of } \mathcal{M}\} -1.
    \end{equation*}
\end{theorem}

The proof and some more detailed discussions will be given in section \ref{s:homdim}. 
For more supporting examples for Conjecture \ref{conj:elso} see also subsections \ref{ss:NNexamples}
and \ref{ss:HOMDIM}.

\bekezdes\textbf{Weight grading.} By the definition of the weight function, for any given CDP realization $\mathcal{D}$ of $N$, the set of weight values of $w_{\mathcal{D}, 0}$ is bounded from below and the minimum is denoted by $m_{w_{\mathcal{D}}}$. Moreover, this value can be read off from the bigraded lattice homology $\mathbb{Z}[U]$-module $\mathbb{H}_*(\mathbb{R}^{|\mathcal{D}|}, w_{\mathcal{D}})$ as $m_{w_{\mathcal{D}}}=-\max\{ n \in \mathbb{Z} \,:\, \mathbb{H}_{0,2n} \neq 0\}$ (see the grading convention in paragraph \ref{9complexb}). Therefore, as a consequence of the Independence Theorem \ref{th:IndepMod} we have the following:

\begin{cor}\label{cor:minw}
    For any CDP realization $\mathcal{D}$ of $N$, the minimum value $m_{w_{\mathcal{D}}}=\min w_{\mathcal{D}}$ of the corresponding weight function is independent of the realization chosen, it is merely an invariant of the realizable submodule $N$. 
    It can be read from the lattice homology module as well:
    \begin{center}
    $m_{w_{\mathcal{D}}}=\min w_{\mathcal{D}}=-\max\{ n \in \mathbb{Z} \,:\, \mathbb{H}_{0,2n}(N \hookrightarrow M) \neq 0\}$.
    \end{center}
\end{cor}

We would like to highlight, however, that the computation of $m_w$  and its `geometric' interpretation seems nontrivial, even in the case of the previously defined lattice (co)homology theories (see e.g., \cite{HofN} for some results in the case of analytic lattice homology of reduced curve singularities). In comparison with the Heegaard Floer theory, it is the analogue of the $d$-invariant (cf. \cite{OSzF}).

On the other end, Theorem \ref{th:properties} \textit{(a)} implies that, for $n \gg 0$ large enough, the $S_n$ spaces are strong deformation retractable to the rectangle $R(0, d)$, i.e., they are all contractible, hence $\mathbb{H}_{{\rm red},-2n} = 0$ for $n\gg 0$. In the case of an algebraically closed field and of a local base ring (which is the case for most applications below) we also have the following strengthening:

\begin{theorem}[Nonpositivity Theorem]\label{th:upperbound}
    Let $k$ be an algebraically closed field, $(\mathcal{O}, \mathfrak{m})$ a local Noetherian $k$-algebra, $M$ a finitely generated $\mathcal{O}$-module and $N \leq M$ a finite $k$-codimensional realizable submodule. Then for every positive integer $n >0$, the $S_n$ space is contractible (this is highlighted in the name of the theorem), hence the weight-grading of the \emph{reduced} lattice homology module is nonnegative,  i.e.,
    $\mathbb{H}_{{\rm red},q, -2n}(N \hookrightarrow M) = 0$
    for all $q\geq 0$ and $n >0$. In other words,
\begin{equation*}
    \mathbb{H}_{{\rm red}, q}(N \hookrightarrow M) = \bigoplus_{0 \geq n \geq m_w}\mathbb{H}_{{\rm red},q, -2n}(N \hookrightarrow M) \hspace{5mm} \text{ for all } q \geq 0.
\end{equation*}
Moreover, if $N \neq M$, the upper bound is sharp: $S_{0}$ is not connected.
\end{theorem}

The proof will be given in section \ref{s:nonpositivity}. The name `Nonpositivity Theorem' comes from the original version \cite[Theorem 6.1.1]{KNS2} presented in the case of the analytic lattice \emph{co}homology of reduced complex analytic curve singularities. In fact, our proof is just a direct generalization of the original argument presented in [loc. cit.].

\section{Symmetric lattice homology of integrally closed ideals}\label{s:ideals}

In this section we will discuss the case when $M=\mathcal{O}$ and $N=\mathcal{I} \triangleleft \mathcal{O}$
is a finite $k$-codimensional ideal in it. In this case, the realizability of $\mathcal{I}$ is equivalent with integral closedness, moreover,  
the weight function corresponding to any realization is symmetric. 

\subsection{Extending valuations to a rank one free module}\label{ss:sym1}\,

Let us fix a field $k$ and a Noetherian $k$-algebra $\mathcal{O}$. If we consider the rank one free module $M=\mathcal{O}$, then  every discrete valuation $\mathfrak{v}:\mathcal{O} \rightarrow \overline{\N}$ extends to $M$  (see Remark \ref{rem:translating} (i)). 
In fact, these extensions can be easily classified. 

\begin{lemma}
    If we consider $\mathcal{O}$, as a finitely generated module $M$ over itself, and  a discrete valuation $\mathfrak{v}:\mathcal{O} \rightarrow \overline{\N}$, then its  extensions $(\mathfrak{v}, \mathfrak{v}^M)$ are precisely pairs of form $(\mathfrak{v}, \mathfrak{v}+n')$ for some $n' \in \mathbb{Z}$. 
\end{lemma}

\begin{proof}
    Clearly $\mathfrak{v}^M(1)=n'$ for some $n' \in \mathbb{Z}$ and part (a) of Definition \ref{def:edv} describes $\mathfrak{v}^M$ uniquely.
\end{proof}

The corresponding filtrations (cf. (\ref{eq:fvMl})) then naturally satisfy $\mathcal{F}_{\mathfrak{v}}^M(\ell) = \mathcal{F}_{\mathfrak{v}}(\ell -n') \triangleleft \mathcal{O}$.

\begin{cor}\label{cor:idealREAL}
    An ideal $\mathcal{I} \triangleleft\mathcal{O}$ is realizable in the sense of Definition \ref{def:REAL} (with $M=\mathcal{O}$ and $N = \mathcal{I}$), if and only if there exists some finite collection $\mathcal{D}:=\{\mathfrak{v}_1, \ldots, \mathfrak{v}_r\}$ of discrete valuations on $\mathcal{O}$  and a lattice point $d \in (\mathbb{Z}_{\geq 0})^r$ such that $\mathcal{I}=\mathcal{F}_{\mathcal{D}}(d)$. (Here the name `collection' refers to the fact that all these valuations respect Assumption \ref{ass:fincodim})
\end{cor}

\begin{remark}
    Notice that the nonnegativity conditions $d_v \geq 0 \ (\text{for all }\ 1 \leq v \leq r)$ on the components of the lattice point $d=\sum_{v=1}^rd_ve_v$ are, in fact, unnecessary; Corollary \ref{cor:idealREAL} remains true without them. Indeed, the negative components $d_w \leq 0$ contribute only trivial $\mathcal{F}_{\mathfrak{v}_w}(d_w)=\mathcal{O}$ terms to the intersection $\mathcal{I}=\cap_{v=1}^r\mathcal{F}_{\mathfrak{v}_v}(d_v)$.
\end{remark}

This reformulation allows us to compare the realizability property with \emph{integral closedness}. This latter one has a much more intrinsic definition for ideals than for modules. Indeed, an element $r \in \mathcal{O}$ is called  \textit{integral over $\mathcal{I} \triangleleft \mathcal{O}$}, if there exists an integer $n \in \mathbb{Z}_{>0}$ and elements $a_i \in \mathcal{I}^{i},\  i \in \{ 1, \ldots, n\}$,
such that
\begin{equation}\label{eq:intdep}
    r^n + a_1r^{n-1} + a_2r^{n-2} + \ldots + a_{n-1} r + a_n = 0.
\end{equation}
The \emph{integral closure of $\mathcal{I}$} is the ideal $\overline{\mathcal{I}}$ containing all such elements, and $\mathcal{I}$ is called \emph{integrally closed} if  $\mathcal{I}=\overline{\mathcal{I}}$. The equation (\ref{eq:intdep}) of integral dependence implies that for any discrete valuation $\mathfrak{v}:\mathcal{O} \rightarrow \overline{\mathbb{N}}$ (in the sense of Definition \ref{def:gdv}) and integer $\ell \in \mathbb{Z}$, the ideal $\mathcal{F}_{\mathfrak{v}}(\ell)$ is integrally closed. Moreover, the intersections of integrally closed ideals or their preimages under ring maps are also integrally closed. 
On the other hand, we do not have a similar property for the images of integrally closed ideals:
\begin{remark}\label{rem:imageofintclosed}
    The image of an integrally closed ideal under a ring map is not (necessarily) integrally closed, not even under surjective maps. For example, the ideal $(t^m)$ is integrally closed in the one-variable polynomial ring $k[t]$, but its image $0\triangleleft k[t]/(t^m)$ under the natural quotient homomorphism $q:k[t]\to k[t]/(t^m)$ is not.
\end{remark}

For more details on integral closure of ideals consult \cite{SH}.

\begin{cor}\label{cor:intersection}
Let $\cO$ be a Noetherian ring with a $k$-algebra structure and ${\mathcal {I}}$ a finite codimensional ideal in it. Then $\mathcal{I}$ is realizable (i.e., by Corollary \ref{cor:idealREAL}, there exists a finite collection 
$\mathcal{D}=\{\frv_1,\ldots, \frv_r\}$ of discrete valuations and a lattice point $d\in\Z^r$ such that $\mathcal{I}=\mathcal{F}_{\mathcal{D}}( d)$) if and only if it is integrally closed.
\end{cor}
\begin{proof} 
 Clearly, if $\mathcal{I}=
\mathcal{F}_{\mathcal{D}}(d)$ for some collection $\mathcal{D}$ and lattice point $d$ as in the statement, 
then $\mathcal{I}=\cap_{v \in \mathcal{V}}\mathcal{F}_{\mathfrak{v}}(d_{v})$ is integrally closed in $\calO$.

For the other direction, if $\mathcal{I}$ is finite codimensional and integrally closed, it is enough to find any (not necessarily finite) set $\mathcal{D}'$ of valuations which together realize the ideal (finite codimensionality already implies finiteness of any minimal collection $\mathcal{D}'$, if there exists such).
However, \cite[Theorem 6.8.3]{SH} (translated to our language) states exactly that for an integrally closed ideal $\mathcal{I}$ the set of all discrete valuations $\mathfrak{v}:\mathcal{O} \rightarrow \mathbb{N}\cup \{ \infty\}$, with $\mathfrak{p}_{\mathfrak{v}}=\{ f \in \mathcal{O}\,:\,\mathfrak{v}(f)=\infty\}$ being a minimal prime ideal of $\mathcal{O}$, satisfies the above requirement. Hence any finite codimensional integrally closed ideal is realizable. [Another, in principle more involved, proof would be to simply point to the existence of Rees valuations in Noetherian rings, see \cite{Rees} or \cite[Theorem 10.2.2]{SH}.]
\end{proof}

\begin{prop}\label{prop:symforideals}
    Let $\mathcal{O}$ be a Noetherian $k$-algebra and $\mathcal{I}$ a finite $k$-codimensional, integrally closed ideal in it. Let us regard $M:= \mathcal{O}$ and $N:=\mathcal{I}$ as finitely generated $\mathcal{O}$ modules. Then, for any realization $\mathcal{D}=\{(\mathfrak{v}_1, \mathfrak{v}_1-d_1), \ldots, (\mathfrak{v}_r, \mathfrak{v}_r-d_r)\}$ of $N$, the corresponding weight function $w_{\mathcal{D}, 0}$ is symmetric with respect to $d=\sum_v d_v e_v\in \mathbb{Z}^r$, i.e., for any $\ell \in \mathbb{Z}^r$ we have
    \begin{equation*}
        w_{\mathcal{D}, 0}(\ell)=w_{\mathcal{D}, 0}(d-\ell).
    \end{equation*}
    In fact, already the corresponding height functions satisfy $\mathfrak{h}_{\mathcal{D}}^\circ(\ell) = \mathfrak{h}_{\mathcal{D}, d}^{sym}(\ell) = \mathfrak{h}_{\mathcal{D}}(d-\ell)$ for every lattice point $\ell \in \mathbb{Z}^r$ (for the notations see (\ref{eq:handhcirc0})). 
\end{prop}

\begin{proof}
    By the form of the extensions of discrete valuations $\mathcal{F}_{\mathcal{D}}^M(\ell) = \mathcal{F}_{\mathcal{D}}(\ell+d)$, and, thus, we have  $\mathfrak{h}^{\circ}_{\mathcal{D}}(\ell)=\dim_kM/\mathcal{F}_{\mathcal{D}}^M(-\ell)=\dim_k \mathcal{O}/\mathcal{F}_{\mathcal{D}}(d-\ell) = \mathfrak{h}_{\mathcal{D}}(d-\ell) $ for every $\ell \in \mathbb{Z}^r$. 
\end{proof}

\begin{prop}\label{prop:symforoffsetideals}
    Let us consider the setting of the previous Proposition \ref{prop:symforideals} and introduce the notation $d_+:=\max\{d, 0\} \in (\mathbb{Z}_{\geq 0})^r$. Then the weight function $w_{\mathcal{D}, 0}$ is symmetric in the rectangle $R(0, d_+)$, i.e., for any $\ell \in R(0, d_+) \cap \mathbb{Z}^r$ we have $w_{\mathcal{D}, 0}(\ell)=w_{\mathcal{D}, 0}(d_+-\ell)$. 
    
    Moreover, if $\mathcal{D}$ additionally satisfies the Combinatorial Duality Property in the rectangle $R(0, d_+)$, then it is a CDP realization of $N$ with the lattice point $d_{\mathcal{D}}$  satisfying $0 \leq d_{\mathcal{D}} \leq d_+$ (cf. Notation \ref{not:d_D}). Hence, by Theorem \ref{th:properties} \textit{(a)}, the rectangle $R(0, d_+)$ with symmetric weight function contains all the data required to compute $\mathbb{H}_*(N \hookrightarrow_{\mathcal{O}}M)$.
\end{prop}

\begin{proof}
    By the previous Proposition \ref{prop:symforideals} we know that for any lattice point $\ell \in R(0, d_+) \cap \mathbb{Z}^r$ we have $w_{\mathcal{D}, 0}(\ell)=w_{\mathcal{D}, 0}(d-\ell)$. Thus, we only need to prove that $w_{\mathcal{D}, 0}(d-\ell)=w_{\mathcal{D}, 0}(\max\{d-\ell, 0\})$, with $\max\{d-\ell, 0\}=d_+-\ell$.
    However, the identities from (\ref{eq:0d-stabilization}), with the substitutions $d^-=0$ and $d^+=d$, imply that $\mathfrak{h}_{\mathcal{D}}(d-\ell)=\mathfrak{h}_{\mathcal{D}}(\max\{d-\ell, 0\})$ and $\mathfrak{h}^{\circ}_{\mathcal{D}}(\max\{d-\ell, 0\})=\mathfrak{h}^{\circ}_{\mathcal{D}}(\min\{\max\{d-\ell, 0\}, d\})$, with $\min\{\max\{d-\ell, 0\}, d\}=d-\ell$.

    If, additionally, $\mathcal{D}$ satisfies the Combinatorial Duality Property in the rectangle $R(0, d_+)$, then, by Proposition \ref{prop:rectCDP->fullCDP} (with $d^-=0$ and $d^+=d_+$), it is a CDP realization of $N$. Moreover, 
    \begin{center} 
    $\mathcal{F}_{\mathcal{D}}^M(-d_+) \supseteq  \mathcal{F}_{\mathcal{D}}^M(-d)=\mathcal{F}_{\mathcal{D}}(0)=\mathcal{O} =M$ with $d_+ \geq 0$,
    \end{center}
    hence $d_+ \geq d_{\mathcal{D}}$, by the latter's definition.
\end{proof}

Therefore, the lattice homology module $\mathbb{H}_*(\mathcal{I} \hookrightarrow_{\mathcal{O}} \mathcal{O})$ (defined in full generality in Definition \ref{def:LCofMOD}) of a finite codimensional integrally closed ideal 
always has symmetric weight function. For the sake of clarity, we will describe its construction in the next subsection from the more natural viewpoint of a collection $\mathcal{D}$ of discrete valuations and a dualizing element $d \in (\mathbb{Z}_{\geq 0})^r$. 

\subsection{Review of the construction in the symmetric case}\,

Assume that $\calO$ is a Noetherian ring with a $k$-algebra structure. 
The construction that we wish to present might have two different starting points. In the first case, we start from a finite collection of discrete valuations on $\mathcal{O}$ and a lattice point and assign to them an ideal and a symmetric lattice homology module (this will be the case, for example, in the construction of the combinatorial lattice homology assigned to a Newton boundary in paragraph \ref{par:NewtonComb}). In the second case, we start from the integrally closed ideal and consider any realization of it to study the resulting lattice homology (this will be the case in all analytic lattice homology constructions from sections \ref{s:deccurves} and \ref{s:dnagy}).

\begin{define} \label{def:SH(O,D,d)}
Consider a finite collection $\mathcal{D}=\{\frv_1,\ldots, \frv_r\}$ of discrete valuations on $\calO$ and set $\mathcal{V}:=\{ 1, \ldots, r\}$. 
They induce the multifiltration  
\begin{equation}\label{eq:F_D(l)}
    \mathbb{Z}^r=\mathbb{Z}\langle\{e_v\}_{v \in \mathcal{V}} \rangle\ni \ell=\sum_v \ell_v e_v \mapsto \mathcal{F}_{\mathcal{D}}(\ell)=\{ f \in \mathcal{O}\,:\, \mathfrak{v}_v(f) \geq \ell_v \ \ \mbox{for all} \  v \in \mathcal{V}\} 
\end{equation}
and Hilbert function 
$\frh_{\mathcal{D}}(\ell)=\dim_k\, \calO/\mathcal{F}_{\mathcal{D}}(\ell)$ as in (\ref{eq:fdMl}) and (\ref{eq:handhM}). Recall, that the values of this function are finite by Assumption \ref{ass:fincodim}. If we now fix some lattice point $d\in (\Z_{\geq 0})^r$, then we can associate to it the integrally closed ideal $\mathcal{I}:=\mathcal{F}_{\mathcal{D}}(d)$. 
We can now consider the rectangle $R(0, d)$ and the following symmetric weight function on it:
\begin{equation*}
    R(0, d) \cap \mathbb{Z}^r \ni \ell \mapsto w_{\mathcal{D}, 0}(\ell) =\mathfrak{h}_{\mathcal{D}}(\ell) + \mathfrak{h}^{sym}_{\mathcal{D},\, d}(\ell) - \mathfrak{h}^{sym}_{\mathcal{D},\, d}(0)=\mathfrak{h}_{\mathcal{D}}(\ell) + \mathfrak{h}_{\mathcal{D}}(d-\ell) - \dim_k \mathcal{O}/\mathcal{I}. 
\end{equation*}
This is then extended to higher dimensional cubes by (\ref{eq:9weight}) and the $S_{n, \mathcal{D}}$ spaces are defined as the sublevel cubical complexes as in paragraph \ref{9complexb}. We denote the resulting $\mathbb{H}_*(R(0, d), w_{\mathcal{D}})$ lattice homology module by $\mathbb{SH}_*(\mathcal{O}, \mathcal{D}, d)$ and call it the {\it symmetric lattice homology} associated with $(\mathcal{D},d)$. 

It supports an additional and  natural $\Z_2$-action induced by the automorphisms of the $S_{n, \mathcal{D}}$ spaces given by the involution 
$\ell\leftrightarrow d-\ell$. 
\end{define}

The Independence Theorem \ref{th:IndepMod} has the following counterpart in this setting:

\begin{theorem}[Symmetric version of the Independence Theorem]\label{th:Indep} Let $k$ be any field and $\mathcal{O}$ a (Noetherian) $k$-al\-geb\-ra. Let $\mathcal{D}=\{\frv_1, \ldots, \frv_r\}$ and $\mathcal{D}'=\{\frv'_1, \ldots, \frv'_{r'}\}$ be two collections of discrete valuations on $\mathcal{O}$ and $d \in (\mathbb{Z}_{\geq 0})^r, \ d'\in (\mathbb{Z}_{\geq 0})^{r'}$ two lattice points.  
        Suppose that the following are true:
        \begin{itemize}
            \item $\mathcal{F}_{\mathcal{D}}(d) = \mathcal{F}_{\mathcal{D}'}(d')$ and
            \item both pairs $(\mathfrak{h}_{\mathcal{D}}, \mathfrak{h}_{\mathcal{D},\, d}^{sym})$ and $(\mathfrak{h}_{\mathcal{D}'}, \mathfrak{h}_{\mathcal{D}',\, d'}^{sym})$ satisfy the Combinatorial Duality Property.
        \end{itemize}
        Then the spaces $S_{n, \mathcal{D}}$ and $S_{n, \mathcal{D}'}$ 
        associated with the corresponding lattices and weight functions 
        are homotopy equivalent for every $n \in \mathbb{Z}$.
        Even more, the homotopy equivalences respect the $\mathbb{Z}_2$-action and the inclusions
        $S_{n, \mathcal{D}}\hookrightarrow S_{n+1, \mathcal{D}}$ and $S_{n, \mathcal{D}'}\hookrightarrow S_{n+1, \mathcal{D}'}$, hence
        \begin{equation*}
        \mathbb{SH}_*(\mathcal{O}, \mathcal{D}, d) \cong  \mathbb{SH}_*(\mathcal{O}, \mathcal{D}', d') \text{ as bigraded } \mathbb{Z}[U]\text{-modules with } \mathbb{Z}_2\text{-action.}
        \end{equation*}
\end{theorem}

 That is, the lattice homology and the homotopy types of the cubical complexes $S_{n, \mathcal{D}}$ are independent of the pairs $(\mathcal{D}, d)$ of collections of discrete valuations and lattice points, they are invariants merely of  
 the finite codimensional ideal $\mathcal{F}_{\mathcal{D}}(d) \triangleleft \mathcal{O}$. In fact, we have the following identification:

 \begin{remark}\label{rem:SH=HforCDP}
    By Corollary \ref{cor:rectCDP->fullCDP}, Theorem \ref{th:properties} \textit{(a)} and Proposition \ref{prop:symforideals}, if  $(\mathfrak{h}_{\mathcal{D}}, \mathfrak{h}_{\mathcal{D}, \,d}^{sym})$ satisfies the Combinatorial Duality Property on the rectangle $R(0, d)$, then 
    \begin{center}$\mathbb{SH}_*(\mathcal{O}, \mathcal{D}, d) = \mathbb{H}_*(\mathcal{F}_{\mathcal{D}}(d) \hookrightarrow_{\mathcal{O}} \mathcal{O})$. 
    \end{center}
\end{remark}

The proof of the symmetric version of the Independence Theorem \ref{th:Indep} will be given in section \ref{s:proof}.

\subsection{The symmetric lattice homology of integrally closed ideals}\label{ss:Indep}\,

 Motivated 
by the Independence Theorem \ref{th:Indep}, as a second starting point, we might fix the Noetherian ring  $\calO$ with a $k$-algebra structure 
and an integrally closed ideal
$\mathcal{I}\triangleleft \calO$ with finite $k$-codimension. 

\begin{define}\label{def:REALforideals} If the ideal $\mathcal{I}\triangleleft \calO$ can be written as $\mathcal{F}_{\mathcal{D}}(d)$
for some collection of discrete valuations $\mathcal{D}=\{\frv_1,\ldots,\frv_r\}$ and lattice point
$d\in (\Z_{\geq 0})^r$, then we say that 
the pair 
      $(\mathcal{D},d)$ is  a \textit{`realization of $\mathcal{I}$'}. By Corollary \ref{cor:idealREAL}, this is equivalent to the collection
      \begin{center}
          $\mathcal{D}_{ext}:=\{ (\mathfrak{v}_1, \mathfrak{v}_1-d_1), \ldots, (\mathfrak{v}_r, \mathfrak{v}_r - d_r)\}$
      \end{center} 
      of extended discrete valuations being a realization of $\mathcal{I}$ as an $\mathcal{O}$-submodule of the rank one, free module $\mathcal{O}$ in the sense of Definition \ref{def:REAL}. If the pair $(\mathfrak{h}_{\mathcal{D}}, \mathfrak{h}_{\mathcal{D},\, d}^{sym})$ of height functions satisfies the Combinatorial Duality Property in the rectangle $R(0, d)$, then we call $(\mathcal{D}, d)$ a \textit{`CDP realization'}.  
\end{define}

Here several comments are in order. We stress that the realization of some ideal by a pair
 $(\mathcal{D}, d)$ is highly non-unique. For example, assume that $\mathcal{I}$ is realized by 
some   $(\mathcal{D}, d)$. Then even with the very same $\mathcal{D}$ it can happen that 
there are several distinct lattice points $d'$ such that for all of them one has 
$\mathcal{I}=\mathcal{F}_{\mathcal{D}}(d')$ (e.g., in the case of the analytic lattice homology of an irreducible Gorenstein curve singularity we have  $\mathcal{F}(\mathbf{c})=\mathcal{F}(\mathbf{c}-1)=\mathcal{C} \triangleleft \mathcal{O}$ --- for notations see section \ref{s:deccurves}; or, see the case of the analytic lattice homology of weighted homogeneous surface singularities with rational homology sphere link in Remark \ref{rem:whomrathomslink}),
hence, a priori,  all the symmetrization procedure might depend on the  
different choices. Moreover, it is also possible
that  $\mathcal{I}$ can be realized by completely different sets of valuations (and lattice points),
see, e.g., the situation discussed in subsection \ref{ss:NNan&comb}. The previous symmetric version of the Independence Theorem \ref{th:Indep} claims exactly the fact that, using the above construction, all the different CDP realizations yield the very same lattice homology module. 

On the other hand, once a realization $(\mathcal{D}, d) $ of $\mathcal{I}$
is found, we can naturally construct CDP realizations from it, similarly to the case of modules. Indeed, the following specialization of Lemma \ref{lem:duplatrukk} is true: 

\begin{lemma}\label{lem:duplatrukkideal}
Let $\cO$ and $\mathcal{I}$ be as above.  Assume that the pair $(\mathcal{D}, d)$ is a realization of $\mathcal{I}$, where
 $\mathcal{D}=\{(\frv_v, \frv^M_v)\}_{v=1}^r$  is a collection of discrete valuations on $\mathcal{O}$ 
 and $d \in (\mathbb{Z}_{\geq 0})^r$ is a lattice point.
 Then the following facts hold:

(1)  the multiset
 $\mathcal{D}^\sharp:=\{ \frv_1,  \frv_2, 
\ldots , \frv_r, 
\frv_{r+1}=\frv_1,
,\ldots , \frv_{2r}=\frv_r\}$ of $2r$ discrete valuations and the lattice point $d^\sharp:=(d, d) \in (\mathbb{Z}_{\geq 0})^{2r}$ give a CDP realization of $N$;

(2)
The collection $\mathcal{D}^\natural:=\{ 2\frv_v\}_{v=1}^r$  
of `doubled' discrete valuations and the `doubled' lattice point $d^\natural=2d \in (\mathbb{Z}_{\geq 0})^r$
    give a CDP realization of $N$.
\end{lemma}

As a consequence of this discussion, 
 we can introduce the well-defined symmetric lattice homology of any finite codimensional integrally closed ideal.

\begin{define}\label{def:LC}
    Let $k$ be a field, $\mathcal{O}$ a Noetherian $k$-algebra and $\mathcal{I} \triangleleft\mathcal{O}$ a finite $k$-codimen\-sional integrally closed ideal. By Corollary \ref{cor:intersection}, it admits a realization, even more, by Lemma \ref{lem:duplatrukkideal} it admits a CDP realization $(\mathcal{D}, d)$. Then we denote the lattice homology module $\mathbb{SH}_*(\mathcal{O}, \mathcal{D}, d)$ associated with this setup by $\mathbb{SH}_*(\mathcal{I} \triangleleft \mathcal{O})$, and call it \textit{`the symmetric lattice homology of the (integrally closed) ideal $\mathcal{I}$'}. According to the Independence Theorem \ref{th:Indep}, it is independent of the CDP realization chosen.
\end{define}

Clearly, by Remark \ref{rem:SH=HforCDP}, the symmetric lattice homology $\mathbb{SH}_*(\mathcal{I} \triangleleft \mathcal{O})$ agrees with $\mathbb{H}_*(\mathcal{I} \hookrightarrow_{\mathcal{O}} \mathcal{O})$.
Even more, as an analogue of Theorems \ref{th:properties} and \ref{th:upperbound} and Proposition \ref{prop:homdegupperbound}, the symmetric lattice homology of integrally closed ideals satisfy the following:

\begin{theorem}\label{th:propertiesideal}
    Let $\mathcal{O}$ be a Noetherian $k$-algebra, $\mathcal{I} \triangleleft \mathcal{O}$ a finite $k$-codimensional, integrally closed ideal. \\
\noindent (a) The Euler characteristic is well-defined and gives $eu(\mathbb{SH}_*(\mathcal{I} \triangleleft \mathcal{O}))=\dim_k (\mathcal{O}/\mathcal{I})$. \\
\noindent (b) The $\mathbb{Z}[U]$-module 
 $\mathbb{SH}_*(\mathcal{I} \triangleleft \mathcal{O})$ ---
and the homotopy types of the $S_n$ spaces ---  only depend on the Artin algebra $\mathcal{O}/\mathcal{I}$. \\
\noindent (c) ${\rm homdim}(\mathcal{I} \triangleleft \mathcal{O})=\begin{cases}
    \max \{ q: \mathbb{SH}_{q}(\mathcal{I} \triangleleft \mathcal{O}) \neq 0\} & \text{ if }\mathcal{I} \neq \mathcal{O};\\
    -1 & \text{ if }\mathcal{I} = \mathcal{O}.
\end{cases}$
\\

\noindent \hspace{31mm} $
     \leq \min \{ |\mathcal{D}|\,:\, (\mathcal{D}, d) \text{ realization of } \mathcal{I}\} -1$.\\
\noindent (d) 
    If $k$ is an algebraically closed field, $(\mathcal{O}, \mathfrak{m})$ a local Noetherian $k$-algebra and $\mathcal{I}\triangleleft \mathcal{O}$ a finite codimensional integrally closed ideal, then for every $n > 0$ the $S_n$ space is contractible, hence the weight-grading of the \emph{reduced} symmetric lattice homology  is nonnegative,
    \begin{center}i.e.,
    $\mathbb{SH}_{{\rm red},q, -2n}(\mathcal{I} \triangleleft \mathcal{O}) = 0$
    for all $q\geq 0$ and $n >0$. 
    \end{center}
Moreover, if $\mathcal{I} \neq \mathcal{O}$, then $S_{0}$ is not connected.
\end{theorem}

\begin{remark}\label{rem:d_I}
{\bf The dualizing element and valuative ideal correspondence.} 
Given a $k$-algebra $\calO$ as above and  a collection $\mathcal{D}=\{\frv_1, \ldots, \frv_r\}$ of discrete valuations on $\cO$,
we can realize a correspondence between finite codimensional ideals of $\cO$ and
lattice points $\ell\in\mathbb{Z}^r$ as follows. Firstly, for any $\ell$ we can consider the ideal 
$\mathcal{F}_{\mathcal{D}}(\ell)$ as in (\ref{eq:F_D(l)}). Conversely, 
for any  finite codimensional ideal $\mathcal{I} \triangleleft \mathcal{O}$ (such that 
 none of the valuations $\frv_1, \ldots, \frv_r$ is identically $\infty$ on the ideal 
 $\mathcal{I}$), 
we can define the \emph{ dualizing  element $d_{\mathcal{I}}=\sum_{v\in \mathcal{V}}d_{\mathcal{I},v} e_v \in (\Z_{\geq 0 })^r$ of $\mathcal{I}$} by
\begin{equation*}
    d_{\mathcal{I},v}:=\mathfrak{v}_v(\mathcal{I})= \min \{ \frv_v(f) \,:\, f \in \cI \}= \max \{ k \in \N \,:\, \cI \subset \cF_{\frv_v}(k) \}<\infty \hbox{ for all }v \in \mathcal{V}.
\end{equation*}
Then the following facts hold.

\vspace{2mm}

(a)
If $\mathcal{I}\subset \mathcal{F}_{\mathcal{D}}(\ell)$ for some lattice point $\ell\in \mathbb{Z}^{|\mathcal{D}|}$, then
$$\mathcal{I}\subset \overline{\mathcal{I}}\subset \mathcal{F}_{\mathcal{D}}(d_{\mathcal{I}})\subset
\mathcal{F}_{\mathcal{D}}(\ell),$$
where $\overline{\mathcal{I}}$ is the
integral closure of $\mathcal{I}$.

Furthermore, we have $d_{\mathcal{L}}= d_{\mathcal{I}}$
for any ideal $\mathcal{L}$ with $\mathcal{I}\subset \mathcal{L}\subset \overline{\mathcal{I}}$.
Therefore, for any such ideal (or, in fact, for any reduction $\mathcal{L}$ --- in the sense of \cite[Subsection 1.2]{SH} --- of the finitely generated ideal $\overline{\mathcal{I}}$), the lattice homology $\mathbb{SH}_\ast(\cO,\mathcal{D}, d_{\mathcal{L}})$ is tautologically the same. 

Note also that equality between  $\overline{\mathcal{I}}$  and 
$\mathcal{F}_{\mathcal{D}}(d_{\mathcal{I}}) $  cannot be expected in full generality. Indeed, in order
to have such, the set $\mathcal{D}$  should be sufficiently `rich',
it must contain sufficiently many well-chosen valuations (e.g., the Rees valuations associated with $\overline{\mathcal{I}}$). 

(b) From the other viewpoint, if we write 
$\mathcal{I}(\ell)$ for  $\mathcal{F}_{\mathcal{D}}(\ell)$ (for some lattice point $\ell\in (\Z_{\geq 0})^r$), then 
\begin{center}
$d_{\mathcal{I}(\ell)}\geq \ell$ 
and $\mathcal{F}_{\mathcal{D}}( d_{\mathcal{I}(\ell)})=\mathcal{I}(\ell)$.
\end{center}
\end{remark}

\begin{remark}
    Let the integrally closed ideal $\mathcal{I} \triangleleft \mathcal{O}$  be realized by a collection $\mathcal{D}=\{\mathfrak{v}_1, \ldots, \mathfrak{v}_r\}$ of discrete valuations, i.e., $\mathcal{I}=\mathcal{F}_{\mathcal{D}}(d_{\mathcal{I}})$. Also, let us consider the collection $$\mathcal{D}^{\rm ext}=\{(\mathfrak{v}_1, \mathfrak{v}_1-d_{\mathcal{I}, 1}), \ldots, (\mathfrak{v}_r, \mathfrak{v}_r-d_{\mathcal{I}, r})\}$$ of extensions as in Proposition \ref{prop:symforideals}, such that $\mathcal{F}_{\mathcal{D}^{\rm ext}}^M(0) = \mathcal{I}$. In this setting we have $d_{\mathcal{D}^{\rm ext}}=d_{\mathcal{I}}$.
\end{remark}

\begin{remark}\label{ex:Id}
Fix any collection of discrete valuations $\mathcal{D}$ of $\cO$. 
Then for any $\ell\in
(\Z_{\geq 0})^r$ we can consider the ideal $\mathcal{I}=\mathcal{I}(\ell):= \mathcal{F}_{\mathcal{D}}(\ell)\triangleleft \cO$ as in the previous paragraph. 
The symmetric lattice homology module 
${\mathbb S}\bH_*(\mathcal{I}(\ell)\triangleleft \cO)={\mathbb S}\bH_*(\cO,\mathcal{D}^\natural, 2\ell)$
is a well-defined invariant of $\mathcal{I}(\ell)$, hence of $\ell$. It is the categorification of the Hilbert function $\frh$ associated with $\mathcal {D}$. Indeed, by part  \textit{(a)} of  Theorem \ref{th:propertiesideal},
$eu({\mathbb S}\bH_*(\mathcal{I}(\ell)\triangleleft \cO))=
\dim_k\, (\cO/\mathcal{I}(\ell))=\frh_{\mathcal{D}}(\ell)$. It depends essentially on the choice of $\ell$. 

Here the following  warning is in order. Although for every lattice point $\ell \in (\mathbb{Z}_{\geq 0})^r$ the ideal $\mathcal{I}(\ell)$ can be realized by the same collection $\mathcal{D}$ of discrete valuations, the fact whether the pair $(\mathfrak{h}_{\mathcal{D}}, \mathfrak{h}_{\mathcal{D}, \, \ell}^{sym})$ satisfies the Combinatorial Duality Property might depend on the lattice point $\ell$. Therefore, the need to use the doubling construction of Lemma \ref{lem:duplatrukkideal} to compute $\mathbb{SH}_*(\mathcal{I}(\ell) \triangleleft \mathcal{O})$ might change from $\ell$ to $\ell$. 
\end{remark}

\subsection{Some examples of valuations and integrally closed ideals} \label{ss:NEWEX}\

In this subsection we list some examples of discrete valuations and integrally closed ideals, which will be important in the following discussions. We also compute some (symmetric) lattice homology modules.

\bekezdes \textbf{Valuations coming from gradings.} If $\mathcal{O}$ is an $\mathbb{N}$-graded $k$-algebra, i.e., $\mathcal{O}=\bigoplus_{n \in \mathbb{N}}\mathcal{O}_{(n)}$ with $\mathcal{O}_{(n)}$ denoting the homogeneous part of degree $n$, then the usual order function
\begin{equation}
    {\rm ord}:\mathcal{O} \rightarrow \overline{\mathbb{N}}, \ f=\sum_{n\in \mathbb{N}} f_{(n)} \ (\text{where } f_{(n)}\in \mathcal{O}_{(n)}) \mapsto \min\{ n \in \overline{\mathbb{N}} \,:\, f_{(n)}\neq 0\}
\end{equation}
gives a discrete valuation. The standard examples of such $\mathbb{N}$-graded $k$-algebras are the polynomial rings (or formal power series rings in the complete case and convergent power series rings in the complex analytic case).  

\begin{nota}
    Throughout this note we will denote by $\mathcal{O}_m$ the polynomial ring $k[x_1, \ldots, x_m]$ (or $k[[x_1, \ldots, x_m]]$ in the complete case, or $\mathbb{C}\{x_1, \ldots, x_m\}$ in the complex analytic setting --- this will always be clear from the context). Moreover, we will denote the maximal ideal $(x_1, \ldots, x_m)$ by $\mathfrak{m}_m$.
\end{nota}

The usual direct sum decomposition (according to the degree) of $\mathcal{O}_m$ yields the so-called \textit{maximal ideal valuation}. The ideals realized by this are just the powers of the maximal ideal.

\begin{example}\label{ex:Ct}
Set $\mathcal{I}=\mathfrak{m}_1^d=(x_1^d)\triangleleft \cO_1=k[x_1]$ (where $d>0$ is an integer).

Here $\mathcal{D}$ consists of a single discrete valuation, the maximal ideal valuation,
and $d$ serves also as the dualizing element of $\mathcal{I}$ (cf. \ref{rem:d_I}). Then $\frh(\ell)=\ell$ for any $0\leq \ell\leq d$, so the pair $(\mathfrak{h},\mathfrak{h}^{sym}_d)$ does not satisfy the CDP in $R(0, d)$, hence we need to consider the `doubled' version (cf. Lemma \ref{lem:duplatrukkideal}).
Now the values of $\frh^\natural(\ell) $ for $0\leq \ell\leq 2d$ are $0, 1,1,2,2,\ldots ,d,d$
and the values of $w_0$ are $0,1,0,\ldots , 0,1,0$. Therefore 
${\mathbb S}\bH_0=\calt^-_0\oplus (\calt_0(1))^{\oplus d}$. The graded root is 

\begin{picture}(200,50)(0,50)
\linethickness{.5pt}

\put(100,80){\circle*{3}}
\put(110,80){\circle*{3}}
\put(150,80){\circle*{3}}
\put(160,80){\circle*{3}}
\put(130,70){\circle*{3}}
\put(130,60){\circle*{3}}
\put(130,77){\makebox(0,0){$\dots$}}
\put(130,50){\makebox(0,0){$\vdots$}}
\qbezier[20](80,80)(130,80)(180,80)
\put(230,80){\makebox(0,0){$({\mathbb S}\bH_0)_0=\Z^{d+1}$}}

\put(100,80){\line(3,-1){30}}
\put(110,80){\line(2,-1){20}}

\put(150,80){\line(-2,-1){20}}
\put(160,80){\line(-3,-1){30}}
\put(130,70){\line(0,-1){15}}

\put(230,55){\makebox(0,0){$eu=d$}}
\end{picture}
\end{example}

\begin{example}\label{ex:CT}
Consider the ideal $\mathcal{I}=\mathfrak{m}_m^d\triangleleft \cO_m=k[x_1,\ldots, x_m]$, for some $d>0$.
 
Then the $\frh$-values corresponding to the maximal ideal valuation are: $\binom{m-1}{m}, \binom{m}{m}, \ldots, \binom{m-1+d}{m}$, while the 
$\frh^\natural$-values satisfy $\frh^\natural(2\ell-1)=\frh^\natural(2\ell)=\frh(\ell)$.
Note also that $\frh(0)=0$ and $\frh(\ell+1)-\frh(\ell)=\binom{m+\ell-1}{m-1}$.
Therefore, $w_0(0)=0$ and the $w_0$-values for $0\leq \ell\leq 2d$ can be obtained inductively as follows:  
\begin{center}
    $w_0(2\ell+1)-w_0(2\ell)=\binom{m+\ell-1}{m-1}$ and 
 $w_0(2\ell)-w_0(2\ell-1)=-\binom{m+d-\ell-1}{m-1}$.
\end{center}

For example, if $m=2$ and $d=5$ the Hilbert function for $0\leq \ell\leq 5$ takes the values $0,1,3, 6, 10, 15$;\ thus the values of 
$\frh^\natural(\ell)  $  $(0\leq \ell\leq 10)$
are $0,1,1,3,3,\dots, 15,15$ and of the weight function $w_0$ are  
\begin{center}
    $0, 1, -4, -2, -6, -3, -6, -2, -4, 1,0$.
\end{center} 
Hence, the graded root is 

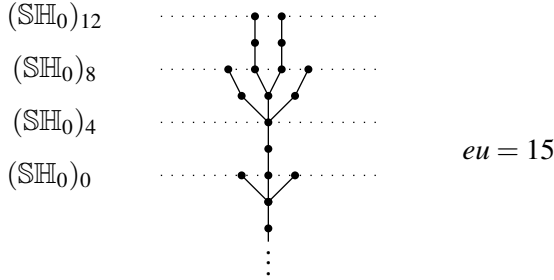
\begin{picture}(200,110)(0,10)
\linethickness{.5pt}

\put(125,100){\circle*{3}}
\put(135,100){\circle*{3}}
\put(125,90){\circle*{3}}
\put(135,90){\circle*{3}}
\put(115,80){\circle*{3}}
\put(125,80){\circle*{3}}
\put(135,80){\circle*{3}}
\put(145,80){\circle*{3}}
\put(120,70){\circle*{3}}
\put(130,70){\circle*{3}}
\put(140,70){\circle*{3}}
\put(130,60){\circle*{3}}
\put(130,50){\circle*{3}}
\put(120,40){\circle*{3}}
\put(130,40){\circle*{3}}
\put(140,40){\circle*{3}}
\put(130,30){\circle*{3}}
\put(130,30){\circle*{3}}
\put(130,20){\circle*{3}}

\put(130,10){\makebox(0,0){$\vdots$}}
\qbezier[20](90,100)(130,100)(170,100)
\qbezier[20](90,80)(130,80)(170,80)
\qbezier[20](90,60)(130,60)(170,60)
\qbezier[20](90,40)(130,40)(170,40)

\put(50,100){\makebox(0,0)[0]{$({\mathbb S}\bH_0)_{12}$}}
\put(50,80){\makebox(0,0)[0]{$({\mathbb S}\bH_0)_{8}$}}
\put(50,60){\makebox(0,0)[0]{$({\mathbb S}\bH_0)_{4}$}}
\put(48,40){\makebox(0,0)[0]{$({\mathbb S}\bH_0)_{0}$}}

\put(125,100){\line(0,-1){20}}
\put(135,100){\line(0,-1){20}}
\put(125,80){\line(1,-2){5}}
\put(135,80){\line(-1,-2){5}}
\put(115,80){\line(1,-2){5}}
\put(145,80){\line(-1,-2){5}}
\put(120,70){\line(1,-1){10}}
\put(140,70){\line(-1,-1){10}}

\put(130,70){\line(0,-1){55}}
\put(120,40){\line(1,-1){10}}
\put(140,40){\line(-1,-1){10}}

\put(220,50){\makebox(0,0){$eu=15$}}
\end{picture}
\end{example}

\bekezdes\textbf{Monomial valuations.} \label{par:monomialval} \ Let us assign some positive integer weights $\{{\rm w}_1, \ldots, {\rm w}_m\}$ to the variables $\{x_1, \ldots, x_m\}$ and consider the corresponding direct sum decomposition of $\mathcal{O}_m$ (i.e., the homogeneous monomials $x_1^{p_1} \cdot \ldots \cdot x_m^{p_m}$ of degree $d$ are those for which $p_1 {\rm w}_1 + \ldots p_m {\rm w}_m=d$). Then the associated valuation is called a \textit{monomial valuation}, and it can, clearly, only realize monomial ideals. On the other hand, these valuations are enough to characterize the integrally closed monomial ideals (using their Newton diagrams). We will discuss these valuations and ideals in section \ref{s:ND}. 

\bekezdes\textbf{Divisorial valuations.}
   Let $(X,o)$ be a normal complex analytic singularity of dimension $\geq 2$.
   Fix a resolution $\phi:\widetilde{X}\to X$ and let $E$ be an irreducible exceptional divisor in $\widetilde{X}$. Then $E$ defines a valuation on the local ring
   $\mathcal{O}_{X,o}$ by the assignment $f\mapsto \frv_E(f):=
   \{ \text{vanishing order of $ f\circ \phi$ along $E$}\}$. 
   It is called the {\it divisorial valuation} associated with $E$. 
If $\phi(E)=o$  then 
   $\frv_{E}$ is called $\mathfrak{m}_{X,o}$-centered.

   In particular, if $\mathcal{D}=\{E_v\}_v $ is a set of irreducible exceptional divisors, and $\ell:=\sum_v \ell_v E_v$ is a divisor supported on them ($\ell_v\in\bZ$ for all $v$), 
   then $\mathcal{F}_{\mathcal{D}}(\ell)=\{f\in\mathcal{O}_{X,o}\,:\, \frv_{E_v}(f)\geq \ell_v \ \text{for all $v$}\}= \phi_*\mathcal{O}_{\widetilde{X}}(-\ell)_o$
   is an integrally closed ideal in $\mathcal{O}_{X,o}$ and 
   ${\rm codim}_{\bC}(\mathcal{F}_{\mathcal{D}}(\ell)\hookrightarrow \mathcal{O}_{X,o})< \infty$. For more details see section \ref{s:dnagy}, especially Remark \ref{rem:filtrationreformulation}.

   If $\dim (X,o)=2$ then, by a theorem of Lipman \cite{Lipman}, any integrally closed ideal of $\mathcal{O}_{X,o}$ can be realized in this way  for a certain choice of the resolution $\phi$ and the divisor $\ell$.

   \bekezdes\textbf{Valuations associated with reduced curves and embedded Weil divisors.}\label{bek:4.10.2} \
   If $(C,o)$ is an irreducible complex analytic curve singularity, then its normalization 
   $n_C:{\rm Spec}(\mathbb{C}\{t\})=\overline{(C,o)}\to (C,o)$ induces a valuation of $\mathcal{O}_{C,o}$ by $f\mapsto \frv_{n_C}(f):={\rm ord} _t(f\circ n_C)$. This is the main ingredient for the definition of the analytic lattice homology of curve singularities.
   For more see section \ref{s:deccurves}.

   Now assume that  $(C,o)\subset (X,o)$ is an irreducible Weil divisor on the normal surface singularity $(X,o)$, and let $\phi:\widetilde{X}\to X$ be an embedded resolution of the pair $(C,o)\subset (X,o)$.  Let $(\widetilde{C},\widetilde{o})$ denote the strict transform of $(C,o)$  and let $E_{\phi,C}$ be the irreducible exceptional curve which intersects $\widetilde{C}$. To such a pair we can associate two natural discrete valuations of 
   $\mathcal{O}_{X,o}$. One of them is the composition  $\frv_{n_C,C}:\mathcal{O}_{X,o}\stackrel{q}{\longrightarrow} \mathcal{O}_{C,o}\stackrel{\frv_{n_c}}{\longrightarrow} \bZ$ associated with the normalization of $(C, o)$;
   while the other one, $\frv _{E_{\phi,C}}$, is the divisorial valuation associated with $E_{\phi,C}$. One has  $\frv_{n_C,C}\geq \frv_{E_{\phi,C}}$, but in general they are not equal (for a comparison of their  ideal filtrations see, e.g., \cite{NV}). A similar phenomenon will be discussed in subsection \ref{ss:NNplanecurves}, see especially Lemma \ref{lem:van>vcomb}.
   
   \noindent[The valuation 
   $\frv_{n_C,C}$ of $\mathcal{O}_{X,o}$ is also called the {\it arc-valuation} associated with $(C,o)$. For some related comments see subsection \ref{ss:arcs}.]

    \bekezdes\textbf{Valuations associated with decorated irreducible Weil divisors.} \label{bek:4.10.3} \ 
Consider  a pair of set-germs $(C,o)\subset (X,o)$ as in the previous paragraph \ref{bek:4.10.2} 
 and let $\phi:\widetilde{X}\to X$ be the minimal good embedded resolution of the pair. Let $P$ denote the intersection point of the strict transform $\widetilde{C}$ with its 
 supporting exceptional curve. Then for any integer $n\in\bZ_{\geq 0}$ --- the \textit{decoration} of $(C,o)$ --- we blow up the 
 intersection point of the strict transform of $C$ with its supporting exceptional curve repeatedly  $n$ times 
 (that is, we blow up $n$ times the corresponding infinitely near point $P$).
 Let $E(n,C)$ be the exceptional curve in the final resolution supporting the last strict transform. Then 
 the divisorial valuation $\frv_{E(n,C)}:\mathcal{O}_{X,o}\to \bZ$ is called the 
 {\it valuation associated with the decorated curve $\{(C,o),n\}$}. 
 [Ideals of type $\mathcal{F}_{\frv_{E(n,C)}}(\ell)$ associated with such a valuation 
 were considered, e.g., in the deformation theory of sandwiched surface singularities, cf. 
 \cite{dJvS}. Note, however, that we use a different convention for the decorations.] One sees that $\mathfrak{v}_{E(n+1, C)} \geq \mathfrak{v}_{E(n, C)}$ for every $n\in \mathbb{Z}_{\geq 0}$. Moreover, for any fixed function $f \in \mathcal{O}_{X, o}$ the limit $\lim_{n\to \infty}\mathfrak{v}_{E(n, C)}(f)=\mathfrak{v}_{n_C, C}(f)$. For the limiting properties of the corresponding filtration see \cite{NV}. 
 
   For example, if $(C,o)=(\{x^2+y^3=0\}, 0)\subset (X,o)=(\bC^2,0)$,
   and $n=2$, then $\mathcal{F}_{\frv_{E(2,C)}}(8)$ is the ideal 
   $(y^4, y^3x, y^3+x^2)$, while for $n=4$  we have 
   $\mathcal{F}_{\frv_{E(4,C)}}(10)=  
   (y^5, y^4x, y^3x^2, y^3+x^2)$. 

By Corollary \ref{cor:intersection}, ideals of type  $\mathcal{F}_{\frv_{E(n,C)}}(\ell)$ are integrally closed, however, if we just look at their generators, e.g., for 
$(y^4, y^3x, y^3+x^2)$ (without knowing its realization  as above),
this property is not immediate  at all (compare with \cite[Example 1.3.3]{SH}). In this way one can construct interesting nontrivial (e.g., non-monomial) integrally closed ideals. 

In fact, we can also present another realization of some of these ideals. Consider again the valuation  $\frv_{E(n,c)}$ as above and let 
$m=\min\{ \frv_{E(n,C)}(f\circ \phi)\,:\, f\in\mathcal{I}_{C,o}\}$, where $\mathcal{I}_{C, o}$ denotes the local ideal of the curve germ $(C, o) \subset (X, o)$.
Now clearly, if $\ell\leq m$, then the ideal $\mathcal{I}_{C,o}$ is contained in 
$\mathcal{F}_{\frv_{E(n,C)}}(\ell)$. Hence such ideals $\{\mathcal{F}_{\frv_{E(n,C)}}(\ell)\}_{\ell\leq m}$ are pullbacks of 
ideals of $\mathcal{O}_{C,o}$ via the quotient map $q:\mathcal{O}_{X,o}\to \mathcal{O}_{C,o}$.
For examples above, with $(C,o)=\{x^2+y^3=0\}\subset (X,o)=(\bC^2,0)$, we have
$\mathcal{F}_{\frv_{E(2,C)}}(8)=q^{-1}(\mathcal{F}_{\frv_{\nu}}(8))$
 and 
$\mathcal{F}_{\frv_{E(4,C)}}(10)=q^{-1}(\mathcal{F}_{\frv_{\nu}}(10))$.
In particular, we see that the very same  integrally closed ideal can be realized  by (sometimes rather) different 
set of valuations (see also how the conductor ideal of a reduced curve singularity 
is realized in multiple different ways in Lemma \ref{lem:agreec}). 

\begin{example}\label{ex:uj1}
Consider the ideal $\mathcal{I}=\mathcal{F}_{E(4,C)}(10)=(y^5, y^4x, y^3x^2, y^3+x^2)$ of $\bC\{x,y\}$ constructed in paragraph \ref{bek:4.10.3}. 
It is realized by a single divisorial valuation $\frv_{E(4,C)}$ of a modification of $(\bC^2,0)$ and has codimension 9. Indeed, the monomials $y^5, y^4x, y^3x^2, y^2x^2, yx^3, x^4$ belong to $\mathcal{I}$, but one also has three other relations $y^3+x^2, \ y^4+yx^2, \ y^3x+x^3$  among the `smaller degree' monomials. Hence, the monomials $1, y, y^2, y^3, y^4, x, xy, xy^2, xy^3$ 
generate $\bC\{x,y\}/\mathcal{I}$, their $\frv_{E(4,C)}$-values 
are 0,2,4,6,8,3,5,7,9. Then the dualizing element is $d_{\mathcal{I}}=10$ and the Hilbert function $\mathfrak{h}(\ell)$-values for $0\leq \ell\leq 10$ are 
0,1,1,2,3,4,5,6,7,8,9. By doubling and symmetrizing this sequence we get the weights and the following graded root (or ${\mathbb S}\bH_0$) 
\begin{center}
\begin{picture}(200,60)(60,35)
\linethickness{.5pt}

\put(100,80){\circle*{3}}
\put(110,80){\circle*{3}}
\put(150,80){\circle*{3}}
\put(120,80){\circle*{3}}
\put(130,80){\circle*{3}}
\put(140,80){\circle*{3}}
\put(160,80){\circle*{3}}
\put(130,60){\circle*{3}}
\put(130,70){\circle*{3}}
\put(130,60){\circle*{3}}

\put(150,70){\circle*{3}}
\put(110,70){\circle*{3}}

\put(130,45){\makebox(0,0){$\vdots$}}
\qbezier[20](80,70)(130,70)(180,70)
\put(60,70){\makebox(0,0){$n=0$}}

\put(100,80){\line(3,-1){30}}
\put(110,80){\line(2,-1){20}}
\put(120,80){\line(1,-1){10}}
\put(110,70){\line(2,-1){20}}
\put(130,80){\line(0,-1){30}}

\put(150,80){\line(-2,-1){20}}
\put(160,80){\line(-3,-1){30}}
\put(140,80){\line(-1,-1){10}}
\put(150,70){\line(-2,-1){20}}

\put(220,55){\makebox(0,0){$eu=9$, \ $\mathbb{S}\bH_{\geq 1}=0$}}
\end{picture}
\end{center}

Recall that $\mathcal{F}_{E(4,C)}(10)=q^{-1}(\mathcal{F}_{\frv_\nu}(10))$, where $q:\bC\{x,y\}\to \mathcal{O}_{C,o}=
\bC\{x,y\}/(y^3+x^2)$ is the natural quotient map (cf. paragraph \ref{bek:4.10.3}). Hence, the above $\mathbb{S}\bH_*(\mathcal{I} \triangleleft \mathbb{C}\{x, y\})$ coincides with the symmetric lattice homology of the ideal $\mathcal{F}_{\frv_{n_C}}(10) \triangleleft \mathcal{O}_{C,o}$ as well. 

It is instructive to compare the above module with the symmetric lattice homology of the ideal $\mathcal{I}'\hookrightarrow \bC\{x,y\}$ given by $(y^5, y^4x, y^3x^2, y^2x^2, yx^3, x^4) $,  realized by the monomial valuation corresponding to the weights $\frv(x)={\rm w}_x=3$, $\frv(y)={\rm w}_y=2$:
$\mathcal{I}'=(\{x^ay^b\,|\, 3a+2b\geq 10\})$ (for the formal definition of monomial valuations see, e.g., Definition 6.1.4 in \cite{HS}). In this case the codimension is 12, 
 the 
Hilbert function $\mathfrak{h}(\ell)$-values for $0 \leq \ell \leq 10$ are 0,1,1,2,3,4,5,7,8,10,12,  and the graded root is 

\begin{center}
\begin{picture}(200,80)(60,40)
\linethickness{.5pt}

\put(130,60){\circle*{3}}

\put(150,70){\circle*{3}}
\put(110,70){\circle*{3}}
\put(130,70){\circle*{3}}

\put(130,80){\circle*{3}}
\put(130,90){\circle*{3}}

\put(110,100){\circle*{3}}
\put(120,100){\circle*{3}}
\put(130,100){\circle*{3}}
\put(140,100){\circle*{3}}
\put(150,100){\circle*{3}}

\put(130,110){\circle*{3}}
\put(120,110){\circle*{3}}
\put(140,110){\circle*{3}}

\put(130,45){\makebox(0,0){$\vdots$}}
\qbezier[20](80,70)(130,70)(180,70)
\put(60,70){\makebox(0,0){$n=0$}}
\qbezier[20](80,100)(130,100)(180,100)
\put(60,100){\makebox(0,0){$n=3$}}

\put(110,100){\line(2,-1){20}}
\put(120,100){\line(1,-1){10}}\put(120,110){\line(1,-1){10}}
\put(110,70){\line(2,-1){20}}
\put(130,110){\line(0,-1){60}}

\put(150,100){\line(-2,-1){20}}
\put(140,100){\line(-1,-1){10}}\put(140,110){\line(-1,-1){10}}
\put(150,70){\line(-2,-1){20}}

\put(220,45){\makebox(0,0){$eu=12$, \ $\mathbb{S}\bH_{\geq 1}=0$}}

\end{picture}
\end{center}
\end{example}

\begin{example}
Let us examine the situation of paragraph \ref{bek:4.10.3}, but at this time with
$(C,o)\subset (X,o)$ being
a non-Cartier irreducible divisor. In our case, let the ambient surface singularity $(X,o)$ be a cyclic quotient singularity  with resolution graph 
\begin{picture}(50,10)(0,0)
    \put(10,0){\line(1,0){40}}
\put(10,0){\circle*{3}}
\put(30,0){\circle*{3}}
\put(50,0){\circle*{3}}
\put(10,5){\makebox(0,0){\tiny{$-2$}}}
\put(30,5){\makebox(0,0){\tiny{$-3$}}}
\put(50,5){\makebox(0,0){\tiny{$-2$}}}
\end{picture} \hspace{3mm}.

In this discussions we will use the language of toric geometry and we also assume certain familiarity with Weil divisors associated with Newton diagrams in this setting (see, e.g.,  \cite{Hir1,NBook,NSigNN}).
In particular, we consider the (real) cone $\sigma$ in $\bR_{\geq 0}^2$ generated by the vectors 
$(1,0)$ and $(5,8)$, then $(X,o)$ is the germ at the singular point of the affine variety
${\rm Spec}(\bC[\sigma \cap \bN^2])$. This algebra $\bC[\sigma \cap \bN^2]$ is generated by 
the semigroup elements $x=(1,0)$, $y=(1,1)$, $z=(2,3)$ and $w=(5,8)$, and the semigroup relations provide the equations of $(X,o)$ in $(\bC^4,0)$ (with coordinates $(x,y,z,w)$):
$xz=y^3$, $xw=y^2z^2$ and $yw=z^3$). 

We fix a {\it  Weil divisor in $(X,o)$ identified by a  (nondegenerate) Newton diagram}:
the segment with endpoints $(4,0)$ and $(3,4)$ (or any shift of it by an integral vector).
In particular, the equations of $(C,o)$ in $(X,o)$ are $zy=x^4$, $z^2=y^2x^3$ and $w=zyx^3$.

The minimal good embedded resolution of the pair $(C,o)\subset (X,o)$ can be read in the dual 
cone $\widehat{\sigma}$ generated by $(0,1)$ and $(8,-5)$ as follows: it is the regular fan generated by the additional vectors 
$(1,1), \ (2,1),\ (3, 1), \ (4,1),\ (1,0), \ (2, -1), \ (5, -3)$. The dual graph 
(where the arrowhead represents the strict transform of $(C,o)$)
is

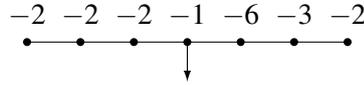
\begin{picture}(220,20)(-130,0)
    \put(10,0){\line(1,0){120}}
\put(10,0){\circle*{3}}
\put(30,0){\circle*{3}}
\put(50,0){\circle*{3}}
\put(70,0){\circle*{3}}
\put(90,0){\circle*{3}}
\put(110,0){\circle*{3}}\put(130,0){\circle*{3}}

\put(10,10){\makebox(0,0){$-2$}}
\put(30,10){\makebox(0,0){$-2$}}
\put(50,10){\makebox(0,0){$-2$}}
\put(70,10){\makebox(0,0){$-1$}}
\put(90,10){\makebox(0,0){$-6$}}
\put(110,10){\makebox(0,0){$-3$}}
\put(130,10){\makebox(0,0){$-2$}}
  \put(70,0){\vector(0,-1){15}}
\end{picture}

\vspace{5mm}

\noindent where the vertices correspond to the above list of vectors of $\widehat{\sigma}$. 

Let $F$ be the irreducible exceptional curve supporting the strict transform. 
Then the divisorial valuation associated with $F$ gives $\frv_F(x)=4$,  $\frv_F(y)=5$,
 $\frv_F(z)=11$,
 $\frv_F(w)=28$, hence the defining ideal $\mathcal{I}_{C,o}$ of $(C,o)$ is contained in the ideal $\mathcal{F}_{\frv_F}(16)$ of $\mathcal{O}_{X,o}$. (Indeed, 
 $\frv_F(zy-x^4)$ is exactly $16$, the other equations of $(C,o)$ have even higher values.)
 In fact, any ideal $\mathcal{F}_{\frv_F}(\ell)\subset \mathcal{O}_{X,o}$
  with $\ell\in \mathbb{Z}_{\geq 0}$ is exactly  the 
  monomial ideal 
  generated by $\{(a,b)\in \sigma\cap \bN^2\,|\, 4a+b\geq \ell\} $ (where $(4,1)$ is the normal vector of the Newton diagram of $(C,o)$).
  
  Let us now blow up the intersection point of $F$ with the strict transform once, in this way we get the resolution associated with the decorated curve $\{(C,o), 1\}$. Let the new irreducible exceptional curve be denoted by $E$. Then the divisorial valuation associated with $E$ in this new resolution  gives   $\frv_E(x)=4$,  $\frv_E(y)=5$,
 $\frv_E(z)=11$,
 $\frv_E(w)=28$, hence $\frv_E(x^4)=\frv_E(zy)=16$, however $\frv_E(zy-x^4)=17$. Thus, the ideal $\mathcal{I}:=\mathcal{F}_{\frv_E}(17)\hookrightarrow \mathcal{O}_{X,o}$ is a non-monomial integrally closed ideal. 
 It is generated by monomials associated with lattice points 
 $\{(a,b)\in \sigma\cap \bN^2\,|\, 4a+b\geq 17\} $ and the additional equation $zy-x^4$.
 (This ideal depends on $(C,o)$ because of the choice of the last blow up: the strict transform of $(C,o)$ is supported on $E$). The ideal $\mathcal{F}_{\frv_E}(17)$ contains $\mathcal{I}_{C,o}$.

The quotient Artin algebra $\bC[\sigma\cap \bN^2]/\mathcal{I}$ has basis
$1, x,y, 
x^2, xy, y^2, z, x^3, x^2y, xy^2, y^3, zy $ (corresponding to the lattice points 
$(0,0), (1,0), (1,1), (2,0), (2,1), (2,2), (2,3), (3,0), \ldots, (3, 4)$)
with $\mathfrak{v}_E$-values $0, 4, 5, 8, 9, 10, 11, 12, 13, 14,  15, 16$.
The Hilbert function  for $0\leq \ell\leq 17$  is 
\begin{center}
$0,1,1,1,1, 2, 3,3, 3, 4,5,6,7,8,9,10,11,12$,
\end{center}
thus, the graded roots is

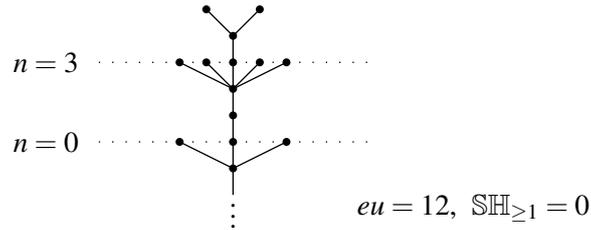
\begin{center}
\begin{picture}(200,80)(50,40)
\linethickness{.5pt}

\put(130,60){\circle*{3}}

\put(150,70){\circle*{3}}
\put(110,70){\circle*{3}}
\put(130,70){\circle*{3}}

\put(130,80){\circle*{3}}
\put(130,90){\circle*{3}}

\put(110,100){\circle*{3}}
\put(120,100){\circle*{3}}
\put(130,100){\circle*{3}}
\put(140,100){\circle*{3}}
\put(150,100){\circle*{3}}

\put(130,110){\circle*{3}}
\put(120,120){\circle*{3}}
\put(140,120){\circle*{3}}

\put(130,45){\makebox(0,0){$\vdots$}}
\qbezier[20](80,70)(130,70)(180,70)
\put(60,70){\makebox(0,0){$n=0$}}
\qbezier[20](80,100)(130,100)(180,100)
\put(60,100){\makebox(0,0){$n=3$}}

\put(110,100){\line(2,-1){20}}
\put(120,100){\line(1,-1){10}}\put(120,120){\line(1,-1){10}}
\put(110,70){\line(2,-1){20}}
\put(130,110){\line(0,-1){60}}

\put(150,100){\line(-2,-1){20}}
\put(140,100){\line(-1,-1){10}}\put(140,120){\line(-1,-1){10}}
\put(150,70){\line(-2,-1){20}}

\put(220,45){\makebox(0,0){$eu=12$, \ $\mathbb{S}\bH_{\geq 1}=0$}}
\end{picture}
\end{center}

Again, this can be compared with the symmetric lattice homology of 
$\mathcal{J}:=\mathcal{F}_{\frv_F}(17)\triangleleft \mathcal{O}_{X,o}$.
(This ideal does not contains $\mathcal{I}_{C,o}$ since $\frv_F(zy-x^4)=16$.)
The corresponding Artin algebra has dimension 13, it  is generated by 
classes of
$1, x,y, 
x^2, xy, y^2, z, x^3, x^2y, xy^2, y^3, zy, x^4 $ (corresponding to the lattice points 
$(0,0), (1,0), (1,1), (2,0), (2,1), (2,2), (2,3), (3,0), \ldots, (3, 4), (4,0)$)
with their $\mathfrak{v}_F$-values being $0, 4, 5, 8, 9, 10, 11, 12, 13, 14,  15, 16, 16$.
The Hilbert function  for $0\leq \ell\leq 17$  is 
\begin{center}
    $0,1,1,1,1, 2, 3,3, 3, 4,5,6,7,8,9,10,11,13$, 
\end{center}
whereas the graded root is

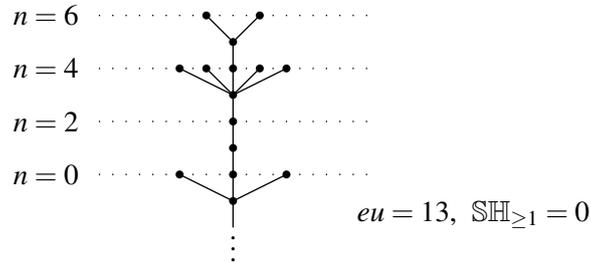
\begin{center}
\begin{picture}(200,105)(50,25)
\linethickness{.5pt}

\put(130,60){\circle*{3}}\put(130,50){\circle*{3}}

\put(150,60){\circle*{3}}
\put(110,60){\circle*{3}}
\put(130,70){\circle*{3}}

\put(130,80){\circle*{3}}
\put(130,90){\circle*{3}}

\put(110,100){\circle*{3}}
\put(120,100){\circle*{3}}
\put(130,100){\circle*{3}}
\put(140,100){\circle*{3}}
\put(150,100){\circle*{3}}

\put(130,110){\circle*{3}}
\put(120,120){\circle*{3}}
\put(140,120){\circle*{3}}

\put(130,35){\makebox(0,0){$\vdots$}}
\qbezier[20](80,60)(130,60)(180,60)
\put(60,60){\makebox(0,0){$n=0$}}
\qbezier[20](80,80)(130,80)(180,80)
\put(60,80){\makebox(0,0){$n=2$}}
\qbezier[20](80,100)(130,100)(180,100)
\put(60,100){\makebox(0,0){$n=4$}}
\qbezier[20](80,120)(130,120)(180,120)
\put(60,120){\makebox(0,0){$n=6$}}

\put(110,100){\line(2,-1){20}}
\put(120,100){\line(1,-1){10}}\put(120,120){\line(1,-1){10}}
\put(110,60){\line(2,-1){20}}
\put(130,110){\line(0,-1){70}}

\put(150,100){\line(-2,-1){20}}
\put(140,100){\line(-1,-1){10}}\put(140,120){\line(-1,-1){10}}
\put(150,60){\line(-2,-1){20}}

\put(220,45){\makebox(0,0){$eu=13$, \ $\mathbb{S}\bH_{\geq 1}=0$}}
\end{picture}
\end{center}

  Finally, we note that the semigroup of $(C,o)$ as an abstract curve 
  (that is, the semigroup of values of $\frv_{n_C}$, cf. subsection \ref{ss:excurves}) is $\mathcal{S}_{C,o}=\langle 4,5 , \geq 8\rangle\subset \bN$.
  Hence $(C,o)$  has  conductor 8,  delta invariant  5,  and 
  it is a non-Gorenstein  curve singularity with embedded dimension 3 (for terminology see section \ref{s:deccurves}). 
\end{example}

For examples with non-vanishing $\bH_1$ (i.e., of higher homological dimension), see subsection \ref{ss:excurves} in the setting of analytic lattice homology of complex analytic reduced curve singularities.

\section{Symmetric lattice homology of `integrally reduced' Artin  $k$-algebras}\label{ss:ARTIN}

In this section we will present some properties of those Artin algebras, which arise as quotients of Noetherian algebras by finite codimensional integrally closed ideals. The reason behind this inquiry is 
part \textit{(b)} of Theorem \ref{th:propertiesideal}, which states that the symmetric lattice homology $\mathbb{SH}_*(\mathcal{I} \triangleleft \mathcal{O})$ of an integ\-rally closed ideal $\mathcal{I}$ is in fact an invariant of the Artin algebra $\mathcal{O}/\mathcal{I}$. We will discuss some characterizations in different settings and provide a Künneth-type formula for computing their symmetric lattice homology with respect to their classical product decomposition.

\subsection{\emph{`Integrally reduced'} Noetherian $k$-algebras}\,

In this subsection we will examine the properties of $k$-algebras $A$ which can be realized as
$\cO/\mathcal{I}$, where $\cO$ is a
 Noetherian  ring with a $k$-algebra structure
 and $\mathcal{I}=\overline{\mathcal{I}}$ is an integrally closed ideal in $\cO$ of finite codimension. Recall also that, by part \textit{(2)} of Proposition \ref{prop:pullbackring}, the Artin $k$-algebra $\mathcal{O}/{\rm Ann}_{\mathcal{O}}(M/N)$ of Theorem \ref{th:properties} \textit{(c)} is also of this type, so, in fact, all the lattice homology theory of realizable submodules concerns the representation theory of such algebras.

 \begin{define}\label{def:intreduced}
     We call a Noetherian $k$-algebra $A$ \emph{`integrally reduced'} if there exists a presentation $A \cong \mathcal{O}/\mathcal{I}$ with $\mathcal{O}$ a Noetherian $k$-algebra and $\mathcal{I} \triangleleft \mathcal{O}$ an integrally closed ideal. 
 \end{define}

 In such a case, since $\cO$ is finitely generated, there exists a surjective morphism $\alpha:\cO_m\to \cO$ for some $m \in \mathbb{N}$,
 where $\cO_m$ denotes the polynomial ring with $k$-coefficients $k[x_1,\ldots, x_m]$ (in the complete\,/\,complex analytic case this could be $k[[x_1,\ldots, x_m]]$ or $\mathbb{C}\{x_1,\ldots, x_m\}$). 
 Since $\mathcal{J}:=\alpha ^{-1}(\mathcal {I})$ is integrally closed in $\cO_m$ as well, the Artin algebra
 $A$ can also be represented as $\cO_m/\mathcal{J}$, where $\overline {\mathcal{J}}=\mathcal{J}\triangleleft \cO_m$.

 On the other hand, clearly not any $k$-algebra is integrally reduced. E.g., in the complex analytic setting  we have the following (valuative) characterization.
 \begin{proposition}\label{prop:VAL}
 (see e.g., \cite[2.1 Théorème]{LT} or \cite{Teissier})\\  Assume $A$ is a complex analytic algebra. Then 
 
  (i)  
  $A$ is of integrally reduced if and only if the following fact holds:
\begin{equation}\label{eq:VAL}\mbox{
 if $a\in A$ and $\gamma(a)=0$
 for any homomorphism $\gamma:A\to \bC\{t\}/(t^m)$, then $a=0$.}
 \end{equation}
 
 (ii)
Moreover, if $A$ satisfies this valuative criterion  (\ref{eq:VAL}), 
and it is written as 
$\cO_m/{\mathcal{I}}$ for some $\mathcal{I}$, then $\mathcal{I}$ is automatically integrally closed.
\end{proposition}

Note that part \textit{(ii)}  does not imply that for any realization of $A$ as 
$\cO/\mathcal{I}$ (with $\cO$ not necessarily of type $\cO_m$)
the ideal $\mathcal{I}$ is integrally closed. E.g., even if $A$ satisfies  (\ref{eq:VAL}), 
its  zero ideal $\mathcal{I}=(0)$ is not necessarily  integrally closed (e.g., consider $A=\C\{t\}/(t^2)$).
Hence   property {\it (ii)} fails for the realization $A=A/(0)$.  [With the slightly more general and somewhat combinatorial definition of \textit{`quasi-valuations'} and \textit{`quasi-realizability'}, this problem can be remedied, see subsection \ref{ss:quasi}, especially Remark \ref{rem:robustness-q}.]

Nevertheless, this property \textit{(ii)} is true also in the general algebraic case. With the aim of proving this statement we introduce the notion of the \emph{`integral reduction'} of an algebra.

\bekezdes \textbf{Integral reduction algebra.}
Given {\it any} finitely generated $k$-algebra $A$, by choosing a generating set of $m$ elements, we can present it as $A = \mathcal{O}_m/\mathcal{I}$. By taking the integral closure $\overline{\mathcal{I}}$ of $\mathcal{I}$ in $\mathcal{O}$, we can naturally consider the {\it 
 `integral reduction algebra'} $\widehat{\mathcal{O}_m/\mathcal{I}}=\mathcal{O}_m/\overline{\mathcal{I}}$. This 
 construction, a priori, depends on the presentation of $A$ as $\cO_m/\mathcal{I}$, 
 however, we show that, in fact, it is independent of that.
 
\begin{lemma} \label{lem:pushintclosed}
    Let $\phi: \mathcal{R} \rightarrow \mathcal{O}$ be a split epimorphism of rings and let $\mathcal{I} \triangleleft \mathcal{O}$ be an ideal. Then the integral closedness of $\phi^{-1}(\mathcal{I})$ implies the same property for $\mathcal{I}$.
\end{lemma}

\begin{proof}
    Let us denote the splitting morphism by $\iota: \mathcal{O} \rightarrow \mathcal{R}$ (i.e., $\phi \circ \iota=id_{\mathcal{O}}$).  
    
    We have to prove the containment $\overline{\mathcal{I}} \subset \mathcal{I}$. Consider some $f \in \overline{\mathcal{I}}$, then, by definition, 
\begin{equation*}
    f^n + a_1 \cdot f^{n-1} + a_2 \cdot f^{n-2} + \ldots + a_n =0 \ \text{for some } a_i \in \mathcal{I}^{i}.
\end{equation*}
Using the splitting, we have $ \iota (f)^n+\iota(a_1)\cdot \iota (f)^{n-1} + \ldots+ \iota(a_n)=0$ with $\iota(a_i) \in \phi^{-1}(\mathcal{I})^i$, i.e., by definition $\iota  (f) \in \overline{\phi^{-1}(\mathcal{I})}=\phi^{-1}(\mathcal{I})$.  Hence, $f=\phi ( \iota (f) )\in \mathcal{I}$, and thus $\overline{\mathcal{I}}=\mathcal{I}$.
\end{proof}

 \begin{prop} \label{prop:reductionwdef}
     Let $A$ be finitely generated $k$-algebra
 and consider two different presentations of it: $A = \mathcal{O}_m/\mathcal{I} = \mathcal{O}_{m'}/\mathcal{I}'$. Then the corresponding integral reduction algebras agree: 
 $\mathcal{O}_m/\overline{\mathcal{I}}\cong \mathcal{O}_{m'}/\overline{\mathcal{I}'}$. 
 Now, $\mathcal{O}_m/\overline{\mathcal{I}}$  being an invariant of the $k$-algebra $A$, we denote it by $\widehat{A}$.  
 \end{prop}

\begin{proof}
    Let us denote the presentation maps by $\phi: \mathcal{O}_m = k[x_1, ..., x_m] \rightarrow A$  with kernel $\mathcal{I}$, and, respectively, $\phi': \mathcal{O}_{m'} = k[x'_1, ..., x'_{m'}] \rightarrow A$, with kernel $\mathcal{I}'$. Then we can construct a third presentation map 
    \begin{equation*}
        \phi^+:\mathcal{O}_{m+m'}=k[x_1, ..., x_m,x'_1, ..., x'_{m'}] \rightarrow A: \ x_i \mapsto \phi(x_i) \ \forall i, \ x'_j \mapsto \phi'(x'_j) \ \forall j,
    \end{equation*} 
    with kernel denoted by $\mathcal{I}^+$. Using the freeness of $\mathcal{O}_{m+m'}$ we can also define the following split epimorphisms:
    \begin{align*}
        \alpha: \mathcal{O}_{m+m'} \rightarrow \mathcal{O}_{m}: x_i \mapsto x_i \ \forall i; \ x'_j \mapsto y'_j \text{ for a certain $y_j'$ with property }\phi(y'_j)=\phi'(x'_j) \ \forall j; \\
        \alpha': \mathcal{O}_{m+m'} \rightarrow \mathcal{O}_{m'}: x_i \mapsto y_i \text{ for a certain $y_i$ with property }\phi'(y_i)=\phi(x_i) \ \forall i;\  x'_j \mapsto x'_j \ \forall j.
    \end{align*}
    They have  natural splittings denoted by $\beta: \mathcal{O}_{m} \rightarrow \mathcal{O}_{m+m'}$, respectively $\beta':\mathcal{O}_{m'} \rightarrow \mathcal{O}_{m+m'}$. Then we have the following commutative diagram:

   \begin{equation}\label{eq:DIAG}
    \begin{tikzcd}[ampersand replacement=\&]
	0 \& {\mathcal{I}} \& {\mathcal{O}_m} \\
	0 \& {\mathcal{I}^+} \& {\mathcal{O}_{m+m'}} \& A \& 0 \\
	0 \& {\mathcal{I}'} \& {\mathcal{O}_{m'}}
	\arrow[from=1-1, to=1-2]
	\arrow[from=1-2, to=1-3]
	\arrow["\beta"', shift right, curve={height=6pt}, dashed, from=1-3, to=2-3]
	\arrow["\phi", from=1-3, to=2-4]
	\arrow[from=2-1, to=2-2]
	\arrow[from=2-2, to=1-2]
	\arrow[from=2-2, to=2-3]
	\arrow[from=2-2, to=3-2]
	\arrow["\alpha"', from=2-3, to=1-3]
	\arrow["{\phi^+}", from=2-3, to=2-4]
	\arrow["{\alpha'}", from=2-3, to=3-3]
	\arrow[from=2-4, to=2-5]
	\arrow[from=3-1, to=3-2]
	\arrow[from=3-2, to=3-3]
	\arrow["{\beta'}", shift left, curve={height=-6pt}, dashed, from=3-3, to=2-3]
	\arrow["{\phi'}"', from=3-3, to=2-4]
\end{tikzcd}
\end{equation}
Notice, that $\alpha^{-1}(\mathcal{I})=\mathcal{I}^+ = \alpha'^{-1}(\mathcal{I}')$. Therefore, by relying on symmetry, it is enough to prove $\mathcal{O}_m/\overline{\mathcal{I}}\cong \mathcal{O}_{m+m'}/\overline{\mathcal{I}^+}$, or, equivalently, that $\alpha^{-1}\big(\overline{\mathcal{I}}\big)=\overline{\mathcal{I}^+}$. From the contraction property of integrally closed ideals, we always have $\alpha^{-1}\big(\overline{\mathcal{I}}\big) \supset \overline{\mathcal{I}^+}$. For the other containment it is enough to notice that $\mathcal{I}^+ \supset {\rm ker} (\alpha)$ and hence, by Lemma \ref{lem:pushintclosed}, $\alpha\big(\overline{\mathcal{I}^+}\big)$ is integrally closed and contains $\mathcal{I}$. Thus $\alpha\big(\overline{\mathcal{I}^+}\big) \supset \overline{\mathcal{I}}$ and finally $\alpha^{-1}\big(\overline{\mathcal{I}}\big)=\overline{\mathcal{I}^+}$.
\end{proof}

\begin{remark}\label{rem:pullbackandintclose}
    In general it is not true that, given a homomorphism of $k$-algebras $\alpha: A' \rightarrow A$ and an ideal $\mathcal{I} \triangleleft A$, the operations of taking contraction and taking integral closure commute, that is, in general $\alpha^{-1}\left(\overline{\mathcal{I}}\right)\neq\overline{\alpha^{-1}(\mathcal{I})}$. 
    See, for example, the map $\alpha: k[t] \rightarrow k[x, y], \ t \mapsto xy$ and the ideal $\mathcal{I}=(x^2, y^2)$: in this case $\alpha^{-1}\left(\overline{\mathcal{I}}\right)=(t) \neq (t^2) = \overline{\alpha^{-1}(\mathcal{I})}$. In fact, not even the condition of $\alpha$ being an epimorphism is enough, as demonstrated by the example from Remark \ref{rem:imageofintclosed}: in the case of the natural quotient map $k[t] \rightarrow k[t]/(t^m)$ and the trivial ideal $\mathcal{I}=0$ we have  $\alpha^{-1}\left(\overline{\mathcal{I}}\right)=(t) \neq (t^m) = \overline{\alpha^{-1}(\mathcal{I})}$. 
\end{remark}

In the case of a complex analytic algebra $A$, following Proposition \ref{prop:VAL} \textit{(i)}, we can give an intrinsic definition of the integral reduction algebra $\widehat{A}$ as follows:

\begin{proposition}\label{prop:analintred}
    Let $A$ be a complex analytic algebra. Then its integral reduction is 
    \begin{equation*}
        \widehat{A}=A/\{a \in A\,:\, \text{ for all } \gamma: A \rightarrow \mathbb{C}\{t\}/(t^m) \text{ we have }\gamma(a)=0\}.
    \end{equation*}
\end{proposition}

The possibility of a similar intrinsic formulation in the general case in the realm of quasi-valuations and quasi-realizability will be discussed around Question \ref{q:qreal=real}.

We can now prove the following generalization of part \textit{(ii)} of Proposition \ref{prop:VAL} to the general algebraic or completed setting:

\begin{cor}\label{cor:intredwdef}
    The property of a Noetherian $k$-algebra $A$ that it is of form $\mathcal{O}_m/\overline{\mathcal{I}}$ is well-defined, i.e., if $A\cong \mathcal{O}_m/\mathcal{I}$ with $\overline{\mathcal{I}}=\mathcal{I}$, then for any presentation $A \cong \mathcal{O}_{m'}/\mathcal{I}'$ we have $\overline{\mathcal{I}'}=\mathcal{I}'$. Thus, given an integrally reduced Noetherian $k$-algebra, in any presentation of it with generators and relations, the ideal of relations is inherently integrally closed. 
    
    \noindent[Therefore, if $A\cong \mathcal{O}_m/\mathcal{I}$ with $\mathcal{I} \neq \overline{\mathcal{I}}$, then $A$ is not integrally reduced.]
\end{cor}

We will need the following folklore statement (the proof is taken from the answer of user325968 to the question https://math.stackexchange.com/questions/1905186):

\begin{lemma}\label{lem:noisowfactor}
    Let $\mathcal{R}$ be a Noetherian ring and $\mathcal{I} \triangleleft \mathcal{R}$ an ideal in it. Then $\mathcal{R}/\mathcal{I} \cong \mathcal{R}$ as rings if and only if $\mathcal{I} = (0)$.
\end{lemma}

\begin{proof}
    The sufficiency is clear. For the necessity part let us assume that we have a concrete ring isomorphism $\phi: \mathcal{R}/\mathcal{I} \rightarrow \mathcal{R}$ and let us denote the natural quotient map by $q: \mathcal{R} \rightarrow \mathcal{R}/\mathcal{I}$. Now, if we assume indirectly that $\mathcal{I} \neq 0$, then we necessarily have $\mathcal{I}_1:=(q \circ \phi)^{-1}(\mathcal{I}) \supsetneq \mathcal{I}$. Similarly, $\mathcal{I}_2:=(q \circ \phi)^{-1}(\mathcal{I}_1) \supsetneq \mathcal{I}_1$ and by repeating this procedure we obtain a strictly increasing sequence of ideals in $\mathcal{R}$. This contradicts the Noetherian property.
\end{proof}

\begin{proof}[Proof of Corollary \ref{cor:intredwdef}]
    In this setting Proposition \ref{prop:reductionwdef} states that $$A\cong \mathcal{O}_m/\overline{\mathcal{I}} \cong\hat{A} \cong \mathcal{O}_{m'}/{\overline{\mathcal{I}'}} \cong (\mathcal{O}_{m'}/\mathcal{I}')/(\overline{\mathcal{I}'}/\mathcal{I}')\cong A/(\overline{\mathcal{I}'}/\mathcal{I}')$$ which, by Lemma \ref{lem:noisowfactor}, implies the integral closedness of $\mathcal{I}'$.
\end{proof}

\subsection{Direct product decomposition in the Artinian case}\,

By the {\it structure theorem for Artinian rings} (see, e.g., \cite[Theorem 8.7]{AM})
any Artinian ring $A$ can be written in a unique way (up to isomorphism) as a direct product of Artinian local rings $\prod_i A_i$. We want to verify this statement in the category of integrally reduced Artin $k$-algebras. Clearly, in the complex analytic case $A$ satisfies (\ref{eq:VAL}) if and only if each factor $A_i$ satisfies it. In general we have the following:

\begin{prop}\label{prop:proddecompintred}
    Let $A$ be an Artin $k$-algebra with $A=\prod_{i=1}^l A_i$ its direct product decomposition to local Artin $k$-algebras. Then $A$ is integrally reduced if and only if all the terms $A_i$ are. 
\end{prop}

\begin{proof}
    Let us present the terms $\{A_i\}_i$ in the form $A_i \cong \mathcal{O}_{m_i}/\mathcal{I}_i$ for all $i$, with $\mathcal{O}_{m_i}=k[x_{i, 1},\, x_{i, 2}, \ldots,\, x_{i, m_i}]$. We can also suppose without loss of generality, that the maximal ideal of $A_i$ is just $\mathfrak{m}_{m_i}/\mathcal{I}$. Then, by Corollary \ref{cor:intredwdef}, every $A_i$ is integrally reduced if and only if $\overline{\mathcal{I}_i}=\mathcal{I}_i$. 

    On the other hand, we can consider the following presentation of $A$:
    \begin{equation*}
        \mathcal{O}_m \xrightarrow{\phi} \mathcal{O}_{m_1} \oplus \ldots \oplus \mathcal{O}_{m_l} \rightarrow \mathcal{O}_{m_1}/\mathcal{I}_1 \oplus \ldots \oplus \mathcal{O}_{m_l}/\mathcal{I}_l \cong A,
    \end{equation*}
    where $\mathcal{O}_{m}=k[x_{1, 0},\, x_{1, 1}, \ldots,\, x_{1, m_1},\, x_{2, 0},\, x_{2, 1}, \ldots,\, x_{2, m_2}, \ldots,\, x_{l, 0},\, x_{l, 1}, \ldots,\, x_{l, m_l}]$ (i.e., $m = \sum_{i}m_i+l$) and $\phi:\mathcal{O}_{m} \rightarrow \oplus_i \mathcal{O}_{m_i}$ is defined by the associations:
    \begin{align*}
        &\,x_{i, j}  \mapsto (0, \ldots, 0,\, x_{i, j},\ \ \, 0, \ldots, 0)  &\, \text{ for all } 1 \leq i \leq l \text{ and } 1 \leq j \leq m_i;\\
        &\, \begin{matrix}
x_{i, 0} \mapsto 
        (0,\ldots,0, &1, &0,\ldots,0)\\
        &\uparrow&     \\
        &i&         \\
\end{matrix}  &\,\text{ for all } 1 \leq i \leq l.
    \end{align*}
    Then ${\rm ker}(\phi) = \big(\{\sum_i x_{i, 0}-1\} \cup \{x_{i, 0}\cdot x_{i', j'} \, :\, 1 \leq i, i' \leq l, \ 1 \leq j' \leq m_{i'}, \text{ such that } i \neq i'\}\big)$. We denote the ideal of relations of this presentation by $\mathcal{I}$, i.e. $\mathcal{I} = \phi^{-1}(\mathcal{I}_1 \oplus \ldots \oplus \mathcal{I}_l)$ and $A \cong \mathcal{O}_m / \mathcal{I}$. Again, by Corollary \ref{cor:intredwdef}, the Artin $k$-algebra $A$ is integrally reduced if and only if $\overline{\mathcal{I}}=\mathcal{I}$. 

    Now, if all of the $\mathcal{I}_i$ are integrally closed, then so is $\mathcal{I}_1 \oplus \ldots \oplus \mathcal{I}_l$ and, by contraction, the ideal $\mathcal{I}$ as well.  
For the other direction let us now suppose that $\mathcal{I}$ is integrally closed. Then for any $1 \leq i \leq l$ we can consider the map $\psi_i: \mathcal{O}_{m_i} \rightarrow \mathcal{O}_m$ defined by the associations $x_j \mapsto x_{i, j}, \ \forall 1 \leq j \leq m_i$. Since $\mathcal{I}_i \subset \mathfrak{m}_{m_i}$, we have $\mathcal{I}_i \subset \psi_i^{-1}(\mathcal{I})$. On the other hand, the composition $\Psi_i:\mathcal{O}_{m_i} \xrightarrow{\psi_i} \mathcal{O}_m \xrightarrow{\phi} \mathcal{O}_{m_1} \oplus \ldots \oplus \mathcal{O}_{m_l} \xrightarrow{{\rm pr}_i} \mathcal{O}_{m_i}$ is the identity map. Hence $\psi_i^{-1}(\mathcal{I})\subset \psi_i^{-1}(\phi^{-1}({\rm pr_i}^{-1}(\mathcal{I}_i)))=\mathcal{I}_i$, and thus $\mathcal{I}_i$ is also integrally closed by contraction.
\end{proof}

\subsection{Lattice homology of integrally reduced Artin $k$-algebras} \,

By unpacking part \textit{(b)} of Theorem \ref{th:propertiesideal} in this setting, we get the following:

\begin{theorem}\label{th:Artin} 
Assume that $A$ is an integrally reduced Artin $k$-algebra, i.e.,
 it can be realized as
$\cO_m/\mathcal{I}$, where $\mathcal{I}=\overline{\mathcal{I}}$.
 Then the lattice homology ${\mathbb S}\bH_*(\mathcal{I}\triangleleft \cO_m )$ is independent of the choice of the realization $\cO_m/\mathcal{I}$ of $A$, it depends merely on the isomorphism type of the $k$-algebra $A$.

 It will be denoted by ${\mathbb S}\bH_*(A)$. By part \textit{(a)} of Theorem \ref{th:propertiesideal}, $eu({\mathbb S}\bH_*(A))=\dim_k(A)$.
 In other words, the matching $A\rightsquigarrow {\mathbb S}\bH_*(A)$ is the categorification of the dimension (of such type of Artin algebras).

 If $k$ is an algebraically closed field and $A$ is local, then the weight-grading of the \emph{reduced} symmetric lattice homology  is nonnegative, i.e.,
    $\mathbb{SH}_{{\rm red},q, -2n}(A) = 0$
    for all $q\geq 0$ and $n >0$. Moreover, if $A$ is nonzero, then so is $\mathbb{SH}_{{\rm red}, 0, 0}$.
\end{theorem}

\begin{question}
    What information does the lattice homology of integrally reduced Artin $k$-algebras contain? Which invariants of the Artin $k$-algebra can we read off from its symmetric lattice homology?
\end{question}

\begin{example}\label{ex:nonisoArtin}
    We produce nonisomorphic integrally reduced Artin $k$-algebras with isomorphic lattice homology. We will rephrase Example 1.2 of \cite{KNS1} in our algebraic language (with $a=10$ specified). Let us consider $k[t]$ with the usual $t$-adic valuations and consider the following subalgebras $\widetilde{A_1}=k\langle t^{10}, t^{13}, t^{14}, t^{17}, t^{18}, t^{21}, t^{22}, t^{25}\rangle$ and $\widetilde{A_2}=k\langle t^{11}, t^{12},t^{14}, t^{16}, t^{19}, t^{20}, t^{21} \rangle$. Then the quotient algebras $A_1:=\widetilde{A_1}/(t^{30})$ and $A_2:=\widetilde{A_2}/(t^{30})$ are integrally reduced local Artin $k$-algebras (with maximal ideals denoted by $\mathfrak{m}_1$ respectively $\mathfrak{m}_2$), for which $\dim_k \mathfrak{m}_1/\mathfrak{m}_1^2 = 8 \neq 7 = \dim_k \mathfrak{m}_2/\mathfrak{m}_2^2$. On the other hand $\mathbb{SH}_*(A_1) \cong \mathbb{SH}_*(A_2)$ with graded root
    
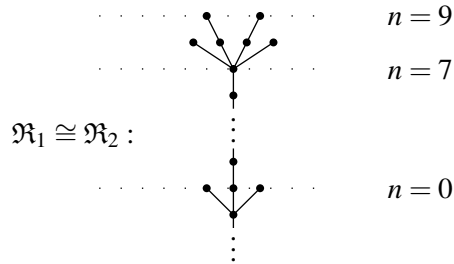
\begin{center}
\begin{picture}(200,110)(15,-115)
\linethickness{.5pt}
\put(100,-80){\circle*{3}}
\put(100,-90){\circle*{3}}
\put(100,-70){\circle*{3}}
\put(100,-35){\circle*{3}}
\put(85,-25){\circle*{3}}
\put(115,-25){\circle*{3}}
\put(90,-80){\circle*{3}}
\put(110,-80){\circle*{3}}
\put(95,-25){\circle*{3}}
\put(90,-15){\circle*{3}}
\put(110,-15){\circle*{3}}
\put(105,-25){\circle*{3}}
\put(100,-45){\circle*{3}}
\put(100,-90){\line(-1,1){10}}
\put(100,-90){\line(1,1){10}}
\put(100,-80){\line(0,1){15}}
\put(100,-80){\line(0,-1){15}}
\put(100,-35){\line(-1,2){10}}
\put(100,-35){\line(1,2){10}}
\put(100,-35){\line(3,2){15}}
\put(100,-35){\line(-3,2){15}}
\put(100,-35){\line(0,-1){15}}
\put(100,-100){\makebox(0,0){$\vdots$}}
\put(100,-55){\makebox(0,0){$\vdots$}}
\put(170,-80){\makebox(0,0)[0]{$n=0$}}
\put(170,-35){\makebox(0,0)[0]{$n=7$}}
\put(170,-15){\makebox(0,0)[0]{$n=9$}}
\put(40,-60){\makebox(0,0){$\mathfrak{R}_1 \cong \mathfrak{R}_2:$}}
\qbezier[10](50,-80)(90,-80)(130,-80)
\qbezier[10](50,-35)(90,-35)(130,-35)
\qbezier[10](50,-15)(90,-15)(130,-15)
\end{picture}
\end{center}
\end{example}
 
Let us recall, that any integrally reduced Artin $k$-algebra is the direct product of finitely many integrally reduced local Artin $k$-algebras (cf. Proposition \ref{prop:proddecompintred}). It is not very surprising, that their symmetric lattice homology modules are related by a Künneth-type formula. Nevertheless, before presenting this result we have to introduce the chain complex description of lattice homology.

\subsection{The chain complex description of lattice homology.}\label{ss:chain}\

In this subsection we present an equivalent definition of the lattice homology modules $\mathbb{H}_{\ast}(\mathbb{R}^r, w)$ associated with the lattice $\mathbb{Z}^r$ and weight function $w_0:\mathbb{Z}^r \rightarrow \mathbb{Z}$.  
Here we restrict ourselves to a very concise description. More details and proofs
can be found in \cite[Section 2]{NFilt}, about the  \emph{co}homological version in \cite[Chapter 11.1]{NBook}, or in greater generality in \cite{AgNeCurves}. This homological description will be used in the next subsection, 
in the proof of the K\"unneth formula. 

We consider the setup and the cubical decomposition of paragraph
\ref{bek:211}. Moreover, we fix an ordering of the basis $\{e_v\}_{v \in \mathcal{V}}$.
Then, on each cube $(\ell, I)$ we get an induced
orientation. 
The (classical)
boundary of an oriented cube 
 $\square_q$ has the form $\sum_k\varepsilon_k \, \square_{q-1}^k$ for some signs
$\varepsilon_k\in \{-1,+1\}$ (depending on the orientations), where the 
$(q-1)$-cubes $\square_{q-1}^k$ are the  faces of $\square_q$. In the presence of any  fixed weight function
$w=\{w_q\}_q$  we can proceed to construct a chain complex as follows. 
Let $\mathcal{L}_q$ be the free $\Z[U]$-module generated by all oriented cubes
$\square_q$ of dimension $q$ ($0\leq q\leq r$).
It is clearly a $\Z[U]$-module by
$U*(U^m\square):= U^{m+1}\square$.
Moreover, we endow $\mathcal{L}_q$ with a $\Z$-grading: by definition, the degree
of $U^m\square$ is ${\rm deg}_{\mathcal{L}}(U^m\square):=-2m-2w(\square)$.
Furthermore, we 
 define the boundary operator 
 $\partial_{w,q}:\mathcal{L}_{q}\to \mathcal{L}_{q-1}$ as follows.  First, consider the classical cubical boundary
$\partial\square_{q}=\sum_k\varepsilon_k \square ^k_{q-1}$, then
set
$$\partial_{w,q}(U^m\square_{q}):=U^m\sum_k\,\varepsilon_k\,
U^{w(\square_{q})-w(\square^k_{q-1})}\, \square^k_{q-1}.$$
Finally,  one shows that
$\partial_w\circ\partial_w=0$, i.e.,
$(\mathcal{L}_*,\partial_{w,*})$ is a chain complex.
It turns out that the homology of this chain complex is the lattice homology 
$\bH_*(\R^r,w)$. The ${\rm deg}_{\mathcal{L}}$-grading agrees with the weight grading coming from the other description, moreover, the $U$-action and the homological grading also match.

\subsection{K\"unneth formula.}\label{ss:KUNNETH}\,

Consider two integrally reduced Artin $k$-algebras of form $A'=\cO'/\mathcal{I'}$ and
 $A''=\cO''/\mathcal{I''}$, where  the corresponding 
  ideals are finite codimensional and integrally closed. Then the direct product algebra $A:=A'\times A''$  
  is also integrally reduced, hence it admits a symmetric lattice homology module. Even more, we can relate it to the lattice homology modules of its summands via Künneth-type formulas.

  We will use the chain complex description of lattice homology presented in subsection \ref{ss:chain} above. For notational clarity, we adopt the convention to denote the objects associated with the algebras $A', \ A''$ and $A' \times A''$ by the superscripts $'$,  $''$ and $^{\times}$. Then we have the following statement:

  \begin{proposition}\label{prop:chainproduct}
      Let $A'\cong \mathcal{O}'/\mathcal{I}'$ and $A''\cong \mathcal{O}''/\mathcal{I}''$ be integrally reduced  Artin $k$-algebras. If $(\mathcal{D}'=\{\frv'_{v'}\}_{v'\in \mathcal{V}'}, d_{\mathcal{I}'})$, respectively $(\mathcal{D}''=\{\frv''_{v''}\}_{v''\in \mathcal{V}''}, \mathcal{I}'')$, are CDP realizations of $\mathcal{I}'$ in $\cO'$, respectively $\mathcal{I}''$ in $\cO''$ (in the sense of Definition \ref{def:REALforideals}), then the finite collection of discrete $\mathcal{O}' \times \mathcal{{O}}''$-valuations $\mathcal{D}^{\times}=\lbrace \{\frv'_{v'}\circ pr'\}_{v'\in \mathcal{V}'}, \{\frv''_{v''}\circ pr''\}_{v''\in \mathcal{V}''}\rbrace$ and the lattice point $d_{\mathcal{I}'\times\mathcal{I}''}=( d_{\mathcal{I}'},  d_{\mathcal{I}''})$ give a CDP realization of $\mathcal{I}' \times \mathcal{I}''$ (where $pr'$ and $pr''$ denote the natural projections to the two summands of $\mathcal{O}' \times \mathcal{O}''$). Moreover, we have the following graded $\mathbb{Z}[U]$-module isomorphism of the corresponding chain complexes
    \begin{equation}\label{eq:chainproduct}
         \left(\mathcal{L}_{\ast}(\cO'/\mathcal{I}', \mathcal{D}' ), \partial' \right) \otimes_{\mathbb{Z}[U]} \left(\mathcal{L}_{\ast}(\cO''/\mathcal{I}'', \mathcal{D}''), \partial'' \right) \cong
         \left(\mathcal{L}_{\ast}(\cO'/\mathcal{I}' \times \cO''/\mathcal{I}'', \mathcal{D}^{\times}), \partial^{\times} \right) 
    \end{equation}
    induced by the chain map: $(\square' \otimes P'(U)) \otimes (\square'' \otimes P''(U)) \mapsto (\square' \times \square'') \otimes (P' \cdot P'')(U)$.
  \end{proposition}
  \begin{proof}

The collection $\mathcal{D}^{\times}$ automatically satisfies 
\begin{center}
    $\mathcal{F}_{\mathcal{D}^{\times}}( d_{\mathcal{I}'\times\mathcal{I}''})=\mathcal{I}'\times \mathcal{I}''$, with $d_{\mathcal{I}'\times\mathcal{I}''}=( d_{\mathcal{I}'},  d_{\mathcal{I}''}) \in \mathbb{Z}^{|\mathcal{V}'|}\times\mathbb{Z}^{|\mathcal{V}''|}$.
\end{center} 
Moreover, for any $(\ell', \ell'') \in \mathbb{Z}^{|\mathcal{V}'|}\times\mathbb{Z}^{|\mathcal{V}''|}$ we also have $\mathfrak{h}^{\times}(\ell', \ell'') = \mathfrak{h}'(\ell') + \mathfrak{h}''(\ell'')$ for the corresponding Hilbert functions. This identity, on one hand, implies $w^{\times}(\square' \times \square'') = w'(\square') + w''(\square'')$ on the level of weight functions for any two cubes of $R(0, d_{\mathcal{I}'})$ and $R(0, d_{\mathcal{I}''})$. On the other hand, it also implies that $\mathcal{D}^\times$ is a CDP realization of $\mathcal{I}' \times \mathcal{I}''$. Now the isomorphism (\ref{eq:chainproduct}) follows through the definitions of subsection \ref{ss:chain}.
\end{proof}

\begin{remark}\label{rem:Snproduct}
    On the level of $S_n$-spaces we have the following relation:
    \begin{equation*}
        S_{n, \mathcal{D}^{\times}} = \bigcup_{i+j=n}S_{i, \mathcal{D}'} \times S_{j, \mathcal{D}''} \subset R(0, d_{\mathcal{I}'}) \times R(0, d_{\mathcal{I}''}).
    \end{equation*}
\end{remark}

\begin{remark}\label{rem:nochainhomot}
    Note that the current form of the Independence Theorem does not imply the well-definedness of the \emph{chain homotopy type} of the chain complex $\left(\mathcal{L}_{\ast}(N \hookrightarrow M, \mathcal{D} ), \partial \right)$ associated with a CDP realization $\mathcal{D}$ of $N$ in $M$, it only asserts the well-definedness of its homology $\mathbb{H}_{\ast}(N \hookrightarrow M)$. Similarly, we do not have a well-defined \emph{filtered homotopy type} 
    \begin{center}
        $\left(\emptyset \subset S_{m_{w_{\mathcal{D}}}, \mathcal{D}} \subset \ldots \subset S_{n, \mathcal{D}} \subset S_{n+1, \mathcal{D}} \subset \ldots \subset \mathbb{R}^{|\mathcal{D}|}\right)$, 
    \end{center}
    instead, we only have well defined $S_n$ homotopy types at each level $n \in \mathbb{Z}$. For more details see subsection \ref{ss:filteredhom}. This is the reason why we can only formulate Proposition \ref{prop:chainproduct} and Remark \ref{rem:Snproduct} with specified CDP realizations. Nevertheless, by turning to the homologies --- via the Künneth spectral sequence --- we can get rid of these uncertainties.
\end{remark}

Given two (flat) chain complexes $\mathcal{L}_*$ and $\mathcal{K}_*$ over a general (commutative) ring $R$, the homologies of the separate complexes and the $\mathcal{L}_* \otimes_{R} \mathcal{K}_*$ total product complex are related by the Künneth spectral sequence (cf., e.g., \cite[Theorem 10.90]{Rotman}). Its second page consists of 
\begin{equation*}
    E^2_{p, q}=\bigoplus_{i+j=q} {\rm Tor}^R_{p}(H_i(\mathcal{L}_*), H_{j}(\mathcal{K}_*))
\end{equation*}
and converges to $H_{p+q}(\mathcal{L}_* \otimes_R \mathcal{K}_*)$. In the case when $R$ is a principal ideal domain --- such as a polynomial ring over a field) --- (or if one of the chain complexes is projective and the ring is hereditary), then this implies the Künneth formula exact sequence (cf., e.g., \cite[Corollary 10.91]{Rotman}, especially the use of the condition ${\rm Tor}^R_{\geq 2}(H_i(\mathcal{L}_*), M)=0$ for all $R$-modules $M$)
\begin{equation*}
    0 \rightarrow \bigoplus_{p+q=n}H_p(\mathcal{L}_*) \otimes_R H_q(\mathcal{K}_*) \rightarrow H_{n}(\mathcal{L}_* \otimes_R \mathcal{K}_*) \rightarrow \bigoplus_{p+q=n-1}{\rm Tor}_1^{R}(H_p(\mathcal{L}_*), H_q(\mathcal{K}_*)) \rightarrow 0.
\end{equation*}

\begin{cor}[Künneth formula]\label{prop:Künneth} Let $A'$ and $A''$ be integrally reduced  Artin $k$-algebras.

(a) If we compute the lattice homology of these algebras (introduced in Theorem \ref{th:Artin}) with coefficients from some field $\mathbb{F}$ (e.g., $\mathbb{Q}$ or $\mathbb{F}_2$), then we have the following isomorphism of $\mathbb{F}[U]$-modules:
    \begin{equation*}
        \mathbb{SH}_q(A' \times A'', \mathbb{F}) = \bigoplus_{i+j=q}\mathbb{SH}_{i}(A', \mathbb{F}) \otimes_{\mathbb{F}[U]}\mathbb{SH}_{j}(A'', \mathbb{F}) \oplus \bigoplus_{i+j=q-1}{\rm Tor}^{\mathbb{F}[U]}_1(\mathbb{SH}_{i}(A', \mathbb{F}), \mathbb{SH}_{j}(A'', \mathbb{F})).
    \end{equation*}

(b) If we use $\mathbb{Z}$-coefficients and have ${\rm Tor}^{\mathbb{Z}[U]}_{\geq 2}(\mathbb{SH}_*(A'), M)=0$, then we get the exact sequence
\begin{equation*}
    0 \rightarrow \bigoplus_{p+q=n}\mathbb{SH}_p(A') \otimes_{\mathbb{Z}[U]} \mathbb{SH}_q(A'') \rightarrow \mathbb{SH}_{n}(A' \times A'') \rightarrow \bigoplus_{p+q=n-1}{\rm Tor}_1^{\mathbb{Z}[U]}(\mathbb{SH}_p(A'), \mathbb{SH}_q(A'')) \rightarrow 0.
\end{equation*}
This happens, e.g., if $\mathbb{SH}_*(A')$ and $\mathbb{SH}_*(A'')$ has no $\mathbb{Z}$-torsion (see Remark \ref{rem:tensorandTor}), which is the case for all main examples.
\end{cor}

\begin{proof}
    The Independece Theorem \ref{th:Indep} and Theorem \ref{th:Artin} imply that for any integrally reduced Artin $k$-algebra $A$ we can choose any presentation $A\cong \mathcal{O}/\mathcal{I}$ (with $\overline{\mathcal{I}}=\mathcal{I}$) and any CDP realization $\mathcal{D}$ to compute its symmetric lattice homology.

    Now, starting from some presentations $A'\cong \mathcal{O}'/\mathcal{I}'$ and $A''\cong \mathcal{O}''/\mathcal{I}''$ we can choose the product presentation for $A' \times A''\cong (\mathcal{O}' \times \mathcal{O}'')/(\mathcal{I}' \times \mathcal{I}'')$ (see also Proposition \ref{prop:proddecompintred}). Also, by Proposition \ref{prop:chainproduct} for any given CDP realizations $\mathcal{D}'$ of $\mathcal{I}'$ in $\mathcal{O}'$ and $\mathcal{D}''$ of $\mathcal{I}''$ in $\mathcal{O}''$ we can choose the CDP realization $\mathcal{D}^\times$ of $\mathcal{I}'\times \mathcal{I}''$ in $\mathcal{O}' \times \mathcal{O}''$ to get the isomorphism (\ref{eq:chainproduct}). Then we can use the algebraic Künneth formula for the \textit{(a)} PID case and the previous exact sequence for the  \textit{(b)} case.
\end{proof}

\begin{remark}\label{rem:tensorandTor}
    Note also the following identities regarding the tensor products of $\mathbb{Z}[U]$-modules (in fact, they will all remain true for coefficients coming from a field $\mathbb{F}$, too): 
    \begin{itemize}
        \item $\mathcal{T}_n^{-}\otimes_{\mathbb{Z}[U]}\mathcal{T}_m^{-}=\mathcal{T}_{n+m}^{-}$,
        $ \mathcal{T}_n^{-}\otimes_{\mathbb{Z}[U]}\mathcal{T}_m(k)=\mathcal{T}_{n+m}(k)$,
        \vspace {1mm}
        \item $ \mathcal{T}_n(k') \otimes_{\mathbb{Z}[U]} \mathcal{T}_m(k) = \mathcal{T}_{n+m}(\min\{k, k'\})$, \vspace {1mm}
        \item ${\rm Tor}^{\mathbb{Z}[U]}_{\geq 1}(\mathcal{T}_n^{-}, T)=0$, ${\rm Tor}^{\mathbb{Z}[U]}_{\geq 2}(\mathcal{T}_m(k), T)=0$ for any $\mathbb{Z}[U]$-module $T$, and \vspace {1mm}
        \item ${\rm Tor}^{\mathbb{Z}[U]}_1(\mathcal{T}_n(k'), \mathcal{T}_m(k))=\mathcal{T}_{n+m -2\max\{k, k'\}}(\min\{k', k\})$
    \end{itemize}
\end{remark}

\begin{remark}
    Similar results are already present in the related literature: e.g., in low-dimensional topology a similar formula is used to compute the Heegaard Floer homology of the connected sum of $3$-manifolds (cf. section 6 in \cite{OSz7}). The analytic lattice homology of decomposable curve singularities has an analogous description, too, see \cite[Theorem 6.3.7]{NFilt2}.
\end{remark}

\begin{example}\label{ex:product} If $A=k[x]/(x^3)\times k[y]/(y^2)$ then ${\mathbb S}\bH_*(A)$ is the tensor product 
\begin{center}
\begin{picture}(300,50)(0,25)
\linethickness{.5pt}

\put(105,60){\circle*{3}}
\put(115,60){\circle*{3}}
\put(125,60){\circle*{3}}
\put(135,60){\circle*{3}}
\put(120,50){\circle*{3}}
\put(120,40){\circle*{3}}
\put(120,30){\makebox(0,0){$\vdots$}}
\put(158,55){\makebox(0,0){$\otimes_{\Z[U]}$}}
\put(75,55){\makebox(0,0){$({\mathbb S}\bH_0$ \ \mbox{of}}}
\put(140,55){\makebox(0,0){$)$}}
\put(193,55){\makebox(0,0){$({\mathbb S}\bH_0$ \ \mbox{of}}}
\put(245,55){\makebox(0,0){$)$}}

\put(105,60){\line(3,-2){15}}
\put(115,60){\line(1,-2){5}}

\put(125,60){\line(-1,-2){5}}
\put(135,60){\line(-3,-2){15}}
\put(120,50){\line(0,-1){15}}

\put(220,60){\circle*{3}}
\put(230,60){\circle*{3}}
\put(240,60){\circle*{3}}
\put(230,50){\circle*{3}}
\put(230,40){\circle*{3}}
\put(230,30){\makebox(0,0){$\vdots$}}

\put(220,60){\line(1,-1){10}}

\put(240,60){\line(-1,-1){10}}
\put(230,60){\line(0,-1){25}}
\end{picture}
\end{center}

Therefore, the cubical subspaces  $S_n^{\times}$ (associated with $\mathcal{D}^{\times, \, \natural}$, where $\mathcal{D}^{\times}$ consists of the pullbacks of the component's maximal ideal valuations, cf. Example \ref{ex:Ct})
are the following: $S^{\times}_{<0}=\emptyset $, 
\begin{center}
\setlength{\unitlength}{0.45mm}
\begin{picture}(200,80)(-20,5)
\linethickness{.5pt}

\multiput(10,40)(5,0){7}{\color{lightgray}\line(0,1){20}}
\multiput(10,40)(0,5){5}{\color{lightgray}\line(1,0){30}}
\multiput(10,40)(0,10){3}{\multiput(0,0)(10,0){4}{\circle*{2}}}
\put(25,72){\makebox(0,0){$S^{\times}_0$}}
\put(25,25){\makebox(0,0){$({\mathbb S}\bH_0)_0=\Z^{4\cdot3}$}}

\multiput(70,40)(5,0){7}{\color{lightgray}\line(0,1){20}}
\multiput(70,40)(0,5){5}{\color{lightgray}\line(1,0){30}}
\multiput(70,40)(0,10){3}{\line(1,0){30}}
\multiput(70,40)(10,0){4}{\line(0,1){20}}
\multiput(70,40)(0,5){5}{\multiput(0,0)(10,0){4}{\circle*{2}}}
\multiput(75,40)(0,10){3}{\multiput(0,0)(10,0){3}{\circle*{2}}}

\put(80,72){\makebox(0,0){$S^{\times}_1$}}
\put(85,25){\makebox(0,0){$({\mathbb S}\bH_0)_{-2}=\Z$}}
\put(87,10){\makebox(0,0){$({\mathbb S}\bH_1)_{-2}=\Z^{3\cdot2}$}}

{\color{darkgray}\polygon*(130,40)(160,40)(160,60)(130,60)}
\multiput(130,40)(0,5){5}{\line(1,0){30}}
\multiput(130,40)(5,0){7}{\line(0,1){20}}
\multiput(130,40)(0,5){5}{\multiput(0,0)(5,0){7}{\circle*{2}}}

\put(145,72){\makebox(0,0){$S^{\times}_n, \ n\geq 2$}}
\put(150,25){\makebox(0,0){$({\mathbb S}\bH_0)_{-2n}=\Z$}}
\put(152,10){\makebox(0,0){$({\mathbb S}\bH_{\geq 1})_{-2n}=0$}}

\end{picture}
\setlength{\unitlength}{1pt}
\end{center}

In fact, if $A=\prod_{i=1}^l k[x_i]/(x_i^{d_i})$, with $d_i\geq 1$,
then the corresponding $S_n$-spaces are the $n$-skeleta of 
$R(0, d)\subset \R^l$ (associated with its cubical decomposition), where $d=(d_1, \ldots, d_l)\in(\Z_{\geq 1})^l$  --- more precisely, they are the double of these with respect to the standard metric. 
In particular, ${\mathbb S}\bH_q(A)\not=0$  for any $0\leq q\leq l-1$.
\end{example}

As a consequence, by the  above K\"unneth formula \ref{prop:Künneth} and Proposition \ref{prop:proddecompintred},  the lattice homology of any integrally reduced Artin 
$k$-algebra is determined by the lattice homologies
of local integrally reduced Artin $k$-algebras (at least clearly over field coeffiecients). We also have the following generalization of the Nonpositivity Theorem:

\begin{theorem}\label{th:NonposArtin}
    Let $k$ be an algebraically closed field and $A$ an integrally reduced Artin $k$-algebra. Then $A$ can be written in a unique way as a direct product of integrally reduced local Artin $k$-algebras $\prod_{i=1}^l A_i$. In this case for every $n \geq l$ the $S_n$ space is contractible, thus the weight grading of the \emph{reduced} symmetric lattice homology is not lower than $2-2l$, i.e., $\mathbb{SH}_{{\rm red}, q, -2n}(A)=0$ for all $q\geq 0$ and $n \geq l$. 
    
    Moreover, if $A$ is nonzero, then so is $\mathbb{SH}_{{\rm red}, l-1, 2-2l}$.
\end{theorem}

The proof will be given in subsection \ref{ss:KünnethNonposproof} as it relies on the deformations constructed for the proof of the Nonpositivity Theorem. One might also be able to give a proof using the Künneth formula of Corollary \ref{prop:Künneth}, though for that one has to take extra care with the coefficient ring used.
\newpage
\section{Analytic lattice homology of reduced curve singularities}\label{s:deccurves}

In this section we will present the first example of categorification via lattice homology of a numerical invariant defined as the codimension of a realizable submodule: the delta invariant of reduced complex curve singularities. Our construction is equivalent with the analytic lattice (co)homology defined by T. \'Agoston and the first author in \cite{AgNeCurves}, however, the presentation from the point of view of lattice homology of realizable modules is new.

In the first part of this section
we review some needed facts about the theory of
complex analytic reduced curve singularities, then we present the categorification of the delta invariant along the lines of the general construction and theory from section \ref{s:4}. For this we need to
consider Rosenlicht's regular differential forms and a filtered version of the classical
Rosenlicht--Serre duality. Finally,
we show that this new construction agrees with the analytic lattice (co)homology of \cite{AgNeCurves} and present some examples.

\subsection{Categorification of the $\delta$ invariant (as an application of the main construction)}\, \label{bek:ANcurves}

Let $(C,o)$ be a reduced complex
curve singularity with local $\C$-algebra $\mathcal{O}=\mathcal{O}_{C,o}$. Let
$\cup_{v\in \mathcal{V}}(C_v,o)$ be the irreducible decomposition of $(C,o)$ (with $\mathcal{V}=\{1, \ldots, r\}$) and denote 
the local $\C$-algebra of
$(C_v,o)$ by $\mathcal{O}_v$. We denote the integral closure of $\mathcal{O}_v$ by $\overline{\mathcal{O}_v}= \mathbb{C}\{t_v\} $,
and we consider $\mathcal{O}_v$ as a subring of $\overline{\mathcal{O}_v}$. Similarly, we denote
 the integral closure of $\mathcal{O}$ by $\overline{\mathcal{O}}= \oplus_v \mathbb{C}\{t_v\}$. 
     The \textit{delta invariant} of $(C,o)$ is defined  as $\delta=\delta(C,o)=\dim _\mathbb{C}\, \overline{\mathcal{O}}/\mathcal{O}$.

Our objective is to establish a  categorification of $\delta(C, o)$ using the construction of lattice homology of realizable submodules from subsection \ref{ss:lathomofmods}. Therefore we need a finitely generated module and a realizable submodule of it with codimension $\delta$. A result of Serre will provide these in the realm of differential forms. In order to present it, we first have to introduce some notations.

The inclusion $\mathcal{O} \hookrightarrow \overline{\mathcal{O}}$ corresponds to the normalization map $n_C:\overline{(C,o)}\to (C,o)$, where $\overline{(C,o)}$ is a multigerm comprised of $r$ smooth germs. 
Let $\Omega^1(*)$ denote the module of germs of meromorphic differential $1$-forms on the normalization
$\overline{(C,0)}$  with possible  poles (of any order) at most in the preimage $\overline {o}:=n_C^{-1}(o)$ of the singular point.
Notice, that, through the ring inclusion $n_C^*: \mathcal{O} \hookrightarrow \overline{\mathcal{O}}$, this
$\Omega^1(*)$
is a torsion-free ${\mathcal{O}}$-module. 
We can consider the following pairing:
\begin{center}
    $\overline{\mathcal{O}} \times \Omega^1(*) \rightarrow \mathbb{C}, \ (f, \alpha) \mapsto \sum_{p\in\overline{o}} \, {\rm res}_p(f\cdot\alpha)$.
\end{center} 
Then {\it Rosenlicht's regular differential forms} are 
 defined as meromorphic forms  $\alpha\in \Omega^1(*)$, such that $\sum _{p\in \overline{o}} {\rm res}_p(n_C^*(f
 )\cdot\alpha)=0 \ 
 \text{ for all } f\in \calO$. Clearly, they form a (torsion-free) $\mathcal{O}$-submodule in $\Omega^1(*)$;
 this module will be denoted by $\omega^R_{C}$. In fact, it is canonically isomorphic with the dualizing module of Grothendieck associated with the germ $(C,o)$ (see e.g., \cite[Chapter VIII]{AltKlei}). 
Let  $\Omega^1 _{\overline{C}}
$ denote the submodule (in $\omega_C^R$) of multigerms of holomorphic differential $1$-forms on $\overline{(C,o)}$. Then, by \cite[Ch 2 Sect 9]{Serre} (see also \cite[Section 1]{BG80} or \cite[Section 2]{Rosenlicht}), 
the above  
 pairing gives a perfect duality
 \begin{equation}\label{eq:ROS}
\overline{\calO}/ \calO\ \times \ \omega^R_{C}/ \Omega^1_{\overline{C}}\to \C,\ \ \
 [f]\times [\alpha]\mapsto \sum_{p\in\overline{o}} \, {\rm res}_p(f\cdot\alpha); 
 \end{equation}
 in particular $\delta = {\rm codim}_{\mathbb{C}}(\Omega^1_{\overline{C}} \hookrightarrow \omega_C^R)$.
 We can also conclude that $\omega_C^R$ is a finitely generated $\mathcal{O}$-module. 

 Moreover, the submodule $\Omega^1_{\overline{C}} \leq \omega_C^R$ is realizable (in the sense of Definition \ref{def:REAL}). Indeed, consider the discrete valuations corresponding to the normalization of the components $(C_v, o)$ of $(C, o)$: 
 \begin{equation}\label{eq:normalizationvaluation}
     \mathfrak{v}_v: \mathcal{O} \rightarrow \overline{\N}, \ g \mapsto {\rm ord}_{t_v}(\pi_v(n_C^*(g))) \hspace{3mm} (\text{with } \pi_v: \overline{\mathcal{O}} \to \overline{\mathcal{O}_v}\cong \mathbb{C}\{t_v\} \text{ the natural projection} ) .
 \end{equation}  
 These can 
  be naturally extended to $\omega_C^R$ in the sense of Definition \ref{def:edv} as 
\begin{equation}\label{eq:v_v^omega}
    \mathfrak{v}_v^{\omega}:\omega_C^R \rightarrow \mathbb{Z} \cup\{\infty\}, \ \alpha=\sum_{u=1}^r a_u(t_u)dt_u  \mapsto {\rm ord}_{t_v}a_v,
\end{equation} 
  where $a_u(t_u)$ is a 
  convergent Laurent series in the variable $t_u$.
Now, for the collection $\{ (\mathfrak{v}_v, \mathfrak{v}_v^{\omega})\}_{v \in \mathcal{V}}$ of extended discrete valuations we can assign the $\mathbb{Z}^r$-indexed multifiltration
\begin{center}
$\ell \mapsto {\mathcal F}^{\omega}(\ell)=\lbrace \alpha \in \omega_C^R\,:\, \mathfrak{v}_v^{\omega}(\alpha) \geq \ell_v \ \mbox{for all} \ v\in \mathcal{V}\}$ \hspace{3mm} (compare with (\ref{eq:fdMl})).
\end{center}
Then clearly 
  $\Omega^1 _{\overline{C}}
=\mathcal{F}^{\omega}(0)$. Therefore we can use the Independence Theorem \ref{th:IndepMod} and Theorem \ref{th:properties} to conclude

\begin{cor}
    For a reduced complex analytic curve singularity $(C, o)$ the lattice homology $\mathbb{Z}[U]$-module $\mathbb{H}_*(\Omega^1 _{\overline{C}} \hookrightarrow_{\mathcal{O}} \omega_C^R)$ is well-defined and gives a categorification of $\delta(C, o)$. Moreover, it can be computed on a finite rectangle depending on the specific CDP realization chosen.
\end{cor}

Moreover, since $\mathcal{O}$ is local, by Theorem \ref{th:upperbound} we also have:

\begin{cor}\label{cor:upperboundcurves}
    For every 
$n \geq 1$ the space $S_n$ is contractible, hence the weight-grading of the \emph{reduced} lattice homology is nonnegative, i.e.,
\begin{center}
    $\mathbb{H}_{{\rm red}, q}(\Omega^1_{\overline{C}} \hookrightarrow_{\mathcal{O}} \omega_C^R) = \bigoplus_{0 \geq n \geq m_w}\mathbb{H}_{{\rm red}, q, -2n}(\Omega^1_{\overline{C}} \hookrightarrow_{\mathcal{O}} \omega_C^R)$ \hspace{5mm} for all $q \geq 0$.
\end{center} 
Moreover, if $(C, o)$ is non-smooth, then the upper bound is sharp: $S_0$ is not connected.
\end{cor}

In fact, $\mathbb{H}_*(\Omega^1 _{\overline{C}} \hookrightarrow_{\mathcal{O}} \omega_C^R)$ is isomorphic to the analytic lattice homology of the curve singularity $(C,o)$ defined by \'Agoston and the first author in \cite{AgNeCurves}. We will devote the following two subsections to present that construction, some of its properties and prove the isomorphism. In that language the above result (i.e., Corollary \ref{cor:upperboundcurves}) was already established in \cite[Theorem 6.1.1]{KNS2} by Kubasch and the authors. In fact, our proof of the general Nonpositivity Theorem \ref{th:upperbound} is just an adapted 
generalization of the one given there.

\subsection{Analytic lattice homology of reduced curve singularities after \cite{AgNeCurves}.} \label{ss:ancurves}\,

We fix a reduced curve singularity $(C, o)$ and will use the notations of the previous subsection.  Consider the lattice $\mathbb{Z}^r$ (where $r$ is the number  irreducible components) with its natural basis $\{e_v\}_{v=1}^r$.
Then
$\overline{\calO}$ has a natural multifiltration indexed by $\ell\in (\Z_{\geq 0})^r$ given by
$\overline{\calF}(\ell):=\{g \in \overline{\mathcal{O}}\,:\, \frakv(g)\geq \ell\}$ (compare with (\ref{eq:fdMl})), where $\mathfrak{v}$ denotes the tuple $(\mathfrak{v}_1, \ldots, \mathfrak{v}_r)$ of discrete valuations corresponding to the normalization of the components (cf. (\ref{eq:normalizationvaluation})). (Under the identification $\overline{\mathcal{O}}\cong \mathbb{C}\{t_1\} \times \ldots \times \mathbb{C}\{t_r\}$ the ideal $\overline{\mathcal{F}}(\ell)$ corresponds to $(t_1^{\ell_1} , \ldots , t_r^{\ell_r}) \overline{\mathcal{O}}$.) This induces an ideal filtration of
$\calO$ by $\calF(\ell):=\overline{\calF}(\ell)\cap \calO$, and also in $\overline {\calO}/\calO$ by
$\calF^\bullet (\ell):=(\, \overline{\calF}(\ell)+\calO\,)/\calO
$.
Then $\ell \mapsto \hh(\ell)=\dim \calO/\calF(\ell)$ is called the \textit{Hilbert function of $(C,o)$}. 
It is increasing and satisfies $\mathfrak{h}(0)=0$ and 
$\mathfrak{h}(\ell+e_v)-\mathfrak{h}(\ell)\in\{0,1\}$ for any $v \in \mathcal{V}$ (since the considered valuations are $0$-dimensional in the sense of  \cite{CutVal}). Furthermore, we also set   $\hh^\bullet (\ell):=\dim \mathcal{F}^{\bullet}(\ell)$.
Then $\hh^\bullet$ is decreasing, $\hh^{\bullet}(0)=\delta$ and 
\begin{equation}\label{eq:dag}
 \hh(\ell)-\hh^\bullet(\ell)+\delta =\dim\, \overline{\calO}/\overline{\calF}(\ell)=|\ell|.\end{equation}

\begin{define}\cite{AgNeCurves}
In the construction of the analytic lattice (co)homology of $(C,o)$ we consider only the first quadrant $R(0, \infty):=(\R_{\geq 0})^r$.
The weight function on the lattice points of $(\Z_{\geq 0})^r$ is defined by
\begin{equation}\label{eq:w_0analytic}
w_{an, 0}(\ell)=\hh(\ell)+\hh^\bullet (\ell)-\hh^\bullet(0)= 2\cdot \hh(\ell)-|\ell|, 
\end{equation}
and $w_{an, q}(\square):=\max\{w_{an, 0}(\ell)\,:\, \mbox{$\ell$ is a vertex of $\square$}\}$ (cf. (\ref{eq:9weight})).
Then corresponding bigraded $\mathbb{Z}[U]$-module $\bH_*(R(0, \infty),w_{an})$ is called the \textit{analytic lattice homology of the reduced curve singularity $(C, o)$} (after \cite{AgNeCurves}) and is denoted by 
 $\bH_{an,*}(C,o)$.
\end{define}

\begin{remark}\label{rem:hh^bulletCDP}
    Since 
$\mathfrak{h}^\bullet (\ell+e_v)-\mathfrak{h}^\bullet (\ell)=\mathfrak{h} (\ell+e_v)-
    \mathfrak{h}(\ell)-1\in \{0, -1\}$ (cf. $(\ref{eq:dag})$), the pair $(\mathfrak{h},\mathfrak{h}^\bullet)$
    satisfies the Combinatorial Duality Property (cf. Definition \ref{def:COMPGOR}). 
\end{remark}

\bekezdes \label{par:conductor} In order to present further properties of this construction, we need the notion of the conductor element. 
Let ${\mathcal C}=(\mathcal{O}:\overline{\mathcal{O}})$  denote the {\it conductor ideal} of $\overline{\mathcal{O}}$, i.e., the largest ideal
of $\mathcal{O}$ which is an ideal of $\overline{\mathcal{O}}$ as well. It has the form
$(t_1^{c_1}, \ldots, t_r^{c_r})\overline{\mathcal{O}} $. The lattice point $\mathbf{c}=(c_1,\ldots , c_r) \in (\mathbb{Z}_{\geq 0})^r$ is called the
{\it conductor} or {\it conductor element} of $(C, o)$.

\begin{remark} \label{rem:condofapc}
    For a \emph{plane} curve singularity germ $(C, 0)\subset (\mathbb{C}^2, 0)$ the components of the conductor are of form $c_v=\mathbf{c}(C_v) + \sum_{w \neq v} i_0(C_w, C_v)$ for all $v \in \mathcal{V}$, 
 where ${\bf c}(C_v)$ is the conductor of the component $(C_v,0)$ and 
 $i_0(\ ,\,) $ denotes the intersection multiplicity at $0\in\C^2$. Furthermore, $\mathbf{c}(C_v)$ equals the Milnor number of the branch $(C_v, 0)$.
For a formula of $c_v$ in the general case see \cite{D'Anna}.
\end{remark} 

From the definitions, it is clear that $\mathfrak{h}^\bullet(\ell) = 0$ for $\ell\geq {\bf c}$
(since $\overline {\calF}(\ell)\subset
\mathcal{C}\subset \calO$ for such $\ell$ --- in fact, for $\ell = \mathbf{c}$ we have $\mathcal{F}(\mathbf{c})=\overline{\mathcal{F}}(\mathbf{c})=\mathcal{C}$). 

 \begin{theorem}\label{cor:EUcurves} \cite{AgNeCurves}
(a) For  any $d\geq {\bf c}$ the natural inclusion $S_n\cap R(0,d)\hookrightarrow S_n$ is a homotopy equi\-valence.

(b)  One has a bigraded $\Z[U]$-module isomorphism  $\bH_{an,*}(C,o)\cong\bH_*(R(0,d),w_{an})$ for any $d\geq \mathbf{c}$
induced by the above inclusion map. Therefore, $\bH_{an,*}(C,o)$ is determined by the weighted  cubes
of the finite rectangle $R(0,\mathbf{c})$. 

(c) $eu(\bH_{an,*}(C,o))=\delta(C,o)$, that is,
$\bH_{an,*}(C,o)$ is a `categorification' of \,$\delta(C,o)$.
\end{theorem}

\begin{remark} \label{rem:AlexesHilbert}
(a)  Ágoston and the first author proved in \cite{AgNeCurves} that a flat deformation $\{(C_t, o)\}_{t \in (\mathbb{C}, 0)}$ of curve singularities induces a graded $\mathbb{Z}[U]$-module morphism $\mathbb{H}_{an, 0}(C_{t \neq 0}, o) \rightarrow \mathbb{H}_{an, 0}(C_{t=0}, o)$. Moreover, the Functoriality Conjecture states that this map can be extended to the full lattice homology $\mathbb{Z}[U]$-module.

(b) The minimal value $m_{w_{an, 0}}$ of the analytic weight function is a `new and mysterious' invariant of the reduced curve singularity $(C, o)$. By work of Hof and the first author this invariant is the main ingredient of the characterization of the Cohen-Macauley type of the germ \cite{HofN}.

(c) By previous work of the authors together with Kubasch \cite[Theorem 4.2.1]{KNS2}, the analytic lattice homology of reduced curve singularities detects the Gorenstein property. In fact, the Gorenstein property corresponds to the weight function being symmteric with respect the conductor element $\mathbf{c}$ which property can also be read off from the ranks of the graded pieces of $\mathbb{H}_{an, 0}(C, o)$. For more on the Gorenstein case see subsection \ref{ss:GorCurves}.  

(d) 
In the case of a {\it plane} curve singularity,   $\hh$ is computable from the partial
multivariable Alexander polynomials  \cite[3.3--3.4]{GorNem2015}, while these Alexander polynomials from the embedded resolution graph of $(C,0)\subset (\bC^2,0)$ \cite{EN}.
In particular, the Hilbert function $\mathfrak{h}$, and, hence, the weight function $w_{an, 0}$ is an embedded topological invariant of 
$(C,0)\subset (\C^2,0)$. In fact, one can read off the multiplicity $m(C, 0)$ from the graded $\mathbb{Z}[U]$-module structure of $\mathbb{H}_{an, 0}(C, 0)$ \cite[Theorem 5.2.5]{KNS2}. Moreover, by work of Kubasch and the second author, for \textit{irreducible} plane curve singularities the analytic lattice homology $\mathbb{H}_{an}(C, 0)=\mathbb{H}_{an, 0}(C, 0)$ is a complete embedded topological invariant  \cite[Theorem 3.1.1]{KS3}. 

(e) The first author introduced in \cite{NFilt} an increasing subspace filtration on the $S_{n}$ spaces, which further refines the analytic lattice homology $\mathbb{H}_{an}(C, o)$. In fact, the corresponding spectral sequence encodes rich geometric information about the singularity: it is a complete topological invariant for \textit{ any (not necessarily irreducible)} plane curve singularity, it recovers the motivic Poincaré series, while the elements of the first page can be associated with the Heegaard Floer Link homology of the link.
\end{remark}

\subsection{Equivalence of the two definitions.}\label{ss:equivforcurves}\,

In this subsection we prove that the two lattice homological categorifications of the delta invariant (presented in the previous two subsections) agree.

\begin{theorem}\label{th:equivforcurves}
    For any reduced curve singularity $(C, o)$ we have: 
    \begin{center}
        $\mathbb{H}_*(\Omega^1_{\overline{C}} \hookrightarrow_{\mathcal{O}} \omega_C^R) \cong \mathbb{H}_{an,*}(C,o)$.
    \end{center}
\end{theorem}

\begin{proof}
    Recall, that the analytic lattice homology of $(C, o)$ (after \cite{AgNeCurves}) is computed on a lattice of dimension $r$ (the number of irreducible components), on the finite rectangle $R(0, \mathbf{c})$ and weight function derived from the valuations corresponding to the normalization of the components (cf. (\ref{eq:w_0analytic})). We compare this setup with the one corresponding to the realization $\mathcal{D}:=\{ (\mathfrak{v}_v, \mathfrak{v}_v^{\omega})\}_{v \in \mathcal{V}}$ (presented in  (\ref{eq:normalizationvaluation}) and (\ref{eq:v_v^omega})) of $\Omega^1_{\overline{C}} \leq \omega_C^R$ (see Definition \ref{def:LCofMOD}). This latter one is also defined on $\mathbb{Z}^r$ and, since we use the same valuations on $\mathcal{O}$, the Hilbert functions $\mathfrak{h}=\mathfrak{h}_{\mathcal{D}}$ agree (on $R(0, \mathbf{c})$) in the two weight function definitions, thus, it is enough to compare the functions $\ell \mapsto \mathfrak{h}_{\mathcal{D}}^\circ(\ell) = \dim_{\mathbb{C}}\omega^R_C\big/ \mathcal{F}_{\mathcal{D}}^\omega(-\ell)$ (defined in (\ref{eq:handhcirc0})) and $\ell \mapsto\mathfrak{h}^\bullet(\ell) = \dim_{\mathbb{C}}\mathcal{F}^\bullet(\ell)$.

    Notice that the extensions $\{\mathfrak{v}_v^{\omega}\}_{v \in \mathcal{V}}$ of the valuations $\{\mathfrak{v}_v\}_{v \in \mathcal{V}}$ can also be defined on the whole $\Omega^1(*)$ with the same formula (\ref{eq:v_v^omega}). We denote the resulting extensions by $\{\mathfrak{v}_v^{\Omega^1}\}_{v \in \mathcal{V}}$ and the corresponding multifiltration on $\Omega^1(*)$ by $\mathcal{F}^{\Omega^1}_{\mathcal{D}}$. Then still $\Omega^1 _{\overline{C}}
=\mathcal{F}_{\mathcal{D}}^{\Omega^1}(0)$ and we have the following perfect duality between $\overline{\mathcal{O}}/\mathcal{C}=\overline{\mathcal{O}}/\overline{\mathcal{F}}(\mathbf{c})$ and 
$\mathcal{F}_{\mathcal{D}}^{\Omega^1}(-\mathbf{c})/\Omega^1_{\overline{C}}$:
 \begin{equation}\label{eq:perfectdualuptoc}
 {\rm Res}:\overline{\calO}/ \mathcal{C}\ \times \ \mathcal{F}_{\mathcal{D}}^{\Omega^1}(-\mathbf{c})/ \Omega^1_{\overline{C}}\ \longrightarrow \ \C,\ \ \
 [f]\times [\alpha]\mapsto \sum_{p\in\overline{o}} \, {\rm res}_p(f\cdot\alpha).
 \end{equation}
 Indeed, both have finite $\mathbb{C}$-dimension $|\mathbf{c}|=\sum_{v=1}^r c_v$ and it is easy to construct for every $f \in \overline{\calO}  \setminus \mathcal{C}$ an $\alpha_f \in \mathcal{F}_{\mathcal{D}}^{\Omega^1}(-\mathbf{c})$ with ${\rm Res}([f], [\alpha_f])\neq 0$ (and similarly, for any $\alpha \in \mathcal{F}_{\mathcal{D}}^{\Omega^1}(-\mathbf{c}) \setminus \Omega^1_{\overline{C}}$ an $f_\alpha \in \overline{\mathcal{O}}$ with ${\rm Res}([f_{\alpha}], [\alpha])\neq 0$).
 Moreover, this duality respects the multifiltrations $\overline{\mathcal{F}}/\mathcal{C}$ and $\mathcal{F}_{\mathcal{D}}^{\Omega^1}/\Omega^1_{\overline{C}}$:
 \begin{align}\label{eq:dualityoffiltrations}
     \forall \ell \in (\mathbb{Z}_{\geq 0})^r\cap R(0,\mathbf{c}): \ \big(\overline{\mathcal{F}}(\ell)/&\,\mathcal{C}\big)^{\perp_{\rm Res}}=\mathcal{F}_{\mathcal{D}}^{\Omega^1}(-\ell)/\Omega^1_{\overline{C}}\,, \nonumber \\ \text{ since }\forall f \in \overline{\mathcal{F}}(\ell),  & \ \alpha \in \mathcal{F}_{\mathcal{D}}^{\Omega^1}(-\ell): \ f \cdot \alpha \in \Omega^{1}_{\overline{C}}\,, \\
      \text{ and } {\rm codim}  & (\overline{\mathcal{F}}(\ell) \hookrightarrow \overline{\mathcal{O}})=\dim \mathcal{F}_{\mathcal{D}}^{\Omega^1}(-\ell)/\Omega^1_{\overline{C}} = |\ell|. \nonumber
 \end{align}
 The dual space of $\mathcal{O}$ under this pairing (\ref{eq:perfectdualuptoc}) is exactly $\omega_C^R$, the module of { Rosenlicht's regular differential forms}, with $\left(\mathcal{O}/\mathcal{C}\right)^{\perp_{\rm Res}}=\omega^R_{C} /
 \Omega^1_{\overline{C}}$. (Indeed, notice that in formula (\ref{eq:perfectdualuptoc})  we can replace $\mathbf{c}$ by any effective lattice point $\ell \in (\mathbb{Z}_{\geq 0})^r$ --- i.e., replace $\mathcal{C}$ by $\overline{\mathcal{F}}(\ell)$ and $\mathcal{F}_{\mathcal{D}}^{\Omega^1}(-\mathbf{c})$ by $\mathcal{F}_{\mathcal{D}}^{\Omega^1}(-\ell)$ --- and and still have a perfect duality with property (\ref{eq:dualityoffiltrations}). Hence, Rosenlicht's regular differential forms, being orthogonal to every element of $\mathcal{C}$, must have valuations $\geq - \mathbf{c}$.)
 This yields the perfect duality (\ref{eq:ROS}) between $\omega^R_{C}/ \Omega^1_{\overline{C}}$ and
 $\overline{\calO}/\calO$.
 Then, from relation (\ref{eq:dualityoffiltrations}) we get \begin{equation*}
     \left(\frac{\overline{\mathcal{F}}(\ell) + \mathcal{O}}{\mathcal{C}}\right)^{\perp_{\rm Res}}=\left({\overline{\mathcal{F}}(\ell)/\mathcal{C}}\right)^{\perp_{\rm Res}} \cap \left(\mathcal{O}/\mathcal{C}\right)^{\perp_{\rm Res}}=\mathcal{F}_{\mathcal{D}}^{\Omega^1}
     (-\ell)/\Omega^1_{\overline{C}} \cap \omega^R_{C}/\Omega^1_{\overline{C}}=
     \frac{\mathcal{F}_{\mathcal{D}}^{\omega}(-\ell)}
     {\Omega^1_{\overline{C}}},
 \end{equation*}
 hence the multifiltrations $\ell \mapsto \mathcal{F}^{\bullet}(\ell) = \frac{\overline{\mathcal{F}}(\ell) + \mathcal{O}}{\mathcal{O}}$ and $\ell \mapsto  \frac{\mathcal{F}_{\mathcal{D}}^{\omega}
     (-\ell) }{\Omega^1_{\overline{C}}}$
 give for any fixed $\ell \in (\mathbb{Z}_{\geq 0})^r$ orthogonal (in the sense of the perfect duality (\ref{eq:ROS})) subspaces of $\overline{\mathcal{O}}/\mathcal{O}$ and $\omega^R_{C}/ \Omega^1_{\overline{C}}$. Therefore, $\mathfrak{h}^\bullet(\ell) = \mathfrak{h}_{\mathcal{D}}^\circ(\ell)$ for any $\ell \in (\mathbb{Z}_{\geq 0})^r$. 
 
 On one hand, it follows that $\mathfrak{h}_{\mathcal{D}}^\circ(\ell) =0$ for all $\ell \geq \mathbf{c}$. On the other hand, we know that $\mathfrak{h}_{\mathcal{D}}(\ell) =0$ for all $\ell \leq 0$. Via Observation \ref{obs:hhMmatroid}, the assumptions of Proposition \ref{prop:rectCDP->fullCDP} are satisfied, hence the CDP property for $\mathfrak{h}$ and $\mathfrak{h}^\bullet$ (see Remark \ref{rem:hh^bulletCDP}) implies the Combinatorial Duality Property for the pair $(\mathfrak{h}_{\mathcal{D}}, \mathfrak{h}^\circ_{\mathcal{D}})$ on the whole lattice $\mathbb{Z}^r$. Therefore, the realization $\mathcal{D}=\{ (\mathfrak{v}_v, \mathfrak{v}_v^{\omega})\}_{v \in \mathcal{V}}$ of $\Omega^1_{\overline{C}} \leq \omega_C^R$ investigated before is in fact a CDP realization of $\Omega^1_{\overline{C}}$, hence the weight function $w_{\mathcal{D}, 0}: \ell \mapsto \mathfrak{h}_{\mathcal{D}}(\ell) + \mathfrak{h}_{\mathcal{D}}^\circ(\ell)-\mathfrak{h}_{\mathcal{D}}^{\circ}(0)$ indeed computes $\mathbb{H}_*(\Omega^1_{\overline{C}} \hookrightarrow_{\mathcal{O}} \omega_C^R)$. However, as we already proved, $w_{\mathcal{D},0}(\ell)=w_{an, 0}(\ell)$ for any lattice point $\ell \in (\mathbb{Z}_{\geq 0})^r$, hence part \textit{(a)} of Theorem \ref{th:properties} implies the desired isomorphism.
\end{proof}

\begin{cor}\label{cor:dD=cforcurves}
    For the realization $\mathcal{D}$  of $\Omega^1_{\overline{C}} \leq \omega_C^R$ presented in (\ref{eq:normalizationvaluation}) and (\ref{eq:v_v^omega}) the lattice point $d_{\mathcal{D}}$ (cf. Notation \ref{not:d_D}) coincides with the conductor element $\mathbf{c}$. Indeed, since $\mathfrak{h}_{\mathcal{D}}^\circ(\ell)=\mathfrak{h}^\bullet(\ell)$ for any $\ell \in (\mathbb{Z}_{\geq 0})^r$, we can rely on the facts that (by the definition of the conductor ideal and the filtration $\mathcal{F}^\bullet$): $\mathfrak{h}^\bullet(\mathbf{c})=0$ and $\mathfrak{h}^\bullet(\mathbf{c}-e_v)=1$ for all $v \in \mathcal{V}$.

    Moreover, by Proposition \ref{prop:pullbackring} \textit{(2)} we get that $\mathcal{C}={\rm Ann}_{\mathcal{O}}\big(\omega^R_C/\Omega^1_{\overline{C}}\big)$.
 \end{cor}

\subsection{Gorenstein case}\label{ss:GorCurves}\,

The original definition of Grothendieck says that a quotient of a regular local ring is \textit{Gorenstein} if it is Cohen--Macaulay and its (canonical) dualizing module is free of rank one (see the first written reference in \cite[Subsection II.6]{SerreGor}). In the case of a reduced complex analytic curve singularity $(C, o)$ this means that the $\calO$-module $\omega^R_{C}$ of Rosenlicht's regular differential forms is free of rank one, that is,
there exists a `Gorenstein form' $\alpha_{Gor}\in  \omega^R_{C}$ such that $ \omega^R_{C}
=\alpha_{Gor} \calO$. Moreover, by work of Rosenlicht \cite{Rosenlicht} and Serre \cite{Serre}, $(C,o)$ is {Gorenstein}  if and only if
$\dim (\overline{\mathcal{O}}/\mathcal{O})=\dim(\mathcal{O}/\mathcal{C})$
(equivalently, using our previous notations from subsection \ref{ss:ancurves}, $\delta=\hh({\bf c})$ or $|{\bf c}|=2\delta$). Thus $\alpha_{Gor}\mathcal{C}=\Omega^1_{\overline{C}}$ and, using the perfect pairing from (\ref{eq:ROS}), we get a perfect pairing
 \begin{equation}\label{eq:ROS2}
D: \overline{\calO}/\calO\ \times \ \calO/\mathcal{C}\to \C,\ \ \
 [f]\times [g]\mapsto \sum_{p\in\overline{o}} \, {\rm res}_p(f \cdot n_C^*(g)\cdot \alpha_{Gor}).\end{equation}

\begin{cor}\label{cor:curveArtin}
    Let $(C, o)$ be a reduced Gorenstein curve singularity. Then $$\mathbb{H}_{an,*}(C, o) \cong \mathbb{H}_*(\Omega_{\overline{C}}^1 \hookrightarrow_{\mathcal{O}} \omega_C^R) \cong \mathbb{SH}_*(\mathcal{C} 
    \hookrightarrow\mathcal{O}).$$
    In particular, the analytic weight function is symmetric with respect to the conductor element, i.e., for all $\ell \in R(0, \mathbf{c})\cap \mathbb{Z}^r:\ w_{an, 0}(\ell) = w_{an, 0}(\mathbf{c}-\ell)$. (In fact, T. Ágoston and the first author already proved in \cite[Subsection 4.4]{AgNeCurves}, that $\frh^\bullet(\ell)=\frh^{sym}_{{\bf c}}(\ell) \ \forall \ell$.)
\end{cor}

\begin{proof}
    Combine Theorem \ref{th:equivforcurves} and and the discussion behind Definition \ref{def:LC}. 
\end{proof}

 However, we warn the reader that if the reduced curve singularity $(C,o)$ is not Gorenstein then usually  $\bH_{an,*}(C, o) \not={\mathbb{SH}}_*(\mathcal{C}\triangleleft \cO_{C,o})$. For example, take the  numerical semigroup $\mathcal{S}=\langle 4,5, \geq 8\rangle \subset \bN$
 and let $(C,o)$ be the singularity of ${\rm Spec}(\bC[\mathcal{S}])$. Then the conductor of $(C,o)$ is 8, 
 the delta invariant is $\delta=5$ and ${\rm codim}(\mathcal{C}\hookrightarrow \mathcal{O})=3$.
 We invite the reader to determine both $\bH_{an,0}(C,o)$ and ${\mathbb S}\bH_0(\mathcal{C}\hookrightarrow \mathcal{O})$:
 they are not isomorphic (the first one has Euler characteristic $\delta=5$, the other one 3). 
 For another example see  Example \ref{ex:345534}.
 
 In fact,  A. A. Kubasch and
the authors proved that the Gorenstein property can be read from the rank of the kernel of the $U$-action in bigrading $(0, 0)$ \cite[Theorem 4.2.1]{KNS2} (see also Remark \ref{rem:AlexesHilbert} (c)).

\begin{remark}
Grothendieck's definition of the Gorenstein property was generalized to all local rings, and thereby all rings, by Bass \cite{Bass} (see also \cite{Huneke}), which gives today's standard commutative algebra formalism.
\end{remark}

\subsection{Computational tools and examples}
\label{ss:excurves}\,

\bekezdes \textbf{Semigroup of values.} In order to compute the analytic lattice homology of reduced curve singularities one is required to obtain the values of the Hilbert function. As this is a very important classical invariant, it has well-understood connections to other invariants which can offer paths towards its computation.  
For example, the Hilbert function $\mathfrak{h}$ contains equivalent information with the \textit{semigroup of values} of the reduced 
 curve singularity $(C,o)$, defined as
\begin{equation*}
    \mathcal{S}_{C,o} = \big\{ \, \big(\mathfrak{v}_1(g), \dots, \mathfrak{v}_r(g)\big)  \in (\bZ_{\geq 0})^r \ : \ g \in \mathcal{O} \text{ is not a zero-divisor} \, \big\}.
\end{equation*}
In general it is a hard problem to classify those numerical semigroups $\mathcal{S}\subset (\bZ_{\geq 0})^r$ which can be 
realized as semigroup of values of curve singularities. However, there is a family of semigroups, defined
by certain combinatorial properties
(listed, e.g, in subsection 2.8 of \cite{KNS2}), called 
\textit{good} semigroups of  $(\bZ_{\geq 0})^r$, which contains all the semigroups of values of curve germs. 
For more see \cite{good,delaMata87,delaMata,Garcia}. The connection of the Hilbert function with the 
semigroup of values is manifested as follows: 
\begin{equation}\label{eq:Hilbertfelcsopbol}
    \text{ for any } \ell \in \mathbb{Z}^r, 1 \leq v \leq r:\ \mathfrak{h}(\ell +e_v)-\mathfrak{h}(\ell)=\begin{cases}
        1, &\ \text{if } \exists\, s \in \mathcal{S}_{C, o}\,:\, s \geq \ell \text{ and }s_v =\ell_v;\\
        0, &\ \text{otherwise;}
    \end{cases}
\end{equation}
(cf. subsection 2.2 of \cite{KNS2}), whereas $\mathcal{S}_{C,o}$ has the description 
\begin{equation}\label{eq:Sfromh}
  \mathcal{S}_{C,o}=\{\, \ell\ \in (\bZ_{\geq 0})^r \, : \,  \hh(\ell+e_v)>\hh(\ell) \text{ for all } v \in \mathcal{V} \, \}.
 \end{equation}

Moreover, the conductor element  $\mathbf{c}=(c_1, \ldots, c_r)$ has the following description from the point of view of the semigroup of values: it is the smallest lattice point satisfying $\mathbf{c}+(\mathbb{Z}_{\geq 0})^r \subset \mathcal{S}_{C, o}$. 

\bekezdes \textbf{The Alexander polynomial and the Poincar\'e series.} \label{par:AlexesHilbert} Another path to finding the Hilbert function could be using the following general formula, which will be our typical approach for plane curve singularities. This is really useful in the case $r \ge 2$, when no other algorithms are known (in the $r =1$ case the semigroup $\mathcal{S}_{C, 0}$ is usually known).

The embedded topological type of a   \emph{plane} curve germ can be completely characterized by several 
invariants, e.g., by its embedded resolution graph, or its splice diagram, or by 
the classical  multivariable  Alexander polynomial $\Delta (t_1, \ldots , t_r)$. For their description and properties see e.g., \cite{BrKn}.
For concrete formula for $\Delta({\bf t}) $ in terms of the splice diagram see \cite{EN}. 

In the next formula, which determines the Hilbert series $$H_{C,0}({\bf t})=H_{C,0}(t_1, \ldots ,t_r)=
\sum_{\ell \in \bN^r} \mathfrak{h}(\ell) {t_1}^{\ell_1} \cdots {t_r}^{\ell_r}$$ of an $r$-component plane curve $(C, 0) \subset (\mathbb{C}^2, 0)$,  it is convenient to replace the Alexander polynomial by the so-called Poincar\'e series. Namely, we set 
 $P(t):= \Delta(t)/(1-t)$ if $r = 1$ and $P({\bf t}) : = \Delta({\bf t})$ for $r > 1$.
 (For  more see,  e.g., \cite{cdg2,cdg}.) Moreover, we  can consider these invariants for any 
   subgerms $(C_J,0) =\cup_{j\in J} (C_j, 0)$, where $J \subseteq \{1, \ldots, r\}$. 

\begin{thm}[{\cite[Theorem 3.4.3]{GorNem2015}}; see also  {\cite[Corollary 4.3]{julioproj}}] \label{thm:reconst}
    With the above notations, $$H_{C, 0}({\bf t}) = \frac{1}{\prod_{i=1}^{r}(1-t_i)} \sum_{\emptyset \ne J = \{i_1, \ldots, i_{|J|}\} \subseteq \{1, \ldots, r\}} (-1)^{|J|-1} t_{i_1} \cdots t_{i_{|J|}} P_{C_J, 0}(t_{i_1},\ldots,t_{i_{|J|}}).$$
[In fact, this formula is also true for not necessarily plane curves as well, however, in the general case the Poincar\'e series does not have such a geometric description.]
\end{thm}

In particular,  the Hilbert function of a reduced \textit{plane} curve germ is an invariant of the embedded topological type, and hence the weight function $w_{an,0}$ is as well. For completeness, we note for plane germs the multivariable Alexander polynomial  is  already  a \emph{complete} invariant for the embedded topological type (cf. \cite{Yamamoto}); in particular, $H_{C,0}({\bf t})$ can be determined merely by $\Delta_{C,0}({\bf t})$, although this correspondence does not 
provide an explicit formula.

\begin{example} \label{ex:hilbert}
    If $r = 1$, then (writing $t$ for $t_1$) we have $H(t)= t\Delta(t)/(1-t)^2$. If $r=2$, then $$H_{C,0}(t_1, t_2) = \frac{1}{(1-t_1)(1-t_2)} \Big(
    \frac{t_1 \Delta_{C_1,0}(t_1)}{1-t_1} + \frac{t_2 \Delta_{C_2,0}(t_2)}{1-t_2} 
    - t_1 t_2 \Delta_{C,0}(t_1, t_2)\Big).$$ 
   Hence, if we write $\Delta_{C,0}(t_1,t_2)=
   \sum_{p=(p_1, p_2) \in {\rm supp} \Delta_{C,0} } {t_1}^{p_1}{t_2}^{p_2}$, 
   then for any $\ell \in (\bZ_{\geq 0})^2$
$$\mathfrak{h}_{C,0}(\ell) = \mathfrak{h}_{C_1,0}(\ell_1) + 
\mathfrak{h}_{C_2,0}(\ell_2) - \#\{p \in {\rm supp} \Delta_{C,0} \,:\, \ell \ge p + (1,1)\}, \ \mbox{and}$$
   $$w_{C,0}(\ell) = w_{C_1,0}(\ell_1) + w_{C_2,0}(\ell_2) -
   2 \cdot \#\{p \in {\rm supp} \Delta_{C,0}\,:\, \ell \ge p + (1,1)\}.$$
 \end{example}

\begin{example}\label{ex:H1}
If $r=1$ then $\bH_{an,\geq 1}(C,o)=0$, cf. paragraph \ref{bek:homdeg}.  On the other hand, if $r=2$, then $\bH_{an,1}(C,o)$ can be nonzero. The last formula from Example \ref{ex:hilbert} suggests the position of certain generators in the case of plane germs. Indeed, if $s_i\in\mathcal{S}_{C_i,0}$, but  $s_i+1\not\in\mathcal{S}_{C_i,0}$  then 
$s_i+1$ is a local maximum of the weight function of $(C_i,0)$: $w_{C_i, 0}(s_i)=w_{C_i,0}(s_i+1)-1=w_{C_i,0}(s_i+2)$ ($i=1,2$).
Therefore, in the context of $(C,0)=(C_1,0)\cup (C_2,0) \subset (\mathbb{C}^2, 0)$, $(s_1+1,s_2+1)$ is a local maximum of the association
$(\ell_1,\ell_2)\mapsto w_{C_1,0}(\ell_1)+ w_{C_2,0}(\ell_2)$, that is, this function is larger at  $(s_1+1,s_2+1)$
than at any other lattice point of $R((s_1,s_2), (s_1+2,s_2+2))$.
Hence, such lattice points give candidates for generators of $\bH_{an, 1}$. \end{example}

\begin{example} \label{ex:locmin}
Similarly, the local minimum points of the analytic weight function $w_{C, o}$, which correspond to certain generators of $\mathbb{H}_{an, 0}$, can be completely characterized by the semigroup of values, see \cite[Section 3]{KNS2}. For another formulation, let us consider the `semigroup module'
\begin{center}
    $\mathcal{M}_{C,o} := \big\{ \, \big(\mathfrak{v}^\omega_1(m), \dots, \mathfrak{v}^\omega_r(m)\big)  \in (\bZ \cup \{\infty\})^r \ : \ m \in M  \, \big\} \cap \mathbb{Z}^r$.
\end{center}
By the definition of extended valuations (cf. Definition \ref{def:edv}), it satisfies that $\mathcal{S}_{C, o}+\mathcal{M}_{C, o} \subset \mathcal{M}_{C, o}$. Moreover, it contains equivalent information to the $\mathfrak{h}^\circ$ height function:
\begin{equation}
    \mathcal{M}_{C,o} = \big\{ \, \ell  \in \bZ^r \ : \ \mathfrak{h}^\circ(-\ell-e_v) > \mathfrak{h}^\circ(-\ell) \text{ for all } v \in \mathcal{V}  \, \big\} 
\end{equation}
and one can compute the decrements of the height function $\mathfrak{h}^\circ$ along the edges of the lattice from $\mathcal{M}_{C,o}$ similarly to the case of the semigroup of values and the Hilbert function (described in \cite[subsection 2.2]{KNS2}). Using this language, a lattice point $\ell \in \mathbb{Z}^r$ is a local minimum point of the analytic weight function $w_{an, 0}$ (and thus gives a one-point connected component of $S_{w_{an, 0}(\ell)}$, cf. \cite[Section 3]{KNS2}) if and only if $ \ell \in \mathcal{S}_{C, o}$ and $-\ell \in \mathcal{M}_{C, o}$.

In the Gorenstein case $\mathcal{M}_{C, o}=(\mathfrak{v}_1^\omega(\alpha_{Gor}), \ldots, \mathfrak{v}_r^{\omega}(\alpha_{Gor}))+\mathcal{S}_{C, o}=-\mathbf{c}+S_{C, o}$, hence the local minimum points of the analytic weight function are the lattice points of the intersection $\mathcal{S}_{C, o} \cap \left(\mathbf{c}-\mathcal{S}_{C, o}\right)$.
\end{example}

\begin{example} \label{ex:2552} \cite[Example 4.5.2]{AgNeCurves} Let us consider the plane (in particular Gorenstein) curve singularity $(C, 0)= (\{ (x^2 - y^5)(x^5-y^2)=0\}, 0) \subset (\mathbb{C}^2, 0)$.
    One can compute the analytic Hilbert function $\mathfrak{h}$ with the help of the Alexander polynomial (c.f. paragraph \ref{par:AlexesHilbert}) or obtain the semigroup of values $\mathcal{S}_{C, 0}$ by direct calculation. The conductor element is $\mathbf{c}=(8,8)$.

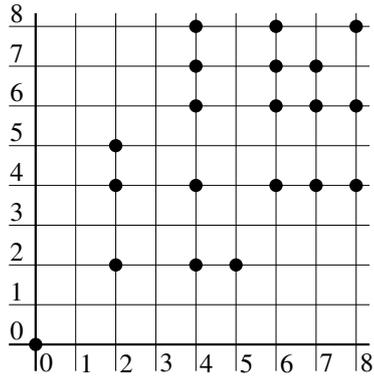
\begin{picture}(400,165)(-250,-40)

\put(-8,-1){\line(1,0){140}}
\put(-1,-8){\line(0,1){117}}

\put(5,6){\makebox(0,0){{$0$}}}
\put(5,18){\makebox(0,0){{$1$}}}
\put(5,30){\makebox(0,0){{$1$}}}
\put(5,42){\makebox(0,0){{$2$}}}
\put(5,54){\makebox(0,0){{$2$}}}
\put(5,66){\makebox(0,0){{$3$}}}
\put(5,78){\makebox(0,0){{$4$}}}
\put(5,90){\makebox(0,0){{$5$}}}
\put(5,102){\makebox(0,0){{$6$}}}

\put(20,6){\makebox(0,0){{$1$}}}
\put(20,18){\makebox(0,0){{$1$}}}
\put(20,30){\makebox(0,0){{$1$}}}
\put(20,42){\makebox(0,0){{$2$}}}
\put(20,54){\makebox(0,0){{$2$}}}
\put(20,66){\makebox(0,0){{$3$}}}
\put(20,78){\makebox(0,0){{$4$}}}
\put(20,90){\makebox(0,0){{$5$}}}
\put(20,102){\makebox(0,0){{$6$}}}

\put(35,6){\makebox(0,0){{$1$}}}
\put(35,18){\makebox(0,0){{$1$}}}
\put(35,30){\makebox(0,0){{$1$}}}
\put(35,42){\makebox(0,0){{$2$}}}
\put(35,54){\makebox(0,0){{$2$}}}
\put(35,66){\makebox(0,0){{$3$}}}
\put(35,78){\makebox(0,0){{$4$}}}
\put(35,90){\makebox(0,0){{$5$}}}
\put(35,102){\makebox(0,0){{$6$}}}

\put(50,6){\makebox(0,0){{$2$}}}
\put(50,18){\makebox(0,0){{$2$}}}
\put(50,30){\makebox(0,0){{$2$}}}
\put(50,42){\makebox(0,0){{$3$}}}
\put(50,54){\makebox(0,0){{$3$}}}
\put(50,66){\makebox(0,0){{$4$}}}
\put(50,78){\makebox(0,0){{$4$}}}
\put(50,90){\makebox(0,0){{$5$}}}
\put(50,102){\makebox(0,0){{$6$}}}

\put(65,6){\makebox(0,0){{$2$}}}
\put(65,18){\makebox(0,0){{$2$}}}
\put(65,30){\makebox(0,0){{$2$}}}
\put(65,42){\makebox(0,0){{$3$}}}
\put(65,54){\makebox(0,0){{$3$}}}
\put(65,66){\makebox(0,0){{$4$}}}
\put(65,78){\makebox(0,0){{$4$}}}
\put(65,90){\makebox(0,0){{$5$}}}
\put(65,102){\makebox(0,0){{$6$}}}

\put(80,6){\makebox(0,0){{$3$}}}
\put(80,18){\makebox(0,0){{$3$}}}
\put(80,30){\makebox(0,0){{$3$}}}
\put(80,42){\makebox(0,0){{$4$}}}
\put(80,54){\makebox(0,0){{$4$}}}
\put(80,66){\makebox(0,0){{$5$}}}
\put(80,78){\makebox(0,0){{$5$}}}
\put(80,90){\makebox(0,0){{$6$}}}
\put(80,102){\makebox(0,0){{$7$}}}

\put(95,6){\makebox(0,0){{$4$}}}
\put(95,18){\makebox(0,0){{$4$}}}
\put(95,30){\makebox(0,0){{$4$}}}
\put(95,42){\makebox(0,0){{$4$}}}
\put(95,54){\makebox(0,0){{$4$}}}
\put(95,66){\makebox(0,0){{$5$}}}
\put(95,78){\makebox(0,0){{$5$}}}
\put(95,90){\makebox(0,0){{$6$}}}
\put(95,102){\makebox(0,0){{$7$}}}

\put(110,6){\makebox(0,0){{$5$}}}
\put(110,18){\makebox(0,0){{$5$}}}
\put(110,30){\makebox(0,0){{$5$}}}
\put(110,42){\makebox(0,0){{$5$}}}
\put(110,54){\makebox(0,0){{$5$}}}
\put(110,66){\makebox(0,0){{$6$}}}
\put(110,78){\makebox(0,0){{$6$}}}
\put(110,90){\makebox(0,0){{$7$}}}
\put(110,102){\makebox(0,0){{$8$}}}

\put(125,6){\makebox(0,0){{$6$}}}
\put(125,18){\makebox(0,0){{$6$}}}
\put(125,30){\makebox(0,0){{$6$}}}
\put(125,42){\makebox(0,0){{$6$}}}
\put(125,54){\makebox(0,0){{$6$}}}
\put(125,66){\makebox(0,0){{$7$}}}
\put(125,78){\makebox(0,0){{$7$}}}
\put(125,90){\makebox(0,0){{$8$}}}
\put(125,102){\makebox(0,0){{$8$}}}

\linethickness{0.3mm}
\put(-210,-5){\line(1,0){135}}
\put(-200,-15){\line(0,1){135}}

\linethickness{0.05mm}
  \multiput(-210,10)(0,15){8}{\line(1,0){135}}
  \multiput(-185,-15)(15,0){8}{\line(0,1){135}}

\put(-200,-5){\circle*{5}}
\put(-170,25){\circle*{5}}
\put(-140,25){\circle*{5}}
\put(-170,55){\circle*{5}}
\put(-140,55){\circle*{5}}
\put(-170,70){\circle*{5}}
\put(-125,25){\circle*{5}}
\put(-110,55){\circle*{5}}
\put(-140,85){\circle*{5}}
\put(-140,100){\circle*{5}}
\put(-140,115){\circle*{5}}
\put(-95,55){\circle*{5}}
\put(-80,55){\circle*{5}}
\put(-110,85){\circle*{5}}
\put(-95,85){\circle*{5}}
\put(-80,85){\circle*{5}}
\put(-110,100){\circle*{5}}
\put(-110,115){\circle*{5}}
\put(-95,100){\circle*{5}}
\put(-80,115){\circle*{5}}

\put(-50,-33){\makebox(0,0){The semigroup of values $\mathcal{S}_{C, 0}$ and the Hilbert function $\mathfrak{h}$}}

\put(-196,-12){\makebox(0,0){\small{$0$}}}
\put(-181,-12){\makebox(0,0){\small{$1$}}}
\put(-166,-12){\makebox(0,0){\small{$2$}}}
\put(-151,-12){\makebox(0,0){\small{$3$}}}
\put(-136,-12){\makebox(0,0){\small{$4$}}}
\put(-121,-12){\makebox(0,0){\small{$5$}}}
\put(-106,-12){\makebox(0,0){\small{$6$}}}
\put(-91, -12){\makebox(0,0){\small{$7$}}}
\put(-76, -12){\makebox(0,0){\small{$8$}}}

\put(-207,0){\makebox(0,0){\small{$0$}}}
\put(-207,15){\makebox(0,0){\small{$1$}}}
\put(-207,30){\makebox(0,0){\small{$2$}}}
\put(-207,45){\makebox(0,0){\small{$3$}}}
\put(-207,60){\makebox(0,0){\small{$4$}}}
\put(-207,75){\makebox(0,0){\small{$5$}}}
\put(-207,90){\makebox(0,0){\small{$6$}}}
\put(-207,105){\makebox(0,0){\small{$7$}}}
\put(-207,120){\makebox(0,0){\small{$8$}}}

\end{picture}

Since the singularity is Gorenstein the `semigroup module' $\mathcal{M}_{C, 0}$ is $-\mathbf{c}+\mathcal{S}_{C, 0}$ and the analytic weight function is $w_{an, 0}=\mathfrak{h} + \mathfrak{h}^{sym}_\mathbf{c}-\mathfrak{h}(\mathbf{c})$. We thus obtain the following weight table and the corresponding graded root:

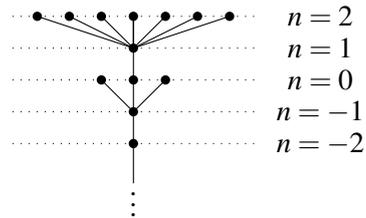
\begin{picture}(400,130)(0,-25)

\qbezier[30](255,44)(300,44)(345,44)
\qbezier[30](255,56)(300,56)(345,56)
\qbezier[30](255,68)(300,68)(345,68)
\qbezier[30](255,32)(300,32)(345,32)
\qbezier[30](255,20)(300,20)(345,20)

\put(370,44){\makebox(0,0){{$n=0$}}} \put(370,56){\makebox(0,0){{$n=1$}}}
\put(370,68){\makebox(0,0){{$n=2$}}}
\put(370,32){\makebox(0,0){{$n=-1$}}}
\put(370,20){\makebox(0,0){{$n=-2$}}}

\put(300,0){\makebox(0,0){$\vdots$}} 
\put(300,20){\circle*{3.5}}
\put(300,32){\circle*{3.5}} 
\put(288,44){\circle*{3.5}}
\put(312,44){\circle*{3.5}} 
\put(300,44){\circle*{3.5}}
\put(300,56){\circle*{3.5}} 
\put(300,5){\line(0,1){63}}
\put(288,44){\line(1,-1){12}} 
\put(300,32){\line(1,1){12}}

\put(264,68){\circle*{3.5}} 
\put(276,68){\circle*{3.5}}
\put(288,68){\circle*{3.5}} 
\put(300,68){\circle*{3.5}}
\put(312,68){\circle*{3.5}} 
\put(324,68){\circle*{3.5}}
\put(336,68){\circle*{3.5}}
\put(300,56){\line(0,1){12}}
\put(300,56){\line(1,1){12}} 
\put(300,56){\line(-1,1){12}}
\put(300,56){\line(2,1){24}} 
\put(300,56){\line(-2,1){24}}
\put(300,56){\line(3,1){36}} 
\put(300,56){\line(-3,1){36}}

\put(24,-18){\line(1,0){160}}
\put(32,-26){\line(0,1){115}}

\put(40,-11){\makebox(0,0){{$0$}}}
\put(40,1){\makebox(0,0){{$1$}}}
\put(40,13){\makebox(0,0){{$0$}}}
\put(40,25){\makebox(0,0){{$1$}}}
\put(40,37){\makebox(0,0){{$0$}}}
\put(40,49){\makebox(0,0){{$1$}}}
\put(40,61){\makebox(0,0){{$2$}}}
\put(40,73){\makebox(0,0){{$3$}}}
\put(40,85){\makebox(0,0){{$4$}}}

\put(57,-11){\makebox(0,0){{$1$}}}
\put(57,1){\makebox(0,0){{$0$}}}
\put(57,13){\makebox(0,0){{$-1$}}}
\put(57,25){\makebox(0,0){{$0$}}}
\put(57,37){\makebox(0,0){{$-1$}}}
\put(57,49){\makebox(0,0){{$0$}}}
\put(57,61){\makebox(0,0){{$1$}}}
\put(57,73){\makebox(0,0){{$2$}}}
\put(57,85){\makebox(0,0){{$3$}}}

\put(74,-11){\makebox(0,0){{$0$}}}
\put(74,1){\makebox(0,0){{$-1$}}}
\put(74,13){\makebox(0,0){{$-2$}}}
\put(74,25){\makebox(0,0){{$-1$}}}
\put(74,37){\makebox(0,0){{$-2$}}}
\put(74,49){\makebox(0,0){{$-1$}}}
\put(74,61){\makebox(0,0){{$0$}}}
\put(74,73){\makebox(0,0){{$1$}}}
\put(74,85){\makebox(0,0){{$2$}}}

\put(91,-11){\makebox(0,0){{$1$}}}
\put(91,1){\makebox(0,0){{$0$}}}
\put(91,13){\makebox(0,0){{$-1$}}}
\put(91,25){\makebox(0,0){{$0$}}}
\put(91,37){\makebox(0,0){{$-1$}}}
\put(91,49){\makebox(0,0){{$0$}}}
\put(91,61){\makebox(0,0){{$-1$}}}
\put(91,73){\makebox(0,0){{$0$}}}
\put(91,85){\makebox(0,0){{$1$}}}

\put(108,-11){\makebox(0,0){{$0$}}}
\put(108,1){\makebox(0,0){{$-1$}}}
\put(108,13){\makebox(0,0){{$-2$}}}
\put(108,25){\makebox(0,0){{$-1$}}}
\put(108,37){\makebox(0,0){{$-2$}}}
\put(108,49){\makebox(0,0){{$-1$}}}
\put(108,61){\makebox(0,0){{$-2$}}}
\put(108,73){\makebox(0,0){{$-1$}}}
\put(108,85){\makebox(0,0){{$0$}}}

\put(125,-11){\makebox(0,0){{$1$}}}
\put(125,1){\makebox(0,0){{$0$}}}
\put(125,13){\makebox(0,0){{$-1$}}}
\put(125,25){\makebox(0,0){{$0$}}}
\put(125,37){\makebox(0,0){{$-1$}}}
\put(125,49){\makebox(0,0){{$0$}}}
\put(125,61){\makebox(0,0){{$-1$}}}
\put(125,73){\makebox(0,0){{$0$}}}
\put(125,85){\makebox(0,0){{$1$}}}

\put(142,-11){\makebox(0,0){{$2$}}}
\put(142,1){\makebox(0,0){{$1$}}}
\put(142,13){\makebox(0,0){{$0$}}}
\put(142,25){\makebox(0,0){{$-1$}}}
\put(142,37){\makebox(0,0){{$-2$}}}
\put(142,49){\makebox(0,0){{$-1$}}}
\put(142,61){\makebox(0,0){{$-2$}}}
\put(142,73){\makebox(0,0){{$-1$}}}
\put(142,85){\makebox(0,0){{$0$}}}

\put(159,-11){\makebox(0,0){{$3$}}}
\put(159,1){\makebox(0,0){{$2$}}}
\put(159,13){\makebox(0,0){{$1$}}}
\put(159,25){\makebox(0,0){{$0$}}}
\put(159,37){\makebox(0,0){{$-1$}}}
\put(159,49){\makebox(0,0){{$0$}}}
\put(159,61){\makebox(0,0){{$-1$}}}
\put(159,73){\makebox(0,0){{$0$}}}
\put(159,85){\makebox(0,0){{$1$}}}

\put(176,-11){\makebox(0,0){{$4$}}}
\put(176,1){\makebox(0,0){{$3$}}}
\put(176,13){\makebox(0,0){{$2$}}}
\put(176,25){\makebox(0,0){{$1$}}}
\put(176,37){\makebox(0,0){{$0$}}}
\put(176,49){\makebox(0,0){{$1$}}}
\put(176,61){\makebox(0,0){{$0$}}}
\put(176,73){\makebox(0,0){{$1$}}}
\put(176,85){\makebox(0,0){{$0$}}}

\put(125,49){\circle{12}}
\put(91,25){\circle{12}}
\put(66,8){\framebox(16,12){}}
\put(100,8){\framebox(16,12){}}
\put(66,32){\framebox(16,12){}}
\put(100,32){\framebox(16,12){}}
\put(100,56){\framebox(16,12){}}
\put(134,32){\framebox(16,12){}}
\put(134,56){\framebox(16,12){}}

\put(33,-17){\framebox(14,12){}}
\put(169,79){\framebox(14,12){}}
\end{picture}

\begin{center}
$\min w_0=-2, \ {\rm rk}\mathbb{H}_{an,red, 0}=8, \ {\rm rk}\mathbb{H}_{an, 1}=2$, thus $eu(\mathbb{H}_{an, *})=8=\delta(C, 0).$
\end{center}

 On the diagram above the circles denote the generators of $\mathbb{H}_{an,1}$ (corresponding to the `circles' around the lattice points $(2+1, 2+1)$ and $(4+1, 4+1)$ --- compare with Example \ref{ex:H1}), while the rectangles denote the local minimum points of the weight function (these correspond to the semigroup elements $s \in \mathcal{S}_{C, 0}$ with $\mathbf{c}-s \in \mathcal{S}_{C, 0}$).
\end{example}

\begin{example}
For the plane curve singularity $\{(x^3+y^5)(x^5+y^3)=0\}$, the semigroup of both components is $\langle 3,5\rangle=\{
1,3,5,6,\geq 8\}$ and the `circles' around $(4,4)$ and $(7,7)$  and their symmentrical points $(10,10)$ and 
$(13,13)$ with respect to ${\bf c}=(17,17) $ generate $\bH_{an, 1}$. 

In general, for $\{(x^2+y^{2k+1})(x^{2l+1}+y^2)=0\}$ ($k,l\geq 2$) the rank of $\bH_{an, 1}$ is $kl-2$.     
\end{example}

\begin{example}\label{ex:345534}
    Let us consider the space curve germ $(C, 0) \subset (\mathbb{C}^3, 0)$ consisting of two irreducible components parametrized by $C_1:t \mapsto (t^3, t^4, t^5)$ and $C_2:s \mapsto (s^5, s^3, s^4)$. One can compute the semigroup $\mathcal{S}_{C, 0}$ (or, equivalently, the Hilbert function $\mathfrak{h}$) and obtain the conductor element $\mathbf{c}=(6,6)$ (e.g., using Lemma 2.4.2 of \cite{KNS2}), therefore $\delta(C, 0)=|\mathbf{c}|-\mathfrak{h}(\mathbf{c})=8$. 
    \begin{center}
        \begin{picture}(360,143)(0, -20)
    
    \linethickness{0.3mm}
\put(10,20){\line(1,0){105}}
\put(20,10){\line(0,1){105}}

\linethickness{0.05mm}
  \multiput(10,35)(0,15){6}{\line(1,0){105}}
  \multiput(35,10)(15,0){6}{\line(0,1){105}}

\put(20,20){\circle*{5}}
\put(65,65){\circle*{5}}
\put(65,80){\circle*{5}}
\put(80,65){\circle*{5}}
\put(65,95){\circle*{5}}
\put(95,80){\circle*{5}}
\put(110,110){\circle*{5}}

\put(180,-11){\makebox(0,0){The semigroup of values $\mathcal{S}_{C, 0}$ and the Hilbert function $\mathfrak{h}$}}

\put(24,13){\makebox(0,0){\small{$0$}}}
\put(39,13){\makebox(0,0){\small{$1$}}}
\put(54,13){\makebox(0,0){\small{$2$}}}
\put(69,13){\makebox(0,0){\small{$3$}}}
\put(84,13){\makebox(0,0){\small{$4$}}}
\put(99,13){\makebox(0,0){\small{$5$}}}
\put(114,13){\makebox(0,0){\small{$6$}}}

\put(13,25){\makebox(0,0){\small{$0$}}}
\put(13,40){\makebox(0,0){\small{$1$}}}
\put(13,55){\makebox(0,0){\small{$2$}}}
\put(13,70){\makebox(0,0){\small{$3$}}}
\put(13,85){\makebox(0,0){\small{$4$}}}
\put(13,100){\makebox(0,0){\small{$5$}}}
\put(13,115){\makebox(0,0){\small{$6$}}}

\put(208,23){\line(1,0){107}}
\put(214,17){\line(0,1){90}}

\put(220,30){\makebox(0,0){$0$}}
\put(235,30){\makebox(0,0){$1$}}
\put(250,30){\makebox(0,0){$1$}}
\put(265,30){\makebox(0,0){$1$}}
\put(280,30){\makebox(0,0){$2$}}
\put(295,30){\makebox(0,0){$3$}}
\put(310,30){\makebox(0,0){$4$}}

\put(220,42){\makebox(0,0){$1$}}
\put(235,42){\makebox(0,0){$1$}}
\put(250,42){\makebox(0,0){$1$}}
\put(265,42){\makebox(0,0){$1$}}
\put(280,42){\makebox(0,0){$2$}}
\put(295,42){\makebox(0,0){$3$}}
\put(310,42){\makebox(0,0){$4$}}

\put(220,54){\makebox(0,0){$1$}}
\put(235,54){\makebox(0,0){$1$}}
\put(250,54){\makebox(0,0){$1$}}
\put(265,54){\makebox(0,0){$1$}}
\put(280,54){\makebox(0,0){$2$}}
\put(295,54){\makebox(0,0){$3$}}
\put(310,54){\makebox(0,0){$4$}}

\put(220,66){\makebox(0,0){$1$}}
\put(235,66){\makebox(0,0){$1$}}
\put(250,66){\makebox(0,0){$1$}}
\put(265,66){\makebox(0,0){$1$}}
\put(280,66){\makebox(0,0){$2$}}
\put(295,66){\makebox(0,0){$3$}}
\put(310,66){\makebox(0,0){$4$}}

\put(220,78){\makebox(0,0){$2$}}
\put(235,78){\makebox(0,0){$2$}}
\put(250,78){\makebox(0,0){$2$}}
\put(265,78){\makebox(0,0){$2$}}
\put(280,78){\makebox(0,0){$3$}}
\put(295,78){\makebox(0,0){$3$}}
\put(310,78){\makebox(0,0){$4$}}

\put(220,90){\makebox(0,0){$3$}}
\put(235,90){\makebox(0,0){$3$}}
\put(250,90){\makebox(0,0){$3$}}
\put(265,90){\makebox(0,0){$3$}}
\put(280,90){\makebox(0,0){$4$}}
\put(295,90){\makebox(0,0){$4$}}
\put(310,90){\makebox(0,0){$4$}}

\put(220,102){\makebox(0,0){$4$}}
\put(235,102){\makebox(0,0){$4$}}
\put(250,102){\makebox(0,0){$4$}}
\put(265,102){\makebox(0,0){$4$}}
\put(280,102){\makebox(0,0){$4$}}
\put(295,102){\makebox(0,0){$4$}}
\put(310,102){\makebox(0,0){$4$}}
    \end{picture}
    \end{center}
    Then $\mathfrak{h}^\circ=8+ \mathfrak{h}-|\,.\,|$ and the `semigroup module' $\mathcal{M}_{C,0}$ associated with the extended valuations are as follows:
    \begin{center}
        \begin{picture}(330,143)(0, -10)
    
    \linethickness{0.3mm}
\put(15,110){\line(1,0){105}}
\put(110,15){\line(0,1){105}}

\linethickness{0.05mm}
  \multiput(15,20)(0,15){6}{\line(1,0){105}}
  \multiput(20,15)(15,0){6}{\line(0,1){105}}

\put(20,20){\circle*{5}}
\put(65,65){\circle*{5}}
\put(65,80){\circle*{5}}
\put(80,65){\circle*{5}}
\put(65,95){\circle*{5}}
\put(95,80){\circle*{5}}
\put(110,110){\circle*{5}}
\put(35,20){\circle*{5}}
\put(50,20){\circle*{5}}
\put(20,35){\circle*{5}}
\put(35,50){\circle*{5}}
\put(80,95){\circle*{5}}
\put(80,80){\circle*{5}}
\put(95,65){\circle*{5}}
\put(110,65){\circle*{5}}
\put(65,110){\circle*{5}}
\put(80,110){\circle*{5}}
\put(95,95){\circle*{5}}
\put(110,80){\circle*{5}}

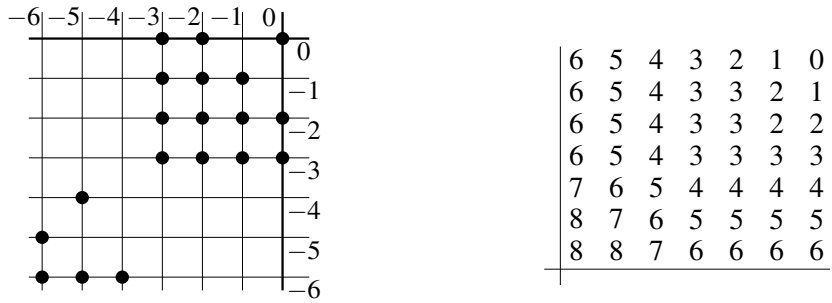
\put(180,-4){\makebox(0,0){The semigroup module $\mathcal{M}_{C, 0}$ and the height function $\mathfrak{h}^\circ$}}

\put(105,118){\makebox(0,0){\small{$0$}}}
\put(89,118){\makebox(0,0){\small{$-1$}}}
\put(73,118){\makebox(0,0){\small{$-2$}}}
\put(58,118){\makebox(0,0){\small{$-3$}}}
\put(43,118){\makebox(0,0){\small{$-4$}}}
\put(28,118){\makebox(0,0){\small{$-5$}}}
\put(13,118){\makebox(0,0){\small{$-6$}}}

\put(118,105){\makebox(0,0){\small{$0$}}}
\put(118,90){\makebox(0,0){\small{$-1$}}}
\put(118,75){\makebox(0,0){\small{$-2$}}}
\put(118,60){\makebox(0,0){\small{$-3$}}}
\put(118,45){\makebox(0,0){\small{$-4$}}}
\put(118,30){\makebox(0,0){\small{$-5$}}}
\put(118,15){\makebox(0,0){\small{$-6$}}}

\put(208,23){\line(1,0){107}}
\put(214,17){\line(0,1){90}}

\put(220,30){\makebox(0,0){$8$}}
\put(235,30){\makebox(0,0){$8$}}
\put(250,30){\makebox(0,0){$7$}}
\put(265,30){\makebox(0,0){$6$}}
\put(280,30){\makebox(0,0){$6$}}
\put(295,30){\makebox(0,0){$6$}}
\put(310,30){\makebox(0,0){$6$}}

\put(220,42){\makebox(0,0){$8$}}
\put(235,42){\makebox(0,0){$7$}}
\put(250,42){\makebox(0,0){$6$}}
\put(265,42){\makebox(0,0){$5$}}
\put(280,42){\makebox(0,0){$5$}}
\put(295,42){\makebox(0,0){$5$}}
\put(310,42){\makebox(0,0){$5$}}

\put(220,54){\makebox(0,0){$7$}}
\put(235,54){\makebox(0,0){$6$}}
\put(250,54){\makebox(0,0){$5$}}
\put(265,54){\makebox(0,0){$4$}}
\put(280,54){\makebox(0,0){$4$}}
\put(295,54){\makebox(0,0){$4$}}
\put(310,54){\makebox(0,0){$4$}}

\put(220,66){\makebox(0,0){$6$}}
\put(235,66){\makebox(0,0){$5$}}
\put(250,66){\makebox(0,0){$4$}}
\put(265,66){\makebox(0,0){$3$}}
\put(280,66){\makebox(0,0){$3$}}
\put(295,66){\makebox(0,0){$3$}}
\put(310,66){\makebox(0,0){$3$}}

\put(220,78){\makebox(0,0){$6$}}
\put(235,78){\makebox(0,0){$5$}}
\put(250,78){\makebox(0,0){$4$}}
\put(265,78){\makebox(0,0){$3$}}
\put(280,78){\makebox(0,0){$3$}}
\put(295,78){\makebox(0,0){$2$}}
\put(310,78){\makebox(0,0){$2$}}

\put(220,90){\makebox(0,0){$6$}}
\put(235,90){\makebox(0,0){$5$}}
\put(250,90){\makebox(0,0){$4$}}
\put(265,90){\makebox(0,0){$3$}}
\put(280,90){\makebox(0,0){$3$}}
\put(295,90){\makebox(0,0){$2$}}
\put(310,90){\makebox(0,0){$1$}}

\put(220,102){\makebox(0,0){$6$}}
\put(235,102){\makebox(0,0){$5$}}
\put(250,102){\makebox(0,0){$4$}}
\put(265,102){\makebox(0,0){$3$}}
\put(280,102){\makebox(0,0){$2$}}
\put(295,102){\makebox(0,0){$1$}}
\put(310,102){\makebox(0,0){$0$}}
    \end{picture}
    \end{center}    
    We can clearly see, that $\mathfrak{h}^\circ$ is not the symmetrized version of $\mathfrak{h}$, showing that this spacial curve $(C, 0)$ is not Gorenstein, equivalently,  $\omega^R_C$ is not free of rank one above $\mathcal{O}$ (compare also with \cite[Theorem 4.1.1]{KNS2}).
    
    We also highlight that, in this case, we can obtain concrete bases $\{f_1, \ldots, f_4\}$ and $\{m_1, \ldots, m_8\}$ of the $\mathbb{C}$-vector spaces $\mathcal{O}/\mathcal{C}$ and $\omega^R_C/\Omega^1_{\overline{C}}$ satisfying that
    \begin{align*}
        \forall f = \sum_{i=1}^{4} a_i f_i \in \mathcal{O}/\mathcal{C}: & \    \mathfrak{v}_1(f)=\min\{\mathfrak{v}_1(f_i)\,:\, a_i \neq 0\} \text{ and } \mathfrak{v}_2(f)=\min\{\mathfrak{v}_2(f_i)\,:\, a_i \neq 0\};\\
        \forall m = \sum_{j=1}^{8} b_j m_j \in \omega^{R}_C/\Omega^1_{\overline{C}}: & \  \mathfrak{v}^M_1(m)=\min\{\mathfrak{v}^M_1(m_j)\,:\, b_j \neq 0\} \text{ and } \mathfrak{v}^M_2(m)=\min\{\mathfrak{v}^M_2(m_j)\,:\, b_j \neq 0\}.
    \end{align*}
    For the significance and existence of these bases see 
    Remarks \ref{rem:constructionofq-val} and \ref{rem:constructionofextq-val}, in our concrete case one such choice is for example:
    
    \begin{align*}
    \mathcal{O}/\mathcal{C} \cong \mathcal{O}_{\mathbb{C}^3, 0}/\mathfrak{m}_{\mathbb{C}^3, 0}^2 =&\, \mathbb{C}\big\langle [1], [x], [y], [z]\big\rangle \text{ with }\\
    \mathfrak{v}([1])=(0,0), \hspace{10mm} \mathfrak{v}([x])=(3,5), \hspace{10mm} &\, \mathfrak{v}([y])=(4,3), \hspace{10mm} \mathfrak{v}([z])=(5,4), 
    \end{align*}
    whereas $\omega^R_C/\Omega^1_{\overline{C}}=$  
    \begin{sizeddisplay}{\small}
    \begin{align*} 
    \mathbb{C}\Bigg\langle\left[\dfrac{1}{t}dt-\dfrac{1}{s}ds\right], \left[\dfrac{1}{t^2}dt\right], \left[\dfrac{1}{s^2}ds\right], \left[\dfrac{1}{t^3}dt\right], \left[\dfrac{1}{s^3}ds\right], &\,  \left[\dfrac{1}{t^4}dt-\dfrac{1}{s^6}ds\right], \left[\dfrac{1}{t^5}dt-\dfrac{1}{s^4}ds\right],
    \left[\dfrac{1}{t^6}dt-\dfrac{1}{s^5}ds\right]\Bigg\rangle \\
    \mathfrak{v}^M\left( \left[\dfrac{1}{t}dt-\dfrac{1}{s}ds\right]\right)=(-1,-1), \hspace{5mm} \mathfrak{v}^M\left( \left[\dfrac{1}{t^2}dt\right]\right)&\,=(-2,0), \hspace{5mm} \mathfrak{v}^M\left( \left[\dfrac{1}{s^2}ds\right]\right)=(0,-2), 
    \end{align*}
    \end{sizeddisplay}
    \begin{align*} 
    \mathfrak{v}^M\left( \left[\dfrac{1}{t^3}dt\right]\right)=(-3,0), \hspace{5mm} \mathfrak{v}^M\left( \left[\dfrac{1}{s^3}ds\right]\right)=&\,(0,-3), \hspace{5mm} \mathfrak{v}^M\left( \left[\dfrac{1}{t^4}dt-\dfrac{1}{s^6}ds\right]\right)=(-4,-6), \\
    \mathfrak{v}^M\left( \left[\dfrac{1}{t^5}dt-\dfrac{1}{s^4}ds\right]\right)=(-5,-4), \hspace{5mm} &\, \mathfrak{v}^M\left( \left[\dfrac{1}{t^6}dt-\dfrac{1}{s^5}ds\right]\right)=(-6,-5).
    \end{align*}
    \begin{center}
        \begin{picture}(360,118)(0, -12)
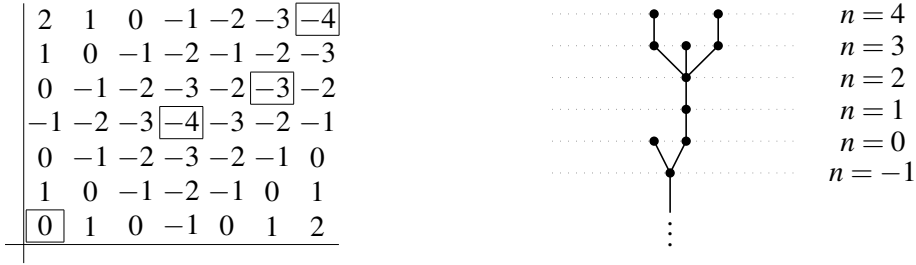
\put(180,-4){\makebox(0,0){The weight function $w_{an,0}$ and the graded root $\mathfrak{R}(C, 0)$}}
\linethickness{0.05mm}
\put(5,13){\line(1,0){125}}
\put(12,6){\line(0,1){98}}

\put(20,20){\makebox(0,0){$0$}}
\put(37,20){\makebox(0,0){$1$}}
\put(54,20){\makebox(0,0){$0$}}
\put(71,20){\makebox(0,0){$-1$}}
\put(88,20){\makebox(0,0){$0$}}
\put(105,20){\makebox(0,0){$1$}}
\put(122,20){\makebox(0,0){$2$}}

\put(20,33){\makebox(0,0){$1$}}
\put(37,33){\makebox(0,0){$0$}}
\put(54,33){\makebox(0,0){$-1$}}
\put(71,33){\makebox(0,0){$-2$}}
\put(88,33){\makebox(0,0){$-1$}}
\put(105,33){\makebox(0,0){$0$}}
\put(122,33){\makebox(0,0){$1$}}

\put(20,46){\makebox(0,0){$0$}}
\put(37,46){\makebox(0,0){$-1$}}
\put(54,46){\makebox(0,0){$-2$}}
\put(71,46){\makebox(0,0){$-3$}}
\put(88,46){\makebox(0,0){$-2$}}
\put(105,46){\makebox(0,0){$-1$}}
\put(122,46){\makebox(0,0){$0$}}

\put(20,59){\makebox(0,0){$-1$}}
\put(37,59){\makebox(0,0){$-2$}}
\put(54,59){\makebox(0,0){$-3$}}
\put(71,59){\makebox(0,0){$-4$}}
\put(88,59){\makebox(0,0){$-3$}}
\put(105,59){\makebox(0,0){$-2$}}
\put(122,59){\makebox(0,0){$-1$}}

\put(20,72){\makebox(0,0){$0$}}
\put(37,72){\makebox(0,0){$-1$}}
\put(54,72){\makebox(0,0){$-2$}}
\put(71,72){\makebox(0,0){$-3$}}
\put(88,72){\makebox(0,0){$-2$}}
\put(105,72){\makebox(0,0){$-3$}}
\put(122,72){\makebox(0,0){$-2$}}

\put(20,85){\makebox(0,0){$1$}}
\put(37,85){\makebox(0,0){$0$}}
\put(54,85){\makebox(0,0){$-1$}}
\put(71,85){\makebox(0,0){$-2$}}
\put(88,85){\makebox(0,0){$-1$}}
\put(105,85){\makebox(0,0){$-2$}}
\put(122,85){\makebox(0,0){$-3$}}

\put(20,98){\makebox(0,0){$2$}}
\put(37,98){\makebox(0,0){$1$}}
\put(54,98){\makebox(0,0){$0$}}
\put(71,98){\makebox(0,0){$-1$}}
\put(88,98){\makebox(0,0){$-2$}}
\put(105,98){\makebox(0,0){$-3$}}
\put(122,98){\makebox(0,0){$-4$}}

\put(13,14){\framebox(14,12)}
\put(63,53){\framebox(16,12)}
\put(97,66){\framebox(16,12)}
\put(114,92){\framebox(16,12)}

\qbezier[30](210,100)(255,100)(300,100)
\qbezier[30](210,88)(255,88)(300,88)
\qbezier[30](210,76)(255,76)(300,76)
\qbezier[30](210,40)(255,40)(300,40)
\qbezier[30](210,64)(255,64)(300,64)
\qbezier[30](210,52)(255,52)(300,52)

\put(330,100){\makebox(0,0){$n=4$}}
\put(330,88){\makebox(0,0){$n=3$}}
\put(330,76){\makebox(0,0){$n=2$}}
\put(330,64){\makebox(0,0){$n=1$}}
\put(330,52){\makebox(0,0){$n=0$}}
\put(330,40){\makebox(0,0){$n=-1$}}

\put(248,100){\circle*{3.5}}
\put(260,88){\circle*{3.5}}
\put(272,100){\circle*{3.5}}
\put(248,88){\circle*{3.5}}
\put(272,88){\circle*{3.5}}
\put(260,76){\circle*{3.5}}
\put(260,64){\circle*{3.5}}
\put(260,52){\circle*{3.5}}
\put(248,52){\circle*{3.5}}
\put(254,40){\circle*{3.5}}
\put(254,20){\makebox(0,0){$\vdots$}}

\linethickness{0.6pt}
\put(248,100){\line(0,-1){12}}
\put(272,100){\line(0,-1){12}}
\put(248,88){\line(1,-1){12}}
\put(260,88){\line(0,-1){36}}
\put(272,88){\line(-1,-1){12}}
\put(248,52){\line(1,-2){6}}
\put(260,52){\line(-1,-2){6}}
\put(254,40){\line(0, -1){15}}
    \end{picture}
    \end{center}

Clearly $eu(\bH_*(C,o))=\delta=8$. However, if in $R(0,{\bf c})$ we symmetrize the Hilbert function $\mathfrak{h}$  with respect to ${\bf c}=(6,6)$, 
then we get $\mathbb{S}\bH_* (\mathcal{C}\hookrightarrow \mathcal{O})$, with Euler characteristic 4.

    \begin{center}
        \begin{picture}(360,128)(0, -17)
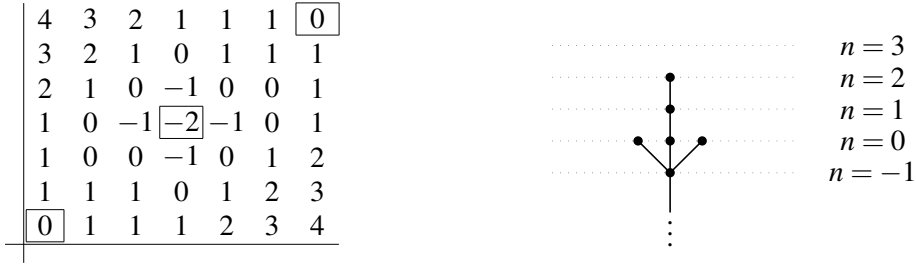
\put(180,-9){\makebox(0,0){The symmetrized weight function $\mathfrak{h} + \mathfrak{h}^{sym}_{\mathbf{c}}-\mathfrak{h}(\mathbf{c})$ and the corresponding graded root}}
\linethickness{0.05mm}
\put(5,13){\line(1,0){125}}
\put(12,6){\line(0,1){98}}

\put(20,20){\makebox(0,0){$0$}}
\put(37,20){\makebox(0,0){$1$}}
\put(54,20){\makebox(0,0){$1$}}
\put(71,20){\makebox(0,0){$1$}}
\put(88,20){\makebox(0,0){$2$}}
\put(105,20){\makebox(0,0){$3$}}
\put(122,20){\makebox(0,0){$4$}}

\put(20,33){\makebox(0,0){$1$}}
\put(37,33){\makebox(0,0){$1$}}
\put(54,33){\makebox(0,0){$1$}}
\put(71,33){\makebox(0,0){$0$}}
\put(88,33){\makebox(0,0){$1$}}
\put(105,33){\makebox(0,0){$2$}}
\put(122,33){\makebox(0,0){$3$}}

\put(20,46){\makebox(0,0){$1$}}
\put(37,46){\makebox(0,0){$0$}}
\put(54,46){\makebox(0,0){$0$}}
\put(71,46){\makebox(0,0){$-1$}}
\put(88,46){\makebox(0,0){$0$}}
\put(105,46){\makebox(0,0){$1$}}
\put(122,46){\makebox(0,0){$2$}}

\put(20,59){\makebox(0,0){$1$}}
\put(37,59){\makebox(0,0){$0$}}
\put(54,59){\makebox(0,0){$-1$}}
\put(71,59){\makebox(0,0){$-2$}}
\put(88,59){\makebox(0,0){$-1$}}
\put(105,59){\makebox(0,0){$0$}}
\put(122,59){\makebox(0,0){$1$}}

\put(20,72){\makebox(0,0){$2$}}
\put(37,72){\makebox(0,0){$1$}}
\put(54,72){\makebox(0,0){$0$}}
\put(71,72){\makebox(0,0){$-1$}}
\put(88,72){\makebox(0,0){$0$}}
\put(105,72){\makebox(0,0){$0$}}
\put(122,72){\makebox(0,0){$1$}}

\put(20,85){\makebox(0,0){$3$}}
\put(37,85){\makebox(0,0){$2$}}
\put(54,85){\makebox(0,0){$1$}}
\put(71,85){\makebox(0,0){$0$}}
\put(88,85){\makebox(0,0){$1$}}
\put(105,85){\makebox(0,0){$1$}}
\put(122,85){\makebox(0,0){$1$}}

\put(20,98){\makebox(0,0){$4$}}
\put(37,98){\makebox(0,0){$3$}}
\put(54,98){\makebox(0,0){$2$}}
\put(71,98){\makebox(0,0){$1$}}
\put(88,98){\makebox(0,0){$1$}}
\put(105,98){\makebox(0,0){$1$}}
\put(122,98){\makebox(0,0){$0$}}

\put(13,14){\framebox(14,12)}
\put(63,53){\framebox(16,12)}
\put(114,92){\framebox(16,12)}

\qbezier[30](210,88)(255,88)(300,88)
\qbezier[30](210,76)(255,76)(300,76)
\qbezier[30](210,40)(255,40)(300,40)
\qbezier[30](210,64)(255,64)(300,64)
\qbezier[30](210,52)(255,52)(300,52)

\put(330,88){\makebox(0,0){$n=3$}}
\put(330,76){\makebox(0,0){$n=2$}}
\put(330,64){\makebox(0,0){$n=1$}}
\put(330,52){\makebox(0,0){$n=0$}}
\put(330,40){\makebox(0,0){$n=-1$}}

\put(254,76){\circle*{3.5}}
\put(254,64){\circle*{3.5}}
\put(266,52){\circle*{3.5}}
\put(254,52){\circle*{3.5}}
\put(242,52){\circle*{3.5}}
\put(254,40){\circle*{3.5}}
\put(254,20){\makebox(0,0){{$\vdots$}}}

\linethickness{0.6pt}

\put(254,40){\line(1,1){12}}

\put(254,40){\line(-1,1){12}}

\put(254,25){\line(0, 1){51}}
    \end{picture}
    \end{center}
\end{example}

\section{Categorifications of analytic invariants of isolated singularities of $\dim\geq 2$}
\label{s:dnagy}

In this section we will present further examples of categorification via lattice homology of an invariant defined as the codimension of a realizable submodule: the geometric genus $p_g$ and irregularity $q$ of normal isolated complex analytic singularities of $\dim\geq 2$ and the plurigenera of normal surface singularities. Our construction in the case of $p_g$ is slightly more general than the analytic lattice (co)homology defined by T. \'Agoston and the first author in \cite{AgNe1, AgNeHigh} (where certain restrictions were imposed), 
the present case is  well-defined for any normal isolated singularity. On the other hand, the categorification of the irregularity is a new construction, whereas the case of plurigenera here is treated differently (with different output) than the authors' first attempt in \cite{NSplurig}.

\subsection{Geometric genus and irregularity of normal isolated singularities}\label{ss:iso}\,

Let $(X,o)$ be a complex analytic normal 
isolated singularity of dimension $n\geq 2$. By a \emph{resolution} of $(X, o)$ we mean the following: fix a sufficiently small Stein representative $X$ of $(X,o)$ and let $\phi:\widetilde{X}\to X$ be a holomorphic map, such that 
\begin{itemize}
    \item $\widetilde{X}$ is smooth,
    \item $\phi$  is an analytic isomorphism over the regular part $U:=X\setminus \{o\}$ and
    \item $E:=\phi^{-1}(o)_{red}$ is a divisor, with each irreducible component of dimension $n-1$ (when considered with reduced structure we call it the \emph{exceptional divisor}).
\end{itemize}  
Similarly to the surface case, we call a resolution \emph{good}, if $E$ is a simple normal crossing divisor, and \emph{minimal (good)} if any (good) resolution dominates it. Good resolutions always exist (at least, in characteristic $0$), indeed they correspond to the `strong resolutions' of Kollár (cf. \cite{Kollar}, where he also provides a proof of this fact). 

In \cite{Yau0, Yau1} S. S.-T. Yau defined the invariants $s^{(p)}=\dim_{\mathbb{C}} H^0(\widetilde{X}\setminus E, \Omega^{p}_{\widetilde{X}})\big/H^0(\widetilde{X}, \Omega^{p}_{\widetilde{X}})$, where $\Omega^{p}_{\widetilde{X}}$ denotes the sheaf of holomorphic $p$-forms on $\widetilde{X}$. We show its independence  from the resolution $\phi$ and the Stein representative $X$ following 
\cite[Proposition 4.3]{Yau0}.  
On one hand, let $\overline{\Omega}^{p}_X$ be the $0^\text{th}$ direct image sheaf $\phi_* \Omega^{p}_{\widetilde{X}}$ of the sheaf of holomorphic $p$-forms on $\widetilde{X}$. By \cite[Proposition (1.1)]{StrStee} this is independent of the resolution $\phi$, since its sections are the holomorphic $p$-forms of `first kind', i.e., holomorphic $p$-forms $\omega$ on $U$ such that for any resolution $\phi':\widetilde{X}' \rightarrow X$ the pullback $\phi'^*(\omega)$ extends holomorphically to the exceptional set $E'=\phi'^{-1}(o)$ (equivalently, for any $\mathcal{C}^{\infty}$ map $\gamma:\triangle^p \rightarrow X$ the integral $\int_\gamma w$ exists). On the other hand, let $\iota: U \hookrightarrow X$ be the inclusion map. Then, again, the $0^\text{th}$ direct image sheaf $\overline{\overline{\Omega}}^p_X=\iota_*\Omega^p_U$ is also coherent, hence the quotient sheaf $\overline{\overline{\Omega}}^p_X\big/\overline{\Omega}^p_X$ is cohe\-rent and supported on $o$. Moreover, since $H^0\big(X, \overline{\overline{\Omega}}^p_X\big)=H^0(U, \Omega^p_U) =H^0(\widetilde{X}\setminus E, \Omega^p_{\widetilde{X}})$ and $H^0\big(X, \overline{\Omega}^p_X\big)=H^0(\widetilde{X}, \Omega^p_{\widetilde{X}})$ we get the complex vector space isomorphism
\begin{equation}\label{eq:indepforirreg&p_g}
    \dfrac{H^0(\widetilde{X}\setminus E, \Omega^{p}_{\widetilde{X}})}{H^0(\widetilde{X}, \Omega^{p}_{\widetilde{X}})}\cong \dfrac{H^0\big(X, \overline{\overline{\Omega}}^p_X\big)}{H^0\big(X, \overline{\Omega}^p_X\big)}\cong\left(\dfrac{\overline{\overline{\Omega}}^p_{X}}{\overline{\Omega}^p_{X}}\right)_o \cong \dfrac{\overline{\overline{\Omega}}^p_{X, o}}{\overline{\Omega}^p_{X, o}}
\end{equation}
by applying Cartan Theorem B for the coherent sheaf $\overline{\Omega}^p_{X}$.
Then $s^{(p)}=\dim_{\mathbb{C}}\overline{\overline{\Omega}}^p_{X, o}\big/\overline{\Omega}^p_{X, o}$, where both stalks are finitely generated modules over the local ring $\mathcal{O}_{X, o}$. This description is clearly independent of the Stein representative $X$ chosen.

By \cite[Theorem (1.3)]{StrStee} (see also \cite{greuel}) the invariants $s^{(p)}$ vanish for $0 \leq p \leq n-2$, whereas for $p=n-1$ 
\begin{equation}\label{eq:irreg}
    q(X, o)=s^{(n-1)}=\dim_{\mathbb{C}}\overline{\overline{\Omega}}^{n-1}_{X, o}\big/\overline{\Omega}^{n-1}_{X, o}=\dim_{\mathbb{C}} H^0(\widetilde{X}\setminus E, \Omega^{n-1}_{\widetilde{X}} ) \big/H^0(\widetilde{X}, \Omega^{n-1}_{\widetilde{X}}) 
\end{equation}
is an interesting invariant called the \textit{`irregularity'} of the isolated singularity $(X, o)$ (first discussed in \cite{wahl, Yau2}). 

For $p=n$, the invariant $s^{(n)}$ agrees with the \textit{`geometric genus'} (originally proved in \cite{Laufer72, Yau77}), defined by Wagreich as
\begin{equation}\label{eq:p_g}
    p_g(X, o)= \dim_{\mathbb{C}}(R^{n-1}\phi_{*}\mathcal{O}_{\widetilde{X}})_o=\dim_{\mathbb{C}}H^{n-1}(\widetilde{X}, \mathcal{O}_{\widetilde{X}}).
\end{equation}
It counts the 
adjunction conditions imposed by the singularity, cf. \cite{MT}. Moreover, the holomorphic $n$-forms of first kind also agree with the $L^2$-integrable $n$-forms on $U$ in the sense of Griffiths \cite[II.(a)]{Griffiths}. 

\begin{remark} A  normal surface singularity is called {\it rational} if $p_g=0$.
This analytic property can be characterized topologically as well, 
 it only depends on the resolution graph $\Gamma_\phi$, which in this case  is  called {\it rational}, see \cite{Artin62,Artin66,Laufer72} or Remark \ref{rem:badvertices} (b). 
 
If $(X,o)$ is a normal surface singularity, then $q\leq p_g-h^1(\calO_E)$, with equality if
$(X,o)$ is weighted homogeneous (cf. \cite{wahl,StrStee}, see also subsection \ref{ex:whq} below). 
In particular, $q=0$ for all rational singularities. 
\end{remark}

Comparing identity (\ref{eq:indepforirreg&p_g}) and Theorem \ref{th:properties} we get the following categorifications of the irregularity and the geometric genus:

\begin{cor}\label{cor:irreg&p_gcategorification} Let $(X, o)$ be a complex analytic normal  isolated singularity of dimension $n \geq 2$. Fix a  Stein representative $X$ 
and $\phi: \widetilde{X}\rightarrow X$, a resolution with exceptional divisor $E$. Then

    (a) $\mathbb{H}_*\big(\overline{\Omega}^{n-1}_{X, o} \hookrightarrow_{\mathcal{O}_{X, o}} \overline{\overline{\Omega}}^{n-1}_{X, o}\big) \cong \mathbb{H}_*\big(H^0(\widetilde{X}, \Omega^{n-1}_{\widetilde{X}})\hookrightarrow_{H^0(\widetilde{X}, \mathcal{O}_{\widetilde{X}})} H^0(\widetilde{X}\setminus E, \Omega^{n-1}_{\widetilde{X}})\big)$ is a well-defined invariant of the singularity and has Euler characteristic the irregularity $q(X, o)$.

    (b) $\mathbb{H}_*\big(\overline{\Omega}^{n}_{X, o} \hookrightarrow_{\mathcal{O}_{X, o}} \overline{\overline{\Omega}}^{n}_{X, o}\big) \cong \mathbb{H}_*\big(H^0(\widetilde{X}, \Omega^{n}_{\widetilde{X}})\hookrightarrow_{H^0(\widetilde{X}, \mathcal{O}_{\widetilde{X}})} H^0(\widetilde{X}\setminus E, \Omega^{n}_{\widetilde{X}})\big)$ is a well-defined invariant of the singularity  and has Euler characteristic the geometric genus $p_g(X, o)$.\\
    \noindent Both can be computed in some finite rectangle (depending on the CDP realization chosen). We call  them
the \emph{`analytic lattice homology of $p$-forms'} associated with the isolated singularity $(X,o)$ and denote them by $\mathbb{H}_*((X, o); \Omega^p)$ with $p \in \{ n-1, n\}$.
\end{cor}

\begin{nota}\label{nota:OMN}
    Let us fix some further notations for this section for the sake of better readability: $\mathcal O:=
H^0(\widetilde{X},\calO_{\widetilde{X}}), \ M^p:=H^0(\widetilde{X}\setminus E, \Omega^{p}_{\widetilde{X}})$ and $N^p:=H^0(\widetilde{X}, \Omega^{p}_{\widetilde{X}})$ for $p \in \{n-1, n\}$. 
\end{nota}

\begin{proof}[Proof of Corollary \ref{cor:irreg&p_gcategorification}]
As we have already seen that the finite dimensional quotient $\mathbb{C}$-vector spaces $M^p/N^p$ are independent of the resolution and representative chosen (cf. (\ref{eq:indepforirreg&p_g})), we have the following two more facts to check to be able to use Theorem \ref{th:properties} \textit{(c)}:
\begin{align}
    \text{ the} \text{  map } &\,\ \ \ \ \mathcal{O}=H^0(\widetilde{X}, \mathcal{O}_{\widetilde{X}}) \xrightarrow{(\phi^*)^{-1}} H^0(X, \mathcal{O}_X) \xrightarrow{\rho_o} \mathcal{O}_{X, o} \ \ \ \ \text{ induces an isomorphism} \label{eq:idekellanormalitas}\\  &\,\mathcal{O}/{\rm Ann}_{\mathcal{O}}(M^p/N^p) \cong \mathcal{O}_{X, o}/{\rm Ann}_{\mathcal{O}_{X, o}}\big(\overline{\overline{\Omega}}_{X, o}^p/\overline{\Omega}^p_{X, o}\big) \text{ for }p \in\{n-1, n\} \text{, and} \nonumber \\
    \text{the} &\, \text{ submodules } H^0(\widetilde{X}, \Omega^{p}_{\widetilde{X}})\leq  H^0(\widetilde{X}\setminus E, \Omega^{p}_{\widetilde{X}}), \ p \in\{n-1, n\} \text{ and} \label{eq:realizalhatosag} \\
    &\ \hspace{30mm} \overline{\Omega}^p_{X, o} \leq \overline{\overline{\Omega}}^p_{X, o}, \hspace{15mm} \ p \in\{n-1, n\} \text{ are realizable.} \nonumber 
\end{align} 

For (\ref{eq:idekellanormalitas}) let us first notice that the morphism $\phi^*: H^0(X, \mathcal{O}_X) \rightarrow H^0(\widetilde{X}, \mathcal{O}_{\widetilde{X}})$ is indeed invertible, since $(X, o)$ is normal. Secondly, let us denote by $\psi: M^p/N^p \xrightarrow{\cong} \overline{\overline{\Omega}}_{X, o}^p/\overline{\Omega}^p_{X, o}$ the vector space isomorphism of (\ref{eq:indepforirreg&p_g}).  Then for any $m + N^p \in M^p/N^p$ and $f \in \mathcal{O}$ we have 
\begin{equation*}
    \psi(f \cdot (m+N^p))=(\rho_o\circ(\phi^*)^{-1})(f)\cdot \psi(m+N^p),
\end{equation*}
i.e., this isomorphism in compatible with the module structures.
Thus, to obtain (\ref{eq:idekellanormalitas}), it is enough to prove for some ideal $I_{\rm Ann} \triangleleft \mathcal{O}_{X, o}, \ I_{\rm Ann} \leq {\rm Ann}_{\mathcal{O}_{X, o}}\big(\overline{\overline{\Omega}}^p_{X, o}/\overline{\Omega}^p_{X, o} \big)$ the isomorphism 
$$\rho_o\circ(\phi^*)^{-1}:\mathcal{O}/\phi^*(\rho_o^{-1}(I_{\rm Ann})) \xrightarrow{\cong}\mathcal{O}_{X, o}/I_{\rm Ann}.$$
Since the quotient $\overline{\overline{\Omega}}^p_{X, o}/\overline{\Omega}^p_{X, o}$ is finite dimensional, for a  large enough lattice point $l \in \mathbb{Z}_{> 0}$ the ideal $\phi_*\mathcal{O}_{\widetilde{X}}(-lE)_o \leq {\rm Ann}_{\mathcal{O}_{X, o}}\big(\overline{\overline{\Omega}}^p_{X, o}/\overline{\Omega}^p_{X, o}\big)$ and we might choose it as $I_{\rm Ann}$. Indeed, since $(X, o)$ is normal, $\phi^*\Big(\rho_o^{-1}\big(\phi_*\mathcal{O}_{\widetilde{X}}(-lE)_o\big)\Big)=H^0(\widetilde{X}, \mathcal{O}_{\widetilde{X}}(-lE))$ and from the short exact sequence of sheaves
\begin{equation*}
    0 \rightarrow \phi_*\mathcal{O}_{\widetilde{X}}(-lE) \rightarrow \mathcal{O}_{X} \rightarrow \mathcal{O}_{X}/\phi_*\mathcal{O}_{\widetilde{X}}(-lE) \rightarrow 0
\end{equation*}
(with the last term being supported only on the singular point $o$) we obtain the long exact sequence of sheaf cohomologies
\begin{center}
\begin{tikzcd}[column sep = 3mm]
 0 \rar & \underbrace{H^0(X, \phi_*\mathcal{O}_{\widetilde{X}}(-lE)) }_{=H^0(\widetilde{X}, \mathcal{O}_{\widetilde{X}}(-lE))} \rar & \underbrace{ H^0(X, \mathcal{O}_X)}_{=H^0(\widetilde{X}, \mathcal{O}_{\widetilde{X}})} \rar & \underbrace{H^0(\mathcal{O}_{X}/\phi_* \mathcal{O}_{\widetilde{X}}(-lE))}_{=\mathcal{O}_{X, o}/I_{\rm Ann}} \rar & \underbrace{H^1(X, \phi_*\mathcal{O}_{\widetilde{X}}(-lE)) }_{=0, \text{ by coherence}} \rar & \ldots,
\end{tikzcd}
\end{center}
which implies the desired isomorphism.

 For (\ref{eq:realizalhatosag}) we will use the extended divisorial valuations corresponding to the components of the exceptional divisor $E$. We will concentrate on the $\mathcal{O}$-submodule $N^p \leq M^p$, the case of $\overline{\Omega}^p_{X, o} \leq \overline{\overline{\Omega}}^p_{X, o}$ can be treated similarly.

    Let $\cup_{v\in\mathcal{V}}E_v$ be the irreducible decomposition of $E$ (with reduced structure). Then each $E_v$ provides a discrete valuation
\begin{equation}\label{eq:v_vdivisorial}
    \frv_v:{\calO}\to \mathbb{N}\cup \{\infty\}\text{ by }f\mapsto {\rm ord}_{E_v}(f), \text{ the order of vanishing of } f \text{ along } E_v.
\end{equation}
[Similarly, on $\mathcal{O}_{X, o}$ we can consider the valuations $g \mapsto {\rm ord}_{E_v}(g \circ \phi)$.]\\
These can be extended to valuations $\frv^{M}_v:M^p\to \Z\cup\{\infty\}$ on $M^p$ ($v \in \mathcal{V}$) as follows.  Let $x$ be a generic point of $E_v$, and let us fix 
some local neighborhood $U_x \subset \widetilde{X}$ of $x$ with local analytic coordinates 
$(u_1,\ldots, u_n)$ such that $U_x\cap E_v=\{u_1=0\}$. 

\textit{(b)} If $p=n$ then any $\omega\in M^n\setminus\{0\}$  has the form 
$adu_1\wedge \cdots \wedge du_n$ locally in $U_x$, where $a$ is meromorphic in $U_x$ with possible poles along $\{u_1=0\}$. Then $a$ can be written as $u_1^{\frv_v(a)}a'$, where 
$a'(x)\in\C^*$. Then we define 
\begin{equation}\label{eq:v_v^Mdivisorial}
    \frv_v^M:M^n\rightarrow \overline{\bZ}:\ \omega \ (\text{with } \omega\big|_{U_x} = adu_1 \wedge \ldots \wedge du_n ) \mapsto \frv_v(a)\in \Z\ \text{ and } \omega =0 \mapsto \infty.  
\end{equation}
It is independent of all the choices. 

\textit{(a)} If $p=n-1$ then, using the same local picture, we define 
\begin{align}
    \frv_v^M:M^{n-1}\rightarrow \overline{\bZ}:\ \omega \Big( \text{with }\omega\big|_{U_x} = \sum_i a_idu_1\wedge \cdots \wedge \widehat{du_i}\wedge \cdots\wedge  du_n \Big) \mapsto &\  \min_i\{\frv_v(a_i)\}\in\Z \nonumber \\ \text{ and } \omega =0 \mapsto &\ \infty.  
\end{align}
Once again, the values are independent of all the choices.

\noindent [Similarly, we can use the same formulas to define the extended valuations of $\overline{\overline{\Omega}}^p_{X, o}$ for $p\in \{n, n-1\}$.]

By a direct verification, if $f\in \calO$ and $\omega\in M^p$ then we have $\frv^M_v(f\omega)=\frv_v(f)+\frv^M_v(\omega)$ for any $v\in\mathcal {V}$, i.e., for any $v\in\mathcal{V}$ and $p \in \{ n-1, n\}$ the pair $(\frv_v, \frv^M_v)$ is an extended discrete valuation (in the sense of Definition \ref{def:edv}).

Denote by $\mathcal{D}=\mathcal{D}_\phi:= \{(\frv_v, \frv^M_v)\}_{v\in \mathcal{V}}$ the collection of extended  discrete valuations corresponding to 
the irreducible exceptional divisors 
$\{E_v\}_{v \in \mathcal{V}}$. We claim, that these give a realization of the finite codimensional submodule $N^p \leq M^p$ ($p\in \{n-1,n\}$). Indeed, the corresponding $\Z^r$-filtration $\mathcal{F}_{\mathcal{D}}^M$ in $M^p$ (cf. (\ref{eq:fdMl})) satisfies $\mathcal{F}_{\mathcal{D}}^M(0)= N^p=H^0(\widetilde{X}, \Omega^p_{\widetilde{X}})$, since the meromorphic $p$-forms $\omega \in H^0(\widetilde{X}\setminus E, \Omega_{\widetilde{X}}^p)$, with  $\mathfrak{v}_v^M(\omega) \geq 0$ for all $v \in \mathcal{V}$, extend to $\widetilde{X}$ (for both $p\in \{n-1, n\}$). 

\noindent [The same argument applies for the extended valuations of $\overline{\overline{\Omega}}^p_{X,o}$ to yield the realizability of $\overline{\Omega}^p_{X, o}$.]

Thus, with both (\ref{eq:idekellanormalitas}) and (\ref{eq:realizalhatosag}) satisfied, we can use Theorem \ref{th:properties} \textit{(c)} to obtain the desired well-definedness and Euler characteristic.
\end{proof}

\begin{remark}\label{rem:filtrationreformulation}
    If we identify the abstract lattice $\Z^r$ ($r=|\mathcal{V}|$) with the free abelian group $\mathbb{Z}\langle \{ E_v\}_{v \in \mathcal{V}}\rangle$ of divisors supported on $E$, then the multifiltrations (defined in (\ref{eq:fdMl})) corresponding to the collection $\mathcal{D}:=\mathcal{D}_\phi= \{(\frv_v, \frv^M_v)\}_{v\in \mathcal{V}}$ from the previous proof satisfy the following: 
    \begin{equation}
        \text{ for any } \ell \in \mathbb{Z}^r: \ \mathcal{F}_{\mathcal{D}}(\ell) =H^0(\widetilde{X},\calO_{\widetilde{X}}(-\ell)) \text{, whereas } \mathcal{F}^M_{\mathcal{D}}(\ell)=H^0(\widetilde{X}, \Omega^p_{\widetilde{X}}(-\ell)).
    \end{equation}
    Therefore, the corresponding $\mathfrak{h}_\mathcal{D}$ and $\mathfrak{h}_{\mathcal{D}}^\circ$ functions (cf. (\ref{eq:handhcirc0})) satisfy for any $\ell \in \mathbb{Z}^r$:
    \begin{equation}\label{eq:highdimHilbert}
        \mathfrak{h}_{\mathcal{D}}(\ell) = \dim_{\mathbb{C}}\dfrac{H^0(\widetilde{X}, \mathcal{O}_{\widetilde{X}})}{H^0(\widetilde{X}, \mathcal{O}_{\widetilde{X}}(-\ell))}, \text{ whereas } \mathfrak{h}^\circ_\mathcal{D}(\ell) = \dim_{\mathbb{C}}\dfrac{H^0(\widetilde{X}\setminus E, \Omega^p_{\widetilde{X}})}{H^0(\widetilde{X}, \Omega^p_{\widetilde{X}}(\ell))}.
    \end{equation}
\end{remark}

\bekezdes \textbf{Nonpositivity Theorem.} From the local sheaf theoretical description of the algebra and modules used in the definition of the analytic lattice homology of $p$-forms we see that the assumptions of the Nonpositivity Theorem \ref{th:upperbound} are satisfied with field $k=\mathbb{C}$ and local $k$-algebra $\mathcal{O}_{X, o}$. Hence we have the following statement:

\begin{cor} \label{cor:upperboundhighdim}
 Let $(X, o)$ be a complex analytic normal isolated singularity of dimension $ \geq 2$. Then for both $p \in \{\dim X-1, \dim X\}$ and for every positive integer $n>0$ the corresponding $S_n$ space is contractible, hence the weight-grading of the \emph{reduced} analytic lattice homology of $p$-forms is nonnegative, i.e.,
 $\mathbb{H}_{{\rm red},q, -2n}((X, o); \Omega^p) = 0$
    for all $q\geq 0$ and $n >0$.  \\ Moreover, if $s^{(p)} \neq 0$ (with $p \in \{\dim X-1, \dim X\}$), then the upper bound is sharp: $S_0$ is not connected.
\end{cor}

\subsection{Comparison with the analytic lattice homology after \cite{AgNe1,AgNeHigh,NBook}}\label{ss:compar}\,

In the case of $p=n$, \'Agoston and the first author defined a categorification of the geometric genus via the so-called analytic lattice homology $\mathbb{H}_{an, *}(X, o)$ of the isolated singularity $(X, o)$ of dimension $n \geq 2$ \cite{AgNeHigh} (see also \cite{AgNe1,NBook} for the normal surface singularity case).  In this subsection we introduce this invariant and compare it with the construction of the previous subsection \ref{ss:iso}. It turns out, that whenever the analytic lattice homology in the sense of 
\cite{AgNe1,AgNeHigh,NBook}
exists (i.e., it is well-defined and categorifies $p_g$), then the two theories agree. 
Therefore, $\mathbb{H}_*\big(({X, o}); \Omega^{n}\big)$, defined for any normal isolated singularity, is a generalization of $\mathbb{H}_{an, *}(X, o)$.  

\begin{define}\cite{AgNeHigh}
Let $(X, o)$ be a complex analytic isolated singularity of dimension  $n \geq 2$. Fix a \textit{good} resolution $\phi:\widetilde{X} \rightarrow X$ with reduced exceptional divisor $E$, having prime decomposition $\cup_{v \in \mathcal{V}}E_v$ (i.e., $E$ is a simple normal crossing divisor). Suppose that $h^{n-1}(\calO_E):=\dim_{\mathbb{C}}\big(H^{n-1}(E, \calO_E)\big)=0$.
In the $n=2$ surface case, this is equivalent with the fact that 
the link $L_X$ of the singularity is a rational homology sphere, equivalently, the dual resolution graph $\Gamma_{\phi}$ is a tree of rational curves, cf. \cite{AgNe1, NBook}. Then the \textit{`analytic lattice homology'} of $(X, o)$, denoted by $\mathbb{H}_{an, *}(X, o)$ (after \cite{AgNe1, AgNeHigh, NBook}), is defined as follows. Consider the `lattice' of effective divisors in $\widetilde{X}$ supported on $E$: we identify it with $(\mathbb{Z}_{\geq 0})^r$ with fixed basis $\{E_v\}_{v \in \mathcal{V}}$ (here $r=|\mathcal{V}|$). Then the weight function is defined as
\begin{align}\label{eq:w_anhighdim}
    w_{{an}, 0}: (\mathbb{Z}_{\geq 0})^r \rightarrow \mathbb{Z}, \ \ell& \, \mapsto \mathfrak{h}(\ell) - h^{n-1}(\mathcal{O}_\ell) \\ \text{ and is extended to higher dimensional cubes } \square_q \text{ as }& w_{an, q}(\square):=\max\{w_{an, 0}(\ell)\,:\, \ell \text{ is a vertex of }\square_q\}, \nonumber
\end{align}
where the function $ (\Z_{\geq 0})^r\ni \ell\mapsto \frh(\ell)=\dim_{\mathbb{C}} \, H^0(\widetilde{X},\calO_{\widetilde{X}})/
H^0(\widetilde{X},\calO_{\widetilde{X}}(-\ell)) $ is the analytic Hilbert function associated with the resolution (studied, e.g., in \cite{CHR}) and $h^{n-1}(\mathcal{O}_\ell)=\dim_{\mathbb{C}}\big(H^{n-1}(\ell, \mathcal{O}_\ell)\big)$. Then $\mathbb{H}_{an, *}(X, o)$ is the lattice homology module $\mathbb{H}_*((\mathbb{R}_{\geq 0})^r, w_{an})$ associated with this setup. According to \cite{AgNeHigh} it is independent of the resolution $\phi$, can be computed on a finite rectangle, moreover, it categorifies the geometric genus $p_g(X, o)$.
\end{define}

\begin{theorem}\label{th:newvsoldhighdim}
$\mathbb{H}_*((X, o); \Omega^n)$ generalizes the analytic lattice homology of  $(X, o)$ after \cite{AgNe1, AgNeHigh, NBook}, as for such a complex analytic normal isolated singularity of dimension $n \geq 2$, having some (or, equivalently, every) good resolution $\phi: \widetilde{X} \rightarrow X$ satisfying $h^{n-1}(\mathcal{O}_E)=0$, we have \begin{center}
    $\mathbb{H}_{an, *}(X, o)\cong   \mathbb{H}_*((X, o); \Omega^n)$.
\end{center}
\end{theorem}

\begin{proof}
    We will proceed similarly to the proof of Theorem \ref{th:equivforcurves}. We fix a good resolution $\phi: \widetilde{X} \rightarrow X$, with exceptional divisor $E=\cup_{v=1}^rE_v$, satisfying the required condition $h^{n-1}(\mathcal{O}_E)=0$. On one hand, $\mathbb{H}_{an, *}(X, o)$ is computed on the lattice $(\mathbb{Z}_{\geq 0})^r$ with weight function $w_{an}$ defined in (\ref{eq:w_anhighdim}). On the other hand, by Corollary \ref{cor:irreg&p_gcategorification} \textit{(b)} we have $\mathbb{H}_*\big(({X, o} );\Omega^{n}\big) \cong \mathbb{H}_*\big(H^0(\widetilde{X}, \Omega^{n}_{\widetilde{X}})\hookrightarrow_{H^0(\widetilde{X}, \mathcal{O}_{\widetilde{X}})} H^0(\widetilde{X}\setminus E, \Omega^{n}_{\widetilde{X}})\big)$. Now, we can consider the realization $\mathcal{D}=\mathcal{D}_\phi=\{(\mathfrak{v}_v, \mathfrak{v}_v^M)\}_{v=1}^r$ of $H^0(\widetilde{X}, \Omega^{n}_{\widetilde{X}})$ using the extended discrete valuations corresponding to the exceptional divisors $\{E_v\}_{v=1}^r$ (see (\ref{eq:v_vdivisorial}) and (\ref{eq:v_v^Mdivisorial})) as before. The corresponding weight function $w_{\mathcal{D}}$ (cf. (\ref{eq:handhcirc0})) is defined on the same lattice $\mathbb{Z}^r$ and we claim, that $w_{\mathcal{D}}\big|_{(\mathbb{Z}_{\geq 0})^r}=w_{an}$. Indeed, we only have to prove that $\mathfrak{h}_{\mathcal{D}}^\circ(0)-\mathfrak{h}_{\mathcal{D}}^{\circ}(\ell) = h^{n-1}(\mathcal{O}_{\ell})$ for every $\ell \in \mathbb{Z}^r, \ \ell \geq 0$ (since $\mathfrak{h}_{\mathcal{D}}\big|_{(\mathbb{Z}_{\geq 0})^r} \equiv \mathfrak{h}$). By Remark \ref{rem:filtrationreformulation} we have for every $\ell \in (\mathbb{Z}_{\geq 0})^r$
\begin{equation}\label{eq:elso}
    \mathfrak{h}_{\mathcal{D}}^\circ(0)-\mathfrak{h}_{\mathcal{D}}^{\circ}(\ell) =\dim_{\mathbb{C}} \, \mathcal{F}_{\mathcal{D}}^M(-\ell)/\mathcal{F}_{\mathcal{D}}^M(0)=\dim_{\mathbb{C}} 
H^0(\widetilde{X},\Omega^n_{\widetilde{X}}(\ell))/
H^0(\widetilde{X},\Omega^n_{\widetilde{X}}).
\end{equation}
Then, using the exact sequence
$$0\to \Omega^n_{\widetilde{X}}\to \Omega^n_{\widetilde{X}}(\ell)\to \Omega^n_{\widetilde{X}}(\ell)\otimes \calO_{\ell}\to 0,$$
the Grauert--Riemenschneider vanishing 
$H^1(\widetilde{X}, \Omega^n_{\widetilde{X}})=0$ (cf. \cite{GrRie,Laufer72,Ram} or \cite[Theorem 6.4.3]{NBook}), and Serre duality (cf., e.g., \cite[Theorem 2.2]{Okuma} and the references therein), we get 
\begin{equation}\label{eq:SerreH} \frh_{\mathcal{D}}^\circ (0)-\frh_{\mathcal{D}}^\circ (\ell)= \dim_{\mathbb{C}}  H^0( \Omega^n_{\widetilde{X}}(\ell)\otimes \calO_\ell)=
\dim_{\mathbb{C}} H^{n-1}(\calO_\ell)^*=h^{n-1}(\calO_\ell)
\end{equation}
for all  $\ell \in (\mathbb{Z}_{\geq 0})^r$; i.e.,
$\mathfrak{h}_{\mathcal{D}}^\circ(\ell)=h^{n-1}(\calO_{\widetilde{X}}) - h^{n-1}(\calO_\ell)$. Moreover, \cite[Lemma 4.4.2]{AgNeHigh} states that the pair $\ell\mapsto (\mathfrak{h}(\ell), h^{n-1}(\calO_{\widetilde{X}}) - h^{n-1}(\cO_{\ell}))$ satisfies the Combinatorial Duality Property on the positive orthant $(\mathbb{Z}_{\geq 0})^r$. Also,  by the Formal Function Theorem 
(see also the discussion in subsection \ref{ss:zcoh}) we obtain:
$h^{n-1}(\calO_\ell)=h^{n-1}(\calO_{\widetilde{X}})=p_g$ for $\ell \gg 0$ large enough. Therefore, by Proposition \ref{prop:rectCDP->fullCDP}, the pair $(\mathfrak{h}_{\mathcal{D}}, \mathfrak{h}_{\mathcal{D}}^\circ)$ corresponding to the realization $\mathcal{D}$ of $H^0(\widetilde{X}, \Omega^n_{\widetilde{X}})$ also satisfies the Combinatorial Duality Property on the whole lattice $\mathbb{Z}^r$. Consequently $\mathcal{D}$ is in fact a CDP realization of $H^0(\widetilde{X}, \Omega^n_{\widetilde{X}})$, hence, $w_{\mathcal{D}}$ computes $\mathbb{H}_*\big(H^0(\widetilde{X}, \Omega^{n}_{\widetilde{X}})\hookrightarrow_{H^0(\widetilde{X}, \mathcal{O}_{\widetilde{X}})} H^0(\widetilde{X}\setminus E, \Omega^{n}_{\widetilde{X}})\big)$ already on the rectangle $R(0, \infty)$ (see part \textit{(a)} of Theorem \ref{th:properties}). But there it agrees with $w_{an}$.
\end{proof}

As a conclusion, we would like to highlight, that although the two definitions give the same lattice homology  module, our new construction has the novelty, that the independence of the resolution statement  follows without any additional assumption (compare with \cite{AgNe1, AgNeHigh, NBook}, where the assumption $h^{n-1}(\calO_E)=0$ --- or, in the $n=2$ case, the equivalent topological assumption that the link is a rational homology sphere --- was necessary). 
Moreover, in [loc.cit.],   the Combinatoral Duality Property of the pair 
$(\frh_{\mathcal{D}_{\phi}},\frh_{\mathcal{D}_{\phi}}^\circ)$ associated with $\mathcal{D}_{\phi}$ was used already in the construction\,/\,definition 
of the lattice (co)homology, 
while this was  proved, again, under the same assumption  $h^{n-1}(\calO_E)=0$. 
The CDP  in the present case is not a priori needed
thanks to `doubling procedure' from Lemma \ref{lem:duplatrukk} (and to the corresponding new, adopted definition). 

\begin{remark}
    In their article \cite{AgNeHigh}, Ágoston and the first author did not assume that the isolated singularity $(X, o)$ of dimension $n \geq 2$ is \emph{normal} as well, nevertheless, their analytic lattice homology $\mathbb{H}_{an,*}(X, o)$ is, in fact, an invariant of the normalization $(\overline{X}, \overline{o})$. Indeed, every resolution $\phi: \widetilde{X} \rightarrow X$ factors through $(\overline{X}, \overline{o})$ and $\mathbb{H}_{*}(X, o)$ is,  computed on the resolution as $\mathbb{H}_*( H^0(\widetilde{X}, \Omega^n_{\widetilde{X}})) \hookrightarrow _{H^0(\widetilde{X}, \mathcal{O}_{\widetilde{X}})} H^0(\widetilde{X} \setminus E, \Omega_{\widetilde{X}}^n))$ \-(cf. the proof of Theorem \ref{th:newvsoldhighdim}).

    On the other hand, for non-normal isolated singularities, the identifications of Corollary \ref{cor:irreg&p_gcategorification} do not hold, hence $\mathbb{H}_*\big(\overline{\Omega}^{p}_{X, o} \hookrightarrow_{\mathcal{O}_{X, o}} \overline{\overline{\Omega}}^{p}_{X, o}\big), \ p \in \{n-1, n\}$ might give different categorifications of $q(X, o)$ and $p_g(X, o)$ from the analytic lattice homology $\mathbb{H}_*((\overline{X}, \overline{o}); \Omega^p)$ of $p$-forms of the normalization $(\overline{X}, \overline{o})$. Using the language of Notation \ref{nota:OMN} the difference between these two theories consists of the choice of the ring $\mathcal{O}:= H^0(\widetilde{X}, \mathcal{O}_{\widetilde{X}})$ versus  $H^0(X, \mathcal{O}_X)$ over which we consider the module $M^p$ and the realizable submodule $N^p$, with $p \in \{n-1, n\}$.

    We will investigate this phenomenon in a forthcoming project.
\end{remark}

\subsection{The cohomological cycle and the Gorenstein property}\label{ss:zcoh}\,

Recall that in subsection \ref{ss:ddep} (in the general setup of a CDP realization $\mathcal{D}$ of some finite codimensional submodule $N\leq M$) we introduced the notation $d_{\mathcal{D}}$ for the unique smallest lattice point, such that the lattice homology was already realized in the finite 
rectangle $R(0,d_{\mathcal{D}})$ (i.e., the smallest one satisfying $\mathcal{F}_{\mathcal{D}}^M(-d_{\mathcal{D}})=M$). [In the case of reduced curve singularities this lattice point was the conductor element $\mathbf{c}$, see Corollary \ref{cor:dD=cforcurves}.]

We wish to identify a concrete cycle which has this property in the 
case of the analytic lattice homology of $n$-forms associated with the normal isolated singularity $(X, o)$ of dimension $n \geq 2$ computed via the realization given in the proof of Corollary \ref{cor:irreg&p_gcategorification}. More precisely, we fix a resolution $\phi: \widetilde{X} \rightarrow X$ and use the (CDP) realization $\mathcal{D}=\{ (\mathfrak{v}_v, \mathfrak{v}_v^M)\}_{v \in \mathcal{V}}$ of $H^0(\widetilde{X}, \Omega_{\widetilde{X}}^{n}) \leq H^0(\widetilde{X} \setminus E, \Omega_{\widetilde{X}}^{n})$ via the extended divisorial valuations (defined in (\ref{eq:v_vdivisorial}) and (\ref{eq:v_v^Mdivisorial})) corresponding to the irreducible components of the exceptional divisor $E$. It turns out that in this setting $d_{\mathcal{D}}$ is the {\it cohomological cycle} $Z_{coh}\in (\Z_{\geq 0})^r$
of the resolution $\phi$. We  recall its definition briefly (for more see, e.g., \cite{MR,Ishii,AgNeHigh,NBook}).
If $p_g=0$ then we set $Z_{coh}=0$. If $p_g>0$ 
then there exists a unique minimal cycle $Z_{coh}>0$ such that
$h^{n-1}(\calO_{\tX})=h^{n-1}(\calO_{Z_{coh}})$. It  has the property that for any $\ell\not\geq Z_{coh}$
one has $h^{n-1}(\calO_\ell)<h^{n-1}(\calO_{\tX})$. 
Then clearly, by the identifications (\ref{eq:SerreH}) and (\ref{eq:elso}), $Z_{coh}$ corresponds to $d_{\mathcal{D}}$ (it satisfies the conditions imposed in subsection \ref{ss:ddep}).
Now, from  Corollary \ref{cor:finrectcomb} we get

\begin{cor}
    In the above setup, for any lattice point $d$ with $2\cdot Z_{coh}\leq d\leq \infty$, the inclusion 
$R(0,d)\hookrightarrow \bR^r$ induces a  bigraded $\Z[U]$-module 
isomorphism 
\begin{center}
    $\bH_*(R(0,d),w_{\mathcal{D}^\natural})\xrightarrow{\cong} \bH_*(\bR^r,w_{\mathcal{D}^\natural}) \cong \mathbb{H}_*((X, o); \Omega^n)$, 
\end{center}
where $\mathcal{D}^\natural$ is the doubling defined in Lemma \ref{lem:duplatrukk} \textit{(2)}.  (We do not need this doubling if $\mathcal{D}$ is already a CDP realization; then $d\geq Z_{coh}$ is enough.) 
\end{cor}

Let  us  fix a  \textit{canonical divisor} $K_{\widetilde{X}}$, i.e., a divisor satisfying  $\Omega^n_{\widetilde{X}}=\calO_{\widetilde{X}}(K_{\widetilde{X}})$. 

\begin{remark}
    Assume that $n=2$ and let $(\,,)$ denote the natural negative definite intersection form on $H_2(\widetilde{X},\Z)=\Z^r=\mathbb{Z}\langle \{ E_v\}_{v \in \mathcal{V}}\rangle$. 
Recall, that in this setup $K_{\widetilde{X}}$ numerically (on the cohomological level) can be represented by a rational cycle $Z_K\in 
\Z^r\otimes \Q$, the \textit{(anti)canonical cycle},  defined by $(-Z_K, E_v)={\rm deg}(K_{\widetilde{X}}\big|_{E_v}) $. Let 
$\lfloor Z_K\rfloor $ denote the integral part of $Z_K$ and  
$\lfloor Z_K\rfloor_+ =\max\{\lfloor Z_K\rfloor,0\} $ be its positive integral part. 
Then, by \cite[Corollary 6.4.7]{NBook} or by \cite[4.1.5]{AgNe1} we have $Z_{coh}\leq \lfloor Z_K\rfloor_+$. (Note that the cohomological cycle usually 
depends on the analytic structure of $(X,o)$, however $Z_K$ 
can be determined from the topological type via the \textit{adjunction formulae} 
${\rm deg}(K_{\widetilde{X}}|_{E_v})=-2+2\cdot \mbox{genus}(E_v)-( E_v,E_v)$.
Recall also, that $Z_K$ can also be identified with $-c_1(\Omega^2_{\widetilde{X}})\in H^2(\widetilde{X},\Z)$.)
\end{remark}

Let us suppose, that the isolated normal complex analytic singularity $(X, o)$ is \textit{Gorenstein}, i.e., for any good resolution $\phi:\tX\to X$, with reduced exceptional divisor $E=\cup_{v \in \mathcal{V}}E_v$, there exists an integral cycle $Z_K\in \Z\langle\{ E_v\}_{v \in \mathcal{V}}\rangle=\Z^r$ such that
$\Omega_{\tX}^n=\calO_{\tX}(-Z_K)$ (i.e., the canonical divisor $K_{\widetilde{X}}$ is linearly equivalent to $-Z_K$). Equivalently, this means that $\Omega^n_{X\setminus \{o\}}\cong \Omega^n_{\widetilde{X} \setminus E}$ is holomorphically trivial. Its nowhere vanishing section (defined up to multiplication by a nonzero constant) is called the \textit{`Gorenstein form'} denoted by
$\omega_G\in H^0(X\setminus \{o\}, \Omega^n _{X\setminus \{o\}})$, while 
its pullback extends meromorphically and satisfies ${\rm div} (\phi^* \omega_G)={\rm div}_E(\phi^* \omega_G)=-Z_K$. Moreover, $Z_K$ is the largest possible pole a meromorphic $n$-form $\omega \in \Omega^n_{\widetilde{X}\setminus E}$ can achieve, as by normality $\omega /\phi ^* \omega_G \in H^0(\widetilde{X}, \mathcal{O}_{\widetilde{X}})$. This also implies that $M^n:=H^0(\widetilde{X}\setminus E, \Omega^{n}_{\widetilde{X}})$ is free of rank one above $\mathcal O:=
H^0(\widetilde{X},\calO_{\widetilde{X}})$, generated by $\omega_G$.

Let us write the positive part of  $Z_K$ as $Z_{K,+}$, as above.  Then  
\begin{center}$Z_{coh}\leq Z_{K,+} $ in any dimension $n\geq 2$ by  \cite[3.7]{Ishii}
\end{center} (for the $n=2$ case see also  \cite[4.20]{MR}). Then using the description $M^n = \omega_{G} \cdot \mathcal{O}$ and Proposition \ref{prop:symforideals} we see that the weight function $w_{\mathcal{D}^{\natural}}$ (corresponding to the doubling $\mathcal{D}^\natural$ from Lemma \ref{lem:duplatrukk} of the  realization $\mathcal{D}=\mathcal{D}_\phi:= \{(\frv_v, \frv^M_v)\}_{v\in \mathcal{V}}$ from the proof of Corollary \ref{cor:irreg&p_gcategorification}) is symmetric with respect to $2\cdot Z_{K, +}$ (cf. Proposition \ref{prop:symforoffsetideals}). More precisely, for any lattice point $\ell \in R(0, 2\cdot Z_{K, +}) \cap \mathbb{Z}^r$ we have $w_{\mathcal{D}^\natural, 0}(\ell) = w_{\mathcal{D}^\natural, 0}(2 \cdot Z_{K, +}-\ell)$, with $R(0, 2\cdot Z_{K, +})$ containing the rectangle $R(0, 2\cdot Z_{coh})$ required for computing the analytic lattice homology $\mathbb{H}_*((X, o); \Omega^n)$ of $n$-forms (cf. Theorem \ref{th:properties} \textit{(a)}).

In fact, we can verify this directly in the following specific situation. A good resolution is called \textit{`essential'} if $Z_K\geq 0$. In the case of normal surface singularities essential good resolutions
exist (e.g., the minimal good resolution is such, see  \cite[Example 6.3.4]{NBook}).
However, in higher dimensions there are singularities
without any essential good resolution (cf. \cite[Remark 3.6]{Ishii}).

\begin{prop}\label{prop:GorSymhighdim}
    For a Gorenstein normal isolated singularity $(X, o)$ with $\phi: \widetilde{X} \rightarrow X$ an essential good resolution, the height functions $\hh_{\mathcal{D}^\natural}$ and $\hh_{\mathcal{D}^\natural}^\circ$ corresponding to the doubling $\mathcal{D}^\natural$ of the realization $\mathcal{D}:=\mathcal{D}_\phi=\{ (\frv_v, \frv_v^M)\}_{v \in \mathcal{V}}$ given by the extended discrete valuations (\ref{eq:v_vdivisorial}) and (\ref{eq:v_v^Mdivisorial}), satisfy $\hh_{\mathcal{D}^\natural}^\circ (\ell)=\hh_{\mathcal{D}^\natural}(2\cdot Z_K-\ell)$ in the rectangle $R(0,2\cdot Z_K)$. 
That is, $\hh_{\mathcal{D}^\natural}^\circ$ is given by symmetrization of $\hh_{\mathcal{D}^\natural}$ with respect to $2\cdot Z_K\geq 0$. (If $\mathcal{D}$ is a CDP realization, then we do not need the doubling and the symmetry is with respect to $Z_K$).
\end{prop}

\begin{proof}
    
Since $H^0(\widetilde{X}, \Omega^n_{\widetilde{X}}(Z_K))\subset H^0(\widetilde{X}\setminus E, \Omega^n_{\widetilde{X}})$, and in both spaces the codimension of the submodule $H^0(\widetilde{X}, \Omega^n_{\widetilde{X}})$ is $p_g$ (use $\Omega_{\tX}^n=\calO_{\tX}(-Z_K)$ and $h^{n-1}(\mathcal{O}_{Z_K}) = p_g$ with the identification (\ref{eq:SerreH}) and (\ref{eq:elso})), we deduce that 
$H^0(\widetilde{X}, \Omega^n_{\widetilde{X}}(Z_K))= H^0(\widetilde{X}\setminus E, \Omega^n_{\widetilde{X}})$.
Hence
 \begin{equation*}
 \begin{split}
\hh_{\mathcal{D}^\natural}^\circ (\ell)&=\dim H^0(\widetilde{X} \setminus E, \Omega^n _{\tX})/H^0( \Omega^n _{\tX}(\lfloor\ell/2\rfloor)) \hspace{30mm} \text{ (use Remark \ref{rem:filtrationreformulation})} \\
 &=
 \dim H^0(\Omega^n _{\tX}(Z_K))/H^0(\Omega^n _{\tX}(\lfloor\ell/2\rfloor))=
\dim H^0(\calO_{\tX})/H^0(\calO_{\tX}(-Z_K+\lfloor\ell/2\rfloor)) \\
&=\hh_{\mathcal{D}^\natural}(2\cdot Z_K-\ell) \hspace{60mm} \text{ for any } 0\leq \ell\leq 2\cdot Z_K.\end{split}
\end{equation*}
\end{proof}

Still in the Gorenstein normal isolated singularity case, under the identification 
$$M^n = H^0(\widetilde{X} \setminus E, \Omega_{\widetilde{X}}^n) = H^0(\widetilde{X}, \Omega^n_{\widetilde{X}}(Z_K)) = \omega_G \cdot H^0(\widetilde{X}, \mathcal{O}_{\widetilde{X}}) =\omega_G \cdot \mathcal{O},$$
the submodule $N^n=H^0(\widetilde{X}, \Omega^n_{\widetilde{X}})$ corresponds to $\omega_G \cdot H^0(\widetilde{X}, \mathcal{O}_{\widetilde{X}}(-Z_K))$. Now, similarly to the original case (i.e., identity  (\ref{eq:indepforirreg&p_g})), we have the $\mathbb{C}$-vector space isomoprhism
\begin{equation}\label{eq:condidindep}
    \frac{M^n}{N^n} =\frac{H^0(\widetilde{X}\setminus E, \Omega^n_{\widetilde{X}})}{H^0(\widetilde{X}, \Omega^n_{\widetilde{X}})} \cong \frac{H^0(\widetilde{X}, \mathcal{O}_{\widetilde{X}})}{H^0(\widetilde{X}, \mathcal{O}_{\widetilde{X}}(-Z_K))} \cong \left(\frac{\phi_*(\mathcal{O}_{\widetilde{X}})}{\phi_*(\mathcal{O}_{\widetilde{X}}(-Z_K))}\right)_o \cong\frac{\mathcal{O}_{X, o}}{\phi_*(\mathcal{O}_{\widetilde{X}}(-Z_K))_o}.
\end{equation}
In fact, the last two identifications are Artin $\mathbb{C}$-algebra isomorphisms. Hence, we have the isomorphism 
\begin{center}
    $\mathcal{O}/{\rm Ann}_{\mathcal{O}}(M^n/N^n) \cong \mathcal{O}_{X, o} / \phi_*(\mathcal{O}_{\widetilde{X}}(-Z_K))_o$.
\end{center}
which, combined with identity (\ref{eq:idekellanormalitas}), implies the resolution-free description 
$$\phi_*(\mathcal{O}_{\widetilde{X}}(-Z_K))_o={\rm Ann}_{\mathcal{O}_{X, o}}\overline{\overline{\Omega}}^n_{X, o}/\overline{\Omega}^n_{X, o}$$
(compare this with \cite[1.6 Définition]{MT} in the hypersurface case).

\begin{define}\label{def:CONDINDEF}
    For Gorenstein normal isolated singularities $(X, o)$ of dimension $n \geq 2$ we call the ideal $\phi_*(\mathcal{O}_{\widetilde{X}}(-Z_K))_o={\rm Ann}_{\mathcal{O}_{X, o}}\overline{\overline{\Omega}}^n_{X, o}/\overline{\Omega}^n_{X, o}\triangleleft\mathcal{O}_{X, o}$ the \emph{conductor ideal} of $\mathcal{O}_{X, o}$. 
\end{define}
\noindent [Compare this definition with the description of the conductor ideal of a curve germ in Corollary \ref{cor:dD=cforcurves}.]

We can also check directly whether the conductor ideal is well-defined, i.e., independent of the resolution $\phi$ chosen:

  \begin{prop}\label{prop:CONDINDEP}
For any resolution $\phi:\widetilde{X}\to X$ set
$\mathcal{C}=\phi_*(\cO_{\widetilde{X}}(-Z_K))_o \subset 
\cO_{X,o}$. 
Then $\mathcal{C}$ is independent of the  choice of the resolution $\phi$.  
  \end{prop}

\begin{proof}
    Consider a resolution $\phi:\widetilde{X}\to X$ and another modification of it 
    $\psi:\widetilde{X}'\to \widetilde{X}$. Then 
    $\psi_*(\Omega^{n}_{\widetilde{X}'})=\Omega^{n}_{\widetilde{X}}$ (cf., e.g., \cite[page 156]{KollarMori}). That is,
    $\psi_*(\cO_{\widetilde{X}'}(-Z_K(\widetilde{X}')))=\cO_{\widetilde{X}}(-Z_K(\widetilde{X}))$.
    Therefore,  $$\phi_*\psi_*(\cO_{\widetilde{X}'}(-Z_K(\widetilde{X}')))=\phi_*\cO_{\widetilde{X}}(-Z_K(\widetilde{X})).$$ Since any two resolutions can be dominated by a third one \cite{Hir, Kollar}, this shows that $\phi_*\cO_{\widetilde{X}}(-Z_K(\widetilde{X}))$ is independent of the resolution. 

    In the case of surface singularities we can argue in the following elementary way as well.
Let $\phi:\widetilde{X}\to X$ be a resolution as above, and let $\psi:\widetilde{X}'\to\widetilde{X} $ be a  blow up of $\widetilde{X}$ at a point of the exceptional curve of $\phi$. Let $E_{new} $ be the exceptional curve of $\psi$.
Then by an elementary computation (see, e.g., \cite[Example 6.3.3 (4)]{NBook})  
$Z_K(\widetilde{X}')=\psi^*(Z_K(\widetilde{X}))-E_{new}$. Hence, for any function germ $g\in \cO_{X,o}$
the inequality $\phi^*(g)\geq Z_K(\widetilde{X})$ holds exactly when the counterpart 
 $\psi^*\phi^*(g)\geq Z_K(\widetilde{X}')$ holds. 
\end{proof}

Now, from identity (\ref{eq:condidindep}) and the discussion after Definition \ref{def:LC} we get the following equivalent description of the analytic lattice homology of $n$-forms in the Gorenstein case:

\begin{corollary}\label{cor:GORCOND} Assume that $(X,o)$ is  a Gorenstein normal isolated singularity of dimension $n \geq 2$. Let 
$\mathcal{C}\subset \cO_{X,o}$  be its  the conductor ideal (cf. Definition \ref{def:CONDINDEF}). 
 Then 
 \begin{equation*}
     \mathbb{SH}_*(\cO_{X,o}/\mathcal{C}) \cong \mathbb{H}_*((X, o); \Omega^n) \ (\cong \mathbb{H}_{an,*}(X, o) \text{ whenever this latter is well-defined}).
 \end{equation*}
\end{corollary}

\subsection{Reduction theorems}\label{ss:RedTh}\,

Let us recall that the lattice homology is determined by certain cubical subspaces $\{S_n\}_n$ embedded in some rectangles of $\Z^r\otimes\R$. In our present context $r$, the rank of the lattice $\Z^r$, 
is the number of valuations we use to `realize' the submodule $N$ in $M$. 

There is a great interest in reducing the rank $r$ of the lattice, if possible. Here are some reasons.

(1) Even theoretically, but definitely in concrete computations, it is hard to
identify the filtrations $\mathcal{F}(\ell)$ and $\mathcal{F}^M(\ell)$, whenever $\ell$ runs over a higher-rank lattice. E.g., in the present case (of analytic lattice cohomology of isolated singularities of dimension $n\geq 2$)
already for the first filtration we have to understand the vanishing orders of all the functions along all the 
irreducible exceptional divisors. It would be highly desirable to reduce  this set
 to some geometrically defined  `essential' set of divisors only. 

(2) By Proposition \ref{prop:homdegupperbound} we have $\bH_q(N\hookrightarrow M)=0$ for $q\geq r$, hence reduction of $r$ implies  additional vanishing. See also section \ref{s:homdim} about the (lattice) homological dimension and its properties. 
\vspace{2mm}

In the case of analytic lattice homology of $p$-forms associated with an isolated singularity $(X,o)$
of dimension $n\geq 2$ 
a possible reduction can occur as follows.

\begin{theorem}[Reduction Theorem] \label{prop:redhigh}
Let $\phi: \widetilde{X}\to X$ be a good resolution with irreducible exceptional divisors 
$\{E_v\}_{v\in\mathcal{V}}$. Consider also the $\mathbb{C}$-algebra $\calO$, modules $M^p$ and $N^p$ (cf. Notation \ref{nota:OMN}) and the 
set of extended discrete valuations $\mathcal{D}_{\phi}=\{(\frv_v, \frv_v^M)\}_{v\in\mathcal{V}}$ as in (\ref{eq:v_vdivisorial}) and (\ref{eq:v_v^Mdivisorial}). 
Recall that we have $N^p=H^0(\widetilde{X}, \Omega^p_{\widetilde{X}})=\mathcal{F}_{\mathcal{D}}^M(0)$.

Assume, that there exists a subset $\overline{\mathcal{V}}\subset \mathcal{V}$, such that 
the collection $\overline{\mathcal{D}}=\{\frv_v\}_{v\in \overline{\mathcal{V}}}$ satisfies 
\begin{equation}\label{eq:dbar}
\mathcal{F}_{\overline{\mathcal{D}}}^M(0)=\cap_{v\in \overline{\mathcal{V}}}\ \mathcal{F}_v^M(0)=N^p.
\end{equation}
Then $\bH_*((X, o);\Omega^p)\cong\bH_*(\bR^{|\mathcal{V}|}, w_{\mathcal{D}^\natural})\cong\bH_*(\bR^{|\overline{\mathcal{V}}|}, w_{\overline{\mathcal{D}}^\natural})$ (where the superscript $^\natural $ denotes the doubling operation of Lemma \ref{lem:duplatrukk} \textit{(b)}, necessary only if $\mathcal{D}$ is not already a CDP realization). In other words, the lattice homology associated with the realization $\mathcal{D}^\natural$ (and the lattice $\Z^{|\mathcal{V}|}$) can be computed on $\Z^{|\overline{\mathcal{V}}|}$ with the weight function corresponding to 
$\overline{\mathcal{D}}^\natural$. 
\end{theorem}
\begin{proof} The statement is a direct consequence of the Independence Theorem \ref{th:IndepMod}.
\end{proof}

In the next discussion it is convenient to write $\mathcal{V}$ as the disjoint union of $\overline {\mathcal{V}}$ and $\mathcal{V}^*$. Let $\bar{r}:=|\overline {\mathcal{V}}|$ 
and  $r^*:=|\mathcal{V}^*|$. Then any lattice point 
 $\ell$ can be written as $(\overline{\ell},\ell^*)\in \Z^{\bar{r}}\times \Z^{r^*}$.
 
We wish to emphasize that when we replace $\mathcal{D}$ by $\overline{\mathcal{D}}$ 
all the functions $\frh,\ \frh^\circ $ and $  w$ will be modified, but rather compatibly with the original situation. Indeed, for any $\overline{\ell}\in \Z^{\bar{r}}$
\begin{equation} \label{eq:reducedfiltr} \mathcal{F}_{\overline{\mathcal{D}}}(\overline{\ell})=\mathcal{F}_{\mathcal{D}}((\overline{\ell},0)) \ \ \mbox{and} \ \ \
\mathcal{F}^M_{\overline{\mathcal{D}}}(\overline{\ell})=
\mathcal{F}^M_{\mathcal{D}}((\overline{\ell},\ell^*_{-\infty})),
\end{equation}
where $\ell^*_{-\infty}$ is a cycle with all entries sufficiently negative. When computing analytic lattice homology of $p$-forms of isolated singularities, the aim is to find the smallest possible set $\overline{\mathcal{V}}$ for which the Reduction Theorem works. It is thus useful to have reformulations of the assumption (\ref{eq:dbar}).

\begin{obs} For any $p\in \{n-1,n\}$, 
    a subset $\overline{\mathcal{V}} \subset \mathcal{V}$ satisfies condition (\ref{eq:dbar}) if and only if it satisfies the following
    \begin{equation}\label{eq:forms} \mbox{
if $\omega\in H^0(\widetilde{X}\setminus E,\Omega^p_{\widetilde{X}})$ has no poles along 
$\cup_{v\in \overline{\mathcal{V}}}E_v$, then $\omega$ has no poles at all. }\end{equation}
\end{obs}
In sheaf cohomological language, condition (\ref{eq:forms}) for $q=n$ reads as 
\begin{center}
    $H^0(\widetilde{X},\Omega^n_{\widetilde{X}}(\ell^*))/
H^0(\widetilde{X},\Omega^n_{\widetilde{X}})=0$ for any cycle $\ell^*>0$ supported on $\mathcal{V}^*$.
\end{center}
 Then, by the combination of 
(\ref{eq:elso}) and (\ref{eq:SerreH}), 
we deduce another reformulation:
\begin{obs} \label{obs:2}
 A subset $\overline{\mathcal{V}} \subset \mathcal{V}$ satisfies condition (\ref{eq:dbar}) for $p=n$ if and only if 
  \begin{equation}\label{eq:hegy} \mbox{
  $h^{n-1}(\cO_{\ell^*})=0$ for any cycle $\ell^*>0$ supported on $\mathcal{V}^*$.  }\end{equation}
\end{obs}

\begin{remark}\label{rem:badvertices}
(a) In \cite[Theorem 4.7.5]{AgNeHigh}, for the case $p=n$, the 
Reduction Theorem \ref{prop:redhigh} was proved under the following assumption: 
  \begin{equation}\label{eq:formsregi} 
  \mbox{
if 
 $\omega\in H^0(\widetilde{X}\setminus E,\Omega^n_{\widetilde{X}})$ 
satisfies ${\rm div}_E(\omega)\big|_{\overline{\mathcal{V}}}\geq -\sum_{v\in \overline{\mathcal{V}}}E_v$, then $\omega$ has no poles. }\end{equation}
Note that our present assumption (\ref{eq:forms}) only requires $\omega$ not to have poles at all, if 
 ${\rm div}_E(\omega)\big|_{\overline{\mathcal{V}}}\geq 0$.
In particular, Theorem \ref{prop:redhigh} is more general than \cite[Theorem 4.7.5]{AgNeHigh}. 

(b) Assume that $n=2$. Though the equivalent conditions  (\ref{eq:forms}) and (\ref{eq:hegy}) look analytical, in this case they can be reformulated in a  topological language 
(hence they can be tested on the resolution graph). 

First, recall that a singularity is rational if $p_g=0$. By the Formal Function Theorem this happens if and only if 
in any resolution with exceptional divisor $E$ one has 
$h^1(\cO_Z)=0$ for any $Z>0$ supported on $E$. By Artin \cite{Artin62,Artin66}
this can be reformulated in terms of the resolution graph as well (namely: $-( Z, Z-Z_K)>0$ for any effective cycle $Z>0$). If a graph satisfies 
the corresponding combinatorial condition we say that it is a rational graph (for a 
different combinatorial criterion see \cite{Laufer72}). 

Then condition (\ref{eq:hegy}) for $n=2$ is equivalent with the fact that
all the connected components of the  full subgraph supported by $\mathcal{V}^*$ are rational. This topological criterion was already used in the literature: by
Definition 7.3.21
of \cite{NBook}, 
a subset  $\overline{\mathcal{V}}$ is a \emph{WR-set} (`weak reduction set')
if it satisfies exactly this condition 
(i.e., by deleting all the vertices  $\overline{\mathcal{V}}$ and the adjacent edges 
we remain with a rational graph). 

Hence, if $p=n=2$ then the Reduction Theorem \ref{prop:redhigh} holds for any WR-set 
$\overline {\mathcal{V}}$.

(c) Recall, that in the case of $n=2$ any normal surface singularity satisfies $q\leq p_g$.
This means that if $\overline{\mathcal{V}}$ is a WR-set, i.e., 
$\mathcal{V}^*$ supports only rational singularities, then condition 
(\ref{eq:forms})  is valid for $p=1$, too. Hence 
the Reduction Theorem \ref{prop:redhigh} holds for any WR-set 
$\overline {\mathcal{V}}$ for $n=2$ and $p=1$ as well. 

(d) 
In the $n=2$ case a  typical example  for a WR-set 
is when we choose  
 $\overline{\mathcal{V}}$ as  the set of nodes (vertices with valency $\geq 3$ or genus $>0$). 
For example, if  we consider the star shaped graph of the minimal good resolution of a weighted homogeneous singularity, then the 
unique central vertex forms a WR-set, hence 
$\bar{r}=1$ and  $\bH_{\geq 1}((X, o); \Omega^p)=0$.  
For a generalization of higher dimensional weighted homogeneous singularities see \cite[section 6]{AgNeHigh}.
\end{remark}

\begin{remark}
    Clearly, in the setting of the Reduction Theorem \ref{prop:redhigh} for every subset $\overline{\mathcal{V}}\subset \mathcal{V}$ satisfying condition (\ref{eq:dbar}) (for $p\in \{n-1, n\}$), by Proposition \ref{prop:homdegupperbound}, the cardinality 
    \begin{center}
        $|\overline{\mathcal{V}}|-1 \geq {\rm homdim}((X, o); \Omega^p):=
            \begin{cases}
                -1 & \text{ if }eu(\mathbb{H}_*)=0;\\
                \max\{q\,:\,\mathbb{H}_{q}((X, o); \Omega^p)\neq 0 \} & \text{ otherwise.}
            \end{cases}$
    \end{center}
\end{remark}

\begin{remark}
    In the case of the analytic lattice homology of reduced curve singularities (presented in section \ref{s:deccurves}), there is no natural stimulus for reducing the rank of the underlying lattice. Indeed, its generators correspond bijectively to irreducible components of the singularity. Nevertheless, in some specific cases one can produce lower dimensional lattices with weight functions yielding the same lattice homology module, e.g., in the Newton nondegenerate plane curve case. In fact, in this case one can use monomial valuations coming from the convex geometry of the Newton polygon to obtain `minimal' realizations of the conductor ideal (cf. Theorem \ref{th:homdimmon}). In a sense, one could consider this result the corresponding Reduction Theorem in the $n=1$ dimensional case.
\end{remark}

\subsection{Example: the weighted homogeneous case}\label{ex:whq}\,

In this subsection we will follow through the construction of the analytic lattice homology of $p$-forms in the more restricted case of Gorenstein weighted homogeneous normal surface singularities ($n=2, \ p \in \{1, 2\}$). We emphasize, however, that  we  impose neither any additional analytic-- (being hypersurface or complete intersection) nor topological restriction (regarding, e.g., the link of the singularity). We will compute the Hilbert functions $\mathfrak{h}$ and $\mathfrak{h}^\circ$ in terms of the coefficients of the Poincaré series and reason that, similarly to $\mathbb{H}_{an, *}(X, o)$, the analytic lattice homology $\mathbb{H}_*((X, o); \Omega^1)$ of $1$-forms can often be nontrivial. Finally, we will provide a concrete example.

\bekezdes \textbf{Weighted homogeneous singularities \cite{OW1, OW2}.}
An isolated singularity $(X,0)\subset (\mathbb{C}^N, 0)$ is called \emph{weighted homogeneous with a good 
$\C^*$-action} if it satisfies the following conditions.
First, we fix an action of $\C^*$ on the affine space $\C^N$ with
positive integer weights  $(w_1,\ldots , w_N)$ by 
$$t\cdot (x_1,\ldots, x_N):=
(t^{w_1}x_1,\ldots , t^{w_N}x_N).$$
As usual, we also assume that ${\rm gcd}\{w_i\}_i=1$. The singularity $(X, 0)$ is weighted homogeneous with these weights, if the ideal of functions vanishing on  $X\subset \C^N$ is generated by weighted homogeneous polynomials. We  assume that 
the origin is the only singular point of the affine variety $X$. 

For such a singularity, the affine coordinate ring, and also the local algebra $\calO=\calO_{X,0}$ is $\Z_{\geq 0}$ graded, i.e., $\calO=\oplus_{\ell\geq 0}R_l$, 
with $R_0=\C$. We will denote the Poincar\'e series of this graded ring by
$P_X(t):=\sum_{\ell\geq 0}\dim_{\mathbb{C}} R_\ell\, t^\ell=:\sum_{\ell\geq 0}d_\ell\, t^\ell.$
In several cases (e.g., if $X$ is a complete intersection or a surface singularity with 
$\mathbb{Q}HS^3$ link) $P_X(t)$ can be combinatorially determined, see, e.g., 
 \cite[paragraph 5.1.26 and Example 9.7.24]{NBook}. 

\bekezdes \textbf{Analytic lattice homology in the Gorenstein surface case.}
Let us now assume for simplicity that the isolated weighted homogeneous normal surface singularity $(X,0)$ is additionally Gorenstein  (however, several steps of the following argument can be  extended  to the non-Gorenstein and higher dimensional cases as well). For some basic results about Gorenstein weighted homogeneous singularities see, e.g., 
\cite{NBook,wahl,Yau2}. We also assume that $(X, o)$ is not rational (otherwise $p_g=q=0$).

Our goal is to determine the analytic lattice homologies of $p$-forms $\mathbb{H}_*((X, 0); \Omega^p)$ in this setting. We use the previous notations 
$\calO=H^0(\widetilde{X}, \mathcal{O}_{\widetilde{X}})$,  $M^p=H^0(\widetilde{X}\setminus E, \Omega^p_{\widetilde{X}})$ and $N^p=H^0(\widetilde{X}, \Omega^p_{\widetilde{X}})$ for the specific minimal good resolution $\phi: \widetilde{X} \to X$ (with reduced exceptional divisor $E$), obtained by resolving the quotient singularities of the exceptional curve of the weighted blow up.
In the case $p=2$ several particular examples of the analytic lattice homology after \cite{AgNe1} are known (consult, e.g., \cite{NBook}), in particular we know that the analytic lattice (co)homology of $2$-forms is usually nontrivial. Our goal is to show the same fact for 
$p=1$ by connecting  the modules $M^2$ and $M^1$ via the Euler vector field and by producing concrete computations for both. 

For the computation of  the analytic lattice homology modules we use the Reduction Theorem   \ref{prop:redhigh}: by construction, the dual resolution graph is star-shaped (with all the legs rational), hence 
we need only one valuation, namely the one associated with the `central' irreducible 
exceptional divisor $E_0$, cf. Remark \ref{rem:badvertices} (d). In fact, $E_0$ is the base of the associated orbifold Seifert bundle, or the exceptional curve of the weighted blow up.
The point is that the direct sum decompositions of the local algebra $\mathcal{O}$ and of the modules $M^1$ and $M^2$ provided by the 
$\C^*$-action 
coincide with  the graded vector spaces associated with the $E_0$-divisorial filtration (see, e.g., \cite[Example 8.6.13]{NBook}). 
More precisely, we have the identity $\mathcal{F}^M_{E_0}(-\ell)=H^0(\Omega^p_{\widetilde{X}}(\ell E_0))$ (where, for the sake of simplicity, the collection consisting of the single extended divisorial valuation $\mathfrak{v}_{E_0}$ corresponding to $E_0$ is also denoted by $E_0$) and
$$R_\ell\simeq \mathcal{F}_{E_0}(\ell)/ \mathcal{F}_{E_0}(\ell+1), 
\ \ \ \
M^p\simeq\oplus _{\ell}\ M^p_{-\ell}=\oplus_{\ell}\ \mathcal{F}^M_{E_0}(-\ell)/ \mathcal{F}^M_{E_0}(-\ell+1).$$
\bekezdes \label{bek:conductorofwhGorss}  Let $c$ be the $E_0$-coefficient of $Z_K=-c_1(\Omega^2_{\widetilde{X}})$ (i.e., $c={\rm mult}_{E_0}(Z_K)$). Recall, that $Z_K$ is the largest pole of a 2-form: it is the pole of the Gorenstein form $\omega_G$. All the other forms of $M^2$ 
can be realized as $f\cdot\omega_G$ for some $f\in\calO$, i.e., $M^2$ is free of rank one (cf. subsection \ref{ss:zcoh}). Since $(X ,o)$ is not rational, we have $c\geq 1$.
Then the `Gorenstein duality' reads as follows:
$$M^2_{-\ell}\simeq R_{c-\ell} \cdot  \omega_G\simeq R_{c-\ell}, \ \ \ \text{ and }  N^2=\oplus _{\ell\leq 0}\ 
M^2_{-\ell}\simeq \oplus _{\ell\geq c}\, R_{\ell}\cdot \omega_G\simeq \oplus _{\ell\geq c}\, R_{\ell}. $$
In particular, $p_g=\dim M^2/N^2=\sum_{\ell=0}^{c-1}d_\ell$ (where $d_\ell=\dim_{\mathbb{C}}R_\ell$).
Therefore,  in the $p=2$ case
\begin{equation*}\begin{split}
\frh(\ell)=&\ \dim\calO/\mathcal{F}_{E_0}(\ell)=\sum_{i=0}^{\ell-1}d_i,\\
\frh^\circ_{p=2} (\ell)=&\ \dim \, M^2/\mathcal{F}_{E_0}^{M^2}(-\ell)=\sum _{\ell'> \ell}\dim\, M^2_{-\ell'}=\sum _{\ell'>\ell}d_{c-\ell'}=
\sum _{i<c-\ell}d_i \\
=&\ \frh(c-\ell)=\frh^{sym}_c(\ell), \text{ hence }
\\
\frh^\circ_{p=2} (\ell)-\frh^\circ_{p=2} (0)=&\ -\sum _{c-\ell\leq i<c}\, d_i.
\end{split}
\end{equation*}
Usually, this pair $(\frh,\frh^\circ_{p=2})$ does not satisfy the CDP, however, if we replace the valuation $\frv_{E_0}$ by $2\frv_{E_0}$ (as in  Lemma \ref{lem:duplatrukk}), we get a pair with CDP, hence a weight function which computes the analytic lattice homology of $2$-forms. (For a concrete computation see Example \ref{ex:555} below.)
By Theorem \ref{th:properties} \textit{(a)}, this $\mathbb{Z}[U]$-module depends only on the weighted cubes of the rectangle $R(0,2c)$.

\bekezdes \textbf{Analytic lattice homology of $1$-forms.}
The connecting bridge between $M^2$ and $M^1$ in the Gorenstein weighted homogeneous normal surface singularity setup is realized by the Euler vector field:
the presence of the $\C^*$-action provides automatically the holomorphic vector field
$\xi:= \sum _i w_ix_i \frac{\partial}{\partial x_i}$
on $(X,0)$.
This is the tangent vector field of the parametrized orbits of the action. (For more details on this see, e.g., \cite{Yau2}.)
 With regards to the minimal good resolution $\phi: \widetilde{X} \rightarrow X$ from the previous paragraph, the Euler vector field $\xi$ lifts to a vector field $\widetilde{\xi}$  on $\widetilde{X}$, with $\phi_*(\widetilde{\xi})=\xi$. Let $x$ be a generic point of $E_0$,
and $U_x$ a sufficiently small neighborhood of $x$ with local analytic coordinates $(u,v)$,
such that the family of orbits intersected with  $U_x$ are given by 
$\{v=\mbox{constant}\}$, whereas $\{u=0\}=E_0\cap U_x$.
Then one sees that $\widetilde{\xi}|_U=u\frac{\partial}{\partial u}$. 

Now let $\iota _{\widetilde{\xi}}:\Omega^2_{\widetilde{X}}\to 
\Omega^1_{\widetilde{X}}$ denote the contraction (inner derivation) by the vector field
$\widetilde{\xi}$. In the neighbourhood $U_x$ of the point $x \in E_0$ in the convenient local coordinate system 
\begin{center}
    $\iota_{\widetilde{\xi}}(u^{-\ell}a(u,v)du\wedge dv)=u^{-\ell+1}a(u,v)dv$.  
\end{center}
Hence it 
induces a homogeneous linear map  $\iota_{\widetilde{\xi}}:M^2_{-\ell}\to 
M^1_{-\ell+1}$. It turns out that this is an isomorphism for $\ell >0$ (compare Yau's Lemma from \cite[pp. 834]{Yau2} with \cite[Corollary 1.11]{wahl}).
Therefore, we have $M^2/M^2_{\geq -1} \cong M^1/N^1$, thus by Theorem \ref{th:properties} \textit{(c)} 
\begin{center}$\mathbb{H}_*((X, 0); \Omega^1) \cong \mathbb{H}_*(M^2_{\geq -1} \hookrightarrow_{\mathcal{O}}M^2)$ and  $q=\dim \, M^1/N^1=\sum_{\ell=0}^{c-2}d_\ell$.
\end{center}

Note that $p_g-q=\dim\, M^2_{-1}=d_{c-1}$. This value can be computed with the help of 2-forms and Serre duality. 
Indeed, by (\ref{eq:SerreH}) we have 
\begin{center}
    $d_{c-1}=\frh(c)-\frh(c-1)=\frh^\circ_{p=2}(0)-\frh^\circ_{p=2}(1)=h^1(\calO_{E_0})=\mbox{genus}(E_0)$.
\end{center}
Therefore: $q=p_g-\mbox{genus}(E_0)$ for Gorenstein weighted homogeneous normal surface singularities, as already proved in \cite{wahl}. 

In conclusion, in this $p=1$ setup, the function $\frh$ clearly is the same, however, 
using the identification  $\iota_{\widetilde{\xi}}:M^2_{-\ell}\to 
M^1_{-\ell+1}$, for $\frh^\circ _{p=1}$ we get 
\begin{equation*}\begin{split}
\frh^\circ_{p=1} (\ell)=&\ \dim \, M^1/\mathcal{F}_{E_0}^{M^1}(-\ell)=\sum _{\ell'> \ell}\dim\, M^2_{-\ell'-1}=\sum _{\ell'>\ell}d_{c-1-\ell'}=
\sum _{i<c-1-\ell}d_i \\
=&\ \frh(c-1-\ell)=\frh^{sym}_{c-1}(\ell), \text{ hence}
\\
 \frh^\circ_{p=1} (\ell)-\frh^\circ_{p=1} (0)=&\ -\sum _{c-1-\ell\leq i<c-1}\, d_i.
\end{split}
\end{equation*}
Again, the  pair $(\frh,\frh^\circ_{p=1})$  associated with $M^1$ 
does not (in general) satisfy the CDP, but if we multiply the valuation by 2, then
we get a pair with the CDP, so the corresponding weight function
computes the analytic lattice homology of $1$-forms. 
Similarly as before, $\mathbb{H}_*((X, 0); \Omega^1)$ only depends on the weights inside the rectangle $R(0,2c-2)$.
It is slightly surprising that the corresponding weight function is also symmetric and that its dualizing element is $c-1$.

\begin{remark}\label{rem:whomrathomslink}
If the link is a rational homology sphere, then $h^1(\calO_E)=0$, hence $p_g=q$ and 
$$\bH_*(N^1\hookrightarrow M^1)=\bH_* (N^2\hookrightarrow M^2) \hspace{5mm} (\text{since }M^2/N^2 \cong M^1/N^1).$$ 
(Notice, however, that the corresponding $\mathfrak{h}^\circ$ functions do not agree, the extensions of $\mathfrak{v}_{E_0}$ to $M^2$ and $M^1$ differ in their degree by 1. For a similar `translation phenomenon' see Example \ref{ex:eltolasplurig} (c). Nevertheless, here it does not cause a problem due to the generality of the Independence Theorem \ref{th:IndepMod}.) Moreover, by Theorem \ref{th:newvsoldhighdim}, they also agree with $\mathbb{H}_{an, *}(X, o)$, which is usually nontrivial. 
\end{remark}

\begin{example}\label{ex:555}
As a concrete example of the Gorenstein weighted homogeneous normal surface singularity let us consider the hypersurface singularity $(\{x_1^d+x_2^d+x_3^d=0\}, 0)$ in $(\C^3,0)$, $d\geq 3$.
Its minimal embedded resolution can be achieved by a single embedded blow-up, with irreducible exceptional divisor $E_0$ having 
self-intersection $-d$ and genus $(d-1)(d-2)/2$.  Besides, 
$Z_K=cE_0=(d-2)E_0$ by the adjunction formula. 
Then the Poincar\'e series $P_X(t)$ is given by $P_X(t)=
(1-t^d)/(1-t)^3$, hence for $i\leq c=d-2$ one has $d_i=\binom{i+2}{2}$. Thus $p_g=d(d-1)(d-2)/6$ and 
$q=(d-1)(d-2)(d-3)/6$.

For an even more concrete example, assume that  $d=5$ and fix $p=2$. Then the  $\frh(\ell)$ values 
are $0, 1,4,10, 20,\ldots$ and  the 
$\frh^\circ_{p=2}(\ell) $ values are  $10, 4, 1, 0, 0, \ldots$ for $\ell=0,1, 2, 3, \ldots$ with $c=3$. Since $p_g=10$, the res\-pective
$\frh^\circ_{p=2}(\ell)-\frh^\circ_{p=2}(0) $ values are  $0, -6, -9, -10, -10, \ldots$. 
If we consider the doubled valuation, then 
the $\frh^\natural$ values are $0, 1,1,4,4,10,10, 20,20, \ldots$ and the 
$\frh^{\circ\,\natural}_{p=2}-\frh^{\circ \, \natural}_{p=2}(0)$ 
va\-lues are $0, 0, -6, -6, -9, -9, -10, -10, -10, \ldots$, hence 
the weights $w_0^\natural $ are in order 
$$ 0, 1, -5, -2, -5, 1, 0, 10, 10, \ldots.$$
The graded root  is shown on the left hand side below. Clearly, $\bH_{\geq 1}=0$.

For $d=5$ and  $p=1$  the $\frh$ values are the same, 
but the $\frh^{\circ}_{p=1}(\ell)$ values are $4,1,0,0,0,\ldots $. Since $q=4$, the 
$\frh^\circ_{p=1}(\ell)-\frh^\circ_{p=1}(0) $ values are  $0, -3, -4, -4,\ldots$. 
If we double the valuation, then 
the corresponding $\frh^\natural$ values are $0, 1,1,4,4,10,10,  \ldots$ and the 
$\frh^{\circ \, \natural}_{p=1}-\frh^{\circ\,\natural}_{p=1}(0)$ 
values are $0, 0, -3, -3, -4, -4, -4,  \ldots$, hence 
the weights $w_0^\natural $ for $p=1$ are 
$$ 0, 1, -2, 1, 0, 6, 6,  \ldots.$$ 
The graded root  is shown on the right  hand side below. 

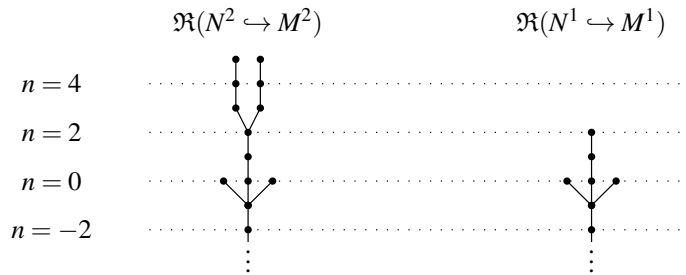
\begin{center}
\resizebox{13cm}{!}{
\begin{picture}(400,110)(0,10)
\linethickness{.5pt}

\put(125,90){\circle*{3}}
\put(135,90){\circle*{3}}
\put(125,80){\circle*{3}}
\put(135,80){\circle*{3}}
\put(125,70){\circle*{3}}
\put(135,70){\circle*{3}}
\put(130,60){\circle*{3}}
\put(130,50){\circle*{3}}
\put(120,40){\circle*{3}}
\put(130,40){\circle*{3}}
\put(140,40){\circle*{3}}
\put(130,30){\circle*{3}}
\put(130,30){\circle*{3}}
\put(130,20){\circle*{3}}

\put(130,10){\makebox(0,0){$\vdots$}}

\qbezier[50](90,80)(200,80)(310,80)
\qbezier[50](90,60)(200,60)(310,60)
\qbezier[50](90,40)(200,40)(310,40)
\qbezier[50](90,20)(200,20)(310,20)

\put(50,80){\makebox(0,0)[0]{$n=4$}}
\put(50,60){\makebox(0,0)[0]{$n=2$}}
\put(50,40){\makebox(0,0)[0]{$n={0}$}}
\put(50,20){\makebox(0,0)[0]{$n=-2$}}

\put(125,90){\line(0,-1){20}}
\put(135,90){\line(0,-1){20}}
\put(125,70){\line(1,-2){5}}
\put(135,70){\line(-1,-2){5}}
\put(130,60){\line(0,-1){45}}
\put(120,40){\line(1,-1){10}}
\put(140,40){\line(-1,-1){10}}

\put(270,60){\circle*{3}}
\put(270,50){\circle*{3}}
\put(260,40){\circle*{3}}
\put(270,40){\circle*{3}}
\put(280,40){\circle*{3}}
\put(270,30){\circle*{3}}
\put(270,30){\circle*{3}}
\put(270,20){\circle*{3}}

\put(270,60){\line(0,-1){45}}
\put(260,40){\line(1,-1){10}}
\put(280,40){\line(-1,-1){10}}

\put(270,10){\makebox(0,0){$\vdots$}}

\put(130,105){\makebox(0,0)[0]{$\mathfrak{R}(N^2 \hookrightarrow M^2)$}}
\put(270,105){\makebox(0,0)[0]{$\mathfrak{R}(N^1 \hookrightarrow M^1)$}}

\end{picture}}
\end{center}
\end{example}

\subsection{Categorification of the $L^2$-plurigenera of isolated singularities}\label{s:pluri}\,

In the following subsections we will present the various versions of  plurigenera of normal surface singularities introduced by 
 Watanabe \cite{Watanabe80}, Knöller \cite{Kno73} and Morales \cite{Mor83}, and show that they can be computed as the codimension of some realizable submodules. As a corollary, we can construct lattice homological categorifications of these invariants via the construction of subsection \ref{ss:lathomofmods}. The authors already introduced a similar concept for Gorenstein normal surface singularities in \cite{NSplurig}, however, the present version is more general and seems better suited to give interesting new invariants (see also part (c) of Example \ref{ex:eltolasplurig}).

Let $(X, o)$ be a complex analytic normal isolated singularity of dimension $n \geq 2$. 
Let $\phi: \tX \rightarrow X$ be a   good resolution with reduced exceptional divisor $E=\phi^{-1}(o)_{red}$. Set also $U:=X \setminus\{o\}$ and consider the sheaf $\overline{\overline{\Omega}}^n_X:=\iota_{\ast}\Omega_U^n$ defined in subsection \ref{ss:iso}. It agrees with the Grothendieck dualizing sheaf and it is divisorial (cf. \cite[pp. 283]{Reid}), hence $\overline{\overline{\Omega}}^n_X\cong \mathcal{O}_{X}(K_X)$ for some `canonical divisor'  $K_X$.
For any positive integer $m\in \mathbb{Z}_{>0}$, we call the sections of the sheaf $\mathcal{O}_{X}(mK_X)$  \emph{`$m$-ple holomorphic $n$-forms'}, since above $U$ we have $(\Omega^n_X)^{\otimes m}\big|_U \cong \mathcal{O}_{X}(mK_X)\big|_U$. An $m$-ple holomorphic $n$-form $\omega$ is said to be \emph{$L^{2/m}$-integrable} if $\int_{W \setminus \{o\}}(\omega \wedge \overline{\omega})^{1/m} <\infty$ for any sufficiently small relatively compact neighbourhood $W$ of $o\in X$. The submodule of $L^{2/m}$-integrable forms on $U$ is denoted by $L^{2/m}(U)\leq H^0(U, \mathcal{O}_X(mK_X))$. Watanabe in \cite{Watanabe80} defined the \emph{$m$-th $L^2$-plurigenus of the singularity $(X, o)$} as the codimension of this submo\-dule $\delta_{m}(X, o)=\dim_{\mathbb{C}} H^0(U, \mathcal{O}_{X}(mK_X))/L^{2/m}(U)$ for every $m \in \mathbb{Z}_{>0}$. It is independent of the choice of the Stein neighbourhood of $o$, since, by a theorem of Sakai \cite[Theorem 2.1]{Sakai} and the argument of Okuma in \cite[pp. 10-11]{Okuma}, we have the vector space isomorphisms
\begin{equation}\label{eq:L2plurig}
     \dfrac{H^0(\widetilde{X}\setminus E, \mathcal{O}_{\widetilde{X}}(mK_{\widetilde{X}}))}{H^0(\widetilde{X}, \mathcal{O}_{\widetilde{X}}(mK_{\widetilde{X}}+(m-1)E))} \cong \dfrac{H^0(U, \mathcal{O}_{X}(mK_X))}{L^{2/m}(U)}
 \cong \left( \dfrac{\mathcal{O}_X(mK_X)}{\phi_\ast \mathcal{O}_{\widetilde{X}}(mK_{\widetilde{X}} + (m-1) E)}\right)_o,\end{equation}
where $K_{\widetilde{X}}$ denotes the canonical divisor of the resolution space $\widetilde{X}$. This identity, together with Theorem \ref{th:properties} implies the following categorification of the $\delta_m(X, o)$ plurigenera:

\begin{cor}\label{cor:delta_mcategorification} Let $(X, o)$ be a complex analytic normal isolated singularity of dimension $n \geq 2$. Fix a  Stein representative $X$ 
and let $\phi: \widetilde{X}\rightarrow X$ denote a good resolution with reduced exceptional divisor $E$. 
Then 
\begin{align*}
\mathbb{H}_*\big( H^0(\widetilde{X}, \mathcal{O}_{\widetilde{X}}(mK_{\widetilde{X}}+(m-1)E)) \hookrightarrow_{H^0(\widetilde{X},\calO_{\widetilde{X}})} &\ H^0(\widetilde{X}\setminus E, \mathcal{O}_{\widetilde{X}}(mK_{\widetilde{X}}))\big) \cong\\
    \cong \mathbb{H}_*\big(L^{2/m}(U)  &\  \hookrightarrow_{H^0(U, \mathcal{O}_X)}H^0(U, \mathcal{O}_X(mK_X)\big) \cong \\ \cong &\ \mathbb{H}_*\left(\phi_*\mathcal{O}_{\widetilde{X}}(mK_{\widetilde{X}} + (m-1)E)_o\hookrightarrow_{\mathcal{O}_{X,o}} \mathcal{O}_{X}(mK_{X})_o\right)
\end{align*} is a well-defined invariant of the singularity and has Euler characteristic the $L^2$-plurigenus $\delta_m(X, o)$ for every $m \in \mathbb{Z}_{>0}$.
It can be computed in some finite rectangle (depending on the CDP realization chosen). We will denote it by  $\mathbb{H}_*((X, o), \delta_m)$.
\end{cor}

\begin{nota}\label{nota:OMNplurig}
    Let us fix some further notations for this subsection for the sake of better readability: $\mathcal O:=
H^0(\widetilde{X},\calO_{\widetilde{X}}),$ $ M:=H^0(\widetilde{X}\setminus E, \mathcal{O}_{\widetilde{X}}(mK_{\widetilde{X}}))$ and $N:=H^0(\widetilde{X}, \mathcal{O}_{\widetilde{X}}(mK_{\widetilde{X}}+(m-1)E))$. Then clearly $\mathbb{H}_*((X, o), \delta_m) \cong \mathbb{H}_*(N \hookrightarrow_\mathcal{O} M)$.
\end{nota}

\begin{proof}[Proof of Corollary \ref{cor:delta_mcategorification}]
    Similarly to the case of identity (\ref{eq:idekellanormalitas}), by normality and finite dimensionality of the vector spaces involved, one can prove that the composition  $$\mathcal{O}=H^0(\widetilde{X}, \mathcal{O}_{\widetilde{X}}) \xrightarrow{(\phi^*)^{-1}} H^0(X, \mathcal{O}_X) \xrightarrow{\rho_o} \mathcal{O}_{X, o}$$ induces $\mathbb{C}$-algebra isomorphisms
    \begin{align*}
    \dfrac{\mathcal{O}}{{\rm Ann}_{\mathcal{O}}(M/N)} \cong &\ \dfrac{H^0(U, \mathcal{O}_{X})}{{\rm Ann}_{H^0(U, \mathcal{O}_X)}\big(H^0(U, \mathcal{O}_{X}(mK_X))/L^{2/m}(U)\big)}\\ \cong &\ \dfrac{\mathcal{O}_{X, o}}{{\rm Ann}_{\mathcal{O}_{X, o}}\big(\mathcal{O}_{X}(mK_X)_o/\phi_*(\mathcal{O}_{\widetilde{X}}(mK_X+(m-1)E)_o\big)}.
    \end{align*}
    Moreover, the vector space isomorphisms of (\ref{eq:L2plurig}) are equivariant with respect to the action of these Artin algebras. Hence, the only remaining part is to prove that the $\mathcal{O}$-submodule $N \leq M$ is realizable (also, the same should be proved about the other two submodules occurring in the other formulations, though, one could just apply the natural identifications to obtain the desired extended valuations). For this we will use extensions of the divisorial valuations $\{\mathfrak{v}_v\}_{v\in \mathcal{V}}$ corresponding to the components $\{ E_v \}_{v \in \mathcal{V}}$ of the exceptional divisor $E$ (see (\ref{eq:v_vdivisorial})). Indeed, if $x$ is a generic point of $E_v$, and 
$U_x \subset \widetilde{X}$ is of some neighborhood of $x$ with local analytic coordinates 
$(u_1,\ldots, u_n)$ such that $U_x\cap E_v=\{u_1=0\}$, then every element $\omega\in M \setminus \{0\}$ can be written as $\omega\big|_{U_x}=a(du_1 \wedge \ldots \wedge du_n)^{\otimes m}$, where $a$ is meromorphic function of form $u_1^{\mathfrak{v}_v(a)}a'$, with $a'(x)\in \mathbb{C}^*$. Then we define 
\begin{equation}\label{eq:v_v^Mdelta_m}
    \frv_v^M:M\rightarrow \overline{\bZ}:\ \omega\ \big(\text{with }\omega\big|_{U_x} = a(du_1 \wedge \ldots \wedge du_n)^{\otimes m}\big) \mapsto \frv_v(a)-m+1 \text{ and } \omega =0 \mapsto \infty.  
\end{equation}
This is independent of all the choices. One can easily verify, that for any $v\in\mathcal {V}$  the pair $(\frv_v, \frv^M_v)$ is an extended discrete valuation.

Denote by $\mathcal{D}=\mathcal{D}_\phi:= \{(\frv_v, \frv^M_v)\}_{v\in \mathcal{V}}$ the collection of extended  discrete valuations corresponding to 
the irreducible exceptional divisors 
$\{E_v\}_{v \in \mathcal{V}}$. We claim that it gives a realization of the finite codimensional submodule $N \leq M$. Indeed, if we identify the abstract lattice $\Z^r$ ($r=|\mathcal{V}|$) with the free abelian group $\mathbb{Z}\langle \{ E_v\}_{v \in \mathcal{V}}\rangle$ of divisors supported on $E$, then the multifiltrations (defined in (\ref{eq:fdMl})) corresponding to the collection $\mathcal{D}$ satisfy that 
    \begin{align*}
        \text{ for any } \ell \in \mathbb{Z}^r: \ \mathcal{F}_{\mathcal{D}}(\ell) =&\ H^0(\widetilde{X},\calO_{\widetilde{X}}(-\ell)), \\ \text{ whereas } \mathcal{F}^M_{\mathcal{D}}(\ell)=&\ H^0(\widetilde{X}, \mathcal{O}_{\widetilde{X}}(mK_{\widetilde{X}} + (m-1)E-\ell)).
    \end{align*}
Then clearly $\mathcal{F}_{\mathcal{D}}^M(0)= N=H^0(\widetilde{X}, \mathcal{O}_{\widetilde{X}}(mK_{\widetilde{X}} + (m-1)E))$, i.e., $\mathcal{D}$ is indeed a realization. Then Lemma \ref{lem:duplatrukk} implies the existence of a CDP realization and the Independence Theorem \ref{th:IndepMod} and Theorem \ref{th:properties} yield the well-definedness and the desired properties of the corresponding lattice homology module.
\end{proof}

Once again, using the local sheaf theoretical description, we can use the Nonpositivity Theorem \ref{th:upperbound} with $k=\mathbb{C}$ and $\mathcal{O}=\mathcal{O}_{X, o}$ to conclude

\begin{cor}\label{cor:nonposfordeltam}
     Let $(X, o)$ be a complex analytic normal isolated singularity of dimension $ \geq 2$. Then for every $m \in \mathbb{Z}_{>0}$ and for every positive integer $n>0$ the corresponding $S_n$ space is contractible, hence the weight-grading of $\mathbb{H}_{{\rm red}, *}((X, o), \delta_m)$ is nonnegative, i.e.,
 $\mathbb{H}_{{\rm red},q,-2n}((X, o), \delta_m) = 0$
    for all $m > 0, \ q\geq 0$ and $n >0$. 

    Moreover, if $\delta_m(X, o)\neq 0$ (i.e., with the notations of \ref{nota:OMNplurig}: $N \neq M$), the upper bound of non-contractibility is sharp: $S_0$ is not even connected.
\end{cor}

\subsection{Categorification of the $\Delta_{m,n}$ plurigenera of normal surface singularities}\label{ss:plurig1}\,

In the normal surface singularity case (i.e., $\dim=n=2$) Knöller and Morales introduced two more plurigenera, which can be generalized even further. They can all be computed as the codimension of some realizable submodule, hence our general construction from subsection \ref{ss:lathomofmods} can be applied to obtain categorifications of these invariants as well. We will use the previous notations and the good resolution $\phi:\widetilde{X} \rightarrow X$.

The {\it plurigenera} of a given normal surface singularity $(X, o)$ are defined for each pair of integers $m \in
\Z_{>0}, n \in \mathbb{Z} \cap [-m, 0]$ as the (finite) dimension $\Delta_{m,n}(X, o)$ of the quotient vector space:
\begin{equation}\label{eq:plurigmn}
\dfrac{H^0(\widetilde{X}\setminus E, \mathcal{O}_{\widetilde{X}}(mK_{\widetilde{X}}))}{H^0(\widetilde{X}, \mathcal{O}_{\widetilde{X}}(m(K_{\widetilde{X}}+E)+nE))} \cong \left( \dfrac{\mathcal{O}_X(mK_X)}{\phi_\ast \mathcal{O}_{\widetilde{X}}(m(K_{\widetilde{X}}+E)+nE))}\right)_o.
\end{equation}
(Compare with \cite[Definition 6.8.57]{NBook}.)
For specific values of $-m \leq n \leq 0$ we get 
\begin{itemize}
    \item the $m$-th $L^2$-plurigenus (Watanabe \cite{Watanabe80}):
    $\delta_m(X, o) = \Delta_{m, -1}(X, o)$;
     \item the $m$-th plurigenus (Knöller \cite{Kno73}): $\gamma_m(X, o) = \Delta_{m, -m}$;
     \item the $m$-th log-plurigenus (Morales \cite{Mor83}):
     $\lambda_m(X, o)= \Delta_{m, 0}$.
\end{itemize}
These definitions are independent of resolutions. Indeed, it is enough to check that for a single blow-up $\Phi:\widetilde{X}' \rightarrow\widetilde{X}$, with center $x \in E$ and exceptional divisor $E_{new}\simeq \mathbb{CP}^1$, the numerator and the denominator of the left hand side of (\ref{eq:plurigmn}) agree for the pairs $(\widetilde{X}', E')$ and $(\widetilde{X}, E)$ (where $E'=\Phi^*(E)_{red}=E+E_{new}$). Since 
\begin{equation*}
    H^0(\widetilde{X}'\setminus E', \mathcal{O}_{\widetilde{X}'}(mK_{\widetilde{X}'}))\cong H^0(X \setminus \{o\}, (\Omega_{X \setminus \{o\}}^2)^{\otimes m})\cong H^0(\widetilde{X}\setminus E, \mathcal{O}_{\widetilde{X}}(mK_{\widetilde{X}}))
\end{equation*}
as $H^0(\widetilde{X}', \mathcal{O}_{\widetilde{X}'}) \cong H^0(X, \mathcal{O}_X) \cong H^0(\widetilde{X}, \mathcal{O}_{\widetilde{X}})$-modules (here we use the fact that $(X, o)$ is normal),
we only have to compare the denominators. Let $\mathcal{L}$ denote the sheaf $\mathcal{O}_{\widetilde{X}}(m(K_{\widetilde{X}}+E)+nE)$. Then we have $\Phi^*(\mathcal{L})=\Phi^*\left(\mathcal{O}_{\widetilde{X}}(m(K_{\widetilde{X}}+E)+nE)\right)=\mathcal{O}_{\widetilde{X}'}(\Phi^*(m(K_{\widetilde{X}}+E)+nE))$ and (since $\Phi^{*}(K_{\widetilde{X}})=K_{\widetilde{X}'}-E_{new}$)
\begin{equation*}
    \Phi^*(m(K_{\widetilde{X}}+E)+nE)=
    \begin{cases}
        m(K_{\widetilde{X}'}-E_{new}+E') + nE' \hspace{5mm} & \text{ if } x \text{ is a smooth point of }E;\\
         m(K_{\widetilde{X}'}+E') + nE' + nE_{new} & \text{ if } x \text{ is a singular point of }E.
    \end{cases}
\end{equation*}
Thus, if we denote by $\mathcal{L}':=\mathcal{O}_{\widetilde{X}'}(m(K_{\widetilde{X}'} + E') +nE')$, then $\Phi^*(\mathcal{L})=\mathcal{L}'(-kE_{new})$ for some $k \geq 0$. The sheaf exact sequence $0 \rightarrow \Phi^*\mathcal{L} \rightarrow \mathcal{L}' \rightarrow  \mathcal{O}_{kE_{new}}\otimes \mathcal{L}' \rightarrow0$ induces the cohomological long exact sequence $$0 \rightarrow H^0(\widetilde{X}', \Phi^*(\mathcal{L})) \rightarrow H^0(\widetilde{X}', \mathcal{L}') \rightarrow H^0(kE_{new}, \mathcal{O}_{kE_{new}}\otimes \mathcal{L}') \rightarrow \ldots,$$ where the last term is $0$, since $E_{new} \simeq \mathbb{CP}^1$ and the intersection $(kE_{new}, m(K_{\widetilde{X}'}+E') + nE')_{\widetilde{X}'}\leq 0$. On the other hand, $H^0(\widetilde{X}', \Phi^*(\mathcal{L})) \cong H^0(\widetilde{X}, \mathcal{L})$ by the properties of the direct image functor, thus
\begin{center}
    $H^0(\widetilde{X}', \mathcal{O}_{\widetilde{X}'}(m(K_{\widetilde{X}'} + E') +nE'))\cong H^0(\widetilde{X}, \mathcal{O}_{\widetilde{X}}(m(K_{\widetilde{X}} + E) +nE))$
\end{center} 
as $H^0(\widetilde{X}', \mathcal{O}_{\widetilde{X}'}) \cong H^0(\widetilde{X}, \mathcal{O}_{\widetilde{X}})$-modules
and, clearly, this identification is compatible with that of the numerators above.

\begin{cor}\label{cor:Delta_m,ncategorification} Let $(X, o)$ be a normal surface singularity, and $\phi: \widetilde{X}\rightarrow X$ a good  resolution with reduced exceptional divisor $E$. 
Then 
\begin{align*}
    \mathbb{H}_*\big(H^0(\widetilde{X}, \mathcal{O}_{\widetilde{X}}(m(K_{\widetilde{X}}+E) + nE)) \hookrightarrow_{H^0(\widetilde{X}, \mathcal{O}_{\widetilde{X}})} & \ H^0(\widetilde{X}\setminus E, \mathcal{O}_{\widetilde{X}}(mK_{\widetilde{X}}))\big) \cong \\
    \cong & \ \mathbb{H}_*\big(\phi_*\mathcal{O}_{\widetilde{X}}(m(K_{\widetilde{X}}+E)+nE)_o \hookrightarrow_{\mathcal{O}_{X, o}}\mathcal{O}_X (mK_X)_o\big)
\end{align*} 
is a well-defined invariant of the singularity and has Euler characteristic the plurigenus $\Delta_{m, n}(X, o)$ for every $m \in \mathbb{Z}_{>0}, n \in \mathbb{Z} \cap [-m, 0]$.
It can be computed in some finite rectangle (depending on the CDP realization chosen). We will denote it by  $\mathbb{H}_*((X, o), \Delta_{m, n})$.
\end{cor}

\begin{proof}
    We already saw the independence of the quotient module of the resolution, whereas identity (\ref{eq:plurigmn}) also shows independence of the Stein representative chosen (it can be lifted to a module isomorphism similarly to the case of identities (\ref{eq:indepforirreg&p_g}) and (\ref{eq:idekellanormalitas})). The fact that the underlying algebra actions are equivalent can be checked analogously to the case of identity (\ref{eq:idekellanormalitas}). The realizability of the submodule can be proved similarly as in the case of Corollary \ref{cor:delta_mcategorification}: by using the divisorial valuations and their natural extensions (cf. (\ref{eq:v_v^Mdelta_m})) translated corresponding to the value of $n \in \mathbb{Z}_{\leq 0}$. Finally, the Independence Theorem \ref{th:IndepMod} and Theorem \ref{th:properties} yield the desired well-definedness and properties.
\end{proof}

Relying on the local sheaf theoretic formulation we obtain the following version of the Nonpositivity Theorem \ref{th:upperbound}:
\begin{cor}\label{cor:nonposforDeltamn}
     Consider the normal surface singularity $(X, o)$ and set $m \in \mathbb{Z}_{>0}$ and $n \in \mathbb{Z} \cap[-m, 0]$. Then the weight-grading of $\,\mathbb{H}_{{\rm red}, *}((X, o), \Delta_{m, n})$ is nonnegative, i.e.,
 $\mathbb{H}_{{\rm red},q, -2n'}((X, o), \Delta_{m, n}) = 0$
    for all $q\geq 0$ and $n' >0$. In fact, even more is true: the corresponding $S_{n'}$-spaces are contractible for every $n' >0$.
    Moreover, if $\Delta_{m,n}(X, o)\neq 0$, the upper bound of non-contractibility is sharp: $S_0$ is not even connected.
\end{cor}

\begin{remark}\label{rem:gorsymforplurig}
    In the Gorenstein case $H^0(\widetilde{X}\setminus E, \mathcal{O}_{\widetilde{X}}(mK_{\widetilde{X}})) \cong H^0(\widetilde{X}, \mathcal{O}_{\widetilde{X}})$ (see  \cite[Lemma 6.1.1]{NSplurig})  and, hence, the weight function of the lattice homology corresponding to any realization is symmetric (cf. Proposition \ref{prop:symforideals}). Moreover, in this case $K_{\widetilde{X}}=-Z_K$ and, hence, $\Delta_{m, n}=\mathfrak{h}(mZ_K-(m+n)E)$, where $\mathfrak{h}$ denotes the `standard' Hilbert function of Remark \ref{rem:filtrationreformulation}. Even more, 
    \begin{equation*}
        \mathbb{H}_*((X, o), \Delta_{m,n})\cong \mathbb{SH}_*\big(H^0(\widetilde{X}, \mathcal{O}_{\widetilde{X}}(-mZ_K + (m+n)E)) \triangleleft H^0(\widetilde{X}, \mathcal{O}_{\widetilde{X}})\big).
    \end{equation*}
\end{remark}

\begin{remark}\label{rem:RedThmforplurig}
    The Reduction Theorem \ref{prop:redhigh} could be formulated for $\mathbb{H}_*((X, o),\Delta_{m,n})$ as well, however, the characterization of the `essential' subsets $\overline{\mathcal{V}}$ of the divisorial valuations (needed to realize $H^0(\widetilde{X}, \mathcal{O}_{\widetilde{X}}(m(K_{\widetilde{X}}+E)+nE))$) seems much harder. In particular, any topological characterization (especially with the notions like WR--, R-- or SR-sets from \cite[Definition 7.3.21]{NBook} or Remark \ref{rem:badvertices}) seems obstructed by the fact that even rational singularities might have nontrivial higher plurigenera: see, e.g., \cite[Example 6.8.64 and 6.8.65]{NBook}. This means that the computation of this new invariant is rather challenging in general. 
\end{remark}

In the sequel we present some calculations for a surface singularity whose minimal resolution graph consists of a single vertex only.

\begin{example} \label{ex:plurig1}
    We go back to Example \ref{ex:555} and consider the homogeneous hypersurface singularity of $X=\{x_1^5+x_2^5+x_3^5=0\}$ at $0$ in $(\C^3,0)$. We use the notations introduced there and will compute $\mathbb{H}_*((X, 0), \delta_2)$. As before, the  $\frh(\ell)$ values for $\ell \geq 0$
are $0, 1,4,10, 20, 35\ldots$.
    
    On the other hand, $mK_{\widetilde{X}} + (m-1)E=-5E_0$, hence, (since $(X, 0)$ is Gorenstein), the 
$\frh^\circ(\ell) $ values for $\ell \geq 0$ are  $35, 20, 10, 4, 1, 0, 0, \ldots$ and $\delta_2=35$.
If we double the valuation, then 
the corresponding $\frh^\natural$ values are $0, 1,1,4,4,10,10, 20,20, 35, 35,  \ldots$ and the 
$\frh^{\circ\,\natural}$ 
values are $35, 35, 20, 20, 10, 10, 4, 4, 1, 1, 0 \ldots$, hence 
the weights $w_0^\natural $ are $$ 0, 1, -14, -11, -21, -15, -21, -11, -14, 1, 0 \ldots.$$ 
The graded root  is shown  below. Clearly, $\bH_{\geq 1}=0$. 

\begin{center}
\begin{picture}(200,130)(0,-20)
\linethickness{.5pt}

\put(125,90){\circle*{2}}
\put(135,90){\circle*{2}}
\put(125,86){\circle*{2}}
\put(135,86){\circle*{2}}
\put(135,82){\circle*{2}}
\put(125,82){\circle*{2}}
\put(135,78){\circle*{2}}
\put(125,78){\circle*{2}}
\put(135,74){\circle*{2}}
\put(125,74){\circle*{2}}
\put(125,70){\circle*{2}}
\put(135,70){\circle*{2}}
\put(130,66){\circle*{2}}
\put(120,62){\circle*{2}}
\put(130,62){\circle*{2}}
\put(140,62){\circle*{2}}
\put(120,58){\circle*{2}}
\put(130,58){\circle*{2}}
\put(140,58){\circle*{2}}
\put(120,54){\circle*{2}}
\put(130,54){\circle*{2}}
\put(140,54){\circle*{2}}
\put(130,50){\circle*{2}}
\put(130,46){\circle*{2}}
\put(130,42){\circle*{2}}
\put(130,38){\circle*{2}}
\put(130,34){\circle*{2}}
\put(130,30){\circle*{2}}
\put(130,26){\circle*{2}}
\put(130,22){\circle*{2}}
\put(130,18){\circle*{2}}
\put(130,14){\circle*{2}}
\put(130,10){\circle*{2}}
\put(120,6){\circle*{2}}
\put(130,6){\circle*{2}}
\put(140,6){\circle*{2}}
\put(130,2){\circle*{2}}
\put(130,-2){\circle*{2}}

\put(130,10){\makebox(0,0){$\vdots$}}

\qbezier[20](90,86)(130,86)(170,86)
\qbezier[20](90,66)(130,66)(170,66)
\qbezier[20](90,46)(130,46)(170,46)
\qbezier[20](90,26)(130,26)(170,26)
\qbezier[20](90,6)(130,6)(170,6)

\put(50,86){\makebox(0,0)[0]{$n=20$}}
\put(50,66){\makebox(0,0)[0]{$n=15$}}
\put(50,46){\makebox(0,0)[0]{$n=10$}}
\put(50,26){\makebox(0,0)[0]{$n=5$}}
\put(50,6){\makebox(0,0)[0]{$n=0$}}

\put(125,90){\line(0,-1){20}}
\put(135,90){\line(0,-1){20}}
\put(125,70){\line(5,-4){5}}
\put(135,70){\line(-5,-4){5}}
\put(130,66){\line(0,-1){72}}
\put(120,6){\line(10,-4){10}}
\put(140,6){\line(-10,-4){10}}
\put(120,62){\line(0,-1){8}}
\put(140,62){\line(0,-1){8}}
\put(120,54){\line(10,-4){10}}
\put(140,54){\line(-10,-4){10}}

\put(130,-10){\makebox(0,0){$\vdots$}}

\put(130,100){\makebox(0,0)[0]{$\mathfrak{R}((X, 0), \delta_2)$}}
\end{picture}
\end{center}
\end{example}

\begin{remark}\label{rem:newplurivsold}
    A priori, this new construction seems to be the `better' categorification of the plurigenera then the one given by the authors in \cite{NSplurig}: we can already see in Example \ref{ex:plurig1} that it deviates far from the analytic lattice homology after \cite{AgNe1} (or, more generally, of $2$-froms); in sharp contrast with the previous construction (cf. \cite[Corollary 6.3.3 and Theorems 8.2.1 and 8.5.3]{NSplurig}, see also the discussion in part (c) of Example \ref{ex:eltolasplurig} in the present note). Moreover, it is also well-defined for every normal surface singularity, independent of the resolution used for its computation and altogether seems to be a finer invariant. Its connections with singularity deformations and geometric reinterpretations of $\delta_m, \gamma_m$ and $\lambda_m$ (cf., e.g.,  \cite{Okuma}) are yet to be explored.
\end{remark}
\newpage
\section{Monomial ideals and Newton diagrams}\label{s:ND}

In this section we will produce a combinatorial lattice homology construction for integral closures and adjoints of monomial ideals in polynomial rings. We compare it to the analytic lattice homo\-logy of Newton nondegenerate hypersurface singularities with the same Newton diagram and  read off consequences regarding deformations and Rees valuations.

\subsection{The integral closure of a monomial ideal. Newton diagram in $(\mathbb{R}_{\geq 0})^m$}\label{ss:ND}\,

Let us consider a monomial ideal $\mathcal{M}$ in the polynomial ring $\cO_m=k[x_1,\ldots, x_m]$ generated by the monomials $\{x^p=x_1^{p_1}\cdots x_m^{p_m}\}_{p\in \mathfrak{S}}$
for a set
$\mathfrak{S}$ of lattice points from  $(\Z_{\geq 0})^m$. 
Any such collection of lattice  points $\mathfrak{S}$ 
defines a Newton diagram as follows. 
 Let  $N\Gamma_+$ denote the convex closure of the set
 $\bigcup_{p\in\mathfrak{S}} (p + (\mathbb{R}_{\geq 0})^m)$ --- we call it the \emph{Newton polytope}. 
 The union of compact boundary faces of $N\Gamma_+$ forms the 
 \textit{Newton boundary}, denoted by $N\Gamma$.
  We also define 
  $N \Gamma_-$ as  the real cone $\bigcup_{x \in N\Gamma}[0, x]$ with base $N\Gamma $ and vertex $0$.  
  If $\mathfrak{S}$ consists of the exponents of  a generating set of the monomial ideal $\mathcal{M}$, then we use the notations  $N\Gamma_\pm(\mathcal{M})$ and  $N\Gamma(\mathcal{M})$. 
 
 We will assume that $\mathcal{M}$ has finite codimension in $\cO_m$. 
 This means that for any index $i$ there exists at least one  monomial generator of type
 $x_i^{p_i}$.  This translates in the language of the Newton boundary into the fact that 
$N\Gamma=N\Gamma(\mathcal{M})$ is  \textit{convenient}: it  
  intersects all the coordinate axes. 

 \begin{nota}
      Below we will write  ${\bf 1}$ for the lattice point $(1,..., 1)\in \Z^m$.
 \end{nota}
 
Associated with the monomial ideal $\mathcal{M}$ we can study two, naturally arising, integrally closed ideals, the integ\-ral closure-- 
 and the adjoint ideal of $\mathcal{M}$. For the definition of the adjoint of an ideal in the general context see, e.g., \cite{LipmanAdj,HS}, 
 for geometric motivations and interpretations see  \cite{DuVal,MT} 
 (this will be used  in our examples provided by analytic germs). 
 
 \begin{proposition}\label{prop:closadj}  (See e.g., subsection 18.4 in \cite{HS})
 
(1) The integral closure  of $\mathcal{M}$ 
is the monomial ideal $\overline{\mathcal{M}}=\big(x^p\,:\,p\in N\Gamma_+(\mathcal{M})\big)\triangleleft\mathcal{O}_m$. 

(2) The adjoint ideal  of $\mathcal{M}$ is the monomial ideal ${\rm adj}(\mathcal{M})=\big(x^p\,:\,p+{\bf 1}\in N\Gamma _+(\mathcal{M})\setminus 
N\Gamma(\mathcal{M})\big)\triangleleft\mathcal{O}_m$. It equals 
${\rm adj}(\overline{\mathcal{M}})$, too. 

Both  $\overline{\mathcal{M}}$ and ${\rm adj}(\mathcal{M})$ are 
integrally closed, in fact, they can be realized --- in the sense of Corollary \ref{cor:idealREAL} or Definition \ref{def:REALforideals} --- by the monomial valuations associated with the facets of 
the Newton boundary, see below in (\ref{eq:c1}) and (\ref{eq:c2}).
\end{proposition}

 \subsection{Combinatorial symmetric lattice homologies of the Newton boundary} \label{ss:NewtonLatCoh}\,

Throughout this subsection we fix a monomial ideal $\mathcal{M} \triangleleft\mathcal{O}_m$ with convenient Newton boundary.
Consider the set  $\{F_\sigma\}_\sigma$  of $(m-1)$-dimensional faces of $N\Gamma(\mathcal{M})$.
Let $r'$ be its cardinality.
Let us denote by $\mathbf{n}_\sigma$ the 
primitive normal vector of $F_\sigma$ with positive coefficients. Then $F_\sigma$ is contained in the affine hyperplane $\{ x \in \R^m\,:\, \langle \mathbf{n}_\sigma, x \rangle = m_\sigma \}$ for some integer $m_\sigma \in \mathbb{N}$ (where $\langle\, ,\rangle$ denotes the Euclidean inner product on $\mathbb{R}^m$). Set ${\bf m}^\Gamma:=(m_1,\ldots, m_{r'})$. 
Then define the monomial  valuations
$\frv^\Gamma_\sigma:\cO_m\to \Z$ ($1\leq \sigma\leq r'$) given by the formula
    $\frv^\Gamma_\sigma\left(\sum_p a_p x^p\right) :=\min\{\langle {\bf n}_\sigma, p\rangle\,:\, a_p \neq 0\}$
(for the definition and main properties of monomial valuations see, e.g., \cite[Definition 6.1.4]{HS}). Set $\mathcal{D}:=\{\mathfrak{v}_\sigma^{\Gamma}\}_{\sigma=1}^{r'}$. 
Define also the integers 
$c^\Gamma_\sigma:=m_\sigma +1 -\langle {\bf n}_\sigma,{\bf 1}\rangle$ for all $1\leq \sigma\leq r'$.
We say that  
 $\mathbf{c}^{\Gamma}:=(c_1^{\Gamma},\ldots , c_{r'}^{\Gamma})$ is 
 the {\it `combinatorial conductor'} of $N\Gamma(\mathcal{M})$.
From Proposition \ref{prop:closadj} we get that
\begin{align}
   \overline{\mathcal{M}} \ \text{is generated by monomials } x^p &\text{ such that
} \frv^\Gamma_\sigma(x^p)\geq m_\sigma \text{ for all } \sigma, \nonumber \\
&\text{i.e., with the notations of (\ref{eq:F_D(l)}), }
\overline{\mathcal{M}}=\mathcal{F}_{\mathcal{D}}({\bf m}^{\Gamma});\label{eq:c1} \\
{\rm adj}({\mathcal{M}}) \ \text{is generated by monomials } x^p &\text{ such that
 } \frv^\Gamma_\sigma(x^p)\geq c^{\Gamma}_\sigma \text{ for all  } \sigma, \nonumber \\ 
 &\text{i.e., with the notations of (\ref{eq:F_D(l)}), }
{\rm adj}(\mathcal{M})=\mathcal{F}_{\mathcal{D}}({\bf c}^\Gamma).\label{eq:c2}
\end{align}

One also verifies that $d_{\overline{\mathcal{M}}}={\bf m}^{\Gamma}$ and 
$ d_{{\rm adj}(\mathcal{M})}={\bf c}^\Gamma$ (for the notation $d_{\mathcal{I}}$ see Remark \ref{rem:d_I}). Thus, we can use the construction of subsection \ref{ss:Indep} to obtain the following symmetric lattice homology theories.

\begin{cor}
    Given a monomial ideal $\mathcal{M}\triangleleft \mathcal{O}_m$ with convenient Newton boundary, we can associate with both integrally closed, finite codimensional ideals $\overline{\mathcal{M}}$ and ${\rm adj}(\mathcal{M})$ symmetric lattice homology modules. Theorem \ref{th:propertiesideal} \textit{(a)} implies that 
\begin{align*}
eu({\mathbb S}\bH_*(\overline{\mathcal{M}}\triangleleft \cO_m))=&\, \big| (\Z_{\geq 0})^{m}\setminus 
N\Gamma _+(\mathcal{M})\big|, \\ 
eu({\mathbb S}\bH_*({\rm adj}({\mathcal{M}})\triangleleft \cO_m))=&\, \big| (\Z_{>0})^{m}\cap  
N\Gamma _-(\mathcal{M})\big|.
\end{align*}
Moreover, since $\overline{\mathcal{M}}$ and ${\rm adj}(\mathcal{M})$ are monomial ideals, thus contained in the homogeneous maximal ideal $\mathfrak{m}_m=(x_1, \ldots, x_m)$, by Theorem \ref{th:propertiesideal} \textit{(b)} we have
\begin{equation*}
    \mathbb{SH}_*(\overline{\mathcal{M}}\triangleleft\mathcal{O}_m) \cong \mathbb{SH}_*(\mathcal{O}_m/\overline{\mathcal{M}}) \cong\mathbb{SH}_*((\mathcal{O}_m)_{\mathfrak{m}_m}/\overline{\mathcal{M}}_{\mathfrak{m}_m}) \cong \mathbb{SH}_*(\overline{\mathcal{M}}_{\mathfrak{m}_m} \triangleleft (\mathcal{O}_m)_{\mathfrak{m}_m}); 
    \end{equation*}
    \begin{equation*}
    \mathbb{SH}_*({\rm adj}(\mathcal{M})\triangleleft\mathcal{O}_m) \cong \mathbb{SH}_*(\mathcal{O}_m/{\rm adj}(\mathcal{M})) \cong\mathbb{SH}_*((\mathcal{O}_m)_{\mathfrak{m}_m}/{\rm adj}(\mathcal{M})_{\mathfrak{m}_m}) \cong  \mathbb{SH}_*({\rm adj}(\mathcal{M})_{\mathfrak{m}_m} \triangleleft (\mathcal{O}_m)_{\mathfrak{m}_m}).
\end{equation*}
Therefore, when $k=\overline{k}$ is algebraically closed we can also use Thoerem \ref{th:propertiesideal} \textit{(d)} to obtain that
\begin{align*}
    \mathbb{SH}_{{\rm red}, q, -2n}(\overline{\mathcal{M}}\triangleleft\mathcal{O}_m)=0 \text{ for all } q\geq 0 \text{ and  }n>0;\\
    \mathbb{SH}_{{\rm red}, q, -2n}({\rm adj}(\mathcal{M})\triangleleft\mathcal{O}_m)=0 \text{ for all } q\geq 0 \text{ and  }n>0.
\end{align*}
\end{cor}

\begin{remark}
  Both symmetric lattice homology modules can be computed using the set of monomial valuations $\{\frv_\sigma^\Gamma\}_\sigma$ corresponding to the $(m-1)$-dimensional faces  of the Newton boundary. Note however, that in general the corresponding height functions will not satisfy Combinatorial Duality Property, hence a priori one has to use, e.g., any of the doubling constructions of Lemma \ref{lem:duplatrukkideal} to get the correct height functions. In practice, however, it is enough to just double the valuations corresponding to the non-primitive facets of the Newton boundary. We will elaborate on this in the following discussion.
\end{remark}

\bekezdes\label{par:NewtonComb} \textbf{Combinatorial description.} The construction of the  symmetric lattice homology modules ${\mathbb S}\bH_*(\overline{\mathcal{M}}\triangleleft \cO_m)$ and ${\mathbb S}\bH_*({\rm adj}({\mathcal{M}})\triangleleft \cO_m)$ turns out to be a completely combinatorial computation depending on a finite set of lattice points determined by  the Newton boundary. 

 We will comment on and make explicit the second case of 
 ${\rm adj}(\mathcal{M})$.

Consider the following set of lattice points: $\mathcal{P}:= 
\big( (\Z_{>0})^{m}\cap  
N\Gamma _-(\mathcal{M})\big)-{\bf 1}$ (this naturally corresponds to the set of exponents of the standard free generating set of the finite dimensional $k$-vector space $\cO_m/{\rm adj}(\mathcal{M})$). In order to define a grading on it, which will have the corresponding height function satisfying the Combinatorial Duality Property, let us consider any primitive facet decomposition of the $(m-1)$-dimensional compact boundary faces of $N\Gamma_+(\mathcal{M})$ (in fact, the differences between possible primitive facet decompositions will not play a role, we will only use the equations of the hypersurfaces containing these faces). Recall that we call a polyhedron, sitting inside $\mathbb{Z}^m \otimes \mathbb{R}$, with vertices in $\mathbb{Z}^m$, {\it primitive} if it does not contain any lattice point in its relative interior. Let us denote these primitive facets by $\{F_v\}_{v=1}^r$, the corresponding primitive positive normal vectors by $\{{\bf n}_v\}_{v=1}^{r}$ and the `right hand side scalars' of the defining equations by $\{m_v\}_{v=1}^r$ (i.e., $F_v \subset \{x \in \mathbb{Z}^m \otimes \mathbb{R}\,:\, \langle x, {\bf n}_v \rangle = m_v\}$). Let us also set $c^\Gamma_v:=m_v +1 -\langle {\bf n}_v,{\bf 1}\rangle$ for all $v$ and then 
 $\mathbf{c}^{\Gamma}:=(c_1^{\Gamma},\ldots , c_{r}^{\Gamma})$.  
 Notice that, due to our primitive decomposition, it can happen that ${\bf n}_v ={\bf n}_w$ and $m_v=m_w$ for $v \neq w$. Now we define the following degree function on $\mathcal{P}$:
\begin{equation} \label{eq:degfctn}
    \deg: \mathcal{P} \rightarrow \mathbb{Z}^r, p \mapsto (\langle p, {\bf n}_1\rangle, \ldots , \langle p, {\bf n}_r\rangle).
\end{equation}
Then, for any lattice point $p\in (\Z_{\geq 0})^{m}$
we have the characterization: $p\in \mathcal{P}$ if and only if there exists a certain index $v$ such that $\langle {\bf n}_v, p\rangle<c^{\Gamma}_v$, i.e., 
$\deg(p)\not\geq {\bf c}^\Gamma$ (hence, this $\deg$ function implies a $\mathbb{Z}^r$-grading of the finite dimensional $k$-vector space $\cO_m/{\rm adj}(\mathcal{M})$ --- compare this situation with Example \ref{ex:Z^rgrading}). Then we define the following height function $\mathfrak{h}^\Gamma : R(0, {\bf c}^{\Gamma}) \cap \mathbb{Z}^r \rightarrow \mathbb{Z}_{\geq 0}$ on the rectangle $R(0, \mathbf{c}^{\Gamma})\cap \mathbb{Z}^r$ (with standard basis $\{e_v\}_{v=1}^r$): 
\begin{equation}\label{eq:hgamma}
\mathfrak{h}^\Gamma(\ell)= \big|\{p\in \mathcal{P}\,:\, \deg(p)\not\geq \ell\}\big| \, \hspace{10mm} (= \dim_k \mathcal{O}_m/\{ g \in \mathcal{O}_m\,:\, \mathfrak{v}^{\Gamma}(g) \geq \ell\},
\end{equation}
where $\mathfrak{v}^\Gamma$ corresponds to the tuple of monomial valuations assigned to the primitive facets considered).\\
Clearly, $\mathfrak{h}^\Gamma(0)=0 $  and  $\mathfrak{h}^\Gamma({\bf c}^\Gamma)=|\mathcal{P}| \ \ \ (=\dim_k \mathcal{O}_m/{\rm adj}(\mathcal{M}))$.

Then in the rectangle $R(0,{\bf c}^\Gamma)$ we can consider $\mathfrak{h}^\Gamma$,
its symmetrized function  $(\mathfrak{h}^\Gamma)^{sym}_{{\bf c}^\Gamma}$, and the weight function
$R(0, \mathbf{c}^{\Gamma})\cap \mathbb{Z}^r \ni \ell \mapsto w_0^\Gamma (\ell)= \mathfrak{h}^\Gamma(\ell)+\mathfrak{h}^\Gamma({\bf c}^\Gamma-\ell)-|\mathcal{P}|$ extended as in (\ref{eq:9weight}). For the sake of convenience, we make the following definition.

\begin{define}\label{def:comblathomofNG}
    The bigraded lattice homology $\mathbb{Z}[U]$-module $\mathbb{H}_*(R(0,{\bf c}^\Gamma), w^\Gamma)$ associated with the rectangle $\mathbb{Z}^r \cap R(0,{\bf c}^\Gamma)$ and the weight function $w^\Gamma$ is called the `{\it combinatorial lattice homology assigned to  $N\Gamma$}'.
\end{define}

\begin{lemma}\label{lem:CDPcomb}
The combinatorial height function $\mathfrak{h}^\Gamma$ satisfies the Combinatorial Duality Property in 
the rectangle $R:=R(0,{\bf c}^{\Gamma})$.
\end{lemma}

\begin{proof} 
A proof when all the faces $F_\sigma$ of the Newton boundary are primitive (hence $r=r'$ and for all indices $1 \leq v \leq r$ we have $F_v=F_{\sigma}$ for different $\sigma$'s) can be found in \cite{AgNeCurves} or in \cite[Example 11.8.13]{NBook}. Here we modify that proof 
to cover the general case, too. 

      Suppose indirectly that there exists some lattice point $\ell \in R$ and some index $1 \leq v \leq r$, such that $\ell+e_{v} \in R$ and both $\mathfrak{h}^{\Gamma}(\ell) \neq \mathfrak{h}^{\Gamma}(\ell+e_{v})$ and $ \mathfrak{h}^{\Gamma}(\mathbf{c}^{\Gamma}-\ell) \neq \mathfrak{h}^{\Gamma}(\mathbf{c}^{\Gamma}-\ell-e_{v})$. This means, by (\ref{eq:hgamma}),
      that there exist
$p, q \in \mathcal{P}$ such that 
\begin{itemize}
    \item $\langle \mathbf{n}_w, p \rangle \geq \ell_w, \ \langle \mathbf{n}_w, q\rangle \geq c^{\Gamma}_w-\ell_w$ for all $w \neq v$ and
    \item $\langle \mathbf{n}_{v}, p \rangle =\ell_{v}, \ \langle \mathbf{n}_{v}, q\rangle = c^{\Gamma}_{v}-\ell_{v}-1$
\end{itemize}  
(where $\ell=\sum_{v=1}^r\ell_v e_v$). Then
\begin{equation}\label{eq:hGammaCDP}
    \langle \mathbf{n}_w, p+q+{\bf 1} \rangle \geq c^{\Gamma}_w + \langle \mathbf{n}_w, {\bf 1} \rangle > m_w \text{ for all }w \neq v, \text{ and } \langle \mathbf{n}_{v}, p+q+ {\bf 1} \rangle =m_{v}.
\end{equation}
 This means that $p+q+{\bf 1}$ is a relative interior point of the facet 
 $\{x\,:\, \langle {\bf n}_v, x\rangle=m_v\}\cap N\Gamma$.
If such an interior point does not exist (i.e., when the facet is a primitive) then we get a contradiction. Otherwise, if there exist at least two 
 primitive facets $F_{v_1} \neq F_{v_0}$ on this face of $N\Gamma_+$, then $\mathbf{n}_{v_1} = \mathbf{n}_{v_0}$ and $c^{\Gamma}_{v_1} = c^{\Gamma}_{v_0}$ and thus the two equations of (\ref{eq:hGammaCDP}) corresponding to the indices $v_0$ and $v_1$ contradict each other.
\end{proof}

 The previous proof shows that in order to obtain the CDP it is enough double the nonprimitive facets of the Newton boundary.

\begin{cor}\label{cor:adj=combforNewton}
    Fix a monomial ideal $\mathcal{M}\triangleleft \mathcal{O}_m$ with convenient Newton boundary. The combinatorial lattice homology assigned to $N\Gamma(\mathcal{M})$ (cf. Definition \ref{def:comblathomofNG}) agrees with the symmetric lattice homology module of the adjoint ideal ${\rm adj}(\mathcal{M})$ in $\mathcal{O}_m$: $\mathbb{H}_*(R(0,{\bf c}^\Gamma), w^\Gamma) \cong \mathbb{SH}_*({\rm adj}(\mathcal{M})\triangleleft \cO_m)$. Its Euler characteristic is $\vert \mathcal{P}\vert$.
\end{cor}

\begin{proof}
    We can use the CDP realization $(\{\mathfrak{v}_v^{\Gamma}\}_{v=1}^r, {\bf c}^{\Gamma})$ of ${\rm adj}(\mathcal{M})$ to compute its lattice homology (see Definition \ref{def:LC}). The reader can easily check that the corresponding lattice, rectangle and weight function agrees with the above defined combinatorial one.
\end{proof}

\begin{example}\label{ex:2552ND} 
Consider the monomial ideal $\mathcal{M}=(x_1^8,\, x_1^5 x_2^2, \, x_1^3x_2^2,\, x_2^4) \triangleleft \mathcal{O}_2$. Clear\-ly, its Newton boundary $N \Gamma (\mathcal{M})$ is convenient, while the $1$-dimensional primitive segments are: $F_1=((3,2), (8, 0))$ and $F_2=((0, 4), (3,2))$, thus $\mathbf{n}_1=(2,5)$ and
$\mathbf{n}_2=(2,3)$. 
Moreover, $\vert \mathcal{P}\vert=9$ and ${\bf c}^{\Gamma}=(10,8)$.
The next diagram shows the support of the given generating set of $\mathcal{M}$ and the Newton boundary in blue, its $(-1,-1)$ translate and the lattice points $\mathcal{P}$ in black, and also their 
$\deg$-values and the combinatorial height function 
$\frh^\Gamma$. 

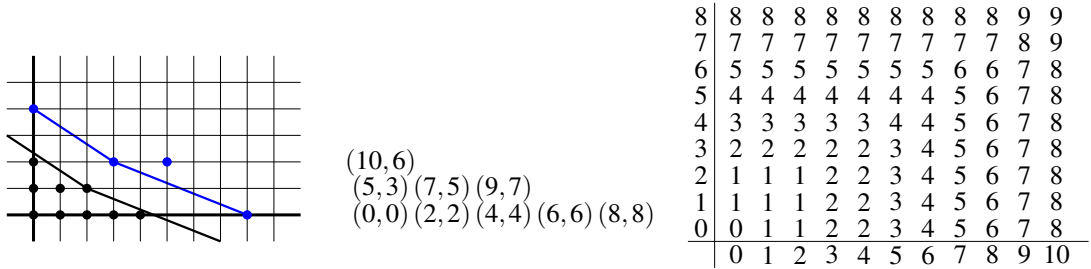
\begin{picture}(420,110)(-20,-25)

\linethickness{0.4mm}
\put(-10,0){\line(1,0){110}}
\put(0,-10){\line(0,1){70}}
\linethickness{0.3mm}
\put(30,20){\color{blue}\line(5,-2){50}}
\put(30,20){\color{blue}\line(-3,2){30}}
\put(20,10){\line(5,-2){50}}
\put(20,10){\line(-3,2){30}}

\linethickness{0.05mm}
  \multiput(-10,10)(0,10){5}{\line(1,0){110}}
  \multiput(10,-10)(10,0){9}{\line(0,1){70}}

\put(10,10){\circle*{3.5}}
\put(20,10){\circle*{3.5}}
\put( 0,10){\circle*{3.5}}
\put( 0,20){\circle*{3.5}}
\put(30,20){\color{blue}\circle*{3.5}}
\put( 0, 0){\circle*{3.5}}
\put(10, 0){\circle*{3.5}}
\put(20, 0){\circle*{3.5}}
\put(30, 0){\circle*{3.5}}
\put(40, 0){\circle*{3.5}}
\put(0, 40){\color{blue}\circle*{3.5}}
\put(80, 0){\color{blue}\circle*{3.5}}
\put(50,20){\color{blue}\circle*{3.5}}

\put(130,00){\makebox(0,0){\small{$(0,0)$}}}
\put(130,10){\makebox(0,0){\small{$(5,3)$}}}
\put(130,20){\makebox(0,0){\small{$(10,6)$}}}
\put(153,00){\makebox(0,0){\small{$(2,2)$}}}
\put(153,10){\makebox(0,0){\small{$(7,5)$}}}
\put(176,00){\makebox(0,0){\small{$(4,4)$}}}
\put(176,10){\makebox(0,0){\small{$(9,7)$}}}
\put(199,00){\makebox(0,0){\small{$(6,6)$}}}
\put(222,00){\makebox(0,0){\small{$(8,8)$}}}

\put(245,-10){\line(1,0){150}}
\put(255,-20){\line(0,1){100}}

\put(263,-15){\makebox(0,0){\small{$0$}}}
\put(275,-15){\makebox(0,0){\small{$1$}}}
\put(287,-15){\makebox(0,0){\small{$2$}}}
\put(299,-15){\makebox(0,0){\small{$3$}}}
\put(311,-15){\makebox(0,0){\small{$4$}}}
\put(323,-15){\makebox(0,0){\small{$5$}}}
\put(335,-15){\makebox(0,0){\small{$6$}}}
\put(347,-15){\makebox(0,0){\small{$7$}}}
\put(359,-15){\makebox(0,0){\small{$8$}}}
\put(371,-15){\makebox(0,0){\small{$9$}}}
\put(383,-15){\makebox(0,0){\small{$10$}}}

\put(250,-5){\makebox(0,0){\small{$0$}}}
\put(250, 5){\makebox(0,0){\small{$1$}}}
\put(250,15){\makebox(0,0){\small{$2$}}}
\put(250,25){\makebox(0,0){\small{$3$}}}
\put(250,35){\makebox(0,0){\small{$4$}}}
\put(250,45){\makebox(0,0){\small{$5$}}}
\put(250,55){\makebox(0,0){\small{$6$}}}
\put(250,65){\makebox(0,0){\small{$7$}}}
\put(250,75){\makebox(0,0){\small{$8$}}}

\put(263,-5){\makebox(0,0){\small{$0$}}}
\put(275,-5){\makebox(0,0){\small{$1$}}}
\put(287,-5){\makebox(0,0){\small{$1$}}}
\put(299,-5){\makebox(0,0){\small{$2$}}}
\put(311,-5){\makebox(0,0){\small{$2$}}}
\put(323,-5){\makebox(0,0){\small{$3$}}}
\put(335,-5){\makebox(0,0){\small{$4$}}}
\put(347,-5){\makebox(0,0){\small{$5$}}}
\put(359,-5){\makebox(0,0){\small{$6$}}}
\put(371,-5){\makebox(0,0){\small{$7$}}}
\put(383,-5){\makebox(0,0){\small{$8$}}}

\put(263, 5){\makebox(0,0){\small{$1$}}}
\put(275, 5){\makebox(0,0){\small{$1$}}}
\put(287, 5){\makebox(0,0){\small{$1$}}}
\put(299, 5){\makebox(0,0){\small{$2$}}}
\put(311, 5){\makebox(0,0){\small{$2$}}}
\put(323, 5){\makebox(0,0){\small{$3$}}}
\put(335, 5){\makebox(0,0){\small{$4$}}}
\put(347, 5){\makebox(0,0){\small{$5$}}}
\put(359, 5){\makebox(0,0){\small{$6$}}}
\put(371, 5){\makebox(0,0){\small{$7$}}}
\put(383, 5){\makebox(0,0){\small{$8$}}}

\put(263,15){\makebox(0,0){\small{$1$}}}
\put(275,15){\makebox(0,0){\small{$1$}}}
\put(287,15){\makebox(0,0){\small{$1$}}}
\put(299,15){\makebox(0,0){\small{$2$}}}
\put(311,15){\makebox(0,0){\small{$2$}}}
\put(323,15){\makebox(0,0){\small{$3$}}}
\put(335,15){\makebox(0,0){\small{$4$}}}
\put(347,15){\makebox(0,0){\small{$5$}}}
\put(359,15){\makebox(0,0){\small{$6$}}}
\put(371,15){\makebox(0,0){\small{$7$}}}
\put(383,15){\makebox(0,0){\small{$8$}}}

\put(263,25){\makebox(0,0){\small{$2$}}}
\put(275,25){\makebox(0,0){\small{$2$}}}
\put(287,25){\makebox(0,0){\small{$2$}}}
\put(299,25){\makebox(0,0){\small{$2$}}}
\put(311,25){\makebox(0,0){\small{$2$}}}
\put(323,25){\makebox(0,0){\small{$3$}}}
\put(335,25){\makebox(0,0){\small{$4$}}}
\put(347,25){\makebox(0,0){\small{$5$}}}
\put(359,25){\makebox(0,0){\small{$6$}}}
\put(371,25){\makebox(0,0){\small{$7$}}}
\put(383,25){\makebox(0,0){\small{$8$}}}

\put(263,35){\makebox(0,0){\small{$3$}}}
\put(275,35){\makebox(0,0){\small{$3$}}}
\put(287,35){\makebox(0,0){\small{$3$}}}
\put(299,35){\makebox(0,0){\small{$3$}}}
\put(311,35){\makebox(0,0){\small{$3$}}}
\put(323,35){\makebox(0,0){\small{$4$}}}
\put(335,35){\makebox(0,0){\small{$4$}}}
\put(347,35){\makebox(0,0){\small{$5$}}}
\put(359,35){\makebox(0,0){\small{$6$}}}
\put(371,35){\makebox(0,0){\small{$7$}}}
\put(383,35){\makebox(0,0){\small{$8$}}}

\put(263,45){\makebox(0,0){\small{$4$}}}
\put(275,45){\makebox(0,0){\small{$4$}}}
\put(287,45){\makebox(0,0){\small{$4$}}}
\put(299,45){\makebox(0,0){\small{$4$}}}
\put(311,45){\makebox(0,0){\small{$4$}}}
\put(323,45){\makebox(0,0){\small{$4$}}}
\put(335,45){\makebox(0,0){\small{$4$}}}
\put(347,45){\makebox(0,0){\small{$5$}}}
\put(359,45){\makebox(0,0){\small{$6$}}}
\put(371,45){\makebox(0,0){\small{$7$}}}
\put(383,45){\makebox(0,0){\small{$8$}}}

\put(263,55){\makebox(0,0){\small{$5$}}}
\put(275,55){\makebox(0,0){\small{$5$}}}
\put(287,55){\makebox(0,0){\small{$5$}}}
\put(299,55){\makebox(0,0){\small{$5$}}}
\put(311,55){\makebox(0,0){\small{$5$}}}
\put(323,55){\makebox(0,0){\small{$5$}}}
\put(335,55){\makebox(0,0){\small{$5$}}}
\put(347,55){\makebox(0,0){\small{$6$}}}
\put(359,55){\makebox(0,0){\small{$6$}}}
\put(371,55){\makebox(0,0){\small{$7$}}}
\put(383,55){\makebox(0,0){\small{$8$}}}

\put(263,65){\makebox(0,0){\small{$7$}}}
\put(275,65){\makebox(0,0){\small{$7$}}}
\put(287,65){\makebox(0,0){\small{$7$}}}
\put(299,65){\makebox(0,0){\small{$7$}}}
\put(311,65){\makebox(0,0){\small{$7$}}}
\put(323,65){\makebox(0,0){\small{$7$}}}
\put(335,65){\makebox(0,0){\small{$7$}}}
\put(347,65){\makebox(0,0){\small{$7$}}}
\put(359,65){\makebox(0,0){\small{$7$}}}
\put(371,65){\makebox(0,0){\small{$8$}}}
\put(383,65){\makebox(0,0){\small{$9$}}}

\put(263,75){\makebox(0,0){\small{$8$}}}
\put(275,75){\makebox(0,0){\small{$8$}}}
\put(287,75){\makebox(0,0){\small{$8$}}}
\put(299,75){\makebox(0,0){\small{$8$}}}
\put(311,75){\makebox(0,0){\small{$8$}}}
\put(323,75){\makebox(0,0){\small{$8$}}}
\put(335,75){\makebox(0,0){\small{$8$}}}
\put(347,75){\makebox(0,0){\small{$8$}}}
\put(359,75){\makebox(0,0){\small{$8$}}}
\put(371,75){\makebox(0,0){\small{$9$}}}
\put(383,75){\makebox(0,0){\small{$9$}}}
\end{picture}

The next diagram shows the $w^\Gamma$-values in the rectangle $R\cap \mathbb{Z}^r$
and the graded root. 

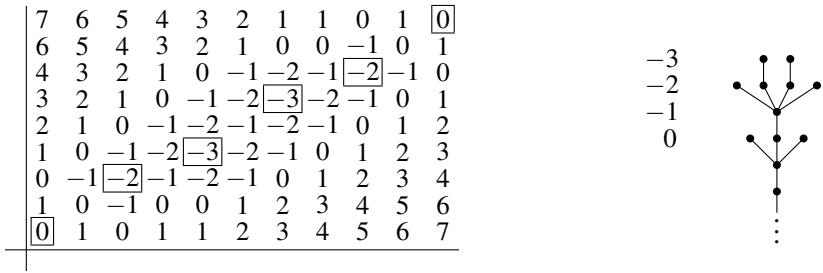
\begin{picture}(420,110)( 5,-20)

\put(51,-12){\line(1,0){170}}
\put(59,-20){\line(0,1){100}}

\put(65,-5){\makebox(0,0){\small{$0$}}}
\put(80,-5){\makebox(0,0){\small{$1$}}}
\put(95,-5){\makebox(0,0){\small{$0$}}}
\put(110,-5){\makebox(0,0){\small{$1$}}}
\put(125,-5){\makebox(0,0){\small{$1$}}}
\put(140,-5){\makebox(0,0){\small{$2$}}}
\put(155,-5){\makebox(0,0){\small{$3$}}}
\put(170,-5){\makebox(0,0){\small{$4$}}}
\put(185,-5){\makebox(0,0){\small{$5$}}}
\put(200,-5){\makebox(0,0){\small{$6$}}}
\put(215,-5){\makebox(0,0){\small{$7$}}}
\put(61,-10){\framebox(8,10)}

\put(65, 5){\makebox(0,0){\small{$1$}}}
\put(80, 5){\makebox(0,0){\small{$0$}}}
\put(95, 5){\makebox(0,0){\small{$-1$}}}
\put(110, 5){\makebox(0,0){\small{$0$}}}
\put(125, 5){\makebox(0,0){\small{$0$}}}
\put(140, 5){\makebox(0,0){\small{$1$}}}
\put(155, 5){\makebox(0,0){\small{$2$}}}
\put(170, 5){\makebox(0,0){\small{$3$}}}
\put(185, 5){\makebox(0,0){\small{$4$}}}
\put(200, 5){\makebox(0,0){\small{$5$}}}
\put(215, 5){\makebox(0,0){\small{$6$}}}

\put(65,15){\makebox(0,0){\small{$0$}}}
\put(80,15){\makebox(0,0){\small{$-1$}}}
\put(95,15){\makebox(0,0){\small{$-2$}}}
\put(110,15){\makebox(0,0){\small{$-1$}}}
\put(125,15){\makebox(0,0){\small{$-2$}}}
\put(140,15){\makebox(0,0){\small{$-1$}}}
\put(155,15){\makebox(0,0){\small{$0$}}}
\put(170,15){\makebox(0,0){\small{$1$}}}
\put(185,15){\makebox(0,0){\small{$2$}}}
\put(200,15){\makebox(0,0){\small{$3$}}}
\put(215,15){\makebox(0,0){\small{$4$}}}
\put(88,10){\framebox(14,10)}

\put(65,25){\makebox(0,0){\small{$1$}}}
\put(80,25){\makebox(0,0){\small{$0$}}}
\put(95,25){\makebox(0,0){\small{$-1$}}}
\put(110,25){\makebox(0,0){\small{$-2$}}}
\put(125,25){\makebox(0,0){\small{$-3$}}}
\put(140,25){\makebox(0,0){\small{$-2$}}}
\put(155,25){\makebox(0,0){\small{$-1$}}}
\put(170,25){\makebox(0,0){\small{$0$}}}
\put(185,25){\makebox(0,0){\small{$1$}}}
\put(200,25){\makebox(0,0){\small{$2$}}}
\put(215,25){\makebox(0,0){\small{$3$}}}
\put(118,20){\framebox(14,10)}

\put(65,35){\makebox(0,0){\small{$2$}}}
\put(80,35){\makebox(0,0){\small{$1$}}}
\put(95,35){\makebox(0,0){\small{$0$}}}
\put(110,35){\makebox(0,0){\small{$-1$}}}
\put(125,35){\makebox(0,0){\small{$-2$}}}
\put(140,35){\makebox(0,0){\small{$-1$}}}
\put(155,35){\makebox(0,0){\small{$-2$}}}
\put(170,35){\makebox(0,0){\small{$-1$}}}
\put(185,35){\makebox(0,0){\small{$0$}}}
\put(200,35){\makebox(0,0){\small{$1$}}}
\put(215,35){\makebox(0,0){\small{$2$}}}

\put(65,45){\makebox(0,0){\small{$3$}}}
\put(80,45){\makebox(0,0){\small{$2$}}}
\put(95,45){\makebox(0,0){\small{$1$}}}
\put(110,45){\makebox(0,0){\small{$0$}}}
\put(125,45){\makebox(0,0){\small{$-1$}}}
\put(140,45){\makebox(0,0){\small{$-2$}}}
\put(155,45){\makebox(0,0){\small{$-3$}}}
\put(170,45){\makebox(0,0){\small{$-2$}}}
\put(185,45){\makebox(0,0){\small{$-1$}}}
\put(200,45){\makebox(0,0){\small{$0$}}}
\put(215,45){\makebox(0,0){\small{$1$}}}
\put(148,40){\framebox(14,10)}

\put(65,55){\makebox(0,0){\small{$4$}}}
\put(80,55){\makebox(0,0){\small{$3$}}}
\put(95,55){\makebox(0,0){\small{$2$}}}
\put(110,55){\makebox(0,0){\small{$1$}}}
\put(125,55){\makebox(0,0){\small{$0$}}}
\put(140,55){\makebox(0,0){\small{$-1$}}}
\put(155,55){\makebox(0,0){\small{$-2$}}}
\put(170,55){\makebox(0,0){\small{$-1$}}}
\put(185,55){\makebox(0,0){\small{$-2$}}}
\put(200,55){\makebox(0,0){\small{$-1$}}}
\put(215,55){\makebox(0,0){\small{$0$}}}
\put(178,50){\framebox(14,10)}

\put(65,65){\makebox(0,0){\small{$6$}}}
\put(80,65){\makebox(0,0){\small{$5$}}}
\put(95,65){\makebox(0,0){\small{$4$}}}
\put(110,65){\makebox(0,0){\small{$3$}}}
\put(125,65){\makebox(0,0){\small{$2$}}}
\put(140,65){\makebox(0,0){\small{$1$}}}
\put(155,65){\makebox(0,0){\small{$0$}}}
\put(170,65){\makebox(0,0){\small{$0$}}}
\put(185,65){\makebox(0,0){\small{$-1$}}}
\put(200,65){\makebox(0,0){\small{$0$}}}
\put(215,65){\makebox(0,0){\small{$1$}}}

\put(65,75){\makebox(0,0){\small{$7$}}}
\put(80,75){\makebox(0,0){\small{$6$}}}
\put(95,75){\makebox(0,0){\small{$5$}}}
\put(110,75){\makebox(0,0){\small{$4$}}}
\put(125,75){\makebox(0,0){\small{$3$}}}
\put(140,75){\makebox(0,0){\small{$2$}}}
\put(155,75){\makebox(0,0){\small{$1$}}}
\put(170,75){\makebox(0,0){\small{$1$}}}
\put(185,75){\makebox(0,0){\small{$0$}}}
\put(200,75){\makebox(0,0){\small{$1$}}}
\put(215,75){\makebox(0,0){\small{$0$}}}
\put(211,70){\framebox(8,10)}

\put(300,30){\makebox(0,0){\small{$0$}}} 
\put(297,40){\makebox(0,0){\small{$-1$}}}
\put(297,50){\makebox(0,0){\small{$-2$}}}
\put(297,60){\makebox(0,0){\small{$-3$}}}
\put(340,-2){\makebox(0,0){$\vdots$}} 
\put(340,2){\line(0,1){38}}

\put(340,10){\circle*{3}}

\put(340,20){\circle*{3}} \put(340,20){\line(1,1){10}} \put(340,20){\line(-1,1){10}}

\put(330,30){\circle*{3}}
\put(350,30){\circle*{3}} 
\put(340,30){\circle*{3}}

\put(340,40){\circle*{3}} \put(340,40){\line(1,2){5}} \put(340,40){\line(3,2){15}} \put(340,40){\line(-1,2){5}} \put(340,40){\line(-3,2){15}}

\put(325,50){\circle*{3}}
\put(335,50){\circle*{3}} \put(335,50){\line(0,1){10}}
\put(345,50){\circle*{3}} \put(345,50){\line(0,1){10}}
\put(355,50){\circle*{3}} 

\put(335,60){\circle*{3}}
\put(345,60){\circle*{3}}
\end{picture}

On the diagram above 
the rectangles denote the local minimum points of the weight function $w^\Gamma$ (cf. subsection \ref{ss:glm}). Notice that 
$\min w_0^{\Gamma}=-3, \ {\rm rank}_{\mathbb{Z}}(\mathbb{H}_{{\rm red},0})=6$, while ${\rm rank}_{\mathbb{Z}}(\mathbb{H}_{\geq 1})=0$, and thus $eu(\mathbb{H}_*)=9=|\mathcal{P}|.$
\end{example}

\subsection{Hypersurface singularities with Newton nondegenerate principal part}\, \label{ss:NNPP}

Given a convenient Newton diagram, we can associate to it a hypersurface singularity, \emph{nondegenerate} in the sense of Kouchnirenko and Varchenko \cite{Kou,Var}, with well-defined embedded topological type. We seek to compare its analytic lattice homology (of $m-1$ forms) with the combinatorial one associated to the corresponding adjoint ideal introduced in the previous subsection. However, we begin the presentation from the viewpoint of a given complex analytic hypersurface singularity.

Let $f:(\C^{m},0)\to (\C,0)$, $f(x)=\sum _{p}a_px^p$, $a_p\in\C \ (\forall p \in (\mathbb{Z}_{>0})^m)$, 
be a complex analytic hypersurface germ. We define its support as $\mathfrak{S}(f):=\{p\in (\mathbb{Z}_{\geq 0})^m\,:\, a_p\not=0\}$. Then, similarly to subsection \ref{ss:ND},
we can construct the Newton boundary $N\Gamma$ and $N\Gamma_{\pm}$ associated to $\mathfrak{S}(f)$; we denote them by 
 $N\Gamma(f)$ and $N\Gamma_{\pm}(f)$.
 In the sequel, we will assume that $N\Gamma(f)$ is convenient. 
 We also write $f_0:=\sum_{p\in N\Gamma(f)}a_px^p$, and we call it the \textit{Newton principal part of $f$}. 

 In order to succesfully connect the combinatorics of 
 $N\Gamma(f)$ and the discrete invariants of the hypersurface singularity $f$,  we need to impose a genericity 
 condition on $f$: for any face $F$ of $N\Gamma(f)$ (of any dimension 
 $0\leq q\leq m-1$) we write $f_F(x)=\sum_{p\in F}a_px^p$ and we  assume that the equation
 $\partial f_F(x)=0$ has no solution in $(\C^*)^m$.  If $f$ satisfies this condition  we say, after Kouchnirenko and Varchenko,  that it \textit{has a Newton nondegenerate principal part}. Such a germ always defines an isolated singularity. We also remark that this nondegeneracy condition is satisfied generically (i.e., for a generic choice of the coefficients of $f_0$).
 For details see  \cite{Kou,Var}. Moreover, Oka proved the following:

\begin{theorem}\cite[Theorem 2.1]{Oka}
    Suppose that $f:(\mathbb{C}^m, 0) \rightarrow (\mathbb{C}, 0)$ is an isolated singularity and it has Newton nondegenerate principal part. Then the Milnor fibration, and hence the embedded topological type  of $(X,0):=(\{f=0\},0) \subset (\C^{m},0)$ is determined by the Newton boundary $N\Gamma(f)$.
\end{theorem}

Similarly, the same is true for the delta invariant in the case $m=2$ and for the geometric genus in the case $m>2$.  The precise statement is provided by the following 
result of Merle and Teissier. In order to formulate it, we recall that a hypersurface singularity is automatically Gorenstein (that is, it admits a nowhere vanishing form 
$\omega_G\in H^0(X\setminus \{0\}, \Omega^{m-1} _{X\setminus \{0\}})$), in fact, 
its Gorenstein form (up to a nonzero constant) is 
$$\omega_G=\dfrac{dx_2\wedge dx_3 \wedge \ldots \wedge dx_m}{\frac{\partial f}{\partial x_1}}\Big|_{X\setminus \{0\}}=\
\ldots\  = (-1)^{m-1}\ 
\dfrac{dx_1\wedge dx_2 \wedge \ldots \wedge dx_{m-1}}{\frac{\partial f}{\partial x_m}}\Big|_{X\setminus \{0\}}.
$$
We denote by $q$ the  natural quotient map
$\mathcal{O}_{\mathbb{C}^m, 0} \rightarrow\mathcal{O}_{\mathbb{C}^m, 0}/(f)=\cO_{X,0}$
and by $\phi: \widetilde{X} \rightarrow X$ any good resolution of $(X,0)$. 
If $m>2$ then the divisor of the  pullback of $\omega_G$  by the resolution $\phi$ is $-Z_K$ (compare with subsection \ref{ss:zcoh}). 

\begin{theorem} \cite[2.1.1 Théorème]{MT} \label{th:MerleTeissier}
    Let $f:(\mathbb{C}^m, 0) \rightarrow (\mathbb{C}, 0)$ be an isolated hypersurface singularity with Newton nondegenerate principal part.
  Let    $\mathcal{C}_f=q^{-1}(\mathcal{C})$ be the 
`adjunction conductor' (conducteur d'adjonction) --- the pullback of the conductor ideal --- and let us
denote by $\mathcal{M}\triangleleft \mathcal{O}_{\mathbb{C}^m, 0}$ the monomial ideal generated by monomials $x^p$ with exponent $p\in\mathfrak{S}(f)$.
    Then we have the following:
    \begin{equation}
        \mathcal{C}_f = {\rm adj}(\mathcal{M}) = \left(\{ x^p\,:\, p+{\bf 1}\in N\Gamma _+(f)\setminus 
        N\Gamma(f)\}\right) \triangleleft \mathcal{O}_{\mathbb{C}^m, 0}.
    \end{equation}
    As a corollary,  $\dim_{\mathbb{C}} \mathcal{O}_{\mathbb{C}^m, 0}/\mathcal{C}_f=\begin{cases}
         \delta(X, 0) & \text{ if } m=2; \\
         p_g(X, 0) & \text{ if } m\geq3;
    \end{cases}$ \ which can be read from 
     the Newton boundary $N \Gamma(f)$ as the cardinality $\vert \mathcal{P}\vert$ of the lattice points $\mathcal{P}= 
\big( (\Z_{>0})^{m}\cap  
N\Gamma _-(\mathcal{M})\big)-{\bf 1}$.
\end{theorem}

\begin{remark} Recall that the conductor ideal was defined for curve singularities in paragraph \ref{par:conductor}, whereas for higher dimensional singularities in Definition \ref{def:CONDINDEF} (see also Proposition \ref{prop:CONDINDEP}). Accordingly, the adjunction conductor $\mathcal{C}_f$ of a hypersurface singularity $(\{f=0\}, 0)\subset (\mathbb{C}^m, 0)$ is
    \begin{itemize}
        \item $m=2$: the pullback  $q^{-1}(\mathcal{C})$ of the conductor ideal $\mathcal{C}$,
       
        \item $m\geq 3$: the ideal $q^{-1}\left((\phi_{\ast}(\mathcal{O}_{\widetilde{X}}(-Z_K)))_0\right)\triangleleft\mathcal{O}_{\mathbb{C}^m, 0}$ of holomorphic functions $\varphi$ satisfying that the pullback along $\phi$ of the meromorphic $(m-1)$-form $\varphi \omega_G$ has a holomorphic extension to the whole of  $\widetilde{X}$.
    \end{itemize} 
\end{remark}

 \subsection{Equivalence of the analytic and the combinatorial lattice homology}\label{ss:NNan&comb}\,

Having  a convenient Newton boundary $N\Gamma$,  we face the following interesting situation: associated with 
$N\Gamma$ we have {\it two naturally defined symmetric lattice homologies}.

The first one, the `analytic one', is constructed as follows. 
Take an isolated complex analytic hypersurface singularity $f$ with nondegenerate principal part supported on $N\Gamma$ and then consider its analytic lattice homology $\bH_{an, *}(\{f=0\},0)$ (in the $m-1\geq 2$ dimensional case we denoted this by $\mathbb{H}_*((\{f=0\}, 0), \Omega^{m-1})$ in section \ref{s:dnagy}, nevertheless, here we will use the $\bH_{an, *}$ notation to stay compatible with the $m=2$ case from section \ref{s:deccurves}).
In the curve case it depends only on the embedded topological type, hence only on the Newton boundary $N\Gamma(f)$. In higher dimensions $m>2$ a priori it is not clear that 
 a similar statement holds, in principle it might be possible that $\bH_{an, *}(\{f=0\},0)$
 depends on the choice of the coefficients  of the monomials of 
 $f$ (i.e., of the `analytic type' of $(\{f=0\},0)$).  

 On the other hand, we have a second module, constructed as in subsection \ref{ss:NewtonLatCoh}.
 Namely, we consider the monomial ideal $\mathcal{M}$ generated by 
 monomials of type $\{x^p\,:\, p\in N\Gamma_+\}$, and then associate to it the symmetric lattice homology
  $\mathbb{S}\bH_*({\rm adj}({\mathcal{M}})\triangleleft \cO_m)$. By Corollary \ref{cor:adj=combforNewton} it agrees with the combinatorial lattice homology assigned to $N\Gamma$, and hence it can be computed combinatorially by counting specific lattice points under the Newton boundary. 

It turns out,  however, that the two lattice homology modules are isomorphic, in particular showing also that the output of the analytic construction depends only on Newton boundary $N\Gamma$.
  In fact this isomorphism question was our original starting point and main motivation for one of the main results, the Independence 
  Theorem \ref{th:Indep},  of the article.

\begin{theorem}\label{th:comparison}
    Let $f:(\mathbb{C}^m, 0) \rightarrow (\mathbb{C}, 0)$ be an isolated hypersurface singularity with Newton nondegenerate principal part. Then, its analytic lattice homology  module (in the $m>2$ case we can think of the analytic lattice homology of $(m-1)$-forms as well) agrees with the combinatorial lattice homology module assigned to its Newton boundary $N\Gamma (f)$:
\begin{equation*}
    \mathbb{H}_{an, \ast}(\{f=0\}, 0) \cong \mathbb{H}_{\ast}(R(0, {\bf c}^{\Gamma}), w^{\Gamma}).
\end{equation*}
\end{theorem}

\begin{proof}
    Let us set the notation $X:=\{ f=0\} \subset \mathbb{C}^m$. By Corollaries \ref{cor:curveArtin} ($m=2$) and \ref{cor:GORCOND} ($m > 2$), we have that $\mathbb{H}_{an, \ast}(X, 0) \cong \mathbb{SH}_{\ast}(\mathcal{C} \triangleleft \mathcal{O}_{X, 0})$, which, by Theorem \ref{th:Artin} {\it (1)}, agrees with $\mathbb{SH}_{\ast}(\mathcal{C}_f \triangleleft \mathcal{O}_{\C^m, 0})$, where $\mathcal{C}_f=q^{-1}(\mathcal{C})$ is the pullback of the conductor ideal via
    the quotient map $q:\mathcal{O}_{\mathbb{C}^m, 0} \rightarrow \mathcal{O}_{X,0}$. On the other hand, by Corollary \ref{cor:adj=combforNewton} we have $\mathbb{H}_*(R(0,{\bf c}^\Gamma), w^\Gamma) \cong \mathbb{SH}_*({\rm adj}(\mathcal{M})\triangleleft \cO_m)$, where $\mathcal{M}$ is the monomial ideal in $\mathcal{O}_{\mathbb{C}^m, 0}$ gene\-rated by monomials of type $\{x^p\,:\, p\in N\Gamma_+(f)\}$. Now,  \cite[2.1.1 Théorème]{MT} of Merle and Teissier (presented here as Theorem \ref{th:MerleTeissier}) identifies $\mathcal{C}_f$ and ${\rm adj}(\mathcal{M})$, and, hence, provides the proof.
\end{proof}

As a corollary, we can answer the  question \cite[Problem 11.9.52]{NBook}, proposed by the first author:

\begin{cor}\label{cor:HanstabforNN}
    The analytic lattice homology (of $(m-1)$-forms) is constant along an equisingular family of hypersurface singularities (of dimension $m-1$) with Newton nondegenerate principal part associated with a fixed Newton boundary $N\Gamma$. The above Theorem \ref{th:comparison}  also provides a combinatorial construction of this module.  
\end{cor}

\subsection{Another approach for the plane curve singularity case}\label{ss:NNplanecurves}\,

For the identification of the ideals part of the proof of Theorem \ref{th:comparison} we also present another, more computation based, but elementary argument.
The motivation to provide this chain of identifications at the level of more 
classical invariants is that this alternative proof can be generalized to some Newton degenerate cases as well, as we will see in Examples \ref{bek:10.6.3} and \ref{bek:10.6.5}.
Along these steps we also see that the height functions $\hh_{an}$ and $\hh^{\Gamma}$ associated with the analytic 
and combinatorial situations  do not agree (see Example \ref{ex:2523}). In particular, 
we face the situation that two weight functions
(both naturally assigned to a geometric picture)  differ but their 
lattice homologies agree. 

\bekezdes \textbf{Newton nondegenerate plane curve singularities.}
The following facts are well known for plane curve singularities with Newton nondegenerate principal part:

\begin{proposition}(See, e.g., in \cite{CNP, wall}.) \label{prop:NNcurves}
    Let $(C, 0) \subset (\mathbb{C}^2, 0)$ be a reduced plane curve singularity, the zero locus of a holomorphic function $f \in \mathbb{C}\{x_1, x_2\}$. Suppose that $f$ has nondegenerate Newton principal part and convenient Newton boundary $N \Gamma(f)$. Then
    
    {\it (i)} given a $1$-dimensional face $F$ of the Newton boundary with positive primitive normal vector ${\bf n}=(n_1, n_2)$, the weighted homogeneous polynomial $f_F(x)=\sum_{p\in F} a_p x^p$ splits as a product
    \begin{center}
        $f_F(x_1, x_2)=c x_1^{a_1}x_2^{a_2}\prod_{i=1}^k(x_1^{n_2}-b_ix_2^{n_1}), \text{ with } b_i \neq b_j$ if $i\neq j$,
    \end{center} where $F$ has endpoints $(a_1+kn_2, a_2)$ and $(a_1, a_2+kn_1)$ and contains $k-1$ interior lattice points (i.e., it consists of $k$ primitive segments);
    
    {\it (ii)} via Newton's algorithm, the quasi-tangents $\{x_1^{n_2}-b_ix_2^{n_1}=0\}_{i=1}^k$ correspond each to different irreducible components (having the corresponding Newton principal part) of the curve singularity $(C, 0)$, hence the total number of primitive segments of $N\Gamma(f)$ equals the number $r$ of components in the irreducible decomposition $(C,0)=\bigcup_{i=1}^r(C_i,0)$;
    
    {\it (iii)} the component corresponding to the quasi-tangent $\{x_1^{n_2}-b_ix_2^{n_1}=0\}$ (we can suppose without loss of generality that $n_1< n_2$) has a single Puiseux pair, specifically $(n_1, n_2)$,  and Puiseux parametri\-zation of form $t \mapsto (t^{n_1}, b_i^{-1/n_1}t^{n_2}+$ higher order terms$)$ (we remark here, that the corres\-ponding map on the local algebras $\mathcal{O}_{C, 0} \rightarrow \mathbb{C}\{t\}$ is the normalization);

    {\it (iv)} the intersection multiplicity at $0$ of the components $C_i$ and $C_j$, $i \neq j$, corresponding to the quasi-tangents  $\{x_1^{n_2}-b_ix_2^{n_1}=0\}$ and $\{x_1^{n'_2}-b_jx_2^{n'_1}=0\}$, is $i_0(C_i, C_j)=\min\{n_1n'_2, n_2n'_1\}$.
\end{proposition}

From part {\it (ii)} we already see, that the lattices underlying the analytic lattice homology of the Newton nondegenerate curve $(C, 0)=(\{f=0\}, 0)$ and the combinatorial lattice homology assigned with the Newton boundary $N\Gamma(f)$ are both identified with $\mathbb{Z}^r$. Let us fix a common indexing of the irreducible components $(C,0)=\bigcup_{v=1}^r(C_v,0)$ and the correspon\-ding primitive edges of the Newton boundary $N\Gamma(f)=\bigcup_{v=1}^rF_v$ 

 Let us recall that  ${\bf c}_{an}:={\bf c}\in (\Z_{\geq 0})^{r}$ denotes the (analytic) conductor element of $(C,0)$ (cf. subsection \ref{ss:ancurves}), while ${\bf c}^\Gamma\in (\Z_{\geq 0})^{r}$ is the combinatorial conductor, whose coefficients
satisfy 
$$c^\Gamma_v=m_v+1-\langle {\bf n}_v, {\bf 1}\rangle = \min\{\frv^\Gamma_v(p)\,:\, x^p\in {\rm adj}(\mathcal{M})\} \text{ for all } 1 \leq v \leq r$$
(cf. (\ref{eq:c2}) and the notations of paragraph \ref{par:NewtonComb}). 

\begin{lemma}\label{lem:COND}
Under the above identification of lattices ${\bf c}_{an}={\bf c}^\Gamma$.
That is, the analytic and combinatorial conductor elements agree. 
(We will denote this common conductor by ${\bf c}$.) Also, for the set of lattice points $\mathcal{P}= 
\big( (\Z_{>0})^{2}\cap  
N\Gamma _-(f)\big)-{\bf 1}$ we have $|\mathcal{P}|=|{\bf c}|/2=\delta(C, 0)$.
\end{lemma}
\begin{proof}
For the irreducible germ $x_1^{n_2}-bx_2^{n_1}$ we have 
\begin{center}
    $\mathbf{c}_{an}=\mathbf{c}^\Gamma=2\delta=(n_1-1)(n_2-1)$ and $|\mathcal{P}|=(n_1-1)(n_2-1)/2$.
\end{center} 
Then one can proceed by induction on the number of components $r$ using the identity 
\begin{center}
$c_{an,v}=\mathbf{c}_{an}(C_v)+i_0(C_v,\cup_{w\not=v}C_w)$ (cf. Remark \ref{rem:condofapc})
\end{center}
and the corresponding part {\it (iv)} of Proposition \ref{prop:NNcurves}. 
\end{proof}

\bekezdes Let $\mathfrak{h}_{an}:=\mathfrak{h}$, respectively $\mathfrak{h}^\Gamma$, be the two Hilbert\,/\,height functions associated with the analytic, respectively combinatorial setup as defined in subsection \ref{ss:ancurves} via the valuations $\mathfrak{v}^{an}_v$ corresponding to the normalization (we add here an extra $_{an}$ or $^{an}$ sub\,/\,superscript to the notation, to differentiate from the combinatorial counterparts), respectively in (\ref{eq:hgamma}) via the monomial valuations $\mathfrak{v}^{\Gamma}_v$ corresponding to the primitive segments of the Newton diagram (defined in the first paragraph of subsection \ref{ss:NewtonLatCoh}). 

\begin{lemma}\label{lem:CDPcurves}
Both height functions $\mathfrak{h}_{an}$ and $\mathfrak{h}^\Gamma$ satisfy the Combinatorial Duality Property in 
the rectangle $R:=R(0,{\bf c})$.
\end{lemma}

\begin{proof}
The analytic case was verified in Remark \ref{rem:hh^bulletCDP}, while the combinatorial case in Lemma \ref{lem:CDPcomb}.
\end{proof}

\begin{lemma}\label{lem:van>vcomb}
    For any lattice point $\ell \in \mathbb{Z}^r \cap R(0, {\bf c})$ we have $\mathfrak{h}_{an}(\ell)\leq \mathfrak{h}^\Gamma(\ell)$, 
hence 
\begin{equation}\label{eq:wanwGamma}
    w_{an}(\ell)= \mathfrak{h}_{an}(\ell)+\mathfrak{h}_{an}({\bf c}-\ell)-\delta \leq \mathfrak{h}^\Gamma(\ell) + \mathfrak{h}^\Gamma({\bf c}-\ell)- |\mathcal{P}|= w^\Gamma(\ell).
\end{equation}
\end{lemma}
\begin{proof}
    Consider again the quotient map $q:\cO_{\mathbb{C}^2,0}\to \cO_{\mathbb{C}^2, 0}/(f)=\cO_{C,0}$, 
    and let us define 
the composition valuations $\tilde {\frv}^{an}_v := \frv^{an}_v\circ q: \cO_{\mathbb{C}^2, 0}\to
\Z_{\geq 0}\cup \{\infty\}$.
Since $q$ is a quotient map, the Hilbert functions associated with 
$\tilde{\frv}^{an}_v$ and $\frv^{an}_v$ agree, we will use the notation 
$\mathfrak{h}_{an}$ for both of them. Also, from part {\it (iii)} of Proposition \ref{prop:NNcurves} we see that $\tilde{\mathfrak{v}}^{an}_v(x^p)=\langle {\bf n}_v, p\rangle=\frv_v^{\Gamma}(x^{p})$ for every monomial. Therefore, by the minimality property of monomial valuations (cf.  \cite[Definition 6.1.4]{HS}), 
 we obtain 
that 
 $   \tilde{\frv}_v^{an}(g) \geq \frv_v^{\Gamma}(g)$ for all $  g \in \cO_{\mathbb{C}^2, 0}$.
\end{proof}

We want to stress, though, that  \emph{in general we do not have equality in these inequalities}, e.g., in equation (\ref{eq:wanwGamma}).

\begin{example}\label{ex:2523} \
(i) Consider 
 the plane curve singularity $(C, 0) \subset (\mathbb{C}^2, 0)$, defined as the zero set of the holomorphic function $f(x_1, x_2)=(x_1^5-x_2^2)(x_1^3-x_2^2)$. The following diagram
 shows the combinatorial and analytic weight functions in $R(0,{\bf c})$. 
 One can see, that they differ in the boxed regions. However, 
 (as a computation based on the diagram shows, too) 
 the lattice homology modules are still the same.

 \begin{picture}(420,145)(35,-40)

\put(51,-12){\line(1,0){170}}
\put(59,-20){\line(0,1){118}}

\put(65,-5){\makebox(0,0){\small{$0$}}}
\put(80,-5){\makebox(0,0){\small{$1$}}}
\put(95,-5){\makebox(0,0){\small{$0$}}}
\put(110,-5){\makebox(0,0){\small{$1$}}}
\put(125,-5){\makebox(0,0){\small{$1$}}}
\put(140,-5){\makebox(0,0){\small{$2$}}}
\put(155,-5){\makebox(0,0){\small{$3$}}}
\put(170,-5){\makebox(0,0){\small{$4$}}}
\put(185,-5){\makebox(0,0){\small{$5$}}}
\put(200,-5){\makebox(0,0){\small{$6$}}}
\put(215,-5){\makebox(0,0){\small{$7$}}}
\put(120,-10){\framebox(100,23)}

\put(65, 7){\makebox(0,0){\small{$1$}}}
\put(80, 7){\makebox(0,0){\small{$0$}}}
\put(95, 7){\makebox(0,0){\small{$-1$}}}
\put(110, 7){\makebox(0,0){\small{$0$}}}
\put(125, 7){\makebox(0,0){\small{$0$}}}
\put(140, 7){\makebox(0,0){\small{$1$}}}
\put(155, 7){\makebox(0,0){\small{$2$}}}
\put(170, 7){\makebox(0,0){\small{$3$}}}
\put(185, 7){\makebox(0,0){\small{$4$}}}
\put(200, 7){\makebox(0,0){\small{$5$}}}
\put(215, 7){\makebox(0,0){\small{$6$}}}

\put(65,19){\makebox(0,0){\small{$0$}}}
\put(80,19){\makebox(0,0){\small{$-1$}}}
\put(95,19){\makebox(0,0){\small{$-2$}}}
\put(110,19){\makebox(0,0){\small{$-1$}}}
\put(125,19){\makebox(0,0){\small{$-2$}}}
\put(140,19){\makebox(0,0){\small{$-1$}}}
\put(155,19){\makebox(0,0){\small{$0$}}}
\put(170,19){\makebox(0,0){\small{$1$}}}
\put(185,19){\makebox(0,0){\small{$2$}}}
\put(200,19){\makebox(0,0){\small{$3$}}}
\put(215,19){\makebox(0,0){\small{$4$}}}

\put(65,31){\makebox(0,0){\small{$1$}}}
\put(80,31){\makebox(0,0){\small{$0$}}}
\put(95,31){\makebox(0,0){\small{$-1$}}}
\put(110,31){\makebox(0,0){\small{$-2$}}}
\put(125,31){\makebox(0,0){\small{$-3$}}}
\put(140,31){\makebox(0,0){\small{$-2$}}}
\put(155,31){\makebox(0,0){\small{$-1$}}}
\put(170,31){\makebox(0,0){\small{$0$}}}
\put(185,31){\makebox(0,0){\small{$1$}}}
\put(200,31){\makebox(0,0){\small{$2$}}}
\put(215,31){\makebox(0,0){\small{$3$}}}

\put(65,43){\makebox(0,0){\small{$2$}}}
\put(80,43){\makebox(0,0){\small{$1$}}}
\put(95,43){\makebox(0,0){\small{$0$}}}
\put(110,43){\makebox(0,0){\small{$-1$}}}
\put(125,43){\makebox(0,0){\small{$-2$}}}
\put(140,43){\makebox(0,0){\small{$-1$}}}
\put(155,43){\makebox(0,0){\small{$-2$}}}
\put(170,43){\makebox(0,0){\small{$-1$}}}
\put(185,43){\makebox(0,0){\small{$0$}}}
\put(200,43){\makebox(0,0){\small{$1$}}}
\put(215,43){\makebox(0,0){\small{$2$}}}

\put(65,55){\makebox(0,0){\small{$3$}}}
\put(80,55){\makebox(0,0){\small{$2$}}}
\put(95,55){\makebox(0,0){\small{$1$}}}
\put(110,55){\makebox(0,0){\small{$0$}}}
\put(125,55){\makebox(0,0){\small{$-1$}}}
\put(140,55){\makebox(0,0){\small{$-2$}}}
\put(155,55){\makebox(0,0){\small{$-3$}}}
\put(170,55){\makebox(0,0){\small{$-2$}}}
\put(185,55){\makebox(0,0){\small{$-1$}}}
\put(200,55){\makebox(0,0){\small{$0$}}}
\put(215,55){\makebox(0,0){\small{$1$}}}

\put(65,67){\makebox(0,0){\small{$4$}}}
\put(80,67){\makebox(0,0){\small{$3$}}}
\put(95,67){\makebox(0,0){\small{$2$}}}
\put(110,67){\makebox(0,0){\small{$1$}}}
\put(125,67){\makebox(0,0){\small{$0$}}}
\put(140,67){\makebox(0,0){\small{$-1$}}}
\put(155,67){\makebox(0,0){\small{$-2$}}}
\put(170,67){\makebox(0,0){\small{$-1$}}}
\put(185,67){\makebox(0,0){\small{$-2$}}}
\put(200,67){\makebox(0,0){\small{$-1$}}}
\put(215,67){\makebox(0,0){\small{$0$}}}

\put(65,79){\makebox(0,0){\small{$6$}}}
\put(80,79){\makebox(0,0){\small{$5$}}}
\put(95,79){\makebox(0,0){\small{$4$}}}
\put(110,79){\makebox(0,0){\small{$3$}}}
\put(125,79){\makebox(0,0){\small{$2$}}}
\put(140,79){\makebox(0,0){\small{$1$}}}
\put(155,79){\makebox(0,0){\small{$0$}}}
\put(170,79){\makebox(0,0){\small{$0$}}}
\put(185,79){\makebox(0,0){\small{$-1$}}}
\put(200,79){\makebox(0,0){\small{$0$}}}
\put(215,79){\makebox(0,0){\small{$1$}}}
\put(61,73){\framebox(100,24)}

\put(65,91){\makebox(0,0){\small{$7$}}}
\put(80,91){\makebox(0,0){\small{$6$}}}
\put(95,91){\makebox(0,0){\small{$5$}}}
\put(110,91){\makebox(0,0){\small{$4$}}}
\put(125,91){\makebox(0,0){\small{$3$}}}
\put(140,91){\makebox(0,0){\small{$2$}}}
\put(155,91){\makebox(0,0){\small{$1$}}}
\put(170,91){\makebox(0,0){\small{$1$}}}
\put(185,91){\makebox(0,0){\small{$0$}}}
\put(200,91){\makebox(0,0){\small{$1$}}}
\put(215,91){\makebox(0,0){\small{$0$}}}

\put(140,-30){\makebox(0,0){\small{The combinatorial weight function}}}

\put(251,-12){\line(1,0){170}}
\put(259,-20){\line(0,1){118}}

\put(265,-5){\makebox(0,0){\small{$0$}}}
\put(280,-5){\makebox(0,0){\small{$1$}}}
\put(295,-5){\makebox(0,0){\small{$0$}}}
\put(310,-5){\makebox(0,0){\small{$1$}}}
\put(325,-5){\makebox(0,0){\small{$0$}}}
\put(340,-5){\makebox(0,0){\small{$1$}}}
\put(355,-5){\makebox(0,0){\small{$2$}}}
\put(370,-5){\makebox(0,0){\small{$3$}}}
\put(385,-5){\makebox(0,0){\small{$4$}}}
\put(400,-5){\makebox(0,0){\small{$5$}}}
\put(415,-5){\makebox(0,0){\small{$6$}}}
\put(319,-10){\framebox(100,23)}

\put(265, 7){\makebox(0,0){\small{$1$}}}
\put(280, 7){\makebox(0,0){\small{$0$}}}
\put(295, 7){\makebox(0,0){\small{$-1$}}}
\put(310, 7){\makebox(0,0){\small{$0$}}}
\put(325, 7){\makebox(0,0){\small{$-1$}}}
\put(340, 7){\makebox(0,0){\small{$0$}}}
\put(355, 7){\makebox(0,0){\small{$1$}}}
\put(370, 7){\makebox(0,0){\small{$2$}}}
\put(385, 7){\makebox(0,0){\small{$3$}}}
\put(400, 7){\makebox(0,0){\small{$4$}}}
\put(415, 7){\makebox(0,0){\small{$5$}}}

\put(265,19){\makebox(0,0){\small{$0$}}}
\put(280,19){\makebox(0,0){\small{$-1$}}}
\put(295,19){\makebox(0,0){\small{$-2$}}}
\put(310,19){\makebox(0,0){\small{$-1$}}}
\put(325,19){\makebox(0,0){\small{$-2$}}}
\put(340,19){\makebox(0,0){\small{$-1$}}}
\put(355,19){\makebox(0,0){\small{$0$}}}
\put(370,19){\makebox(0,0){\small{$1$}}}
\put(385,19){\makebox(0,0){\small{$2$}}}
\put(400,19){\makebox(0,0){\small{$3$}}}
\put(415,19){\makebox(0,0){\small{$4$}}}

\put(265,31){\makebox(0,0){\small{$1$}}}
\put(280,31){\makebox(0,0){\small{$0$}}}
\put(295,31){\makebox(0,0){\small{$-1$}}}
\put(310,31){\makebox(0,0){\small{$-2$}}}
\put(325,31){\makebox(0,0){\small{$-3$}}}
\put(340,31){\makebox(0,0){\small{$-2$}}}
\put(355,31){\makebox(0,0){\small{$-1$}}}
\put(370,31){\makebox(0,0){\small{$0$}}}
\put(385,31){\makebox(0,0){\small{$1$}}}
\put(400,31){\makebox(0,0){\small{$2$}}}
\put(415,31){\makebox(0,0){\small{$3$}}}

\put(265,43){\makebox(0,0){\small{$2$}}}
\put(280,43){\makebox(0,0){\small{$1$}}}
\put(295,43){\makebox(0,0){\small{$0$}}}
\put(310,43){\makebox(0,0){\small{$-1$}}}
\put(325,43){\makebox(0,0){\small{$-2$}}}
\put(340,43){\makebox(0,0){\small{$-1$}}}
\put(355,43){\makebox(0,0){\small{$-2$}}}
\put(370,43){\makebox(0,0){\small{$-1$}}}
\put(385,43){\makebox(0,0){\small{$0$}}}
\put(400,43){\makebox(0,0){\small{$1$}}}
\put(415,43){\makebox(0,0){\small{$2$}}}

\put(265,55){\makebox(0,0){\small{$3$}}}
\put(280,55){\makebox(0,0){\small{$2$}}}
\put(295,55){\makebox(0,0){\small{$1$}}}
\put(310,55){\makebox(0,0){\small{$0$}}}
\put(325,55){\makebox(0,0){\small{$-1$}}}
\put(340,55){\makebox(0,0){\small{$-2$}}}
\put(355,55){\makebox(0,0){\small{$-3$}}}
\put(370,55){\makebox(0,0){\small{$-2$}}}
\put(385,55){\makebox(0,0){\small{$-1$}}}
\put(400,55){\makebox(0,0){\small{$0$}}}
\put(415,55){\makebox(0,0){\small{$1$}}}

\put(265,67){\makebox(0,0){\small{$4$}}}
\put(280,67){\makebox(0,0){\small{$3$}}}
\put(295,67){\makebox(0,0){\small{$2$}}}
\put(310,67){\makebox(0,0){\small{$1$}}}
\put(325,67){\makebox(0,0){\small{$0$}}}
\put(340,67){\makebox(0,0){\small{$-1$}}}
\put(355,67){\makebox(0,0){\small{$-2$}}}
\put(370,67){\makebox(0,0){\small{$-1$}}}
\put(385,67){\makebox(0,0){\small{$-2$}}}
\put(400,67){\makebox(0,0){\small{$-1$}}}
\put(415,67){\makebox(0,0){\small{$0$}}}

\put(265,79){\makebox(0,0){\small{$5$}}}
\put(280,79){\makebox(0,0){\small{$4$}}}
\put(295,79){\makebox(0,0){\small{$3$}}}
\put(310,79){\makebox(0,0){\small{$2$}}}
\put(325,79){\makebox(0,0){\small{$1$}}}
\put(340,79){\makebox(0,0){\small{$0$}}}
\put(355,79){\makebox(0,0){\small{$-1$}}}
\put(370,79){\makebox(0,0){\small{$0$}}}
\put(385,79){\makebox(0,0){\small{$-1$}}}
\put(400,79){\makebox(0,0){\small{$0$}}}
\put(415,79){\makebox(0,0){\small{$1$}}}
\put(261,73){\framebox(100,24)}

\put(265,91){\makebox(0,0){\small{$6$}}}
\put(280,91){\makebox(0,0){\small{$5$}}}
\put(295,91){\makebox(0,0){\small{$4$}}}
\put(310,91){\makebox(0,0){\small{$3$}}}
\put(325,91){\makebox(0,0){\small{$2$}}}
\put(340,91){\makebox(0,0){\small{$1$}}}
\put(355,91){\makebox(0,0){\small{$0$}}}
\put(370,91){\makebox(0,0){\small{$1$}}}
\put(385,91){\makebox(0,0){\small{$0$}}}
\put(400,91){\makebox(0,0){\small{$1$}}}
\put(415,91){\makebox(0,0){\small{$0$}}}
\put(340,-30){\makebox(0,0){\small{The analytic weight function}}}
\end{picture}

(ii) In some examples the analytic and combinatorial  height and weight functions might agree, e.g., for 
$f(x_1,x_2)=(x_1^2+x_2^5)(x_1^5+x_2^2)$ this is the case. In fact, the differences correspond to the cases when we can connect two lattice points of $\mathcal{P}$ by some translate of a primitive edge of the Newton boundary. In such a case, on some edges of the underlying lattice the combinatorial height function jumps by more than one, whereas the analytic Hilbert function by at most one (compare with formula (\ref{eq:Hilbertfelcsopbol})). 
\end{example}

\bekezdes Let $\cF^{an}(\ell)$ and $ \cF^{\Gamma}(\ell)$, $\ell\in 
R(0,{\bf c})\cap \Z^r$, be the corresponding filtrations determined by the 
collections of analytic and combinatorial valuations $\{\tilde{\frv}^{an}_v\}_v$ and 
$\{\frv^\Gamma_v\}_v$ (cf. (\ref{eq:F_D(l)})). 
From the previous discussion, 
$\cF^{an}(\ell) \supset \cF^{\Gamma}(\ell)$ for every lattice point $\ell$. 
However, most importantly, they agree for the conductor element.

\begin{lemma}\label{lem:agreec} In $\cO_{\mathbb{C}^2, 0}$ the following ideals agree:
$\cF^{an}({\bf c})=\cF^\Gamma({\bf c}).$
\end{lemma}
\begin{proof} On one hand $\cF^{an}({\bf c}) \supset \cF^{\Gamma}({\bf c})$, on the other hand, by Lemma \ref{lem:COND}, both of them have the same codimension  $\delta(C,0)=|\mathcal{P}|$, hence they should agree. 
\end{proof}

In conclusion, we got another proof for the main theorem of the section: 

\begin{corollary}[=Theorem \ref{th:comparison}]\label{cor:egyenlo}
For a plane curve singularity $(C, 0) = (\{f=0\}, 0)$ with Newton nondegenerate principal part and convenient Newton diagram
 $$\bH_{an, *}(C,0)\cong\mathbb{SH}_*(\cO_{\mathbb{C}^2, 0}/{\rm adj}(\mathcal{M}))\cong \mathbb{H}_{\ast}(R(0, {\bf c}), w^{\Gamma}),$$
 where $\mathcal{M}$ is the monomial ideal $(\{x^p\,:\, p\in (\mathbb{Z}_{>0})^r\cap N\Gamma_+(f)\}) \triangleleft \mathcal{O}_{\mathbb{C}^2, 0}$.
 \end{corollary}
\begin{proof}
    The ideal $q^{-1}({\mathcal{C}})=\cF^{an}({\bf c})=\cF^\Gamma({\bf c})={\rm adj}(\mathcal{M})$ is realized by two different collection of discrete valuations, namely by  $\{\tilde{\frv}^{an}_v\}_v$ and 
$\{\frv^\Gamma_v\}_v$. They satisfy the CDP condition (see Lemma \ref{lem:CDPcurves}), hence the corresponding homologies agree by Theorem \ref{th:Indep}.
\end{proof}

\subsection{Some more examples}\label{ss:NNexamples}\,

The Independence Theorem \ref{th:Indep}
implies, that the only really important part for the proof of the previous Corollary is the fact that the two set of discrete valuations realize the same ideal, that is, in the terminology of 
Theorem \ref{th:comparison}:  $q^{-1}(\mathcal{C})={\rm adj}(\mathcal{M})$. This can also happen in other, more general situations, e.g., when the germ is \emph{not} Newton nondegenerate, or when 
we compare two germs with different Newton boundaries. We will give some examples below.

\bekezdes\label{bek:10.2.1} \textbf{Example: different  Newton boundary shapes but the same conductor ideal.} \,
We can have different monomial ideals $\mathcal{M}\neq \mathcal{M}'$ such that their adjoints agree. Equivalently, in the plane curve singularity language, we can have two different Newton nondegenerate curve singularities (with different Newton boundaries and different embedded topological types) but with agreeing conductor ideals (more precisely, their pullbacks in  $\cO_{\mathbb{C}^2,0}$ along the natural projection maps $\mathcal{O}_{\mathbb{C}^2, 0} \rightarrow \mathcal{O}_{C,0}$ agree). Indeed, this is just the case if the set $(\mathbb{Z}_{> 0})^2 \setminus \mathcal{P}$ of positive coordinate lattice points lying strictly above the $-\mathbf{1}$-translate of the Newton boundary agrees for the two singularities, i.e., if the corresponding $\mathcal{P}$ sets, defined in paragraph \ref{par:NewtonComb}, agree.

Let us see a concrete example: ${\color{blue} (C_1, 0):=(V_{(x_1^3-x_2^2)(x_1^2-x_2^3)}, 0)}$  and  ${\color{red} (C_2, 0):=(V_{x_1^5-x_2^4}, 0)}$, where $V_f$ denotes the vanishing locus $\{f=0\}$ of the germ $f \in \mathcal{O}_{\mathbb{C}^2, 0}$. They have different number of components, but the corresponding $\mathcal{P}$ sets agree. We can easily compute their combinatorial weight functions and see that the lattice homologies agree:

\begin{picture}(420,180)(-20,-70)

\linethickness{0.05mm}
  \multiput(-10,15)(0,15){6}{\line(1,0){105}}
  \multiput(15,-10)(15,0){6}{\line(0,1){105}}

\linethickness{0.4mm}
\put(-10,0){\line(1,0){105}}
\put(0,-10){\line(0,1){105}}
\linethickness{0.4mm}
\put(30,30){\color{blue}\line(3,-2){45}}
\put(30,30){\color{blue}\line(-2,3){30}}
\put(0,75){\color{blue}\circle*{4}}
\put(30,30){\color{blue}\circle*{4}}
\put(75, 0){\color{blue}\circle*{4}}
\put(45,45){\color{blue}\circle*{4}}

\put(15,15){\color{blue}\line(3,-2){37.5}}
\put(15,15){\color{blue}\line(-2,3){25}}

\put(75, 0){\color{red}\line(-5,4){75}} 
\put(53.75,-10){\color{red}\line(-5,4){63.75}} 
\put(0,60){\color{red}\circle*{4}}
\put(75, 0){\color{red}\circle*{4}}

  \linethickness{0.05mm}

\put(15,15){\circle*{4}}
\put( 0,15){\circle*{4}}
\put( 0,30){\circle*{4}}

\put( 0, 0){\circle*{4}}
\put(15, 0){\circle*{4}}
\put(30, 0){\circle*{4}}

\color{blue}
\put(120,00){\makebox(0,0){\small{$(0,0)$}}}
\put(120,15){\makebox(0,0){\small{$(3,2)$}}}
\put(120,30){\makebox(0,0){\small{$(6,4)$}}}
\put(143,00){\makebox(0,0){\small{$(2,3)$}}}
\put(143,15){\makebox(0,0){\small{$(5,5)$}}}
\put(166,00){\makebox(0,0){\small{$(4,6)$}}}

\put(212, 0){\line(1,0){116}}
\put(220,-8){\line(0,1){80}}

\put(228,5){\makebox(0,0){\small{$0$}}}
\put(243,5){\makebox(0,0){\small{$1$}}}
\put(258,5){\makebox(0,0){\small{$0$}}}
\put(273,5){\makebox(0,0){\small{$1$}}}
\put(288,5){\makebox(0,0){\small{$2$}}}
\put(303,5){\makebox(0,0){\small{$3$}}}
\put(318,5){\makebox(0,0){\small{$4$}}}

\put(228,15){\makebox(0,0){\small{$1$}}}
\put(243,15){\makebox(0,0){\small{$0$}}}
\put(258,15){\makebox(0,0){\small{$-1$}}}
\put(273,15){\makebox(0,0){\small{$0$}}}
\put(288,15){\makebox(0,0){\small{$1$}}}
\put(303,15){\makebox(0,0){\small{$2$}}}
\put(318,15){\makebox(0,0){\small{$3$}}}

\put(228,25){\makebox(0,0){\small{$0$}}}
\put(243,25){\makebox(0,0){\small{$-1$}}}
\put(258,25){\makebox(0,0){\small{$-2$}}}
\put(273,25){\makebox(0,0){\small{$-1$}}}
\put(288,25){\makebox(0,0){\small{$0$}}}
\put(303,25){\makebox(0,0){\small{$1$}}}
\put(318,25){\makebox(0,0){\small{$2$}}}

\put(228,35){\makebox(0,0){\small{$1$}}}
\put(243,35){\makebox(0,0){\small{$0$}}}
\put(258,35){\makebox(0,0){\small{$-1$}}}
\put(273,35){\makebox(0,0){\small{$0$}}}
\put(288,35){\makebox(0,0){\small{$-1$}}}
\put(303,35){\makebox(0,0){\small{$0$}}}
\put(318,35){\makebox(0,0){\small{$1$}}}

\put(228,45){\makebox(0,0){\small{$2$}}}
\put(243,45){\makebox(0,0){\small{$1$}}}
\put(258,45){\makebox(0,0){\small{$0$}}}
\put(273,45){\makebox(0,0){\small{$-1$}}}
\put(288,45){\makebox(0,0){\small{$-2$}}}
\put(303,45){\makebox(0,0){\small{$-1$}}}
\put(318,45){\makebox(0,0){\small{$0$}}}

\put(228,55){\makebox(0,0){\small{$3$}}}
\put(243,55){\makebox(0,0){\small{$2$}}}
\put(258,55){\makebox(0,0){\small{$1$}}}
\put(273,55){\makebox(0,0){\small{$0$}}}
\put(288,55){\makebox(0,0){\small{$-1$}}}
\put(303,55){\makebox(0,0){\small{$0$}}}
\put(318,55){\makebox(0,0){\small{$1$}}}

\put(228,65){\makebox(0,0){\small{$4$}}}
\put(243,65){\makebox(0,0){\small{$3$}}}
\put(258,65){\makebox(0,0){\small{$2$}}}
\put(273,65){\makebox(0,0){\small{$1$}}}
\put(288,65){\makebox(0,0){\small{$0$}}}
\put(303,65){\makebox(0,0){\small{$1$}}}
\put(318,65){\makebox(0,0){\small{$0$}}}

\put(221,1){\framebox(14,9)}
\put(251,21){\framebox(14,9)}
\put(281,41){\framebox(14,9)}
\put(311,61){\framebox(14,9)}

\color{red}
\put(120,60){\makebox(0,0){\small{$(0)$}}}
\put(120,75){\makebox(0,0){\small{$(5)$}}}
\put(120,90){\makebox(0,0){\small{$(10)$}}}
\put(143,60){\makebox(0,0){\small{$(4)$}}}
\put(143,75){\makebox(0,0){\small{$(9)$}}}
\put(166,60){\makebox(0,0){\small{$(8)$}}}

\put(172,90){\line(1,0){206}}
\put(180,82){\line(0,1){20}}

\put(188,95){\makebox(0,0){\small{$0$}}}
\put(203,95){\makebox(0,0){\small{$1$}}}
\put(218,95){\makebox(0,0){\small{$0$}}}
\put(233,95){\makebox(0,0){\small{$-1$}}}
\put(248,95){\makebox(0,0){\small{$-2$}}}
\put(263,95){\makebox(0,0){\small{$-1$}}}
\put(278,95){\makebox(0,0){\small{$0$}}}
\put(293,95){\makebox(0,0){\small{$-1$}}}
\put(308,95){\makebox(0,0){\small{$-2$}}}
\put(323,95){\makebox(0,0){\small{$-1$}}}
\put(338,95){\makebox(0,0){\small{$0$}}}
\put(353,95){\makebox(0,0){\small{$1$}}}
\put(368,95){\makebox(0,0){\small{$0$}}}
\put(181,91){\framebox(14,9)}
\put(241,91){\framebox(14,9)}
\put(301,91){\framebox(14,9)}
\put(361,91){\framebox(14,9)}

\color{black}

\put(150,-45){\makebox(0,0){\small{\color{blue}$\mathbb{H}_{\ast}(C_1, 0)$\color{black} $\cong$ \color{red} $\mathbb{H}_{\ast}(C_2, 0)$ \color{black} $\cong \mathbb{H}_0$ of}}}

\put(260,-45){\makebox(0,0){\small{$0$}}} 
\put(257,-35){\makebox(0,0){\small{$-1$}}}
\put(257,-25){\makebox(0,0){\small{$-2$}}}

\put(230,-65){\makebox(0,0){$\vdots$}} 
\put(230,-61){\line(0,1){16}}

\put(230,-55){\circle*{3}} 
\put(230,-55){\line(1,1){10}} 
\put(230,-55){\line(-1,1){10}}

\put(220,-45){\circle*{3}}
\put(230,-45){\circle*{3}} 
\put(240,-45){\circle*{3}}

\put(230,-45){\line(1,2){5}} 
\put(230,-45){\line(-1,2){5}}

\put(225,-35){\circle*{3}} \put(225,-35){\line(0,1){10}}
\put(235,-35){\circle*{3}} \put(235,-35){\line(0,1){10}}

\put(225,-25){\circle*{3}}
\put(235,-25){\circle*{3}}
\end{picture}

\bekezdes\label{bek:10.2.2} \textbf{Example: monomial ideals and Rees valuations.} Given a Noetherian ring $\mathcal{O}$ and an ideal $\mathcal{I} \triangleleft\mathcal{O}$, Rees proved the existence of a finite minimal set of discrete valuations $\mathcal{D}^R=\{\mathfrak{v}_1^R, \ldots, \mathfrak{v}_r^R\}$ on $\mathcal{O}$, such that there exists some lattice point $d^R_{\mathcal{I}}\in (\mathbb{Z}_{\geq 0})^r$ with the property that for every $n \in \mathbb{N}$: the pair $(\mathcal{D}^R, n\cdot d^R_{\mathcal{I}})$ is a realization of the integral closure $\overline{\mathcal{I}^n}$ of the power $\mathcal{I}^n$ \cite{Rees}. For a monomial ideal $\mathcal{M} \triangleleft \mathcal{O}_m$ in a polynomial ring, these, so-called `Rees valuations' are the monomial valuations corresponding to the facets of the Newton diagram of $\mathcal{M}$. Although these give a good starting point to compute the symmetric lattice homology of $\overline{\mathcal{M}}$, we can easily produce examples when the integrally closed monomial ideal can be realized by strictly less valuations than the number of Rees valuations. This shows that ${\rm homdim}(\overline{\mathcal{I}} \triangleleft \mathcal{O})$ (for the notation see Theorem \ref{th:propertiesideal} \textit{(c)}) can be strictly smaller than the number of Rees valuations. It is an interesting new invariant of the finite codimensional ideal. In fact, for monomial ideals we can characterize it via geometric properties of the Newton polygon, see section \ref{s:homdim}.

For a concrete example we can go back to Example \ref{ex:2552ND} discussing the integrally closed monomial ideal ${\rm adj}(\mathcal{M}) = \overline{(x_1^5,\, x_1 x_2^2,\, x_2^3)} \triangleleft \mathcal{O}_2$. Hence the monomial Rees valuations correspond to the primitive normal vectors $(1, 1)$ and $(1, 2)$. On the other hand, ${\rm adj}(\mathcal{M}) = (\{ x^p\,:\, \langle (4, 7), p\rangle \geq 18\})$. This fact corresponds to the isomorphism
$$ \mathbb{H}_{an, *}({\color{blue}V_{(x_1^5-x_2^2)(x^3_1 - x_2^2)}, 0}) \cong \mathbb{SH}_*({\rm adj}(\mathcal{M}) \triangleleft \mathcal{O}_2) \cong \mathbb{H}_{an,*}({\color{red}V_{x_1^7-x_2^4}, 0}).$$
One can also verify this identity directly by comparing the following weight function with the one in Example \ref{ex:2552ND}

\begin{picture}(420,125)(-20,-20)

\linethickness{0.5mm}
\put(-15,0){\line(1,0){165}}
\put(0,-15){\line(0,1){105}}
\linethickness{0.35mm}
\put(45,30){\color{blue}\line(5,-2){75}}
\put(45,30){\color{blue}\line(-3,2){45}}
\put( 0,60){\color{red}\line(7,-4){105}}
\put(-15,45){\line(7,-4){105}}
\put(30,15){\line(5,-2){75}}
\put(30,15){\line(-3,2){45}}

\linethickness{0.05mm}
  \multiput(-15,15)(0,15){5}{\line(1,0){165}}
  \multiput(15,-15)(15,0){9}{\line(0,1){105}}

\put(15,30){\circle*{4}}
\put( 0,45){\circle*{4}}
\put(45,30){\color{blue}\circle*{4}}
\put(75, 0){\circle*{4}}
\put(0, 60){\color{blue}\circle*{4}}
\put(120, 0){\color{blue}\circle*{4}}
\put(75,30){\color{blue}\circle*{4}}
\put(105, 0){\color{red}\circle*{4}}
\put(0, 60){\color{red}\circle*{4}}

\color{red}
\put(180,00){\makebox(0,0){\small{$(0)$}}}
\put(180,15){\makebox(0,0){\small{$(7)$}}}
\put(180,30){\makebox(0,0){\small{$(14)$}}}
\put(203,00){\makebox(0,0){\small{$(4)$}}}
\put(203,15){\makebox(0,0){\small{$(11)$}}}
\put(226,00){\makebox(0,0){\small{$(8)$}}}
\put(226,15){\makebox(0,0){\small{$(15)$}}}
\put(249,00){\makebox(0,0){\small{$(12)$}}}
\put(272,00){\makebox(0,0){\small{$(16)$}}}

\put(172,75){\line(1,0){206}}
\put(180,67){\line(0,1){20}}

\put(188,80){\makebox(0,0){\small{$0$}}}
\put(203,80){\makebox(0,0){\small{$1$}}}
\put(218,80){\makebox(0,0){\small{$0$}}}
\put(233,80){\makebox(0,0){\small{$-1$}}}
\put(248,80){\makebox(0,0){\small{$-2$}}}
\put(263,80){\makebox(0,0){\small{$-1$}}}
\put(278,80){\makebox(0,0){\small{$-2$}}}
\put(293,80){\makebox(0,0){\small{$-3$}}}
\put(308,80){\makebox(0,0){\small{$-2$}}}
\put(323,80){\makebox(0,0){\small{$-1$}}}
\put(338,80){\makebox(0,0){\small{$-2$}}}
\put(353,80){\makebox(0,0){\small{$-3$}}}
\put(368,80){\makebox(0,0){\small{$-2$}}}
\put(181,76){\framebox(14,9)}
\put(241,76){\framebox(14,9)}
\put(287,76){\framebox(14,9)}
\put(347,76){\framebox(14,9)}

\put(283,55){\line(1,0){95}}
\put(293,60){\makebox(0,0){\small{$-1$}}}
\put(308,60){\makebox(0,0){\small{$-2$}}}
\put(323,60){\makebox(0,0){\small{$-1$}}}
\put(338,60){\makebox(0,0){\small{$0$}}}
\put(353,60){\makebox(0,0){\small{$1$}}}
\put(368,60){\makebox(0,0){\small{$0$}}}
\put(302,56){\framebox(14,9)}
\put(362,56){\framebox(14,9)}
\end{picture}

Notice, that in the previous two examples the symmetric lattice homology of the ideals has trivial $\mathbb{SH}_{\geq 1}$ and, in compliance with Conjecture \ref{conj:elso} (or Theorem \ref{th:homdim}), we indeed managed to find a single monomial valuation realizing the respective ideals. For more discussions, results and examples on this topic see section \ref{s:homdim} and subsection \ref{ss:HOMDIM}.

\bekezdes\label{bek:10.6.3} \textbf{Example: Newton degenerate case still determined by the Newton diagram.}
Consider the irreducible  plane curve singularity $(C, 0)\subset (\mathbb{C}^2, 0)$ with Puiseux parametrisation $x_1=t^4, \ x_2=t^6 + t^7$ and equation $f(x_1,x_2)=(x_1^3-x_2^2)^2-4x_1^5x_2-x_1^7=0$. This is clearly a Newton degenerate curve, being irreducible with two Puiseux pairs.
Its semigroup of values is $\mathcal{S}_{C, 0}=\langle 4, 6, 13 \rangle$ with analytic weight function $w_{an}$ in $R(0,\mathbf{c})=R(0,16)$:

\begin{center}
\begin{picture}(320,25)(-20,-10)
\put(-8,-2){\line(1,0){268}}
\put(0,-8){\line(0,1){20}}

\put(10,5){\makebox(0,0){\small{$0$}}}
\put(25,5){\makebox(0,0){\small{$1$}}}
\put(40,5){\makebox(0,0){\small{$0$}}}
\put(55,5){\makebox(0,0){\small{$-1$}}}
\put(70,5){\makebox(0,0){\small{$-2$}}}
\put(85,5){\makebox(0,0){\small{$-1$}}}
\put(100,5){\makebox(0,0){\small{$-2$}}}
\put(115,5){\makebox(0,0){\small{$-1$}}}
\put(130,5){\makebox(0,0){\small{$-2$}}}
\put(145,5){\makebox(0,0){\small{$-1$}}}
\put(160,5){\makebox(0,0){\small{$-2$}}}
\put(175,5){\makebox(0,0){\small{$-1$}}}
\put(190,5){\makebox(0,0){\small{$-2$}}}
\put(205,5){\makebox(0,0){\small{$-1$}}}
\put(220,5){\makebox(0,0){\small{$0$}}}
\put(235,5){\makebox(0,0){\small{$1$}}}
\put(250,5){\makebox(0,0){\small{$0$}}}
\put(3,0){\framebox(14,10)}
\put(63,0){\framebox(14,10)}
\put(93,0){\framebox(14,10)}
\put(123,0){\framebox(14,10)}
\put(153,0){\framebox(14,10)}
\put(183,0){\framebox(14,10)}
\put(243,0){\framebox(14,10)}
\end{picture}
\end{center}
On the other hand, we could start from its Newton boundary as well, coloured in blue in the following diagram. If we associate to it the monomial ideal $\mathcal{M}$, as before, we can compute $\mathbb{SH}_*({\rm adj}(\mathcal{M})\triangleleft\mathcal{O}_2)$ by  the combinatorial degree function (cf. (\ref{eq:degfctn})):
$${\rm deg}:\mathbb{Z}^2 \rightarrow \mathbb{Z}^2, \  {\rm deg}(p)=(\, \langle (2, 3), p-{\bf 1}\rangle, 
 \langle (2, 3), p-{\bf 1}\rangle \, )\in \mathbb{Z}^2$$
 (corresponding to the primitive edge decomposition of the boundary) and combinatorial conductor $\mathbf{c}^{{\Gamma}}=(8,8)$. We obtain the following weight function:
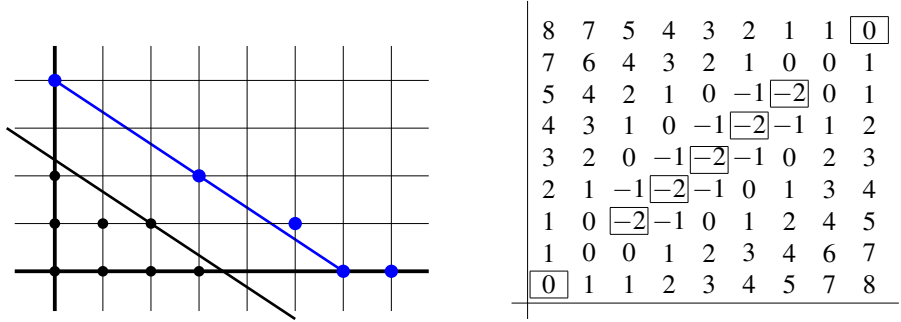
\begin{center}
    \begin{picture}(340,125)(-40,-20)

\linethickness{0.5mm}
\put(-35,0){\line(1,0){155}}
\put(-20,-15){\line(0,1){100}}
\linethickness{0.35mm}
\put(34,36){\color{blue}\line(3,-2){54}}
\put(34,36){\color{blue}\line(-3,2){54}}

\put(16,18){\line(3,-2){54}}
\put(16,18){\line(-3,2){54}}

\linethickness{0.05mm}
  \multiput(-35,18)(0,18){4}{\line(1,0){155}}
  \multiput(-2,-15)(18,0){7}{\line(0,1){100}}

\put(34,36){\color{blue}\circle*{5}}
\put(-20,72){\color{blue}\circle*{5}}
\put(88, 0){\color{blue}\circle*{5}}
\put(70,18){\color{blue}\circle*{5}}
\put(106, 0){\color{blue}\circle*{5}}
\put(-20,18){\circle*{4}}
\put(-20,36){\circle*{4}}
\put(-20, 0){\circle*{4}}
\put(16,18){\circle*{4}}
\put(-2,18){\circle*{4}}
\put(16,0){\circle*{4}}
\put(34,0){\circle*{4}}
\put(-2,0){\circle*{4}}

\put(151,-12){\line(1,0){145}}
\put(157,-18){\line(0,1){120}}

\put(165,-5){\makebox(0,0){\small{$0$}}}
\put(180,-5){\makebox(0,0){\small{$1$}}}
\put(195,-5){\makebox(0,0){\small{$1$}}}
\put(210,-5){\makebox(0,0){\small{$2$}}}
\put(225,-5){\makebox(0,0){\small{$3$}}}
\put(240,-5){\makebox(0,0){\small{$4$}}}
\put(255,-5){\makebox(0,0){\small{$5$}}}
\put(270,-5){\makebox(0,0){\small{$7$}}}
\put(285,-5){\makebox(0,0){\small{$8$}}}
\put(158,-10){\framebox(14,10)}

\put(165, 7){\makebox(0,0){\small{$1$}}}
\put(180, 7){\makebox(0,0){\small{$0$}}}
\put(195, 7){\makebox(0,0){\small{$0$}}}
\put(210, 7){\makebox(0,0){\small{$1$}}}
\put(225, 7){\makebox(0,0){\small{$2$}}}
\put(240, 7){\makebox(0,0){\small{$3$}}}
\put(255, 7){\makebox(0,0){\small{$4$}}}
\put(270, 7){\makebox(0,0){\small{$6$}}}
\put(285, 7){\makebox(0,0){\small{$7$}}}

\put(165,19){\makebox(0,0){\small{$1$}}}
\put(180,19){\makebox(0,0){\small{$0$}}}
\put(195,19){\makebox(0,0){\small{$-2$}}}
\put(210,19){\makebox(0,0){\small{$-1$}}}
\put(225,19){\makebox(0,0){\small{$0$}}}
\put(240,19){\makebox(0,0){\small{$1$}}}
\put(255,19){\makebox(0,0){\small{$2$}}}
\put(270,19){\makebox(0,0){\small{$4$}}}
\put(285,19){\makebox(0,0){\small{$5$}}}
\put(188,14){\framebox(14,10)}

\put(165,31){\makebox(0,0){\small{$2$}}}
\put(180,31){\makebox(0,0){\small{$1$}}}
\put(195,31){\makebox(0,0){\small{$-1$}}}
\put(210,31){\makebox(0,0){\small{$-2$}}}
\put(225,31){\makebox(0,0){\small{$-1$}}}
\put(240,31){\makebox(0,0){\small{$0$}}}
\put(255,31){\makebox(0,0){\small{$1$}}}
\put(270,31){\makebox(0,0){\small{$3$}}}
\put(285,31){\makebox(0,0){\small{$4$}}}
\put(203,26){\framebox(14,10)}

\put(165,43){\makebox(0,0){\small{$3$}}}
\put(180,43){\makebox(0,0){\small{$2$}}}
\put(195,43){\makebox(0,0){\small{$0$}}}
\put(210,43){\makebox(0,0){\small{$-1$}}}
\put(225,43){\makebox(0,0){\small{$-2$}}}
\put(240,43){\makebox(0,0){\small{$-1$}}}
\put(255,43){\makebox(0,0){\small{$0$}}}
\put(270,43){\makebox(0,0){\small{$2$}}}
\put(285,43){\makebox(0,0){\small{$3$}}}
\put(218,37){\framebox(14,10)}

\put(165,55){\makebox(0,0){\small{$4$}}}
\put(180,55){\makebox(0,0){\small{$3$}}}
\put(195,55){\makebox(0,0){\small{$1$}}}
\put(210,55){\makebox(0,0){\small{$0$}}}
\put(225,55){\makebox(0,0){\small{$-1$}}}
\put(240,55){\makebox(0,0){\small{$-2$}}}
\put(255,55){\makebox(0,0){\small{$-1$}}}
\put(270,55){\makebox(0,0){\small{$1$}}}
\put(285,55){\makebox(0,0){\small{$2$}}}
\put(233,50){\framebox(14,10)}

\put(165,67){\makebox(0,0){\small{$5$}}}
\put(180,67){\makebox(0,0){\small{$4$}}}
\put(195,67){\makebox(0,0){\small{$2$}}}
\put(210,67){\makebox(0,0){\small{$1$}}}
\put(225,67){\makebox(0,0){\small{$0$}}}
\put(240,67){\makebox(0,0){\small{$-1$}}}
\put(255,67){\makebox(0,0){\small{$-2$}}}
\put(270,67){\makebox(0,0){\small{$0$}}}
\put(285,67){\makebox(0,0){\small{$1$}}}
\put(248,62){\framebox(14,10)}

\put(165,79){\makebox(0,0){\small{$7$}}}
\put(180,79){\makebox(0,0){\small{$6$}}}
\put(195,79){\makebox(0,0){\small{$4$}}}
\put(210,79){\makebox(0,0){\small{$3$}}}
\put(225,79){\makebox(0,0){\small{$2$}}}
\put(240,79){\makebox(0,0){\small{$1$}}}
\put(255,79){\makebox(0,0){\small{$0$}}}
\put(270,79){\makebox(0,0){\small{$0$}}}
\put(285,79){\makebox(0,0){\small{$1$}}}

\put(165,91){\makebox(0,0){\small{$8$}}}
\put(180,91){\makebox(0,0){\small{$7$}}}
\put(195,91){\makebox(0,0){\small{$5$}}}
\put(210,91){\makebox(0,0){\small{$4$}}}
\put(225,91){\makebox(0,0){\small{$3$}}}
\put(240,91){\makebox(0,0){\small{$2$}}}
\put(255,91){\makebox(0,0){\small{$1$}}}
\put(270,91){\makebox(0,0){\small{$1$}}}
\put(285,91){\makebox(0,0){\small{$0$}}}
\put(278,86){\framebox(14,10)}
\end{picture}
\end{center}

It turns out that the analytical lattice cohomology agrees with this combinatorial one:
\begin{equation*}
    \mathbb{H}_{an, \ast}(C, 0) \cong \mathbb{H}_{\ast}(R(0, \mathbf{c}^{\Gamma}), w^{\Gamma}) \cong \mathcal{T}_4^- \oplus (\mathcal{T}_4(1))^{\oplus 4} \oplus (\mathcal{T}_0(1))^{\oplus 2}.
\end{equation*}
Indeed, one can also show that the corresponding discrete valuations (and conductor elements) realize the same ideal $\mathcal{I} \triangleleft \mathbb{C}\{ x_1, x_2\}$, so the above isomorphism is, in fact, a corollary of the Independence Theorem. For a formal proof, one would have to use once again the inequality between the analytic and monomial valuations and the identification of the codimensions. 

Moreover, 
whenever $(C', 0):=(V_{x_1^6+ax_1^3x_2^2+x_2^4}, 0)$ is a Newton nondegenerate plane curve singularity (i.e., if $a \neq \pm2$), one also has 
\begin{equation}\label{eq:ISOUJ}
    \mathbb{H}_{an,\ast}(C', 0) \cong \mathbb{H}_{\ast}(R(0, \mathbf{c}^{\Gamma}), w^{\Gamma}) \cong \mathbb{H}_{an, \ast}(C, 0).
\end{equation}

\bekezdes  {\bf Example (Example \ref{bek:10.6.3} generalized).}
\label{bek:10.6.5} 
Let $(C,0)$ be the irreducible Newton degenerate plane curve singularity with Puiseux parametrization 
$x_1=t^{p_1p_2}$, $x_2=t^{p_2q_1}(1+t)$ (whose semigroup is generated by $p_1p_2, 
\ p_2q_1,\ p_1p_2q_1+1 $), where $2\leq p_1<q_1$, $(p_1,q_1)=1$, $2\leq p_2$.
Then the analytic lattice homology of $(C,0)$ agrees with the  combinatorial lattice homology associated with the Newton diagram generated by the lattice points $(p_1p_2,0)$ and $(0, p_2q_1)$. 

This is the only case (up to equisingularity) of an irreducible plane curve singularity with two Puiseux pairs when such a 
coincidence holds. In fact, this is due to the fact, that this is only case when we have the coincidence of the delta invariant with the cardinality $\vert \mathcal{P}\vert$ of positive coordinate lattice points lying strictly above the $-\mathbf{1}$-translate of the Newton boundary.

For some examples in the Newton nondegenerate surface singularity case, see subsection \ref{ss:CondIdConst}.

\section{Lattice homological dimension and monomial ideals}\label{s:homdim}

Given a Noetherian algebra $\mathcal{O}$ over a field $k$ and a finitely generated $\mathcal{O}$-module $M$, for any realizable submodule $N \leq M$ we defined its (lattice) homological dimension in (\ref{eq:defhomdim}) as 
\begin{equation*}
    {\rm homdim}(N \hookrightarrow M):=\begin{cases}
        \max \{ q\,:\, \mathbb{H}_{q}(N \hookrightarrow M)\not=0\} & \text{ if }N \neq M;\\
        -1 & \text{ if }N=M.
    \end{cases}
\end{equation*} 
By construction, the cardinality of any realization of $N$ gives a strict upper bound, more precisely 
\begin{center}
    ${\rm homdim}(N \hookrightarrow M) \leq \min \{|\mathcal{D}|\,:\, \mathcal{D} \text{ a realization of } N\} -1$ \hspace{10mm} (cf. Proposition \ref{prop:homdegupperbound}).
\end{center}

In this subsection we present some conditions under which the above inequality becomes an identity. More concretely, we assign to certain combinations of ring and module elements non-vanishing maximal dimensional homological cycles. In the case when $M=\mathcal{O}=k[x_1, x_2]$ and $N=\mathcal{M}$ is an integrally closed monomial ideal we can produce such ring elements --- in fact, monomials --- using the geometry of the corresponding Newton polygon. Thus, in this case we can understand the homological dimension by means of convex geometry. By the results of the previous section (see especially Corollary \ref{cor:egyenlo}) this also yields a Newton diagrammatic characterization of the \emph{analytic} lattice homological dimension of plane curve singularities with Newton nondegenerate principal part.

\subsection{A condition for the non-vanishing of the top homology}\label{ss:tophomnon0}\,

In this subsection we present how the existence of a certain combination of ring and module elements can bound both the cardinality of any realization and the homological dimension. We will use the following general setting: let $k$ be any field, $\mathcal{O}$ a Noetherian $k$-algebra, $M$ a finitely generated $\mathcal{O}$-module and $N \leq M$ a realizable submodule in the sense of Definition \ref{def:REAL}. 

\begin{prop}\label{prop:htopnem0}
    Fix a nonnegative integer $r\in \mathbb{Z}_{\geq 0}$. Suppose that there exist certain ring elements $f_1, f_2, \ldots, f_r \in \mathcal{O} \setminus{\rm Ann}_{\mathcal{O}}(M/N)$ and module elements $m_1, m_2, \ldots, m_r \in M\setminus N$ satisfying the following:
    \begin{itemize}
        \item $f_im_i \notin N$, for all $i \in \{1, \ldots, r\}$, while
        \item $f_im_j\in N$, for all $i, j \in \{1, \ldots, r\}, \ i \neq j$.
    \end{itemize}
    Then $N$ cannot be realized by less than $r$ extended valuations. Moreover, if there exists a CDP realization of cardinality $r$, then, by Proposition \ref{prop:homdegupperbound}, we get that
    \begin{equation*}
        {\rm homdim}(N \hookrightarrow M)=r-1=\min\{|\mathcal{D}|\,:\, \mathcal{D} \text{ a realization of }N\}-1.
    \end{equation*}
    More precisely, we prove that $\mathbb{H}_{{\rm red},{r-1}}(N \hookrightarrow M)\neq 0$ in the case when $r \geq 1$.
\end{prop} 

For the sake of completeness and simpler referencing we also present the symmetric version of the previous statement (for notations see subsection \ref{ss:Indep}):

\begin{cor}\label{cor:htopnem0sym}
    Let $\mathcal{I} \triangleleft \mathcal{O}$ be a finite codimensional integrally closed ideal and suppose that there exist $f_1, f_2, \ldots, f_r, g_1, g_2, \ldots, g_r \in \mathcal{O} \setminus\mathcal{I}$ ring elements (for some $r\in \mathbb{Z}_{\geq 0}$) satisfying the following:
    \begin{itemize}
        \item  $f_ig_i \notin \mathcal{I}$, for all $i \in \{1, \ldots, r\}$, while
        \item $f_ig_j\in \mathcal{I}$, for all $i, j \in \{1, \ldots, r\}, \ i \neq j$.
    \end{itemize}
    Then $\mathcal{I}$ cannot be realized by less than $r$ discrete valuations. Moreover, if there exists a CDP realization of cardinality $r$, then, by Theorem \ref{th:propertiesideal} \textit{(c)}, we get that
    \begin{equation*}
        {\rm homdim}(\mathcal{I} \triangleleft \mathcal{O})=r-1=\min\{|\mathcal{D}|\,:\, (\mathcal{D}, d) \text{ a realization of }\mathcal{I}\}-1.
    \end{equation*}
\end{cor}

\begin{example}
    In the setting of Example \ref{ex:product} (i.e. $\mathcal{O}=k[x]\times k[y]$ and $\mathcal{I}=(x^3)\times (y^2)$) such combination of ring elements can be chosen, e.g., as $f_1=g_1:=(1,0)$ and $f_2=g_2:=(0,1)$. Also, $\mathcal{I}$ can be realized by $2$ valuations and ${\rm homdim}(\mathcal{I}\triangleleft\mathcal{O})=1$.
\end{example}

\begin{example}\label{ex:2552H1} 
Consider the plane curve singularity $(C, 0)= (\{ (x^2 - y^5)(x^5-y^2)=0\}, 0) \subset (\mathbb{C}^2, 0)$ with Newton nondegenerate principal part from Example \ref{ex:2552} (see also \cite[Example 4.5.2]{AgNeCurves}). Its analytic lattice homology $\mathbb{H}_{an,*}(C, 0)$ agrees with the symmetric lattice homology of the monomial adjunction conductor ideal $\mathcal{C}_f=(x^4, x^2y, xy^2, y^4) \triangleleft\mathbb{C}[x,y]=\mathcal{O}_2$ (see Theorem \ref{th:MerleTeissier} and Corollary \ref{cor:egyenlo}). In this setting we can choose the following combination of ring elements $$f_1:=x, \ f_2:=y, \ g_1:=x^2, \ g_2:=y^2,$$ which satisfy $f_1g_1, \ f_2g_2 \notin \mathcal{C}_f$ and $f_1g_2, \ f_2g_1\in \mathcal{C}_f$. The two discrete valuations coming from the normalization of the components (cf. (\ref{eq:normalizationvaluation})) give a CDP realization of $\mathcal{C}_f$ and we can clearly see from the weight table in Example \ref{ex:2552} that $\mathbb{H}_{an, 1}(C, 0) \neq 0$.
\end{example}

\begin{proof}[Proof of Proposition \ref{prop:htopnem0}]
    Let $\mathcal{D}'=\{ \mathfrak{v}_1, \ldots, \mathfrak{v}_{r'}\}$ be a CDP realization of $N$. For a given $i \in \{1, \ldots, r\}$ the non-containment relation $f_im_i \notin N=\mathcal{F}_{\mathcal{D}'}^M(0)$ implies that there exists some index $v(i)\in \{ 1, \ldots , r'\}$ such that 
    \begin{equation}\label{eq:vv(fvmv)}
        \mathfrak{v}_{v(i)}^M(f_im_i)=\mathfrak{v}_{v(i)}(f_i)+\mathfrak{v}_{v(i)}^M(m_i)<0.
    \end{equation}
    Moreover, since for every $j\in \{1, \ldots, r\}, \ i \neq j$ we have $f_im_j, f_jm_i\in N$, we also get that 
    \begin{align*}
        \forall j\in \{1, \ldots, r\}, \ i \neq j:&\ \mathfrak{v}_{v(i)}(f_i)+\mathfrak{v}_{v(i)}^M(m_j)\geq 0 \\
        &\Rightarrow \mathfrak{v}^{M}_{v(i)}(m_i)<\min\{ \mathfrak{v}^M_{v(i)}(m_j)\,:\,i \neq
        j \in\mathcal{V}\};\\
        \forall j\in \{1, \ldots, r\}, \ i \neq j:&\ \mathfrak{v}_{v(i)}(f_j)+\mathfrak{v}_{v(i)}^M(m_i)\geq 0\\
        &\Rightarrow \mathfrak{v}_{v(i)}(f_i)<\min\{ \mathfrak{v}_{v(i)}(f_j)\,:\,i \neq
        j \in\mathcal{V}\}.
    \end{align*}
    Thus $v(i) \neq v(j)$ if $i \neq j$, hence, any such realization $\mathcal{D}'$ must satisfy the identity $|\mathcal{D}'|=r' \geq r$.

    For the second part let us suppose that $\mathcal{D}=\{\mathfrak{v}_1, \ldots, \mathfrak{v}_r\}$ is a CDP realization of $N$ (with cardinality $|\mathcal{D}|=r$) and denote the tuples $(\mathfrak{v}_1, \ldots, \mathfrak{v}_r)$ and $(\mathfrak{v}_1^M, \ldots, \mathfrak{v}^M_r)$ by $\mathfrak{v}$ and $\mathfrak{v}^M$ respectively. Also, let us denote the index set for simplicity by $\mathcal{V}:=\{1, \ldots, r\}$. For $r=0$ we have $N=M$ and the statement is automatic by definition, so we can turn to the case of $r\geq 1$.

    By simple reordering, we can also assume that $v(i)=i$ for all indices $1\leq i \leq r$, i.e., 
    \begin{itemize}
        \item for all $v \in \mathcal{V}$: $a^-_v:=\mathfrak{v}_v(f_v)<\min\{ \mathfrak{v}_v(f_w)\,:\,v \neq
        w \in\mathcal{V}\}$;
        \item for all $v \in \mathcal{V}$: $-a^+_v:=\mathfrak{v}^M_v(m_v)<\min\{ \mathfrak{v}^M_v(m_w)\,:\,v \neq
        w \in\mathcal{V}\}$.
    \end{itemize}
    Now let us consider the lattice points $a^-=(a^-_1, \ldots, a^-_r)$ and $a^+=(a^+_1, \ldots, a^+_r)$ and the rectangle $R(a^-, a^+)$ bounded by them. We claim the following:
    \begin{lemma}
        In the setting just described $a^+_v \geq a^-_v + 2$ for all $v \in \mathcal{V}$. Moreover, the rectangle $R(a^-, a^+)$ contains an interior lattice point with weight larger than $w\big|_{\partial R(a^-,a^+)}$ (i.e., the weights of all the cubes contained in $\partial R(a^-,a^+)$), inducing a nontrivial element of $\mathbb{H}_{r-1}(N \hookrightarrow M)$.
    \end{lemma}

    \begin{proof}
    First of all, notice that $a^-_v-a^+_v=\mathfrak{v}_v(f_vm_v)<0$ (cf. (\ref{eq:vv(fvmv)})), hence, $a^+_v >a^-_v$ for all $v \in \mathcal{V}$. Moreover, we will prove the following:
        \begin{enumerate}        
            \item if a lattice point $\ell \in R(a^-, a^+)$ has $\ell_v=a^-_v$ for some index $v \in \mathcal{V}$, then $w(\ell+e_v)>w(\ell)$;
            \item if a lattice point $\ell \in R(a^-, a^+)$ has $\ell_v=a^+_v$ for some index $v \in \mathcal{V}$, then $w(\ell-e_v)>w(\ell)$.
        \end{enumerate}
        These two properties
        imply that $a^+_v$ cannot be $a^-_v+1$, hence, $a^+_v\geq a^-_v+2$  for all $v \in \mathcal{V}$. On the other hand, they also imply that $\max \big(w\big|_{\partial R(a^-,a^+)}\big)=\max\{w(\ell)\,:\,\ell \in \partial R(a^-, a^+)\cap \mathbb{Z}^r\}$ is obtained in some lattice point $\ell^+$ on the interior of an $(r-1)$-dimensional facet (i.e., for some $v \in \mathcal{V}:\ \ell^+_v=a^\pm_v$ and for all $v \neq w \in \mathcal{V}:\, a^-_w <\ell^+_w<a^+_w$) and its single neighbour $\ell^+ \mp e_v$ in the interior of the rectangle $R(a^-, a^+)$ satisfies $w(\ell^+ \mp e_v)>w(\ell^+)= \max\big(w\big|_{\partial R(a^-,a^+)}\big)$. Thus $\partial R(a^-, a^+)$ is contained, but is not contractible in $S_{w(\ell^+)}$, hence, $\mathbb{H}_{{\rm red},r-1, -2w(\ell^+)}(N\hookrightarrow M) \neq 0$. Therefore, there only remains the verification of properties $(1)$ and $(2)$.

        Let us now consider a lattice point  $\ell \in R(a^-, a^+)$ with $\ell_v=a^-_v$ for some index $v \in \mathcal{V}$ (i.e., $\ell_v=a^-_v$ and $a^-_w \leq  \ell_w \leq a^+_w$ for all $v \neq w \in \mathcal{V}$). By the Combinatorial Duality Property (cf. Definition \ref{def:COMPGOR}), it is enough to prove that $\mathfrak{h}_{\mathcal{D}}(\ell+e_v)>\mathfrak{h}_{\mathcal{D}}(\ell)$. For this we check that $f_v\in \mathcal{F}_{\mathcal{D}}(\ell) \setminus \mathcal{F}_{\mathcal{D}}(\ell+e_v)$ (for the notations see section \ref{s:4}). First, notice that $\mathfrak{v}_v(f_v)=a^-_v=\ell_v$, hence $f_v \notin \mathcal{F}_{\mathcal{D}}(\ell+e_v)$. Secondly, we also know that for any $v \neq w \in \mathcal{V}:\ f_vm_w\in N$, hence $0 \leq \mathfrak{v}_w^M(f_vm_w)=\mathfrak{v}_w(f) -a^+_w$. This finally implies that 
        \begin{center}
            $\mathfrak{v}_w(f_v)\geq a^+_w\geq \ell_w$ for all $w \in \mathcal{V}$, hence $f_v \in \mathcal{F}_{\mathcal{D}}(\ell)$, and thus $w(\ell+e_v)>w(\ell)$.
        \end{center}

        The proof of property $(2)$ is completely analogous: by the Combinatorial Duality Property it is enough to check that $m_v \in \mathcal{F}_\mathcal{D}^M(-\ell) \setminus \mathcal{F}_\mathcal{D}^M(-\ell+e_v)$, which is the direct consequence of our assumption that $f_wm_v\in N$ for all $v \neq w \in \mathcal{V}$. 
    \end{proof}

    This also finishes the proof of Proposition \ref{prop:htopnem0}.
\end{proof}

\begin{example}[Continuation of Example \ref{ex:2552H1}]
    For the choice of $f_1, f_2, g_1, g_2 $ as above and for the valuations coming from the normalization of the components we get $a^-=(2,2)$ and $a^+=(4,4)$. On the weight table in Example \ref{ex:2552} we clearly see that $R\big((2,2), (4,4)\big)$ contains the local maximum point $(3,3)$ of the weight function $w_{an}$ with the $1$-cycle $\partial R\big((2,2), (4,4)\big)$ giving a nontrivial generator of $\mathbb{H}_{an, 1, 2}(C, 0)=\mathbb{H}_{an,1}(C, 0)=\mathbb{Z}^2$. Another generator (around the local maximum $(5,5)$) corresponds, e.g., to the combination
    $f_1:=x^2, \ f_2:=y^2$, $ g_1:=x, \ g_2:=y.$
\end{example}

In order to use Corollary \ref{cor:htopnem0sym} (or Proposition \ref{prop:htopnem0}) to its full strength, we have to understand what the minimal cardinality $r$ of a realizing collection of (extended) discrete valuations is and then find ring (and module) elements satisfying the required conditions. In the case of realizable monomial ideals in $k[x_1, x_2]$ we can fulfill this plan using the tool of Newton polytopes and their geometric realizations. 

\subsection{Geometric realizations of monomial ideals in $k[x_1,x_2]$}\label{ss:geomrelmonid}\,

Let us consider a finite codimensional integrally closed monomial ideal $\mathcal{M}=\overline{\mathcal{M}}$ in the two variable polynomial ring $\mathcal{O}_2=k[x_1, x_2]$. By our previous results, its symmetric lattice homology $\mathbb{SH}_*(\mathcal{M} \triangleleft \mathcal{O}_2)$ is well-defined and can be computed using any CDP realization of $\mathcal{M}$. [In fact, since the lattice homology only depends on the quotient Artin algebra $\mathcal{O}_2/\mathcal{M}$, this whole argument remains true in $k[[x_1, x_2]]$ or $k[x_1, x_2]_{(x_1,x_2)}$ --- or $\mathbb{C}\{x_1, x_2\}$ in the complex analytic case --- as well.] A particularly nice case to handle is when we choose a CDP realization consisting of only monomial valuations, since then all the computations become combinatorial in terms of the lattice points under the Newton polytope (see subsection \ref{ss:ND} for the definitions and notations). We present this process in this subsection (see also subsection \ref{ss:NewtonLatCoh} for the case of the monomial ideal ${\rm adj}(\mathcal{M})$).

\bekezdes \textbf{The realization induced by the Newton boundary.}\label{par:NGammafelbontas}
Since $\mathcal{M}$ is finite $k$-codimensional and monomial, its Newton boundary $N\Gamma(\mathcal{M})$ is convenient. Its edges $\{F_\sigma\}_{\sigma=1}^s$ lie on affine lines determined by positive primitive normal vectors $\{\mathbf{n}_{F_\sigma}\}_{\sigma=1}^s$ and positive integers $\{d_{F_\sigma}\}_{\sigma=1}^s$ for which $F_\sigma \subset \{p\in \mathbb{R}^2\,:\,\langle p, \mathbf{n}_{F_\sigma}\rangle =d_{F_\sigma}\}$ for all $1 \leq \sigma \leq s$. Thus 
\begin{equation}\label{eq:NG+kivagva}
N\Gamma_+(\mathcal{M})=\{p \in (\mathbb{R}_{\geq 0})^2\,:\,\langle p, \mathbf{n}_{F_\sigma} \rangle \geq d_{F_\sigma}  \text{ for all } 1 \leq \sigma \leq s\}
\end{equation}
and the corresponding collection $\mathcal{D}_{N\Gamma}=\{\mathfrak{v}_{\mathbf{n}_{F_1}}, \ldots, \mathfrak{v}_{\mathbf{n}_{F_s}}\}$ of monomial valuations (i.e., we assign to the edge $F_\sigma$ the valuation $\mathfrak{v}_{\mathbf{n}_{F_\sigma}}\big(\sum_{p\in (\mathbb{Z}_{\geq 0})^2} a_p x^p\big)=\min\{ \langle p, \mathbf{n}_{F_\sigma} \rangle\,:\, a_p \neq 0\}$) and tuple $d_{N\Gamma}=(d_{F_1}, \ldots, d_{F_s})$ give a realization of $\mathcal{M}$ in $\mathcal{O}_2$ in the sense of Definition \ref{def:REALforideals}.

On the other hand, if we only want to obtain a realization of $\mathcal{M}$ it is enough to require from the monomial valuations and corresponding tuple the weaker (than 
(\ref{eq:NG+kivagva})) condition  to separate only the lattice points of $N\Gamma_+(\mathcal{M})$. 

\begin{nota}\label{nota:RB}
    Before that, let us introduce the following notations for the sets of lattice points
    \begin{itemize}
        \item  lying in the support of the integrally closed monomial ideal $\mathcal{M}$, $\mathcal{R}:=N\Gamma_+(\mathcal{M}) \cap (\mathbb{Z}_{\geq 0})^2$;
        \item  lying under the Newton boundary, $\mathcal{B}:=(\mathbb{Z}_{\geq 0})^2 \setminus \mathcal{R}$.
    \end{itemize}
\end{nota}

\begin{obs}\label{obs:geomreal}
    Monomial realizations of a finite codimensional integrally closed monomial ideal $\mathcal{M} \triangleleft \mathcal{O}_2$ (i.e., realizations $\big(\mathcal{D}=\{\mathfrak{v}_1, \ldots, \mathfrak{v}_r\},\, d_{\mathcal{D}}=(d_1, \ldots, d_r)\big)$ using monomial valuations) corres\-pond to certain (multi)sets of lines $\{L_v\}_{v=1}^r$ given by equations $L_v=\{(p_1, p_2)\in \mathbb{R}^2\,:\,p_2=m^v p_1+b^v\}$ having negative rational slope and positive rational $y$-intercept (i.e., $m^v=-\mu_v/\nu_v$ and $b^v=\beta_v/\gamma_v$ with $\mu_v, \nu_v, \beta_v, \gamma_v \in \mathbb{Z}_{>0}$ and $\gcd(\mu_v, \nu_v)=\gcd(\beta_v, \gamma_v)=1$) for all $1 \leq v \leq r$. Indeed, we can assign to a pair $(\mathfrak{v}_v, d_v)$ the line 
    \begin{equation*}
        L_{(\mathfrak{v}_v, d_v)}=\{ \langle (p_1, p_2), (\mathfrak{v}_v(x_1), \mathfrak{v}_v(x_2))\rangle =d_v\} =\Big\lbrace p_2=-\dfrac{\mathfrak{v}_v(x_1)}{\mathfrak{v}_v(x_2)}p_1 + \dfrac{d_v}{\mathfrak{v}_v(x_2)}\Big\rbrace,
    \end{equation*}
    and, vice versa, to a line $L_v=\big\{p_2=-\frac{\mu_v}{\nu_v} p_1+\frac{\beta_v}{\gamma_v}\big\}$ the pair $(\mathfrak{v}_{(\mu_v \gamma_v, \nu_v\gamma_v)}, \beta_v \nu_v)$, where $\mathfrak{v}_{(\mu_v \gamma_v, \nu_v\gamma_v)}$ is the monomial valuation $\sum_{p\in (\mathbb{Z}_{\geq 0})^2} a_p x^p\mapsto \min\{ \langle p,  (\mu_v \gamma_v, \nu_v\gamma_v) \rangle\,:\, a_p \neq 0\}$.

    Moreover, such a collection $\mathcal{D}$ of monomial valuations and tuple $d_{\mathcal{D}}$ of positive integers give a realization of $\mathcal{M}$ if and only if for the corresponding (multi)set of lines $\big\{L_v=\{p_2=m^vp_1+b^v\}\big\}_{v=1}^r$ we have 
    \begin{equation}\label{eq:geomreal}
        \{ (p_1, p_2)\in (\mathbb{Z}_{\geq 0})^2\,:\, p_2-m^vp_1 \geq b^v \text{ for all }1 \leq v \leq r\}=\mathcal{R},
    \end{equation}
    or, equivalently, 
    \begin{equation*}
        \{ (p_1, p_2)\in (\mathbb{Z}_{\geq 0})^2\,:\exists v\, (1 \leq v \leq r):\, p_2-m^vp_1 < b^v\}=\mathcal{B}.
    \end{equation*}
\end{obs}

\begin{define}\label{def:geomreal}
    Given the support $\mathcal{R}$ of a finite codimensional integrally closed monomial ideal $\mathcal{M} \triangleleft \mathcal{O}_2$, we call a set $\mathfrak{L}=\big\{L_v=\{p_2=m^vp_1+b^v\}\}_{v=1}^r$ of lines with negative rational slope and positive rational $y$-intercept a \emph{`geometric realization'} of $\mathcal{R}$ if it satisfies condition (\ref{eq:geomreal}).

    To any such geometric realization $\mathfrak{L}=\big\{L_v=\{p_2=-\frac{\mu_v}{\nu_v}p_1+\frac{\beta_v}{\gamma_v}\}\big\}_{v=1}^r$ of $\mathcal{R}$ we associate the monomial CDP realization $\big(\mathcal{D}_{\mathfrak{L}}=\{\mathfrak{v}_{(2\mu_1\gamma_1, 2\nu_1\gamma_1)}, \ldots, \mathfrak{v}_{(2\mu_r\gamma_r, 2\nu_r\gamma_r)}\}, d_{\mathcal{D}_{\mathfrak{L}}}=(2\nu_1\beta_1, \ldots, 2\nu_r\beta_r)\big)$ of $\mathcal{M}$ in the sense of Definition \ref{def:REALforideals}. (Notice that we already inserted the doubling operation of Lemma \ref{lem:duplatrukkideal} into this definition to obtain the Combinatorial Duality Property for the associated height function.)
\end{define}

Analogously to paragraph \ref{par:NewtonComb} we can then compute the symmetric lattice homology $\mathbb{SH}_*(\mathcal{M} \triangleleft\mathcal{O}_2)$ combinatorially. More concretely, if we consider the degree function \begin{equation*}
    \deg_{\mathfrak{L}}:\mathcal{B} \rightarrow (\mathbb{Z}_{\geq 0})^r, (p_1, p_2) \mapsto \big(\langle (p_1, p_2) , (2\mu_1\gamma_1, 2\nu_1\gamma_1) \rangle, \ldots, \langle (p_1, p_2), (2\mu_r\gamma_r, 2\nu_r\gamma_r)\rangle\big),
\end{equation*}
then the height function on the rectangle $R(0, d_{\mathcal{D}_{\mathfrak{L}}})\cap \mathbb{Z}^r$ is $\mathfrak{h}_{\mathcal{D}_{\mathfrak{L}}}(\ell)=\big|\{p\in \mathcal{B}\,:\,\deg_{\mathfrak{L}}(p) \ngeq \ell\}\big|$ with the weight function being $w_{\mathcal{D}_{\mathfrak{L}}}(\ell) = \mathfrak{h}_{\mathcal{D}_{\mathfrak{L}}}(\ell) + \mathfrak{h}_{\mathcal{D}_{\mathfrak{L}}}(d_{\mathcal{D}_{\mathfrak{L}}}-\ell)-\big|\mathcal{B}\big|$. Since $\mathfrak{h}_{\mathcal{D}_{\mathfrak{L}}}$ satisfies the CDP, we have $\mathbb{SH}_*(\mathcal{M} \triangleleft\mathcal{O}_2)=\mathbb{H}_*(R(0, d_{\mathcal{D}_\mathfrak{L}}), w_{\mathcal{D}_\mathfrak{L}})$.

\subsection{Lattice homological dimension of monomial ideals in $k[x_1,x_2]$}\label{ss:homdimmonid}\,

Recall that in subsection \ref{ss:bounds} we defined the homological dimension of realizable submodules, which in the setting of a proper finite codimensional integrally closed monomial ideal $\mathcal{M} \triangleleft\mathcal{O}_2=k[x_1, x_2]$, with $M\neq \mathcal{O}_2$, reads as ${\rm homdim}(\mathcal{M} \triangleleft \mathcal{O}_2)=\max\{q\,:\,\mathbb{SH}_{q}(\mathcal{M} \triangleleft \mathcal{O}_2)\neq 0\}$. Moreover, by Proposition \ref{prop:homdegupperbound} and Observation \ref{obs:geomreal} we already have the inequality
\begin{align*}
    {\rm homdim}(\mathcal{M} \triangleleft \mathcal{O}_2) \leq \min \{|\mathcal{D}|\,:\, (\mathcal{D}, d_{\mathcal{D}}) & \text{ a realization of } \mathcal{M}\} -1\\ &\leq \min \{|\mathfrak{L}|\,:\, \mathfrak{L} \text{ a geometric realization of } \mathcal{R}\} -1.
\end{align*}
We claim that, in our current setting, this are, in fact, all equalities. In particular, in this case Conjecture \ref{conj:elso} also holds.

\begin{theorem}[= Theorem \ref{th:homdimmon}]\label{th:homdim}
    Let $\mathcal{M}$ be a finite codimensional integrally closed monomial ideal in $\mathcal{O}_2$ (which could mean $k[x_1, x_2], \ k[x_1, x_2]_{(x_1, x_2)}$ or $k[[x_1, x_2]]$ depending on the context). Using the notation $\mathcal{R}$ for its support in the exponent lattice $(\mathbb{Z}_{\geq 0})^2$, we have the following characterization of the (lattice) homological dimension of $\mathcal{M}$:
    \begin{align*}
    {\rm homdim}(\mathcal{M} \triangleleft \mathcal{O}_2) =&\, \min \{|\mathcal{D}|\,:\, (\mathcal{D}, d_{\mathcal{D}}) \text{ a realization of } \mathcal{M}\} -1\\ =&\, \min \{|\mathcal{D}|\,:\, (\mathcal{D}, d_{\mathcal{D}}) \text{ a monomial realization of } \mathcal{M}\} -1\\ =&\, \min \{|\mathfrak{L}|\,:\, \mathfrak{L} \text{ a geometric realization of } \mathcal{R}\} -1.
\end{align*}
\end{theorem}

Translated to the language of plane curve singularities with Newton nondegenerate principal part (via Theorem \ref{th:comparison}) this reads as:
 
 \begin{cor}\label{cor:htopforNNcurves}
 Let $(C,0)=(\{f=0\}, 0)\subset (\mathbb{C}^2, 0)$ be a plane curve singularity with Newton nondegenerate principal part and convenient Newton boundary $N\Gamma(f)$.
Then we have the following:
   \begin{align*}
      {\rm homdim}_{an}&(C, 0):=\begin{cases}
      -1 & \text{ if }(C, 0) \text{ is regular}\,;\\
          \max\{s\,:\,\mathbb{H}_{an, s}(C, 0)\neq 0\} & \text{ if }(C, 0) \text{ is singular}.
      \end{cases}\\
      &=\min\{|\mathcal{D}|\,:\, (\mathcal{D}, d) \text{ is a realization of }q^{-1}(\mathcal{C})\}-1\\
      &=\min\{|\mathcal{D}|\,:\, (\mathcal{D}, d) \text{ is a monomial realization of }q^{-1}(\mathcal{C})\}-1\\
      &=\min\{|\mathfrak{L}|\,:\, \mathfrak{L} \text{ is a geometric realization of }\big((N\Gamma_+(f) \setminus N\Gamma(f))-(1,1)\big)\cap(\mathbb{Z}_{\geq 0})^2\}-1,
  \end{align*}
  where $q^{-1}(\mathcal{C})$ is the pullback of the conductor ideal $\mathcal{C}$ of $(C,o)$ (which is monomial --- see Theorem \ref{th:MerleTeissier}) along the natural projection map $q:\mathcal{O}_{\mathbb{C}^2,0} \rightarrow \mathcal{O}_{C, 0}$.
\end{cor}

The proof of Theorem \ref{th:homdim} consists of two parts: a coordinate geometric study of the minimal cardinality of a geometric realization of $\mathcal{R}$, and a careful choice of ring elements allowing the application of Corollary \ref{cor:htopnem0sym}.
These two parts are connected by the so-called `kerb configurations'.

\begin{define}\label{def:kerbconf}
    Consider the exponent lattice $(\mathbb{Z}_{\geq 0})^2 \subset \mathbb{R}^2$ of $\mathcal{O}_2$ and the support $\mathcal{R}$ of a finite codimensional integrally closed mononial ideal $\mathcal{M}$. $\mathcal{B}:=(\mathbb{Z}_{\geq 0})^2 \setminus \mathcal{R}$ (cf. Notation \ref{nota:RB}). A sequence of lattice points $\{p^l=(p_1^l, p_2^l)\}_{l=0}^{2t} \subset (\mathbb{Z}_{\geq 0})^2$ is called a \emph{`weak kerb configuration (associated with $\mathcal{R}$)'} if it satisfies the following conditions:
    \begin{itemize}
        \item for every $0 \leq l \leq 2t-1: p_1^{l+1}>p_1^l$ and $p_2^{l+1}<p_2^l$;
        \item for every \emph{even} index $0 \leq l\leq 2t$ the lattice point $p^l \in \mathcal{B}$ (`blue'), while for every \emph{odd} index $1 \leq l \leq 2t-1: p^l \in \mathcal{R}$ (`red');
        \item if we denote the slope of the line through the lattice points $p^l$ and $p^{l+1}$ by $m(\overline{p^l p^{l+1}})$, then we require the following inequalities:
        \begin{equation}\label{eq:kerbconfslopesor}
            m(\overline{p^0 p^{1}})\leq m(\overline{p^1 p^{2}})<m(\overline{p^2 p^{3}}) \leq m(\overline{p^3 p^{4}}) < m(\overline{p^4 p^{5}}) \leq \ldots \leq m(\overline{p^{2t-1} p^{2t}}).
        \end{equation}
    \end{itemize}

    Let us denote the $y$-intercept of the line $\overline{p^l p^{l+1}}$ by $b(\overline{p^l p^{l+1}})$ for every $0 \leq l \leq 2t-1$. We call the sequence $\{p^l=(p_1^l, p_2^l)\}_{l=0}^{2t} \subset (\mathbb{Z}_{\geq 0})^2$ a \emph{`strong kerb configuration (associated with $\mathcal{R}$)'} if it additionally satisfies the following property:
    \begin{align*}
        \{(p_1, p_2)\in \mathbb{Z}_{\geq p^0_1} \times \mathbb{Z}_{\geq p^{2t}_2}\,:\,p_2-m(\overline{p^l p^{l+1}})p_1\geq b(\overline{p^l p^{l+1}}) \text{ for all } 0\leq l \leq 2t-1\} \setminus \{p^{2l}\}_{l=0}^{t} \subset \mathcal{R},
    \end{align*}
    i.e., if any lattice point $(p_1, p_2)\in \mathbb{Z}^2$ with $p_1\geq p^0_1$ and $p_2\geq p_2^{2t}$ lying on or above the lines $\{\overline{p^lp^{l+1}}\}_{l=0}^{2t-1}$ must either be `red' or an element of the sequence.
    
    A concrete example of a strong kerb configuration can be seen in Figure \ref{fig:kerbconf}.
\end{define}

\begin{center}
\begin{figure}
\begin{picture}(220,200)(-20,-20)

\linethickness{0.35mm}
\put(-15,0){\line(1,0){210}}
\put(0,-15){\line(0,1){190}}

\linethickness{0.05mm}
  \multiput(-15,40)(0,20){5}{\line(1,0){210}}
  \multiput(40,-15)(20,0){4}{\line(0,1){190}}

  \multiput(120,-15)(20,0){4}{\line(0,1){145}}
  \multiput(120,150)(20,0){4}{\line(0,1){25}}
  \put(-15,140){\line(1,0){125}}
  \put(-15,160){\line(1,0){210}}
  \put(-15,20){\line(1,0){25}}
  \put(30,20){\line(1,0){165}}
  \put(20,-15){\line(0,1){25}}
  \put(20,30){\line(0,1){145}}

\put(0,160){\tikz \coordinate (p1);}
\put(20,100){\tikz \coordinate (p2);}
\put(80,40){\tikz \coordinate (p3);}
\put(120,20){\tikz \coordinate (p4);}
\put(180,0){\tikz \coordinate (p5);}
\put(195,0){\tikz \coordinate (p6);}
\put(195,175){\tikz \coordinate (p7);}
\put(0,175){\tikz \coordinate (p8);}
\put(0,0){\tikz[overlay] \fill[red, nearly transparent] (p1) -- (p2) -- (p3) -- (p4) -- (p5) -- (p6) -- (p7) -- (p8) -- cycle;}

\multiput(0,160)(20,0){10}{\color{red}\circle*{4}}
\multiput(20,140)(20,0){5}{\color{red}\circle*{4}}
\multiput(20,120)(20,0){9}{\color{red}\circle*{4}}
\multiput(20,100)(20,0){9}{\color{red}\circle*{4}}
\multiput(40,80)(20,0){8}{\color{red}\circle*{4}}
\multiput(60,60)(20,0){7}{\color{red}\circle*{4}}
\multiput(80,40)(20,0){6}{\color{red}\circle*{4}}
\multiput(120,20)(20,0){4}{\color{red}\circle*{4}}
\multiput(180,0)(20,0){1}{\color{red}\circle*{4}}

\multiput(0,140)(20,0){1}{\color{blue}\circle*{4}}
\multiput(0,120)(20,0){1}{\color{blue}\circle*{4}}
\multiput(0,100)(20,0){1}{\color{blue}\circle*{4}}
\multiput(0,80)(20,0){2}{\color{blue}\circle*{4}}
\multiput(0,60)(20,0){3}{\color{blue}\circle*{4}}
\multiput(0,40)(20,0){4}{\color{blue}\circle*{4}}
\multiput(40,20)(20,0){4}{\color{blue}\circle*{4}}
\multiput(0,0)(20,0){9}{\color{blue}\circle*{4}}
\put(0,20){\color{blue}\circle*{4}}

\linethickness{0.5mm}
\put(0,140){\line(1,-2){40}}
\put(40,60){\line(2,-1){120}}

\put(0,140){\color{blue}\circle*{6}}
\put(160,0){\color{blue}\circle*{6}}
\put(40,60){\color{blue}\circle*{6}}
\put(20,100){\color{red}\circle*{6}}
\put(80,40){\color{red}\circle*{6}}

\put(4,-10){\makebox(0,0){\small{$0$}}}
\put(24,-10){\makebox(0,0){\small{$1$}}}
\put(44,-10){\makebox(0,0){\small{$2$}}}
\put(64,-10){\makebox(0,0){\small{$3$}}}
\put(84,-10){\makebox(0,0){\small{$4$}}}
\put(104,-10){\makebox(0,0){\small{$5$}}}
\put(124,-10){\makebox(0,0){\small{$6$}}}
\put(144, -10){\makebox(0,0){\small{$7$}}}
\put(164, -10){\makebox(0,0){\small{$8$}}}
\put(184, -10){\makebox(0,0){\small{$9$}}}

\put(-10,5){\makebox(0,0){\small{$0$}}}
\put(-10,25){\makebox(0,0){\small{$1$}}}
\put(-10,45){\makebox(0,0){\small{$2$}}}
\put(-10,65){\makebox(0,0){\small{$3$}}}
\put(-10,85){\makebox(0,0){\small{$4$}}}
\put(-10,105){\makebox(0,0){\small{$5$}}}
\put(-10,125){\makebox(0,0){\small{$6$}}}
\put(-10,145){\makebox(0,0){\small{$7$}}}
\put(-10,165){\makebox(0,0){\small{$8$}}}

\put(110,130){\tikz \coordinate (p9);}
\put(190,130){\tikz \coordinate (p10);}
\put(190,150){\tikz \coordinate (p11);}
\put(110,150){\tikz \coordinate (p12);}

\put(150,140){\makebox(0,0){${\color{red} \mathcal{R}} \subset N\Gamma_+(\mathcal{M})$}}
\put(20,20){\makebox(0,0){\color{blue}$\mathcal{B}$}}
\end{picture}
\caption{In this example we consider the finite codimensional integrally closed monomial ideal $\mathcal{M}=(x_1^9, x_1^6x_2, x_1^4x_2^2, x_1^3 x_2^3, x_1^2 x_2^4, x_1x_2^5, x_2^8) \triangleleft k[x_1, x_2]$. The sequence $\{{\color{blue}(0,7)},\, {\color{red}(1, 5)},\, {\color{blue}(2,3)},\, {\color{red}(4, 2)},\, {\color{blue}(8,0)}\}$ of lattice points is a strong kerb configuration (similarly to the sequence $\{{\color{blue}(0,7)},\, {\color{red}(1, 5)},\, {\color{blue}(2,3)},\, {\color{red}(6, 1)},\, {\color{blue}(8,0)}\}$).} \label{fig:kerbconf}
\end{figure}
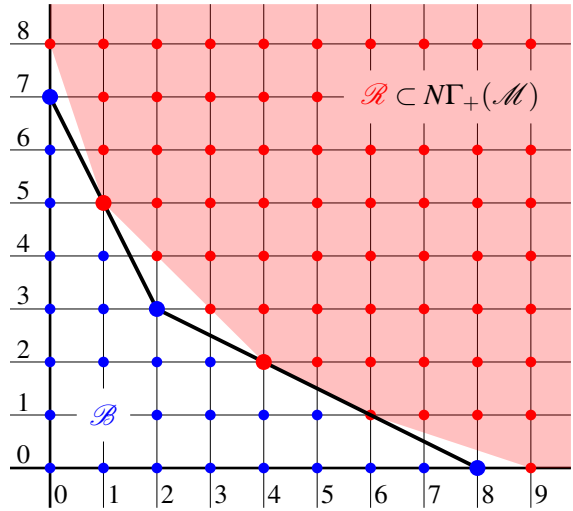
\end{center}
\vspace{-10mm}

 \begin{remark}
        Notice that, since $\mathcal{M}$ is integrally closed, $\mathcal{R}$ contains all the lattice points in its convex closure (cf. Proposition \ref{prop:closadj} \textit{(1)}), hence, already the inequality $m(\overline{p^{2l-1} p^{2l}})\leq m(\overline{p^{2l} p^{2l+1}})$ implies $m(\overline{p^{2l-1} p^{2l}})< m(\overline{p^{2l} p^{2l+1}})$ in (\ref{eq:kerbconfslopesor}). Moreover, by the first condition, all these slopes are negative rational numbers.
    \end{remark}

\begin{prop}\label{prop:linedim-1=kerbdim} 
    For a finite codimensional integrally closed monomial ideal $\mathcal{M} \triangleleft\mathcal{O}_2$ and its support $\mathcal{R}={\rm supp}(\mathcal{M}) \subset (\mathbb{Z}_{\geq 0})^2$ we have the following identity:
    \begin{align}
     \min \{|\mathfrak{L}|\,:\, \mathfrak{L} &\text{ a geometric realization of } \mathcal{R} \text{ in the sense of Definition \ref{def:geomreal}} \} -1 = \nonumber\\
     &\,=\max\big(\{t\in \mathbb{Z}_{\geq 0}\,:\,\exists \text{ a weak kerb configuration } \{p^l\}_{l=0}^{2t} \text{ associated with }\mathcal{R}\}\cup\{-1\}\big);
     \label{eq:linedim-1=kerbdim}\\
     &\,=\max\big(\{t\in \mathbb{Z}_{\geq 0}\,:\,\exists \text{ a strong kerb configuration } \{p^l\}_{l=0}^{2t} \text{ associated with }\mathcal{R}\}\cup \{-1\}\big). \nonumber
\end{align}
\end{prop}

\begin{nota}\label{nota:linekerbdim}
    For the sake of better readability, let us use the following notations:
    \begin{itemize}
        \item ${\rm linedim}(\mathcal{R}):=\min \{|\mathfrak{L}|\,:\, \mathfrak{L} \text{ a geometric realization of } \mathcal{R} \text{ in the sense of Definition \ref{def:geomreal}}\}$ ;
        \item ${\rm wkerbdim}(\mathcal{R}):=\max\big(\{t\in \mathbb{Z}_{\geq 0}\,:\,\exists\text{ a weak kerb configuration } \{p^l\}_{l=0}^{2t} \text{ ass. with }\mathcal{R}\}\cup\{-1\}\big)$;
        \item ${\rm skerbdim}(\mathcal{R}):=\max\big(\{t\in \mathbb{Z}_{\geq 0}\,:\,\exists\text{ a strong kerb configuration } \{p^l\}_{l=0}^{2t} \text{ ass. with }\mathcal{R}\}\cup\{-1\}\big)$.
    \end{itemize}
    Therefore, Proposition \ref{prop:linedim-1=kerbdim} claims that ${\rm linedim}(\mathcal{R})-1={\rm wkerbdim}(\mathcal{R})={\rm skerbdim}(\mathcal{R})$.
\end{nota}
    
\begin{proof}[Proof of Proposition \ref{prop:linedim-1=kerbdim}]\,
Clearly, since every strong kerb configuration is also a weak kerb configuration, we automatically have that ${\rm wkerbdim}(\mathcal{R})\geq{\rm skerbdim}(\mathcal{R})$. This implies that it is enough to prove that  ${\rm linedim}(\mathcal{R})-1\geq{\rm wkerbdim}(\mathcal{R})$ and ${\rm skerbdim}(\mathcal{R})\geq {\rm linedim}(\mathcal{R})-1$. In the case when $\mathcal{M}=\mathcal{O}_2$ these are satisfied, since all these values are $(-1)$. Therefore, in the sequel we suppose  that $\mathcal{M}$ is a proper ideal with ${\rm linedim}(\mathcal{R})-1, \ {\rm wkerbdim}(\mathcal{R})$, and ${\rm skerbdim}(\mathcal{R})$ being non-negative.

First, let us show that ${\rm linedim}(\mathcal{R})-1\geq{\rm wkerbdim}(\mathcal{R})$.

\noindent This is equivalent to the fact, that if there exists a weak kerb configuration $\{p^l\}_{l=0}^{2t}$ consisting of $(2t+1)$ lattice points, then $\mathcal{R}$
    cannot be geometrically realized by $t$ or less lines (in the sense of Definition \ref{def:geomreal}), i.e., every geometric realization $\mathfrak{L}=\{L_v\}_{v=1}^r$ of $\mathcal{R}$ must contain at least $t+1$ lines. This, however, is a consequence of the following fact: in order for $\mathfrak{L}$ to satisfy the identity (\ref{eq:geomreal}), it must contain some lines $\{L_{v_i}\}_{i=0}^t$ with their slopes $\{m(L_{v_i})\}_{i=0}^t$ satisfying the following series of inequalities:
    \begin{equation*}
            m(L_{v_0})<m(\overline{p^0 p^{1}})\leq \ldots  \leq m(\overline{p^{2i-1} p^{2i}})<m(L_{v_i})<m(\overline{p^{2i} p^{2i+1}}) \leq \ldots \leq m(\overline{p^{2t-1} p^{2t}})<m(L_{v_t}).
        \end{equation*}
    Indeed, without any $L_v \in \mathfrak{L}$ having slope between $m(\overline{p^{2i-1} p^{2i}})$ and $m(\overline{p^{2i} p^{2i+1}})$, the `blue' lattice point $p^{2i} \in \mathcal{B}$ could not be separated from the `red' lattice points $p^{2i-1}, p^{2i+1} \in \mathcal{R}$.  Thus, for any geometric realization $\mathfrak{L}$ of $\mathcal{R}$ and any kerb configuration $\{p^l\}_{l=0}^{2t}$ we have $\big|\mathfrak{L}\big| \geq t+1$.

    Secondly, we show that ${\rm linedim}(\mathcal{R})-1\leq{\rm skerbdim}(\mathcal{R})$.

    \noindent More precisely, we will prove that, since $\mathcal{R}$ cannot be geometrically realized by ${\rm linedim}(\mathcal{R})-1$ lines, there exists a strong kerb configuration $\{p^l\}_{l=0}^{2t}$ with $t={\rm linedim}(\mathcal{R})-1$. 
    In fact, we give a constructive argument.

 First let us consider the following family of `tangential' lines of slope $-\infty \leq m \leq0$:
    \begin{equation}\label{eq:tanline(m)}
        L(m):=\big\{p=(p_1, p_2)\in \mathbb{R}^2\,:\,p_2-mp_1=\langle p, (-m,1)\rangle= \max\{p'_2-mp'_1\,:\,(p'_1, p'_2) \in \mathcal{R}\}\big\}.
    \end{equation}
    For the sake of shorter formulas let us denote the right hand side of these linear equations by 
    \begin{equation*}
        b(L(m)):=\max\{\langle p', (-m,1)\rangle =p'_2-mp'_1\,:\,p'=(p'_1, p'_2) \in \mathcal{R}\}.
    \end{equation*}
    Now, for every negative slope $m$ let us consider the set $\mathcal{R}(m)=\mathcal{R} \cap L(m)$ of `red' lattice points lying on the `tangential' line $L(m)$ with slope $m$. 
    We will be looking for `blue' lattice points lying on or above the `tangential' line $L(m)$ and on the left side of $\mathcal{R}(m)$. We will construct inductively the sets $\mathcal{B}_i$ and $\mathcal{B}_i(m)$ (with $i\geq 0$ and $m\in [-\infty, 0]$ a slope). We start with $\mathcal{B}_0=\mathcal{B}$ and
    \begin{equation*}
        \mathcal{B}_0(m)=\big\{(p_1, p_2)\in \mathcal{B}_0\,:\, p_2-mp_1 \geq b(L(m)), \ p_1 < \min\{p_1'\,:\, (p_1', p_2') \in \mathcal{R}(m)\}\big\}.
    \end{equation*}
    (For an illustrating example see Figure \ref{fig:L(m)R(m)B(m)}.) 
    \begin{center}
\begin{figure}
\begin{picture}(300,220)(-60,-40)

\linethickness{0.35mm}
\put(-15,0){\line(1,0){210}}
\put(0,-15){\line(0,1){190}}

\linethickness{0.05mm}
  \multiput(-15,20)(0,20){2}{\line(1,0){210}}
  \multiput(-15,80)(0,20){5}{\line(1,0){210}}
  \multiput(20,-15)(20,0){5}{\line(0,1){190}}
  \multiput(160,-15)(20,0){2}{\line(0,1){190}}
  \multiput(120,-15)(20,0){2}{\line(0,1){65}}
  \multiput(120,70)(20,0){2}{\line(0,1){105}}
 
  \put(-15,60){\line(1,0){125}}
  \put(150,60){\line(1,0){45}}

\put(0,160){\tikz \coordinate (p1);}
\put(20,100){\tikz \coordinate (p2);}
\put(80,40){\tikz \coordinate (p3);}
\put(120,20){\tikz \coordinate (p4);}
\put(180,0){\tikz \coordinate (p5);}
\put(195,0){\tikz \coordinate (p6);}
\put(195,175){\tikz \coordinate (p7);}
\put(0,175){\tikz \coordinate (p8);}
\put(0,0){\tikz[overlay] \fill[red, nearly transparent] (p1) -- (p2) -- (p3) -- (p4) -- (p5) -- (p6) -- (p7) -- (p8) -- cycle;}

\linethickness{0.5mm}
\put(40,60){\line(-2,1){80}}
\put(40,60){\line(2,-1){180}}

\linethickness{0.2mm}
\put(0,140){\color{blue}\circle*{6}}
\put(0,140){\color{blue}\line(-2,-1){26}}
\put(0,120){\color{blue}\circle*{6}}
\put(0,120){\color{blue}\line(-1,0){24}}
\put(0,100){\color{blue}\circle*{6}}
\put(0,100){\color{blue}\line(-2,1){26}}
\put(0,80){\color{blue}\circle*{6}}
\put(0, 80){\color{blue}\line(-1,1){31}}
\put(20,80){\color{blue}\circle*{6}}
\put(20, 80){\color{blue}\line(-3,2){48}}
\put(40,60){\color{blue}\circle*{6}}
\put(40,60){\color{blue}\line(-4,3){69}}

\put(120,20){\color{red}\circle*{6}}
\put(120,20){\color{red}\line(1,4){8}}
\put(80,40){\color{red}\circle*{6}}
\put(80,40){\color{red}\line(5,2){36}}

\put(4,-10){\makebox(0,0){\small{$0$}}}
\put(24,-10){\makebox(0,0){\small{$1$}}}
\put(44,-10){\makebox(0,0){\small{$2$}}}
\put(64,-10){\makebox(0,0){\small{$3$}}}
\put(84,-10){\makebox(0,0){\small{$4$}}}
\put(104,-10){\makebox(0,0){\small{$5$}}}
\put(124,-10){\makebox(0,0){\small{$6$}}}
\put(144, -10){\makebox(0,0){\small{$7$}}}
\put(164, -10){\makebox(0,0){\small{$8$}}}

\put(-10,5){\makebox(0,0){\small{$0$}}}
\put(-10,25){\makebox(0,0){\small{$1$}}}
\put(-10,45){\makebox(0,0){\small{$2$}}}
\put(-10,65){\makebox(0,0){\small{$3$}}}

\put(110,130){\tikz \coordinate (p9);}
\put(190,130){\tikz \coordinate (p10);}
\put(190,150){\tikz \coordinate (p11);}
\put(110,150){\tikz \coordinate (p12);}

\put(130,60){\makebox(0,0){${\color{red} \mathcal{R}(m)} $}}
\put(-40,120){\makebox(0,0){\color{blue}$\mathcal{B}_0(m)$}}
\put(220,-15){\makebox(0,0){$L(m)$}}
\end{picture}
\caption{Continuing with the previous example ideal from Figure \ref{fig:kerbconf}, for the slope $m=-1/2$ we have $L(-1/2)=\{(x_1, x_2)\,:\,x_2=-\frac{1}{2}x_1 + 4\}, \ {\color{red} \mathcal{R}(-1/2)}=\{{\color{red}(4,2)}, {\color{red}(6,1)}\}$ and ${\color{blue}\mathcal{B}_0(-1/2)}=\{{\color{blue}(0, 7)}, {\color{blue}(0, 6)}, {\color{blue}(0, 5)}, {\color{blue}(0,4)}, {\color{blue}(1, 4)}, {\color{blue}(2, 3)}\}$. } \label{fig:L(m)R(m)B(m)}
\end{figure}
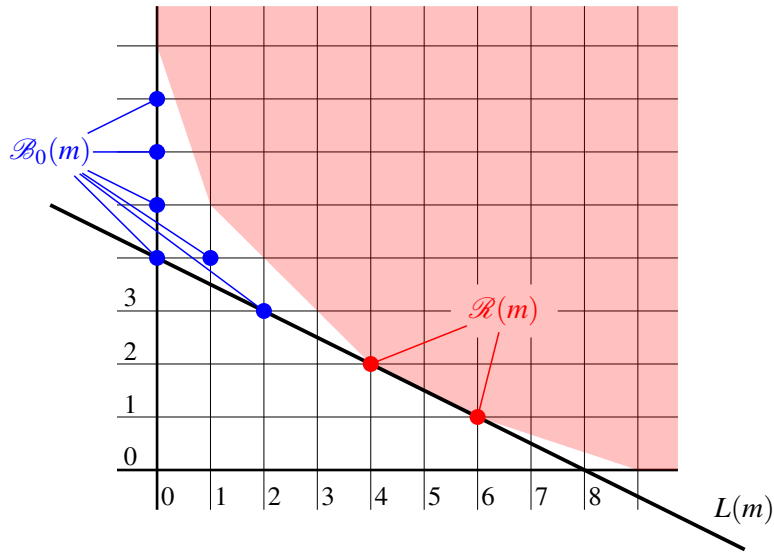
\end{center}

\vspace{-10mm}
    
    Now clearly
    \begin{itemize}
        \item $\mathcal{B}_0(-\infty)=\emptyset$, since $\min\{p'_1\,:\,(p'_1, p'_2)\in \mathcal{R}(-\infty)\}=0$, and
        \item $\mathcal{B}_0(0)=\mathcal{B}\neq \emptyset$, since $\mathcal{M}$ is proper.
    \end{itemize}
    Thus, by continuity of the function $m \mapsto L(m)$ and the finiteness of the set $\mathcal{B}_0=\mathcal{B}$, there exists the minimal slope $m^0 \in \mathbb{Q}_{<0}, \ -\infty < m^0 \leq 0$ such that $\mathcal{B}_0(m^0) \neq \emptyset$,
    \begin{center}
        i.e., $m^0=\min \{ m > -\infty\,:\,\mathcal{B}_0(m) \neq \emptyset\}$.
    \end{center}
    In fact, by the minimality of $m^0$, we know that $\emptyset \neq \mathcal{B}_0(m^0)\subset L(m^0)$. Therefore, there exist lattice points $p^0, \ p^1 \in L(m^0)$ ($p^1$ not necessarily unique) such that:
    \begin{itemize}
        \item $p^0\in \mathcal{B}_0(m^0)$ with $p^0_1=\max\{p_1\,:\,\text{there exists } (p_1, p_2)\in \mathcal{B}_0(m^0)\}$ and $p^1 \in \mathcal{R}(m^0)$;
        \item $p^0_1<p^1_1$ and $p^0_2>p^1_2$;
        \item $m(\overline{p^0 p^1})=m^0\in \mathbb{Q}_{<0}$ and
        \item $\{(p_1, p_2)\in \mathbb{Z}_{\geq p^0_1}\times \mathbb{Z}_{\geq p^1_2}\,:\,p_2-m^0p_1\geq b(L(m^0))\} \setminus \{p^0\} \subset \mathcal{R}$.
    \end{itemize}

    Now, we set $\mathcal{B}_1:=\{(p_1, p_2)\in \mathcal{B}_0 \setminus \mathcal{B}_0(m^0)\,:\,p_2-m^0 p_1 \geq b(L(m^0)) \big\}$. (For a concrete example see Figure \ref{fig:m0B1}.) In fact, this is the same as $\{(p_1, p_2)\in \mathcal{B}_0\,:\,p_2-(m^0-\varepsilon) p_1 \geq b(L(m^0-\varepsilon)) \big\}$ for some $\varepsilon \in \mathbb{Q}_{>0}$ small enough. Thus, $\mathcal{B}_1=\emptyset$ if and only if ${\rm linedim}(\mathcal{R})=1$. Indeed, if $\mathcal{B}_1$ is empty, then $L(m^0-\varepsilon)$ gives a geometric realization of $\mathcal{R}$, whereas if there exists some $q^2\in \mathcal{B}_1$, then $\{{\color{blue}p^0},\, {\color{red}p^1},\, {\color{blue}q^2}\}$ gives a weak kerb configuration implying (by the first part of this Proposition) that ${\rm linedim}(\mathcal{R})\geq 2$.

\begin{center}
\begin{figure}
\begin{picture}(220,200)(-20,-20)

\linethickness{0.35mm}
\put(-15,0){\line(1,0){210}}
\put(0,-15){\line(0,1){190}}

\linethickness{0.05mm}
  \put(-15, 20){\line(1,0){85}}
  \put(90,20){\line(1,0){105}}
  \multiput(-15,40)(0,20){3}{\line(1,0){210}}
  \put(-15,100){\line(1,0){45}}
  \put(50,100){\line(1,0){145}}
  \put(-15,120){\line(1,0){210}}
  \put(-10, 140){\line(1,0){205}}
  \put(-15,160){\line(1,0){210}}
  \put(20,-15){\line(0,1){190}}
  \put(40,-15){\line(0,1){105}}
  \put(40,110){\line(0,1){65}}
  \put(60,-15){\line(0,1){190}}
  \put(80,-15){\line(0,1){25}}
  \put(80,30){\line(0,1){145}}
  \multiput(100,-15)(20,0){5}{\line(0,1){190}}

\put(0,160){\tikz \coordinate (p1);}
\put(20,100){\tikz \coordinate (p2);}
\put(80,40){\tikz \coordinate (p3);}
\put(120,20){\tikz \coordinate (p4);}
\put(180,0){\tikz \coordinate (p5);}
\put(195,0){\tikz \coordinate (p6);}
\put(195,175){\tikz \coordinate (p7);}
\put(0,175){\tikz \coordinate (p8);}
\put(0,0){\tikz[overlay] \fill[red, nearly transparent] (p1) -- (p2) -- (p3) -- (p4) -- (p5) -- (p6) -- (p7) -- (p8) -- cycle;}

\linethickness{0.5mm}
\put(40,60){\line(-1,2){55}}
\put(40,60){\line(1,-2){38}}

\put(0,140){\color{blue}\circle*{6}}
\put(20,100){\color{red}\circle*{6}}

\put(40,60){\color{blue}\circle*{4}}
\put(60,40){\color{blue}\circle*{4}}
\put(60,20){\color{blue}\circle*{4}}
\put(100,20){\color{blue}\circle*{4}}
\multiput(80,0)(20,0){5}{\color{blue}\circle*{4}}

\put(4,-10){\makebox(0,0){\small{$0$}}}
\put(24,-10){\makebox(0,0){\small{$1$}}}
\put(44,-10){\makebox(0,0){\small{$2$}}}
\put(64,-10){\makebox(0,0){\small{$3$}}}
\put(84,-10){\makebox(0,0){\small{$4$}}}
\put(104,-10){\makebox(0,0){\small{$5$}}}
\put(124,-10){\makebox(0,0){\small{$6$}}}
\put(144, -10){\makebox(0,0){\small{$7$}}}
\put(164, -10){\makebox(0,0){\small{$8$}}}
\put(184, -10){\makebox(0,0){\small{$9$}}}

\put(-10,5){\makebox(0,0){\small{$0$}}}
\put(-10,25){\makebox(0,0){\small{$1$}}}
\put(-10,45){\makebox(0,0){\small{$2$}}}
\put(-10,65){\makebox(0,0){\small{$3$}}}
\put(-10,85){\makebox(0,0){\small{$4$}}}
\put(-10,105){\makebox(0,0){\small{$5$}}}

\put(110,130){\tikz \coordinate (p9);}
\put(190,130){\tikz \coordinate (p10);}
\put(190,150){\tikz \coordinate (p11);}
\put(110,150){\tikz \coordinate (p12);}

\put(-20,140){\makebox(0,0){$p^0$}}
\put(40,100){\makebox(0,0){$p^1$}}
\put(80,20){\makebox(0,0){\color{blue}$\mathcal{B}_1$}}
\end{picture}
\caption{In our running example we have $m^0=-2$, $L(-2)=\{x_2=-2x_1+7\}$,  $\mathcal{B}_0(-2)=\{p^0\}=\{{\color{blue}(0,7)}\}$ and $\mathcal{R}(-2)=\{p^1\}=\{{\color{red}(1,5)}\}$. Using these we obtain that $\mathcal{B}_1=\{{\color{blue}(2, 3)}, {\color{blue}(3, 2)}, {\color{blue}(3, 1)}, {\color{blue}(4, 1)}, {\color{blue}(5, 1)}, {\color{blue}(4, 0)}, {\color{blue}(5, 0)}, {\color{blue}(6, 0)}, {\color{blue}(7, 0)}, {\color{blue}(8, 0)}\}$.} \label{fig:m0B1}
\end{figure}
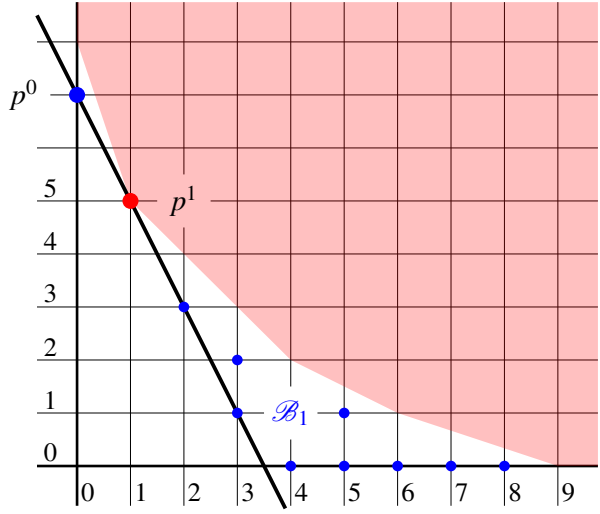
\end{center}

\vspace{-10mm}

If ${\rm linedim}(\mathcal{R})=1$, then the `blue' lattice point $p^0$ gives a one-point strong kerb configuration, otherwise ${\rm linedim}(\mathcal{R})\geq 2$ and we proceed similarly as before, using $\mathcal{B}_1$ in place of $\mathcal{B}_0$:  for any slope $m^0 \leq m \leq 0$ we consider the `blue' lattice point sets
    \begin{equation*}
        \mathcal{B}_1(m)=\big\{(p_1, p_2)\in \mathcal{B}_1\,:\, p_2-mp_1 \geq b(L(m)), p_1 < \min \{p'_1\,:\, (p'_1, p'_2)\in \mathcal{R}(m)\} \big\}.
    \end{equation*}
    For $m=m^0$ we have $\mathcal{B}_1(m^0)=\emptyset$, whereas for $m=0$ we have $\mathcal{B}_1(0)=\mathcal{B}_1\neq \emptyset$. 

    Now, similarly as before, there exists a minimal slope $m^1=\min\{m>m^0\,:\, \mathcal{B}_1(m) \neq \emptyset\}$, satisfying $\mathcal{B}_1(m^1)\subset L(m^1)$ and providing lattice points $p^2, \ p^3 \in L(m^1)$ ($p^3$ not necessarily unique) such that:
    \begin{itemize}
        \item $p^2\in \mathcal{B}_1(m^1)$ with $p^2_1=\max\{p_1\,:\,\text{there exists } (p_1, p_2)\in \mathcal{B}_1(m^1)\}$ and $p^3 \in \mathcal{R}(m^1)$;
        \item $p^2_1>p^1_1$ and $p^2_2<p^1_2$;
        \item $p^3_1>p^2_1$ and $p^3_2<p^2_2$;
        \item $m(\overline{p^1p^2}), \ m(\overline{p^2, p^3}) \in \mathbb{Q}_{<0}$ with 
        \begin{equation*}
            m(\overline{p^0p^1}) =m^0 \leq m(\overline{p^1p^2}) < m(\overline{p^2p^3})=m^1,
            \end{equation*}  
         where the strict inequality is implied by the convexity of $N\Gamma_+(\mathcal{M})$; and
         \item $\{(p_1, p_2)\in \mathbb{Z}_{\geq p^0_1}\times \mathbb{Z}_{\geq p^3_2}\,:\,p_2-m(\overline{p^l p^{l+1}})p_1\geq b(\overline{p^l p^{l+1}}) \text{ for all } 0\leq l \leq 2\} \setminus \{p^0, p^2\} \subset \mathcal{R}$.
    \end{itemize}
    Once again, we can define $\mathcal{B}_2:=\{(p_1, p_2)\in \mathcal{B}_1 \setminus \mathcal{B}_1(L(m^1))\,:\,p_2-m^1p_1\geq b(m^1)\}$ and see that it is empty if and only if ${\rm linedim}(\mathcal{R})=2$.
    In this latter case the sequence $\{{\color{blue}p^0},\, {\color{red}p^1},\, {\color{blue}p^2}\}$ gives a strong kerb configuration, since 
    \begin{align*}
        \{(p_1, p_2)\in &\ \mathbb{Z}_{\geq p^0_1}\times \mathbb{Z}_{\geq p^2_2}\,:\,p_2-m(\overline{p^l p^{l+1}})p_1\geq b(\overline{p^l p^{l+1}}) \text{ for all } 0\leq l \leq 1\} \subset \\
        &\subset \{(p_1, p_2)\in \mathbb{Z}_{\geq p^0_1}\times \mathbb{Z}_{\geq p^3_2}\,:\,p_2-m(\overline{p^l p^{l+1}})p_1\geq b(\overline{p^l p^{l+1}}) \text{ for all } 0\leq l \leq 2\}.
    \end{align*}
    If, however, ${\rm linedim}(\mathcal{R})\geq3$ we can  repeat the arguments for $\mathcal{B}_2, \mathcal{B}_3,$ etc. (obtaining during this process the slopes $m^2, m^3,$ etc., and the 
    alternatingly coloured lattice points ${\color{blue} p^4}, {\color{red} p^5}, {\color{blue} p^6}, {\color{red} p^7},$ etc.), until we arrive to $\mathcal{B}_{{\rm linedim}(\mathcal{R})}= \emptyset$.  
    Then, similarly as before, the sequence $\{p^l\}_{l=0}^{2t}$ with $t={\rm linedim}(\mathcal{R})-1$ is a strong kerb configuration associated with $\mathcal{R}$. Thus ${\rm skerbdim}(\mathcal{R}) \geq {\rm linedim}(\mathcal{R})-1$. 

    \begin{center}
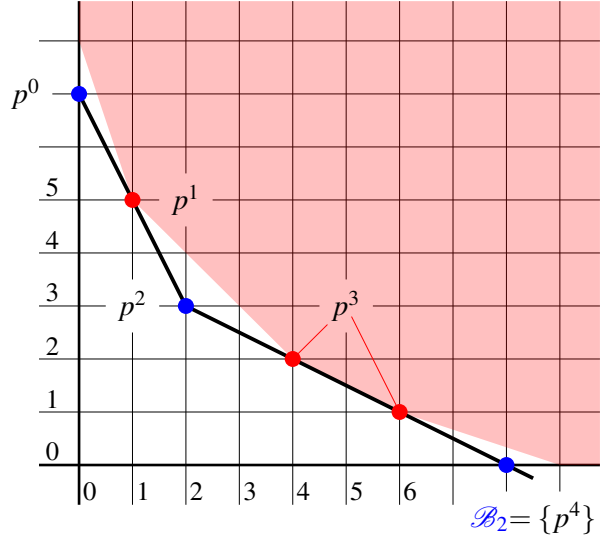
\begin{figure}
\begin{picture}(220,210)(-20,-30)

\linethickness{0.35mm}
\put(-15,0){\line(1,0){210}}
\put(0,-15){\line(0,1){190}}

\linethickness{0.05mm}
  \put(-15, 20){\line(1,0){210}}
  \put(-15, 40){\line(1,0){210}}
  \put(-15,60){\line(1,0){25}}
  \put(30,60){\line(1,0){60}}
  \put(110,60){\line(1,0){85}}
  \multiput(-15,80)(0,20){1}{\line(1,0){210}}
  \put(-15,100){\line(1,0){45}}
  \put(50,100){\line(1,0){145}}
  \put(-15,120){\line(1,0){210}}
  \put(-10, 140){\line(1,0){205}}
  \put(-15,160){\line(1,0){210}}
  \put(20,-15){\line(0,1){65}}
  \put(20,70){\line(0,1){105}}
  \put(40,-15){\line(0,1){105}}
  \put(40,110){\line(0,1){65}}
  \put(60,-15){\line(0,1){190}}
  \put(80,-15){\line(0,1){190}}
  \put(100,-15){\line(0,1){65}}
  \put(100,70){\line(0,1){105}}
  \multiput(120,-15)(20,0){2}{\line(0,1){190}}
  \multiput(160,-10)(20,0){2}{\line(0,1){185}}

\put(0,160){\tikz \coordinate (p1);}
\put(20,100){\tikz \coordinate (p2);}
\put(80,40){\tikz \coordinate (p3);}
\put(120,20){\tikz \coordinate (p4);}
\put(180,0){\tikz \coordinate (p5);}
\put(195,0){\tikz \coordinate (p6);}
\put(195,175){\tikz \coordinate (p7);}
\put(0,175){\tikz \coordinate (p8);}
\put(0,0){\tikz[overlay] \fill[red, nearly transparent] (p1) -- (p2) -- (p3) -- (p4) -- (p5) -- (p6) -- (p7) -- (p8) -- cycle;}

\put(80,40){\color{red}\line(1,1){13}}
\put(120,20){\color{red}\line(-1,2){18}}

\linethickness{0.5mm}
\put(40,60){\line(-1,2){40}}
\put(40,60){\line(2,-1){130}}

\put(0,140){\color{blue}\circle*{6}}
\put(20,100){\color{red}\circle*{6}}
\put(40,60){\color{blue}\circle*{6}}

\put(80,40){\color{red}\circle*{6}}
\put(120,20){\color{red}\circle*{6}}
\put(160,0){\color{blue}\circle*{6}}

\put(4,-10){\makebox(0,0){\small{$0$}}}
\put(24,-10){\makebox(0,0){\small{$1$}}}
\put(44,-10){\makebox(0,0){\small{$2$}}}
\put(64,-10){\makebox(0,0){\small{$3$}}}
\put(84,-10){\makebox(0,0){\small{$4$}}}
\put(104,-10){\makebox(0,0){\small{$5$}}}
\put(124,-10){\makebox(0,0){\small{$6$}}}

\put(-10,5){\makebox(0,0){\small{$0$}}}
\put(-10,25){\makebox(0,0){\small{$1$}}}
\put(-10,45){\makebox(0,0){\small{$2$}}}
\put(-10,65){\makebox(0,0){\small{$3$}}}
\put(-10,85){\makebox(0,0){\small{$4$}}}
\put(-10,105){\makebox(0,0){\small{$5$}}}

\put(110,130){\tikz \coordinate (p9);}
\put(190,130){\tikz \coordinate (p10);}
\put(190,150){\tikz \coordinate (p11);}
\put(110,150){\tikz \coordinate (p12);}

\put(-20,140){\makebox(0,0){$p^0$}}
\put(40,100){\makebox(0,0){$p^1$}}
\put(20,60){\makebox(0,0){$p^2$}}
\put(100,60){\makebox(0,0){$p^3$}}
\put(170,-20){\makebox(0,0){{\color{blue} $\mathcal{B}_2$}$=\{p^4\}$}}
\end{picture}
\caption{In our particular case $m^1=-1/2$ inducing that $p^2={\color{blue} (2,3)}$ and $p^3={\color{red} (4,2)}$ or ${\color{red} (6,1)}$. Then $\mathcal{B}_2=\{{\color{blue}(0,8)}\}$, so $m^2=0$ and $\mathcal{B}_3=\emptyset$, implying that $p^4={\color{blue}(0,8)}$ comp\-letes the strong kerb configuration ($p^5$ could be any {\color{red}$(p_1,0)$}, with $p_1\geq 9$).} \label{fig:m1B2}
\end{figure}
\end{center}
\vspace{-10mm}
\end{proof}

\begin{remark}
    The inductive nature of the proof suggests that Proposition \ref{prop:linedim-1=kerbdim} would also hold true in the more general setting of an affine lattice $\Lambda \simeq \mathbb{Z}^2, \ \Lambda \subset \mathbb{R}^2$ and an unbounded convex polytope defined by finitely many linear inequalities (in our specific case these are the standard lattice $\mathbb{Z}^2 \leq \mathbb{R}^2$ and the Newton polytope $N\Gamma_+(\mathcal{M})$). However, for the sake of simplicity we omit the description of this more general case here.
\end{remark}

We can now turn to the proof of Theorem \ref{th:homdim}.

 \begin{proof}[Proof of Theorem \ref{th:homdim}]
 By Propositions \ref{prop:homdegupperbound} and \ref{prop:linedim-1=kerbdim} it is enough to prove that for a finite codimensional integrally closed monomial ideal $\mathcal{M} \triangleleft \mathcal{O}_2$ we have
 \begin{equation*}
     {\rm homdim}(\mathcal{M}\triangleleft\mathcal{O}_2)\geq {\rm skerbdim}(\mathcal{R}),
 \end{equation*}
 where $\mathcal{R}={\rm supp}(\mathcal{M}) \subset (\mathbb{Z}_{\geq 0})^2$ (cf. Notations \ref{nota:RB} and \ref{nota:linekerbdim}). In practice, we will construct from a strong kerb configuration $\{p^l\}_{l=0}^{2t}$, with $t={\rm skerbdim}(\mathcal{R})={\rm linedim}(\mathcal{R})-1$, two sets $\{f_0, \ldots, f_t\}$ and $\{g_0, \ldots, g_t\}$ of ring elements satisfying the conditions of Corollary \ref{cor:htopnem0sym}. Since $\mathcal{R}$ can be geometrically realized by ${\rm linedim}(\mathcal{R})=t+1$ lines (i.e., $\mathcal{M}$ can be realized by $t+1$ monomial valuations --- cf. Observation \ref{obs:geomreal}) we can apply the second part of this Corollary \ref{cor:htopnem0sym} to ensure that 
 ${\rm homdim}(\mathcal{M} \triangleleft \mathcal{O}_2)\geq t$. 
 
 In fact, all the above ring elements will be monomials of form $f_i=x^{F_i}$ ($0\leq i \leq t$) and $g_j=x^{G_j}$ ($0\leq j \leq t$) for `blue' lattice points $F_0, \ldots, F_t, G_0, \ldots, G_t \in \mathcal{B}=(\mathbb{Z}_{\geq 0})^2\setminus \mathcal{R}$, for 
 which we require 
 the following properties:
 \begin{align}
     & \bullet \ F_i + G_i \in \mathcal{B} && \text{ for all } i\in \{0, \ldots, t\}; \label{eq:F+GinB}\\
     & \bullet \ F_i+G_j \in \mathcal{R} && \text{ for all } i, j \in \{0, \ldots, t\},\, i \neq j. \label{eq:F+GinR}
 \end{align}
Indeed, in the case of our integrally closed monomial ideal $\mathcal{M}\triangleleft \mathcal{O}_2 $, these conditions on the lattice point level are equivalent to those of Corollary \ref{cor:htopnem0sym} on the ring element\,/\,monomial level.

Now, given the strong kerb configuration $\{p^l\}_{l=0}^{2t}$, with $t={\rm skerbdim}(\mathcal{R})$, we define these lattice points as follows:
\begin{align}
    F_i:=&\ \left(\sum_{l=0}^{2i-1}(-1)^{l-1}p^l_1\,, \ \sum_{l=2i}^{2t}(-1)^lp^l_2\right) &\text{ for all } 0 \leq i \leq t;\\
    G_j:=&\ \left(\sum_{l=0}^{2j}(-1)^{l}p^l_1\,, \ \sum_{l=2j+1}^{2t}(-1)^{l+1}p^l_2\right) &\text{ for all } 0 \leq j \leq t.
\end{align}
(For a concrete example see Figure \ref{fig:summands}.)

\begin{center}
\begin{figure}
\begin{picture}(220,200)(-20,-20)

\linethickness{0.35mm}
\put(-15,0){\line(1,0){210}}
\put(0,-15){\line(0,1){190}}

\linethickness{0.05mm}
  \multiput(-15,20)(0,20){8}{\line(1,0){210}}
  \multiput(20,-15)(20,0){9}{\line(0,1){190}}

\put(0,160){\tikz \coordinate (p1);}
\put(20,100){\tikz \coordinate (p2);}
\put(80,40){\tikz \coordinate (p3);}
\put(120,20){\tikz \coordinate (p4);}
\put(180,0){\tikz \coordinate (p5);}
\put(195,0){\tikz \coordinate (p6);}
\put(195,175){\tikz \coordinate (p7);}
\put(0,175){\tikz \coordinate (p8);}
\put(0,0){\tikz[overlay] \fill[red, nearly transparent] (p1) -- (p2) -- (p3) -- (p4) -- (p5) -- (p6) -- (p7) -- (p8) -- cycle;}

\linethickness{0.4mm}
\put(0,140){\line(1,-2){40}}
\put(40,60){\line(2,-1){120}}

\put(0,140){\color{blue}\circle*{5}}
\put(160,0){\color{blue}\circle*{5}}
\put(40,60){\color{blue}\circle*{5}}
\put(20,100){\color{red}\circle*{5}}
\put(80,40){\color{red}\circle*{5}}

\linethickness{0.5mm}
\put(20,20){\line(-1,2){20}}
\put(20,20){\line(2,-1){40}}
\put(20,40){\line(-1,2){20}}
\put(20,40){\line(2,-1){80}}

\put(0,60){\color{blue}\circle*{6}}
\put(20,20){\color{blue}\circle*{6}}
\put(60,0){\color{blue}\circle*{6}}

\put(0,80){\color{blue}\circle*{6}}
\put(20,40){\color{blue}\circle*{6}}
\put(100,0){\color{blue}\circle*{6}}

\put(4,-10){\makebox(0,0){\small{$0$}}}
\put(24,-10){\makebox(0,0){\small{$1$}}}
\put(44,-10){\makebox(0,0){\small{$2$}}}
\put(64,-10){\makebox(0,0){\small{$3$}}}
\put(84,-10){\makebox(0,0){\small{$4$}}}
\put(104,-10){\makebox(0,0){\small{$5$}}}
\put(124,-10){\makebox(0,0){\small{$6$}}}
\put(144, -10){\makebox(0,0){\small{$7$}}}
\put(164, -10){\makebox(0,0){\small{$8$}}}
\put(184, -10){\makebox(0,0){\small{$9$}}}

\put(-10,5){\makebox(0,0){\small{$0$}}}
\put(-10,25){\makebox(0,0){\small{$1$}}}
\put(-10,45){\makebox(0,0){\small{$2$}}}
\put(-10,65){\makebox(0,0){\small{$3$}}}
\put(-10,85){\makebox(0,0){\small{$4$}}}
\put(-10,105){\makebox(0,0){\small{$5$}}}
\put(-10,125){\makebox(0,0){\small{$6$}}}
\put(-10,145){\makebox(0,0){\small{$7$}}}
\put(-10,165){\makebox(0,0){\small{$8$}}}

\put(110,130){\tikz \coordinate (p9);}
\put(190,130){\tikz \coordinate (p10);}
\put(190,150){\tikz \coordinate (p11);}
\put(110,150){\tikz \coordinate (p12);}
\end{picture}
\caption{In the case of our running example we consider the strong kerb configuration $\{{\color{blue}(0,7)},\,{\color{red}(1, 5)},\, {\color{blue}(2,3)},\,{\color{red}(4, 2)},\, {\color{blue}(8,0)}\} $, and assign to it the `blue' lattice points $F_0={\color{blue}(0,3)}, \ F_1={\color{blue}(1,1)},\ F_2={\color{blue}(3,0)}$ and $G_0={\color{blue}(0,4)},\ G_1={\color{blue}(1,2)},\ G_2={\color{blue}(5,0)}$. One can easily check that these satisfy conditions (\ref{eq:F+GinB}) and (\ref{eq:F+GinR}).} \label{fig:summands}
\end{figure}
\end{center}
\vspace{-15mm}

Then we automatically have that 
\begin{align}
    &\bullet \ F_i+G_i=(p^{2i}_1, p^{2i}_2)=p^{2i} \in \mathcal{B} &&\text{for all }0 \leq i \leq t;\label{eq:Fi+Gi}\\
    &\bullet \ F_{i+1} + G_{i}=p^{2i+1}\in \mathcal{R} && \text{for all }0 \leq i \leq t-1;\label{eq:Fi+1+Gi}\\
    &\bullet \ F_{i+1}-F_i=p^{2i+1}-p^{2i} &&\text{for all }0 \leq i \leq t-1;\label{eq:Fi+1-Fi}\\
    &\bullet \ G_{j+1}-G_j=p^{2j+2}-p^{2j+1} &&\text{for all }0 \leq i \leq t-1.\label{eq:Gj+1-Gj}
\end{align}
Thus, there only remains to prove property \ref{eq:F+GinR}: $F_i+G_j\in \mathcal{R}$ for all $i \notin \{j, j-1\}$. For this, notice that from the identities (\ref{eq:Fi+1-Fi}) and the (\ref{eq:Gj+1-Gj}) and sequence of inequalities (\ref{eq:kerbconfslopesor}) we get
\begin{equation}\label{eq:FGarevertices}
\begin{split}
    -\infty<\, m(\overline{F_0,F_1})\ <\ m(\overline{F_1,F_2})\ &\,< \ldots <m(\overline{F_{t-1},F_t})\ <0;\\
    -\infty <m(\overline{G_0,G_1})<m(\overline{G_1,G_2})&\,< \ldots <m(\overline{G_{t-1},F_t})<0.
\end{split}
\end{equation}
Hence, we can consider the following unbounded convex `Newton' polytopes:
\begin{align*}
    P_F:=&\,N\Gamma_+\big(\{F_i\}_{i=0}^t\big)={\rm ConvHull}\Big(\bigcup_{i=0}^t(F_i+\big(\mathbb{R}_{\geq 0})^2\big)\Big)\\=&\,\{p\in \mathbb{R}_{\geq 0}\times\mathbb{R}_{\geq p^{2t}_2}\,:\,\langle p, (-1,m(\overline{F_i, F_{i+1}})\rangle \geq \langle F_i, (-1,m(\overline{F_i, F_{i+1}})\rangle \ \forall\, 0\leq i <t\}; \\
    P_G:=&\,N\Gamma_+\big(\{G_j\}_{j=0}^t\big)={\rm ConvHull}\Big(\bigcup_{j=0}^t(G_j+\big(\mathbb{R}_{\geq 0})^2\big)\Big)\\=&\,\{p\in \mathbb{R}_{\geq p^0_1}\times\mathbb{R}_{\geq 0}\,:\,\langle p, (-1,m(\overline{G_j, G_{j+1}})\rangle \geq \langle G_j, (-1,m(\overline{G_j, G_{j+1}})\rangle \ \forall\, 0\leq j <t\};\\
    P_\mathcal{K}:=&\,N\Gamma_+\big(\{p^l\}_{l=0}^{2t}\big)={\rm ConvHull}\Big(\bigcup_{l=0}^{2t}(p^l+\big(\mathbb{R}_{\geq 0})^2\big)\Big)={\rm ConvHull}\big(\mathcal{K}\big)\\=&\,\{p\in \mathbb{R}_{\geq p^0_1}\times\mathbb{R}_{\geq p^{2t}_2}\,:\,\,\langle p, (-m(\overline{p^l p^{l+1}}), 1)\rangle \geq b(\overline{p^l p^{l+1}}) \ \forall l:0\leq l\leq 2t-1\}.
\end{align*}
We claim that the Minkowski sum $P_F+P_G=\{ p+p'\in \mathbb{R}^2\,:\, p \in P_F, \ p' \in P_G\}$ agrees with $P_{\mathcal{K}}$. Indeed, by the inequalities from (\ref{eq:FGarevertices}), the lattice points $\{F_i\}_{i=0}^t$ are exactly the vertices of $P_F$ (respectively, the lattice points $\{G_j\}_{j=0}^t$ are exactly the vertices of $P_G$). Therefore, the algorithm which produces the Minkowski sum of convex planar polytopes by gluing their edges in the order of their slopes (see, e.g., \cite[Section 3]{MinkSum}) yields $P_{\mathcal{K}}$ via the identities (\ref{eq:Fi+Gi})--(\ref{eq:Gj+1-Gj}) and property (\ref{eq:kerbconfslopesor}) of kerb configurations. 

Consequently, for any $0\leq i,j \leq t$ we have 
\begin{center}$F_i+G_j \in (P_F +P_G)\cap \mathbb{Z}^2=P_{\mathcal{K}}\cap \mathbb{Z}^2=\mathcal{K}$.
\end{center}
Moreover, by \cite[Proposition 12.1]{Fukuda}, the description of any vertex $p^l$ of $P_{\mathcal{K}}$ as a sum $p^l=F_{i}+G_{j}$ of vertices of $P_F$ and $P_G$ is unique, in particular, the decomposition of $p^{2l}=F_l+G_l$, for all $0\leq l \leq t$, is also unique (indeed, by property (\ref{eq:kerbconfslopesor}) of kerb configurations these blue lattice points $\{p^{2l}\}_{l=0}^t$ are vertices of $P_{\mathcal{K}}$). Thus, if $i \neq j$, then $F_i+G_j \in \mathcal{K}\setminus \{p^{2l}\}_{l=0}^t$, which lattice point set, by the distinguishing property of \emph{strong} kerb configurations, is contained in $\mathcal{R}$. Hence, we have proved both properties (\ref{eq:F+GinB}) and (\ref{eq:F+GinR}), thus, we can use the second part of Corollary \ref{cor:htopnem0sym}, to obtain the identity
    \begin{equation*}
    {\rm homdim}(\mathcal{M}\triangleleft \mathcal{O}_2)=t={\rm skerbdim}(\mathcal{R})={\rm linedim}(\mathcal{R})-1. \qedhere
    \end{equation*}
 \end{proof}

\begin{remark} \label{rem:minrealmonid}
    Notice that in the proof of Proposition \ref{prop:linedim-1=kerbdim} we actually produced a strong kerb configuration $\{p^l\}_{l=0}^{2t}$ with maximal length $t={\rm homdim}(\mathcal{M}\triangleleft\mathcal{O}_2)$ for any a realizable monomial ideal $\mathcal{M}$ of $\mathcal{O}_2=k[x_1, x_2]$. Moreover, by construction, for a small enough rational $\varepsilon\in\mathbb{Q}_{>0}$ the set of lines $\mathfrak{L}:=\{L(m^{i}-\varepsilon)\,:\,0\leq i\leq t\}$ gives a minimal geometric realization of $\mathcal{R}={\rm supp}(\mathcal{M}) \subset (\mathbb{Z}_{\geq0})^2$ in the sense of Definition \ref{def:geomreal} (here we use the notations of the proof of Proposition \ref{prop:linedim-1=kerbdim}, i.e., $m^{i}$ is the slope $m(\overline{p^{2i} p^{2i+1}})$ of the line $\overline{p^{2i} p^{2i+1}}$ and the `tangential' line $L(m)$ with slope $m$ is as defined in  (\ref{eq:tanline(m)})). Thus, via Observation \ref{obs:geomreal} and Definition \ref{def:geomreal}, we can produce `by hand' a monomial realization $\mathcal{D}_{\mathfrak{L}}$ of $\mathcal{M}$ of minimal cardinality $|\mathcal{D}_{\mathfrak{L}}|={\rm homdim}(\mathcal{M}\triangleleft \mathcal{O}_2)$ using this recipe.
\end{remark}

\begin{remark}\label{rem:homdimhighdim}
    We expect Theorem \ref{th:homdim} to also be true in the higher dimensional setting, i.e., for  finite codimensional integrally closed monomial ideals in $\mathcal{O}_n=k[x_1, \ldots, x_n]$. In such a case, however, one has to find the correct generalization of the (strong) kerb configurations. Another difficulty is that one cannot sweep through the `tangential hyperplanes' just by moving a single real valued (thus totally ordered) parameter (in the $n=2$ case this was the \emph{slope} of the `tangential line'). 

    Nevertheless, we point out the encouraging examples of $\mathcal{M}(l)=\overline{(x_1^l, x_2^l, x_3^l, x_1x_2x_3)} \triangleleft \mathcal{O}_3$, where $l\in\{3,4,5\}$. For $l=3$, we  have $\mathcal{M}(3)=\mathfrak{m}_3^3=(x_1, x_2, x_3)^3$ which can trivially be realized by a single valuation, thus ${\rm homdim}(\mathcal{M}(3)\triangleleft\mathcal{O}_3)=0$. On the other hand, for $l=5$ we can use the three facets of the Newton boundary to provide a geometric realization, and we can consider the following combination of ring elements:
    $$f_1=g_1:=x_1x_2, \ f_2=g_2:=x_1x_3, \ f_3=g_3:=x_2x_3.$$
    Clearly, $f_1g_1, \ f_2g_2, \ f_3g_3 \notin \mathcal{M}(5)$ and $f_1g_2, \ f_1g_3, \ f_2g_1, \ f_2g_3, \ f_3g_1, \ f_3g_2 \in \mathcal{M}(5)$, thus, by Corollary \ref{cor:htopnem0sym}, we have ${\rm homdim}(\mathcal{M}(5)\triangleleft\mathcal{O}_3)=2$. Finally, one can check that the support of $\mathcal{M}(4)$ can be geometrically realized by two planes, but not with a single one, and that, correspondingly, it has $\mathbb{SH}_{1,12}(\mathcal{M}(4)\triangleleft\mathcal{O}_3)=\mathbb{Z}^4$ (thus ${\rm homdim}(\mathcal{M}(4)\triangleleft\mathcal{O}_3)=1$). Its graded root is the following:
    
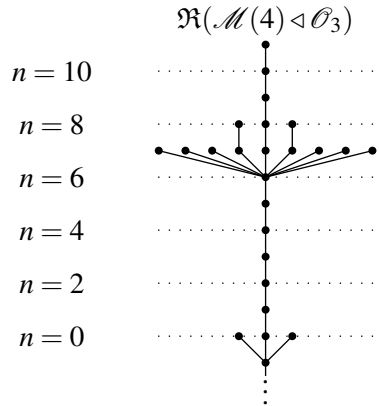
\begin{center}
    \begin{picture}(200,155)(0,-20)
\linethickness{0.5pt}

\put(130,116){\circle*{3}}
\put(130,106){\circle*{3}}
\put(130,96){\circle*{3}}
\put(130,86){\circle*{3}}
\put(140,86){\circle*{3}}
\put(120,86){\circle*{3}}
\put(90,76){\circle*{3}}
\put(100,76){\circle*{3}}
\put(110,76){\circle*{3}}
\put(120,76){\circle*{3}}
\put(130,76){\circle*{3}}
\put(140,76){\circle*{3}}
\put(150,76){\circle*{3}}
\put(160,76){\circle*{3}}
\put(170,76){\circle*{3}}
\put(130,66){\circle*{3}}
\put(130,56){\circle*{3}}
\put(130,46){\circle*{3}}
\put(130,36){\circle*{3}}
\put(130,26){\circle*{3}}
\put(130,16){\circle*{3}}
\put(130,06){\circle*{3}}
\put(120,06){\circle*{3}}
\put(140,06){\circle*{3}}
\put(130,-4){\circle*{3}}

\qbezier[20](90,106)(130,106)(170,106)
\qbezier[20](90,86)(130,86)(170,86)
\qbezier[20](90,66)(130,66)(170,66)
\qbezier[20](90,46)(130,46)(170,46)
\qbezier[20](90,26)(130,26)(170,26)
\qbezier[20](90,6)(130,6)(170,6)

\put(50,106){\makebox(0,0)[0]{$n=10$}}
\put(50,86){\makebox(0,0)[0]{$n=8$}}
\put(50,66){\makebox(0,0)[0]{$n=6$}}
\put(50,46){\makebox(0,0)[0]{$n=4$}}
\put(50,26){\makebox(0,0)[0]{$n=2$}}
\put(50,6){\makebox(0,0)[0]{$n=0$}}

\put(120,86){\line(0,-1){10}}
\put(140,86){\line(0,-1){10}}
\put(130,116){\line(0,-1){125}}
\put(120,76){\line(1,-1){10}}
\put(110,76){\line(2,-1){20}}
\put(100,76){\line(3,-1){30}}
\put(90,76){\line(4,-1){40}}
\put(140,76){\line(-1,-1){10}}
\put(150,76){\line(-2,-1){20}}
\put(160,76){\line(-3,-1){30}}
\put(170,76){\line(-4,-1){40}}
\put(120,06){\line(1,-1){10}}
\put(140,06){\line(-1,-1){10}}

\put(130,-12){\makebox(0,0){$\vdots$}}

\put(130,125){\makebox(0,0)[0]{$\mathfrak{R}(\mathcal{M}(4)\triangleleft\mathcal{O}_3)$}}
\end{picture}
\end{center}
\end{remark}

\section{Deformations and 
superisolated surface singularities}\label{s:NNSI}

 In this section we analyze the behaviour of the different analytic lattice homology theories  with respect to certain complex analytic singularity deformations. Here we will focus mostly on the curve and surface cases, and in the surface case we will include in the picture the topological lattice homology as well.  

We will analyze several types of deformations.
The introduction of those deformations  which are defined by the property that they preserve the 
conductor ideal
is a novelty  in the literature  (as far as 
we know),  and it is imposed by the results of the present manuscript. Our guiding question is: for which deformations is the lattice homology stable?

A particularly interesting class of singularities with deformations worth analyzing
is the family of superisolated surface singularities 
and their (canonical) deformations. We start this section with a short review of their main properties.

\subsection{Superisolated singularities}\label{ss:superiso} \cite{ArtalJordan,ALM,Ignacio,NBook}\, 

A hypersurface singularity $(X,0)\subset (\C^3,0)$ is called \textit{superisolated} if the strict transform $\widetilde{X}_{0}$ of $(X,0)$ under 
the blow-up at the origin $0\in \C^3$ (or, in other words, the blow-up of the maximal ideal) 
is smooth. This condition already implies that 
$(X,0)$ must have an isolated singularity. Moreover, it turns out, that if $(X,0)$ is not already smooth then $\widetilde{X}_{0}\to X$ is the minimal resolution $\widetilde{X}_{min}\to X$ of $(X,0)$. 

A different characterization  is the following. Assume that $(X,0)$ is given as
$(\{f=0\},0)$, where $f:(\C^3,0)\to (\C,0)$ is a holomorphic map germ. Let us write 
$f$ as $f_d+f_{d+1}+\cdots$, where $\{f_e\}_e$ are its homogeneous components of degree  $e$, with $f_d\not\equiv 0$. Then $(X,0)$ is superisolated if and only if the projective curve 
$C:=\{f_d=0\}$ in $\mathbb{P}^2$ is reduced,   and 
its isolated singular points $\{p_i\}_i$   are not situated on $\{f_{d+1}=0\}$.
It turns out that the embedded topological type 
(hence, in particular, the link $L_X$ of $(X,0)$) is independent of the choice of the higher degree homogeneous parts
$\{f_{e}\}_{e\geq d+1}$ (under the assumption that $\{f_{d+1}=0\}\cap {\rm Sing}(C)=\emptyset$).
Hence, it basically depends only on $C$. 
Via this identification we can `embed' the theory of 
reduced projective plane curves into the theory of isolated hypersuface surface singularities. In fact, this 
connection  turns out to bear fruit in the classification of projective plane curves. 
For additional  connections and relations
 with lattice (co)homology theory see,  e.g., \cite{Spany, BLMN2}.

The family of superisolated germs plays a key role in the theory of normal surface singularities as well: besides that it creates a bridge between singular projective plane curves and surface singularities, it is one of the best `testing families' for different 
conjectures, properties, see e.g., \cite{ArtalJordan,ALM,Ignacio,LMNsi,NBook} and references therein. 

\bekezdes {\bf Some numerical invariants.}
In the minimal resolution $\widetilde{X}_{min}\to X$ the exceptional curve can be identified with $C$.
Even more, for any singular point $p_i$
of $C$ the embedded topological types of $(C,p_i)\subset (\mathbb{P}^2,p_i)$
and of $(C,p_i)\subset (\widetilde{X}_{min},p_i)$ agree. 
Hence the minimal \emph{good} resolution $\widetilde{X}\to X$ of $(X,0)$ can be obtained 
from $\widetilde{X}_{min}$ by resolving these embedded local plane curve singularities (following the resolution steps of plane curve singularities). 

Let $C_1\cup\cdots\cup  C_r$ be the irreducible decomposition of 
the projective curve $C$, and let $d_i$ be the degree of $C_i$ in  $\mathbb{P}^2$. Clearly $\sum_id_i=d$. Then the irreducible curves $\{C_i\}_{i=1}^r$, 
considered in  $\widetilde{X}_{min}$ as exceptional divisors of the minimal resolution, have the following intersection properties:
$(C_i,C_j)_{\widetilde{X}_{min}}= d_id_j$ for $i\not=j$ and 
$(C_i,C_i)_{\widetilde{X}_{min}}= - d_i(d-d_i+1)$. In particular, 
$(C,C)_{\widetilde{X}_{min}}= -d$. 
 Moreover, in $\widetilde{X}_{min}$ 
one also has $Z_K=(d-2)C$, denoted by $Z_{K, min}$. Therefore, 
$Z_{K,min}^2=-(d-2)^2d$. 
This shows that the numerical link invariant, 
computed in any resolution as  
$Z_K^2+|\mathcal{V}| $ (where $\mathcal{V}$ is the set of irreducible exceptional divisors of that resolution), equals 
$Z_{K,min}^2+r=-(d-2)^2d+r$. (For more about this invariant consult \cite[Section 6.3.A]{NBook} and the references therein.)

It is also known that the geometric genus $p_g$ of a superisolated germ  $(X,0)$ is 
$d(d-1)(d-2)/6$, i.e., it is independent of the concrete form of the equation, it depends only on $d$.

If we denote the dual resolution graph of a {\it good} resolution by $\Gamma$, and $\{v_{C_1}, \ldots , v_{C_r}\}$ 
are the  vertices of $\Gamma$ which corresponds to the strict transforms of the irreducible  exceptional curves $\{C_i\}_{i=1}^r$ of $\widetilde{X}_{min}$, then 
$\Gamma\setminus\{v_{C_1}, \ldots , v_{C_r}\}$ is a collection of connected rational graphs (they are exactly the 
minimal embedded resolution graphs of the plane curve singularities $\{(C,p_i)\subset (\widetilde{X}_{min},p_i)\}_i$). Hence, $\{v_{C_i}\}_{i=1}^r$ is a WR-set, it satisfies 
(\ref{eq:hegy}) and (\ref{eq:dbar}). In particular, the corresponding Reduction 
Theorem \ref{prop:redhigh} for the analytic lattice homology holds.

\begin{example} \label{ex:superisolated} (a) 
If $C$ is smooth then $(X,0)$ belongs to a special  family of 
superisolated singularities. In this case the minimal resolution
$\widetilde{X}_{min}$ is already good, and the embedded topological type 
is independent of the choices of $\{f_e\}_{e\geq d+1}$, and even of the choice of 
$f_d$ (under the assumption that $C$ remains smooth). Hence this embedded topological type coincides with  the type of $(\{x_1^d+x_2^d+x_3^d=0\}, 0)$. 
In the minimal resolution the exceptional curve is irreducible of genus 
$(d-1)(d-2)/2$ and it has self-intersection $-d$.  Also, 
$Z_{K}^2+|\mathcal{V}|=-(d-2)^2d+1$ (and $p_g=d(d-1)(d-2)/6$).

(b)  
An important family of superisolated singularities consists of those germs for which 
$C$ is irreducible and all of its singular points are cusps (i.e., locally irreducible). If $C$ is additionally rational, then the link $L_X$ of $(X,0)$ 
is a rational homology  sphere (its resolution graph is a tree with all genus decorations being $0$).  In fact,  it has the following surgery description: $L_X=S^3_{-d}(\#_iK_i)$, where 
$\{K_i\subset S^3\}_i$
are the links of the local embedded topological types at $p_i\in C$. 
\end{example}

\subsection{Deformations of superisolated singularities} 
\label{bek:SIdef}\,

Superisolated surface singularities come with two natural deformations, in both cases the deformation happens at the homogeneous level of $f_d$. We recall these constructions below.

\bekezdes \label{bek:NAIV} {\bf The `naive' deformation.}
Assume that $f$ defines a superisolated singularity with $C$ not smooth. 
Let $\{f_{d,t}\}_{t\in(\C,0)}$ be a deformation of $f_d$ such that $f_{d,t}$ is homogeneous of degree $d$ for every $t \in (\C, 0)$, $f_{d,t=0}=f_d$, and 
the projective curve $\{f_{d,t\not=0}=0\}$ is smooth. 
(E.g., $f_{d,t}$ is a small generic homogeneous deformation of $f_d$.)
Then set the deformation $f_t=f_{d,t}+f_{d+1}+\cdots $ \ of $f$ with $t \in (\mathbb{C}, 0)$. By the previous arguments, for $t\not=0$ the map germ  
$f_t:(\bC^3,0)\to (\bC,0)$ defines an isolated hypersurface singularity whose embedded topological type 
agrees with that of  $(\{x_1^d+x_2^d+x_3^d=0\}, 0)$.

Along this deformation the embedded topological types, the links  and the Milnor numbers are not 
stable, however, $p_g=d(d-1)(d-2)/6$ and $Z_{K,min}^2=-d(d-2)^2$ remain constant. 
Note also, that unless $d\leq 2$, the link of $(X_{t\not=0},0)$ is not a rational homology sphere.
(As we already saw in Example \ref{ex:superisolated} (a), in the minimal resolution the exceptional curve is irreducible of genus 
$(d-1)(d-2)/2$ and self-intersection $-d$.) 

\bekezdes \label{bek:nonNAIV} {\bf The `non-naive' deformation.}
Assume that $f(x_1, x_2, x_3)=0$ defines a superisolated surface singularity with the coordinates $(x_1,x_2,x_3)$ chosen in such a way that $x_3=0$
intersects the projective curve $C$ transversely along smooth points of $C$ and $f_{d+1}(0,0,1)\not=0$.
Then consider the deformation $f_t(x_1,x_2,x_3):=f_d(x_1,x_2,tx_3)+\sum_{e\geq d+1}f_e(x_1,x_2,x_3)$.

Then for each $t$ with $|t|$ sufficiently small 
$f_t$ defines a superisolated singularity. 
The deformation does not preserve the embedded topological types, the links or the Milnor numbers, however, it is $p_g$-constant and $Z_{K,min}$-constant. 
Moreover, the embedded topological type (hence the link, too) 
of $(X_{t=0},0)$ coincides with that of $(\{x_1^d+x_2^d+x_3^{d+1}=0\}, 0)$, whose link  is a rational homology sphere. 
In particular, if $C=\{f_d=0\}$ is a rational cuspidal curve then for any 
$t$ the links are 
$\Q HS^3$'s. 
This has the advantage that along such deformation we can also analyze the  stability
of those invariants which are defined only for singularities with rational homology sphere links (e.g., the 
topological lattice homology). 

\vspace{2mm}

Note that if we start with a superisolated singularity $(X,o)$, then in the 
`naive' deformation this germ appears as the central fiber of the deformation, 
while in the `non-naive' deformation the nearby fiber has this topological type. 

\subsection{Certain types of deformations}\label{ss:TYPE}\,

In this subsection we introduce several constraints on general singularity deformations, which make them somehow accessible or interesting, especially from the point of view of lattice homology. A novelty of this manuscript is the definition of `conductor ideal constant deformations', which come directly from the previous results (in particular, under these types of deformations the analytic lattice homologies remain constant). We also compare them with more classical notions, such as $p_g$-constant or $\delta$-constant deformations. In order to discuss these properly we have to fix the setting and some notations.

Assume that $F:(\mathcal{X},0)\to (\C,0)$ is a deformation of
$F^{-1}(0)=(X,0)$, which factors through an inclusion 
$({\mathcal{X},0})\hookrightarrow (\C^m\times \C,(0,0))$ and the projection to the second factor $ (\C^m\times \C,(0,0))\to (\C,0)$,
such that each fiber $(F^{-1}(t),0)=(X_t,0)\subset (\C^m,0)$ is reduced with an isolated singularity.
Let us denote by $q_t:\cO_{\C^m,0}\to \cO_{X_t,0}$ the natural quotient morphism. 

In the next discussions we will assume for simplicity  that each
singularity $(X_t,0)$ is Gorenstein. Moreover, conforming to sections \ref{s:deccurves} and \ref{s:dnagy}, we also require the following:

\begin{ass}\label{ass:singsallowed}
    We assume that all of our studied complex analytic curve singularities are reduced, whereas all higher dimensional singularities are isolated and normal.
\end{ass}

\begin{define}\label{def:deformations} We say that 

(1) $F$ is a {\it `conductor ideal constant deformation'} if the pullback
$q_t^{-1}(\mathcal{C}_{X_t,0})\subset \cO_{\C^m,0}$ is independent of the parameter $t$,
where $\mathcal{C}_{X_t,0}$ is the conductor ideal of $(X_t,0)$ (for singularities of dimension $\geq 2$ this was defined in Definition \ref{def:CONDINDEF}, for curve singularities in paragraph \ref{par:conductor}--- see also the adjunction conductor ideal of Merle and Teissier presented in Theorem \ref{th:MerleTeissier} ); 

(2) (a) ($\dim_{\mathbb{C}} (X_t,0)=1$ case) $F$ is a {\it $\delta$-constant deformation} if the delta invariant $\delta(X_t,0)$ is independent of the parameter $t$;

(2) (b) ($\dim_{\mathbb{C}} (X_t,0)\geq 2$ case) $F$ is a {\it $p_g$-constant deformation} if the geometric genus $p_g(X_t,0)$ is independent of the parameter $t$;

(3)   ($\dim_{\mathbb{C}} (X_t,0)= 2$ case) $F$ is a {\it $Z_{K,min}^2$-constant deformation} if the numerical invariant  $Z_{K,min}^2$ of 
$(X_t,0)$ (i.e., $Z_K^2=(Z_K, Z_K)_{\widetilde{X}_{t, min}}$  computed in its minimal resolution 
$\widetilde{X}_{t, min}\to X_t$)  is independent of the parameter $t$. 
\end{define}

The conductor ideal constant deformations are motivated by the previous theorems of this manuscript, see, e.g.,
 Corollaries \ref{cor:curveArtin} and  \ref{cor:GORCOND}. In particular, along such deformations the analytic lattice homology remains constant, moreover, as its Euler characteristic (or the codimension of the conductor ideal), so does the $\delta$ invariant, or, in the higher dimensional case, the geometric genus $p_g$.

On the other hand, the interest in  $Z_{K,min}^2$-constant deformations in our context is motivated by the following fact. 
Let  $F:(\mathcal{X},0)\to (\C,0)$ be  the germ of a flat deformation of a normal 
Gorenstein two-dimensional singularity $(X,0)$. 
Then, by a result of Laufer \cite[Theorem 5.7]{LauferWeak}, 
$F$ admits a very weak simultaneous resolution (possibly after a finite base change)
if and only if $Z_{K,min}^2$ is constant along $F$. 
(This can be compared with  a theorem of Teissier, which says that a deformation of reduced curve singularities is equinormalizable if and only if it is $\delta$-constant \cite{Teissiersimultan}.)
Having such a simultaneous resolution implies that one can find a natural linear map between 
the lattices
associated with 
the resolutions of  $(X_{t\not=0},0)$ and $(X_{t=0},0)$ \cite{Teissiersimultan}, which might help to define  a connecting functor (a natural bigraded $\Z[U]$-module homomorphism)
between the analytic (or topological) lattice homologies of the central and the generic fiber. 

\subsection{Conductor ideal constant deformations}\, \label{ss:CondIdConst}

In this subsection we will study, through numerous examples, the main properties of conductor ideal constant deformations. E.g., we will prove that both the naive and non-naive deformations of superisolated singularities are of this kind and compare this new notion with the other types of deformations introduced in the previous subsection.

First, though, we formalize the following consequence of our previous theorems (especially Corollaries \ref{cor:curveArtin} and  \ref{cor:GORCOND}): 

\begin{theorem}\label{th:DEFCOND} 
Let $F:(\mathcal{X}, 0)\rightarrow (\mathbb{C}, 0)$ be a conductor ideal constant deformation of the $n$-dimen\-sional Gorenstein germs $(X_t,0)$ satisfying Assumption \ref{ass:singsallowed} as above. Then
\begin{itemize}
    \item if $n=1$, then the analytic lattice homology $\mathbb{H}_{an,*}(X_t, 0)$ stays constant;
    \item if $n\geq 2$, then  the analytic 
lattice homology of $n$-forms 
$\bH_{*}((X_t,0);\Omega^n)$  stays constant.
\end{itemize}
\end{theorem}

Before showing some concrete examples, let us recall also the following corollary, which will be our main technical tool in this section:

\begin{cor}[= Corollary \ref{cor:HanstabforNN}]
    The analytic lattice homology (of $n$-forms) is constant along an equisingular deformation of isolated hypersurface singularities with Newton nondegenerate principal part associated with a fixed Newton boundary $N\Gamma$. Moreover, it agrees with the symmetric lattice homology of the adjoint of the monomial ideal supported on $N\Gamma_+$ (for notations and terminology see section \ref{s:ND}).
\end{cor}

\begin{example} The isomorphism of the lattice homologies of the pair of curves
$C_1$ and $C_2$ from Example \ref{bek:10.2.1} can be seen from the viewpoint of Theorem \ref{th:DEFCOND}  as well. Indeed, one can consider the deformation 
$(x_1^2-x_2^3-tx_2^2)(x_1^3-x_2^2)=0$ with parameter $t \in (\mathbb{C}, 0)$ connecting $C_1 $ and an intermediary curve
$C_3:=\{(x_1^2-x_2^2-x_2^3)(x_1^3-x_2^2)=0\}$ and also a similar deformation  
  connecting $C_2$ and $C_3$. 
Both deformations preserve the pullbacks of the conductor ideals: 
$q^{-1}(\mathcal{C})$ is $\mathfrak{m}_{\C^2,0}^3$ in all these cases (where $\mathfrak{m}_{\C^2,0}$ denotes the maximal ideal of the local ring $\mathcal{O}_{\mathbb{C}^2, 0}$).
\end{example}

\begin{example}
Consider the curves $(C,0)$ and $(C', 0)$ from Example \ref{bek:10.6.3}.
A computation shows that the pullback of their conductor ideals
is the same $q^{-1}(\mathcal{C})={\rm adj}(\mathcal{\overline{M}})$, where $\mathcal{M}$ is the monomial ideal generated by $x_1^6$ and $x_2^4$. 

Consider now the one parameter family of plane curve singularities
$$f_a(x_1,x_2)=(x_1^3-x_2^2)(x_1^3-x_2^2+ax_2^2) -4x_1^5x_2-x_1^7=0,$$ where the parameter $a\in(\mathbb{C},0)$.
Then, along this deformation the adjunction conductor ideal $\mathcal{C}_{f_a}$ in $\cO_{\mathbb{C}^2,0}$ is constant.
Hence, the lattice homology modules must stay constant as well. Thus, the isomorphism (\ref{eq:ISOUJ}) exemplifies Theorem \ref{th:DEFCOND}, too.
\end{example}

\begin{example}\label{ex:NNSurfaces} (a)
 Similar examples can be given for surface singularities as well.    Let us consider the germs 
$\{x^{13}+y^{13}+x^3y^2+x^2y^3+z^3=0\}$,  $\{x^{14}+y^{14}+x^2y^2+z^3=0\}$ and
$\{x^{13}+y^{9}z+x^2y^2+z^3=0\}$. All of them are  conductor ideal constant deformations of
the germ $\{x^{13}+y^{13}+x^2y^2+z^3=0\}$. Indeed, they all have Newton nondegenerate principal parts and the corresponding $\mathcal{P}$-sets (defined in paragraph \ref{par:NewtonComb}) agree. Then we can use the Theorem \ref{th:MerleTeissier} of Merle and Teissier.

  (b) The deformation  $\{x_1^d+x_2^d+tx_3^d+x_3^{d+1}=0\}$ is a conductor ideal  constant deformation (once again, by Newton nondegeneracy). The two germs corresponding to 
$t\not=0$ and $t=0$ appear in both the 
naive and non-naive deformations of superisolated singularities  from \ref{bek:NAIV} and \ref{bek:nonNAIV}.  (In fact, this is the naive deformation of $x_1^d+x_2^d+x_3^{d+1}$.)
\end{example}

\begin{proposition}\label{prop:11.5.5}
Assume that $(X,o) $ is a superisolated surface singularity.
Then both the naive and non-naive deformations from \ref{bek:NAIV} and \ref{bek:nonNAIV} are conductor ideal constant. In fact, the pullback of the corresponding conductor ideals is $\mathfrak{m}_{\bC^3,0}^{d-2}\triangleleft \cO_{\bC^3,0}$. 

Therefore, the analytic lattice homology  of a superisolated surface singularity  given by the equation
$f=f_d+\cdots=0$ agrees with the analytic lattice homology of the germ
 $\{x_1^d+x_2^d+x_3^{d+1}=0\}$ (or, equivalently,  of  $\{x_1^d+x_2^d+x_3^d=0\}$), and it also agrees with the 
 symmetric lattice homology of the integrally closed monomial ideal $\mathfrak{m}_{\mathbb{C}^3, 0}^{d-2}\triangleleft \cO_{\mathbb{C}^3, 0}$.

This also shows that the analytic lattice homology of $\{f=0\}$ depends merely on the degree  $d$, it
does not depend on the choice of the particular  
equation.
\end{proposition}
\begin{proof}
For singularities of type  $\{x_1^d+x_2^d+x_3^{d+1}=0\}$  or  $\{x_1^d+x_2^d+x_3^d=0\}$,
the identity $\mathcal{C}_f=\mathfrak{m}_{\mathbb{C}^3, 0}^{d-2}$ already follows from a  Newton diagram argument. 
Let us verify this identity for an arbitrary  superisolated germ $(X, 0)$ as well. 
By Proposition \ref{prop:CONDINDEP} we can do the computations in the minimal resolution
$\phi_{min}:\widetilde{X}_{min}\to X$. 
In this case, the Gorenstein form 
$\omega_G$ (cf. its formula in subsection \ref{ss:NNPP})
has a pole of order $(d-2)$ along $C$. That is, $Z_K=(d-2)C$. Hence, the conductor ideal consists 
of functions $\varphi$ such that $\phi_{min}^*\varphi$ has a vanishing order $\geq (d-2)$ along $C$
(cf. Theorem \ref{th:MerleTeissier}). Since this vanishing order in the case of a 
monomial of degree $e$ is $e$, we get that  $\mathcal{C}_f\supset\mathfrak{m}_{\mathbb{C}^3, 0}^{d-2}$. On the other hand, ${\rm codim}(\mathfrak{m}_{\mathbb{C}^3, 0}^{d-2} \hookrightarrow\mathcal{O}_{\mathbb{C}^3, 0})=d(d-1)(d-2)/6=p_g(X, 0)={\rm codim}(\mathcal{C}_f \hookrightarrow\mathcal{O}_{\mathbb{C}^3, 0})$ for every superisolated singularity.
\end{proof}

In the following example we compare the different types of deformations (cf. Definition \ref{def:deformations}):

\begin{example}\label{ex:LIST}
    (a)
    Since the codimension of the conductor ideal is the 
delta invariant, res\-pectively the geometric genus, a conductor ideal constant deformation is always a $\delta$--, respectively $p_g$-constant deformation.

(b) If $\dim_{\mathbb{C}} (X_t,0)=2$, then along a Gorenstein 
 $Z_{K,min}^2$-constant deformation the geometric genus stays constant as well \cite[Theorem 5.3]{LauferWeak}.

\vspace{2mm} 

 On the other hand, in general,  no other implication holds  between these properties.

 \vspace{2mm} 

 (c) Consider the following deformation of hypersurface singularities 
\begin{center}
$X_t=\{x_1^2+tx_2^3+ x_2^4+ x_3^r+x_1x_2x_3=0\}$, $\ r\geq 7, \ t\in (\mathbb{C}, 0)$. 
\end{center}
For every $t\in (\mathbb{C}, 0)$ the surface singularity $(X_t,0)$ is a 
cusp singularity with $p_g=1$ and conductor ideal the maximal ideal. 
Hence the deformation is $p_g$-- and 
conductor ideal constant, but not  $Z_{K,min}^2$-constant.
Indeed, $Z_{K,min}^2=-1$ for $t\not=0$,    while 
 $Z_{K,min}^2=-2$  for $t=0$. \\ (For the minimal resolution graphs consult \cite[pp. 60--62]{Dimca}.)

 (d) Consider the following deformation of isolated 
 plane curve singularities $x_1^3+(x_2^2+tx_1)^2$. Then, by 
 \cite[Example 4.7.5]{AgNeCurves}, this is a delta constant deformation but along this deformation the lattice homology is not constant, hence it
 is not a conductor ideal constant deformation.

 Indeed, this fact can be checked by direct calculation as well. In the $t=0$ case we have the $E_6$ singularity $C_{t=0}=\{x_1^3+x_2^4=0\}$ having conductor element $\mathbf{c}_{t=0}=6$, delta invariant $\delta(C_{t=0}, 0)=3$ and pulled back conductor ideal $q_{t=0}^{-1}(\mathcal{C}_{t=0})=\mathfrak{m}^2=(x_1^2, x_1x_2, x_2^2) \triangleleft\mathbb{C}\{x_1, x_2\}$. In the $t\neq 0$ case we can consider the new holomorphic coordinates $(x_2, y)$ with $y=x^2_2 + tx_1$, hence the equation has the form $y^2+\big(\frac{y-x_2^2}{t}\big)^3=0$. It is Newton nondegenerate with support and Newton boundary presented in the following diagram. 

 
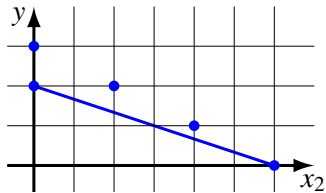
\begin{center}
 \begin{picture}(140,90)(-20,-20)

\linethickness{0.4mm}
\put(-10,0){\vector(1,0){115}}
\put(0,-10){\vector(0,1){70}}
\linethickness{0.4mm}
\put(45,15){\color{blue}\line(3,-1){45}}
\put(45,15){\color{blue}\line(-3,1){45}}

\put(100,-8){\makebox{$x_2$}}
\put(-8,55){\makebox{$y$}}

\linethickness{0.05mm}
  \multiput(-10,15)(0,15){3}{\line(1,0){115}}
  \multiput(15,-10)(15,0){6}{\line(0,1){70}}

\put( 0, 30){\color{blue}\circle*{4}}
\put(60, 15){\color{blue}\circle*{4}}
\put(90, 0){\color{blue}\circle*{4}}
\put( 0,45){\color{blue}\circle*{4}}
\put(30,30){\color{blue}\circle*{4}}
 \end{picture}
 \end{center}
 
 Thus, $C_{t\neq 0}$ is in fact the $A_5$ singularity with $\mathbf{c}_{t \neq 0}=(3,3)$ and $\delta(C_{t\neq 0}, 0)=3$. Moreover, the pullback $q_{t \neq 0}^{-1}(\mathcal{C}_{t\neq 0}) \triangleleft \mathbb{C}\{x_1, x_2\}$ agrees with the combinatorial conductor ideal (cf. Lemma \ref{lem:agreec} --- or, more generally, with the adjoint ideal from Theorem \ref{th:MerleTeissier}) $(y, x^3_2)=(x_2^3, x_2^2 + tx_1) \neq \mathfrak{m}^2$. Finally, the graded roots corresponding to the $t=0$ and $t \neq 0$ cases are the following:

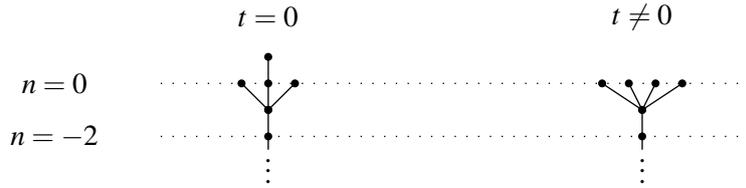
\begin{center}
\begin{picture}(350,65)(0,5)
\linethickness{.5pt}

\put(130,50){\circle*{3}}
\put(120,40){\circle*{3}}
\put(130,40){\circle*{3}}
\put(140,40){\circle*{3}}
\put(130,30){\circle*{3}}
\put(130,30){\circle*{3}}
\put(130,20){\circle*{3}}

\put(130,10){\makebox(0,0){$\vdots$}}

\qbezier[50](90,40)(200,40)(310,40)
\qbezier[50](90,20)(200,20)(310,20)

\put(50,40){\makebox(0,0)[0]{$n={0}$}}
\put(50,20){\makebox(0,0)[0]{$n=-2$}}

\put(130,50){\line(0,-1){35}}
\put(120,40){\line(1,-1){10}}
\put(140,40){\line(-1,-1){10}}

\put(255,40){\circle*{3}}
\put(265,40){\circle*{3}}
\put(275,40){\circle*{3}}
\put(285,40){\circle*{3}}
\put(270,30){\circle*{3}}
\put(270,30){\circle*{3}}
\put(270,20){\circle*{3}}

\put(270,30){\line(-1,2){5}}
\put(270,30){\line(-3,2){15}}
\put(270,30){\line(1,2){5}}
\put(270,30){\line(3,2){15}}
\put(270,30){\line(0,-1){15}}

\put(270,10){\makebox(0,0){$\vdots$}}

\put(130,65){\makebox(0,0)[0]{$t=0$}}
\put(270,65){\makebox(0,0)[0]{$t \neq 0$}}
\end{picture}
\end{center}

 (e) Next we construct a surface singularity deformation along which the geometric genus 
 and  $Z_{K,min}^2$ are constant but the lattice homology and the conductor ideal not. 
 
 Consider the deformation of isolated surface singularities 
 $x_1^4+(x_2^3+tx_1)^3+x_3^3=0$. For $t=0$, this gives the Brieskorn-type singularity $X_{t=0}=\{x_1^4+x_2^9+x_3^3\}$. This is Newton nondegenerate with its Newton boundary consisting of a single facet having primitive normal vector $\mathbf{n}=(9, 4, 12)$. Hence, by the results of section \ref{s:ND} (especially Theorems \ref{th:MerleTeissier} and \ref{th:comparison}), one can easily compute the geometric genus $p_g(X_{t=0})=4$, the `conductor coefficient' $c=12$ (cf.  \ref{bek:conductorofwhGorss}) and the (combinatorial) Hilbert function inducing (by Gorenstein symmetry) the weight function of the lattice homology:

 \begin{center}
     \begin{tabular}{c|ccccccccccccc}
          $\ell$ & $0$ & $1$ & $2$ & $3$ & $4$ & $5$ & $6$ & $7$ & $8$ & $9$ & $10$ & $11$ & $12$ \\
          \hline 
          $\mathfrak{h}(\ell)$ & $0$ & $1$ & $1$ & $1$ & $1$ & $2$ & $2$ & $2$ & $2$ & $3$ & $4$ & $4$ & $4$ \\
          $\mathfrak{h}^{sym}_{c}(\ell)$ & $4$ & $4$ & $4$ & $3$ & $2$ & $2$ & $2$ & $2$ & $1$ & $1$ & $1$ & $1$ & $0$ \\
           \hline 
          $w(\ell)$ & $0$ & $1$ & $1$ & $0$ & $-1$ & $0$ & $0$ & $0$ & $-1$ & $0$ & $1$ & $1$ & $0$
     \end{tabular}
 \end{center}
 The minimal (good) resolution graph is the following:
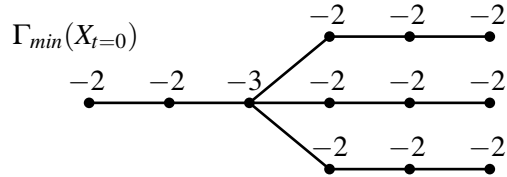
\begin{center}
\begin{picture}(200,70)(0,10)
\linethickness{1pt}

\put(85,40){\circle*{4}}
\put(55,40){\circle*{4}}
\put(25,40){\circle*{4}}
\put(115,40){\circle*{4}}
\put(145,40){\circle*{4}}
\put(175,40){\circle*{4}}
\put(115,15){\circle*{4}}
\put(145,15){\circle*{4}}
\put(175,15){\circle*{4}}
\put(115,65){\circle*{4}}
\put(145,65){\circle*{4}}
\put(175,65){\circle*{4}}

\put(24,48){\makebox(0,0)[0]{$-2$}}
\put(54,48){\makebox(0,0)[0]{$-2$}}
\put(83,48){\makebox(0,0)[0]{$-3$}}
\put(114,48){\makebox(0,0)[0]{$-2$}}
\put(144,48){\makebox(0,0)[0]{$-2$}}
\put(174,48){\makebox(0,0)[0]{$-2$}}
\put(115,23){\makebox(0,0)[0]{$-2$}}
\put(144,23){\makebox(0,0)[0]{$-2$}}
\put(174,23){\makebox(0,0)[0]{$-2$}}
\put(114,73){\makebox(0,0)[0]{$-2$}}
\put(144,73){\makebox(0,0)[0]{$-2$}}
\put(174,73){\makebox(0,0)[0]{$-2$}}

\put(25,40){\line(1,0){150}}
\put(85,40){\line(6,5){30}}
\put(85,40){\line(6,-5){30}}
\put(115,15){\line(1,0){60}}
\put(115,65){\line(1,0){60}}

\put(20,65){\makebox(0,0)[0]{$\Gamma_{min}(X_{t=0})$}}
\end{picture}
\end{center}
This can be computed with cyclic covering techniques (cf. \cite{NCCover}). Then, by direct calculation or by Laufer's formula for smoothings of Gorenstein normal surface singularities (\cite{Lauferform}, see also Durfee's formula \cite{Durfeeform}), $Z_{K,min}^2=-12$.

In the $t\neq 0$ case we can consider the new variable $x_1':= x_2^3+tx_1$, so that in the holomorphic coordinates $(x_1', x_2, x_3)$ the singularity $X_{t \neq 0}$ is given by the equation: $\big(\frac{x_1'-x_2^3}{t}\big)^4+x_1'^3 + x_3^3=0$. This is once again Newton nondegenerate, with the Newton boundary consisting of a single facet with primitive normal vector $\mathbf{n}=(4,1,4)$. One can easily compute, that $p_g(X_{t \neq 0})=4$ with `conductor coefficient' $c=4$ and (combinatorial and doubled) weight function of the lattice homology:

 \begin{center}
     \begin{tabular}{c|ccccccccc}
          $\ell$ & $0$ & $1$ & $2$ & $3$ & $4$ & $5$ & $6$ & $7$ & $8$ \\
          \hline 
          $\mathfrak{h}^{\natural}(\ell)$ & $0$ & $1$ & $1$ & $2$ & $2$ & $3$ & $3$ & $4$ & $4$ \\
          $(\mathfrak{h}^{sym}_{c})^{\natural}(\ell)$ & $4$ & $4$ & $3$ & $3$ & $2$ & $2$ & $1$ & $1$ & $0$ \\
           \hline 
          $w(\ell)$ & $0$ & $1$ & $0$ & $1$ & $0$ & $1$ & $0$ & $1$ & $0$
     \end{tabular}
 \end{center}

 The minimal (good) resolution is the following (we stress that the link $L_{X_{t\neq 0}}$ is not a rational homology sphere with the central node having genus $1$), with $Z_{K, min}^2=-12$.

 
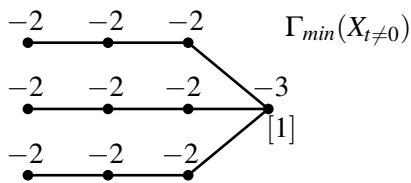
\begin{center}
\begin{picture}(150,70)(0,10)
\linethickness{1pt}

\put(85,40){\circle*{4}}
\put(55,40){\circle*{4}}
\put(25,40){\circle*{4}}
\put(115,40){\circle*{4}}
\put(25,15){\circle*{4}}
\put(55,15){\circle*{4}}
\put(85,15){\circle*{4}}
\put(25,65){\circle*{4}}
\put(55,65){\circle*{4}}
\put(85,65){\circle*{4}}

\put(24,48){\makebox(0,0)[0]{$-2$}}
\put(54,48){\makebox(0,0)[0]{$-2$}}
\put(83,48){\makebox(0,0)[0]{$-2$}}
\put(116,48){\makebox(0,0)[0]{$-3$}}
\put(24,23){\makebox(0,0)[0]{$-2$}}
\put(54,23){\makebox(0,0)[0]{$-2$}}
\put(82,23){\makebox(0,0)[0]{$-2$}}
\put(24,73){\makebox(0,0)[0]{$-2$}}
\put(54,73){\makebox(0,0)[0]{$-2$}}
\put(84,73){\makebox(0,0)[0]{$-2$}}
\put(120,32){\makebox(0,0)[0]{$[1]$}}

\put(25,40){\line(1,0){90}}
\put(85,15){\line(6,5){30}}
\put(85,65){\line(6,-5){30}}
\put(25,15){\line(1,0){60}}
\put(25,65){\line(1,0){60}}

\put(145,70){\makebox(0,0)[0]{$\Gamma_{min}(X_{t\neq0})$}}
\end{picture}
\end{center}
 
 In conclusion, along this deformation the geometric genus is constant with $p_g=4$, and  $Z_{K,min}^2$ is also constant with value $-12$. However, as one can compute, the analytic lattice homology is not 
 constant (see below), hence the conductor ideal is non-constant as well. 

 
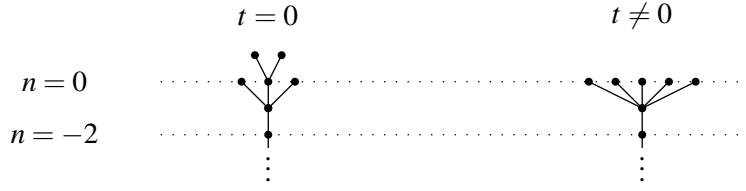
\begin{center}
\begin{picture}(350,70)(0,5)
\linethickness{.5pt}

\put(125,50){\circle*{3}}
\put(135,50){\circle*{3}}
\put(120,40){\circle*{3}}
\put(130,40){\circle*{3}}
\put(140,40){\circle*{3}}
\put(130,30){\circle*{3}}
\put(130,30){\circle*{3}}
\put(130,20){\circle*{3}}

\put(130,10){\makebox(0,0){$\vdots$}}

\qbezier[50](90,40)(200,40)(310,40)
\qbezier[50](90,20)(200,20)(310,20)

\put(50,40){\makebox(0,0)[0]{$n={0}$}}
\put(50,20){\makebox(0,0)[0]{$n=-2$}}

\put(130,40){\line(1,2){5}}
\put(130,40){\line(-1,2){5}}
\put(130,40){\line(0,-1){25}}
\put(120,40){\line(1,-1){10}}
\put(140,40){\line(-1,-1){10}}

\put(250,40){\circle*{3}}
\put(260,40){\circle*{3}}
\put(270,40){\circle*{3}}
\put(280,40){\circle*{3}}
\put(290,40){\circle*{3}}
\put(270,30){\circle*{3}}
\put(270,30){\circle*{3}}
\put(270,20){\circle*{3}}

\put(270,30){\line(0,1){10}}
\put(270,30){\line(1,1){10}}
\put(270,30){\line(2,1){20}}
\put(270,30){\line(-1,1){10}}
\put(270,30){\line(-2,1){20}}
\put(270,30){\line(0,-1){15}}

\put(270,10){\makebox(0,0){$\vdots$}}

\put(130,65){\makebox(0,0)[0]{$t=0$}}
\put(270,65){\makebox(0,0)[0]{$t \neq 0$}}
\end{picture}
\end{center}

\vspace{2mm}

Hence,  by part (c), the stability of the  conductor ideal and $p_g$ does not imply the stability of   $Z_{K,min}^2$, and, by 
 part (e), the stability of the $p_g$ and   $Z_{K,min}^2$
does not imply the stability of the conductor ideal. 
\end{example}

\begin{remark} On another note, it is appropriate to make the following observation. 
Consider a Gorenstein normal surface singularity $(X, o)$ and its minimal resolution $\phi:\widetilde{X}_{min}\to X$ and take the short 
exact sequence of sheaves $0\to \cO_{\widetilde{X}_{min}}(-mZ_{K,min})\to
 \cO_{\widetilde{X}_{min}}\to \cO_{mZ_{K,min}}\to 0$ for any $m\geq 1$.
Since $h^1( \cO_{\widetilde{X}_{min}}(-mZ_{K,min}))=0$ by the Grauert--Riemenschneider Vanishing Theorem (see e.g. \cite[Theorem 6.4.3]{NBook}), from the cohomological long exact sequence we deduce that 
\begin{align*}
    \dim\, \cO_{X,o}/(\phi_*(\cO_{\widetilde{X}_{min}}(-mZ_{K,min})))_o=&\, \dim H^0(\widetilde{X}_{min}, \mathcal{O}_{\widetilde{X}_{min}})\big/  H^0(\widetilde{X}_{min},\cO_{\widetilde{X}_{min}}(-mZ_{K,min})) \\
    =&\,\chi(mZ_{K,min})+p_g\\
    =&\,-Z_{K,min}^2\cdot m(m-1)/2+p_g.
\end{align*}
Now, if it happens that $(\dag)$ $(\phi_*(\cO_{\widetilde{X}_{min}}(-mZ_{K,min})))_o =
\big((\phi_*(\cO_{\widetilde{X}_{min}}(-Z_{K,min})))_0\big)^m$, then the above identity 
connects $\dim \cO_{X, o}/
\mathcal{C}^m$ (that is, the Hilbert function / polynomial  of $\mathcal{C}$) with 
$Z_{K,min}^2$. Usually $(\dag)$ is not valid (due to the existence of certain base points),
however still $\mathcal{C}$ and $Z_{K,min}^2$ (and $p_g$) are not totally independent. 
\end{remark}

\subsection{Further connections between deformations and lattice homology theories}\label{ss:DefsinLH}\,

In this subsection we wish to review some further facts and conjectures regarding the behaviour of both the analytic and topological
lattice (co)homology modules along some surface singularity deformations. 
We start this  discussion by recalling the construction of the natural morphism 
which connects the above two lattice homology theories.

Let us fix a normal surface singularity  $(X,o)$ whose link $L_X$ is a rational homology sphere. 
In particular, the topological lattice homology is well defined (cf. Example \ref{toplc}), whereas the analytic lattice homology can also be computed through the original definition of \cite{AgNe1} (see subsection \ref{ss:compar}).

Then for any good resolution $\phi:\widetilde{X} \rightarrow X$ and any effective cycle $\ell \in L=H_2(\widetilde{X}, \mathbb{Z})$,
from the exact sequence $0\to \cO_{\widetilde{X}}(-\ell)\to \cO_{\widetilde{X}}\to \cO_\ell\to 0$,
we obtain that $\hh(\ell)\leq h^0(\cO_{\ell})$, where $\mathfrak{h}$ is the analytic Hilbert function associated with the divisorial valuations of the resolution (cf. (\ref{eq:highdimHilbert}) or (\ref{eq:w_anhighdim})). In particular, 
$$w_{an}(\ell)=\hh(\ell)-h^1(\cO_\ell)\leq h^0(\cO_\ell)-h^1(\cO_\ell)=\chi(\ell)=w_{top}(\ell) \ \ \forall \ell \in L_{\geq 0}.$$
This shows that we automatically have an inclusion of the filtered spaces:
$S_{n, top}\subset S_{n, an}$
 for any $n\in\Z$. In particular, we have a bigraded $\Z[U]$-module morphism 
 $$\mathfrak{H}_*:\bH_{top,*}(X,o)= \bH_*((\R_{\geq 0})^r,w_{top})\to 
 \bH_*((\R_{\geq 0})^r, w_{an})=\bH_{an,*}(X,o).$$
By \cite[Example 11.9.29]{NBook} $\mathfrak{H}_*$ is an isomorphism if $(X,o)$ is contained in one of the
following families of singularities: 
rational, Gorenstein elliptic, weighted  homogeneous, superisolated associated with a unicuspidal irreducible rational curve, or splice quotient singularities with almost rational graphs.

If for a certain normal surface singularity $(X,o)$ the morphism $\mathfrak{H}_*$ is an isomorphism, then the Euler characterisitics of the two lattice homologies must agree, too. In particular, the geometric genus $p_g$ 
(which is, a priori, an analytic invariant of $(X,o)$) must agree with the normalized Seiberg--Witten invariant of the link $L_X$ (associated with the canonical spin$^c$-structure). This provides a topological description of $p_g$  (compare with the Seiberg-Witten Invariant Conjecture from \cite{NBook}).

On the other hand,  $\mathfrak{H}_*$ is not
always an isomorphism, e.g., for a superisolated singularity associated with a rational multicuspidal curve usually it is not. For example, if $(X,o)$ is the superisolated singularity 
associated with a degree 5 irreducible rational curve $C$ which has two cusps, each with a single Puiseux pair of type 
$(3,4)$ and $(2,7)$,  then $\bH_{top, 1}\not=0$ while $\bH_{an, 1}=0$, cf.
\cite[Example 11.9.49]{NBook}. 

\begin{corollary}\label{cor:topstab}
Consider a deformation of Gorenstein normal surface singularities, such that the link of 
each $(X_t,o)$ is a rational homology sphere and for each $(X_t,o)$ the morphism 
$\mathfrak{H}_*$ is an isomorphism. If this deformation is conductor ideal constant, then 
the topological lattice homology stays constant as well (even if the link is not constant). 
\end{corollary}
 
Parallel with the early results in the theory of topological lattice (co)homology, there appeared an expectation that along a Gorenstein $Z_{K,min}^2$-constant deformations this invariant remains stable, too. However, it turned out that this is not the case, counterexamples might appear already for the non-naive deformation of  superisolated germs associated with curves having more cusps. We see this fact also 
via Proposition \ref{prop:11.5.5}: whenever $\mathfrak{H}_*$ is not an isomorphism for a certain superisolated germ $(X, 0)$, 
$\bH_{top, *}$ cannot be stable along its non-naive deformation (indeed, for the central fiber $\mathfrak{H}_*$ is an isomorphism as $\bH_{top, *}(X_{t=0}, 0)\cong\mathbb{H}_{top, *}(\{x_1^d+x_2^d+x_3^{d+1}=0\}, 0)$ since $X_{t=0}$ is embedded topologically equivalent to $\{x_1^d+x_2^d+x_3^{d+1}=0\}$, and, since this is weighted homogeneous, 
$$\mathbb{H}_{top, *}(\{x_1^d+x_2^d+x_3^{d+1}=0\}, 0)\cong \mathbb{H}_{an, *}(\{x_1^d+x_2^d+x_3^{d+1}=0\}, 0) \cong \mathbb{H}_{an, *}(X_{t=0}, 0)$$
by Proposition \ref{prop:11.5.5}). We raise the attention, however, to the fact that in all known examples the non-stability appears in 
$\bH_{top, \geq 1}$ (see , e.g., \cite[Example 11.4.12]{NBook}). 
Therefore, a weaker version is still open (see \cite[Conjecture 11.3.39]{NBook}):
along a Gorenstein $Z_{K,min}^2$-constant deformation of surface singularities (with $\Q HS^3$ link) $\bH_{top, 0}$ is constant, too.

In \cite[Conjecture 11.9.48]{NBook} there is formulated another conjecture: $\bH_{an, *}$ is constant along a $p_g$-constant deformation of normal surface singularities with $\Q HS^3$ links.

Recall that in this work we extended the analytic lattice homology for any normal surface singularity, without any restriction about their links. Hence we might wonder about the extended version of this last conjecture: is $\bH_{an, *}$ constant along any  $p_g$-constant deformation?
The answer is negative, as Example \ref{ex:LIST} (e) shows, there are Gorenstein deformations with both  $Z_{K,min}^2$
and $p_g$ constant but $\bH_{an, 0}$ non-constant.   In the mentioned example $X_{t\not=0}$ has the embedded topological type of $\{x_1^3+x_2^3+x_3^{12}=0\}$, whose link is not a rational homology sphere, so this example does not contradict the aforementioned conjecture \cite[Conjecture 11.9.48]{NBook}. 
\newpage
\section{The proof of the Independence Theorem and a more general version}\label{s:proof}

In this section we present the proof of the main technical result of this paper: the Independence Theorem. The argument will involve the split weight function techniques pioneered by T. Ágoston and the first author, relying heavily on Combinatorial Duality Property and matroid rank inequality of the height functions. Our key combinatorial result is the Path Ordering Lemma, which replaces the case-by-case analysis used in the original proofs of well-definedness of analytic lattice homology theories and reduction-type theorems. To deal with the final homotopy theoretical question in our cubical setting we use Quillen's Fiber Theorem, or alternatively the technique of quasifibrations by Dold and Thom.  We will also present a more general formulation of the Independence Theorem and give an application for the equivariant analytic lattice homology of normal surface singularities.

\subsection{The Independence Theorem }\label{ss:dualtoN}\,

Let us recall the statement (we use the notations and definitions introduced in section \ref{s:4}):

\begin{theorem}[= Theorem \ref{th:IndepMod}] Let $k$ be a field, $\mathcal{O}$ a (Noetherian) $k$-algebra and $M$ a finitely generated module over it. Let $\mathcal{D}=\{\frv_1, \ldots, \frv_r\}$ and $\mathcal{D}'=\{\frv'_1, \ldots, \frv'_{r'}\}$ be two collections of extended discrete valuations (satisfying Assumption \ref{ass:fincodim}).  
        Suppose that the following conditions hold:
        \begin{itemize}
            \item $\mathcal{F}_{\mathcal{D}}^M(0) = \mathcal{F}_{\mathcal{D}'}^M(0)$ and
            \item both pairs $(\mathfrak{h}_{\mathcal{D}}, \mathfrak{h}_{\mathcal{D}}^\circ)$ and $(\mathfrak{h}_{\mathcal{D}'}, \mathfrak{h}_{\mathcal{D}'}^{\circ})$ satisfy the Combinatorial Duality Property.
        \end{itemize}
        Then the spaces $S_{n, \mathcal{D}}$ and $S_{n, \mathcal{D}'}$ 
        associated with the corresponding lattices and weight functions 
        are homotopy equivalent for every $n \in \mathbb{Z}$.
        Even more, the homotopy equivalences are compatible with the inclusions
        $S_{n, \mathcal{D}}\hookrightarrow S_{n+1, \mathcal{D}}$ and $S_{n, \mathcal{D}'}\hookrightarrow S_{n+1, \mathcal{D}'}\,$ such as to yield an isomorphism
        \begin{equation}
        \mathbb{H}_*(\mathbb{R}^r, w_{\mathcal{D}}) \cong  \mathbb{H}_*(\mathbb{R}^{r'}, w_{\mathcal{D}'}) \text{ of bigraded } \mathbb{Z}[U]\text{-modules.}
        \end{equation}
\end{theorem}

\begin{proof}[Proof of the module version of the Independence Theorem (Theorem \ref{th:IndepMod})]
The idea of the proof is, that if we get two finite collections of extended valuations ($\mathcal{D}:=\{\mathfrak{v}_1, \ldots, \mathfrak{v}_r\}$ and $ \mathcal{D}':=\{\mathfrak{v}'_1, \ldots, \mathfrak{v}'_{r'}\}$) with the required properties, we associate them one by one to each other, to get a common extension, and prove that along any single step the homotopy types of the $S_n$-spaces are preserved. This last part is the statement of the Adjoining Theorem, which will be discussed and proved in the succeeding subsection. Here, instead, we discuss the common extension strategy in detail. 

Let us introduce the notation $\mathcal{D}^{i, j}$ for the following collection of extended discrete va\-luations 
\begin{center}
    $\mathcal{D}^{i, j}=\{\frv_1, ..., \frv_i, \frv'_1, ..., \frv'_j\}$ (e.g., $\mathcal{D}^{r, 0}=\mathcal{D}$ and $\mathcal{D}^{0, r'}=\mathcal{D}'$).
\end{center} To any such collection we can associate the lattice $\Z^{i} \times \mathbb{Z}^{j} \subset \mathbb{Z}^{r}\times \mathbb{Z}^{r'}$ with basis $\{ e_1, ..., e_i, e'_1, ..., e'_j \}$, the multifiltrations $\cF_{\mathcal{D}^{i, j}}$ and $\cF^M_{\mathcal{D}^{i, j}}$ (cf. (\ref{eq:fdMl})),  the  height functions $\frh^{i, j}$  and $\frh^{\circ \ i, j }$(cf. (\ref{eq:handhcirc0})), the weight function $w^{i, j} = \frh^{i, j} + \frh^{\circ \ i, j } - \frh^{\circ \ i, j }(0) $ and for every $n \in \Z$ the cubical complex $S_{n}^{i, j}$. Let us denote the natural projections
\begin{align*}
    \text{by } \ \ \ \ & pr_1^{i, j}: \mathbb{Z}^{i} \times \mathbb{Z}^{j} = \Z \langle e_1, ..., e_i, e'_1, ..., e'_j \rangle \rightarrow  \mathbb{Z}^{i-1} \times \mathbb{Z}^{j} = \Z \langle e_1, ..., e_{i-1}, e'_1, ..., e'_j \rangle \\
    \text{and } \ \  & pr_2^{i, j}:  \mathbb{Z}^{i} \times \mathbb{Z}^{j} = \Z \langle e_1, ..., e_i, e'_1, ..., e'_j \rangle \rightarrow  \mathbb{Z}^{i} \times \mathbb{Z}^{j-1} = \Z \langle e_1, ..., e_i, e'_1, ..., e'_{j-1} \rangle.
\end{align*}

Then we consider the following sequence of maps:
$$
\resizebox{15cm}{!}{
\xymatrix{
S_{n, \mathcal{D}} = S_n^{r, 0}  &   S_n^{r, 1} \ar[l]_{\ \ \ \ \ \ pr_2} &   ...  \ar[l]_{\ \ pr_2} & S_n^{r, r'-1} \ar[l]_{pr_2} & S_n^{r, r'} \ar[l]_{\ \ pr_2} \ar[r]^{pr_1 \ \ } & S_{n}^{r-1, r'} \ar[r]^{\ \ pr_1} & ... \ar[r]^{pr_1 \ } & S_{n}^{1, r'} \ar[r]^{pr_1\ \ \ \ } & S_n^{0, r'}=S_{n, \mathcal{D}'}  \\
}
}
$$
Then the Adjoining Theorem \ref{thm:main} and Remark \ref{rm:F'isjo} from the following subsections ensure that the maps  $pr_1^{i, r'} \big|_{S_{n}^{i, r'}}: S_n^{i, r'} \rightarrow S_n^{i-1, r'} $ and $pr_2^{r, j} \big|_{S_n^{r, j}}: S_n^{r, j} \rightarrow S_n^{r, j-1}$ are all well-defined surjective homotopy equivalences $\big($and the pairs $(\frh^{i, r'}, \frh^{\circ \ i, r'})$ and $(\frh^{r, j}, \frh^{\circ \ r, j})$ all satisfy the Combinatorial Duality Property$\big)$. Clearly, these projections commute with the inclusions $S_n^{i,j} \hookrightarrow S_{n+1}^{i,j}$. 
\end{proof}

\begin{proof}[Proof of the symmetric version of the Independence Theorem (Theorem \ref{th:Indep})]
    By Remark \ref{rem:SH=HforCDP}, if in Theorem \ref{th:IndepMod} we specialize to $M:= \mathcal{O}$ and $N:= \mathcal{I}$ and make the translation to nonnegative extensions of the discrete valuations, then we get back exactly the well-definedness statement for the symmetric lattice homology of finite codimensional integrally closed ideals. In more details, the general strategy of the proof of Theorem \ref{th:IndepMod} in this setting can be presented (using the previous notations in addition with $d:=(d_1, \ldots, d_r)$, $d:'=(d_1', \ldots, d_{r'}')$ and $d^{i, j}=:(d_1, \ldots, d_i, d'_1, \ldots, d'_j)$) as the sequence of isomorphisms
$$
        \resizebox{15cm}{!}{
\xymatrix{
\mathbb{SH}_*(\mathcal{O}, \mathcal{D}, d) & \mathbb{SH}_*(\mathcal{O}, \mathcal{D}^{r, j}, d^{r, j}) \ar[l]^{\ldots}_{\cong} &  \mathbb{SH}_*(\mathcal{O}, \mathcal{D}^{r, r'}, d^{r, r'}) \ar[l]^{\ldots}_{\cong} \ar[r]^{\cong}_{\ldots} & \mathbb{SH}_*(\mathcal{O}, \mathcal{D}^{i, r'}, d^{i, r'}) \ar[r]^{\cong}_{\ldots} & \mathbb{SH}_*(\mathcal{O}, \mathcal{D}', d')  \\
}
}
$$
with $j$ increasing and $i$ decreasing. Moreover, in this case the projection maps  $pr_1^{i, r'} \big|_{S_{n}^{i, r'}}: S_n^{i, r'} \rightarrow S_n^{i-1, r'} $ and $pr_2^{r, j} \big|_{S_n^{r, j}}: S_n^{r, j} \rightarrow S_n^{r, j-1}$ are automatically $\mathbb{Z}_2$-equivariant (since $S_n^{i, j}$ is symmetric with respect to $d^{i, j}$, cf. Definition \ref{def:SH(O,D,d)}).
\end{proof}

\subsection{The Adjoining Theorem}\label{ss:mainthm}\,

Let $\cO$ denote a (Noetherian) $k$-algebra, $M$ a finitely generated $\mathcal{O}$-module and $N \leq M$ a realizable submodule (recall Definition \ref{def:REAL} --- $N$ will play the role of  $\mathcal{F}_{\mathcal{D}}^M(0) = \mathcal{F}_{\mathcal{D}'}^M(0) \leq M$ in the Independence Theorem \ref{th:IndepMod}).

Let $\frv_1, \frv_2, ..., \frv_r, \frv_{r+1}$ be extended discrete valuations on $\cO$ and $M$ in the sense of Definition \ref{def:edv}, translated such that $\frv^M_i(N) \geq 0\ $ for all $i\in \{1, \ldots, r+1\}$. We also suppose that they satisfy the conditions of Assumption \ref{ass:fincodim}. Define the sets $\mathcal{D}:=\{\frv_1, \frv_2, ..., \frv_r\} $ and $\mathcal{D}^+:=\{\frv_1, \frv_2, ..., \frv_r, \frv_{r+1}\} $ with corresponding index sets $\mathcal{V}:=\{1, \ldots, r \}$ and $\mathcal{V}^+:=\{1, \ldots, r, r+1 \}$. We consider the multifiltrations $\mathcal{F}_{\mathcal{D}}$ and $\mathcal{F}_{\mathcal{D}^+}$ (with rank $r$ and $r+1$ respectively) and their extensions  $\mathcal{F}_{\mathcal{D}}^M$ and $\mathcal{F}_{\mathcal{D}^+}^M$ as in (\ref{eq:fdMl}).  (For the sake of notational simplicity, we think of $\mathbb{Z}^r=\mathbb{Z}\langle \{e_v \}_{v \in \mathcal{V}} \rangle$  as a sublattice of $\mathbb{Z}^{r+1}=\mathbb{Z}\langle \{e_v \}_{v \in \mathcal{V}^+} \rangle$ of corank $1$.)

We also recall the notations $\frh$ and $\mathfrak{h}^{\circ}$ (respectively $\frh^+$ and $\mathfrak{h}^{\circ \, +}$) for the height functions (cf. (\ref{eq:handhcirc0})), $w$ (respectively  $w^+$) for the weight function and $S_n$ (respectively $S_n^+$) for the sublevel cubical complexes for every $n \in \Z$ assigned to the collection of extended discrete valuations $\mathcal{D}$ (respectively $\mathcal{D}^+$). 

\begin{theorem}[Adjoining Theorem] \label{thm:main}
Let $\mathcal{D}$ and $\mathcal{D}^+$ be collections of extended discrete valuations as introduced above. Suppose, moreover, that ${\mathcal{D}}$ satisfies the following properties:
\begin{itemize}
    \item $\cF_{\mathcal{D}}^M(0) =N$;
    \item the pair $(\frh, \frh^{\circ})$ satisfies the CDP, i.e., $\frh(\ell+e_v) - \frh(\ell)$ and $\frh^{\circ}(\ell+e_v) - \frh^{\circ}(\ell)$ cannot be simultaneously nonzero for every $\ell \in \mathbb{Z}^r$ and $v \in \mathcal{V}$.
\end{itemize}
Then the natural projection  
    $\pi_{\R}: \R^{r+1} \rightarrow \R^{r}, \ \sum_{v \in \mathcal{V}^+} {a_v}e_v \mapsto 
     \sum_{v \in \mathcal{V}} {a_v}e_v$ (along the last coordinate axis)
induces a well-defined surjective homotopy equivalence between $S^+_n$ and $S_n$ for every $n \in \Z$. Thus $\ \mathbb{H}_*(\mathbb{R}^r,w) \cong \mathbb{H}_*(\mathbb{R}^{r+1},w^+)$ as $\bZ[U]$-modules.
\end{theorem}

\begin{remark}\label{rem:D^+realizes}
    The identity $\cF_{\mathcal{D}}^M(0) = N$ already implies $ \cF_{\mathcal{D}^+}^M(0) =N$, since, by construction,
    \begin{center}$\cF_{\mathcal{D}^+}^M(0) = \cF_{\mathcal{D}}^M(0) \cap \mathcal{F}_{\mathfrak{v}_{r+1}}^M(0)$ with $\mathcal{F}_{\mathfrak{v}_{r+1}}^M(0) \supset N$.
    \end{center}
    We will also prove that, with this realization property assumed, the CDP of the pair $(\mathfrak{h}, \mathfrak{h}^\circ)$ implies the CDP for  $(\mathfrak{h}^+, \mathfrak{h}^{\circ\,+})$. 
\end{remark}

Before the proof of the Adjoining Theorem we will go through some necessary claims and lemmas. 

Let us fix an arbitrary $n \in \Z$. First, we would like to prove that the natural projection $\pi_{\R}$ 
induces a surjective map from $S^+_n$ to $S_n$. Indeed, at least on the level of lattice points, we have the following:

\begin{claim} \label{clm:b'_l}
  For any given lattice point $\ell \in \mathbb{Z}^r$, there exists an integer $0 \leq b_\ell < \infty$, satisfying  that $w^+(\ell + b_\ell e_{r+1})=w(\ell)$. Furthermore, it is also true, that:
\begin{itemize}
    \item for any $-\infty < b < b_\ell: \ w^+(\ell + be_{r+1}) \geq w^+(\ell + (b+1)e_{r+1})$;
    \item for any $b_\ell \leq b < \infty: \ w^+(\ell + be_{r+1})\leq w^+(\ell + (b+1)e_{r+1})$.
\end{itemize}
\end{claim}

\begin{proof}
From $\cF_{\mathcal{D}}^M(0) = \cF_{\mathcal{D}^+}^M(0) =N$ we have  $w(\ell) = \frh(\ell) + \frh^{\circ}(\ell) - \dim_k(M/N)$ whereas
\begin{center}
    $w^+(\ell+b_\ell e_{r+1}) = \frh^+(\ell +b_\ell e_{r+1}) + \frh^{\circ \, +}(\ell +b_\ell e_{r+1})  - \dim_k(M/N)$.
\end{center}
Since any $\mathfrak{h}$ is increasing and $\mathfrak{h}^\circ$ is decreasing, it is enough to prove that for every $\ell \in \Z^{r}$ there exist integers $0 \leq b_\ell$ and $b_\ell^\circ $, such that
\begin{itemize}
    \item for every $b \leq b_\ell: \ \frh^+(\ell + b e_{r+1}) = \frh(\ell)$;
    \item for every $b^\circ_\ell \leq b: \ \frh^{\circ \, +}(\ell + b e_{r+1}) = \frh^{\circ}(\ell)$;
    \item and $b^{\circ}_\ell \leq b_\ell$.
\end{itemize} 
The first assumption can easily be satisfied, as $\mathcal{F}_{\mathcal{D}^+}(\ell + b e_{r+1})=\mathcal{F}_{\mathcal{D}}(\ell) \cap \mathcal{F}_{\mathfrak{v}_{r+1}}(b)$.
We can choose $b_{\ell}$ to be the maximal such $b$, for which $\frh^+(\ell + b e_{r+1}) = \frh(\ell)$. Clearly, we have $0 \leq b_{\ell} < \infty$ from the properties of the discrete valuations (firstly, $\mathcal{F}_{\mathfrak{v}_{r+1}}(0) = \mathcal{O}$, whereas for the second inequality use the non-triviality condition (c) of Definition \ref{def:gdv} to obtain $\lim_{b \to \infty}\mathfrak{h}^+(\ell + be_{r+1})=\infty$). Also, by the maximality of $b_\ell$ we get that
\begin{center}
    there exists an $f \in \cF_{\mathcal{D}}(\ell) \cap \cF_{\frv_{r+1}}(b_\ell)$, such that $f \notin \cF_{\frv_{r+1}}(b_\ell+1) \ \Leftrightarrow \ \frv_{r+1}(f) = b_\ell$.
\end{center} 
Similarly, as $\mathcal{F}_{\mathcal{D}^+}^M(-\ell-b e_{r+1}) = \mathcal{F}_{\mathcal{D}}^M(-\ell)  \cap \mathcal{F}_{\mathfrak{v}_{r+1}}^M(-b) $, we can choose $b^{\circ}_\ell$ as the minimal such $b$, for which $\frh^{\circ \, +}(\ell + b e_{r+1}) = \frh^{\circ}(\ell)$ (once again, we have $b^\circ_{\ell} > - \infty$, since $\lim_{b \to -\infty}\mathfrak{h}^{\circ\, +}(\ell + be_{r+1})=\infty$ from the non-triviality condition (c) of Definition \ref{def:edv}). By its minimality, there exists an $m \in \mathcal{F}_{\mathcal{D}}^M(-\ell)$ with $\mathfrak{v}_{r+1}^M(m)=-b^{\circ}_\ell$.
 But then 
 \begin{center}$f\cdot m \in \cF_{\mathcal{D}}(\ell) \cdot \cF_{\mathcal{D}}^M(-\ell) \subset \cF_{\mathcal{D}}^M(0) = N$, so $\frv_{r+1}^M(f\cdot m) \geq 0 $
 \end{center}
 by the property $\frv_{r+1}^M(N) \geq 0$. However $\frv_{r+1}^M(f\cdot m) = \frv_{r+1}(f) + \frv_{r+1}^M(m)=b_{\ell} -b^{\circ}_{\ell}$.
\end{proof}

 \begin{remark}\label{rem:goingtoinfty}
    Notice that for this previous Claim \ref{clm:b'_l} we do not actually need the non-triviality conditions (c) of Definitions \ref{def:gdv} and \ref{def:edv} (these were, in fact, required for the finiteness property of the weight function, cf. the proof of Lemma \ref{lem:w0matroid}). Indeed, it is even easier to find a desired $b_\ell$ with the stated properties if $\frh^+(\ell + b e_{r+1}) = \frh(\ell)$ or $\frh^{\circ \, +}(\ell + b e_{r+1}) = \frh^{\circ}(\ell)$ for all $b \in \mathbb{Z}$. This observation is important when shifting from discrete valuations to the quasi-valuations of subsection \ref{ss:quasi}.
\end{remark}

\begin{cor}
    The projection $\pi_{\mathbb{R}}$ maps  $S_n^+$ into $S_n$. Indeed, since a real point $p^+ \in S_n^+$ if and only if the vertices of the unique open cube containing it have $w^+$-weight  $\leq n$, whereas, by Claim \ref{clm:b'_l}, for any $\ell^+ \in \mathbb{Z}^{r+1}$ we have $w^+(\ell^+) \geq w(\pi_\mathbb{R}(\ell^+))$, hence $w(\pi_{\mathbb{R}}(p^+))\leq w^+(p^+) \leq n$ and $\pi_\mathbb{R}(p^+) \in S_n$. 
\end{cor}

Now, in order to prove that $\pi_{\R}$ induces an ${S^+_n} \rightarrow S_n$ surjection, we would like to show that for every $q$-cube $\square_q= (\ell, I)$ of $S_n$ the $q$-cube $(\ell + b_\ell e_{r+1}, I)$ is in $S^+_n$ (here we use the notations of paragraph \ref{bek:211}). For this, it is enough to prove that
\begin{center}
if $w \big|_{\{ \ell+ e_J\,:\, J \subset I \} } \leq n \Rightarrow w^+\big|_{\{ \ell + b_\ell e_{r+1} + e_J\,:\, J \subset I \}} \leq n$.
\end{center} 
This will be done via the Path Ordering Lemma.

\subsection{The Path Ordering Lemma}\,

We will state and prove this lemma slightly more generally than needed for the Adjoining Theorem. First we need to introduce some convenient notations.

Fix some lattice $\Z^r$ with basis $\{e_v\}_{v\in \mathcal{V}}$ and weight function $w: \mathbb{Z}^r \rightarrow \mathbb{Z} $ on it.
Let $\ell$  and $\ell'$ be two lattice points $\ell, \ell' \in \mathbb{Z}^r$, with $\ell' \geq \ell$. Following part (iii) of paragraph \ref{9SSP} we call a finite sequence $\gamma=\{ x_i \}_{i=0}^{t}$ of lattice points an \textit{`increasing path'} from $\ell$ to $\ell'$, if $x_0=\ell$,  $x_t = \ell'$ and $x_{i+1} = x_i + e_{v(i)}$ for some index $v(i) \in \mathcal{V}$ for all $0 \leq i <t$. Denote by $s_+(\gamma):= \sum_{i=0}^{t-1} \max \{ 0, w(x_{i+1})-w(x_i) \}$ and $s_-(\gamma):= \sum_{i=0}^{t-1} \min \{ 0, w(x_{i+1})-w(x_i) \}$. Obviously $w(\ell')-w(\ell) = s_+(\gamma) + s_-(\gamma)$.

For ease of reference, let us call an  edge $(x, \{e_v\})$ (also denoted, for simplicity, as $(x, x+e_v)$)  from $x$ to $x+e_v$ `positive, `neutral' or `negative' (and mark them with $+, 0$ or $-$ accordingly) depending on the sign of the difference $(w(x+e_v)-w(x))$. 

\begin{lemma}[Path Ordering Lemma]\label{lem:POL}
Suppose that the weight function $w$ satisfies the matroid rank inequality, i.e., 
\begin{equation*}
    w(x+\bar{x}+e_v)-w(x+\bar{x}) \leq w(x+e_v)-w(x) \hspace{5mm} \forall x \in \mathbb{Z}^r, \bar{x} \in (\mathbb{Z}_{\geq 0})^r, v \in \mathcal{V}, e_v \notin {\rm supp}(\bar{x}).
\end{equation*}
Then, for any two lattice points $\ell, \ell' \in \mathbb{Z}^r$, with $\ell' \geq \ell$ being a vertex of the $r$-cube $(\ell,\mathcal{V})$, there exists an increasing path $\widetilde{\gamma}$ from $\ell$ to $\ell'$, such that in $\widetilde{\gamma}$ all the positive edges precede the neutral and negative ones, and all the negative ones follow the positive and neutral ones (i.e., the order of the edges is: positive, neutral, negative), with $s_+(\widetilde{\gamma})$ maximal and $s_-(\widetilde{\gamma})$ minimal (i.e., for any other increasing path $\gamma$ from $\ell$ to $\ell'$ we have  $s_+(\widetilde{\gamma}) \geq s_+(\gamma)$ and $s_-(\widetilde{\gamma}) \leq s_-(\gamma)$).
\end{lemma}

\begin{proof}
It is enough to prove the following: given an increasing path $\gamma$ from $\ell$ to $\ell'$, we can construct another one, $\gamma^{ord}$, having the prescribed edge order and also $s_+(\gamma^{ord}) \geq s_+(\gamma)$ and $s_-(\gamma^{ord})\leq s_-(\gamma)$.

The main observation following from the matroid rank inequality for $w$ is that
\begin{align}
    \text{if the edge} (x, x+e_v) \text{ is positive,} &\nonumber \\
    \text{then } \forall w \neq v: \ (x- e_w, &\ x-e_w+e_v) \text{ is still positive, with} \label{eq:CDP+SP}\\
    w(x-e_w+&\,e_v)-w(x-e_w) \geq w(x+e_v)-w(x). \nonumber
\end{align}
The opposite statement is true for the negative edges as well. 

Now the idea of the proof is to make simple switches to the increasing path until we arrive to a path satisfying the edge order conditions. A simple switch is a single change of lattice points in the path:
\begin{center}
\begin{tabular}{cccccccc}
    from & $\color{red}\g \color{black}$ & $=$ & $\ldots$ & $x_i$ & $x_{i+1}= x_{i}+E_{v(i)}$ & $x_{i+2}=x_{i}+E_{v(i)} + E_{v(i+1)}$ & $\ldots$ \\
    to & $\color{blue}\widehat{\g}\color{black}$ & $=$ & $\ldots$ & $\widehat{x_i}=x_i$ & $\widehat{x_{i+1}}= x_{i} + E_{v(i+1)}$ & $\widehat{x_{i+2}}=x_{i}+ E_{v(i+1)} + E_{v(i)}= x_{i+2}$ & $\ldots$
\end{tabular}
\end{center}
(Notice that, since $\ell'$ is a vertex of $(\ell, \mathcal{V})$, we have the inequality $v(i) \neq v(i+1)$.)\\
We will make this simple switch in the following three cases:
\begin{enumerate}
    \item if the first edge is negative and the second is positive, then, after the switch, the first edge will be positive and the second negative;
    \item if the first edge is neutral and the second is positive, then, after the switch, the first edge will be positive and the second might be neutral or negative;
    \item if the first edge is negative and the second is neutral, then, after the switch, the first edge might be positive or neutral and the second will be negative.
\end{enumerate}
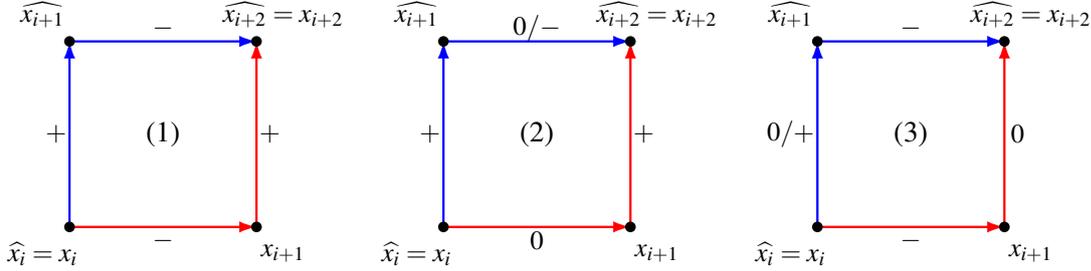
\begin{picture}(420,110)(-5,0)

\multiput(0, 0)(140, 0){3}{\makebox(40, 20){\small{$\widehat{x_i} = x_i$}}}
\multiput(90, 0)(140, 0){3}{\makebox(40, 20){\small{$x_{i+1}$}}}
\multiput(0, 90)(140, 0){3}{\makebox(40, 20){\small{$\widehat{x_{i+1}}$}}}
\multiput(90, 90)(140, 0){3}{\makebox(40, 20){\small{$\widehat{x_{i+2}}=x_{i+2}$}}}

\put(60, 50){\makebox(10, 10){(1)}}
\put(200, 50){\makebox(10, 10){(2)}}
\put(340, 50){\makebox(10, 10){(3)}}

\put(60, 10){\makebox(10, 10){\small{$-$}}}
\put(100, 50){\makebox(10, 10){\small{$+$}}}
\put(20, 50){\makebox(10, 10){\small{$+$}}}
\put(60, 90){\makebox(10, 10){\small{$-$}}}

\put(200, 10){\makebox(10, 10){\small{$0$}}}
\put(240, 50){\makebox(10, 10){\small{$+$}}}
\put(160, 50){\makebox(10, 10){\small{$+$}}}
\put(200, 90){\makebox(10, 10){\small{$0/-$}}}

\put(340, 10){\makebox(10, 10){\small{$-$}}}
\put(380, 50){\makebox(10, 10){\small{$0$}}}
\put(295, 50){\makebox(10, 10){\small{$0/+$}}}
\put(340, 90){\makebox(10, 10){\small{$-$}}}

\linethickness{0.3mm}

\multiput(30,20)(140, 0){3}{\color{red}\vector(1, 0){70}}
\multiput(100,20)(140, 0){3}{\color{red}\vector(0, 1){70}}
\linethickness{0.3mm}
\multiput(30,90)(140, 0){3}{\color{blue}\vector(1, 0){70}}
\multiput(30,20)(140, 0){3}{\color{blue}\vector(0, 1){70}}
\multiput(100,90)(140, 0){3}{\circle*{4}}
\multiput(30,20)(140, 0){3}{\circle*{4}}
\multiput(100,20)(140, 0){3}{\circle*{4}}
\multiput(30,90)(140, 0){3}{\circle*{4}}
\end{picture}

We can also see from (\ref{eq:CDP+SP}) that through these three types of changes, the $s_+$ function cannot decrease, and $s_-$ cannot increase.

Due to the finiteness of the path, the following easy algorithm gives us an increasing path that satisfies the prescribed edge order conditions: first, move all the existing positive edges to the front (through simple switches of type (1) and (2)), then move every negative edge to the back (through simple switches of type (1) and (3)), then, move the newly created positive edges in front of the neutrals, then move the newly created negative edges behind the neutrals etc. repeat this process until you get to the right order.
\end{proof}

\begin{remark}
    The statement of the lemma remains valid if the matroid rank inequality is only satisfied in a (finite) rectangle and we only use lattice points contained in that. On the other hand, the conclusion is no longer true if $\ell'$ is not a vertex of the cube $(\ell, \mathcal{V})$, e.g., if $\ell'=\ell + 2e_v$ there exists only a single increasing path, hence no room for switches.
\end{remark}

\begin{example}\label{ex:h&hcircmatroidimpliesROLemma}
    Let us consider the combinatorial setup and notations of paragraph \ref{bek:comblattice}. Suppose that both $h$ and $h^{\circ}$  satisfy the matroid rank inequality and the pair $(h, h^{\circ})$ satisfies the Combinatorial Duality Property. Then the weight function $w=h+h^{\circ}-h^{\circ}(0)$ satisfies the matroid rank inequality, and, hence, the Path Ordering Lemma holds true for it. In fact, besides the CDP for the pair $(h, h^{\circ})$, it is enough to require the matroid rank inequality for $h$. Indeed, that already prescribes the shape of switches (1) and (2), while in the case of switch (3), although the first new edge can become negative, it will not ruin the order and $s_+({\color{red}\gamma}\big|_{[x_i, x_{i+2}]})=0$ and $s_-({\color{red}\gamma}\big|_{[x_i, x_{i+2}]})=w(x_{i+2})-w(x_i)$ still cannot change in the wrong direction.
\end{example}

Going back to the proof of the Adjoining Theorem, observe that, by Example \ref{ex:h&hcircmatroidimpliesROLemma}, in order to apply the Path Ordering Lemma for $w=\mathfrak{h} + \mathfrak{h}^\circ-\dim_{k}(M/N)$,  all we need is the fact that the height functions $\ell \mapsto \dim_k \mathcal{O}/\mathcal{F}_{\mathcal{D}}(\ell), \ \ell \mapsto \dim_k M/\mathcal{F}_{\mathcal{D}}^M(\ell)$ satisfy the matroid rank inequality (cf. Observation \ref{obs:hhMmatroid}) and the Combinatorial Duality Property. Indeed, the involution $\ell \mapsto -\ell$ used in the definition of $\mathfrak{h}^{\circ}$ keeps the form of inequality (\ref{eq:matroid}). 

Notice, however, that $\frh^+$ and $\frh^{\circ \, +}$ also have these properties: they are Hilbert functions of certain valuative multifiltrations, so they satisfy the matroid rank inequality, thus we only have to prove the CDP. This will follow from the subsequent remark.

\begin{remark} \label{rm:h'CDP}
Let us suppose that for some $v \in \cV$ and lattice point $\ell \in \mathbb{Z}^r$ we have $\frh(\ell)=\frh(\ell+e_v)$ (which is equivalent to $\cF_{\mathcal{D}}(\ell)  =  \cF_{\mathcal{D}}(\ell+E_v) $, 
since $\cF_{\mathcal{D}}$ is a descending multifiltration). But then, for any integer $c\in \mathbb{Z}$:
\begin{center}
    $\cF_{\mathcal{D}}(\ell) \cap \cF_{\frv_{r+1}}(c)=\cF_{\mathcal{D}}(\ell+e_v) \cap \cF_{\frv_{r+1}}(c) \Rightarrow \frh^+(\ell+ce_{r+1}) = \frh^+(\ell+ce_{r+1} + e_v)$.
\end{center}
Similarly, $\frh^{\circ}(\ell) = \frh^{\circ}(\ell +e_v)$ implies $\frh^{\circ\, +}(\ell+ce_{r+1}) = \frh^{\circ\, +}(\ell+ce_{r+1} + e_v)$ for any $c \in \mathbb{Z}$ (and $v \in \cV$).
Hence from the CDP of the pair $(\frh, \frh^{\circ})$ we get the CDP  for $(\frh^+, \frh^{\circ\, +})$ for any $v \in \cV$ and $\ell^+ \in \mathbb{Z}^{r+1}$. The case of $v=r+1$ results from the proof of Claim \ref{clm:b'_l}. Hence the pair $(\frh^+, \frh^{\circ \, +})$ satisfies the Combinatorial Duality Property, so the Path Ordering Lemma can be used for $w^+$ as well.
\end{remark}

\begin{proposition}\label{prop:pisurj}
The natural projection $\pi_{\R}: \R^{r+1} \rightarrow \R^{r}, \ \sum_{v \in \mathcal{V}^+} {a_v}e_v \mapsto 
     \sum_{v \in \mathcal{V}} {a_v}e_v$ induces, for every $n \in \Z$, an $S^+_n \rightarrow S_n$ surjection.
\end{proposition}

\begin{proof}
We will prove that for every $q$-cube $\square_q= (\ell, I)$  of $S_n$ the $q$-cube $\square^+_q=(\ell + b_\ell e_{r+1}, I)$ is in $S^+_n$ (where $b_{\ell}$ is the integer given by Claim \ref{clm:b'_l}). Let us denote $\ell + b_{\ell}e_{r+1}$ by $\ell^+$. Then $w(\ell) = w^+(\ell^+)$ by Claim \ref{clm:b'_l}. Therefore, if $q=0$, we already have the lift. Suppose then that $q > 0$.

Recall that, by their definitions, for any subset $J \subset I$:
\begin{center}$w(\ell+e_J) = w(\ell)+s_+(\g)+s_-(\g)$ for any increasing path $\g$ from $\ell$ to $\ell+e_J$. 
\end{center}
The same is true for the weight function $w^+$ and any increasing path from $\ell^+$ to $\ell^++e_J$.

Recall also that $w^+(\square_q^+)$ is the maximal  $w^+$ weight obtained by the vertices of $\square_q^+$. Then there exists a vertex $\ell'^+ \in (\ell^+, I)$  with $w^+(\ell'^+)=w^+(\square_q^+)$ together with an increasing path $\g^+$ from $\ell^+$ to $\ell'^+$ having only positive edges. Indeed, if we take any $w^+$-maximal lattice point in $(\ell^+, I)$, starting from any increasing path the Path Ordering Lemma \ref{lem:POL} provides a well ordered one, which cannot have negative edges, and after eliminating the neutral edges we still get a $w^+$-maximal lattice point. Thus $w^+(\square_q^+)= w(\ell) + s_+(\g^+)$. As $I \subset \cV$, we can project the path $\g^+$ to $(\ell, I)$, it will give us an increasing path $\g := \pi_{\R}(\g^+)$ from $\ell$ to $\pi_{\R}(\ell'^+) =: \ell'$. We would like to prove, that $\g$ has  only positive edges and $s_+(\g) \geq s_+(\g^+)$.

However, this is just a corollary of the Combinatorial Duality Property  and the matroid rank inequality for the Hilbert function $\mathfrak{h}^+$. Indeed, for any positive edge  $(\ell^++e_J, \ell^++e_J +e_v)$ contained in $(\ell^+, I)$ (with $J \sqcup \{v\} \subset I$) we have
\begin{align*}
    0<w^+(\ell^++e_J & +e_v) - w^+(\ell^+ +e_J) = \\ 
    =\ & \frh^+(\ell^++e_J +e_v) - \frh^+(\ell^+ +e_J)  \hspace{5mm} \text{by edge positivity and CDP}; \\
    = \ & \frh^+(\ell+e_J +e_v + b_{\ell} e_{r+1}) - \frh^+(\ell +e_J+ b_{\ell}e_{r+1}) \\
    \leq \ & \frh^+(\ell+e_J +e_v) - \frh^+(\ell +e_J) \hspace{10mm} \text{by matroid rank inequality;} \\
    = \ & \frh(\ell+e_J +e_v) - \frh(\ell +e_J) \hspace{15mm} \text{by their definitions;}\\
    = \ & w(\ell+e_J +e_v) - w(\ell +e_J) \hspace{14mm} \text{by CDP for } \frh \text{ and }\frh^{\circ}.
\end{align*}
So, clearly, $\gamma$ also has only positive edges and summing over the differences we get the desired inequality $s_+(\g) \geq s_+(\g^+)$ as well. Now
\begin{center}
    $w^+(\square_q^+)= w^+(\ell'^+) = w^+(\ell^+) + s_+(\g^+) \leq w(\ell) + s_+(\g) = w(\ell') \leq w(\square_q)$.
\end{center} 
Hence, if $\square_q \subset S_n \Leftrightarrow w(\square_q) \leq n$, then $w^+(\square_q^+) \leq n \Leftrightarrow \square_q^+ \subset S_n^+$, so we produced a lift for every cube in the cubical decomposition of $S_n$.
\end{proof}

\begin{cor} \label{cor:invimofpts}
The preimage of every point $p \in S_n$ with respect to the projection $\pi_{\R} \big|_{S^+_n}: S^+_n \rightarrow S_n$ is a nonempty closed (maybe degenerate) vertical interval. Even more, for any $q$-cube  $\square_q$  from the cubical decomposition of $S_n$, there exists a nonempty closed (maybe degenerate) vertical interval $\I \subset \mathbb{R}\langle e_{r+1}\rangle$ (with integral endpoints), such that
\begin{center}
$\left( \pi_{\R} \big|_{S^+_{n}} \right)^{-1} \left(\square_q^\circ \right) = \square_q^\circ \times \I$,
\end{center}
where $\square_q^\circ$ denotes the interior of $\square_q$.
\end{cor}
\begin{proof}
For the lattice points $\ell \in S_n$ this results from Claim \ref{clm:b'_l}. Any other point is an interior point of some unique cube $\square_q=(\ell, I)$. Proposition \ref{prop:pisurj} (or, more precisely, its proof) implies that 
\begin{equation}\label{eq:osszelet}
    \square_q^\circ \times \{ b_{\ell}\} \subset\left( \pi_{\R} \big|_{S^+_{n}} \right)^{-1} \left(\square_q^\circ \right).
\end{equation}
for the integer $b_{\ell} \geq 0$.
If the same is true for another integer $b>b_{\ell}$ (respectively $b<b_{\ell}$), then for every vertex $\ell'$ of $\square_q$ both $w^+(\ell'+b_{\ell} e_{r+1}) \leq n$ and $w^+(\ell'+b e_{r+1}) \leq n$. This, in turn, implies, by Claim \ref{clm:b'_l}, that for every intermediate integer $b_{\ell} < c < b$ (respectively $b_{\ell} > c > b)$ and every vertex $\ell'$ of the cube $\square_q$ we have $ w^+(\ell'+c e_{r+1}) \leq n$. Thus the lift of $ \square_q$ to any of these heights $c$ is in $S^+_n$ and so are the $(q+1)$-cubes connecting them. So for any such $b$ we have 
$\square_q \times [ b_{\ell}, b] \subset S^+_n $.
The union of these intervals gives $\I$. \\
\big[We remark, that if the containment (\ref{eq:osszelet}) is true for some non-integer $b$, then, by formula (\ref{eq:9weight}), we already have $\square_q^\circ \times [\lfloor b \rfloor, \lceil b \rceil] \subset\big( \pi_{\R} \big|_{S^+_{n}} \big)^{-1} \left(\square_q^\circ \right).$\big]
\end{proof}

\begin{cor}\label{cor:preimagecontract}
    The inverse image of every closed $q$-cube $\square_q$ of $S_n$ with respect to the projection $\pi_{\R} \big|_{S^+_n}: S^+_n \rightarrow S_n$ is contractible. Indeed, if $\square_q=(\ell, I)$, then $\left( \pi_{\R} \big|_{S^+_{n}} \right)^{-1} \left(\square_q \right)$ deformation retracts along the $e_{r+1}$ axis to its lift $(\ell+b_{\ell}e_{r+1}, I)$. 
\end{cor}

We are now ready to prove the Adjoining Theorem:

\begin{proof}[Proof of Theorem \ref{thm:main}.]
We would like to prove that the natural projection $\pi_{\R}$ restricted to $S^+_n$ is a homotopy equivalence. Notice that $\pi_{\R} \big|_{S^+_n}: S^+_n \rightarrow S_n$ can be thought of as a simplicial map between simplicial complexes, and, by Corollary \ref{cor:preimagecontract}, each simplex has contractible inverse image. Quillen's Theorem A (also known as Quillen's Fiber Theorem) asserts that in this case the simplicial map is a homotopy equivalence (c.f. \cite{Quillen}). 

We can also give a not much longer proof using the theory of quasifibrations.
By Whitehead's Theorem, it is enough to prove, that $\pi_{\R} \big|_{S^+_n}: S^+_n \rightarrow S_n$ is a weak homotopy equivalence. As we already know that the fibers are contractible, it suffices to prove that the map $\pi_{\R} \big|_{S^+_n}: S^+_n \rightarrow S_n$ is a quasifibration, then the homotopical long exact sequence will yield the result. (For the definition and basic properties of quasifibrations see \cite{DoldThom})

\cite[Theorem 6.1.5]{DadNem} provides us simple criteria to check that this projection is indeed a quasifibration:
\begin{enumerate}[label=(\alph*)]
    \item $S_n$ is locally contractible as a closed cubical complex.
    \item The linear homotopies from the contractible neighbourhoods can be lifted to $S^+_n$ 
    as we know from Corollary \ref{cor:invimofpts} how the preimages of the open $q$-cubes look like.
    \item The retractions given by the lifted homotopies induce weak homotopy equivalences on the fibers, as they are all contractible.
    \item For all $x \in S_n,\ \pi_{\R}^{-1}(x) \cap S^+_n$ is path connected, as it is a closed interval.
\end{enumerate}
Thus, $\pi_{\R} \big|_{S^+_n} S^+_n \rightarrow S_n$ is a quasifibration with contractible fibers, so it turns out to be a homotopy equivalence for every $n \in \Z$. 
\end{proof}

\begin{remark} \label{rm:F'isjo}
Notice that Remarks \ref{rem:D^+realizes} and \ref{rm:h'CDP} imply that the multifiltration $\cF_{\mathcal{D}^+}$ satisfies the conditions of the Adjoining Theorem \ref{thm:main}, so if we take another extended discrete valuation $\frv_{r+2}$ on $\cO$ and $M$, such that $\frv_{r+2}^M(N) \geq 0$, then the height functions corresponding to the new multifiltration associated to these $r+2$ valuations gives the same lattice homology. So we can indeed run the argument proving the Independence Theorem.
\end{remark}

\subsection{Filtered homotopy type}\label{ss:filteredhom}

Notice, that we do not state the results neither in the Independence Theorem \ref{th:IndepMod} nor in the Adjoining Theorem \ref{thm:main} in terms of the filtered homotopy type of 
$$\ldots \subset S_{n-1} \subset S_{n} \subset S_{n+1} \subset \ldots \subset \mathbb{R}^n$$
or the chain homotopy type of $(\mathcal{L}_*, \partial_{w,*})$ (in the sense, e.g., of \cite{Seppo}).Indeed, we only have the following weaker equivalence relation on the level of the $\{S_n\}_n$ filtration (which only induces isomorphisms on the lattice homology modules):

\begin{theorem}[Topological version of the Independence Theorem]\label{th:IndepModTop}
    Let $k$ be a field, $\mathcal{O}$ a (Noetherian) $k$-algebra and $M$ a finitely generated module over it. Let $\mathcal{D}=\{\frv_1, \ldots, \frv_r\}$ and $\mathcal{D}'=\{\frv'_1, \ldots, \frv'_{r'}\}$ be two collections of extended discrete valuations (satisfying Assumption \ref{ass:fincodim}) and suppose that the following conditions hold:
        \begin{itemize}
            \item $\mathcal{F}_{\mathcal{D}}^M(0) = \mathcal{F}_{\mathcal{D}'}^M(0)$ and
            \item both pairs $(\mathfrak{h}_{\mathcal{D}}, \mathfrak{h}_{\mathcal{D}}^\circ)$ and $(\mathfrak{h}_{\mathcal{D}'}, \mathfrak{h}_{\mathcal{D}'}^{\circ})$ satisfy the Combinatorial Duality Property.
        \end{itemize}
        Then for every $n \in \mathbb{Z}$ there exist homotopy equivalences $f_n:S_{n, \mathcal{D}} \rightarrow  S_{n, \mathcal{D}'}$ (with homotopy inverses $g_n:S_{n, \mathcal{D}'} \rightarrow  S_{n, \mathcal{D}}$) between the spaces $S_{n, \mathcal{D}}$ and $S_{n, \mathcal{D}'}$ 
        commuting up to homotopy with the natural inclusions $\iota_{n, \mathcal{D}}:S_{n, \mathcal{D}} \hookrightarrow S_{n+1, \mathcal{D}}$ and $\iota_{n, \mathcal{D}'}:S_{n, \mathcal{D}'} \hookrightarrow S_{n+1, \mathcal{D}'}$. More explicitely, these maps satisfy the following homotopy relations:
        \begin{itemize}
            \item $g_{n} \circ f_n \simeq {\rm id}_{S_{n, \mathcal{D}}}$ and $f_{n} \circ g_n \simeq {\rm id}_{S_{n, \mathcal{D}'}}$ \ for all $n \in \mathbb{Z}$;
            \item $f_{n+1} \circ \iota_{n, \mathcal{D}} \simeq \iota_{n, \mathcal{D}'} \circ f_{n}$ \ for all $n \in \mathbb{Z}$;
            \item $g_{n+1} \circ \iota_{n, \mathcal{D}'} \simeq \iota_{n, \mathcal{D}} \circ g_{n}$ \ for all $n \in \mathbb{Z}$.
        \end{itemize}
\end{theorem}

\begin{proof}
    We check that the homotopy equivalences constructed (implicitly) in the original proof of the Independence Theorem \ref{th:IndepMod} (given as the composition of homotopy equivalences $S_n^{r, j} \longleftrightarrow S_{n}^{r, j+1}$ and $S_n^{i, r'} \longleftrightarrow S_{n}^{i-1, r'}$ produced by the Adjoining Theorem \ref{thm:main}) satisfy the required commutativity properties. In fact, due to the step-by-step nature of that construction, it is enough to prove the corresponding topological enhancement of the Adjoining Theorem: let $\pi_n:S_n^+ \rightarrow S_n$ abridge the homotopy equivalence $\pi_{\mathbb{R}}\big\vert_{S_n^+}$, let $s_n:S_n \rightarrow S_{n}^+$ denote its homotopy inverse (existing due to Quillen's Fiber Theorem) and let $\iota_n:S_n \hookrightarrow S_{n+1}$ and $\iota_n^+:S_n^+ \hookrightarrow S_{n+1}^+$ denote the natural inclusions; then we have to prove that
    \begin{align}
        & \pi_{n+1}\circ \iota_n^+\simeq \iota_n \circ \pi_n: S_n^+ \rightarrow S_{n+1} & \text{for all }n \in \mathbb{Z}, \label{eq:commuting}\\
        & s_{n+1} \circ \iota_n \simeq \iota_{n}^+ \circ s_n:S_n \rightarrow S_{n+1}^+ & \text{for all }n \in \mathbb{Z}. \label{eq:commutinguptohom}
    \end{align}
    In fact, identity (\ref{eq:commuting}) is already an equality on the set-map level, due to the explicit description of the $\{\pi_n\}_n$ projection maps. For identity (\ref{eq:commutinguptohom}) we use the fact the $s_n$ is the homotopy inverse of $\pi_n$ for every $n \in \mathbb{Z}$. This implies that
    \begin{equation*}
        s_{n+1} \circ \iota_n =s_{n+1} \circ \iota_n \circ {\rm id}_{S_n} \simeq s_{n+1} \circ \iota_n \circ \pi_n \circ s_n =s_{n+1} \circ \pi_{n+1} \circ \iota_n^+ \circ s_n \simeq {\rm id}_{S_{n+1}^+}\circ \iota_n^+ \circ s_n = \iota_n^+ \circ s_n. 
    \end{equation*}
    [These homotopy equivalences induce the isomorphism of the lattice homology $\mathbb{Z}[U]$-modules.]
\end{proof}

\begin{remark}
    Using the same techniques, one can also prove that in the symmetric setting of Theorem \ref{th:Indep} the homotopy equivalences $\pi_n:S_n^+ \rightarrow S_n$ and $s_n:S_n \rightarrow S_n^+$ commute with the underlying $\mathbb{Z}_2$-symmetry up to homotopy.
\end{remark}

We firmly believe that using more advanced techniques (e.g., homotopy colimits such as in \cite{Seppo}) one could construct the homotopy inverses $\{s_n\}_{n\in \mathbb{Z}}$ such that they would genuinely commute with the inclusion maps, i.e., to have set-map equalities in identity (\ref{eq:commutinguptohom}) (and also in the last two rows of Theorem \ref{th:IndepModTop}). That would mean, that already the filtered homotopy type of $\ldots \subset S_{n-1} \subset S_{n} \subset S_{n+1} \subset \ldots \subset \mathbb{R}^n$ associated with a realizable submodule $N$ could be considered a well-defined invariant. Similarly, we also believe in the well-definedness of the lattice chain homotopy type (see again \cite{Seppo}). This would also help in the case of direct products of integrally reduced Artin $k$-algebras and Künneth-type relations, cf. subsection \ref{ss:KUNNETH}, especially Remark \ref{rem:nochainhomot}. However, this investigation will take place in a subsequent work.

\subsection{A more general case of the Independence Theorem}\label{ss:MOSTGEN}\,

In this section we present a slightly more general setting, where our argument for the Independence Theorem works.

Let $k$ be a field, $M_1, M_{-1}$ and $M_0$ vector spaces with a $k$-bilinear `multiplication' operation
\begin{equation*}
    \bullet: M_1 \times M_{-1} \rightarrow M_0, \ (m_1, m_{-1}) \mapsto m_1 \bullet m_{-1}.
\end{equation*}
[In Section \ref{s:4} these vector spaces are $M_1=\mathcal{O}$, $M_{-1}=M_0=M$ and the ring action on the module gives the $\bullet$ multiplication.]\\
Let $N_0 \leq M_0$ be a subspace. 

\begin{define}\label{def:compfilt}
    A \textit{`compatible filtration'} on the tuple $(M_1, M_{-1}, N_0\leq M_0) $ is a triple, denoted by  $(\mathcal{F}^{M_1}, \mathcal{F}^{M_{-1}}, \mathcal{F}^{M_0})$, of nonconstant decreasing subspace filtrations (i.e., if $\ell \leq \ell'$ in $\mathbb{Z}$, then the corres\-ponding subspaces $\mathcal{F}^{M_\varepsilon}(\ell) \geq \mathcal{F}^{M_\varepsilon}(\ell')$ for all $\varepsilon\in \{\pm 1, 0\}$) such that
    \begin{align*}
        \forall \ell \in \mathbb{Z}: \hspace{5mm} &\mathcal{F}^{M_1}(\ell) \leq M_1, \hspace{20mm} \mathcal{F}^{M_{-1}}(\ell) \leq M_{-1}, \hspace{20mm} \mathcal{F}^{M_0}(\ell) \leq M_0, \hspace{4mm} \mbox{ with } \nonumber \\
        &{\rm codim}_k(\mathcal{F}^{M_1}(\ell) ) < \infty , \hspace{7.5mm} {\rm codim}_k(\mathcal{F}^{M_{-1}}(\ell) ) < \infty; \nonumber\\
       \text{whereas} \hspace{5mm} &\mathcal{F}^{M_1}(0)=M_1, \hspace{10mm} \exists\, d \in \mathbb{Z}_{\leq 0}:\  \mathcal{F}^{M_{-1}}(d) = M_{-1},   \hspace{10mm} \mbox{and }  \mathcal{F}^{M_0}(0) \geq N_0. \nonumber     
    \end{align*}
    We also require the following compatibility relations with the $\bullet$ multiplication operation:
    \begin{align*}
        &\forall\, \ell, \ell' \in \mathbb{Z}: \mathcal{F}^{M_1}(\ell) \bullet \mathcal{F}^{M_{-1}}(\ell') \leq \mathcal{F}^{M_0}(\ell+\ell') \mbox{, but} \\
        \mbox{for any } m_1 \notin &\mathcal{F}^{M_1}(\ell+1) \mbox{ and } m_{-1} \notin \mathcal{F}^{M_{-1}}(\ell'+1): \ m_1 \bullet m_{-1} \notin \mathcal{F}^{M_0}(\ell+\ell'+1) \nonumber
    \end{align*}
    For technical reasons (compare with Remark \ref{rem:goingtoinfty}) we also require that 
    \begin{equation*}
        \lim_{\ell \to \infty} {\rm codim}_k(\mathcal{F}^{M_1}(\ell) ) =\infty \text{ and } \lim_{\ell \to \infty} {\rm codim}_k(\mathcal{F}^{M_{-1}}(\ell) ) =\infty.
    \end{equation*}
    We will often abridge the notation $(\mathcal{F}^{M_1}, \mathcal{F}^{M_{-1}}, \mathcal{F}^{M_0})$ to $\mathcal{F}^*$.
\end{define}

\begin{remark}\label{rem:valfunctions}
    Given such a compatible filtration, we can associate with it a triple of functions (analogues of the discrete valuations) $(\mathfrak{v}^{M_1}, \mathfrak{v}^{M_{-1}}, \mathfrak{v}^{M_0})$ defined by
    \begin{align}
    &\mathfrak{v}^{M_{\pm1}}: M_{\pm1} \rightarrow \mathbb{Z} \cup \{\infty\}, \hspace{30mm} \mathfrak{v}^{M_0}: M_0 \rightarrow \mathbb{Z} \cup \{\pm\infty\}, \hspace{20mm} \mbox{with} \nonumber\\
    & 
    m_{\pm1} \mapsto \sup\{ \ell \in \mathbb{Z}\,:\, m_{\pm1} \in \mathcal{F}^{M_{\pm1}}(\ell)\} \hspace{7mm}
    m_0 \mapsto \sup\{ \ell \in \mathbb{Z}\,:\, m_0 \in \mathcal{F}^{M_0}(\ell)\}. \nonumber
    \end{align}
    These clearly satisfy the following properties (for all $\varepsilon \in \{ \pm1, 0\}$):
    \begin{itemize}
    \item[(a)] $\frv^{M_{\varepsilon}}(am_{\varepsilon}) = \frv^{M_\varepsilon}(m_\varepsilon)$ \ for all  $0 \neq a \in k, m_{\varepsilon} \in M_{\varepsilon}$;
    \item[(b)] $\frv^{M_{\varepsilon}}(m_{\varepsilon}+m'_{\varepsilon}) \geq \min \{ \frv^{M_{\varepsilon}}(m_\varepsilon), \frv^{M_\varepsilon}(m'_{\varepsilon}) \}$ \ for all $m_\varepsilon, m'_\varepsilon \in  M_\varepsilon$;
    \item[(c)]  $\frv^{M_\varepsilon}$ is nonconstant,
    \item[(d)] and $\mathfrak{v}^{M_0}(m_1 \bullet m_{-1})= \mathfrak{v}^{M_1}(m_1) + \mathfrak{v}^{M_{-1}}(m_{-1})$ for all $m_1 \in M_1$ and $m_{-1} \in M_{-1}$.
\end{itemize}
They naturally encode the same information as the filtrations. Indeed, one has the following identity: 
\begin{center}$\mathcal{F}^{M_{\varepsilon}}(\ell) = \{ m_{\varepsilon} \in M_{\varepsilon} \,:\, \mathfrak{v}^{M_{\varepsilon}}(m_\varepsilon) \geq \ell\}$ for all $\varepsilon \in \{ \pm1, 0\}$.
\end{center}
[Although this formulation with the triple of functions seems simpler, we wanted to express the minimal conditions for our proof of the Independence Theorem in the language of filtrations first, as those are absolutely essential for all underlying notions. This equivalent reformulations shows, however, that the whole concept of our construction  relies on a multiplication-type operation and such valuative functions respecting it.]
\end{remark}

Given a \emph{multiset} of such compatible filtrations $\mathcal{D}:=\{ \mathcal{F}_1^*, \ldots , \mathcal{F}_r^*\}$, we can associate to it the \emph{`compatible multifiltration'} $\mathcal{F}_{\mathcal{D}}^*=\big(\mathcal{F}_{\mathcal{D}}^{M_{1}}, \mathcal{F}_{\mathcal{D}}^{M_{-1}}, \mathcal{F}_{\mathcal{D}}^{M_{0}}\big)$ defined as
\begin{equation*}
    \mathbb{Z}^r \ni \ell=\sum_{v} \ell_v e_v \mapsto \mathcal{F}_{\mathcal{D}}^{M_{\varepsilon}}(\ell) = \cap_{v=1}^r\mathcal{F}_v^{M_{\varepsilon}}(\ell_v) \leq M_{\varepsilon};
\end{equation*}
and corresponding height functions: 
\begin{equation*}
    \mathfrak{h}_{\mathcal{D}}, \mathfrak{h}_{\mathcal{D}}^{\circ}: \mathbb{Z}^r \rightarrow \mathbb{Z}: \  \mathfrak{h}_{\mathcal{D}}(\ell) = {\rm codim}_k\big(\mathcal{F}_{\mathcal{D}}^{M_1}(\ell) \hookrightarrow M_1\big), \ \mathfrak{h}^{\circ}_{\mathcal{D}}(\ell) = {\rm codim}_k\big(\mathcal{F}_{\mathcal{D}}^{M_{-1}}(-\ell) \hookrightarrow M_{-1}\big).
\end{equation*}
Clearly, by Lemma \ref{lemma:matroid}, they both satisfy the matroid rank inequality. Also, by our conventions $\mathcal{F}_{\mathcal{D}}^{M_0}(0) \supset N_0$.

\begin{define}
    Given a tuple $(M_1, M_{-1}, N_0 \leq M_0)$. We call a finite (multi)set of compatible filtrations $\mathcal{D}:=\{ \mathcal{F}_1^*, \ldots, \mathcal{F}_r^*\}$ a `\textit{weak realization}' of  $N_0$, if the corres\-ponding multifiltration satisfies $\mathcal{F}_{\mathcal{D}}^{M_0}(0) = N_0$.

    Similarly as before, a realization $\mathcal{D}$ is a \textit{`weak CDP realization'}, if the pair of height functions $(\mathfrak{h}_{\mathcal{D}}, \mathfrak{h}_{\mathcal{D}}^{\circ})$ satisfies the Combinatorial Duality Property. 
\end{define}

\begin{remark}\label{rem:weakreeltolas}
    The term `weak' reflects the fact that such weak realizations $\mathcal{D}:=\{ \mathcal{F}_1^*, \ldots, \mathcal{F}_r^*\}$ can be easily translated with effective lattice points $c=(c_1, \ldots, c_r)\in (\mathbb{Z}_{\geq 0})^r$ without them loosing their realizing property:
    \begin{center}
        in fact, for any compatible filtration $\mathcal{F}^*_v=(\mathcal{F}_v^{M_1}, \mathcal{F}_v^{M_{-1}}, \mathcal{F}_v^{M_0})$ and integer $c_v \in \mathbb{Z}_{\geq 0}$:\\
        the translation $\mathcal{F}^*_v[c_v]=(\mathcal{F}_v^{M_1}[c_v], \mathcal{F}_v^{M_{-1}}[-c_v], \mathcal{F}_v^{M_0})$ is again a compatible filtration, where\\
        for all $l \in \mathbb{Z}$: $\mathcal{F}_v^{M_1}[c_v](l):=\mathcal{F}_v^{M_1}(l+c_V)$ and $\mathcal{F}_v^{M_{-1}}[-c_v](l):=\mathcal{F}_v^{M_{-1}}(l-c_V)$;\\
        moreover, $\mathcal{F}_{\mathcal{D}[c]}^{M_0}(0)=\mathcal{F}_{\mathcal{D}}^{M_0}(0)=N_0$, where $\mathcal{D}[c]:=\{ \mathcal{F}_1^*[c_1], \ldots, \mathcal{F}_r^*[c_r]\}$.
    \end{center}
    
    Phrased differently, these weak realizations come in families, where the subspace $\mathcal{F}_\mathcal{D}^{M_{-1}}(0)\leq M_{-1}$ can be considered as a running parameter, which can be (and will be) fixed later (in Definition \ref{def:REAlnormalized}). Nevertheless, we will need this more `loose' setup and the corresponding Independence Theorem \ref{th:Indep-MGUF} for the proof of the Reduction Theorem \ref{th:REDANh} for the equivariant analytic lattice homology (see also part (c) of Example \ref{ex:eltolasplurig}).
\end{remark}

\begin{remark}
    Just as before, given a weak realization $\mathcal{D}:=\{ \mathcal{F}_1^*, \ldots \mathcal{F}_r^*\}$, we can always construct a weak CDP realization by considering every compatible filtration twice, or changing them to the compatible filtrations assigned to the (valuative) functions $\{2\cdot \mathfrak{v}_1^*, \ldots, 2\cdot \mathfrak{v}_r^*\}$ (see Remark \ref{rem:valfunctions} for the relation between the compatible filtrations and these valuative functions).
\end{remark}

Having a weak CDP realization $\mathcal{D}$, we can consider the lattice $\mathbb{Z}^{|\mathcal{D}|}$ and the following \textit{unnormalized} weight function
\begin{equation}\label{eq:wunnormalized}
    \mathbb{Z}^{|\mathcal{D}|} \ni \ell \mapsto w^u_{\mathcal{D}}(\ell)=\mathfrak{h}_{\mathcal{D}}(\ell) + \mathfrak{h}_{\mathcal{D}}^\circ(\ell)
\end{equation}
to obtain $\mathbb{H}_*(\mathbb{R}^{|\mathcal{D}|}, w^u_{\mathcal{D}})$ (we call $w^u_{\mathcal{D}}$ \textit{unnormalized} since it does not (necessarily) satisfies the usual $w(0)=0$ identity).
In this more general, unnormalized setup, we can apply (an adapted version of) the previously presented proof to get the following theorem.

\begin{theorem}[Independece Theorem --- the general unnormalized form]\label{th:Indep-MGUF}
    Let us fix a tuple of $k$-vec\-tor spaces $(M_1,  M_{-1}, N_0 \leq M_0)$ as before, and suppose, that there exists some weak realization of $N_0$. Then every weak CDP realization will give the same unnormalized lattice homology module, denoted by $\mathbb{UH}_*(M_1, M_{-1}, N_0 \leq M_0)$. 
    
    For a given weak CDP realization $\mathcal{D}$, it can be computed on the finite rectangle $R(d^-, d^+)$, where $d^-=\max\{ \ell \in (\mathbb{Z}_{\geq 0})^{|\mathcal{D}|}\,:\,\mathcal{F}_{\mathcal{D}}^{M_1}(\ell)=M_1\}$ and $d^+=-\max\{ \ell \in (\mathbb{Z}_{\leq 0})^{|\mathcal{D}|}\,:\,\mathcal{F}_{\mathcal{D}}^{M_{-1}}(\ell)=M_{-1}\}$. Moreover, it has Euler characteristic $0$. Even more, it only depends on the equivalence class (under the equivalence defined in the following Definition \ref{def:isotype}) of the $\bullet$ multiplication 
    $$\bullet: \dfrac{M_1}{(N_0\,:\, M_{-1})} \times \dfrac{M_{-1}}{(N_{0}:M_1)} \rightarrow \dfrac{M_\bullet}{N_0\cap M_\bullet},$$
    where the notation $(N_0\,:\, M_{-\epsilon})$ is used for the subspace $\{ m_{\epsilon} \in M_{\epsilon}\,:\, m_{\epsilon} \bullet M_{-\epsilon} \subset N_0\}$ for $\epsilon \in \{\pm 1\}$, moreover, $M_\bullet = M_1 \bullet M_{-1} \leq M_0$, i.e., the image of the $\bullet$ multiplication map in $M_0$.
    \end{theorem}

    \begin{define}\label{def:isotype}
        On the set of all $\bullet$ multiplication maps we consider the equivalence relation generated by the following relation $\mathbf{R}$: given two tuples $(M_1, M_{-1}, N_0\leq M_0)$ and $(M_1', M'_{-1}, N'_0\leq M'_0)$ of $k$-vector spaces, the $k$-bilinear `multiplication' maps
    \begin{align*}
        \bullet: \dfrac{M_1}{(N_0\,:\, M_{-1})} \times \dfrac{M_{-1}}{(N_{0}:M_1)} \rightarrow \dfrac{M_\bullet}{N_0\cap M_{\bullet}} \hspace{5mm} \text{and} \hspace{5mm } \bullet': \dfrac{M'_1}{(N'_0\,:\, M'_{-1})} \times \dfrac{M'_{-1}}{(N'_{0}:M'_1)} \rightarrow \dfrac{M'_{\bullet'}}{N'_0\cap M'_{\bullet'}}
    \end{align*}
    belong to relation $\mathbf{R}$ (i.e., $(\bullet,\bullet')\in \mathbf{R}$) if there exist $k$-linear maps $L_1:M_1 \rightarrow M_1', \ L_{-1}:M_{-1} \rightarrow M'_{-1}$ and $L_0:M_\bullet \rightarrow M'_{\bullet'}$ satisfying that 
    \begin{enumerate}
        \item $L_0^{-1}(N'_0\cap M'_{\bullet'})= N_0\cap M_{\bullet}$;
        \item for all $m_1 \in M_1$ and $m_{-1} \in M_{-1}$ we have $L_0(m_1 \bullet m_{-1})=L_1(m_1)\bullet'L_{-1}(m_{-1})$;
        \item and the induced maps on the quotients
    \begin{equation*}
        L_1^*: \dfrac{M_1}{(N_0\,:\, M_{-1})} \rightarrow \dfrac{M'_1}{(N'_0\,:\, M'_{-1})}, \hspace{5mm} L_{-1}^*:\dfrac{M_{-1}}{(N_{0}:M_1)} \rightarrow \dfrac{M'_{-1}}{(N'_{0}:M'_1)}, \hspace{5mm} L_0^*:  \dfrac{M_\bullet}{N_0\cap M_\bullet} \rightarrow \dfrac{M'_{\bullet'}}{N'_0\cap M'_{\bullet'}}
    \end{equation*}
    are $k$-vector space isomorphisms.
    \end{enumerate}
    \end{define}

    \begin{remark}
        The equivalence relation from the above Definition \ref{def:isotype} is rather cumbersome. This is due to the fact that in the absence of some stronger algebraic structure, we have to make sure artificially that the isomorphisms on the quotient level can be lifted to linear maps on the dividend vector spaces. Nevertheless, if one is satisfied with a weight function and $S_n$-spaces defined only on a finite rectangle, then one can revert to a simpler equivalence relation (cf. Remark \ref{rem:MOSTGEN}) using the  quasi-valuations from section \ref{ss:quasi} instead of the unbounded valuative functions from Remark \ref{rem:valfunctions}.
    \end{remark}

\bekezdes \textbf{Normalized version.}\label{par:Normalized}\,

If we want a normalized version (i.e., a weight function $w$ with $w(0)=0$), we have to subtract $\mathfrak{h}^{\circ}_{\mathcal{D}}(0)=\dim_k M_{-1}/\mathcal{F}^{M_{-1}}_{\mathcal{D}}(0)$ from the unnormalized weight function $w_{\mathcal{D}}^u$ (defined in (\ref{eq:wunnormalized})) to obtain the more standard-looking normalized weight function:
\begin{equation}\label{eq:wmostgeneral}
    \mathbb{Z}^{|\mathcal{D}|} \ni \ell \mapsto w_{\mathcal{D}}(\ell)=\mathfrak{h}_{\mathcal{D}}(\ell) + \mathfrak{h}_{\mathcal{D}}^\circ(\ell) - \mathfrak{h}_{\mathcal{D}}^\circ(0)= \mathfrak{h}_{\mathcal{D}}(\ell) + \mathfrak{h}_{\mathcal{D}}^\circ(\ell) - \dim_k(M_{-1}/\mathcal{F}^{M_{-1}}_{\mathcal{D}}(0)).
\end{equation}
In order for this to make sense independently of the realization, in this normalized setting we have to fix the subspace $N_{-1}:=\mathcal{F}^{M_{-1}}_{\mathcal{D}}(0) \leq M_{-1}$ as well to get a corresponding version of the Independece Theorem (cf. Remark \ref{rem:weakreeltolas}). [In fact, it would already be enough to fix the codimension of the subspace $\mathcal{F}^{M_{-1}}_{\mathcal{D}}(0) \leq M_{-1}$. However, for the sake of simplicity we omit this level of generality here.]
\begin{obs}
Consider a $k$-bilinear `multiplication' map $\bullet:M_1\times M_{-1} \rightarrow M_0$ as before and let $\mathcal{D}=\{ \mathcal{F}_1^*, \ldots, \mathcal{F}_r^*\}$ be a weak realization of some subspace $N_0\leq M_0$. We claim that if $0$ is the largest lattice point $\ell\in \mathbb{Z}^{|\mathcal{D}|}$ such that $M_1=\mathcal{F}_{\mathcal{D}}^{M_{1}}(\ell)$, i.e., for every $1 \leq v \leq r$ there exists some $ \ m_{1,v} \in M_1$ with $m_{1,v} \in \mathcal{F}^{M_1}_{v}(0) \setminus  \mathcal{F}^{M_1}_{v}(1)$ (equivalently, $\mathfrak{v}_{v}^{M_1}(m_{1,v})=0$),  then $N_{-1}=\mathcal{F}^{M_{-1}}_{\mathcal{D}}(0)=(N_0\,:\,M_1)$.
 Indeed, we can use the same argument as in the proof of part \textit{(2)} of Proposition \ref{prop:pullbackring}. In general, (without the assumption on $0$) we only have the containment $N_{-1} \subset (N_0\,:\,M_1)$.
\end{obs} 

Similarly as before, we introduce the following terminology:

\begin{define}\label{def:REAlnormalized}
    Given a tuple $(M_1, N_{-1} \leq M_{-1}, N_0 \leq M_0)$. We call a (multi)set of compatible filtrations $\mathcal{D}:=\{ \mathcal{F}_1^*, \ldots, \mathcal{F}_r^*\}$ a `\textit{strong realization}' of the pair $(N_{-1}, N_0)$, if the corres\-ponding multifiltrations satisfy $\mathcal{F}_{\mathcal{D}}^{M_{-1}}(0) = N_{-1}$ and $\mathcal{F}_{\mathcal{D}}^{M_0}(0) = N_0$. A realization $\mathcal{D}$ is a \textit{`strong CDP realization'}, if the pair of Hilbert functions $(\mathfrak{h}_{\mathcal{D}}, \mathfrak{h}_{\mathcal{D}}^{\circ})$ satisfies the Combinatorial Duality Property. 
\end{define}

\begin{theorem}[Independece Theorem --- the general normalized form]\label{th:Indep-MGF}
    Let us fix a tuple of $k$-vec\-tor spaces $(M_1, N_{-1} \leq M_{-1}, N_0 \leq M_0)$ as before, and suppose that there exists some strong realization of $(N_{-1}, N_0)$. Then every strong CDP realization $\mathcal{D}$ will give the same lattice homology module $\mathbb{H}_*(M_1, N_{-1} \leq M_{-1}, N_0 \leq M_0):=\mathbb{H}_*(\mathbb{R}^{|\mathcal{D}|}, w_{\mathcal{D}})$ (where $w_{\mathcal{D}}$ is the normalized weight function of (\ref{eq:wmostgeneral})). 
    
    For a given weak CDP realization $\mathcal{D}$, it can be computed on the finite rectangle $R(d^-, d^+)$, where $d^-=\max\{ \ell \in (\mathbb{Z}_{\geq 0})^{|\mathcal{D}|}\,:\,\mathcal{F}_{\mathcal{D}}^{M_1}(\ell)=M_1\}$ and $d^+=-\max\{ \ell \in (\mathbb{Z}_{\leq 0})^{|\mathcal{D}|}\,:\,\mathcal{F}_{\mathcal{D}}^{M_{-1}}(\ell)=M_{-1}\}$. Moreover, it has Euler characteristic ${\rm codim}_k(N_{-1} \hookrightarrow M_{-1})$. Even more, it only depends on the equivalence class (under the equivalence defined in Definition \ref{def:isotype}) of the $\bullet$ multiplication map
    $$\bullet: \dfrac{M_1}{(N_0\,:\, M_{-1})} \times \dfrac{M_{-1}}{(N_{0}:M_1)} \rightarrow \dfrac{M_\bullet}{N_0\cap M_{\bullet}},$$
    and the codimension ${\rm codim}_k(N_{-1} \hookrightarrow M_{-1})$ (where the notation $(N_0\,:\, M_{-\epsilon})$ is used for the subspace $\{ m_{\epsilon} \in M_{\epsilon}\,:\, m_{\epsilon} \bullet M_{-\epsilon} \subset N_0\}$ for $\epsilon \in \{\pm 1\}$,  moreover, $M_\bullet = M_1 \bullet M_{-1} \leq M_0$, i.e., the image of the $\bullet$ multiplication map in $M_0$.).
    \end{theorem}

\begin{example}\label{ex:eltolasplurig}
    (a) We get back the setting of Section \ref{s:4} if we choose $M_1:=\mathcal{O}$, $M_{-1}=M_0=M$ and $N_{-1}=N_0=N$, then the ring action on the module will give the $\bullet$ multiplication and the extended valuations will give the compatible filtrations, thus the Independence Theorem \ref{th:IndepMod} is a special case of the above, more general one.

    (b) A natural context in which one could use the more general version is the case of a finite Abelian group action on the algebra and the module. More concretely, if $\mathcal{O}$ is a $k$-algebra (with $k=\overline{k}$ algebraically closed), $M$ is a finitely generated $\mathcal{O}$-module and $G$ is a finite Abelian group acting on both of them equivariantly. Then we can consider the eigenspace decompositions (corresponding to the characters of $G$): $\mathcal{O}=\oplus_{\chi \in \hat{G}}\mathcal{O}_{\chi}$ and $M = \oplus_{\chi \in \hat{G}}M_{\chi}$. In this case, we have a natural multiplication map induced by the ring action: $\bullet: \mathcal{O}_\chi \times M_{\chi^{-1}} \rightarrow M_{\chi_0}$, where $\chi_0$ denotes the trivial character. If for some finite codimensional realizable submodule $N \leq M$, with $N=\oplus_{\chi\in \hat{G}}N_\chi$, we can define equivariant compatible filtrations giving strong CDP realizations for any pair $(N_{\chi^{-1}}, N_{\chi_0})$, then these will give a well defined `\textit{equivariant}' lattice homology of $N \leq M$. We will exemplify this construction in the following subsection. 

    (c) Notice that if we are given a tuple $(M_1, M_{-1}, N_0 \leq M_0)$ of $k$-vector spaces and some compatible filtration $\mathcal{F}^*$ with associated triple of functions $(\mathfrak{v}^{M_1}, \mathfrak{v}^{M_{-1}}, \mathfrak{v}^{M_0})$, then for any $c \in \mathbb{Z}, c \geq -\mathfrak{v}^{M_{1}}(M_1)$ the valuative functions $(\mathfrak{v}^{M_1}+c, \mathfrak{v}^{M_{-1}}-c, \mathfrak{v}^{M_0})$ satisfy the requirements of Remark \ref{rem:valfunctions}.  Let us denote the corresponding filtration by $\mathcal{F}^*[c]$ (compare with Remark \ref{rem:weakreeltolas}). Then we have the following result: suppose that $\mathcal{D}=\{ \mathcal{F}_1^*, \ldots, \mathcal{F}_r^*\}$ is a CDP realization of the pair $(\mathcal{F}_{\mathcal{D}}^{M_{-1}}(0),N_0)$, then for any $d =(d_1, \ldots, d_r) \in \mathbb{Z}^r, d_v \geq -\mathfrak{v}_v^{M_1}(M_1)$ for all $1 \leq v \leq r$, the (multi)set $\mathcal{D}[d]:=\{ \mathcal{F}_1^*[d_1], \ldots, \mathcal{F}_r^*[d_r]\}$ is a CDP realization of the pair $(\mathcal{F}_{\mathcal{D}}^{M_{-1}}(d), N_0)$. Moreover, the value tables of the height functions $\mathfrak{h}_{\mathcal{D}}$ (respectively $\mathfrak{h}_{\mathcal{D}}^\circ$) and $\mathfrak{h}_{\mathcal{D}[d]}$ (respectively $\mathfrak{h}_{\mathcal{D}[d]}^\circ$) are just translates of each other by the vector $d$, whereas the difference between the Euler characteristics of the corres\-ponding normalized lattice homologies comes only from the difference between the normalizing constants $\dim_k(M_{-1}/\mathcal{F}_{\mathcal{D}}^{M_{-1}}(0))$ and $\dim_k(M_{-1}/\mathcal{F}_{\mathcal{D}}^{M_{-1}}(d))$ present in the definitions of the weight functions $w_{\mathcal{D}}$ and $w_{\mathcal{D}[d]}$ (cf. (\ref{eq:wmostgeneral})). Hence $$\mathbb{H}_*(\mathbb{R}^r, w_{\mathcal{D}[d]}) \cong \mathbb{H}_*(\mathbb{R}^r, w_{\mathcal{D}})\big[-2\big(\dim_k(M_{-1}/\mathcal{F}_{\mathcal{D}}^{M_{-1}}(0)) - \dim_k(M_{-1}/\mathcal{F}_{\mathcal{D}}^{M_{-1}}(d)\big)\big] .$$
    In fact, retrospectively we see that this construction was used to obtain the  first categorification of the plurigenera of Gorenstein normal surface singularities in \cite{NSplurig} from the analytic lattice cohomology of \cite{AgNe1}, see in particular \cite[Definition 6.2.1]{NSplurig}.  
\end{example}

\begin{remark}
Here a warning is in order: the Nonpositivity Theorem \ref{th:upperbound} does not generalize (at least using the proof discussed in section \ref{s:nonpositivity}) to the setting of this current subsection (neither to the unnormalized, nor to the normalized setup). For more details see Remark \ref{rem:nononposthm}.
\end{remark}

\subsection{The equivariant analytic lattice homology of surface singularities}\,

T. Ágoston and  the first author defined the equivariant analytic lattice (co)homology of surface singularities in \cite{AgNeIII}. It is an invariant of singularities with rational homology sphere link and categorifies the equivariant geometric genus. It was a priori associated with a resolution $\phi:\widetilde{X}\rightarrow X$ and later proved to be independent of it.  In this subsection we give an equivalent description of it through our more general construction from subsection \ref{ss:MOSTGEN}. Our definition has the advantage that by Theorem \ref{th:Indep-MGF} it is immediately well-defined (i.e., it is independent of the choice of $\phi$) and comes with a natural Reduction Theorem. 

However, contrary to the case of section \ref{s:dnagy}, here we cannot get rid of the rational homology sphere assumption on the link, due to the fact that the universal abelian covering (needed for the construction itself) is well-defined only in this case.

Throughout this section let us fix a normal surface singularity $(X, o)$ with rational homology sphere link $L_X$. 
Before giving the construction we have to introduce several notions and notations. We follow \cite[Section 5]{AgNeIII} and \cite[Sections 6.2 and 6.8]{NBook}. 

\bekezdes \textbf{The combinatorics  of a resolution \cite{NOSz,NGr}.}

Let $\phi:\widetilde{X}\to X$ be a good   resolution   of $(X,o)$ with
 exceptional curve $E:=\phi^{-1}(0)$,  and  let $\cup_{v\in\calv}E_v$ be
the irreducible decomposition of $E$. 
 Let $\Gamma$ be the dual resolution graph of $\phi$.

The lattice $L:=H_2(\widetilde{X},\mathbb{Z})\cong H_2(E, \mathbb{Z})$ is
freely generated by the fundamental homology classes of  $\{E_v\}_{v\in\mathcal{V}}$ (denoted in the very same way) and is  endowed
with the natural  negative definite intersection form  $(\,,\,)$.
 The dual lattice is $L'={\rm Hom}_\Z(L,\Z) \simeq\{
l'\in L\otimes \Q\,:\, (l',L)\in\Z\}$.
It  is generated
by the (anti)dual classes $\{E^*_v\}_{v\in\mathcal{V}}$ defined
by the identity $(E^{*}_{v},E_{w})=-\delta_{vw}$ (where $\delta_{vw}$ stays for the  Kronecker symbol).
$L'$ can also be identified with $H^2(\tX,\Z)\cong H_2(\widetilde{X}, \partial \widetilde{X})$ via the Universal Coefficient Formula (since $H_1(\widetilde{X}, \mathbb{Z}) \cong H_1(E, \mathbb{Z})$ is free) and the Lefschetz duality.

Then the intersection form embeds $L$ into $L'$ with
\begin{center}
 $L'/L\cong {\rm Tors}( H_1(L_X,\mathbb{Z}))={\rm ker}\left(H_1(L_X, \mathbb{Z})\rightarrow H_1(\widetilde{X}, \mathbb{Z})\right)$ \hspace{5mm} (since $L_X \simeq \partial \widetilde{X}$),
 \end{center}
 which will be abridged by $H$.
 The class of $l'$ in $H$ is denoted by $[l']$.

There is a natural partial ordering of $L'$ and $L$: we write $l_1'\geq l_2'$ if
$l_1'-l_2'=\sum _v r_vE_v$ with every  $r_v\geq 0$. We set $L_{\geq 0}=\{l\in L\,:\, l\geq 0\}$ and
$L_{>0}=L_{\geq 0}\setminus \{0\}$.
The support of a cycle $l=\sum n_vE_v$ is defined as  ${\rm supp}(l)=\cup_{n_v\not=0}E_v$.

Since  $L_X$ is a rational homology sphere,
each $E_v$ is rational and the dual graph of any good resolution is a tree. In this case $H_1(L_X,\Z)=H$ is finite and
we denote by $\widehat{H}$ the Pontrjagin dual $\mathrm{Hom}(H,S^1)$ of the abelian group $H$.
Let $\theta:H\to \widehat{H}$ be the isomorphism $[l']\mapsto
e^{2\pi i(l',\,\cdot)}$ of $H$ with  $\widehat{H}$.

\bekezdes \textbf{Natural line bundles.}\label{ss:UAC}

As previously noted, $L'={\rm Hom}_{\mathbb{Z}}(L, \mathbb{Z})$ is isomorphic to $H^2(\widetilde{X},\bZ)$, which is the
target of the first Chern class $c_1:{\rm Pic}(\widetilde{X})\to H^2(\widetilde{X},\bZ)$.
This morphism appears in the following exact sequence (induced by the exponential exact sequence
of sheaves):
\begin{equation}\label{eq:PIC}
0\to
 H^1(\tX,\cO_{\tX}) \longrightarrow 
  {\rm Pic}(\tX)\stackrel{c_1}{\longrightarrow} H^2(\tX,\bZ)\to 0.
\end{equation}

\noindent
In this sequence  $c_1$
admits  a natural group section $s_L$  over the integral cycles
$L\subset L'$. Namely, for any $l\in L$ one takes $\cO_{\tX}(l)\in {\rm Pic}(\tX)$
with $c_1(\cO(l))=l$. By \cite{NGr} (see also \cite[Section 6.2]{NBook})
this $s_L$ can be extended in a unique way  to a natural group section
$s:L'\to {\rm Pic}(\tX)$. Its existence is basically guaranteed by the
facts that $H=L'/L$ is finite, while ${\rm Pic}^0(\tX):=H^1(\tX,\calO_{\tX})$ is torsion free.
The line bundles $s(l')$, indexed by $l'\in L'$, and denoted also by $\cO_{\tX}(l'):=s(l')$,
are called {\it natural line bundles}.
In fact, $\cL\in {\rm Pic}(\tX)$ is natural if and only if
some power of it has the form $\cO_{\tX}(l)$ for an integral cycle $l\in L$.

\bekezdes\label{bek:univabcov}{\bf The universal abelian covering.}
Let $c:(X_a,o)\to (X,o)$ denote the universal abelian covering  \ix{$c:(X_a,o)\to (X,o)$}
of $(X,o)$: $(X_a,o)$  is  the  unique normal singular germ such that
$X_a\setminus \{o\}$ is
 the regular covering of $X\setminus \{o\}$ associated with the kernel of the Hurewicz map
$\pi_1(X\setminus \{o\})\to H$ (cf. \cite{Stein}).

 Since  $\tX\setminus
E$  is isomorphic to $X\setminus \{o\}$, the map $\pi_1(\tX\setminus E)=\pi_1(X\setminus \{o\})\to H$ defines a
regular Galois covering of $\tX\setminus E$ as well. This has a unique
extension $\widetilde{c}:Z\to \tX$ with $Z$ normal and $\widetilde{c}$ finite.
(In other words,  $\widetilde{c}:Z\to \tX$ is the normalized pullback of $c$ via $\phi$.)
The (reduced) branch locus of $\widetilde{c}$ is included in $E$, and the
Galois action of $H$ extends to $Z$ as well.
Since $E$ is a normal
crossing divisor, the only singularities that $Z$ might have are
cyclic quotient singularities (which are, in particular, rational). Let
$r:\widetilde{Z}\to Z$ be  a resolution of these singular points
such that $(\widetilde{c}\circ r)^{-1}(E)_{red}$ is a normal crossing divisor. Set $p:=\widetilde{c}\circ r$.
\begin{equation}\label{eq:diagramUAC}
\begin{tikzcd}[ampersand replacement=\&, row sep=large, column sep=large]
\widetilde{Z} \arrow{r}{r} \ar{dr}{p} \& Z \arrow{d}{\widetilde{c}} \arrow{r}{\psi_a} \& (X_a,o) \arrow{d}{c}\\
 \& (\tX,E)  \arrow{r}{\phi} \&  (X,o)
\end{tikzcd}
 \end{equation}

\begin{theorem}\label{th:3.9} (\cite[Theorem 6.2.9, Corollary 6.2.11 and paragraph 6.8.20]{NBook}, see also the references therein)
$\widetilde{c}_*\cO_Z$ is a vector bundle and its $H$-eigensheaf decomposition
is:
 \begin{equation}\label{eq:3.9SUM}
\widetilde{c}_*\cO_Z = \oplus _{\alpha\in \widehat{H}}\widetilde{\cL}_\alpha,\end{equation}
where $\widetilde{\cL}_{\theta(h)}=\cO_{\tX}(-r_h)$  for any $h\in H$, where $r_h \in L'_{\geq 0}$ is the smallest effective cycle with $[r_h]=h$ (for more on these representatives see \cite[subsection 6.6.A]{NBook})
In particular, $\widetilde{c}_*\cO_Z\cong\oplus_{l'\in Q}\cO_{\tX}(-l')$, where $Q=\{ l' \in L'\,:\, \lfloor l' \rfloor =0\}$. As a corollary, $c_*\mathcal{O}_{X_a}$ also has the eigen-decomposition
\begin{equation}\label{eq:c*OXadecomp}
    c_*\mathcal{O}_{X_a} = \oplus_{\alpha \in \widehat{H}}\mathcal{L}_{\alpha} \cong \oplus_{l' \in Q}\phi_*(\mathcal{O}_{\widetilde{X}}(-l')),
\end{equation}
where the $\theta (h)$-eigenspace $\mathcal{L}_{\theta(h)}$ is indexed by $l'=r_h$. For a fixed $h \in H$, the correspondence between a $\theta(h)$-eigenfunction $f_a$ of $\mathcal{O}_{X_a}$ and a section $s_a$ of $\mathcal{O}_{\widetilde{X}}(-r_h)$ is realized via the pullback map $p^*$: if the restriction of the divisor of $s_a$ to $E$ is $l \in L$, then the restriction of the divisor of $f_a$ (lifted to $\widetilde{Z}$) to the exceptional locus is $p^*(l+r_h)$, which turns out to be an integral cycle on $\widetilde{Z}$.

Similarly, the direct image of the dualizing sheaf $\, \overline{\overline{\Omega}}^2_{Z}$ ($\, \cong r_*\Omega^2_{\widetilde{Z}}$ since $r:\widetilde{Z}\to Z$ is the resolution of cyclic quotient --- hence rational --- singularities) of $Z$ has the eigen-decomposition
\begin{equation} \label{eq:Omegabarbardecomp}
    \widetilde{c}_*\left(\overline{\overline{\Omega}}^2_{Z}\right) \cong \oplus_{l' \in Q} \mathcal{O}_{\widetilde{X}}(l') \otimes_{\mathcal{O}_{\widetilde{X}}}\Omega^2_{\widetilde{X}},
\end{equation}
where the $\theta (h)^{-1}$-eigenspace is indexed by $l'=r_h$. 
\end{theorem}

\bekezdes\label{bek:UACPG} {\bf The geometric genus of the universal abelian covering.}

 Let $(X_a,o)\to (X,o)$ be the universal abelian covering of $(X,o)$,
 and consider the notations of diagram (\ref{eq:diagramUAC}).
By definition, the  geometric genus $p_g(X_a,o)$ is
$h^1(\widetilde{Z},\cO_{\widetilde{Z}})$ (cf. (\ref{eq:p_g})). Since  $r:\widetilde{Z}\to Z$ is the resolution of the cyclic quotient singularities of $Z$, we have  $p_g(X_a,o)=h^1(\cO_Z)$. Since $\widetilde{c}$ is finite, $h^1(\cO_Z)$  equals
$\dim (R^1(\phi \circ\widetilde{c})_*\cO_Z)_o$,
 and this  has eigenspace decomposition
$\oplus_{h\in H}(R^1(\phi \circ \widetilde{c})_*\cO_Z)_{o,\theta(h)}$. By Theorem \ref{th:3.9}
the dimension of the $\theta(h)$-eigenspace is
$$p_g(X_a,o)_{\theta(h)}:= \dim\,(R^1(\phi\circ \widetilde{c})_*\cO_Z)_{o,\theta(h)}=\dim\,(R^1\phi_*\mathcal{O}_{\widetilde{X}}(-r_h))_o =h^1(\tX,\co_{\tX}(-r_h)).$$
It is called the \textit{equivariant geometric genus} of $(X,o)$ associated with $h\in H$. Sometimes we abridge  it by $p_{g,h}=p_{g,h}(X,o)$. Clearly, for $h=0$ we get $p_g(X_a,o)_{\theta(0)}=p_g(X,o)$. Also, by summation:
$$p_g(X_a,o)=\sum_{h\in H}h^1(\tX,\co_{\tX}(-r_h))=\sum_{h\in H}p_{g, h}(X, o).$$
Following Laufer's work \cite{Laufer72} on the geometric genus, in this equivariant setting we get 
\begin{equation*}
    p_g(X_a, o)_{\theta(h)}=H^0(\widetilde{X}\setminus E, \Omega^2_{\widetilde{X}}(r_h)) / H^0(\widetilde{X}, \Omega^2_{\widetilde{X}}(r_h)) \hspace{5mm} \text{(cf. \cite[ Proposition 6.8.15]{NBook}).}
\end{equation*}
We also give a third formulation, which immediately shows that it is an invariant of the singularity $(X, o)$, independent of the Stein neighbourhood and the resolution chosen.

\begin{prop} \label{prop:equivariant}
    Let $(X, o)$ be a normal surface singularity with rational homology sphere link. Let $c:(X_a, o) \rightarrow (X, o)$ denote its universal abelian cover. Consider the sheaves $\overline{\overline{\Omega}}^2_{X_a}$ and $\overline{\Omega}^2_{X_a}$ defined in subsection \ref{ss:iso} and let us denote the eigen-decompositions (with respect to the $H$-action) of their direct images by 
    \begin{equation}\label{eq:Omegabarandbarbardecomp}
c_*\overline{\overline{\Omega}}^2_{X_a}\cong \oplus _{\alpha\in \widehat{H}}\overline{\overline{\cL}}_\alpha, \hspace{10mm} c_*\overline{\Omega}^2_{X_a}\cong \oplus _{\alpha\in \widehat{H}}\overline{\cL}_\alpha.
\end{equation}
Then we claim, that for any $h \in H$ the equivariant geometric genus $p_{g, h}$ is the dimension of the following complex vector space (independent of the good resolution $\phi: (\widetilde{X}, E) \rightarrow (X, o)$ and  the Stein neighbourhood chosen)
\begin{equation}\label{eq:p_g,hmodulus}
    \dfrac{H^0(\widetilde{X}\setminus E, \Omega^2_{\widetilde{X}}(r_h))}{H^0(\widetilde{X}, \Omega^2_{\widetilde{X}}(r_h))} \cong \dfrac{\left(\overline{\overline{\mathcal{L}}}_{\theta(h)^{-1}}\right)_o}{\left(\overline{\mathcal{L}}_{\theta(h)^{-1}}\right)_o} \cong \left(\left( \dfrac{c_* \overline{\overline{\Omega}}^2_{X_a}}{c_* \overline{\Omega}^2_{X_a}}\right)_o \right)_{\theta(h)^{-1}}. 
\end{equation} 
Moreover, the equivalence classes (under the equivalence defined in Definition \ref{def:isotype}) of the induced $\bullet$ multiplication maps (natural multiplication of sections) 
\begin{center}
    $\bullet_\phi:\dfrac{H^0(\widetilde{X}, \mathcal{O}_{\widetilde{X}}(-r_h))}{\left(H^0(\widetilde{X}, \Omega^2_{\widetilde{X}})\,:\,H^0(\widetilde{X}\setminus E, \Omega^2_{\widetilde{X}}(r_h))\right)} \times \dfrac{H^0(\widetilde{X}\setminus E, \Omega^2_{\widetilde{X}}(r_h))}{H^0(\widetilde{X}, \Omega^2_{\widetilde{X}}(r_h))}  \rightarrow \dfrac{H^0(\widetilde{X}\setminus E, \Omega^2_{\widetilde{X}})}{H^0(\widetilde{X}, \Omega^2_{\widetilde{X}})}$
\end{center} and 
$\ \ \bullet_{\mathcal{L}}:\dfrac{(\mathcal{L}_{\theta(h)})_o}{\left((\overline{\Omega}^2_X)_o \,:\, \Big(\overline{\overline{\mathcal{L}}}_{\theta(h)^{-1}}\Big)_o\right)} \times \dfrac{\left(\overline{\overline{\mathcal{L}}}_{\theta(h)^{-1}}\right)_o}{\left(\overline{\mathcal{L}}_{\theta(h)^{-1}}\right)_o} \rightarrow \dfrac{\left(\overline{\overline{\Omega}}^2_{X}\right)_o}{\left(\overline{\Omega}^2_{X}\right)_o}\ \ $ also agree.
\end{prop}

\begin{proof}  
    We will examine the sheaves $c_*\overline{\overline{\Omega}}^2_{X_a}$ and $c_*\overline{\Omega}^2_{X_a}$ more carefully. First notice, that 
    $$\overline{\Omega}^2_{X_a}\cong (\psi_a \circ r)_*\Omega^2_{\widetilde{Z}}\cong (\psi_a)_* \overline{\overline{\Omega}}^2_Z\,$$
    since $r:\widetilde{Z}\to Z$ is the resolution of rational singularities. Therefore, by Theorem \ref{th:3.9} and (\ref{eq:Omegabarbardecomp}),
    \begin{center}
        $c_*\overline{\Omega}^2_{X_a}\cong \phi_* \widetilde{c}_*\overline{\overline{\Omega}}^2_Z \cong \phi_*\left(\oplus_{h\in H} \mathcal{O}_{\widetilde{X}}(r_{-h}) \otimes_{\mathcal{O}_{\widetilde{X}}}\Omega^2_{\widetilde{X}}\right)\cong \oplus_{h \in H}\phi_*\left((\mathcal{O}_{\widetilde{X}}(r_{-h})) \otimes_{\mathcal{O}_{\widetilde{X}}} \Omega^2_{\widetilde{X}} \right)$.
    \end{center} Thus $\overline{\mathcal{L}}_{\theta(h)^{-1}}=\phi_*\left((\mathcal{O}_{\widetilde{X}}(r_h)) \otimes_{\mathcal{O}_{\widetilde{X}}} \Omega^2_{\widetilde{X}}\right)$ and for $h=0$ we have $\overline{\mathcal{L}}_{\theta(0)}=\overline{\Omega}^2_X$. 

    For the second part, we introduce the common notation $\iota$ for each of the inclusion maps 
    \begin{center}
        $X \setminus o \hookrightarrow X, \hspace{5mm} X_a \setminus o \hookrightarrow X_a, \hspace{5mm} \widetilde{X} \setminus E \hookrightarrow \widetilde{X} \hspace{5mm}$ and $Z \setminus E_a \hookrightarrow Z$,
    \end{center} where $E_a$ is the exceptional set of the partial resolution $\psi_a: Z \rightarrow X_a$, or, equivalently, the reduced preimage of $E$ under $\widetilde{c}$. Now \begin{align*}
        c_*\overline{\overline{\Omega}}^2_{X_a} \cong & \ c_* \iota_* \Omega^2_{X_a \setminus o} \cong c_* \iota_* (\psi_a)_*\Omega^2_{Z \setminus E_a}\cong \phi_* \iota_* \widetilde{c}_*\left(\overline{\overline{\Omega}}^2_Z \Big|_{Z \setminus E_a}\right) \\
        \cong & \  \iota_* \phi_* \oplus_{h \in H}  \left(\mathcal{O}_{\widetilde{X}}(r_{-h})\Big|_{\widetilde{X} \setminus E} \otimes_{\mathcal{O}_{\widetilde{X}\setminus E}} \Omega^2_{\widetilde{X}\setminus E}\right) \hspace{20mm} \text{by (\ref{eq:Omegabarbardecomp})} \\
        \cong & \ \iota_* \oplus_{h \in H} \left(    \phi_* \Big(\mathcal{O}_{\widetilde{X}}(r_{-h})\Big|_{\widetilde{X} \setminus E}\Big) \otimes_{\mathcal{O}_{X\setminus o}} \Omega^2_{X\setminus o} \right)\\
        \cong & \  \oplus_{h \in H}  \iota_*  \left( \phi_* \Big(\mathcal{O}_{\widetilde{X}}(r_{-h})\Big|_{\widetilde{X} \setminus E}\Big) \otimes_{\mathcal{O}_{X \setminus o}} {{\Omega}}^2_{X\setminus o}\right).
    \end{align*}
Hence $\overline{\overline{\mathcal{L}}}_{\theta(h)^{-1}}=\iota_*\left(\phi_*\Big(\mathcal{O}_{\widetilde{X}}(r_h)\Big|_{\widetilde{X}\setminus E}\Big) \otimes_{\mathcal{O}_{X\setminus o}} \Omega^2_{X\setminus o}\right)$ and for $h=0$ we get $\overline{\overline{\mathcal{L}}}_{\theta(0)}=\overline{\overline{\Omega}}^2_X$. 

Let us now take a closer look at the exact sequence
\begin{equation} \label{eq:Lexactseq}
    0 \rightarrow \overline{\mathcal{L}}_{\theta(h)^{-1}} \rightarrow \overline{\overline{\mathcal{L}}}_{\theta(h)^{-1}} \rightarrow \overline{\overline{\mathcal{L}}}_{\theta(h)^{-1}}/\overline{{\mathcal{L}}}_{\theta(h)^{-1}} \rightarrow 0.
\end{equation}
The right hand term is clearly supported in $o$, whereas $H^0(X, \overline{\mathcal{L}}_{\theta(h)^{-1}})=H^0(\widetilde{X}, \Omega^2_{\widetilde{X}}(r_h))$ and 
$$H^0(X, \overline{\overline{\mathcal{L}}}_{\theta(h)^{-1}})=H^0\Big(\widetilde{X} \setminus E, \Omega^2_{\widetilde{X}\setminus E}\otimes _{\mathcal{O}_{\widetilde{X}\setminus E}}\mathcal{O}_{\widetilde{X}}(r_h)\Big|_{\widetilde{X}\setminus E}\Big)=H^0(\widetilde{X}\setminus E, \Omega^2_{\widetilde{X}}(r_h)).$$
Now, $H^1(\widetilde{X}, \Omega^2_{\widetilde{X}}(r_h))=0$ by the generalized Grauert--Riemenschneider Vanishing (cf. \cite[Theorem 6.4.3]{NBook}), hence we get $H^1(X, \overline{\mathcal{L}}_{\theta(h)^{-1}})=0$ from the Leray spectral sequence. Thus, from the cohomological long exact sequence of (\ref{eq:Lexactseq}) we get the desired identity (\ref{eq:p_g,hmodulus}). The equivalence of the $\bullet$ multiplication maps (under the equivalence defined in Definition \ref{def:isotype} --- in fact, the pair $(\bullet_\phi, \bullet_{\mathcal{L}})$ belongs to the relation $\mathbf{R}$) now follows similarly to, e.g., the proof of Corollary \ref{cor:irreg&p_gcategorification}, we omit its detailed discussion here.
\end{proof}

\begin{cor}\label{cor:equivariant}
    Let $(X, o)$ be a normal surface singularity with rational homology sphere link. Fix a  Stein representative $X$ 
and let $\phi: \widetilde{X}\rightarrow X$ denote a resolution with reduced exceptional divisor $E$. Then, the for the fixed tuple 
\begin{center}$\left(H^0(\widetilde{X}, \mathcal{O}_{\widetilde{X}}(-r_h)), H^0(\widetilde{X}, \Omega^2_{\widetilde{X}}(r_h)) \leq H^0(\widetilde{X}\setminus E, \Omega^2_{\widetilde{X}}(r_h)),  {H^0(\widetilde{X}, \Omega^2_{\widetilde{X}})} \leq  H^0(\widetilde{X}\setminus E, \Omega^2_{\widetilde{X}})\right)$
\end{center}
there exists a realization of the pair $(H^0(\widetilde{X}, \Omega^2_{\widetilde{X}}(r_h)), H^0(\widetilde{X}, \Omega^2_{\widetilde{X}}))$, hence, its lattice homology is well-defined. It is independent of the Stein representative and good resolution chosen, as it agrees with $$\mathbb{H}_*\left((\mathcal{L}_{\theta(h)})_o, \left(\overline{\mathcal{L}}_{\theta(h)^{-1}}\right)_o \leq \left(\overline{\overline{\mathcal{L}}}_{\theta(h)^{-1}}\right)_o,  \left(\overline{\Omega}^2_{X}\right)_o \leq \left(\overline{\overline{\Omega}}^2_{X}\right)_o\right).$$ 
It is an invariant of the singularity and the class $h \in H$ (or the corresponding spin$^c$-structure on $L_X$) and has Euler characteristic the equivariant geometric genus $p_g(X_a, o)_{\theta(h)}$. In fact, this module  is identical with $\mathbb{H}_{an, *}((X, o), h)$, the equivariant analytic lattice homology of $(X, o)$ associated with $h \in H$ from \cite{AgNeIII}  of Ágoston and the first author. 
\end{cor}

\begin{proof}
    We prove the statements by checking that the multifiltrations used in \cite{AgNeIII} give compatible multifiltrations in the sense of Definition \ref{def:compfilt}, even more, together they give a strong CDP realization of the pair $(H^0(\widetilde{X}, \Omega^2_{\widetilde{X}}(r_h)), H^0(\widetilde{X}, \Omega^2_{\widetilde{X}}))$. This will prove the second equivalence, whereas the first one will then just be the consequence of Theorem \ref{th:Indep-MGF} and Proposition \ref{prop:equivariant}. 

    Let us then consider the multifiltrations from \cite{AgNeIII}: for any $\ell \in L=\mathbb{Z}\langle E_v\rangle_{v \in \mathcal{V}}$ 
    \begin{align}
        \ell \mapsto \mathcal{F}^{M_1}(\ell) =  H^0(\widetilde{X}, \calO_{\tX}(- \ell-r_h)) \leq H^0(\widetilde{X}, \calO_{\tX}(-r_h)) & \nonumber\\ \ell \mapsto \mathcal{F}^{M_{-1}}(\ell)=H^0(\widetilde{X}, \Omega^2_{\tX}(-\ell+r_h))&\, \leq H^0(\widetilde{X} \setminus E, \Omega^2_{\tX}(r_h)) \label{eq:equivfilt} \\
        \ell \mapsto \mathcal{F}^{M_0}(\ell) =&\, H^0(\widetilde{X}, \Omega^2_{\widetilde{X}}(-\ell)) \leq H^0(\widetilde{X} \setminus E, \Omega^2_{\widetilde{X}}) \nonumber
    \end{align}
    and these are clearly compatible with the natural multiplication map of sections
    $$H^0(\widetilde{X}, \calO_{\tX}(-r_h)) \times H^0(\widetilde{X} \setminus E, \Omega^2_{\tX}(r_h)) \rightarrow H^0(\widetilde{X} \setminus E, \Omega^2_{\widetilde{X}});$$
    and also $\mathcal{F}^{M_{-1}}(0) = H^0(\widetilde{X}, \Omega^2_{\widetilde{X}}(r_h)), \ \mathcal{F}^{M_0}(0)=H^0(\widetilde{X}, \Omega^2_{\widetilde{X}})$. By \cite[Lemma 6.3.1]{AgNeIII}, the corresponding height functions 
    \begin{equation}\label{eq:handhcircequiv}
        \mathfrak{h}(\ell)=\dim \frac{H^0(\widetilde{X}, \calO_{\tX}(-r_h))}{H^0(\widetilde{X}, \calO_{\tX}(- \ell-r_h))} \hspace{10mm} \text{and} \hspace{10mm} \mathfrak{h}^\circ(\ell)= \dim \frac{H^0(\widetilde{X} \setminus E, \Omega^2_{\tX}(r_h))}{H^0(\widetilde{X}, \Omega^2_{\tX}(\ell+r_h))}
    \end{equation}
    satisfy the Combinatorial Duality Property, which implies that $(\mathcal{F}^{M_1}, \mathcal{F}^{M_{-1}}, \mathcal{F}^{M_0})$ is a CDP realization (in the sense of Definition \ref{def:REAlnormalized}) of the submodules and the weight function corresponding to this setup (cf. formula (\ref{eq:wmostgeneral})) indeed computes the lattice homology module. So, in fact, this is well-defined and isomorphic to $\mathbb{H}_{an, *}((X, o), h)$.
\end{proof}

Our construction has the advantage that the more general forms of the Independence Theorem imply its invariance directly. Moreover, the Reduction Theorem \cite[Theorem 9.1.6]{AgNeIII} is also a direct consequence of these, similarly to the analytic ($h=0$) case discussed in subsection \ref{ss:RedTh}. Here though we have to make an important observation: the Reduction Theorem \cite[Theorem 9.1.6]{AgNeIII} uses, in fact, an a priori unnormalized weight function (i.e., $\overline{w}(0) \neq 0$), hence when proving it we have to rely on the unnormalized version of the Independence Theorem \ref{th:Indep-MGUF}. Before presenting the detailed statement we introduce some notations and terminology.

\begin{nota}
    Recall the notion of the cohomological cycle $Z_{coh}\in (\mathbb{Z}_{\geq 0})^r$ from subsection \ref{ss:zcoh}, which corresponds to the general $d_{\mathcal{D}}$ lattice point (cf. Notation \ref{not:d_D}) in the case of the analytic lattice homology (of $2$-forms) of normal surface singularities. Now, for any class $h\in H$ we write $Z_{coh,h}$ for $Z_{coh}(\calO_{\tX}(-r_h))$ (for its existence see \cite[Proposition 5.4.1]{AgNeIII}), i.e., $Z_{coh,h}$ is the smallest effective cycle $d \in L_{\geq 0}$ satisfying
    \begin{equation}\label{eq:h^1}
p_{g,h}=h^1(\tX,\calO_{\tX}(-r_h))=h^1(d,\calO_{\tX}(-r_h))=\dim \frac{H^0(\widetilde{X}, \Omega^2_{\widetilde{X}}(r_h+d))}{H^0(\widetilde{X}, \Omega^2_{\widetilde{X}}(r_h))}.
\end{equation}
Equivalently, any section of $\Omega^2_{\widetilde{X}}(r_h)$ has pole along an exceptional divisor $E_v$ of order at most the corresponding coefficient $(Z_{coh, h})_v$. 
Then we can compute $\mathbb{H}_{an,*}((X, o), h)$ on the finite rectangle $R(0, d)$ for any lattice point $d \geq Z_{coh,h}$ (cf. \cite[Lemma 6.1.9]{AgNeIII}).
\end{nota}

Let us now suppose, that $\overline{\mathcal{V}} \subset \mathcal{V}$ is a WR-set, i.e., it satisfies the condition that if $\omega\in H^0(\widetilde{X}\setminus E,\Omega^2_{\widetilde{X}})$ has no poles along 
$\cup_{v\in \overline{\mathcal{V}}}E_v$, then $\omega$ has no poles at all (see conditions (\ref{eq:forms}) and (\ref{eq:hegy}) and Remark \ref{rem:badvertices} (b)). Hence, if we use the disjoint  decomposition  $\calv=\overline{\calv}\sqcup \calv^*$  and write
any  $\ell\in L$
as $\overline{\ell}+\ell^*$, or $(\overline{\ell}, \ell^*)$,
where $\overline {\ell}$ and  $\ell^*$ are  supported on $\overline{\calv}$ and
 $\calv^*$ respectively, then 
 \begin{equation}\label{eq:d*}
     H^0(\widetilde{X},\Omega^2_{\widetilde{X}}(d^*))\cong
H^0(\widetilde{X},\Omega^2_{\widetilde{X}}) \text{ for any effective cycle } d \in L_{\geq 0}.
 \end{equation}

 Now, the Reduction Theorem  \cite[Theorem 9.1.6]{AgNeIII} of Ágoston and the first author states the following:

\begin{theorem}[Reduction Theorem for the equivariant analytic lattice homology] \label{th:REDANh} 
Let $\overline{\mathcal{V}} \subset \mathcal{V}$ be a WR-set, $d\geq Z_{coh,h}$ and denote by $\overline{R}$ the rectangle $R(0, \overline{d})$,  the $\overline{\calv}$-projection of $R(0, d)$.
For any $\overline{\ell}\in \overline {R}\cap \mathbb{Z}^r$ define  the weight function (using the notations of (\ref{eq:handhcircequiv})):
$$\overline{w}_0(\overline{\ell})=\hh(\overline{\ell})+\hh^\circ (\overline{\ell}+d^*)-p_{g,h}
=\hh(\overline{\ell})-h^1(\calO_{\overline{\ell}+d^*}(-r_h)),$$
and extend it to higher dimensional cubes via the formula (\ref{eq:9weight}). 
Then there exists a bigraded $\Z[U]$-module isomorphism
$$\bH_{an, *}((X, o), h)\cong \bH_{*}(\overline{R}, \overline{w}).$$
\end{theorem}

\begin{remark}
    We want to stress, that in this setup $\overline{w}_0(0)$ is not necessarily $0$ (i.e., $\overline{w}_0$ is not normalized in the sense of paragraph \ref{par:Normalized}), it actually depends on the analytic properties of the line bundle $\mathcal{O}_{\widetilde{X}}(-r_h)$. This corresponds to the fact that the filtrations (introduced in (\ref{eq:equivfilt})) corresponding to the subset $\overline{\mathcal{V}} \subset \mathcal{V}$ might not give a strong realization of the pair $(H^0(\widetilde{X}, \Omega^2_{\widetilde{X}}(r_h)), H^0(\widetilde{X}, \Omega^2_{\widetilde{X}}))$.
\end{remark}

\begin{proof}[A new conceputal proof of the Reduction Theorem \ref{th:REDANh}]
    By identity (\ref{eq:d*}) we get that the compatible filtrations (introduced in (\ref{eq:equivfilt})) corresponding to the subset $\overline{\mathcal{V}} \subset \mathcal{V}$ give a weak realization of the subspace $H^0(\widetilde{X}, \Omega^2_{\widetilde{X}})\leq H^0(\widetilde{X} \setminus E, \Omega^2_{\widetilde{X}})$. Hence, if we denote by $\mathfrak{h}_{\overline{\mathcal{V}}}$ and $\mathfrak{h}_{\overline{\mathcal{V}}}^\circ$ the corresponding height functions, then, by the unnormalized Independence Theorem \ref{th:Indep-MGUF}, we get:
    \begin{equation*}
        \bH_{an, *}((X, o), h)=\bH_{*}(\mathbb{R}^{|\mathcal{V}|}, \mathfrak{h}+ \mathfrak{h}^\circ-p_{g, h})\cong \bH_{*}(\mathbb{R}^{|\overline{\mathcal{V}}|}, \mathfrak{h}_{\overline{\mathcal{V}}}+ \mathfrak{h}^\circ_{\overline{\mathcal{V}}}-p_{g, h}) \cong \bH_{*}(\overline{R}, \mathfrak{h}_{\overline{\mathcal{V}}}+ \mathfrak{h}^\circ_{\overline{\mathcal{V}}}-p_{g, h}). 
    \end{equation*}
    Moreover, by the definition of $Z_{coh, h}$ (see also the description after (\ref{eq:h^1})) for any  $\overline{\ell} \in \overline{R}\cap \mathbb{Z}^{|\overline{\mathcal{V}}|}$  we have $\mathfrak{h}_{\overline{\mathcal{V}}}(\overline{\ell})=\mathfrak{h}(\overline{\ell})$ and $\mathfrak{h}^\circ_{\overline{\mathcal{V}}}(\overline{\ell})=\mathfrak{h}^\circ(\overline{\ell} + d^*)$, hence $\bH_{*}(\overline{R}, \mathfrak{h}_{\overline{\mathcal{V}}}+ \mathfrak{h}^\circ_{\overline{\mathcal{V}}}-p_{g, h}) \cong \bH_{*}(\overline{R}, \overline{w}).$
\end{proof}

\section{The proof of the Nonpositivity Theorem}\label{s:nonpositivity}

In this section we will state and prove the Nonpositivity Theorem \ref{th:upperbound} in a slightly more general form. We follow the strategy discussed in \cite[Section 6]{KNS2} (with the required modifications), where the analogue of this theorem was presented in the case of the analytic lattice \emph{co}homology of reduced curve singularities. 

\subsection{The general statement}\,

We use the notations and notions defined in sections \ref{ss:CombLattice} and \ref{s:4}. 

\begin{theorem} \label{th:contr}
    Let $k$ be an infinite field, $\mathcal{O}$ a $k$-algebra, $M$ an $\mathcal{O}$-module and $N \leq M$ a realizable submodule. Suppose that there exists a CDP realization $\mathcal{D}=\{ \mathfrak{v}_1, \ldots, \mathfrak{v}_r\}$ of $N$ (in the sense of Definition \ref{def:REAL}), such that there exists a prime ideal $\mathfrak{p} \triangleleft \mathcal{O}$ with  $\mathfrak{p}=\{f \in \mathcal{O}\,:\, \mathfrak{v}_v(f)\geq 1\}$ for all $v \in \mathcal{V}=\{1, \ldots, r\}$ (usually $\mathfrak{p}$ is a maximal ideal). Let us denote $\mathfrak{d}:=\dim_k(\mathcal{O}/\mathfrak{p})< \infty$. Then the $S_{n, \mathcal{D}}$ spaces corresponding to this setup are contractible for $n \geq \mathfrak{d}$. Hence, the weight-grading of the \emph{reduced} lattice homology is not lower than $2-2\mathfrak{d}$, i.e., $\mathbb{H}_{{\rm red}, q, -2n}(N\hookrightarrow_{\mathcal{O}} M)=0$ for all $q \geq 0$ and $n\geq \mathfrak{d}$.
    
    Moreover, if $N \neq M$, then $S_{\mathfrak{d}-1, \mathcal{D}}$ is in fact non-connected. 
\end{theorem}

\begin{proof}[Proof of the Nonpositivity Theorem \ref{th:upperbound} using Theorem \ref{th:contr}]
    If $k$ is algebraically closed, then it is automatically infinite. Moreover, for valuations $\mathfrak{v}_v$ satisfying  Assumption \ref{ass:fincodim}  with $(\mathcal{O}, \mathfrak{m})$ local, the prime ideal $\{f \in \mathcal{O}\,:\, \mathfrak{v}_v(f)\geq 1\}$ must be the unique, necessarily finite $k$-codimensional, maximal ideal $\mathfrak{m}$. In fact, since $k = \overline{k}$, we have $\dim_k(\mathcal{O}/\mathfrak{m})=\mathfrak{d}=1$. Hence we can apply Theorem \ref{th:contr} with $\mathfrak{d}=1$ to obtain the desired result.
\end{proof}

\begin{example}
    Both $\mathbb{H}_{an, *}(X, o)$ of a normal isolated singularity of $\dim_{\mathbb{C}}\geq 2$ and $\mathbb{H}_{an, *}(C, o)$ of a reduced curve singularity satisfy the above  assumptions with $\mathfrak{d}=1$, hence, have the nonpositivity property (cf. Corollaries \ref{cor:upperboundhighdim} and \ref{cor:upperboundcurves}).
\end{example}

\begin{example}
    $\mathfrak{d}>1$ can happen, for example, if $K/k$ is a finite field extension and $\mathcal{O}$ is in fact a $K$-algebra. E.g., $\mathbb{H}_*((t)\hookrightarrow_{\mathbb{R}}\mathbb{C}[t])\cong \mathcal{T}_0^- \oplus \mathcal{T}_0(2)$.
\end{example}

\begin{remark}\label{rem:Dbar}
In order to prove Theorem \ref{th:contr} it is useful to notice that the homotopy type of the $S_{n, \mathcal{D}}$-spaces are well-defined: by the Independence Theorem \ref{th:IndepMod} they are independent of the CDP realization chosen (see also Theorem \ref{th:IndepModTop}). Thus it is enough to prove their contractibility for a convenient modification $\overline{\mathcal{D}}$ of the original\,/\,given realization $\mathcal{D}$. For the sake of simpler arguments we will make the following operations on the collection of extended valuations:
\begin{itemize}
    \item first consider the unique largest lattice point $a =  \sum_{v \in \mathcal{V}}a_v e_v \in \mathbb{Z}^r$ satisfying that $\mathcal{F}_{\mathcal{D}}^M(a)=N$ (i.e., $a_v=\mathfrak{v}_v^M(N) \geq 0$ for all $v \in \mathcal{V}$) and translate the extended valuations by the corresponding coordinates: $(\mathfrak{v}_v, \mathfrak{v}_v^M) \mapsto (\mathfrak{v}_v, \mathfrak{v}_v^M-a_v)$ for all $v \in \mathcal{V}$;
    \item then double the valuations to obtain
    \begin{equation}\label{eq:modifiedD}
        \overline{\mathcal{D}}:=\{(\overline{\mathfrak{v}}_1, \overline{\mathfrak{v}}_1^M)=(2\mathfrak{v}_1, 2\mathfrak{v}_1^M-2a_1), \ldots , (\overline{\mathfrak{v}}_r, \overline{\mathfrak{v}}_r^M)=(2\mathfrak{v}_r, 2\mathfrak{v}_r^M-2a_r)\}.
    \end{equation}
\end{itemize}
Now $\overline{\mathcal{D}}$ is still a CDP realization of $N$, satisfying the conditions of Theorem \ref{th:contr} and has $\overline{\mathfrak{v}}^M_v(N)=0$ for all $v \in \mathcal{V}$. Also, by the Independece Theorem, $S_{n, \mathcal{D}}$ and $S_{n, \overline{\mathcal{D}}}$ are homotopy equivalent for every $n \in \mathbb{Z}$. Due to their additional good properties, for the proof of Theorem \ref{th:contr} we will use the modified collection $\overline{\mathcal{D}}$ and the corresponding $S_{n, \overline{\mathcal{D}}}$-spaces.
\end{remark}

\begin{remark}\label{rem:w(R(0,E)<kappa}
    Notice that, since $\mathfrak{p}=\{f \in \mathcal{O}\,:\, \overline{\mathfrak{v}}_v(f)\geq 1\}$ for all index $v \in \mathcal{V}$, for all nonempty subsets $\emptyset \neq I \subset \mathcal{V}$ we have $\mathcal{F}_{\overline{\mathcal{D}}}(e_{I}) = \mathfrak{p}$ and, hence, $\mathfrak{h}_{\overline{\mathcal{D}}}(e_I)=\mathfrak{d}$ (where $e_I=\sum_{v\in I}e_v$). Now, by CDP for the pair $(\mathfrak{h}_{\overline{\mathcal{D}}}, \mathfrak{h}_{\overline{\mathcal{D}}}^{\circ})$ we get that $w_{\overline{\mathcal{D}}}(e_v)=\mathfrak{d}$ for all $v\in \mathcal{V}$ and $w_{\overline{\mathcal{D}}}(e_I) \leq \mathfrak{d}$ for all $I \subset \mathcal{V}$.
\end{remark}

The proof of Theorem \ref{th:contr} follows the next strategy: 
we construct a strong deformation  retraction from
the (contractible) space $R(0,d)$ ($d \geq d_{\overline{\mathcal{D}}}$, see subsection \ref{ss:ddep}) to its subspace $S_{n, \overline{\mathcal{D}}}\cap R(0,d)$
 for any  $n\geq \dim_k (\mathcal{O}/\mathfrak{p})=\mathfrak{d}$. This is done inductively
 via strong deformation retractions  from  $S_{n, \overline{\mathcal{D}}}\cap R(0,d)$
 to  $S_{n-1, \overline{\mathcal{D}}}\cap R(0,d)$ for descending values of $n$, which will, in turn, be the composition of deformations eliminating
 lattice points
 $\ell \in S_{n, \overline{\mathcal{D}}}\cap R(0,d)$ of maximal weight $w_{ \overline{\mathcal{D}}}(\ell) =n>\mathfrak{d}$  in a prescribed order. In order to describe the process we need to introduce some relevant notions and their properties.

\subsection{Generalized local minima}\label{ss:glm}\,

The question of contractibility, or, even more, of connectedness of the $S_n$ spaces is strongly related to the local minimum points of the weight function (i.e., those lattice points $\ell \in \mathbb{Z}^r$ for which  for any $v \in \mathcal{V}$ and $\epsilon_v\in\{\pm1\}$ we have $w(\ell)<w(\ell + \epsilon_v e_v)$). Indeed, any such lattice point $\ell$ corresponds to a distinct connected component of $S_{w(\ell)}$ consisting only of itself. For more relations between the local minimum points of the weight function and the corresponding lattice homology\,/\,graded root consult \cite[Subsection 2.3]{KNS1} or \cite[Subsection 3.1]{KNS2}. 
It turns out that in order to understand the local minima of the weight function, it is convenient to use a slightly more general concept, namely that of \emph{`generalized local minima'}. In this subsection we give a detailed discussion of this notion and its main properties. These form the technical backbone of the proof of Theorem \ref{th:contr}.

\begin{define}
    A lattice point $p \in (\bZ_{\geq 0})^r$ satisfying the condition
      $  w(p)< w(p-e_v)$ for any $v \in \mathcal{V}$
     is said to be a \textit{`generalized local minimum point'} of the weight function.
\end{define}

Let us consider the previous setting: let $k$ be an infinite field, $\mathcal{O}$ a $k$-algebra, $M$ an $\mathcal{O}$-module and $N \leq M$ a realizable submodule, and let $w$ denote the weight function corresponding to a CDP realization of $N$. Then we have the following results:

\begin{lemma}\label{lem:existsm_p}
    $p\in (\mathbb{Z}_{\geq 0})^r$ is a generalized local minimum point of the weight function if and only if there exists some $m_p \in M$ such that $\overline{\mathfrak{v}}^M(m_p)=-p$.
\end{lemma}

\begin{proof}
The sufficiency part is clear: if $\overline{\mathfrak{v}}^M(m_p)=-p$, then $\mathcal{F}_{\overline{\mathcal{D}}}^M(-p+e_v) \subsetneq \mathcal{F}_{\overline{\mathcal{D}}}^M(-p)$ for all $v \in \mathcal{V}$, hence, $\mathfrak{h}^\circ(p-e_v) > \mathfrak{h}^\circ(p)$ and we can use the Combinatorial Duality Property. 

Let us prove necessity. 
     The inequality $w(p-e_v) > w(p)$ implies that $\mathcal{F}_{\overline{\mathcal{D}}}^M(-p+e_v) \neq \mathcal{F}_{\overline{\mathcal{D}}}^M(-p)$ for every index $v\in \mathcal{V}$, i.e., there exists some $m^{v} \in \mathcal{F}_{\overline{\mathcal{D}}}^M(-p)$ with $\overline{\mathfrak{v}}_v^M(m^{v})=-p_v$. 
    
    Then we claim that $m_p$ can be chosen as a `generic' $k$-linear combination of these $m^{v}$-s (since $|k|=\infty$ we have an infinite number of distinct nonzero field elements, so genericity makes sense). 
    Indeed, suppose inductively that we have already constructed a linear combination $\overline{m}^{w}$ satisfying the properties
    \begin{equation}\label{eq:m^wprop}
        \overline{m}^{w} \in \mathcal{F}^M_{\overline{\mathcal{D}}}(-p) \text{ and }\overline{\mathfrak{v}}_v^M(\overline{m}^w)=-p_v \text{ for all }v\in\{1, \ldots, w\}.
    \end{equation} We prove that we can find some $b \in k$, such that $\overline{m}^w+b\cdot m^{w+1}$ can be chosen as $\overline{m}^{w+1}$ (satisfying the conditions (\ref{eq:m^wprop})), because for every index $1 \leq v \leq w+1$ there can not be more than one $b_v \in k$, such that $\overline{\mathfrak{v}}_{v}^M(\overline{m}^w+b_v\cdot m^{w+1}) > -p_v$. Indeed, if we assume indirectly that for some index $1 \leq v \leq w+1$ and different elements $b_{v,1} \neq b_{v,2} \in k$ we have  both
    \begin{center}
    $\overline{\mathfrak{v}}_{v}^M(\overline{m}^w+b_{v,1}\cdot m^{w+1}) > -p_v$ and $\overline{\mathfrak{v}}_{v}^M(\overline{m}^w+b_{v,2}\cdot m^{w+1}) > -p_v$,
    \end{center} 
    then these require 
    \begin{align*}
        -p_v= \overline{\mathfrak{v}}_{v}^M(\overline{m}^w) = \overline{\mathfrak{v}}_{v}^M(m^{w+1})=  &\ \overline{\mathfrak{v}}_{v}^M\big( (\overline{m}^w+b_{v,1}\cdot m^{w+1}) - (\overline{m}^w+b_{v,2}\cdot m^{w+1})  \big) \geq \nonumber\\
        \geq&\ \min\{\overline{\mathfrak{v}}_{v}^M(\overline{m}^w+b_{v,1}\cdot m^{w+1}), \overline{\mathfrak{v}}_{v}^M(\overline{m}^w+b_{v,2}\cdot m^{w+1}) \} > -p_v,
    \end{align*}
    a contradiction. Hence, since $k$ is infinite there exists $b\in k$ with $\overline{m}^w+b\cdot m^{w+1}$ satisfying the conditions (\ref{eq:m^wprop}) with index $w+1$, thus, by induction, we get the desired $m_p$.
\end{proof}

\begin{example}
Using the modified collection of extended valuations $\overline{\mathcal{D}}$ of Remark \ref{rem:Dbar}, there exists a `canonical' generalized local minimum point: $0\in (\mathbb{Z}_{\geq 0})^r$. Indeed, this is a consequence of the first step of the modification procedure to obtain $\overline{\mathcal{D}}$ from $\mathcal{D}$.
\end{example}

From the finite $k$-codimensionality of $N$ in $M$ and Lemma \ref{lem:existsm_p} we already know that there are only finitely many generalized local minimum points. 
Moreover, we know that there exists $d \in (\mathbb{Z}_{\geq0})^r$ such that $\mathcal{F}_{\overline{\mathcal{D}}}^M(-d)=M$ and all the lattice homological information is already present on the rectangle $R(0, d)$ (cf. Theorem \ref{th:properties} \textit{(a)} and subsection \ref{ss:ddep}). Therefore we have the following:
 
\begin{corollary}\label{cor:kisebb}
If  $p$ is a generalized local minimum point then $0\leq p\leq d$.
\end{corollary}

\begin{lemma} \label{szimmutak}
    Let $p \in (\bZ_{\geq 0})^r$ be a generalized local minimum point and $v\in\mathcal{V}$. Then
    \begin{equation*}
    w(\ell) < w(\ell + e_v) \ \ \text{ implies  } \ \ w(p-e_v-\ell) > w(p-\ell)
    \end{equation*}
    for any lattice point $\ell\in (\bZ_{\geq 0})^r$ satisfying  $\ell+e_v\leq p$. Even more, in this case 
    \begin{equation*}
        w(p-e_v-\ell) - w(p-\ell) \geq  w(\ell + e_v)-w(\ell).
    \end{equation*}
\end{lemma}

\begin{proof}
    By the Combinatorial Duality Property, the inequality $w(\ell) < w(\ell+e_v)$ is equivalent with 
    $$\dim_k\big(\mathcal{F}_{\overline{\mathcal{D}}}(\ell)/\mathcal{F}_{\overline{\mathcal{D}}}(\ell+e_v)\big)=w(\ell+e_v) - w(\ell) >0.$$
    Hence, for $g:= w(\ell+e_v) - w(\ell)$, there exists a $g$-dimensional complementary complex subspace $\mathcal{G}$ of $\mathcal{F}_{\overline{\mathcal{D}}}(\ell+e_v)$  in $\mathcal{F}_{\overline{\mathcal{D}}}(\ell)$, 
    i.e., $\mathcal{G} \oplus \mathcal{F}_{\overline{\mathcal{D}}}(\ell+e_v) = \mathcal{F}_{\overline{\mathcal{D}}}(\ell)$, with every $f \in \mathcal{G}\setminus\{ 0\}$ satisfying $\overline{\mathfrak{v}}(f)\geq \ell$ and $\overline{\mathfrak{v}}_v(f) = \ell_v$.

    But then, if we choose $m_p$ as in Lemma \ref{lem:existsm_p}, then for every $f \in \mathcal{G}\setminus \{0\}$ we have 
    \begin{center}
        $f \cdot m_p \in \mathcal{F}_{\overline{\mathcal{D}}}^M(-p+\ell)$ with $\overline{\mathfrak{v}}_v^M(f \cdot m_p )=-p_v + \ell_v$.
    \end{center} 
    Notice that this implies that the multiplication-with-$m_p$ map $\mathcal{G} \rightarrow \mathcal{F}_{\overline{\mathcal{D}}}^M(-p+\ell), \ f \mapsto f \cdot m_p$ is linear and injective, as $\overline{\mathfrak{v}}_v^M(0) = \infty$. Hence $\mathcal{F}_{\overline{\mathcal{D}}}^M(-p+\ell+e_v)$ has codimension at least $g$ in $\mathcal{F}_{\overline{\mathcal{D}}}^M(-p+\ell)$, thus, by the Combinatorial Duality Property, 
    \begin{equation*}
        \dim_k\big(\mathcal{F}_{\overline{\mathcal{D}}}^M(-p+\ell)/\mathcal{F}_{\mathcal{D}}^M(-p+\ell+e_v)\big) = w(p-e_v-\ell) - w(p-\ell) \geq g. \qedhere
    \end{equation*}
\end{proof}

\begin{remark}\label{rem:nononposthm} The previous Lemma \ref{szimmutak} and its proof cannot be directly generalized to the setting of subsection \ref{ss:MOSTGEN}, where we used a `multiplication' operation $\bullet:M_1\times M_{-1} \rightarrow M_0$ with possibly different sources and target. In that case, the existence of a linear subspace of $M_{1}$ analogous to $\mathcal{G}$ in the previous proof, would not necessarily imply anything similar in $M_{-1}$. This is causing that the Nonpositivity Theorem \ref{th:upperbound} cannot be directly generalized to that setting.
\end{remark}

\begin{proposition}\label{cl:glm}
Let $p\in (\mathbb{Z}_{\geq 0})^r$ be a generalized local minimum point and $\ell\in R(p/2, p) \cap\mathbb{Z}^r$ a lattice point (i.e., $\ell \geq p-\ell$). Then $w(\ell) \leq w(p-\ell)$. Specifically, if $\ell = p$ is a generalized local minimum point, then $w(p) \leq w(0) =0$.
\end{proposition}

\begin{proof}
Let us choose an increasing path $\gamma = \{x^j\}_{j=0}^t$ from $p-\ell$ to $\ell$ (in the sense of \ref{9SSP} (iii)), with $t=|p|-2|\ell|$, and define the dual path $\overline{\gamma}=\{ \overline{x}^{j}\}_{j=0}^t$ given by $\overline{x}^j := p-x^{t-j}$. Clearly, this is also an increasing path going from $\overline{x}^0=p - x^t= p-\ell$ to $\overline{x}^t=p - x^0=\ell$. By telescoping summation we get that
    \begin{align*}
        w(\ell)-w(p-\ell) &= \frac{1}{2}\bigg( \sum_{j=0}^{t-1}\big( w(x^{j+1})-w(x^j)\big) + \sum_{j=0}^{t-1}\big( w(\overline{x}^{t-j})-w(\overline{x}^{t-j-1})\big)\bigg).
    \end{align*}
    Now, by Lemma \ref{szimmutak}, if $w(x^{j+1})-w(x^j)>0$, then $w(\overline{x}^{t-j})-w(\overline{x}^{t-j-1})\leq -(w(x^{j+1})-w(x^j))$ (and, similarly, if $w(\overline{x}^{t-j})-w(\overline{x}^{t-j-1})>0$, then $w(x^{j+1})-w(x^j)\leq - (w(\overline{x}^{t-j})-w(\overline{x}^{t-j-1}))$), whence the total sum in the right hand side (hence the difference $w(\ell)-w(p-\ell)$ itself) cannot be positive.
    \end{proof}

The following is also true as in \cite[Theorem 3.2.10]{KNS2}:

\begin{theorem}\label{lem:w=0glm}
If $p \in (\bZ_{\geq 0})^r$ is a generalized local minimum point with $w(p)=0$, then
    \begin{equation}\label{eq:sym}
        w(p-\ell) = w(\ell) \hbox{ for any lattice point } \ell \in R(0, p) \cap \bZ^r,
    \end{equation}
    i.e.,
    the weight function restricted to the rectangle $R(0,p)$ is symmetric.
\end{theorem}

Let us now consider the slightly stricter setting of Theorem \ref{th:contr} and use also the modification procedure of Remark \ref{rem:Dbar}. I.e., consider the weight function $\overline{w}:=w_{\overline{\mathcal{D}}}$ corresponding to the CDP realization $\overline{\mathcal{D}}=\{\overline{\mathfrak{v}}_1, \ldots, \overline{\mathfrak{v}}_r\}$ of $N \leq M$ satisfying $\overline{\mathfrak{v}}_v^M(M)\subset 2\mathbb{Z}, \ \overline{\mathfrak{v}}_v^M(N)=0$ and $\{f \in \mathcal{O}\,:\,\overline{\mathfrak{v}}_v(f)\geq 1\}\equiv\mathfrak{p}$ for all $v \in \mathcal{V}$. Then we get the following corollary of Proposition \ref{cl:glm}:

\begin{corollary}\label{cor:wl+cube}
    Let $\ell \geq 0$ be a lattice point and consider a  cube $(\ell, J) \subset (\bR_{\geq 0})^r$ (for some index set $J \subset \mathcal{V}$). Suppose that 
        $\ell^+:=\ell + e_J$ is a generalized local minimum point.
    Then $\overline{w}(\ell) \leq \mathfrak{d}$.
\end{corollary}

\begin{proof}
    For any $\ell$ with $\ell \geq \ell^+-\ell=e_J$ 
    this follows from Proposition \ref{cl:glm} and Remark \ref{rem:w(R(0,E)<kappa}. In every other case we must have some index $v \in J$ such that $\ell_v=0$. But then $\ell^+_v=1$ odd with $\ell^+$ a generalized local minimum point, which contradicts Lemma \ref{lem:existsm_p}, since $\overline{\mathfrak{v}}_v^M(M) \subset 2\mathbb{Z}$ by our modification (\ref{eq:modifiedD}).
\end{proof}

\subsection{The M-vertices of cubes}\label{ss:6.2}\,

For cubes we use several notations, namely $(\ell,J)$, or $R(\ell, \ell+e_J)$, where $e_J:=\sum_{v\in J}e_v$. In this way we basically distinguish the vertices $\ell$ and $\ell+e_J$. However, in certain discussions it is convenient to distinguish some other vertices, e.g., a  vertex $\ell+e_{K}$ ($K\subset J$)  identified by having the maximal weight
$w(\ell+e_{K})=\max\{w(\ell')\,:\, \ell' \ \text{is a vertex of } \ (\ell,J)\}$.
Then we can rename the cube $(\ell, J)$, viewed from $\ell+e_K$ in `positive and negative' coordinate directions, see below.

\begin{nota}
For a lattice point $\ell \in \bZ^r$ and index sets $J^+, J^- \subset \mathcal{V}$ with $ J^+ \cap J^- = \emptyset$, let $(\ell, J^+, J^-)$ denote the cube with vertices $ \ell + e_{K^+} - e_{ K^-}$, where $K^+$ and  $K^-$ run over all subsets of $ J^+$ and $J^-$ respectively (i.e.,  $(\ell-e_{J^-}, J^+ \cup J^-)$ or $R(\ell-e_{J^-}, \ell+e_{J^+})$ using the previous notations). In this case the dimension of the cube is $q=\big|\{J^+ \cup J^-\}\big|$.
\end{nota}

\begin{define}\label{def:goodstarcube}
     The lattice point  $\ell$  is called an  \textit{`M-vertex'} of the cube  $\square$, if it is a vertex of it, i.e.,  $\square=(\ell, J^+, J^-)$ for some $J^+, J^- \subset \mathcal{V}$, and the weight function satisfies the following: 
    \begin{itemize}
        \item[(i)] $w(\ell+e_v) \leq w(\ell)$ for all $v \in J^+$,  and
        \item[(ii)] $w(\ell-e_v) \leq w(\ell)$ for all $v \in J^-$
    \end{itemize}
    (i.e., its weight is no less than any of its neighbours').
\end{define}

A lattice point having maximal weight in a cube is an $M$-vertex of that cube (hence a cube can have multiple M-vertices). The lemma below states the converse in the following combinatorial setting:

\begin{lemma} \label{lem:goodstarcube}
 Consider a weight function with form $w=h+h^\circ - h^{\circ}(0)$ as in (\ref{eq:w_0def}) on a cube $(\ell, J^+, J^-) \subset \mathbb{R}^r$, where $h, h^\circ: \mathbb{Z}^r \cap (\ell, J^+, J^-) \rightarrow \mathbb{Z}$ are an increasing and, respectively, 
    a decreasing function as in paragraph \ref{bek:comblattice}, both satisfying the matroid rank inequality and the pair $(h, h^{\circ})$ satisfying the CDP.\\
    Now, if $\ell$ is an M-vertex of $(\ell, J^+, J^-)$ then $w(\ell') \leq w(\ell)$
     for any vertex $\ell'$ of $(\ell, J^+, J^-)$.
\end{lemma}

\begin{proof}
        (I) \ Part (i) of Definition \ref{def:goodstarcube}
        is equivalent with $h(\ell)=h(\ell+e_v)$ for any  $v \in J^+$. Then induction on the cardinality $|K^+|$, for $K^+\subset J^+$, and
            inequality (\ref{eq:matroid}) imply $h(\ell)=h(\ell+e_{J^+})$.

        (II) \   Similarly, (ii) reads as
         $h^\circ(\ell-e_v)=h^\circ(\ell)-1$ for any  $v \in J^-$. Then, induction on the cardinality $|K^-|$, for
         $K^-\subset J^-$, and
            inequality (\ref{eq:matroid}) show that $h^\circ(\ell-e_{J^-})=h^\circ(\ell)$.

          (III) \
        Now for any vertex $\ell'$ of  $(\ell, J^+, J^-)$, we have $\ell-e_{J^-} \leq \ell' \leq \ell+e_{J^+}$. Firstly, by the increasing property of $h$, we have $h(\ell-e_{J^-}) \leq h(\ell') \leq h(\ell+e_{J^+})= h(\ell)$. Secondly, by the decreasing property of $h^\circ$, we have $h^{\circ}(\ell+e_{J^+}) \leq h^{\circ}(\ell') \leq h^{\circ}(\ell-e_{J^-})= h^{\circ}(\ell)$.
        Hence, by the definition of $w$, we have $w(\ell') \leq w(\ell)$.
\end{proof}

Notice, that the assumptions of Lemma \ref{lem:goodstarcube} are always satisfied when using the construction of section \ref{s:4} (cf. Observation \ref{obs:hhMmatroid}), more specifically, in the setting of Theorem \ref{th:contr} as well.

\begin{remark}
Notice, that Lemma \ref{lem:goodstarcube} is no longer true if we only require $w$ to satisfy the matroid rank inequality.
\end{remark}

\subsection{The deformation theorem}\,

In this subsection we present the proof of the Deformation Theorem and show how it implies Theorem \ref{th:contr}. Throughout this subsection we will use exclusively the modified realization $\overline{\mathcal{D}}$ from Remark \ref{rem:Dbar} and, hence, (without the danger of confusion) we will denote the corresponding height and weight functions and cubical complexes simply by $\mathfrak{h}, \ \mathfrak{h}^\circ, \ w$ and $S_n$ ($n\in \mathbb{Z}$). Recall, that due to the Independence Theorem \ref{th:IndepModTop}, these $S_n$-spaces are homotopy equivalent to those corresponding to any other CDP realization of the submodule $N \leq M$.

\begin{theorem}[Deformation Theorem] \label{th:DefThm}
    In the setting of Theorem \ref{th:contr}  for any integer $n\in \mathbb{Z}$, with $n \geq \mathfrak{d}:=\dim_k(\mathcal{O}/\mathfrak{p})$,  and lattice point $d \geq d_{\overline{\mathcal{D}}}$ (cf. subsection \ref{ss:ddep}) there exists a strong deformation retraction
    \begin{equation*}
        R(0, d) \searrow S_n \cap R(0, d).
    \end{equation*}
In particular,   $ R(0, d) $ and
$ S_n \cap R(0, d)$ are homotopy equivalent and, thus, both are contractible.
\end{theorem}

This clearly implies the Nonpositivity Theorem:

\begin{proof}[Proof of Theorem \ref{th:contr}]
    Theorems \ref{th:properties} \textit{(a)}  and the Deformation Theorem \ref{th:DefThm} imply that for any $n \geq \mathfrak{d}$, the cubical complex $S_n$ is contractible. 

    Let us also check that if $N \neq M$, then $S_{\mathfrak{d}-1}$ is not connected. Since we work over the finite rectangle $R(0, d)$, the lattice point $0$, by Remark \ref{rem:w(R(0,E)<kappa}, gives a distinct connected component of $S_{\mathfrak{d}-1}$ --- all of its neighbours have weight $\mathfrak{d}$. Thus, it is enough to show the existence of another (nonzero) lattice point $\ell \in (\mathbb{Z}_{\geq 0})^r$ with weight $w( \ell) < \mathfrak{d}$. For this, notice that $\{0\}$ is, in fact, a distinct contractible component of every $S_n$ with $0 \leq n \leq \mathfrak{d}-1$. Now, if we supposed indirectly that there was no other lattice point with weight strictly less than $\mathfrak{d}$, then we would get that 
    \begin{equation*}
        S_n \text{ is } \begin{cases}
            \text{contractible for every } n \geq \mathfrak{d},  & \text{ by the first part of the theorem;} \\
            \{0\}, \text{ hence contractible for every } 0 \leq n < \mathfrak{d}, & \text{ by the indirect assumption;} \\
            \text{empty for every }n<0, & \text{ by the indirect assumption.}
        \end{cases} 
        \end{equation*}
        But this would contradict the fact that $eu(\mathbb{H}_*(N \hookrightarrow M))=\dim_k(M/N) \neq 0$, hence, in fact, $S_{\mathfrak{d}-1}$ cannot be connected.
\end{proof}

\begin{proof}[Proof of the Deformation Theorem \ref{th:DefThm}]

We will give a sequence of strong deformation retractions
\begin{equation*}
\resizebox{15cm}{!}{
    $R(0, d) \overset{r_{j}}{\searrow} S_{j-1} \cap  R(0, d) \overset{r_{j-1}}{\searrow} \ldots S_{n} \cap R(0, d) \overset{r_n}{\searrow } S_{n-1} \cap R(0, d) \ldots  \overset{r_{\mathfrak{d}+2}}{\searrow} S_{\mathfrak{d}+1} \cap R(0, d) \overset{r_{\mathfrak{d}+1}}{\searrow } S_{\mathfrak{d}} \cap R(0, d),$}
\end{equation*}
where $j=\max\{ w(\ell) \,: \, \ell \in R(0, d) \cap \bZ^r\}$ (thus
$R(0, d) \subset S_j$).

Even more, each $r_n$ above will itself be the composition of strong deformation retractions, one for each lattice point in $(S_{n} \setminus S_{n-1}) \cap R(0, d)\cap \mathbb{Z}^r$. Let us denote this set by
\begin{equation}\label{eq:W_n}
    W_{n}:=\{ \ell \in \bZ^r\,: \,\ell \in  (S_{n} \setminus S_{n-1}) \cap R(0, d)\} = \{ \ell \in \bZ^r \cap R(0, d)\,: \, w(\ell) = n\}.
\end{equation}
Also, let us order this set $W_n=\{ \ell^1, \ldots, \ell^{|W_n|}\}$ such that for all $i \leq j$ we have $|\ell^{i}| \leq |\ell^j|$. Denote the cubical complex spanned by $(S_{n-1} \cap R(0, d)) \cup \{\ell^1, \ell^2, \ldots, \ell^m\}$ by $R_n^m$ (i.e., the union of those closed cubes, whose every vertex is contained in the set $(S_{n-1} \cap R(0, d) \cap \bZ^r) \cup \{\ell^1, \ell^2, \ldots, \ell^m\}$). Then $r_n$ will be constructed as the composition 
\begin{equation*}
S_{n} \cap  R(0, d) = R_n^{|W_n|} \overset{r^{|W_n|}_{n}}{\searrow} R_n^{|W_n|-1} \overset{r^{|W_n|-1}_{n}}{\searrow} \ldots R_n^m \overset{r^m_n}{\searrow } R_n^{m-1} \ldots  \overset{r^2_n}{\searrow} R_n^1 \overset{r_n^1}{\searrow } R_n^0= S_{n-1} \cap R(0, d)
\end{equation*}
of strong deformation retractions, where $r_n^m$ shall eliminate the lattice point $\ell^m$ and every closed cube containing it from $R_n^m$. 

Let us denote the closed cubical complex consisting of those cubes which have $\ell^m$ as their vertex
(the `closed star' of $\ell^m$) by
\begin{equation}\label{eq:St_ldef}
    St_{\ell^m}:= \bigcup_{\substack{\square_q \subset R_n^m \\ \ell^m \text{ vertex of }\square_q}}\square_q.\end{equation}
    (More precisely, we want to use every cube contained in the topological realization of this space, so as a cubical complex $St_{\ell^m}=\bigcup_{\square_{q'}' \subset \square_q \subset R^m_n, \ \ell^m \in \square_q} \square_{q'}'$.) Then $ St_{\ell^m} \cup R_n^{m-1} = R_n^m$.
Note that the `open star' of $\ell^m $ is 
\begin{center}
    $St_{\ell^m}^{\circ}= \bigcup \{\, \square_q^{\circ} \ : \ \square_q \subset R_n^m \text{ and } \ell^m \text{ is a} $ $\text{vertex of }\square_q \, \}$,
\end{center}
with $ St_{\ell^m}^{\circ} =  R_n^m \setminus R_n^{m-1} $. 

Let us denote by $T_{\ell^m}$ the subset of $St_{\ell^m}\subset\mathbb{R}^r$ consisting of those points which are contained in such closed faces of $St_{\ell^m}$, which do not have $\ell^m$ as a vertex (i.e., $T_{\ell^m}=\bigcup_{ \square_{q'}' \subset St_{\ell^m}, \ \ell^m \notin \square_{q'}'} \square_{q'}'$). Then clearly $St_{\ell^m} \cap R_n^{m-1} = T_{\ell^m}$. Therefore, to obtain $r_n^m: R_n^m \searrow R_n^{m-1}$ we  need a strong deformation retraction $St_{\ell^m} \searrow T_{\ell^m}$, which then could be extended to $R_n^m$  with the identity map on $R_n^{m-1}$. 
A construction for this strong deformation retraction was presented in \cite[subsection 6.5]{KNS2}. In order to being able to apply it here as well, we need the following observation:

\begin{prop}\label{prop:Stesell-Ei}
    For every lattice point $\ell^m$, with $w(\ell^m)=n>\mathfrak{d}$, there exists a `good' coordinate direction $v$, meaning that
    \begin{equation}\label{eq:Stesell-Ei} 
    St_{\ell^m}=\bigcup_{\substack{(\ell^m, J^+, J^-) \subset R_n^m \\ v \in J^-}}(\ell^m, J^+, J^-)
\end{equation}
 as a topological subspace in $\mathbb{R}^r$ (i.e., as a topological subspace $St_{\ell^m}$ can be written as the union of cubes $(\ell^m, J^+, J^-)$ such that $v \in J^-$; in particular,  $\ell^m+e_v \notin St_{\ell^m}$).
\end{prop}

\begin{proof}
As $R_n^m \subset S_n$ and $w(\ell^m) =n$, we get that
$\ell^m$ is an M-vertex of every 
 $q$-cube $\square_q \subset St_{\ell^m}$ containing $\ell^m$ (compare with Lemma \ref{lem:goodstarcube}).

 Set $I^+:=\{ v \in \mathcal{V}\,: \, w(\ell^m + e_v) < w(\ell^m)\}$. We will argue that in the absence of a good coordinate direction predicted  by the theorem, we would have at  $\ell^+:= \ell^m + e_{I^+}$ a gene\-ralized local minimum point and could use Corollary \ref{cor:wl+cube} for the cube $(\ell^m, I^+)$ to get $w(\ell^m) \leq \mathfrak{d}$, which contradicts the assumption.

 First notice, that, by the matroid rank inequality valid for the weight function $w$ (cf. Lemma \ref{lem:w0matroid}), from $w(\ell^m + e_v) < w(\ell^m)$ we get
 \begin{equation}\label{eq:I}
 w(\ell^+)<w(\ell^+-e_v) \ \ \mbox{ for every index $v \in I^+$.}
 \end{equation}
 On the other hand, if we set $I^-:=\{ v \in \mathcal{V}\,: \, w(\ell^m - e_v) > w(\ell^m)\}$, then, using again the matroid rank inequality,
 we get
 \begin{equation}\label{eq:II}
 w(\ell^+)<w(\ell^+-e_v) \ \ \mbox{ for every index $v \in I^-$. }
 \end{equation}
 Now $I^+\cup I^-=\mathcal{V}$ cannot happen since in that case $\ell^+$ would be a generalized local minimum point with $w(\ell^m)> \mathfrak{d}$, contradicting Corollary \ref{cor:wl+cube} used for the cube $(\ell, I^+)$. On the other hand, for any
  $v \in \mathcal{V} \setminus (I^+ \cup I^-)$  we automatically have
    \begin{equation*}
        w(\ell-e_v) \leq w(\ell) \leq w(\ell+e_v).
    \end{equation*}
    Therefore, by the order in which we eliminate the higher-than-$\mathfrak{d}$ weight vertices (discussed right after (\ref{eq:W_n})) we see that $\ell +e_v \notin R_n^m$, hence for every cube $(\ell^m, J^+, J^-) \subset R_n^m$ we have $v \notin J^+$. 
    
    To finish the proof, it is enough to find among $\mathcal{V} \setminus (I^+ \cup I^-)$ a `good' coordinate $v$ satisfying that for each $(\ell^m, J^+, J^-) \subset R_n^m$, the (possibly) larger $(\ell^m, J^+, J^-\cup\{v\})$ is also contained in $R_n^m$.
    Assume on the contrary that for every index  $v \in \mathcal{V} \setminus (I^+ \cup I^-)$
    there exists a  cube $\square_{q,v}=(\ell^m, J^+, J^-) \subset R_n^m$ with $v \notin J^+ \cup J^-$, admitting $\ell$ as an M-vertex and having the property that $(\ell^m, J^+, J^- \cup \{v\}) \not\subset R_n^m$. This implies, that there exists a  vertex $\ell' \in \square_{q,v}$ such that $\ell'-e_v \notin R_n^m$, i.e.,
     \begin{equation}\label{eq:6.3}
     w(\ell'-e_v)> w(\ell').\end{equation} 
    Indeed, as $S_{n-1} \subset R_{n}^m$, we know that $w(\ell'-e_v) \geq n$, and, by $\square_{q,v}\subset R_n^m$, $w(\ell') \leq n$. Also, if $w(\ell') = n$, once again by the ordering of the $w=n$ lattice points, $|\ell^m| \geq|\ell'| > |\ell'-e_i|$, and thus, if $w(\ell'-e_v) = n$, then $\ell'-e_v$ would be contained in $R_n^m$. Therefore, in the sequence $w(\ell'-e_v) \geq n \geq w(\ell')$ at least one of the inequalities must be strict.

    On the other hand, $J^+ \subset I^+$ (again by the ordering of the weight $=n$ lattice points),
    hence we have 
    $\ell'\leq \ell^m+e_{J^+}\leq \ell^m+e_{I^+}= \ell^+$.
    Then, applying the matroid rank inequality to (\ref{eq:6.3}),
\begin{equation}\label{eq:III}
 w(\ell')<w(\ell'-e_v)\ \ \Rightarrow\  \ w(\ell^+)<w(\ell^+-e_v).
    \end{equation}
    But if this is true for every remaining coordinate $v \in \mathcal{V} \setminus (I^+ \cup I^-)$, then via (\ref{eq:I}), (\ref{eq:II}) and (\ref{eq:III}) the lattice point $\ell^+$ is a generalized local minimum point, a contradiction. Therefore, there must exist a `good' coordinate direction $v \in \mathcal{V} \setminus (I^+ \cup I^-)$, which  satisfies (\ref{eq:Stesell-Ei}).
Therefore, as a topological space $St_{\ell^m}$ can be given as a union of closed cubes containing both $\ell^m$ and $\ell^m-e_v$.
\end{proof}

Let us now turn back to the proof of the Deformation Theorem \ref{th:DefThm}. By Proposition \ref{prop:Stesell-Ei}, for every integer $n>\mathfrak{d}$ and $1 \leq m \leq |W_n|$ the space $St_{\ell^m}$ can be described as the union of closed cubes containing both $\ell^{m}$ and $\ell^{m}-e_v$ for the `good' direction $v \in \mathcal{V}$. Therefore, we can apply the construction of \cite[Claim 6.5.8]{KNS2} to obtain the desired strong deformation retraction $St_{\ell^m} \searrow T_{\ell^m}$ through simple `homotheties' with center $\ell^m+e_v$ for every $n>\mathfrak{d}$, $1 \leq m \leq |W_n|$. Their composition will give the deformations $R(0, d) \searrow S_n \cap R(0, d)$ announced by the Deformation Theorem.
 \end{proof}

\subsection{Nonpositivity Theorem for direct products of integrally reduced local Artin algebras}\label{ss:KünnethNonposproof}\,

In this subsection we will prove the Nonpositivity Theorem \ref{th:NonposArtin} for direct products of integrally reduced local Artin algebras. We will use the notions and notations of section \ref{ss:ARTIN}. Let us recall the statement:

\begin{theorem}[= Theorem \ref{th:NonposArtin}]
    Let $k$ be an algebraically closed field and $A$ an integrally reduced Artin $k$-algebra. Then $A$ can be written in a unique way as a direct product of integrally reduced \emph{local} Artin $k$-algebras $\prod_{i=1}^l A_i$. In this case for every $n \geq l$ the $S_n$ space is contractible, thus the weight grading of the \emph{reduced} symmetric lattice homology is not lower than $2-2l$, i.e., $\mathbb{SH}_{{\rm red}, q, -2n}(A)=0$ for all $q\geq 0$ and $n > l-1$. Moreover, if $A$ is nonzero, then so is $\mathbb{SH}_{{\rm red}, l-1, 2-2l}$.
\end{theorem}

\begin{proof}
    We will discuss the case when $l=2$, i.e., $A=A' \times A''$, with $A'$ and $A''$ local and integrally reduced, for more summands one just needs to repeat this proof inductively. Consider some descriptions of the summands as $A'\cong \mathcal{O}'/\mathcal{I}'$ and $A''\cong \mathcal{O}''/\mathcal{I}''$, with $\mathcal{I}'$ and $\mathcal{I}''$ integrally closed, and choose CDP realizations $(\mathcal{D}', d_{\mathcal{I}'})$ for $\mathcal{I}'$ and $(\mathcal{D}'', d_{\mathcal{I}''})$ for $\mathcal{I}''$. Then the collection $\mathcal{D}^\times$ and lattice point $d_{\mathcal{I}'\times \mathcal{I}''}$ of Proposition \ref{prop:chainproduct} will give a CDP realization of $\mathcal{I}'\times\mathcal{I}''\triangleleft\mathcal{O}' \times \mathcal{O}''$ and, by Remark \ref{rem:Snproduct}, the $S_n$-spaces are as follows: $S_{n, \mathcal{D}^{\times}}=\bigcup_{i+j=n} S_{i, \mathcal{D}'} \times S_{j, \mathcal{D}''}$. Recall, that by the Independence Theorem \ref{th:IndepModTop} (using also the fact that $\mathbb{SH}_*(A) \cong \mathbb{SH}_*(\mathcal{I}' \times \mathcal{I}'' \triangleleft \mathcal{O}' \times \mathcal{O}'')\cong \mathbb{H}_*(\mathcal{I}' \times \mathcal{I}'' \hookrightarrow \mathcal{O}' \times \mathcal{O}'')$ also on the $S_n$-space level and part \textit{(a)} of Theorem \ref{th:properties}) these $S_n$-spaces are homotopy equivalent to those corresponding to any other CDP realization of $\mathcal{I}'\times \mathcal{I}''$, hence, we can use these distinguished ones for our proof.  

    First, we will prove that $S_{n, \mathcal{D}^{\times}}$ is contractible for any $n \geq 2$. Indeed, $\bigcup_{i+j=n} S_{i, \mathcal{D}'} \times S_{j, \mathcal{D}''}$ for $n\geq 2$ is homotopy equivalent to $S_{1, \mathcal{D}'} \times S_{n-1, \mathcal{D}''}$, since we can use the (proof of the) Deformation Theorem \ref{th:DefThm} (relying also on part \textit{(a)} of Theorem \ref{th:properties}) to obtain the strong deformation retractions:
    \begin{itemize}
        \item for any $1 \leq i \leq n-\min w_{\mathcal{D}''}$:
        \begin{equation*}
        S_{i, \mathcal{D}'}\times S_{n-i, \mathcal{D}''} \overset{r_{i,  \mathcal{D}'} \times {\rm id}}{\searrow} S_{i-1, \mathcal{D}'}\times S_{n-i, \mathcal{D}''} \overset{r_{i-1,  \mathcal{D}'} \times {\rm id}}{\searrow} \ldots \overset{r_{2,  \mathcal{D}'} \times {\rm id}}{\searrow} S_{1, \mathcal{D}'} \times S_{n-i, \mathcal{D}''} \subset S_{1, \mathcal{D}'} \times S_{n-1, \mathcal{D}''};
    \end{equation*}
    \item for any $n-1 \leq j \leq n-\min w_{\mathcal{D}'}$:
    \begin{equation*}
        S_{n-j, \mathcal{D}'}\times S_{j, \mathcal{D}''} \overset{ {\rm id} \times r_{j,  \mathcal{D}''} }{\searrow} S_{n-j, \mathcal{D}'}\times S_{j, \mathcal{D}''} \overset{ {\rm id} \times r_{j-1,  \mathcal{D}''} }{\searrow} \ldots \overset{ {\rm id} \times r_{n,  \mathcal{D}''} }{\searrow} S_{n-j, \mathcal{D}'} \times S_{n-1, \mathcal{D}''} \subset S_{1, \mathcal{D}'}\times S_{n-1, \mathcal{D}''}.
    \end{equation*}
    \end{itemize}
    Moreover, these maps coincide on the intersections, hence, we indeed get a strong deformation retraction $\bigcup_{i+j=n} S_{i, \mathcal{D}'} \times S_{j, \mathcal{D}''} \searrow S_{1, \mathcal{D}'}\times S_{n-1, \mathcal{D}''}$, with  contractible target (cf. Theorem \ref{th:propertiesideal} \textit{(d)}). 

    Secondly, we have to check that if $A= A'\times A''$ with $A', A'' \neq 0$ local, then $S_{1, \mathcal{D}^{\times}}$ has nontrivial first homology. Similarly to the previous case, by the Deformation Theorem \ref{th:DefThm}, $S_{1, \mathcal{D}^{\times}}$ is homotopy equivalent to $S_{0, \mathcal{D}'}\times S_{1, \mathcal{D}''}\cup S_{1, \mathcal{D}'}\times S_{0, \mathcal{D}''}$, since the other terms of form $S_{i, \mathcal{D}'}\times S_{1-i, \mathcal{D}''}$ can be deformation retracted into  either one of them. On the other hand, since the $S_1$-terms are contractible by the previously proved Nonpositivity Theorem (the specific version used here is Theorem \ref{th:propertiesideal} \textit{(d)}), the cubical complex $S_{1, \mathcal{D}^{\times}}$ is homotopy equivalent to the join of  $S_{0, \mathcal{D}'}$ and $S_{0, \mathcal{D}''}$, which are themselves non-connected (again by Theorem \ref{th:propertiesideal} \textit{(d)}), thus yield nontrivial $H_1(S_{1, \mathcal{D}^{\times}}, \mathbb{Z})$ (for a concrete computation see, e.g., Example \ref{ex:product}). 
    
    For the higher component number cases (i.e., when $l>2$) one has to use the same arguments multiple times inductively.  
\end{proof}

\section{Quasi-valuations and quasi-realizability}\label{ss:quasi}

In this section we will present yet another --- in some sense more general --- approach to the main results of this manuscript: the Independence Theorem and the well-definedness of the lattice homology of realizable submodules. Our starting point is the following deficiency in Definition \ref{def:REAL} of realizability: the fact that the submodule $N\leq M$ is realizable does not (necessarily) imply that $N/N$ is realizable in $M/N$ (see, e.g., the example in Remark \ref{rem:imageofintclosed}), 
even though the lattice homology $\mathbb{H}_*(N \hookrightarrow M)$ only depends on this quotient $M/N$. This weakness follows from the Definition \ref{def:gdv} of discrete valuations, which, by their properties, assign the value $\infty$ to every nilpotent element. Moreover, it is also unfavorable that,  although the output depends only on the quotient $M/N$, the construction of the lattice homology module requires the existence of valuations defined on the whole $M$. In order to remedy these issues, we introduce the notion of \emph{`quasi-valuations'}, which are, in some sense, quotiented\,/\,restricted discrete valuations, crafted to handle nilpotent and torsion elements as well.

\subsection{Quasi-valuations}\,

Consider the following setting: let $k$ be a field and let $\cO$ be a Noetherian $k$-algebra.

\begin{define}\label{def:Nd}
Let us define a commutative monoid structure on the set $\{0, 1, ..., d-1, d \} \subset \mathbb{N}$ (with the induced ordering and $d\in \mathbb{Z}_{>0}$) as follows:
\begin{center}
    \item $\forall a, b \in \{0, 1, ..., d-1, d \}: \  a \ast b = \begin{cases}
    a+b\ & \text{ if } a+b \leq d ;\\
    d\ & \text{ if } a+b \geq d.\\
    \end{cases}$
\end{center}
This operation clearly preserves the ordering. Denote this monoid with $\N_d$.
\end{define}

\begin{define}\label{def:quasi-val}
We call a function $\mathfrak{q}:\mathcal{O} \rightarrow \N_d$ a \emph{`discrete quasi-valuation'} if it has the following properties:
\begin{itemize}
    \item $\mathfrak{q}(0)=d$;
    \item $\mathfrak{q}(fg) = \mathfrak{q}(f) \ast \mathfrak{q}(g)$ for every $f, g \in \cO$;
    \item $\mathfrak{q}(af) = \mathfrak{q}(f)$ for every $f \in \cO, \ a \in k^{\ast}= k \setminus \{ 0 \}$.
    \item $\mathfrak{q}(f+g) \geq \min \{ \mathfrak{q}(f), \mathfrak{q}(g) \}$ for every $f, g \in \cO$.
\end{itemize}
 Given a quasi-valuation $\mathfrak{q}$ on $\mathcal{O}$ we can define the ideals 
     \begin{equation*}
    \cF_{\mathfrak{q}}(n) = \{ f \in \cO \,:\, \mathfrak{q}(f) \geq n \}\triangleleft\mathcal{O} \text{ for any }n \in \mathbb{N}_d.
\end{equation*}
\end{define}

Similarly as before, if $\mathfrak{q}(f) \neq \mathfrak{q}(g)$, then automatically $\mathfrak{q}(f+g) = \min \{ \mathfrak{q}(f), \mathfrak{q}(g) \}$.

\begin{example}\label{ex:levagottv}
    For any discrete valuation $\mathfrak{v}:\mathcal{O} \rightarrow \overline{\mathbb{N}}$ (in the sense of Definition \ref{def:gdv}) and any positive integer $d\in \mathbb{Z}_{>0}$, the composition $\mathfrak{v}_{\leq d}: \mathcal{O} \xrightarrow{\mathfrak{v}}\overline{\mathbb{N}}\xrightarrow{\min_d}\mathbb{N}_d$ is a quasi-valuation, where $\min_d:\overline{\mathbb{N}} \rightarrow \mathbb{N}_d$ is the natural forgetting map: $\min_d(n)=\min\{d, n\}$ (which is a commutative monoid morphism).

    Moreover, for any $n \in \mathbb{N}_d$ we have that $\mathcal{F}_{\mathfrak{v}_{\leq d}}(n)=\mathcal{F}_{\mathfrak{v}}(n) \triangleleft \mathcal{O}$ (cf. (\ref{eq:fvMl})).
\end{example}

\begin{example}\label{ex:quasirest}
    Let us consider a discrete valuation $\mathfrak{v}:\mathcal{O} \rightarrow \mathbb{N}\cup \{ \infty\}$ in the sense of Definition \ref{def:gdv} and consider some ideal $\mathcal{I} \triangleleft \mathcal{O}$ satisfying $\mathfrak{v}(\mathcal{I})\geq d$ (i.e., for all $f \in \mathcal{I}$ we have $\mathfrak{v}(f) \geq d$). Then the association 
    \begin{equation*}
        \mathfrak{q}: \mathcal{O}/\mathcal{I} \rightarrow \mathbb{N}_d, \ f+ \mathcal{I} \mapsto \begin{cases}
            \mathfrak{v}(f) & \text{ if $\mathfrak{v}(f) < d$}\\
            d & \text{ if $\mathfrak{v}(f) \geq d$}
        \end{cases}
    \end{equation*}
    is well-defined (by property (b) of Definition \ref{def:gdv}) and is, in fact, a discrete quasi-valuation. We will sometimes denote it by $\mathfrak{v}\big|_{\mathcal{O}/\mathcal{I}}$.

    Moreover, if we denote by $q:\mathcal{O} \rightarrow \mathcal{O}/\mathcal{I}$ the standard quotient map, then for any $n \in \mathbb{N}_d$ we have that $q^{-1} \big(\mathcal{F}_{\mathfrak{v}\big|_{\mathcal{O}/\mathcal{I}}}(n)\big)=\mathcal{F}_{\mathfrak{v}}(n)\triangleleft \mathcal{O}$. 
\end{example}

On the other hand, not all quasi-valuations $\mathfrak{q}:\mathcal{O}/\mathcal{I}\rightarrow \mathbb{N}_d$ can be extended to a discrete valuation on $\mathcal{O}$, see Example \ref{ex:quasinonext} later. Thus, the following question comes up naturally:

\begin{question} \label{q:whenextendquasi}
    Consider a finite codimensional ideal $\mathcal{I}$ in the polynomial ring $\mathcal{O}_m=k[x_1, \ldots, x_m]$ (respectively the convergent or formal power series ring)  and a discrete quasi-valuation $\mathfrak{q}:\mathcal{O} \rightarrow \mathbb{N}_d$ on $\mathcal{O}=\mathcal{O}_m/\mathcal{I}$. When can we extend this to a discrete valuation $\mathfrak{v}:\mathcal{O}_m \rightarrow \mathbb{N} \cup \{\infty\}$ (in the sense of Definition \ref{def:gdv}) satisfying $\mathfrak{v}(\mathcal{I})\geq d$ and $\mathfrak{q} = \mathfrak{v}\big\vert_{\mathcal{O}_m/\mathcal{I}}$?
\end{question}

One main advantage of quasi-valuations over valuations is the fact that over the finite dimensional quotient $\mathcal{O}/{\rm Ann}(M/N)$ they are finitely determined:

\begin{remark}\label{rem:constructionofq-val}
    Given a \textit{finite dimensional} $k$-algebra $\mathcal{O}$, any discrete quasi-valuation $\mathfrak{q}:\mathcal{O} \rightarrow \mathbb{N}_d$ can be specified by the following `finite' data: for given $\mathfrak{q}:\mathcal{O} \rightarrow \mathbb{N}_d$ there exists a $k$-vector space basis $\{f_1, f_2, \ldots, f_{o}\}$ of $\mathcal{O}$ (with $\dim_k \mathcal{O}=o$) with $\mathfrak{q}(f_i)=:q_i \in \{1, \ldots, d\}$  for all $1 \leq i \leq o$, such that 
    \begin{equation} \label{eq:quasivalwithbasis}
        \text{for any } f = \sum_{i=1}^{o}a_if_i \in \mathcal{O}  \text{ we have } \mathfrak{q}\big(\sum_{i=1}^{o}a_if_i\big)=\min\{q_i\,:\, a_i \in k\setminus\{0\}\}.
    \end{equation}
    Indeed, this basis can be obtained by subsequently extending along the increasing filtration
    \begin{center}
        $\mathcal{F}_{\mathfrak{q}}(d) \subset \mathcal{F}_{\mathfrak{q}}(d-1) \subset \ldots \subset \mathcal{F}_{\mathfrak{q}}(0)=\mathcal{O}$. 
    \end{center}
    
    For multiple quasi-valuations on the same finite dimensional algebra we can find a similar basis and assigned tuples respecting  identity (\ref{eq:quasivalwithbasis}) for every individual quasi-valuation via studying the respective multifiltration. An example for this `description with basis' construction was already given in Example \ref{ex:345534}.
\end{remark}
    On the other hand, (returning to the single quasi-valuation case) we can also describe when such basis and assigned values produces a quasi-valuation:
    
    \begin{prop}\label{prop:quasivalcond}
        Let $\mathcal{O}$ be a finite dimensional $k$-algebra and consider a $k$-basis $\{f_1, \ldots, f_o\}$ of  $\mathcal{O}$ and an assigned set of values $\{q_1, \ldots, q_o\}$, with $q_i \in \mathbb{N}_d$ for all $1\leq i \leq o$. 
    The mapping
    \begin{equation} \label{eq:quasiwithbasis}
        \mathfrak{q}:\mathcal{O} \rightarrow \mathbb{N}_d, f = \sum_{i=1}^{o}a_if_i\mapsto\min\{q_i\,:\, a_i \in k\setminus\{0\}\}
    \end{equation}
    gives a quasi-valuation if and only if the following compatibility condition is satisfied:
        \begin{equation*}
        \begin{split}
\big({\rm Span}_k(\{ f_i\, :\, q_i =a \})\setminus&\, \{0\}\big) \cdot\big( {\rm Span}_k(\{ f_i\, :\, q_i =b \})\setminus\{0\}\big) \subset \\ &\, \subset {\rm Span}_k(\{ f_i\, :\, q_i \geq a\ast b \})\setminus {\rm Span}_k(\{ f_i\, :\, q_i \geq a\ast b +1\})
\end{split}
        \end{equation*}
        \noindent for any $a, b \in \mathbb{N}_d$, i.e., for any $0 \neq f \in {\rm Span}_k(\{ f_i\, :\, q_i =a \})$ and $0 \neq f' \in {\rm Span}_k(\{ f_i\, :\, q_i =b \})$, the product $ff' \in {\rm Span}_k(\{ f_i\, :\, q_i \geq a\ast b \})\setminus {\rm Span}_k(\{ f_i\, :\, q_i \geq a\ast b +1\})$.
    \end{prop}

\begin{proof}
    Straightforward computation.
\end{proof}

The compatibility condition of Proposition \ref{prop:quasivalcond} clearly implies that $\mathfrak{q}(f_jf_{j'})=\mathfrak{q}(f_{j})*\mathfrak{q}(f_{j'})$ (or, equivalently, $f_jf_{j'} \in {\rm Span}_k(\{f_i\,:\,q_i \geq q_j\ast q_{j'}\}) \setminus {\rm Span}_k(\{f_i\,:\,q_i \geq q_j\ast q_{j'}+1\})$) for every pair $(f_j, f_{j'})$ of basis elements. On the other hand, this latter condition is not sufficient for $\mathfrak{q}$ to be a quasi-valuation:

\begin{nonexample} \label{ex:noneasyq-val}
    Consider the Artin $k$-algebra $\mathcal{O}:=k[x, y, z, t]/(\mathfrak{m}^3+(xt-zy,\, xz-yt))$, where $\mathfrak{m}$ denotes the maximal ideal $\mathfrak{m}=(x, y, z, t)$. (For simplicity we will use the same notation for a polynomial $p(x, y, z, t)\in k[x, y, z, t]$ and its class $[p(x, y, z, t)]=p(x, y, z, t) + (\mathfrak{m}^3+(xt-yz,\, xz-yt))$ in $\mathcal{O}$.) Consider the monomial $k$-basis $\mathcal{B}:=\{1, x, y, z, t, x^2, xy, xz=yt, xt=yz, y^2, z^2, zt, t^2\}$ of $\mathcal{O}$ and assign the ordinary degree values ($f_i=x^{a_x}y^{a_y}z^{a_z}t^{a_t} \rightsquigarrow a_x+a_y+a_z+a_t=q_i \in \mathbb{N}_3$) to these basis elements. Then, for any pair $f_j, f_{j'} \in \mathcal{B}$, we have $f_jf_{j'} \in {\rm Span}_k(\{f_i\,:\,q_i \geq q_j\ast q_{j'}\}) \setminus {\rm Span}_k(\{f_i\,:\,q_i \geq q_j\ast q_{j'}+1\})$, however, the formula $\mathfrak{q}: \mathcal{O} \rightarrow \mathbb{N}_3, \ f = \sum_{i=1}^{o}a_if_i \mapsto \min\{q_i\,:\, a_i \in k^\ast\}$ does not define a quasi-valuation on $\mathcal{O}$. Indeed, 
        $1+1=\mathfrak{q}(x+y)\ast\mathfrak{q}(z-t) \neq \mathfrak{q}((x+y)(z-t))=\mathfrak{q}(0)=3$.
\end{nonexample} 

Proposition \ref{prop:quasivalcond} can be best used to construct quasi-valuations when the set of assigned values $\{q_1, \ldots, q_o\}$ contains little to no repetitions, since then we can work with lower dimensional vector spaces to verify the compatibility condition. Using this construction we are now able to  show that quasi-valuations of a quotient algebra does not (necessarily) extend to discrete valuations of its dividend:

\begin{example}\label{ex:quasinonext}
    Let us consider the polynomial ring $k[x,y,z,t]$ and its finite dimensional quotient $\mathcal{O}:=k[x,y,z,t]/(x,y,z,t)^3$. (For simplicity we will use the same notation for $p(x, y, z, t)\in k[x, y, z, t]$ and its class $[p(x, y, z, t)]=p(x, y, z, t) + (x,y,z,t)^3$ in $\mathcal{O}$.) We construct a quasi-valuation $\mathfrak{q}:\mathcal{O}\rightarrow\mathbb{N}_{13}$ by Proposition \ref{prop:quasivalcond} using the following $k$-basis and associated values: 
    \begin{equation*}
\begin{array}{ccccccccccccccccc}
\{ & 1,  & x,  & y,  & z, & t, & x^2, & xy, & xz, & xt+yz, & y^2, & yt+z^2, & zt, & t^2, & xt-yz, &  yt-z^2 & \}\, \\
& \downarrow \,& \downarrow \, & \downarrow \,& \downarrow \,& \downarrow \,& \downarrow \,& \downarrow \,& \downarrow \,& \downarrow \,& \downarrow \,& \downarrow \,& \downarrow \,& \downarrow \,& \downarrow \,& \downarrow & \\
\{ &0, & 5, & 5, & 5, & 5, & 10, & 10, & 10, & 10, & 10, & 10, & 10, &  10, & 11, & 12 & \}.
\end{array}
 \end{equation*} 
 One can check by direct computation that these satisfy the compatibility condition of Proposition \ref{prop:quasivalcond}, thus $\mathfrak{q}:\mathcal{O}\rightarrow\mathbb{N}_{13}$ defined by formula (\ref{eq:quasiwithbasis}) is indeed a quasi-valuation.

On the other hand, there cannot exist a discrete valuation $\mathfrak{v}:k[x,y,z,t]\rightarrow \overline{\mathbb{N}}$ such that its restriction $\mathfrak{v}\big|_{k[x,y,z,t]/(x,y,z,t)^3}$ with $d=13$ in the sense of Example \ref{ex:quasirest} is $\mathfrak{q}$. Indeed, such a valuation would have
\begin{itemize}
    \item $\mathfrak{v}(t(xz-y^2))=\mathfrak{v}(t) + \mathfrak{v}(xz-y^2)=\mathfrak{q}(t) + \mathfrak{q}(xz-y^2)=15$;
    \item $\mathfrak{v}(z(xt-yz))=\mathfrak{v}(z) + \mathfrak{v}(xt-yz)=\mathfrak{q}(z) + \mathfrak{q}(xt-yz)=16$;
    \item $\mathfrak{v}(y(z^2-yt))=\mathfrak{v}(y) + \mathfrak{v}(z^2-yt)=\mathfrak{q}(y) + \mathfrak{q}(z^2-yt)=17$; and thus
    \item $15 = \mathfrak{v}(t(xz-y^2))=\mathfrak{v}(z(xt-yz)+y(z^2-yt))<\min\{ \mathfrak{v}(z(xt-yz)), \mathfrak{v}(y(z^2-yt))\}=16$,
\end{itemize}
which contradicts condition (b) of Definition \ref{def:gdv} of discrete valuations.
\end{example}

\subsection{Extended quasi-valuations}\,

We also define the extensions of quasi-valuations to modules analogously to the case of extended discrete valuations (cf. Definition \ref{def:edv}). Let $\mathcal{O}$ be a Noetherian $k$ algebra and $M$ be a finitely generated $\mathcal{O}$-module.
\begin{define}\label{def:extquasi}
    We call a pair $(\mathfrak{q}, \mathfrak{q}^M)$ an \textit{`extended discrete quasi-valuation'} if the following properties hold:
    \begin{itemize}
        \item $\mathfrak{q}:\cO \rightarrow \N_d$ is a discrete quasi-valuation in the sense of Definition \ref{def:quasi-val} for some $d \in \mathbb{Z}_{>0}$;
        \item $\mathfrak{q}^M:M \rightarrow \mathbb{M}_d:=\{ -d, -d+1,\ldots,  -1, 0\}$ satisfies the following properties:
        \begin{itemize}
    \item[(a)] $\mathfrak{q}^M(fm) = \mathfrak{q}(f) \ast \mathfrak{q}^M(m)$ \ for all  $f \in \cO, m \in M$,
      where, similarly as above,  
      \begin{equation} \label{eq:quasiaction}
     \forall a \in \mathbb{N}_d,\ b\in \mathbb{M}_d: \  a \ast b = \begin{cases}
    a+b\ & \text{ if } a+b \leq 0 ;\\
    0\ & \text{ if } a+b \geq 0.
    \end{cases}
\end{equation}
    \item[(b)] $\mathfrak{q}^M(m+m') \geq \min \{ \mathfrak{q}^M(m), \mathfrak{q}^M(m') \}$ \ for all $m, m' \in  M$;
\end{itemize}
    \end{itemize}
    Given an extended quasi-valuation $(\mathfrak{q}, \mathfrak{q}^M)$ we can define the submodules 
     \begin{equation*}
    \cF_{\mathfrak{q}}(n) = \{ f \in \cO \,:\, \mathfrak{q}(f) \geq n \}\triangleleft\mathcal{O} \ \ \text{ and } \ \ \cF^M_{\mathfrak{q}}(-n) = \{ m \in M \,:\, \mathfrak{q}^M(m) \geq -n \}\leq M \text{ as in (\ref{eq:fvMl})},
\end{equation*}
 for any $n \in \{ 0, 1, \ldots, d\}$. 
\end{define}

\begin{remark}
    Similarly as before, if $\mathfrak{q}^M(m) \neq \mathfrak{q}^M(m')$, then $\mathfrak{q}^M(m+m') = \min \{ \mathfrak{q}^M(m), \mathfrak{q}^M(m') \}$.
\end{remark}

\begin{remark}
    The convention of Definition \ref{def:extquasi} regarding the choice that extensions of quasi-valuations map to $\mathbb{M}_d=\{-d, -d+1, \ldots, -1,0\}$ fits perfectly with our construction and geometric pictures, yet, it is not so natural as mapping to $\mathbb{N}_d=\{0, 1, \ldots, d-1, d\}$ would be. For other applications one might find this latter convention to be more convenient.
\end{remark}

\begin{example}\label{ex:levagottv^M}
    For any extended discrete valuation $(\mathfrak{v}, \mathfrak{v}^M)$ on the $\mathcal{O}$-module $M$ (in the sense of Definition \ref{def:edv}) and any positive integer $d \in \mathbb{Z}_{\geq d_{\{\mathfrak{v}\}}}$ (for $d_{\{\mathfrak{v}\}}$ see Notation \ref{not:d_D}) the compositions $\mathfrak{v}_{\leq d}:\mathcal{O} \xrightarrow{\mathfrak{v}}\overline{\mathbb{N}} \xrightarrow{\min_d} \mathbb{N}_d$ and $\mathfrak{v}_{\leq d}^M:M \xrightarrow{\mathfrak{v}^M}\overline{\mathbb{Z}} \xrightarrow{\min_0} \mathbb{Z}_{\leq 0}$ give an extended quasi-valuation $(\mathfrak{v}_{\leq d}, \mathfrak{v}_{\leq d}^M)$ on $M$, where $\min_d:\overline{\mathbb{N}} \rightarrow \mathbb{N}_d$ and $\min_0:\overline{\mathbb{Z}} \rightarrow \mathbb{Z}_{\geq 0}$ are the natural forgetting  maps: $\min_d(n)=\min\{d, n\}$ (with $d=0$ in the second case).

     Moreover, for any $n \in \mathbb{N}_d:\ \mathcal{F}_{\mathfrak{v}_{\leq d}}(n)=\mathcal{F}_{\mathfrak{v}}(n) \triangleleft \mathcal{O}$ and $\mathcal{F}^M_{\mathfrak{v}_{\leq d}}(-n)=\mathcal{F}^M_{\mathfrak{v}}(-n) \leq M$ (cf. (\ref{eq:fvMl})).
\end{example}

\begin{example}\label{ex:quasi-redM}
Given an extended discrete valuation $(\mathfrak{v}, \mathfrak{v}^M)$ on the $\mathcal{O}$-module $M$ and a submodule $N \leq M$ with $\mathfrak{v}^M(N) \geq 0$ (i.e., for all $m \in N$ we have $\mathfrak{v}^M(m) \geq 0)$ we can consider the pair $(\mathfrak{q}, \mathfrak{q}^M)$ consisting of the restricted quasi-valuation $$\mathfrak{q}=\mathfrak{v}\big\vert_{\mathcal{O}/{\rm Ann}(M/N)}:\mathcal{O}/{\rm Ann}(M/N) \rightarrow \mathbb{N}_{d_{\{\mathfrak{v}\}}}$$ (for $d_{\{\mathfrak{v}\}}$ see Notation \ref{not:d_D} and use that $\mathcal{F}_\mathfrak{v}(d_{\{\mathfrak{v}\}})\supset {\rm Ann}(M/N)$ by Proposition \ref{prop:pullbackring} \textit{(2)}) and the association 
\begin{equation*}
         \mathfrak{q}^{M/N}: M/N \rightarrow \mathbb{M}_d, \ m+ N\mapsto \begin{cases}
            \mathfrak{v}^M(m) & \text{ if $\mathfrak{v}^M(m) < 0$}\\
            0  & \text{ if $\mathfrak{v}^M(m) \geq 0$},
        \end{cases}
    \end{equation*}
    which is once again well-defined (by property (b) of Definition \ref{def:edv}). Then $(\mathfrak{q}, \mathfrak{q}^{M/N})$ is an extended discrete quasi-valuation. We sometimes denote it by $\mathfrak{v}\big\vert_{M/N}$ and call it the \emph{`quasi-reduction'} of $\mathfrak{v}$.

    Moreover, if we denote by $q:\mathcal{O} \rightarrow \mathcal{O}/{\rm Ann}(M/N)$ and $q^M:M\rightarrow M/N$ the standard quotient maps, then for any $n \in \mathbb{N}_d$: $q^{-1} \big(\mathcal{F}_{\mathfrak{v}\big|_{\mathcal{O}/{\rm Ann}(M/N)}}(n)\big)=\mathcal{F}_{\mathfrak{v}}(n)\triangleleft \mathcal{O}$ and $(q^M)^{-1} \big(\mathcal{F}^M_{\mathfrak{v}\big|_{M/N}}(-n)\big)=\mathcal{F}^M_{\mathfrak{v}}(-n)\leq M$. 
\end{example}

Similarly to Question \ref{q:whenextendquasi} we can ask the following:

\begin{question}\label{q:extquasi}
    Let $\mathcal{O}$ be a (Noetherian) $k$-algebra, $M$ a finitely generated $\mathcal{O}$-module and $N \leq M$ a finite $k$-codimensional submodule. Suppose that we have a discrete valuation $\mathfrak{v}:\mathcal{O} \rightarrow \overline{\mathbb{N}}$ with its quasi-reduction $\mathfrak{q}=\mathfrak{v}\big\vert_{\mathcal{O}/{\rm Ann}(M/N)}:\mathcal{O}/{\rm Ann}(M/N) \rightarrow \mathbb{N}_d$ and also a quasi-extension $\mathfrak{q}^{M/N}: M/N \rightarrow \mathbb{M}_d$ of $\mathfrak{q}$. Under what conditions is there an extension $\mathfrak{v}^M:M \rightarrow \overline{\mathbb{Z}}$ of the discrete valuation $\mathfrak{v}$ in the sense of Definition \ref{def:edv}, such that $(\mathfrak{v}, \mathfrak{v}^M)\big\vert_{M/N}=(\mathfrak{q}, \mathfrak{q}^{M/N})$?
\end{question}

Analogously to the case of quasi-valuations (cf. Remark \ref{rem:constructionofq-val}), over the finite $k$-dimensional module $M/N$ and finite $k$-dimensional algebra $\mathcal{O}/{\rm Ann}(M/N)$ all extended quasi-valuations can be `finitely' described:

\begin{remark}\label{rem:constructionofextq-val}
    Consider a finite $k$-dimensional module $M$ over the finite $k$-dimensional algebra $\mathcal{O}$. By  Remark \ref{rem:constructionofq-val}, for any quasi-valuation $\mathfrak{q}:\mathcal{O} \rightarrow \mathbb{N}_d$ there exists a $k$-basis $\{f_1, \ldots, f_o\}$ of $\mathcal{O}$ (with $\dim_k(\mathcal{O})=o$) and a multiset of assigned values $\{q_1, \ldots, q_o\}$ ($q_i \in \mathbb{N}_d$ for all $1 \leq i \leq o$), such that
    \begin{equation*}
        \text{for any } f = \sum_{i=1}^{o}a_if_i \in \mathcal{O}  \text{ we have } \mathfrak{q}\big(\sum_{i=1}^{o}a_if_i\big)=\min\{q_i\,:\, a_i \in k\setminus\{0\}\}.
    \end{equation*}
Similarly,  for every extension $\mathfrak{q}^M$ of the quasi-valuation $\mathfrak{q}:\mathcal{O} \rightarrow \mathbb{N}_{d}$ we can assign a $k$-vector space basis $\{m_1, m_2, \ldots, m_u\}$ of $M$ (with $\dim_k M=u$) and values $\mathfrak{q}^M(m_j) = s_j\in\mathbb{M}_d$ for all $1 \leq j \leq u$, such that 
for every $m=\sum_{j=1}^u b_jm_j$ we have $\mathfrak{q}^M(m) = \min\{s_j \,:\, b_j \neq 0\}$.
\end{remark}
We also have the following analogue of Proposition \ref{prop:quasivalcond}:
\begin{prop}\label{prop:extquasivalcond}
    Let $\mathcal{O}$ be a finite dimensional $k$-algebra and $M$ a finite $k$-dimensional module over it. Suppose that $\mathfrak{q}:\mathcal{O} \rightarrow \mathbb{N}_d$ is a quasi-valuation given by formula (\ref{eq:quasiwithbasis}) using the $k$-basis $\{f_1, \ldots, f_o\}$ of $\mathcal{O}$ and the multiset $\{q_1, \ldots, q_o\}$ of assigned values from $\mathbb{N}_d$. Consider now a $k$-basis $\{m_1, \ldots, m_u\}$ of $M$ (with $\dim_k(M)=u$) and a multiset $\{s_1, \ldots, s_u\}$ of assigned values ($s_j\in \mathbb{M}_d$ for all $1 \leq j \leq u$). The mapping 
    \begin{equation}
        \mathfrak{q}^M:M \rightarrow \mathbb{M}_d, \ m=\sum_{j=1}^ub_jm_j \mapsto \min\{ s_j \,:\, b_j \neq 0\}
    \end{equation}
    gives an extension of $\mathfrak{q}$ if and only if it satisfies the following compatibility condition:
            \begin{equation*}
            \begin{split}
\big({\rm Span}_k(\{ f_i\, :\, q_i =a \})&\,\setminus\{0\}\big)\cdot \big({\rm Span}_k(\{ m_i\, :\, s_i =b \}) \setminus \{0\}\big) \subset \\ &\,\subset  {\rm Span}_k(\{ m_i\, :\, s_i \geq a\ast b \})\setminus {\rm Span}_k(\{ m_i\, :\, s_i \geq a\ast b +1\})
\end{split}
        \end{equation*}
        for any $a \in \mathbb{N}_d$ and $b \in \mathbb{M}_d$, i.e., for any $0 \neq f \in {\rm Span}_k(\{ f_i\, :\, q_i =a \})$ and $0 \neq m \in {\rm Span}_k(\{ m_i\, :\, m_i =b \})$, the product $fm \in  {\rm Span}_k(\{ m_i\, :\, s_i \geq a\ast b \})\setminus {\rm Span}_k(\{ m_i\, :\, s_i \geq a\ast b +1\})$.
\end{prop}

We remark that these conditions are easier to check when the partition sets $\{\{f_i\,:\,q_i=a\}\}_{a \in \mathbb{N}_d}$ and  $\{\{m_i\,:\,s_i=b\}\}_{b \in \{-d, \ldots, 0\}}$ each have smaller cardinality. We will see such a concrete construction in Example \ref{ex:curvewith1qv}.

\subsection{Quasi-realizability and lattice homology}\,

Remarkably, all the results of section \ref{s:4} remain true in the setting of (extended) quasi-valuations. For the sake of completeness we present them in this new language here.

\begin{define}
For any finite (multi)set $\mathcal{D}=\{\mathfrak{q}_1, \ldots, \mathfrak{q}_r\}$
of extended discrete quasi-valuations (with value monoids $\{\mathbb{N}_{d_1}, \mathbb{N}_{d_2}, \ldots, \mathbb{N}_{d_r}\}$) and tuple $\ell=(\ell_1,\ldots, \ell_r)\in \mathbb{Z}^r$, with $0 \leq \ell_v \leq d_v$ for all $v \in \{1;\ldots, r\}$, we define
\begin{equation}\label{eq:fdMl-q}
    \mathcal{F}_{\mathcal{D}}(\ell)= \bigcap_{v =1}^{r}\mathcal{F}_{\mathfrak{q}_v}(\ell_v) \triangleleft \mathcal{O} \ \ \text{ and }\ \ \mathcal{F}^M_{\mathcal{D}}(-\ell) = \bigcap_{v =1}^{r}\mathcal{F}^M_{\mathfrak{q}_v}(-\ell_v) \leq M.
\end{equation}
\end{define}
\begin{ass}\label{ass:fincodim-q}
    We also make the assumption (compare with Assumption \ref{ass:fincodim}) that the ideals $\mathcal{F}_{\mathcal{D}}(\ell)$ (respectively the submodules $\mathcal{F}_{\mathcal{D}}^M(\ell)$) have finite $k$-codimension in $\cO$ (respectively in $M$). A fi\-nite (multi)set of extended discrete quasi-valuations satisfying this finite codimensionality condition will be called  a \textit{`finite collection of extended (discrete) quasi-valuations'}.
\end{ass}

\begin{define}
    For any finite collection $\mathcal{D}=\{\mathfrak{q}_1, \ldots, \mathfrak{q}_r\}$
of extended discrete quasi-valuations (with value monoids $\{\mathbb{N}_{d_1}, \mathbb{N}_{d_2}, \ldots, \mathbb{N}_{d_r}\}$) we denote the tuple $(d_1, \ldots, d_r) \in (\mathbb{Z}_{>0})^r$ by $d_{\mathcal{D}}$ (compare with Notation \ref{not:d_D} --- these two, in fact, agree for an extended valuation and its quasi-reduction, see Example \ref{ex:quasi-redM}) and we introduce 
the following height functions on the rectangle $R(0, d_{\mathcal{D}})$:
\begin{equation}\label{eq:qhandhM}
    \mathfrak{h}_{\mathcal{D}}, \mathfrak{h}^\circ_{\mathcal{D}}: \mathbb{Z}^r\cap R(0, d_{\mathcal{D}}) \rightarrow \mathbb{Z}, \hspace{5mm} \mathfrak{h}_{\mathcal{D}}(\ell) = \dim_k\mathcal{O}/\mathcal{F}_{\mathcal{D}}(\ell) \hspace{5mm} \text{and} \hspace{5mm} \mathfrak{h}^\circ_{\mathcal{D}}(\ell) = \dim_kM/\mathcal{F}_{\mathcal{D}}^M(-\ell).
\end{equation}
Since $\mathcal{F}_\mathcal{D}(0)=\mathcal{O}$, we have $\mathfrak{h}_{\mathcal{D}}(0)=0$ and  $\mathfrak{h}_{\mathcal{D}}$ increasing with $\mathfrak{h}^M_{\mathcal{D}}$ decreasing. Therefore, we can set the weight function as $w_{\mathcal{D}, 0}(\ell)=\frh_{\mathcal{D}}(\ell)+\frh^\circ_{\mathcal{D}} (\ell)-\frh^\circ_{\mathcal{D}} (0)$ for every $\ell \in \mathbb{Z}^r \cap R(0, d_{\mathcal{D}})$ (compare with (\ref{eq:w_0def})) and extend it to higher dimensional cubes by formula (\ref{eq:9weight}). Thus, we obtain the lattice homology module $\mathbb{H}_*(R(0, d_{\mathcal{D}}), w_{\mathcal{D}})$.
\end{define}

The point is that the Independence Theorem \ref{th:IndepMod} and all the results of sections \ref{s:4}--\ref{s:ideals} are also valid in this quasi-valuated setup, too, since during the proofs we make every computation in the finite dimensional quotients $M/N$ and $\mathcal{O}/{\rm Ann}(M/N)$, where $N$ is the submodule $N=\mathcal{F}_\mathcal{D}^M(0)$ (and, by part \textit{(a)} of Theorem \ref{th:properties}, in the finite rectangle $R(0, d_{\mathcal{D}})$). We list the results here, the detailed proofs are left to the reader.

\begin{theorem}[Independence Theorem]\label{th:IndepMod-q} Let $k$ be any field, $\mathcal{O}$ a (Noetherian) $k$-algebra and $M$ a finitely generated module over it. Let $\mathcal{D}=\{\mathfrak{q}_1, \ldots, \mathfrak{q}_r\}$ and $\mathcal{D}'=\{\mathfrak{q}'_1, \ldots, \mathfrak{q}'_{r'}\}$ be two collections of extended discrete quasi-valuations on $M$.  
        Suppose that $\mathcal{F}_{\mathcal{D}}^M(0) = \mathcal{F}_{\mathcal{D}'}^M(0)\leq M$ and
both pairs $(\mathfrak{h}_{\mathcal{D}}, \mathfrak{h}_{\mathcal{D}}^\circ)$ and $(\mathfrak{h}_{\mathcal{D}'}, \mathfrak{h}_{\mathcal{D}'}^{\circ})$ satisfy the Combinatorial Duality Property.
        Then the cubical complexes $S_{n, \mathcal{D}}\subset R(0, d_{\mathcal{D}})$ and $S_{n, \mathcal{D}'}\subset R(0, d_{\mathcal{D}'})$ 
        associated with the corresponding lattices and weight functions 
        are homotopy equivalent for every $n \in \mathbb{Z}$.
        Even more, the homotopy equivalences respect the inclusions
        $S_{n, \mathcal{D}}\hookrightarrow S_{n+1, \mathcal{D}}$ and $S_{n, \mathcal{D}'}\hookrightarrow S_{n+1, \mathcal{D}'}$ (in the precise sense of Theorem \ref{th:IndepModTop}), hence
        \begin{equation*}
        \mathbb{H}_*(R(0, d_{\mathcal{D}}), w_{\mathcal{D}}) \cong  \mathbb{H}_*(R(0, d_{\mathcal{D}'}, w_{\mathcal{D}'}) \text{ as bigraded } \mathbb{Z}[U]\text{-modules.}
        \end{equation*}
\end{theorem}

Therefore, if we set some field $k$, a Noetherian $k$-algebra $\mathcal{O}$ and a finitely generated $\mathcal{O}$-module $M$, we get well-defined lattice homology $\mathbb{Z}[U]$-modules for \textit{`quasi-realizable'} submodules of $M$:

\begin{define}\label{def:REAL-q}
    Let $N$ be a finite codimensional submodule of $M$. $N$ is called \textit{`quasi-realizable'} if some finite collection $\mathcal{D}$ of 
    extended discrete quasi-valuations 
    satisfies 
    $N=\mathcal{F}^M_{\mathcal{D}}(0)$. In this case $\mathcal{D}$ is called a \textit{`quasi-realization of $N$'}. 
      We say that the realization is a `{\it CDP quasi-realization}'
if the associated pair of functions $(\mathfrak{h}_{\mathcal{D}}, \mathfrak{h}^\circ_{\mathcal{D}})$
satisfies the Combinatorial Duality Property. Similarly to the previous cases, from any quasi-realization $\mathcal{D}$ of $N$ we can construct a CDP quasi-realization $\mathcal{D}^{\natural}$ just by postcomposing every extended discrete quasi-valuation with the `doubling' map $\overline{\mathbb{Z}} \xrightarrow{\cdot2}\overline{\mathbb{Z}}$ (the new quasi-valuations are denoted simply by $\{2\mathfrak{q}_v\}_v$).  
\end{define}

Thus, using the Independence Theorem \ref{th:IndepMod-q}
for any  {\it quasi-realizable} submodule $N\leq M$  we can define its lattice homology $\mathbb{H}_*(N \hookrightarrow_{\mathcal{O}} M)$ analogously to Definition \ref{def:LCofMOD}:

\begin{define}
    Let $k$ be any field, $\mathcal{O}$ a (Noetherian) $k$-algebra and $M$ a finitely generated module over it. Let $N \leq M$ be a quasi-realizable submodule and $\mathcal{D}$ a CDP quasi-realization of it. Then consider the multifiltrations $\mathcal{F}_\mathcal{D}$ and $\mathcal{F}_{\mathcal{D}}^M$ of (\ref{eq:fdMl-q}) and the height functions $\mathfrak{h}_\mathcal{D}$ and $\mathfrak{h}_\mathcal{D}^\circ$ on the rectangle $R(0, d_\mathcal{D})$ of (\ref{eq:qhandhM}) associated with this collection. The lattice homology of $N \leq M$ is defined as $\mathbb{H}_*(N \hookrightarrow_{\mathcal{O}} M):=\mathbb{H}_*(R(0, d_{\mathcal{D}}), w_{\mathcal{D}})$, where $w_{\mathcal{D}, 0}(\ell)=\frh_{\mathcal{D}}(\ell)+\frh^\circ_{\mathcal{D}} (\ell)-\frh^\circ_{\mathcal{D}} (0)$ for every $\ell \in \mathbb{Z}^r \cap R(0, d_{\mathcal{D}})$.
\end{define}

\begin{remark} We use the same notation as before, since if $N\leq M$ is realizable in the sense of Definition \ref{def:REAL}, then for any realization $\mathcal{D}=\{ \mathfrak{v}_1, \ldots, \mathfrak{v}_r\}$ of it the (multi)set of the quasi-valuations $\mathcal{D}_{\leq d_{\mathcal{D}}}=\{(\mathfrak{v}_1)_{\leq d_{\{\mathfrak{v}_1\}}}, \ldots,(\mathfrak{v}_1)_{\leq d_{\{\mathfrak{v}_r\}}}\}$ from Example \ref{ex:levagottv^M} will give a quasi-realization of $N$. Moreover, the corresponding height functions, and hence the weight functions agree in the rectangle $R(0, d_{\mathcal{D}})$, which contains all the lattice homological information (see Theorem \ref{th:properties} \textit{(a)}). So, in fact, for realizable a submodule $N \leq M$ both definitions give the same $\mathbb{H}_*(N \hookrightarrow_{\mathcal{O}} M)$. 

Even more, the restricted quasi-valuations $\mathcal{D}\big\vert_{M/N}=\{\mathfrak{v}_1\big\vert_{M/N}, \ldots,\mathfrak{v}_r\big\vert_{M/N}\}$ from Example \ref{ex:quasi-redM} give a quasi-realization of $N/N$ in $M/N$, with the very same height and weight functions on the very same rectangle. Thus $\mathbb{H}_*(N \hookrightarrow_{\mathcal{O}} M)\cong \mathbb{H}_*(N/N \hookrightarrow_{\mathcal{O}/{\rm Ann}_{\mathcal{O}}(M/N)} M/N)$ as well.
\end{remark}

On the other hand, notice that, even though $N$ is realizable in $M$, $N/N$ is not (necessarily) so in $M/N$ (see in particular the example in Remark \ref{rem:imageofintclosed}). This means, that especially in torsion modules over nilpotent algebras, the class of quasi-realizable submodules might be strictly larger than that of realizable submodules. It is not clear, however, whether this phenomena can occur in the free and symmetric  case as well, i.e.,  when $M=\mathcal{O}=\mathcal{O}_m=k[x_1, \ldots, x_m]$ (or $k[x_1, \ldots, x_m]_{(x_1, \ldots, x_m)}, k[[x_1, \ldots , x_m]]$, etc., depending on the context --- for the complex analytic case see Corollary \ref{cor:intcl=real=qreal}). 

\begin{question}\label{q:qreal=real}
    Are there quasi-realizable ideals in $\mathcal{O}_m$ which are not integrally closed? (Here we also consider those ideals, which have infinite $k$-codimension in $\mathcal{O}_m$ but can be obtained as $\mathcal{F}_{\mathcal{D}}(d_{\mathcal{D}})$ for some finite set $\mathcal{D}$ of extended quasi-valuations --- not necessarily respecting Assumption \ref{ass:fincodim-q}.) More generally, under what conditions do quasi-realizable and integrally closed ideals agree in a Noetherian $k$-algebra?
\end{question}

This question naturally holds significance in the study of integrally reduced algebras, i.e., quotients of $\mathcal{O}_m$ by integrally closed ideals (cf. section \ref{ss:ARTIN}). Indeed, an affirmative answer to Question \ref{q:qreal=real} would imply the following characterization:
\begin{align}
    \text{a Noetherian $k$-algebra $A$  is integrally reduced}& \Leftrightarrow \nonumber \\
     \Leftrightarrow \text{ the ideal $0 \triangleleft A$ is quasi-realizable in} & \text{ the wider sense of Question \ref{q:qreal=real} }  \Leftrightarrow \nonumber\\
     \Leftrightarrow\text{ if }a \in A \text{ and }\mathfrak{q}(a)=0 &\text{ for all quasi-valuations }\mathfrak{q}:A \rightarrow\mathbb{N}_d, \text{ then }a=0 .\label{eq:intred}
\end{align}
(Once again, compare with the complex analytic case in Proposition \ref{prop:VAL} and Corollary \ref{cor:intcl=real=qreal}). 

This would also imply the following intrinsic characterization of the integral reduction of a Noetherian $k$-algebra $A$:
\begin{equation*}
        \widehat{A}=A/\{a \in A\,:\, \text{ for all } \mathfrak{q}: A \rightarrow \mathbb{N}_d \text{ we have }\mathfrak{q}(a)=0\}.
    \end{equation*}
    (Compare with the complex analytic case in Proposition \ref{prop:analintred}).)

\bekezdes Regarding the lattice homology of quasi-realizable submodules we also have the following properties inherited from sections \ref{s:4} and \ref{s:ideals}: 

\begin{theorem}\label{th:properties-quasi}
    Let $\mathcal{O}$ be a Noetherian $k$-algebra, $N \leq M$ finitely generated modules with $N$ quasi-realizable and $\dim_k(M/N) < \infty$. \\
\noindent (a) The Euler characteristic is well-defined and is $eu(\mathbb{H}_*(N \hookrightarrow M))=\dim_k (M/N)$, i.e., $\mathbb{H}_*(N \hookrightarrow M)$ categorifies the codimension ${\rm codim}_k(N \hookrightarrow M)$.  \\
\noindent (b) The $\mathbb{Z}[U]$-module (and the homotopy type of the $S_n$ spaces) only depends on the quotient module $M/N$ over the algebra $\mathcal{O}/{\rm Ann}_{\mathcal{O}}(M/N)$. \\
\noindent (c) ${\rm homdim}(N \hookrightarrow M)= \begin{cases}
    \max \{ q: \mathbb{H}_{q}(N \hookrightarrow M) \neq 0\} &  \text{ if }N\neq M;\\
    -1 &  \text{ if }N= M.
\end{cases}$\\

\noindent \hspace{33mm} $\leq \min \{ |\mathcal{D}|\,:\, \mathcal{D} \text{ quasi-realization of } N\} -1$.\\
\noindent (d) 
    If $k$ is an algebraically closed field, $(\mathcal{O}, \mathfrak{m})$ a local Noetherian $k$-algebra, $M$ a finitely generated $\mathcal{O}$-algebra and $N \leq M$ a finite codimensional quasi-realizable submodule, then for every $n > 0$ the $S_n$ space is contractible, hence the weight-grading of the \emph{reduced} lattice homology  is nonnegative.\\
    \noindent (e) If  $M=\mathcal{O}$ and $N=\mathcal{I} \triangleleft \mathcal{O}$ a quasi-realizable ideal, then every weight function $w_{\mathcal{D}}$ correspoding to a quasi-realization $\mathcal{D}$ (in the sense of Definition \ref{def:REAL-q}) is symmetric with respect to $d_{\mathcal{D}}$, such that the $S_n$ spaces (up to homotopy) and the lattice homology $\mathbb{H}_*(\mathcal{I} \hookrightarrow \mathcal{O})$ inherits a $\mathbb{Z}_2$-symmetry. We will denote the latter with $\mathbb{SH}_*(\mathcal{I} \triangleleft \mathcal{O})$.
\end{theorem}

\begin{remark}\label{rem:robustness-q}
    One advantage of the quasi-valuative viewpoint is that quasi-realizability is more robust with respect to module operations: if $N \leq M$ is quasi-realizable, then so is
    \begin{itemize}
        \item $\phi^{-1}(N)$ for any $\mathcal{O}$-module homomorphism $\phi: M' \rightarrow M$, since the precomposition of a quasi-valuation with a homomorphism is still a quasi-valuation;
        \item $\psi(N)$ for any quotient map $\psi: M \rightarrow M/N'$ with $N' \leq N$, since restrictions of quasi-valuations with these parameters are quasi-valuations.
    \end{itemize}
    So, in this setting $N \leq M$ is quasi-realizable if and only if so is $N/N$ in $M/N$ (which is not true for the integral closedness and realizability properties, see, e.g., the example in Remark \ref{rem:imageofintclosed}) and $\mathbb{H}_*(N \hookrightarrow M) \cong \mathbb{H}_*(N/N \hookrightarrow M/N)$.

    Moreover, these quasi-valuations offer a more flexible setting in which to investigate the homological dimension Conjecture \ref{conj:elso}. Indeed, this is due to the fact, that by Propositions \ref{prop:quasivalcond} and \ref{prop:extquasivalcond} all quasi-valuations on finite $k$-dimensional algebras and modules can be specified by `finite data'. We give here a concrete example, providing a single quasi-valuation for the quasi-realization of the submodule $\Omega^1_{\overline{C}}\leq \omega^R_C$ assigned to the space curve singularity of Example \ref{ex:345534} with vanishing $\mathbb{H}_{an, \geq 1}$.
\end{remark}

 \begin{example}\label{ex:curvewith1qv}
     In the non-symmetric (non-Gorenstein) case of Example \ref{ex:345534} we can use a single extended discrete \textit{quasi}-valuation to obtain a quasi-realization of $\Omega^1_{\overline{C}}$. In fact, using Propositions \ref{prop:quasivalcond} and \ref{prop:extquasivalcond}, it is enough to specify $\mathbb{C}$-bases of $\mathcal{O}/\mathcal{C}$ and $\omega^R_C/\Omega^1_{\overline{C}}$, assign to them values and check whether these satisfy the compatibility conditions described there. In our case we make the following choices:
    \begin{align*}
    \mathcal{O}/\mathcal{C} \cong \mathcal{O}_{\mathbb{C}^3, 0}/\mathfrak{m}_{\mathbb{C}^3, 0}^2 = \mathbb{C}\big\langle f_1, f_2, f_3, &\, f_4\big\rangle= \mathbb{C}\big\langle [1], [x], [y], [z]\big\rangle \text{ with }\\
    \mathfrak{q}([1])=0, \hspace{10mm} \mathfrak{q}([x])=8, \hspace{10mm} &\, \mathfrak{q}([y])=7, \hspace{10mm} \mathfrak{q}([z])=9, \\
    \text{with } \mathfrak{q}\left(\sum_{i=1}^4 a_i f_i\right)= \min\{ \mathfrak{q}(f_i)\,:\,&\, a_i \neq 0\} \text{ and } d_{\{\mathfrak{q}\}}=11,
    \end{align*}
    whereas $\omega^R_C/\Omega^1_{\overline{C}}=\mathbb{C}\langle m_1, \ldots, m_8\rangle=$
    \begin{sizeddisplay}{\small}
    \begin{align*}
    \mathbb{C}\Bigg\langle\left[\dfrac{1}{t}dt-\dfrac{1}{s}ds\right], \left[\dfrac{1}{t^2}dt\right], \left[\dfrac{1}{t^2}dt-\dfrac{1}{s^2}ds\right], \left[\dfrac{1}{t^3}dt\right], \left[\dfrac{1}{s^3}ds\right], &\, \left[\dfrac{1}{t^4}dt-\dfrac{1}{s^6}ds\right], \left[\dfrac{1}{t^5}dt-\dfrac{1}{s^4}ds\right],
    \left[\dfrac{1}{t^6}dt-\dfrac{1}{s^5}ds\right]\Bigg\rangle \\
    \mathfrak{q}^M\left( \left[\dfrac{1}{t}dt-\dfrac{1}{s}ds\right]\right)=-2, \hspace{5mm} \mathfrak{q}^M\left( \left[\dfrac{1}{t^2}dt\right]\right)&\,=-1, \hspace{5mm} \mathfrak{q}^M\left( \left[\dfrac{1}{t^2}dt-\dfrac{1}{s^2}ds\right]\right)=-4, \\
    \mathfrak{q}^M\left( \left[\dfrac{1}{t^3}dt\right]\right)=-3, \hspace{5mm} \mathfrak{q}^M\left( \left[\dfrac{1}{s^3}ds\right]\right)=&\,-3, \hspace{5mm} \mathfrak{q}^M\left( \left[\dfrac{1}{t^4}dt-\dfrac{1}{s^6}ds\right]\right)=-10, \\
    \mathfrak{q}^M\left( \left[\dfrac{1}{t^5}dt-\dfrac{1}{s^4}ds\right]\right)=-9, \hspace{5mm} &\, \mathfrak{q}^M\left( \left[\dfrac{1}{t^6}dt-\dfrac{1}{s^5}ds\right]\right)=-11, \\
    \text{with } \mathfrak{q}^M\left(\sum_{j=1}^8 b_j m_j\right)=\min\{ &\,  \mathfrak{q}(m_j)\,:\, b_j \neq 0\} \text{ and } d_{\{\mathfrak{q}\}}=11.
    \end{align*}
    \end{sizeddisplay}
    One can check by direct computations that the compatibility conditions of Propositions \ref{prop:quasivalcond} and \ref{prop:extquasivalcond} are satisfied in this case, hence these assignments indeed give an extended quasi-valuation on $\omega_C^R/\Omega^1_{\overline{C}}$. The doubled height functions and the weight function are the following:
    \vspace{2mm}
    \begin{center}
    \resizebox{15.5cm}{!}{
        $\begin{array}{c|ccccccccccccccccccccccc}
        \ell & 0 & 1 & 2 & 3 & 4 & 5 & 6 & 7 & 8 & 9 & 10 & 11 & 12 & 13 & 14 & 15 & 16 & 17 & 18 & 19 & 20 & 21 & 22 \\
        \hline
        \mathfrak{h}_{\{2\mathfrak{q}\}}(\ell)        & 0 & 1 & 1 & 1 & 1 & 1 & 1 & 1 & 1 & 1 & 1 & 1 & 1 & 1 & 1 & 2 & 2 & 3 & 3 & 4 & 4 & 4 & 4 \\
        \mathfrak{h}_{\{2\mathfrak{q}\}}^\circ(\ell) & 8 & 8 & 7 & 7 & 6 & 6 & 4 & 4 & 3 & 3 & 3 & 3 & 3 & 3 & 3 & 3 & 3 & 3 & 2 & 2 & 1 & 1 & 0 \\
        
        \hline w_{\{2\mathfrak{q}\}}(\ell) & 0 & 1 & 0 & 0 & -1 & -1 & -3 & -3 & -4 & -4 & -4 & -4 & -4 & -4 & -4 & -3 & -3 & -2 & -3 & -2& -3 & -3 & -4  
        \end{array}$}
    \end{center}
    \vspace{2mm}
    Clearly, the lattice homology agrees with the one computed in Example \ref{ex:345534}. However, it is rather unclear whether this quasi-valuation can be extended to a discrete valuation on the whole of $\mathcal{O}$ and $\omega^R_C$, or how else can we find a realization of $\Omega^1_{\overline{C}}$ with a single extended discrete valuation.
 \end{example}

\begin{remark}\label{rem:MOSTGEN}
    One could also reformulate with this quasi-valuative language the generalizations \ref{th:Indep-MGF} and \ref{th:Indep-MGUF} of the Independence Theorem, provided one is satisfied with working only on a finite rectangle. In this case the analogues of the valuative functions of Remark \ref{rem:valfunctions} corresponding to compatible filtrations should be mappings $\mathfrak{q}^{M_1}:M_1\rightarrow \mathbb{N}_d, \ \mathfrak{q}^{M_{-1}}:M_{-1} \rightarrow \mathbb{M}_d$ and $\mathfrak{q}^{M_0}:M_0\rightarrow \mathbb{M}_d$ satisfying the analogues of properties (a)--(d)  in loc. cit. for some $d \in \mathbb{Z}_{>0}$. Even more, using these conventions, one can prove that the resulting lattice homology modules ($\mathbb{UH}_*(M_1, M_{-1}, N_0 \leq M_0)$ or $\mathbb{H}_*(M_1, N_{-1} \leq M_{-1}, N_0 \leq M_0)$ computed in a finite rectangle) only depend on the isomorphism type of the $\bullet$ multiplication under the following weaker (than Definition \ref{def:isotype}) isomorphism relation:
    given two tuples $(M_1, M_{-1}, N_0\leq M_0)$ and $(M_1', M'_{-1}, N'_0\leq M'_0)$ of $k$-vector spaces, the $k$-bilinear `multiplication' maps
    \begin{align*}
        \bullet: \dfrac{M_1}{(N_0\,:\, M_{-1})} \times \dfrac{M_{-1}}{(N_{0}:M_1)} \rightarrow \dfrac{M_\bullet}{N_0\cap M_{\bullet}} \hspace{5mm} \text{and} \hspace{5mm } \bullet': \dfrac{M'_1}{(N'_0\,:\, M'_{-1})} \times \dfrac{M'_{-1}}{(N'_{0}:M'_1)} \rightarrow \dfrac{M'_{\bullet'}}{N'_0\cap M'_{\bullet'}}
    \end{align*}
    are considered isomorphic if there exist $k$-linear isomorphisms 
    \begin{equation*}
        L_1: \dfrac{M_1}{(N_0\,:\, M_{-1})} \xrightarrow{\cong} \dfrac{M'_1}{(N'_0\,:\, M'_{-1})}, \hspace{5mm} L_{-1}:\dfrac{M_{-1}}{(N_{0}:M_1)} \xrightarrow{\cong} \dfrac{M'_{-1}}{(N'_{0}:M'_1)}, \hspace{5mm} L_0:  \dfrac{M_{\bullet}}{N_0\cap M_{\bullet}} \xrightarrow{\cong} \dfrac{M'_{\bullet'}}{N'_0\cap M'_{\bullet'}}
    \end{equation*}
  satisfying that 
 for all $m_1 \in \frac{M_1}{(N_0\,:\, M_{-1})}$ and $m_{-1} \in\frac{M_{-1}}{(N_{0}:M_1)}$ we have $L_0(m_1 \bullet m_{-1})=L_1(m_1)\bullet'L_{-1}(m_{-1})$.
\end{remark}

\subsection{Analogies with arcs and jets}\label{ss:arcs}\,

Given a field $k$ and a (Noetherian) $k$-algebra $\mathcal{O}$, \textit{(formal) arc}s of ${\rm Spec}\,\mathcal{O}$ are defined as morphisms ${\rm Spec}\, K[[t]] \rightarrow {\rm Spec}\,\mathcal{O}$, where $K\supset k$ is a field extension (for more see, e.g., \cite{IshiiArc} --- compare also with paragraph \ref{bek:4.10.2}). On the ring level, these correspond to maps $\alpha: \mathcal{O} \rightarrow K[[t]]$, which, post-composed with the $ord_{t}: K[[t]] \rightarrow \mathbb{N}\cup \{ \infty \}$ order function, give rise to valuations $\mathfrak{v}_{\alpha}=ord_t \circ \alpha$. In fact, the valuations corresponding to the normalization of components of reduced curve singularities (cf. (\ref{eq:normalizationvaluation})) are also of this form with $k=K=\mathbb{C}$.

\begin{remark}\label{rem:arcvaluations}
    Clearly, for an arc $\alpha:\mathcal{O} \rightarrow k[[t]]$, the corresponding Hilbert function $\mathfrak{h}_{\mathfrak{v}_{\alpha}}$, defined as in (\ref{eq:qhandhM}), satisfies $\mathfrak{h}_{\mathfrak{v}_{\alpha}}(\ell+1)-\mathfrak{h}_{\mathfrak{v}_{\alpha}}(\ell) \in \{0, 1\}$ for any $\ell \in \mathbb{Z}$. Indeed, for any two functions $f, f' \in \mathcal{F}_{\mathfrak{v}_{\alpha}}(\ell) \setminus \mathcal{F}_{\mathfrak{v}_{\alpha}}(\ell+1)$, i.e., satisfying 
    \begin{center}
        $\alpha(f)=at^\ell+$ \textit{higher order terms}, $\alpha(f')=a't^\ell+$ \textit{higher order terms} for some $a, a' \in k \setminus \{0\}$,
    \end{center} 
    their linear combination $a'f-af' \in \mathcal{F}_{\mathfrak{v}_{\alpha}}(\ell+1)$. The reverse implication is an interesting question.
\end{remark}

\begin{question}\label{q:hilbertjump01}
    Let  $\mathfrak{v}: \mathcal{O} \rightarrow \mathbb{N} \cup \{\infty\}$ be a discrete valuation with its corresponding Hilbert function satisfying $\mathfrak{h}_{\mathfrak{v}}(\ell+1)-\mathfrak{h}_{\mathfrak{v}}(\ell) \in \{0, 1\}$ for any $\ell \in \mathbb{Z}$. Under what conditions is there an arc $\alpha_{\mathfrak{v}}:\mathcal{O} \rightarrow k[[t]]$ satisfying $\mathfrak{v}=\mathfrak{v}_{\alpha_{\mathfrak{v}}}$?
\end{question}

This question is also related to the theory of maximally valued fields. Indeed, if $\mathcal{O}$ is local, then the quotient domain $\mathcal{O}/\mathfrak{p}_{\mathfrak{v}}$ has Krull dimension $1$ by \cite[Theorem 4.4.6]{Bruns-Herzog}, where $\mathfrak{p}_{\mathfrak{v}}$ denotes the so-called `core' $\{f \in \mathcal{O}\,:\,\mathfrak{v}(f)=\infty\}$. Therefore the natural extension of $\mathfrak{v}$ to ${\rm Frac}(\mathcal{O}/\mathfrak{p}_{\mathfrak{v}})$ is automatically $0$-dimensional, hence we can apply Kaplansky's Embedding Theorem \cite[Theorem 6]{Kap}: if the residue field $k$ of $\mathfrak{v}$ is algebraically closed of characteristic $0$, then there is an embedding ${\rm Frac}(\mathcal{O}/\mathfrak{p}_{\mathfrak{v}}) \subset k((t))$ in the field of formal Laurent series (see also \cite[section 4]{CutVal} and an `embedded' version in \cite{TeissierKaplansky}).

\bekezdes In the sense in which discrete valuations correspond to (formal) arcs, the analogues of quasi-valuations are the $m$-jets. Indeed, these latter objects are defined as follows: for $m \in \mathbb{N}$, a $k$-morphism ${\rm Spec}\,K[t]/(t^{m+1}) \rightarrow {\rm Spec} \, \mathcal{O}$ (i.e., a homomorphism $\mathcal{O} \rightarrow K[t]/(t^{m+1})$) is called an \textit{$m$-jet}. These ring maps post-composed with the suitable order function $ord_t:K[t]/(t^{m+1}) \rightarrow \mathbb{N}_{m+1}$ (with $ord_t(0)=m+1$) give quasi-valuations on $\mathcal{O}$ in the sense of Definition \ref{def:quasi-val}.

This connection helps us to reformulate  Proposition \ref{prop:VAL} 
 (or \cite[2.1 Théorème]{LT}) and answer Question \ref{q:qreal=real} affirmatively in the complex analytic case as follows:

\begin{cor}\label{cor:intcl=real=qreal}
   Let $\mathcal{I}$ be an ideal in the convergent power series ring $\mathcal{O}_m=\mathbb{C}\{x_1, \ldots, x_m\}$. Then
 \begin{center}
     $\mathcal{I}$ is integrally closed $\Leftrightarrow$ $\mathcal{I}$ is realizable $\Leftrightarrow$ $\mathcal{I}$ is quasi-realizable. 
 \end{center}
\end{cor}

Note however, that not all quasi-valuations come from jets:

\begin{example}
    There exists no jet $\mathfrak{j}:k[x,y,z,t]/(x,y,z,t)^{3} \rightarrow K[\tau]/(\tau^{13})$ such that post-composing with the order function $ord_{\tau}:K[\tau]/(\tau^{13})\rightarrow \mathbb{N}_{13}$ gives the quasi-valuation $\mathfrak{q}$ from Example \ref{ex:quasinonext}. Indeed, any such mapping must necessarily have
    \begin{align*}
        & \mathfrak{j}(x)=a_x \tau^5 + \text{ higher order terms } & \mathfrak{j}(y)=a_y \tau^5 + \text{ higher order terms } \\
        & \mathfrak{j}(z)=a_z \tau^5 + \text{ higher order terms } & \mathfrak{j}(t)=a_t \tau^5 + \text{ higher order terms } 
    \end{align*}
    with $a_x, a_y, a_z, a_t \in K\setminus \{0\}$ and
    \begin{align*}
        &ord_\tau(\mathfrak{j}(xt-yz))=\mathfrak{q}(xt-yz)=11, \text{ thus } a_xa_t=a_ya_z \text{ and }\\
        &ord_\tau(\mathfrak{j}(yt-z^2))=\mathfrak{q}(yt-z^2)=12, \text{ thus } a_ya_t=a^2_z.
    \end{align*}
    But these imply that $a_xa_z=a_y^2$, since $a_x a_z a_t^2= a_y a_z^2 a_t = a_y^2 a_t^2$, and hence
    \begin{equation*}
        ord_\tau(\mathfrak{j}(xz-y^2)) > 10 = \mathfrak{q}(xz-y^2), \text{ contradiction}.
    \end{equation*}
\end{example}

\begin{question}\label{q:quasi-jet}
    How can we characterize those quasi-valuations which come from jets? Can all quasi-realizable ideals be quasi-realized using only quasi-valuations coming from jets?
\end{question}

\section{Appendix}

\subsection{Realizable and integrally closed submodules}\,

We recall here Rees's original definition of integral dependence over a submodule. First we describe the setting and fix some notations. Let $\mathcal{O}$ be a Noetherian ring and $N \leq M$ an
inclusion of finitely generated $\mathcal{O}$-modules. 
For any prime ideal $\mathfrak{p} \triangleleft\mathcal{O}$ let us denote the fraction field of $\mathcal{O}/\mathfrak{p}$ by $\kappa(\mathfrak{p})$, the $\kappa(\mathfrak{p})$-vector space $M \otimes_{\mathcal{O}}\kappa(\mathfrak{p}) \cong M_{\mathfrak{p}}/\mathfrak{p} M_{\mathfrak{p}}$ by $W(\mathfrak{p})$ and the image of $N \otimes_{\mathcal{O}} \kappa(\mathfrak{p})$ in $W(\mathfrak{p})$ by $V(\mathfrak{p})$. We also denote by $\iota_\mathfrak{p}$ the canonical composition $M \rightarrow M/\mathfrak{p}M \rightarrow W(\mathfrak{p})$.

\begin{define}\label{def:intdep}(\cite{ReesMod}, see also \cite[Section 16]{SH})
An element $m \in M$ is said to be
\textit{integral over $N$} if for every \textit{minimal} prime ideal $\mathfrak{p}$ in $\mathcal{O}$ and every discrete valuation
ring $\mathcal{O}_{\mathfrak{v}}$ between $\mathcal{O}/\mathfrak{p}$ and $\kappa(\mathfrak{p})$ (i.e., the corresponding discrete valuation $\mathfrak{v}$ satisfies $\mathfrak{v} \geq 0$ on $\mathcal{O}/\mathfrak{p}$ and $\mathfrak{v}(a)=\infty$ if and only if $a=0$ in $\mathcal{O}/\mathfrak{p}$) we have 
\begin{center}$\iota_{\mathfrak{p}}(m) \in \iota_{\mathfrak{p}}(N)\mathcal{O}_{\mathfrak{v}}$,
\end{center}
i.e., the image of $m$ in $M_{\mathfrak{p}}/\mathfrak{p} M_{\mathfrak{p}}$ can be written as $\iota_{\mathfrak{p}}(m)= \sum_{i}v_i n_i$, where $v_i \in \mathcal{O}_{\mathfrak{v}}$ and $n_i$ lies in the canonical image of $N$
inside $W(\mathfrak{p})$ for each $i$.

The submodule $N$ is called \textit{integrally closed} in $M$, if every element $m \in M$, which is integral over $N$, already lies in $N$.
\end{define}

We want to translate this definition to our language of extended discrete valuations on the finitely generated $\mathcal{O}$-module $M$. 

\begin{lemma}
    Consider a finitely generated module $M$ over the Noetherian ring  $\mathcal{O}$ and a discrete valuation $\mathfrak{v}: \mathcal{O} \rightarrow \overline{\mathbb{N}}=\mathbb{N}\cup\{ \infty \}$ (in the sense of Definition \ref{def:gdv}) with $W(\mathfrak{p}_{\mathfrak{v}}) \neq 0$ (where $\mathfrak{p}_{\mathfrak{v}}$ denotes the core $\{f \in \mathcal{O}\,:\,\mathfrak{v}(f)=\infty\}$ of the valuation). Then there exists some (in fact, there exist infinitely many) map $\mathfrak{v}^M:M \rightarrow \overline{\mathbb{Z}}=\mathbb{Z} \cup \{\infty\}$, such that the pair $(\mathfrak{v}, \mathfrak{v}^M)$ is an extended discrete valuation (in the sense of Definition \ref{def:edv}).
\end{lemma}

\begin{proof}
    Any such valuation $\mathfrak{v}: \mathcal{O} \rightarrow \overline{\mathbb{N}}$ with core $\mathfrak{p}_{\mathfrak{v}} \triangleleft \mathcal{O}$ descends to a discrete valuation of $\mathcal{O}/\mathfrak{p}_{\mathfrak{v}}$ which in turn extends uniquely as $\mathfrak{v}^\kappa$ to its field of fractions $\kappa(\mathfrak{p}_{\mathfrak{v}})$. We can now consider the finite dimensional $\kappa(\mathfrak{p}_{\mathfrak{v}})$-vector space $W(\mathfrak{p}_{\mathfrak{v}})\neq 0$. One can then extend this valuation to this vector space to give it a valued vector space structure in the sense of \cite[Definition 3.3]{valued-VS} (this corresponds to our Definition \ref{def:edv} if we consider $W(\mathfrak{p}_{\mathfrak{v}})$ as a module over $\kappa(\mathfrak{p}_{\mathfrak{v}})$). Indeed, we can consider any basis $\{ w_1, \ldots, w_s\}$ of $W(\mathfrak{p}_{\mathfrak{v}})$, assign any values $w_i \mapsto a_i \in \mathbb{Z}$ for all $1 \leq i \leq s$ and then set $\mathfrak{v}^W: W(\mathfrak{p}_{\mathfrak{v}}) \rightarrow \overline{\mathbb{Z}}$ as 
    \begin{equation}\label{eq:minextension}
        w=\sum_i \alpha_i w_i \mapsto \mathfrak{v}^W(w)=\min \{ \mathfrak{v}^\kappa(\alpha_i)+a_i\,:\,1 \leq i\leq s\} \hspace{3mm} (\text{where } \alpha_i \in \kappa(\mathfrak{p}_{\mathfrak{v}}) \text{ for all } 1 \leq i \leq s).
   \end{equation}
    By direct computation (see also \cite[Lemma 3.5]{valued-VS}) this will satisfy the required properties and than we can define the extension $\mathfrak{v}^M: M \rightarrow \overline{\mathbb{Z}}$ as $m \mapsto \mathfrak{v}^W(\iota_{\mathfrak{p}_{\mathfrak{v}}}(m))$. On can easily check, that (independently from the chosen values $\{ a_i\}_{i=1}^s$) the pair $(\mathfrak{v}, \mathfrak{v}^M)$ gives indeed an extended discrete valuation.
\end{proof}

Even more, we also have the following:

\begin{remark}\label{rem:exttoMandW}
    The extensions $\mathfrak{v}^M:M \rightarrow \overline{\mathbb{Z}}$ of a discrete valuation $\mathfrak{v}: \mathcal{O} \rightarrow \overline{\mathbb{N}}$ correspond bijectively to the extensions $\mathfrak{v}^W: W(\mathfrak{p}_{\mathfrak{v}}) \rightarrow \overline{\mathbb{Z}}$ of $\mathfrak{v}^\kappa: \kappa(\mathfrak{p}_{\mathfrak{v}}) \rightarrow \overline{\mathbb{Z}}$, where $\mathfrak{p}_{\mathfrak{v}}$ is the core of the valuation.
    Indeed, any extension $\mathfrak{v}^W$ can be pulled back to $M$ via the canonical map $\iota_{\mathfrak{p}_{\mathfrak{v}}}$, whereas any extension $\mathfrak{v}^M$ descends to $\iota_{\mathfrak{p}_{\mathfrak{v}}}(M)$ (cf. part  (iii) of Remark \ref{rem:translating}), from which we get  $W(\mathfrak{p}_{\mathfrak{v}})$ by localization with the multiplicative set $(\mathcal{O}/\mathfrak{p}_{\mathfrak{v}}) \setminus (\mathfrak{p}_{\mathfrak{v}}/\mathfrak{p}_{\mathfrak{v}})$. The well-definedness of the extension can be checked by direct computation.
\end{remark}

We can now state the following reformulation:

\begin{prop}\label{prop:intdepwithedvs}
    For any fixed discrete valuation $\mathfrak{v}: \mathcal{O} \rightarrow \overline{\mathbb{N}}$ we have the following equivalence: for any $m \in M$ we have 
    \begin{center}
$\iota_{\mathfrak{p}_{\mathfrak{v}}}(m) \in \iota_{\mathfrak{p}_{\mathfrak{v}}}(N) \mathcal{O}_{\mathfrak{v}}$ if and only if $\mathfrak{v}^M(m) \geq \mathfrak{v}^M(N)$ for every extension $\mathfrak{v}^M:M \rightarrow \overline{\mathbb{Z}}$
    \end{center}
    (in the sense of Definition \ref{def:edv}), where $\mathcal{O}_v=\{ \alpha \in \kappa(\mathfrak{p}_{\mathfrak{v}})\,:\, \mathfrak{v}^{\kappa}(\alpha)\geq 0\}$.
\end{prop}

\begin{remark}
    Notice that, since $N$ is finitely generated over the Noetherian ring $\mathcal{O}$, we have 
    \begin{center}
    $\mathfrak{v}^M(N)=\min\{\mathfrak{v}^M(n)\,:\, n \in N\}>-\infty$.
    \end{center}
\end{remark}

\begin{proof}[Proof of Proposition \ref{prop:intdepwithedvs}]
    If $\iota_{\mathfrak{p}_{\mathfrak{v}}}(m) \in \iota_{\mathfrak{p}_{\mathfrak{v}}}(N) \mathcal{O}_{\mathfrak{v}}$, i.e., $\iota_{\mathfrak{p}_{\mathfrak{v}}}(m)= \sum_{i=1}^{k}v_i \iota_{\mathfrak{p}_\mathfrak{v}}(n_i)$, with $v_i \in \mathcal{O}_{\mathfrak{v}}$ and $n_i\in N$, then \begin{equation*}
        \mathfrak{v}^M(m)= \mathfrak{v}^W(\iota_{\mathfrak{p}_{\mathfrak{v}}}(m))\geq \min\{ \mathfrak{v}^\kappa(v_i) + \mathfrak{v}^M(n_i)\, :\, 1\leq i \leq k\} \geq \mathfrak{v}^M(N),
    \end{equation*}
    as $\mathcal{O}_{\mathfrak{v}}=\{ \alpha \in \kappa(\mathfrak{p}_{\mathfrak{v}})\,:\, \mathfrak{v}^\kappa(\alpha)\geq 0\}$ (where $\mathfrak{v}^W:W(\mathfrak{p}_{\mathfrak{v}})\rightarrow \overline{\mathbb{Z}}$ is the corresponding extension of $\mathfrak{v}^M$ from Remark \ref{rem:exttoMandW}).

    Now suppose indirectly that $\iota_{\mathfrak{p}_{\mathfrak{v}}}(m) \notin \iota_{\mathfrak{p}_{\mathfrak{v}}}(N)\mathcal{O}_\mathfrak{v}$. We will construct an extension $\mathfrak{v}^M: M \rightarrow \overline{\mathbb{Z}}$ of $\mathfrak{v}:\mathcal{O} \rightarrow \overline{\mathbb{N}}$, such that $\mathfrak{v}^M(m) < 0 =\mathfrak{v}^M(N)$. The main technical ingredient is the following:

    \begin{lemma} \label{lem:nicebasis}
        Consider a discrete valuation $\mathfrak{v}: \mathcal{O} \rightarrow \overline{\mathbb{N}}$ with corresponding discrete va\-luation ring $\mathcal{O}/\mathfrak{p}_{\mathfrak{v}} \leq \mathcal{O}_{\mathfrak{v}} \leq \kappa(\mathfrak{p}_{\mathfrak{v}})$ for the core ideal $\mathfrak{p}_{\mathfrak{v}}$ and set $\dim_{\kappa(\mathfrak{p}_{\mathfrak{v}})}V(\mathfrak{p}_{\mathfrak{v}}):=s_N$. Then there exists a $\kappa(\mathfrak{p}_{\mathfrak{v}})$-basis of $V(\mathfrak{p}_{\mathfrak{v}})$ of the form $\{ \iota_{\mathfrak{p}_{\mathfrak{v}}}(n_1), \ldots, \iota_{\mathfrak{p}_{\mathfrak{v}}}(n_{s_N})\}$ with $n_i \in N$ for all indices $i\in \{1, \ldots, s_N\}$, such that $\iota_{\mathfrak{p}_{\mathfrak{v}}}(N) \subset \sum_{i=1}^{s_N} \big(\iota_{\mathfrak{p}_{\mathfrak{v}}}(n_i)\mathcal{O}_{\mathfrak{v}}\big)$, i.e., for every $n \in N$
        \begin{equation*}
            \iota_{\mathfrak{p}_{\mathfrak{v}}}(n) = \sum_{i=1}^{s_N} v_i\, \iota_{\mathfrak{p}_{\mathfrak{v}}}(n_i) \text{ with } v_i \in \mathcal{O}_{\mathfrak{v}} \ (\forall i \in \{ 1, \ldots, s_N\}).
        \end{equation*}
    \end{lemma}
    
    We postpone the proof of this lemma in order to finish the proof of Proposition \ref{prop:intdepwithedvs} first. 

    Consider thus the basis $\{ \iota_{\mathfrak{p}_{\mathfrak{v}}}(n_1), \ldots, \iota_{\mathfrak{p}_{\mathfrak{v}}}(n_{s_N})\}$ of $V(\mathfrak{p}_{\mathfrak{v}})$ given by Lemma \ref{lem:nicebasis} and extend it to some $\kappa(\mathfrak{p}_\mathfrak{v})$-basis $\{ \iota_{\mathfrak{p}_{\mathfrak{v}}}(n_1), \ldots, \iota_{\mathfrak{p}_{\mathfrak{v}}}(n_{s_N}), w_{s_N+1}, \ldots, w_s\}$ of $W(\mathfrak{p}_{\mathfrak{v}})$. Then
    \begin{equation*}
        \iota_{\mathfrak{p}_{\mathfrak{v}}}(m)=\sum_{i=1}^{s_N} \alpha_i \, \iota_{\mathfrak{p}_{\mathfrak{v}}}(n_i) + \sum_{i=s_N+1}^{s}\beta_{i}\,w_i \text{ for some } \alpha_i, \beta_i \in \kappa(\mathfrak{p}_{\mathfrak{v}}) \text{ for all }1 \leq i \leq s.
    \end{equation*}
    
    Suppose first, that $\iota_{\mathfrak{p}_{\mathfrak{v}}}(m) \notin V(\mathfrak{p}_{\mathfrak{v}})$. Then there exists some $i_0 \in \{ s_N+1, \ldots, s\}$ such that $ \beta_{i_0} \neq 0$, or equivalently, $\mathfrak{v}^\kappa(\beta_{i_0}) \neq \infty$. Then we can construct an extension $\mathfrak{v}^W:W(\mathfrak{p}_{\mathfrak{v}}) \rightarrow \overline{\mathbb{Z}}$ as follows. We fix the following values on the basis elements:
    \begin{align*}
        \forall \,1 \leq i \leq s_N:& \ \mathfrak{v}^W(\iota_{\mathfrak{p}_{\mathfrak{v}}}(n_i))=0, \\ 
        \forall\, s_N+1 \leq i \leq s, i \neq i_0:& \ \mathfrak{v}^W(w_i)=0 \\ \text{ and } &\ \mathfrak{v}^W(w_{i_0})< -\mathfrak{v}^\kappa(\beta_{i_0});
    \end{align*}
    and then extend it to any $\kappa(\mathfrak{p}_{\mathfrak{v}})$-linear combination as in (\ref{eq:minextension}). Clearly, by Lemma \ref{lem:nicebasis}, in this case $\mathfrak{v}^M(m) < 0 =\mathfrak{v}^M(N)$ for the composition $\mathfrak{v}^M=\mathfrak{v}^W\circ \iota_{\mathfrak{p}_{\mathfrak{v}}}$. 

    Suppose now, that $\iota_{\mathfrak{p}_{\mathfrak{v}}}(m) \in V(\mathfrak{p}_{\mathfrak{v}}) \setminus \iota_{\mathfrak{p}_{\mathfrak{v}}}(N) \mathcal{O}_{\mathfrak{v}}$. Then  $\iota_{\mathfrak{p}_{\mathfrak{v}}}(m)=\sum_{i=1}^{s_N} \alpha_i \, \iota_{\mathfrak{p}_{\mathfrak{v}}}(n_i)$, with $\alpha_i \in \kappa(\mathfrak{p}_{\mathfrak{v}})$ for all $i$, and there exists some index $1 \leq i_0 \leq s_N$ such that $\alpha_{i_0} \notin \mathcal{O}_{\mathfrak{v}}$ (or, equivalently, $\mathfrak{v}^\kappa(\alpha_{i_0})<0$). 
    Then we can construct the extension $\mathfrak{v}^W:W(\mathfrak{p}_{\mathfrak{v}}) \rightarrow \overline{\mathbb{Z}}$ by fixing the value $0$ on each basis element of $\{ \iota_{\mathfrak{p}_{\mathfrak{v}}}(n_1), \ldots, \iota_{\mathfrak{p}_{\mathfrak{v}}}(n_{s_N}), w_{s_N+1}, \ldots, w_s\}$
    and then extending it to any $\kappa(\mathfrak{p}_{\mathfrak{v}})$-linear combination as in (\ref{eq:minextension}). Once again, by Lemma \ref{lem:nicebasis}, $\mathfrak{v}^M(m) < 0 =\mathfrak{v}^M(N)$ for the composition $\mathfrak{v}^M=\mathfrak{v}^W\circ \iota_{\mathfrak{p}}$. 
\end{proof}

\begin{proof}[Proof of Lemma \ref{lem:nicebasis}]
    We can start with any generating set $\{ n'_1, \ldots, n'_{r_N}\}$ of $N$ above $\mathcal{O}$. We can also suppose without loss of generality, that $\{ \iota_{\mathfrak{p}_{\mathfrak{v}}}(n'_1), \ldots, \iota_{\mathfrak{p}_{\mathfrak{v}}}(n'_{s_N})\}$ give a $\kappa(\mathfrak{p}_{\mathfrak{v}})$-basis of $V(\mathfrak{p}_{\mathfrak{v}})$. 
    
    Then $\iota_{\mathfrak{p}_{\mathfrak{v}}}(n'_{r_N})=\sum_{i=1}^{s_N} \alpha_i^{r_N}\, \iota_{\mathfrak{p}_{\mathfrak{v}}}(n'_i)$ for some  $\alpha_i^{r_N} \in \kappa(\mathfrak{p}_{\mathfrak{v}}), \ 1 \leq i \leq s_N$. Set also $\alpha_{r_N}^{r_N}:=-1$. Consider now the index $i_{r_N}\in\{1, \ldots, s_N, r_N\}$ distinguished by the property that $\mathfrak{v}^\kappa(\alpha_{i_{r_N}}^{r_N}) = \min \{ \mathfrak{v}^\kappa(\alpha_{i}^{r_N})\}_i \neq \infty$. Automatically $\alpha_{i_{r_N}}^{r_N} \neq 0 \in \kappa(\mathfrak{p}_{\mathfrak{v}})$, and hence 
    \begin{equation*}
        \iota_{\mathfrak{p}_{\mathfrak{v}}}(n'_{i_{r_N}})= \sum_{i \neq i_{r_N}} - \frac{\alpha_i^{r_N}}{\alpha_{i_{r_N}}^{r_N}} \iota_{\mathfrak{p}_{\mathfrak{v}}}(n'_i), \text{ with } \mathfrak{v}^\kappa \left(- \frac{\alpha_i^{r_N}}{\alpha_{i_{r_N}}^{r_N}}  \right) \geq 0, \text{ hence } \frac{\alpha_i^{r_N}}{\alpha_{i_{r_N}}^{r_N}}  \in \mathcal{O}_{\mathfrak{v}}.
    \end{equation*}
    Let us then fix $n_{r_N}:=n'_{i_{r_N}}$ and swap the indices of $n'_{{r_N}}$ and $n'_{i_{r_N}}$. Note that in this way we obtain that $\iota_{\mathfrak{p}_{\mathfrak{v}}}(n_{r_N}) \in \sum_{i=1}^{r_N-1} \big(\iota_{\mathfrak{p}_{\mathfrak{v}}}(n'_{i}) \mathcal{O}_{\mathfrak{v}}\big)$ and $\{ \iota_{\mathfrak{p}_{\mathfrak{v}}}(n'_1), \ldots, \iota_{\mathfrak{p}_{\mathfrak{v}}}(n'_{s_N})\}$ is still a $\kappa(\mathfrak{p}_{\mathfrak{v}})$-basis of $V(\mathfrak{p}_{\mathfrak{v}})$. 

    Next, we repeat the above process with $n'_{r_N-1}$, obtaining $n_{r_{N}-1}:=n'_{i_{r_N-1}}$ and a new ordering of the subset $\{ n'_1, \ldots, n'_{r_N-1}\}$ of the generators of $N$, such that 
    \begin{center}
        $\iota_{\mathfrak{p}_{\mathfrak{v}}}(n_{r_N-1})\in \sum_{i=1}^{r_N-2} \big(\iota_{\mathfrak{p}_{\mathfrak{v}}}(n'_i)\mathcal{O}_{\mathfrak{v}}\big)$ and $\{ \iota_{\mathfrak{p}_{\mathfrak{v}}}(n'_1), \ldots, \iota_{\mathfrak{p}_{\mathfrak{v}}}(n'_{s_N})\}$ is a $\kappa(\mathfrak{p}_{\mathfrak{v}})$-basis of $V(\mathfrak{p})$\\
        \vspace{2mm}
        (and it is still true that $\iota_{\mathfrak{p}_{\mathfrak{v}}}(n_{r_N}) \in \sum_{i=1}^{r_N-1} \big(\iota_{\mathfrak{p}_{\mathfrak{v}}}(n'_{i}) \mathcal{O}_{\mathfrak{v}}\big)$).
    \end{center} We repeat this process until we get to $n_{s_N+1}$ and then set $n_i:=n'_i$ for all $1 \leq i \leq s_N$. Thus we arrived to a new ordering $\{ n_1, \ldots, n_{r_N}\}$ of the original generating set $\{n'_1, \ldots, n'_{r_N}\}$ of $N$ such that $\{ n_1, \ldots, n_{s_N}\}$ is a  $\kappa(\mathfrak{p}_{\mathfrak{v}})$-basis of $V(\mathfrak{p}_{\mathfrak{v}})$ and for every $s_{N}+1 \leq i \leq r_N$ we have the property that $\iota_{\mathfrak{p}_{\mathfrak{v}}}(n_{i}) \in \sum_{i=1}^{i-1} \big(\iota_{\mathfrak{p}_{\mathfrak{v}}}(n_i) \mathcal{O}_\mathfrak{v}\big)$. As $\mathcal{O}_{\mathfrak{v}} \geq \mathcal{O}/ \mathfrak{p}_{\mathfrak{v}}$, the $\iota_{\mathfrak{p}_{\mathfrak{v}}}$-image of every $\mathcal{O}$-linear combination of these $n_i$-s (i.e., the $\iota_{\mathfrak{p}_{\mathfrak{v}}}$-image of every element of $N$) lies in the $\mathcal{O}_{\mathfrak{v}}$-submodule $\sum_{i=1}^{s_N} \big(\iota_{\mathfrak{p}_{\mathfrak{v}}}(n_i) \mathcal{O}_\mathfrak{v}\big)$.
\end{proof}

Proposition \ref{prop:intdepwithedvs} combined with the original Definition \ref{def:intdep} of integral dependence yields the following:

\begin{cor}\label{cor:integralviaextendedval}
    An element $m \in M$ is integral over $N$ if and only if for every discrete valuation $\mathfrak{v}: \mathcal{O} \rightarrow \overline{\mathbb{N}}$, with its core $\mathfrak{p}_{\mathfrak{v}}:=\{ f \in \mathcal{O}\,:\, \mathfrak{v}(f)=\infty\}$ being a minimal prime ideal in $\mathcal{O}$, and every extension $\mathfrak{v}^M: M \rightarrow \overline{\mathbb{Z}}$ in the sense of Definition \ref{def:edv}, we have $\mathfrak{v}^M(m) \geq \mathfrak{v}^M(N)$. 
\end{cor}

\begin{cor}[= Theorem \ref{th:REES}]\label{cor:Rees=Real}
    If $k$ is a field, $\mathcal{O}$ a Noetherian $k$-algebra, $N \leq M$ an inclusion of finitely generated $\mathcal{O}$-modules with ${\rm codim}_k(N\hookrightarrow M)<\infty$, then $N$ is integrally closed in $M$ if and only if there exists a (finite) collection $\mathcal{D}$ of extended discrete valuations \emph{with their corresponding discrete valuation rings lying between $\mathcal{O}/\mathfrak{p}$ and $\kappa(\mathfrak{p})$ for some minimal prime ideal $\mathfrak{p}\triangleleft \mathcal{O}$}, such that $N = \mathcal{F}_{\mathcal{D}}^M(0)$ (for the definition of $\mathcal{F}_{\mathcal{D}}^M$ see (\ref{eq:fdMl})).
\end{cor}

\begin{proof}
    By Corollary \ref{cor:integralviaextendedval},  $N$ is integrally closed if and only if for every $m \in M \setminus N$ there exists a discrete valuation $\mathfrak{v}: \mathcal{O} \rightarrow \overline{\mathbb{N}}$, with its core $\mathfrak{p}_{\mathfrak{v}}$ a minimal prime ideal in $\mathcal{O}$, and an extension $\mathfrak{v}^M: M \rightarrow \overline{\mathbb{Z}}$ such that $\mathfrak{v}^M(m) < \mathfrak{v}^M(N)$. However, by finite codimensionality, this is equivalent with the existence of a finite collection $\mathcal{D}=\{ (\mathfrak{v}_1, \mathfrak{v}_1^M),  \ldots, (\mathfrak{v}_r, \mathfrak{v}_r^M)\}$ of extended discrete valuations, with their cores $\{\mathfrak{p}_{\mathfrak{v}_v}\}_{v=1}^r$ being minimal prime ideals, and satisfying $N=\{m \in M\,:\, \mathfrak{v}_v^M(m) \geq \mathfrak{v}_v^M(N) \text{ for all }1 \leq v \leq r\}$. The ability to translate the extensions $\{\mathfrak{v}_v^M\}_v$ by the values $\{\mathfrak{v}_v^M(N)\}_v$ (cf. Remark \ref{rem:translating} (i)) implies the final form of the equivalence.
\end{proof}

However, we do not know the answer for the following question:

\begin{question}
    Suppose that $N \leq M$ is a finite $k$-codimensional realizable submodule, i.e., there exists a finite collection $\mathcal{D}$ of extended discrete valuations (without any restriction on their core) with $N=\mathcal{F}_{\mathcal{D}}^M(0)$. Is then $N$ integrally closed in $M$? What if $M$ is a free 
    $\mathcal{O}$-module?
\end{question}

\bekezdes \textbf{The local complex analytic case.}\label{par:locan}\,

We can also consider the local analytic setting of paragraph \ref{bek:locan}, with $\cO=\cO_{X,o}$, the local $\C$-algebra of a complex analytic spacegerm $(X, o)$, and $M=\cO^p$ for some $p\geq 1$. Recall, that --- by the definition of Gaffney in \cite{Gaff,GaffKl} --- an element $x \in \cO_{X,o}^p$ is integral (in the sense of Gaffney) over a submodule $N\leq M$ if for all holomorphic map germs $\phi:(\C,0)\to (X,o)$ the induced morphism $\phi^{*}_M:M=\mathcal{O}^p \rightarrow \mathcal{O}_{\mathbb{C}, 0}^p=\mathbb{C}\{t\}^p$ maps $x$ to the submodule $(\phi_M^*(N))\cO_{\C,0}$ generated over $\cO_{\C,0}$ by the image of $N$.
Using our language of extended discrete valuations we can describe the notion of integral closedness (in the sense of Gaffney) in the following way:

\begin{theorem}[= Theorem \ref{th:Gaffney}]
    Let $\mathcal{O}=\mathcal{O}_{X, o}$, $M=\mathcal{O}^p$ for some $p \geq 1$ and $N \leq M$ a finite codimensional finitely generated submodule. Then $N$ is integrally closed (in the sense of Gaffney) in $M$ if and only if there exists a finite set of holomorphic map germs $\{\phi_{v}:(\mathbb{C}, 0) \rightarrow (X, o)\}_{v=1}^r$ and extensions $\{{\rm ord}_t^v:\mathcal{O}_{\mathbb{C} ,0}^p\cong \mathbb{C}\{t\}^p\rightarrow\overline{\mathbb{Z}}\}_{v=1}^r$ of the standard order function ${\rm ord}_t:\mathcal{O}_{\mathbb{C} ,0}\cong\mathbb{C}\{t\}\rightarrow\overline{\mathbb{N}}$, such that the collection $\mathcal{D}=\{ (\mathfrak{v}_{\phi_v}:={\rm ord}_t \circ \phi^*, \mathfrak{v}_{\phi_v}^M:={\rm ord}_t^v\circ \phi^*_M)\}_{v =1}^r$ satisfies $N = \mathcal{F}_{\mathcal{D}}^M(0)$.
\end{theorem}

\begin{proof}
In fact, we only need the following local analytic counterpart of Corollary \ref{cor:integralviaextendedval}:

\begin{cor}\label{cor:integralGviaextendedval}
An element $x \in \mathcal{O}_{X, o}^p$ is integral over $N$ if and only if for every holomorphic map germ $\phi:(\mathbb{C}, 0) \rightarrow (X, o)$, and every extension ${\rm ord}_t^{\mathbb{C}\{t\}^p}:  \mathbb{C}\{t\}^p\rightarrow \overline{\mathbb{Z}}$ of the standard order function ${\rm ord}_t:\mathcal{O}_{\mathbb{C} ,0}\cong\mathbb{C}\{t\}\rightarrow\overline{\mathbb{N}}$, we have ${\rm ord}_t^{\mathbb{C}\{t\}^p}\circ\phi^*_M(x) \geq {\rm ord}_t^{\mathbb{C}\{t\}^p}\circ \phi^*_M(N)$. 
\end{cor}

\begin{proof} One just has to apply Proposition \ref{prop:intdepwithedvs} to the local complex analytic algebra $\mathcal{O}_{\mathbb{C}, 0}$, the standard order valuation ${\rm ord}_t:\mathcal{O}_{\mathbb{C},0} \rightarrow \overline{\mathbb{N}}$, the free module $\mathcal{O}_{\mathbb{C}, 0}^p$ and the submodule $\phi^*_M(N)\mathcal{O}_{\mathbb{C}, 0}$ to get that for any $x \in M=\mathcal{O}_{X, o}^p$:  \begin{center}$\phi^*_M(x) \in \phi^*_M(N)\mathcal{O}_{\mathbb{C}, 0}$ if and only if ${\rm ord}_t^{\mathbb{C}\{t\}^p}(\phi^*_M(x)) \geq {\rm ord}_t^{\mathbb{C}\{t\}^p}( \phi^*_M(N)\mathcal{O}_{\mathbb{C}, 0})$
\end{center}
 for every extension ${\rm ord}_t^{\mathbb{C}\{t\}^p}:  \mathbb{C}\{t\}^p\rightarrow \overline{\mathbb{Z}}$. Indeed, in this case $\mathfrak{p}_{{\rm ord}_t}=(0) \triangleleft \mathcal{O}_{\mathbb{C}, 0}$.
\end{proof}

To finish the proof of Theorem \ref{th:Gaffney}, one just has to use the same argument as in the case of Corollary \ref{cor:Rees=Real} $=$ Theorem \ref{th:REES}, relying on Corollary \ref{cor:integralGviaextendedval} instead of Corollary \ref{cor:integralviaextendedval}.
\end{proof}

Notice, that in order to be able to state the above Theorem \ref{th:Gaffney} and Corollary \ref{cor:integralGviaextendedval} with more general extensions $\mathfrak{v}_\phi^{M}:M \rightarrow \overline{\mathbb{Z}}$ of $\mathfrak{v}_\phi={\rm ord}_t \circ \phi^*:\mathcal{O} \rightarrow \overline{N}$ (not necessarily factoring through the pullback map $\phi^*_M:M=\mathcal{O}_{X, o}^p \rightarrow \mathbb{C}\{t\}^p$ by design) we would need the following statement (see also the discussion after formula (\ref{eq:curvemapval})):

\begin{question}\label{q:Gaffney}
    Consider an extension of valued fields $(K, \mathfrak{v}) \leq (L, \mathfrak{v}_L)$ (i.e., $\mathfrak{v}_L: L \rightarrow \overline{\mathbb{Z}}$ is a discrete valuation with trivial core, such that $\mathfrak{v}_L \big|_K=\mathfrak{v}$) and a (finite dimensional) $K$-vector space $V$ with an extension $\mathfrak{v}^V:V \rightarrow \overline{\mathbb{Z}}$ of $\mathfrak{v}$ in the sense of Definition \ref{def:edv} (see also \cite[Definition 3.3]{valued-VS}). (Under what conditions) is it true that  there exists an extension $\mathfrak{v}_L^{V \otimes_K L}:V\otimes_K L \rightarrow \overline{\mathbb{Z}}$ of $\mathfrak{v}_L$ (in the sense of Definition \ref{def:edv}), which satisfies that $\mathfrak{v}_L^{V \otimes_K L}\big|_V=\mathfrak{v}^V$?
\end{question}

Although this question seems rather natural to ask, the authors haven't managed to find any related literature. ű

\subsection{Further questions, open problems}\label{s:QOP}\,

In this subsection we list some (further) open problems, questions and conjectures related with the different lattice homological constructions presented in the previous sections, which might be guiding objectives for future research.

\bekezdes \textbf{Homological dimension.}\label{ss:HOMDIM}

Let $\mathcal{O}$ be a Noetherian $k$-algebra  and let $N$ be a finite codimensional realizable submodule in the finitely generated $\mathcal{O}$-module $M$. In this context we can introduce two numerical invariants (see also subsection
\ref{ss:bounds} and section \ref{s:homdim}): the `(lattice) homological dimension' 
\begin{equation*}
    {\rm homdim}(N \hookrightarrow M):=\begin{cases}
        \max \{ q\,:\, \mathbb{H}_{q}(N \hookrightarrow M)\not=0\} & \text{if }N\neq M;\\
        -1 & \text{if }N= M.
    \end{cases}
\end{equation*}
and the `realizability cardinality'
${\rm realiz}(N\hookrightarrow M):= \min\{|\mathcal{D}|\,:\, \text{$\mathcal{D}$ is a realization of $N$}\} $
of $N$. Proposition \ref{prop:homdegupperbound} states the inequality 
 ${\rm homdim}(N\hookrightarrow M)\leq 
{\rm realiz}(N\hookrightarrow M)-1$ as a direct consequence of our construction.

\begin{question}\label{que:1}
Is this inequality sharp? Are there examples when this is not an equality?
\end{question}

In answering this question we are obstructed by two main difficulties: in general it is not easy to compute $\bH_*$ or to
establish vanishing bounds for $\bH_q$ (they are mostly obtained exactly by the previous inequality, or by direct computations); and also  the classification of the corresponding (extended) valuations is not (always) known.

The simplest case seems to be when  ${\rm homdim}(N\hookrightarrow M)=0$. In Conjecture \ref{conj:elso} we formulated the expectation that any integrally closed finite codimensional ideal of $k[x,y]$ with homological dimension zero  can be realized  by a single valuation. In the body of the manuscript  we supported this statement by several examples, see, e.g., Examples \ref{bek:10.2.1}--\ref{bek:10.2.2}, or the case of the ideal  
 $(y^4, y^3x, y^3+x^2)\subset \bC[x,y]$ from paragraph \ref{bek:4.10.3}
 (which can be realized as  $\mathcal{F}_{\frv_{E(2,c)}}(8)$).
 
 Moreover, using the language of quasi-valuations (cf. subsection \ref{ss:quasi}), one could also study the slightly weaker conjecture, whether the vanishing of higher homologies imply \emph{quasi}-realizability by a single quasi-valuation.  We gave a supporting example for this statement in the setting of the analytic lattice homology of reduced curve singularities in Example \ref{ex:curvewith1qv}.
 
 Returning to the more general Question \ref{que:1}, in section \ref{s:homdim} we presented sufficient conditions for the equality ${\rm homdim}(N\hookrightarrow M)= 
{\rm realiz}(N\hookrightarrow M)-1$ and proved that in the case of realizable monomial ideals of $k[x,y]$ these conditions hold. We expect the same to be true in every  polynomial ring, for supporting examples see Remark \ref{rem:homdimhighdim}.

In fact, we formulated the above conjecture only for the algebra $k[x,y]$ since in this case we have a better picture of the discrete valuations. It can happen that the conjecture holds true in more general algebras as well.

 \bekezdes \textbf{Constructing minimal realizations.} 
 
 Let $\mathcal{O}$ be a Noetherian $k$-algebra  and let $N$ be a finite codimensional realizable submodule in the finitely generated $\mathcal{O}$-module $M$ as above. Related  to the previous paragraph \ref{ss:HOMDIM}, we formulate  the following question as well.

 \begin{question}\label{que:2} Is there a natural construction which provides a {\bf minimal 
     set of valuations} which realizes $N\hookrightarrow M$ ?
 \end{question}

In the case when $M=\mathcal{O}=k[x_1, x_2]$ and $N=\mathcal{M}$ is a realizable monomial ideal we gave a recipe to produce such a minimal set of realizing monomial valuations in section \ref{s:homdim} (see especially Remark \ref{rem:minrealmonid}). This depends on the strong connection between the homological dimension and the convex geometry of the Newton polytope (cf. Theorem \ref{th:homdim}). We emphasize, however, that such realizations are highly non-unique, e.g., our exact constructions in section \ref{s:homdim} depend on a choice of a small enough rational number $\varepsilon \in \mathbb{Q}_{>0}$.

We also note that  if $\mathcal{I}\subset \mathcal{O}$ is an integrally closed ideal of $\mathcal{O}$, then the very same set of 
 Rees valuations realizes all the integral closures $\overline{\mathcal{I}^n}$ of its powers $\mathcal{I}^n\subset \mathcal{O}$. 
Moreover, the set of Rees valuations is unique \cite{Rees}. 
However, if we wish to realize only  $\mathcal{I}$, 
then it turns out that it might be realized by fewer valuations than the number of Rees valuations, see, e.g., Example \ref{bek:10.2.2}.

 \bekezdes \textbf{Hilbert--Samuel-type regularity.}\label{bek:HS-type}
 
 Again, consider an integrally closed ideal $\mathcal{I}\triangleleft \mathcal{O}$ of finite codimension, and suppose that all of its powers $\mathcal{I}^n=\overline{\mathcal{I}^n}$ are also integrally closed (such ideals are called \emph{normal} in the literature and are extensively studied.)
 Then, for every $\mathcal{I}^n$ we can consider the symmetric lattice homology module 
 $\mathbb{S}\bH_*(\mathcal{I}^n\triangleleft \mathcal{O}) $. Since all of them are realized by the same set of Rees valuations \cite{Rees}, we can run the construction of the homologies in the very same lattice. Also, since the numerical sequence 
 $eu(\mathbb{S}\bH_*(\mathcal{I}^n\triangleleft \mathcal{O}))=\dim_k \mathcal{O}/\mathcal{I}^n $ shows certain regularity (at least for $n\gg 0$)
 via the Hilbert--Samuel polynomial of $\mathcal{I}$, we expect some regularity in the sequence of modules  $\mathbb{S}\bH_*(\mathcal{I}^n\triangleleft \mathcal{O}) $ as well. 

 \begin{problem} Understand the behavior of the sequence 
 of modules  $\mathbb{S}\bH_*(\mathcal{I}^n\triangleleft \mathcal{O}) $
 (at least for $n\gg 0$). 
 \end{problem}

 \bekezdes \textbf{Stability along a given multifiltration.}

More generally, we can consider the symmetric lattice homologies of all the ideals realized by a fixed collection of discrete valuations $\mathcal{D}$ of a $k$-algebra $\cO$. 
In Remark \ref{ex:Id} we introduced the notation $\mathcal{I}=\mathcal{I}(\ell):= \mathcal{F}_{\mathcal{D}}(\ell)$ for any $\ell\in
(\Z_{\geq 0})^r$. Then the modules  
${\mathbb S}\bH_*(\mathcal{I}(\ell)\triangleleft \cO)=\mathbb{SH}_*(\mathcal{O}, \mathcal{D}^{\natural}, 2\ell)$ give a categorification of the Hilbert function $\frh_{\mathcal{D}}$ associated with $\mathcal {D}$, where $\mathcal{D}^{\natural}$ is the doubling of $\mathcal {D}$ from Lemma \ref{lem:duplatrukkideal}. 

 \begin{problem} Understand the behavior of the  modules  $\ell \mapsto \mathbb{S}\bH_*(\mathcal{I}(\ell)\triangleleft \mathcal{O}) $
 (at least for $\ell\gg 0$). 
 \end{problem}

 \noindent [Notice that, using these notations, in the previous paragraph \ref{bek:HS-type} we were asking a question about the ideals $\mathcal{I}(n\cdot d_{\mathcal{I}})$, where $\mathcal{I}(d_{\mathcal{I}})=\mathcal{I}$.]

 For example, we claim that if $\mathcal{O}=\mathcal{O}_{\mathbb{C}^2, 0}$ and the valuations of $\mathcal{D}$ correspond to normalizations of embedded irreducible curve singularities $\{(C_v, 0)\}_{v=1}^r$ (cf. (\ref{eq:normalizationvaluation})), then the modules  $\mathbb{S}\bH_*(\mathcal{I}(\ell)\triangleleft \mathcal{O})$ for $\ell \geq 2\mathbf{c}(\cup_v C_v)$ can be rather explicitly described, where $\mathbf{c}(\cup_v C_v)$ is the conductor element. [A similar situation was discussed in the last part of paragraph \ref{bek:4.10.3} without the assumption on the conductor.]  The precise results in this case will be presented in a forthcoming manuscript.

\bekezdes \textbf{Connections with the theory of integrally closed ideals and modules over $2$-dimensional regular local algebras}\,

The theory described in the title is rather rich, it  was initiated by Zariski \cite{ZarComplete} (see also \cite{LipmanComplete} for a modern treatment), who proved the following  unique factorization result: the product of any two integrally closed  ideals  in a $2$-dimensional regular local ring is again integrally closed, moreover any such non-zero
ideal can be expressed uniquely (except for ordering) as a product of simple integrally closed ideals.

Notice that, using our language, if $\mathcal{I}_1$ and $\mathcal{I}_2$ are finite codimensional integrally closed ideals in the $2$-dimensional regular local $k$-algebra $\mathcal{O}$, then they are realizable. Moreover, if we choose such a collection $\mathcal{D}$ of discrete valuations on $\mathcal{O}$ which realizes simultaneously $\mathcal{I}_1, \ \mathcal{I}_2$ and $\mathcal{I}_1\mathcal{I}_2$, then we also have the fact that $\mathcal{F}_{\mathcal{D}}(d_{\mathcal{I}_1}+d_{\mathcal{I}_2})=\mathcal{I}_1\mathcal{I}_2$, where $\mathcal{F}_{\mathcal{D}}(d_{\mathcal{I}_1})=\mathcal{I}_1$ and $\mathcal{F}_{\mathcal{D}}(d_{\mathcal{I}_2})=\mathcal{I}_2$  (for the notations see section \ref{s:deccurves}).

\begin{question}
    Can we connect the symmetric lattice homology of $\mathcal{I}_1$ and $\mathcal{I}_2$ with that of $\mathcal{I}_1\mathcal{I}_2$? Does lattice homology see simpleness in this context?
\end{question}
     
 Although in higher dimensions the product of integrally closed ideals might not be so, Zariski's results were extended to some higher dimensional cases by Lipman \cite{LipmanComplete} (see also the case of the Goto-class $\mathcal{G}^*$ of integrally closed ideals in \cite{CNR}). Kodiyalam and Hayasaka extended the study to integrally closed modules \cite{Kodiyalam, KodHay1, KodHay2}. We hope that our lattice homology construction can become a useful tool in further exploration of this area.

 \bekezdes \textbf{Integral reduction and lattice homology.}\,

In section \ref{ss:ARTIN} we introduced `integrally reduced algebras', i.e., quotients of polynomial rings with integrally closed ideals. Although the notion is well defined, an intrinsic characterization (similar to that of Lejeune-Jalabert and Teissier in \cite[2.1 Théorème]{LT} in the analytic case --- see our formulation in Proposition \ref{prop:VAL}) is missing (see especially Question \ref{q:qreal=real} heading into this direction and its possible consequence (\ref{eq:intred})).  

We also introduced the `integral reduction' operation, which produces from any source algebra an integrally reduced one. The connections between the input and the outcome, though, need further study.

\begin{question}
    Which (algebraic) invariants and properties remain stable under such an integral reduction operation? Does it have any geometric significance?
\end{question}

In the case of  integrally reduced \emph{Artin} algebras, we can consider their well-defined lattice homology modules, which categorify their $k$-dimension and depend mysteriously on their multiplicative structure.

\begin{problem}
    Connect the lattice homology with other algebraic invariants of integrally reduced Artin algebras.
\end{problem}

\bekezdes \textbf{Quasi-valuations, arcs and jets}\,

In section \ref{ss:quasi} we introduced the notion of `quasi-valuations', which are generalizations of discrete valuations configured to better handle nilpotent elements and zero divisors of the base $k$-algebra. Already their definition brings up a number of interesting questions (cf. Questions \ref{q:whenextendquasi} and \ref{q:extquasi}) regarding their relations with discrete valuations and quotient operations. We also claimed that the Independence Theorem remains valid in their context, which raised the question whether this generalization really extends our playing field, i.e.,

\begin{question}
    Are all Artin algebras, arising as quotients with quasi-realizable ideals, integrally reduced?
\end{question}

In subsection \ref{ss:arcs} we also presented some analogies between discrete valuations and formal arcs, respectively between quasi-valuations and jets (see Questions \ref{q:hilbertjump01} and \ref{q:quasi-jet} arising naturally in this context). The search for connections between the symmetric lattice homology modules of integrally closed ideals of a Noetherian $k$-algebra and the formal arc-- and jet-spaces of the same algebra seems to be a fascinating project.

In the complex analytic case the authors expect some deep relation between the cohomology of the (restricted) contact loci associated with an isolated hypersurface singularity $f:(\mathbb{C}^m,0)\rightarrow(\mathbb{C}, 0)$ and its analytic lattice homology (of $m-1$-forms) $\mathbb{H}_{*}((\{f=0\},0), \Omega^{m-1}).$ This would create an additional connection with the Floer homology of compactly supported Milnor monodromy (via the results and arc-Floer conjecture of \cite{arcFloer}).

\bekezdes \textbf{Filtered version.}

The analytic lattice homology of reduced curve singularities has a filtered enhancement\,/\,version  introduced by the first author in \cite{NFilt2}, which is a much stronger invariant (e.g., it is a complete embedded topological invariant for plane curves). In its construction an additional filtration by the norm $|\,.\,|$ of the lattice points is taken into account to produce a rich spectral sequence. 

A similar approach could be taken in the case of the lattice homology of realizable submodules of section \ref{s:4} by considering the following `generalized norm'
\begin{equation*}
    \mathfrak{n}:(\mathbb{Z}_{\geq 0})^r \rightarrow \mathbb{Z}_{\geq 0},\ \ell \mapsto \dim(\mathcal{O}/\mathcal{F}(\ell)) + \dim(\mathcal{F}^M(-\ell)/N) =  \mathfrak{h}(\ell) - \mathfrak{h}^{\circ}(\ell)+{\rm codim}(N \hookrightarrow M) 
\end{equation*}
on the lattice points (compare with (\ref{eq:dag}) in the curve singularity case). One can prove, that the Independence Theorem \ref{th:IndepMod} can be strengthened to also respect the $\mathfrak{n}$-value of the local minimum and local maximum points of the weight function. On the other hand, the spectral sequence corresponding to this filtration will not necessarily be invariant. The consequences and applications of this enhancement will be discussed in a forthcoming paper. 

\bekezdes \textbf{Functoriality and exact sequences.}\,

The theory of different lattice homology constructions in its current form lacks the incorporation of functors\,/\,morphisms. In the case of the lattice homology modules associated with complex analytic singularities a lot of effort is put into assigning nontrivial graded $\mathbb{Z}[U]$-module morphisms to certain singularity deformations. The best result so far concerns the analytic lattice homology of reduced curve singularities: \cite[Theorem 6.4.1]{AgNeCurves} claims that a flat deformation $\{(C_t, o)\}_{t \in (\mathbb{C}, 0)}$ induces a degree zero graded $\mathbb{Z}[U]$-module morphism $\mathbb{H}_{an, 0}(C_{t \neq 0}, o) \rightarrow \mathbb{H}_{an, 0}(C_{t=0}, o)$ on the homological degree $0$ part. In fact, Ágoston and the first author constructed compatible maps $F_{0, n}: \pi_0(S_{n, t\neq 0}) \rightarrow \pi_0(S_{n, t=0})$ between the connected components of the $S_n$ spaces for every level $n \in \mathbb{Z}$. The `Functoriality Conjecture' claims that these $\{F_{0, n}\}_n$ maps can be extended to the $S_n$-spaces up to homotopy, such that these would induce a degree zero graded $\mathbb{Z}[U]$-module morphism $\mathbb{H}_{an, *}(C_{t \neq 0}, o) \rightarrow \mathbb{H}_{an, *}(C_{t=0}, o)$ on the whole lattice homology modules.

\begin{problem}
    Find nontrivial graded $\mathbb{Z}[U]$-module morphisms between the lattice homologies of certain realizable submodules. What information can they encode? Can we produce exact sequences?
\end{problem}

One advantage of the theory of lattice homology of realizable submodules in this regard is that  we can use different realizations inducing different lattices and weight functions to compute the same lattice homology module. In this way, if we want to compare, for example, the $S_n$-space filtrations of two realizable submodules of the same ambient module, then we can use the same collection of extended valuations realizing both. 
This might be especially helpful, since we do not really understand yet the combinatorial requirements on the weight\,/\,height functions to obtain nontrivial maps (e.g., not tautological inclusions) between different $S_n$-space filtrations, let alone exact sequences on the level of the lattice homology  modules. 

\bekezdes \textbf{Stability under singularity deformations}

Already at the introduction of the topological lattice homology associated with complex analytic singularities there arose the question of which deformations preserve this new invariant. The addition of new analytic variants just added more layers and complexity to this question. In this manuscript we introduced the notion of \emph{conductor ideal constant deformations} of Gorenstein singularities under which the analytic lattice homology (of top forms) is tautologically stable.   
We initiated their study in subsection \ref{ss:CondIdConst} through numerous examples comparing them with $\delta$--, $p_g$-- and $Z_{K,min}^2$-constant deformations, though, much is yet to be uncovered.

In subsection \ref{ss:DefsinLH} we also gave a short review of the main conjectures and open problems relating the analytic and topological lattice homologies with singularity deformations. Of particular interest is \cite[Conjecture 11.9.48]{NBook}, which claims that along flat $p_g$-constant deformations of normal surface singularities (with $\mathbb{Q}HS^3$ links) the analytic lattice homology is stable. Notice that if we omit the topological condition on the link (and use the more general analytic lattice homology of $2$-forms), then this statement becomes false: see the $p_g$-- (and $Z_{K,min}^2$--), but not $\mathbb{H}_*(\,.\,,\Omega^2)$-constant deformation in Example \ref{ex:LIST} (e) with the generic fiber having non-$\mathbb{Q}HS^3$ link. Nevertheless, \cite[Conjecture 11.9.48]{NBook} in its original form is still open.

\bekezdes \textbf{Characterization of $\bH_{\geq 1}$ and $m_w=\min w_0$ in the case of complex analytic singularities.}

In general, the presence and ranks of higher homological degree summands and the numerical invariant $m_w=\min w_0$ are rather mysterious.

\begin{problem} Give a characterization of the graded $\mathbb{Z}[U]$-modules $\{\bH_q(X,o)\}_{q\geq 1}$
and of $\min w_0$ (in the case of  the analytic lattice homology theories) associated with isolated (hypersurface) singularities of dimension $1$ and $2$ in terms of their classical singularity invariants. 
\end{problem}

Recall that in Proposition \ref{prop:htopnem0}, in the general algebraic case, we assign top dimensional homological cycles to some combinations of ring and module elements. Even though this result can be further strengthened for plane curve singularities with Newton non-degenerate principal part and convenient Newton boundary (cf. Corollary \ref{cor:htopforNNcurves}), in general, for complex analytic singularities the geometric interpretations of $\{\bH_q(X,o)\}_{q\geq 1}$ with classical invariants are still painfully missing. 

On the other hand, the invariant $m_w=\min w_0$ of reduced curve singularities was already successfully linked with their Cohen Macaulay type by A. Hof and the first author in \cite{HofN}. The higher dimensional case, however, is yet to be studied.


\begin{thebibliography}{30}

\bibitem{AgArr} \'Agoston, T.: Lattice cohomology and subspace arrangements: the topological and analytic cases. \emph{Period Math Hung} {\bf 89} (2024), 298–-312.

\bibitem{AgNe1} \'Agoston, T. and N\'emethi, A.: Analytic lattice cohomology of surface singularities, arXiv:2108.12294 (2021). 

\bibitem{AgNeIII} \'Agoston, T. and N\'emethi, A.: Analytic lattice cohomology of surface singularities, II (the equivariant case), arXiv:2108.12429 (2021).

\bibitem{AgNeHigh} \'Agoston, T. and N\'emethi, A.: The analytic lattice cohomology of isolated singularities,
arXiv:2109.11266
(2021).

\bibitem{AgNeCurves} \'Agoston, T. and N\'emethi, A.: The analytic lattice cohomology of isolated curve singularities,\\ arXiv:2301.08981 (2023).

\bibitem{AgNePoset} \'Agoston, T. and N\'emethi, A.: Lattice Cohomology of Partially Ordered Sets, \emph{Studia Sci. Math. Hung.} {\bf 61} (2) (2024), 185--202.

\bibitem{AltKlei} Altman, A. and Kleiman, S.: Introduction to Grothendieck duality theory, \emph{Lecture Notes in Mathematics} {\bf 146}, Springer (2006).
	

\bibitem{ArtalJordan} Artal  Bartolo, E.: Forme de Jordan de la monodrmomie  des
singularit\'es superisol\'ees de surfaces,
{\em Mem. Amer. Math. Soc.} {\bf 525} (1994).

\bibitem{ALM} Artal Bartolo, E., Luengo, I. and Melle Hern\'andez, A. :
Superisolated Surface Singularities, \emph{Singularities and Computer Algebra} (Lossen C., Pfister G., eds.), London Mathematical Society Lecture Note Series, Cambridge University Press (2006), 13--40.

\bibitem{Artin62} Artin, M.:
Some numerical criteria for contractability of curves on algebraic surfaces.
{\em  Amer. J. of Math.}, {\bf 84} (1962), 485--496. 

\bibitem{Artin66} Artin, M.:
On isolated rational singularities of surfaces.
{\em Amer. J. of Math.}, {\bf 88} (1966), 129--136.

\bibitem{AM} Atiyah, M.F. and MacDonald, I.G.: Introduction to Commutative Algebra, \emph{Addison-Wesley Series in Math.}, CRC Press (1969).

\bibitem{Bass} Bass, H.:
		On the ubiquity of Gorenstein rings,
		{\it Mathematische Zeitschrift},
		{\bf 82}
		(1) (1963),
		8--28.

\bibitem{BLZ} Borodzik, M., Liu, B. and Zemke, I.: Lattice homology, formality, and plumbed L-space links. \emph{J. Eur. Math. Soc.} (2024), published online first. 

\bibitem{BrKn} Brieskorn, E. and   Kn\"orrer, H.:  Plane Algebraic Curves,
{\it Birkh\"auser}, Boston, (1986).

\bibitem{Bruns-Herzog} Bruns, W. and Herzog, J.: Cohen Macaulay Rings, \emph{Cambridge studies in advanced mathematics} {\bf 39}, Cambridge University Press (1993).

\bibitem{BG80} Buchweitz, R.-O. and Greuel, G.-M.: The Milnor Number and Deformations of Complex Curve
Singularities, {\it Inventiones Math.} {\bf 58} (1980), 241--281.

\bibitem{arcFloer} Budur, N., de Bobadilla, J.F., Lê, Q.T. and Nguyen, H.D.: Cohomology of contact loci, \emph{J. Differential Geom.} {\bf 120} (2022), 389--409.

 \bibitem{cdg2} Campillo, A., Delgado de la Mata, F. and Gusein-Zade, S. M.: Integrals with respect to the Euler characteristic over spaces of functions and the Alexander polynomial, Proc. of the Steklov Inst. of Math. {\bf 238} (2002), 134--147.
    
    \bibitem{cdg} Campillo, A., Delgado de la Mata, F. and Gusein-Zade, S. M.: The Alexander polynomial of a plane curve singularity via the ring of functions on it, Duke Math. Journal {\bf 117} (1)  (2003), 125--156.

\bibitem{CNR} Conca, A., De Negri, E. and Rossi, M.E. Integrally closed and componentwise linear ideals. \emph{Math. Z.} {\bf 265} (2010), 715–-734.

\bibitem{CNP} Cassou-Nogues, P. and Płoski, A.: Invariants of plane curve singularities and Newton diagrams, {\it Universitatis Iagellonicae Acta Mathematica},  {\bf 49} (2011), 9--34.

\bibitem{CutVal} Cutkosky, S.D.: Valuations in Algebra and Geometry, \emph{Contemporary Mathematics}, {\bf 266} (2000), 45--64.

\bibitem{CHR} Cutkosky, S.D., Herzog, J. and Reguera, A.: Poincar\'e series of resolutions
of surface singularities,  \emph{Trans. of AMS}, {\bf 356} (5) (2003), 1833--1874.

\bibitem{DadNem} D\u{a}d\u{a}rlat, M. and  N\'emethi, A.:
  Shape theory and (connective) K-theory,
\emph{J. Operator Theory}, {\bf 23} (2) (1990), 207--291. 

\bibitem{good} D'Anna, M., Garc\'\i a-S\'anchez, P. A., Micale, V. and  Tozzo, L.:
Good subsemigroups of $ \mathbb{N}^n$,
{\it Internat. J. Algebra Comput.} {\bf 28} (2) (2018), 179–-206.

\bibitem{D'Anna} D'Anna, M.:
			The canonical module of a one-dimensional reduced local ring,
	{\it Communications in Algebra},
		{\bf 25} (9) (1997),
	2939--2965.

\bibitem{MinkSum} Das, S., and Sarvottamananda, S.:  Computing the Minkowski sum of convex polytopes in $\Re^ d$.  arXiv:1811.05812 (2018).

\bibitem{Spany}  de Bobadilla, J.F., Luengo, I., Melle-Hern\'andez, A. and
N\'emethi, A.:  On rational cuspidal projective plane curves,
\emph{Proc. London Math. Soc.} {\bf 92} (2006), 99--138.


\bibitem{BLMN2}
de Bobadilla, J.F., Luengo, I., Melle-Hern\'andez, A. and
N\'emethi, A.: On rational cuspidal curves, open surfaces and
local singularities,
{\em Singularity theory,
Dedicated to Jean-Paul Brasselet on His 60th Birthday,}
Proc.~of the 2005 Marseille Singularity School and Conference
(2007), 411--442.

\bibitem{delaMata87} Delgado de la Mata, F.:
The semigroup of values of a curve singularity with several branches, \emph{Manuscripta Math.} {\bf 59} (1987), 347--374.

\bibitem{delaMata} Delgado de la Mata, F.:
			Gorenstein curves and symmetry of the semigroup of values,
		{\it Manuscripta Mathematica}
		 {\bf 61}
		(3) (1988),
		285--296.

\bibitem{dJvS} de Jong, T. and van Straten, D.: Deformation theory of sandwiched singularities. \emph{Duke Math. Journal}, {\bf 95} (3) (1998), 451–-522.

\bibitem{Dimca} Dimca, A.: Singularities and Topology of Hypersurfaces, \emph{Universitext}, Springer-Verlag (1992).

\bibitem{DoldThom} Dold, A. and  Thom, R. :
 Quasifaserungen und unendliche symmetrische Produkte,
 \emph{Annals of Mathematics, Second Series}, {\bf 67} (1958),  239--281.

\bibitem{Durfeeform} Durfee, A.H.: The signature of smoothings of complex surface singularities,
  \emph{Math. Ann.} {\bf 232} (1) (1978), 85--98.

  \bibitem{DuVal} Du Val, P.: On isolated singularities of surfaces which do not affect
the condition of adjunction, \emph{Proc. Camb. Phil. Soc.} {\bf 30} (4) (1934),  453--459. 

\bibitem{EN} Eisenbud, D. and Neumann, W.: Three-dimensional link theory and invariants of plane curve singularities.
{\em Annals of Mathematics Studies}, Princeton Univ. Press (1985).


\bibitem{FJ} Favre, Ch. and Jonsson, M.: The valuative tree,
\emph{Lecture Notes in Mathematics} {\bf 1853}, Springer (2004).

\bibitem{Fukuda} Fukuda, K: Lecture notes: Polyhedral Computation, \url{https://people.inf.ethz.ch/fukudak/lect/pclect/notes2015/PolyComp2015.pdf}

\bibitem{Gaff} Gaffney, T.: Integral closure of modules and Whitney equisingularity, 
\emph{Invent. math.}, {\bf 107} (1992), 301--322.

\bibitem{GaffKl} Gaffney, T. and Kleiman S.L.: Specialization of integral dependence for  modules, 
\emph{Invent. math.} {\bf 137} (1999), 541--574.

\bibitem{Garcia} Garcia, A.:
Semigroups associated to singular points of plane curves,
\emph{J. für die reine und angewandte Math.} {\bf 336} (1982),
165--184.



\bibitem{GorNem2015}  Gorsky, E. and N\'emethi, A.:
Lattice and Heegaard Floer Homologies of Algebraic Links,
\emph{Int. Math. Res. Notices} {\bf 23} (2015), 12737--12780.

\bibitem{GrRie} Grauert, H. and Riemenschneider, O.: Verschwindungss\"atze f\"ur analytische
kohomologiegruppen auf komplexen R\"aumen, \emph{ Inventiones math.} {\bf 11} (1970), 263--292.

\bibitem{greuel} Greuel, G.-M.: 
Dualit\"at in der lokalen Kohomologie isolierter Singularit\"aten, \emph{Math Annalen} 
{\bf 250} (1980), 157--173.

\bibitem{Griffiths} Griffiths, Ph. A.: Variations on a Theorem of Abel, \emph{Invent. Math.}
{\bf 35} (1976), 321--390.

\bibitem{Hada} Hada, Y.: Optimal algebraic tangent cone of torsion-free sheaves via valuations, arXiv:2602.01112 (2026).

\bibitem{KodHay1} Hayasaka, F. and Kodiyalam, V.: Note on indecomposable integrally closed modules of rank 2 over two-dimensional regular local rings. \emph{Journal of Commutative Algebra} \textbf{15} (4) (2023), 513--518.

\bibitem{KodHay2} Hayasaka, F. and Kodiyalam, V.: Indecomposable integrally closed modules of rank 3 over two-dimensional regular local rings. \emph{Journal of Pure and Applied Algebra} \textbf{228} (6) (2024), 107612.

\bibitem{Hir} Hironaka, H.: Resolution of singularities of an algebraic variety over a field of characteristic zero, \emph{Annals of Mathematics} {\bf 79} (1964), 109--326.

\bibitem{Hir1} Hirzebruch, F.:  \"Uber vierdimensionale Riemannsche
Fla\"chen mehrdeutiger analytischer Functionen von zwei complexen
Ver\"anderlichen,  {\it Math. Ann.} {\bf 126} (1953), 1-22.


\bibitem{HofN} Hof, A. and Némethi, A.: Cohen-Macaulay Type via Lattice Homology and the Motivic Poincaré Series, arXiv:2509.11858 (2025).

\bibitem{Huneke} Huneke, C. L.:
		Hyman Bass and ubiquity: Gorenstein rings, \emph{Algebra, K-theory, Groups, and Education (New York, 1997)}, Contemp. Math. {\bf 243}, Amer. Math. Soc. (1999). 

\bibitem{HS} H\"ubl, R. and Swanson, I.: Adjoints of ideals, 
\emph{Michigan Math. J.} {\bf 57} (2008), 447--462. 

\bibitem{valued-VS} Irwan, S. E., Garminia, H. and Astuti, P.: On valuations of vector spaces. \emph{JP Journal of Algebra, Number Theory and Applications}, {\bf 38} (6) (2016), 633--646.

\bibitem{Ishii} Ishii, S.:  On Isolated Gorenstein Singularities,
{\it Math. Ann.} {\bf 270} (1985), 541--554.

\bibitem{IshiiArc} Ishii, S.: Jet schemes, arc spaces and the Nash problem, \emph{C. R. Math. Rep. Acad. Sci. Canada}, {\bf 29} (1) (2007), 1-–21.

\bibitem{Kap} Kaplansky, I.: Maximal fields with valuations I, \emph{Duke Mathematical Journal} {\bf 9} (1942), 303--321. 

\bibitem{Kho} Khovanov, M.:  A categorification of the Jones polynomial,
{\it Duke Mathematical Journal} {\bf 101} (3) (2000),  359–426.

\bibitem{Kno73} Kn\"oller, F.-W.:
2-dimensionale Singularit\"aten und Differentialformen,
{\it Mathematische Annalen}, {\bf 206} (3) (1973), 205--213.

\bibitem{Kodiyalam} Kodiyalam, V.: Integrally closed modules over two-dimensional regular local rings. \emph{Transactions of the American Mathematical Society} \textbf{347} (9) (1995), 3551--3573.

\bibitem{Kollar} Koll\'ar, J.:  Lectures on resolution of singularities. \emph{Annals of Mathematics Studies} {\bf 166}, Princeton University Press (2009).

\bibitem{KollarMori} Koll\'ar, J. and Mori, S.: Birational geometry of algebraic varieties, 
\emph{Cambridge Tracts in Mathematics}  {\bf 134}, Cambridge University Press (1998). 

\bibitem{Kou} Kouchnirenko, A. G.: Poly\`edres de Newton et nombres de Milnor,
\emph{Invent. Math.} {\bf 32} (1976), 1--31.

\bibitem{KNS1} Kubasch, A. A., N\'emethi, A. and Schefler, G.: Multiplicity and lattice cohomology of plane curve singularities, 
{\it Rev. Roumaine Math. Pures Appl.} {\bf 69} (2) (2024), 191–-234.
    
\bibitem{KNS2} Kubasch, A. A., N\'emethi, A. and Schefler, G.: Structural properties of the lattice cohomology of curve singularities, \textit{Selecta Mathematica, New Series} {\bf 31} (78) (2025).


\bibitem{KS3} Kubasch, A. A. and Schefler, G.: Lattice cohomology and the embedded topological type of plane curve singularities, 	arXiv:2504.13366 (2025).

\bibitem{LN1} L\'aszl\'o, T. and N\'emethi, A.:
Reduction theorem for lattice cohomology,
{\it Int. Math. Res. Notices } {\bf 11} (2015), 2938--2985.

\bibitem{Laufer72} Laufer, H.B.: On rational singularities,
{\em Amer. J. of Math.} {\bf 94} (1972), 597--608.

\bibitem{Lauferform} Laufer, H.B.: On $\mu$ for surface singularities, {\it Several Complex Variables, Part 1} (R. O.,  Wells Jr., ed.):  Proceedings of Symposia in Pure Mathematics {\bf 30} (1), American Mathematical Society, Providence (1977), 45--49.

\bibitem{LauferWeak} Laufer, H.B.: Weak simultaneous resolution for deformations of Gorenstein surface singularities,
{\it Proc. of Symposia in Pure Math.} {\bf 40} (2) (1983), 1--29.

\bibitem{LT} Lejeune-Jalabert, M. and Teissier, B.: Cl\^{o}ture int\'egrale des id\'eaux et \'euisingulait\'e,
\emph{Annales de la Facult\'e des sciences de Toulouse : Math\'ematiques}, serie 6, {\bf 17} (4) (2008), 781--859. 

\bibitem{Lipman} Lipman, J.: Rational singularities, with applications to algebraic surfaces
and unique factorization, {\it Inst. Hautes \'Etudes Sci. Publ. Math.} {\bf 36} (1969), 195-279.

\bibitem{LipmanComplete} Lipman, J.: On complete ideals in regular local rings, \emph{Algebraic geometry and commutative algebra,} Academic Press (1988), 203--231.

\bibitem{LipmanAdj} Lipman, J.: Adjoints of ideals in regular local rings, \emph{Math. Research letters} {\bf 1} (1994), 739--755.


\bibitem{Ignacio} Luengo, I.: The $\mu$-constant stratum is not
smooth, {\em Invent. Math.} {\bf 90} (1) (1987), 139--152.

\bibitem{LMNsi} Luengo, I.,  Melle-Hern\'andez, A. and  N\'emethi, A.:
Links and analytic invariants of superisolated singularities,
 {\em Journal of Algebraic Geometry} {\bf 14} (2005), 543--565.

\bibitem{MT} Merle, M. and Teissier, B.: Conditions d'adjonction, d'apr\`es Du Val, 
S\'eminaire sur les singularit\'es des surfaces (\'Ecole Polytechnique) {\bf 14} (1976--1977), 1--16.

\bibitem{Mor83}
Morales, M.:
Calcul de quelques invariants des singularit\'es de surface normale, in: Weber, C. (ed.) {\it Noeuds, Tresses et Singularités}, 
Monographies de l'Enseignement Math\'ematique {\bf 31} (1983), 191--203.

 \bibitem{julioproj} Moyano-Fern\'andez,  J. J.: Poincar\'e series for plane curve singularities and their behaviour under projections, {\it Journal of Pure and Applied Algebra}  {\bf 219} (6) (2015), 2449--2462.

\bibitem{NCCover} N\'emethi, A.:  Resolution graphs of some surface singularities. I. Cyclic coverings, {\it Singularities in Algebraic and Analytic Geometry} eds.: Melles, C. G. and Michler, R. I., Contemporary Mathematics {\bf 266}, American Mathematical Society, Providence  (2000), 89--128.

\bibitem{NOSz} N\'emethi, A.: On the Ozsv\'ath--Szab\'o invariant of negative definite plumbed 3-manifolds,
{\it Geometry \& Topology} {\bf 9} (2) (2005), 991--1042.

\bibitem{NGr} N\'emethi, A.: Graded roots and singularities,
 Proc. {\em Advanced School and Workshop} on {\em Singularities in Geometry and
Topology} ICTP (Trieste, Italy), World Sci. Publ., Hackensack, NJ (2007), 394--463.

\bibitem{Nlattice} N\'emethi, A.: Lattice cohomology of normal surface
  singularities.   {\em Publ. Res. Inst. Math. Sci.}  {\bf 44} (2)
  (2008),  507--543.

\bibitem{NJEMS} N\'emethi, A.: The Seiberg--Witten invariants of negative definite plumbed 3-manifolds,
{\em J. Eur. Math. Soc.} {\bf 13} (2011), 959--974.

\bibitem{NBook} N\'emethi, A.: Normal Surface Singularities, {\em Ergebnisse der Mathematik und ihrer
Grenzgebiete. 3. Folge} {\bf 74}, Springer Nature (2022).


\bibitem{NFilt2}  N\'emethi, A.:
Filtered lattice homology of surface singularities,
arXiv:2307.16581 (2023).

\bibitem{NFilt}  N\'emethi, A.:
Filtered lattice homology of curve singularities,
arXiv:2306.13889 (2023).

\bibitem{NSplurig} Némethi, A. and Schefler, G.: Categorification of the plurigenera of Gorenstein normal surface singularities. \emph{Mathematische Zeitschrift} {\bf 307} (4) (2024), 68.

\bibitem{NSigNN} N\'emethi, A. and Sigur{\dh}sson, B.:
Local Newton nondegenerate Weil divisors in toric varieties, {\it Michigan Math. J. Advance Publication} (2026), 1--58.


\bibitem{NV} N\'emethi, A. and Veys, W.: Filtrations associated with singularities, 
arXiv:2405.10898 (2024).

\bibitem{Seppo} Niemi-Colvin, S.: Invariance and naturality of knot lattice homology and homotopy, arXiv:2202.08941v4 (2024).

\bibitem{Oka} Oka, M.: On the bifurcation of the multiplicity and topology of the Newton boundary, \textit{Journal of the Mathematical Society of Japan} {\bf 31} (3) (1979), 435--450.

\bibitem{OW1} Orlik, P. and Wagreich P.: Isolated singularities of algebraic surfaces with $\mathbb{C}^*$-action, \emph{Ann. Math.} {\bf 93} (2) (1971), 205--228.

\bibitem{OW2} Orlik, P. and Wagreich P.: Singularities of algebraic surfaces with $\mathbb{C}^*$-action, \emph{Math. Ann.} {\bf 193} (2) (1971), 121--135.

\bibitem{Okuma} Okuma, T.: Plurigenera of surface singularities, \emph{Nova Science Pub. Inc.} (2000).

\bibitem{OSSzSpectral}   Ozsváth, P. S., Stipsicz, A. I. and Szabó, Z.: A spectral sequence on lattice homology. \textit{Quantum Topol.} {\bf 5} (4) (2014), 487–-521.

\bibitem{OSSzKnots} Ozsváth, P. S., Stipsicz, A. I. and Szabó, Z.: Knots in lattice homology. \emph{Comment. Math. Helv.} {\bf 89} (4) (2014), 783-–818.

\bibitem{OSzF}  Ozsv\'ath, P. S. and  Szab\'o, Z.: Absolutely graded Floer homologies and intersection
forms for four-manifolds with boundary, \emph{Advances in Mathematics} {\bf 173} (2003),
179–-261.

\bibitem{OSz} Ozsv\'ath, P.S. and  Szab\'o, Z.:
 Holomorphic discs
and topological invariants for closed  three-manifolds, \emph{Ann. of
Math.} (2) {\bf 159} (3) (2004), 1027--1158.

\bibitem{OSz7}  Ozsv\'ath, P.S. and  Szab\'o, Z.: Holomorphic discs
and three-manifold invariants: properties and applications, \emph{Ann.
of Math.} (2) {\bf 159} (3) (2004), 1159--1245.

\bibitem{OSzLink} Ozsváth, P. S. and
Szabó, Z.: Holomorphic disks, link invariants and the multi-variable
Alexander polynomial. {\it Algebr. Geom. Topol.} {\bf 8} (2) (2008),
615--692.

\bibitem{Quillen} Quillen, D.: Higher algebraic K-theory I, in {\em Lecture Notes in Mathematics},
 Springer-Verlag, Berlin, {\bf 341} (1973), 85--147. 


\bibitem{Ram} Ramanujam, C.P.: Remarks on Kodaira vanishing theorem, {\it J. Indian Math. Soc.} {\bf 36}
(1972), 41--51.

\bibitem{Rees} Rees, D.: Valuations associated with ideals II., \emph{J. London Math. Soc.} {\bf 31} (1956), 221-–228.


\bibitem{ReesMod} Rees, D.: Reduction of modules, in \emph{Mathematical Proceedings of the Cambridge Philosophical Society}, Cambridge University Press, {\bf 101} (3) (1987), 431--449.

\bibitem{Reid} Reid, M.: Canonical 3-folds, \emph{Journées de géometrie algébrique d'Angers} (A. Beauville, ed.), Sijthoff and Noordhoff (1980), 273--310.

\bibitem{MR}  Reid, M.: Chapters on Algebraic Surfaces,
in: \emph{Complex Algebraic Geometry},
IAS/Park City Mathematical Series,  {\bf 3}  (J. Koll\'ar, ed.) (1997),
3--159.

\bibitem{Rosenlicht} Rosenlicht, M.: Equivalence Relations on Algebraic Curves, \emph{Annals of Mathematics} {\bf 56} (1) (1952), 169--191.

\bibitem{Rotman} Rotman, J. J.: An Introduction to Homological Algebra, Second edition, \emph{Springer Universitext}, New York (2009). 

\bibitem{Sakai} Sakai, F.: Kodaira dimension of complements of divisors, \emph{Comples Analysis and Algebraic Geometry} (W. L. Baily and T. Shioda, ed.), Iwanami Shoten, Cambridge University Press (1977),  239--257.

\bibitem{SerreGor} Serre, J.--P.: Sur les modules projectifs, \emph{Séminaire Albert Châtelet et Paul Dubreil}, {\bf 14} (1)  (1960-1961).

\bibitem{Serre} Serre, J.--P.: Algebraic Groups and Class Fields,
  {\em Graduate Texts in Mathematics} {\bf 117}, Springer New York, NY (1988).

\bibitem{Stein} Stein, K.: Analytische Zerlegungen komplexer Räume, \emph{Math. Ann.} {\bf 132} (1) (1956), 63--93. 

\bibitem{StrStee} van Straten, D. and Steenbrink, J.: Extendibility of holomorphic differential forms near isolated hypersurface singularities, 
\emph{Abh. Math. Sem. Univ. Hamburg} {\bf 55} (1985), 97--110.

\bibitem{S} Swanson, I.: Rees valuations,  
in: Fontana, M., Kabbaj, SE., Olberding, B., Swanson, I. (eds) \emph{Commutative Algebra}, Springer, New York, NY (2010), 421--440.

\bibitem{SH} Swanson, I., Huneke, C.: Integral Closure of Ideals, Rings, and Modules,
\emph{London Mathematical Society Lecture Note Series} {\bf 336}, Cambridge University Press (2006).

\bibitem{Teissiersimultan} Teissier, B.: Résolution simultanée, I et II, \emph{Séminaire sur les singularités
des surfaces}, Palaiseau 1976--77, Demazure, H. Pinkham, B. Teissier, Editeurs, Springer \emph{Lecture Notes in Math.} {\bf 777} (1976), 71--146.

\bibitem{Teissier} Teissier, B.:
Mon\^{o}mes, volumes et multiplicit\'es, in: \emph{Introduction \`a la th\'eorie des singularités II}, Travaux en Cours {\bf 37} (1998), 127-–141.

\bibitem{TeissierKaplansky} Teissier, B.: On Kaplansky's Embedding Theorem, 	arXiv:2509.06379 (2025).


\bibitem{Var} Varchenko, A. N.: Zeta-function of monodromy and Newton's diagram,
\emph{Invent. Math.} {\bf 37} (1976), 253--262.

\bibitem{Watanabe80} Watanabe, K.:
On plurigenera of normal isolated singularities I,
{\it Mathematische Annalen} {\bf 250} (1) (1980), 65--94.

\bibitem{wahl} Wahl, J.: A characterization of quasi-homogeneous Gorenstein surface singularities, \emph{Composition Math.} {\bf 55} (3) (1985), 269--288.

\bibitem{wall} Wall, C.T.C.: Singular Points of Plane Curves,
{\em London Math. Soc., Student Texts} {\bf 63}, Cambridge University Press (2004).
   
\bibitem{Yamamoto} Yamamoto, M.: Classification of isolated algebraic singularities by their Alexander polynomials, Topology {\bf 23} (3), 277--287 (1984).
    
\bibitem{Yau77} Yau, S. S.-T.: Two theorems on higher dimensional singularities, \emph{Math. Ann.}
{\bf 231} (1977), 55--59.

\bibitem{Yau0} Yau, S. S.-T.: Various numerical invariants for isolated singularities, \emph{American Journal of Mathematics} {\bf 104} (5) (1982), 1063--1100.

\bibitem{Yau1} Yau, S. S.-T.: On irregularity and geometric genus of isolated singularities, 
\emph{Proc. Symp. Pure Math.} {\bf 40} (2) (1983), 653--662.

\bibitem{Yau2} Yau, S. S.-T.:
$s^{n-1}$  invariant for isolated $n$ dimensional singularities and
 its application to moduli problem, \emph{Amer. J. Math.} {\bf 104} (4) (1982), 829--841.

\bibitem{ZarComplete} Zariski, O.: Polynomial ideals defined by infinitely near base points, \emph{Amer. J. Math.} \textbf{60} (1938), 151--204.

\bibitem{Zemke} Zemke, I.: The equivalence of lattice and Heegaard Floer homology, \emph{Duke Math. J.} {\bf 174} (5) (2025), 857--910. 

\end{thebibliography}
\end{document}